\documentclass[12pt,twoside,leqno,openany]{amsbook}
\usepackage{amssymb,amsbsy,amsmath,amsfonts,amssymb,amscd,times,
graphics,color,xy,xypic,footmisc,fancyhdr,multicol,fancybox,graphicx,mathrsfs}
\usepackage[french]{babel}
\usepackage[T1]{fontenc} 
\definecolor{blue}{cmyk}{1.,1.,0.,0.63}
\definecolor{red}{cmyk}{0.,1.,1.,0.63}
\definecolor{green}{cmyk}{1.,0.,1.,0.63}
\let\mathcal\mathscr

\newcommand{\C}{\mathbb{C}}\newcommand{\K}{\mathbb{K}}
\newcommand{\N}{\mathbb{N}}\renewcommand{\P}{\mathbb{P}}
\newcommand{\Q}{\mathbb{Q}}\newcommand{\R}{\mathbb{R}}

\newcommand{\blue}{\textcolor{blue}}
\newcommand{\green}{\textcolor{green}}
\newcommand{\red}{\textcolor{red}}

\makeatletter
\renewcommand{\@fnsymbol}[1]
{\ensuremath{\ifcase#1\or $*$\or $**$\or $***$\or $****$\or $*****$
\else\@ctrerr\fi}}
\makeatother

\newcommand{\cf}{{\em cf.}}

\newcommand{\deutsch}[1]{[{\sl \!\!\;#1}]}

\newcommand{\deutschplain}[1]{{\sl #1}}

\newcommand{\eg}{{\em e.g.}}

\newcommand{\emphasis}[1]{{\em #1}}

\newcommand{\etc}{{\em etc.} }

\newcommand{\etcv}{{\em etc.,} }

\newcommand{\ie}{{\em i.e.} }

\newcommand{\voir}{{\em voir}}

\newcommand{\HEAD}[2]{%
\pagestyle{fancy}
\fancyhead[RO]{\scriptsize\sf\thepage}
\fancyhead[CO]{{\scriptsize\sf #1}}
\fancyhead[LE]{\scriptsize\sf\thepage}
\fancyhead[CE]{{\scriptsize\sf #2}}
\fancyfoot{}}

\newcommand{\sectiondritterV}[1]{%
\medskip
\nopagebreak
\begin{center}
{\sf #1}
\end{center}
\nopagebreak
\smallskip}

\newcommand{\CITATION}[1]{\smallskip\hfill
\begin{minipage}[t]{11cm}\baselineskip=0.43cm\parindent=0.91cm
{\footnotesize{\sf #1}}\end{minipage}\medskip}
\newcommand{\REFERENCE}[1]{\hfill \green{#1}}

\def\boiteepaisseavecuntitre#1{%
  \def\thickhrulefill{\leavevmode \leaders \hrule height 1pt \hfill \kern \z@}%
  \def\bkvz@before@breakbox{\ifhmode\par\fi\vskip\breakboxskip\relax}%
  \fboxrule=1pt
  \def\bkvz@set@linewidth{\advance\linewidth -2\fboxrule
                          \advance\linewidth -2\fboxsep}%
  \def\bkvz@left{\vrule \@width\fboxrule\hskip\fboxsep}%
  \def\bkvz@right{\hskip\fboxsep\vrule \@width\fboxrule}%
  \def\bkvz@top{\hbox to \hsize{%
      \vrule\@width\fboxrule\@height 1.2pt %%% D'où vient ce 0.2pt ????
      \thickhrulefill{#1}\thickhrulefill
      \vrule\@width\fboxrule\@height 1.2pt}}%
  \def\bkvz@bottom{\hrule\@height\fboxrule}%
  \breakbox}

\newcommand{\engellie}{\smallskip\begin{boiteepaisseavecuntitre}{}
\normalsize\baselineskip=0.53cm\fboxrule=0.77pt}
\newcommand{\stopengellie}{\end{boiteepaisseavecuntitre}}

\begin{document}

\pagenumbering{roman}

\thispagestyle{empty}
%\pagenumbering{roman}

\begin{center}

{\large

Joël\,\, {\sc M\,e\,r\,k\,e\,r}

\smallskip

{\large\footnotesize\sf \'Ecole Normale Supérieure}

{\large\footnotesize\sf D\'epartement de Math\'ematiques et Applications}

{\large\footnotesize\sf 45 rue d'Ulm, F-75230 Paris Cedex 05}

{\large\footnotesize\bf www.dma.ens.fr/$\sim$merker/index.html}

{\large\footnotesize\sf merker@dma.ens.fr}

}

\bigskip

\bigskip\bigskip\bigskip\bigskip\bigskip

\bigskip\bigskip

{\LARGE\bf

Sophus Lie, Friedrich Engel

\bigskip

et le problème de Riemann-Helmholtz

}

\bigskip\bigskip
{\sf\em (Titre actuel)}

\end{center}

\vfill
{\large\sf Monographie}
\hfill
{\bf 5 février 2010}

{\sf (ANR Physique et Géométrie}
\hfill
{\bf viii+311 pages}

{\sf à la charnière 
19\textsuperscript{ème}/20\textsuperscript{ème})}

{\large\sf à paraître aux \'Editions Hermann}

{\large\sf Paris, rue de la Sorbonne}

\newpage

\thispagestyle{empty}

$\:$

\begin{center}

{\LARGE\bf
THEORIE}

\bigskip\bigskip

{\large\sc der}

\bigskip\bigskip

{\LARGE\bf
TRANSFORMATIONSGRUPPEN}

\bigskip\bigskip

{\bf ----------------}

\bigskip\bigskip

{\large
DRITTER UND LETZTER ABSCHNITT}

\bigskip\bigskip

{\bf ----------------}

\bigskip\bigskip

UNTER MITWIRKUNG

\bigskip\bigskip

{\sc von}

\bigskip\bigskip

{\large Prof. Dr. FRIEDRICH ENGEL}

\bigskip\bigskip

{\sc bearbeitet}

\bigskip\bigskip

{\sc von}

\bigskip\bigskip

{\bf SOPHUS LIE,}

\medskip

{\footnotesize\sc professor der geometrie and der
universität leipzig}

\vfill

{\large\bf 
B.G. TEUBNER, LEIPZIG}

\medskip

{\large\bf 
1893}

\medskip

{\bf ----------------}

\medskip

{\large\bf 
CHELSEA PUBLISHING COMPANY}

\medskip

{\large\bf 
NEW YORK, N. Y.}

\medskip

{\large\bf 
1970}

\end{center}

\newpage

\bigskip\bigskip
\bigskip\bigskip
\bigskip

\centerline{\Large\bf Table des matières}

\bigskip

\centerline{\sf\em (actuelle)}

\bigskip\bigskip\bigskip

\noindent
{\large\bf Partie~I}\ \ \ \ \ \ \ \ \ \ \ \ \
{\large\bf\em Introduction philosophique générale}

%%%%%%%%%%%%%%%%%%%%%%%%%  CHAPITRE 1 %%%%%%%%%%%%%%%%%%%%%%%%%%%%%%%

\medskip\noindent
\parbox[t]{12.75cm}{
\parbox[t]{1.5cm}{\hfill{\footnotesize\bf Chap.~1.}}
\parbox[t]{0.25cm}{}
\parbox[t]{10cm}{
\baselineskip=0.43cm
{\footnotesize\bf L'ouverture riemannienne\dotfill}}\hfill
\parbox[t]{1cm}{\scriptsize\sf\hfill\pageref{Chapitre-1}}}

\noindent
\parbox[t]{12.75cm}{
\parbox[t]{1.5cm}{\hfill{\footnotesize\bf }}
\parbox[t]{0.25cm}{}
\parbox[t]{10cm}{
\baselineskip=0.43cm
{\scriptsize\sf 1.1.~Circonstances historiques\dotfill}}\hfill
\parbox[t]{1cm}{
\scriptsize\sf\hfill\pageref{circonstances-historiques}}}

\noindent
\parbox[t]{12.75cm}{
\parbox[t]{1.5cm}{\hfill{\footnotesize\bf }}
\parbox[t]{0.25cm}{}
\parbox[t]{10cm}{
\baselineskip=0.43cm
{\scriptsize\sf 1.2.~Appréciations d'universalité\dotfill}}\hfill
\parbox[t]{1cm}{
\scriptsize\sf\hfill\pageref{appreciations-d-universalite}}}

\noindent
\parbox[t]{12.75cm}{
\parbox[t]{1.5cm}{\hfill{\footnotesize\bf }}
\parbox[t]{0.25cm}{}
\parbox[t]{10cm}{
\baselineskip=0.43cm
{\scriptsize\sf 1.3.~Assembler l'inachevé\dotfill}}\hfill
\parbox[t]{1cm}{
\scriptsize\sf\hfill\pageref{assembler-l-inacheve}}}

\noindent
\parbox[t]{12.75cm}{
\parbox[t]{1.5cm}{\hfill{\footnotesize\bf }}
\parbox[t]{0.25cm}{}
\parbox[t]{10cm}{
\baselineskip=0.43cm
{\scriptsize\sf 1.4.~Le myst\`ere des
notions primitives de la g\'eom\'etrie\dotfill}}\hfill
\parbox[t]{1cm}{
\scriptsize\sf\hfill\pageref{mystere-notions-primitives}}}

\noindent
\parbox[t]{12.75cm}{
\parbox[t]{1.5cm}{\hfill{\footnotesize\bf }}
\parbox[t]{0.25cm}{}
\parbox[t]{10cm}{
\baselineskip=0.43cm
{\scriptsize\sf 1.5.~Fondements de la g\'eom\'etrie\dotfill}}\hfill
\parbox[t]{1cm}{
\scriptsize\sf\hfill\pageref{fondements-de-la-geometrie}}}

\noindent
\parbox[t]{12.75cm}{
\parbox[t]{1.5cm}{\hfill{\footnotesize\bf }}
\parbox[t]{0.25cm}{}
\parbox[t]{10cm}{
\baselineskip=0.43cm
{\scriptsize\sf 1.6.~Le renversement riemannien\dotfill}}\hfill
\parbox[t]{1cm}{
\scriptsize\sf\hfill\pageref{le-renversement-riemannien}}}

\noindent
\parbox[t]{12.75cm}{
\parbox[t]{1.5cm}{\hfill{\footnotesize\bf }}
\parbox[t]{0.25cm}{}
\parbox[t]{10cm}{
\baselineskip=0.43cm
{\scriptsize\sf 1.7.~Décider 
l'ouverture problématique du conceptuel\dotfill}}\hfill
\parbox[t]{1cm}{
\scriptsize\sf\hfill\pageref{decider-l-ouverture}}}

\noindent
\parbox[t]{12.75cm}{
\parbox[t]{1.5cm}{\hfill{\footnotesize\bf }}
\parbox[t]{0.25cm}{}
\parbox[t]{10cm}{
\baselineskip=0.43cm
{\scriptsize\sf 1.8.~Nécessité, suffisance, bifurcation\dotfill}}\hfill
\parbox[t]{1cm}{
\scriptsize\sf\hfill\pageref{necessite-suffisance-bifurcation}}}

\noindent
\parbox[t]{12.75cm}{
\parbox[t]{1.5cm}{\hfill{\footnotesize\bf }}
\parbox[t]{0.25cm}{}
\parbox[t]{10cm}{
\baselineskip=0.43cm
{\scriptsize\sf 1.9.~L'influence épistémologique de
Herbart sur Riemann\dotfill}}\hfill
\parbox[t]{1cm}{
\scriptsize\sf\hfill\pageref{influence-de-Herbart}}}

\noindent
\parbox[t]{12.75cm}{
\parbox[t]{1.5cm}{\hfill{\footnotesize\bf }}
\parbox[t]{0.25cm}{}
\parbox[t]{10cm}{
\baselineskip=0.43cm
{\scriptsize\sf 1.10.~Le réalisme 
dialectique modéré de Herbart\dotfill}}\hfill
\parbox[t]{1cm}{
\scriptsize\sf\hfill\pageref{realisme-dialectique-Herbart}}}

\noindent
\parbox[t]{12.75cm}{
\parbox[t]{1.5cm}{\hfill{\footnotesize\bf }}
\parbox[t]{0.25cm}{}
\parbox[t]{10cm}{
\baselineskip=0.43cm
{\scriptsize\sf 1.11.~La méthode des relations\dotfill}}\hfill
\parbox[t]{1cm}{
\scriptsize\sf\hfill\pageref{methode-des-relations}}}

\noindent
\parbox[t]{12.75cm}{
\parbox[t]{1.5cm}{\hfill{\footnotesize\bf }}
\parbox[t]{0.25cm}{}
\parbox[t]{10cm}{
\baselineskip=0.43cm
{\scriptsize\sf 1.12.~La méthode de la spéculation et les
graphes de concepts\dotfill}}\hfill
\parbox[t]{1cm}{
\scriptsize\sf\hfill\pageref{methode-de-la-speculation}}}

\noindent
\parbox[t]{12.75cm}{
\parbox[t]{1.5cm}{\hfill{\footnotesize\bf }}
\parbox[t]{0.25cm}{}
\parbox[t]{10cm}{
\baselineskip=0.43cm
{\scriptsize\sf 1.13.~La méthode des plus petits
changements conceptuels\dotfill}}\hfill
\parbox[t]{1cm}{
\scriptsize\sf\hfill\pageref{methode-plus-petits-changements}}}

\noindent
\parbox[t]{12.75cm}{
\parbox[t]{1.5cm}{\hfill{\footnotesize\bf }}
\parbox[t]{0.25cm}{}
\parbox[t]{10cm}{
\baselineskip=0.43cm
{\scriptsize\sf 1.14.~Modes amétriques de détermination\dotfill}}\hfill
\parbox[t]{1cm}{
\scriptsize\sf\hfill\pageref{modes-ametriques}}}

\noindent
\parbox[t]{12.75cm}{
\parbox[t]{1.5cm}{\hfill{\footnotesize\bf }}
\parbox[t]{0.25cm}{}
\parbox[t]{10cm}{
\baselineskip=0.43cm
{\scriptsize\sf 1.15.~Genèse 
du multidimensionnel\dotfill}}\hfill
\parbox[t]{1cm}{
\scriptsize\sf\hfill\pageref{genese-du-multidimensionnel}}}

\noindent
\parbox[t]{12.75cm}{
\parbox[t]{1.5cm}{\hfill{\footnotesize\bf }}
\parbox[t]{0.25cm}{}
\parbox[t]{10cm}{
\baselineskip=0.43cm
{\scriptsize\sf 1.16.~Conditions pour la détermination des
rapports métriques\dotfill}}\hfill
\parbox[t]{1cm}{
\scriptsize\sf\hfill\pageref{conditions-suffisantes-metrique}}}

\noindent
\parbox[t]{12.75cm}{
\parbox[t]{1.5cm}{\hfill{\footnotesize\bf }}
\parbox[t]{0.25cm}{}
\parbox[t]{10cm}{
\baselineskip=0.43cm
{\scriptsize\sf 1.17.~Genèse des métriques riemanniennes\dotfill}}\hfill
\parbox[t]{1cm}{
\scriptsize\sf\hfill\pageref{genese-des-metriques-riemanniennes}}}

\noindent
\parbox[t]{12.75cm}{
\parbox[t]{1.5cm}{\hfill{\footnotesize\bf }}
\parbox[t]{0.25cm}{}
\parbox[t]{10cm}{
\baselineskip=0.43cm
{\scriptsize\sf 1.18.~Surfaces de courbure constante\dotfill}}\hfill
\parbox[t]{1cm}{
\scriptsize\sf\hfill\pageref{surfaces-de-courbure-constante}}}

\noindent
\parbox[t]{12.75cm}{
\parbox[t]{1.5cm}{\hfill{\footnotesize\bf }}
\parbox[t]{0.25cm}{}
\parbox[t]{10cm}{
\baselineskip=0.43cm
{\scriptsize\sf 1.19.~Courbure sectionnelle de 
Riemann-Christoffel-Lipschitz\dotfill}}\hfill
\parbox[t]{1cm}{
\scriptsize\sf\hfill\pageref{courbure-Riemann-Christoffel-Lipschitz}}}

\noindent
\parbox[t]{12.75cm}{
\parbox[t]{1.5cm}{\hfill{\footnotesize\bf }}
\parbox[t]{0.25cm}{}
\parbox[t]{10cm}{
\baselineskip=0.43cm
{\scriptsize\sf 1.20.~Caractérisation de l'euclidéanité 
par annulation de la
courbure\dotfill}}\hfill
\parbox[t]{1cm}{
\scriptsize\sf\hfill\pageref{caracterisation-par-courbure-nulle}}}

%%%%%%%%%%%%%%%%%%%%%%%%%  CHAPITRE 2 %%%%%%%%%%%%%%%%%%%%%%%%%%%%%%%

\medskip\noindent
\parbox[t]{12.75cm}{
\parbox[t]{1.5cm}{\hfill{\footnotesize\bf Chap.~2.}}
\parbox[t]{0.25cm}{}
\parbox[t]{10cm}{
\baselineskip=0.43cm
{\footnotesize\bf La mobilité helmholtzienne 
de la rigidité\dotfill}}\hfill
\parbox[t]{1cm}{\scriptsize\sf\hfill\pageref{Chapitre-2}}}

\noindent
\parbox[t]{12.75cm}{
\parbox[t]{1.5cm}{\hfill{\footnotesize\bf }}
\parbox[t]{0.25cm}{}
\parbox[t]{10cm}{
\baselineskip=0.43cm
{\scriptsize\sf 2.1.~Le problème de Riemann-Helmholtz\dotfill}}\hfill
\parbox[t]{1cm}{
\scriptsize\sf\hfill\pageref{probleme-riemann-helmholtz}}}

\noindent
\parbox[t]{12.75cm}{
\parbox[t]{1.5cm}{\hfill{\footnotesize\bf }}
\parbox[t]{0.25cm}{}
\parbox[t]{10cm}{
\baselineskip=0.43cm
{\scriptsize\sf 2.2.~Incomplétudes riemanniennes\dotfill}}\hfill
\parbox[t]{1cm}{
\scriptsize\sf\hfill\pageref{incompletudes-riemanniennes}}}

\noindent
\parbox[t]{12.75cm}{
\parbox[t]{1.5cm}{\hfill{\footnotesize\bf }}
\parbox[t]{0.25cm}{}
\parbox[t]{10cm}{
\baselineskip=0.43cm
{\scriptsize\sf 2.3.~Rendre objectives les propositions
de la géométrie\dotfill}}\hfill
\parbox[t]{1cm}{
\scriptsize\sf\hfill\pageref{objectiver-la-geometrie}}}

\noindent
\parbox[t]{12.75cm}{
\parbox[t]{1.5cm}{\hfill{\footnotesize\bf }}
\parbox[t]{0.25cm}{}
\parbox[t]{10cm}{
\baselineskip=0.43cm
{\scriptsize\sf 2.4.~Les quatre axiomes de Helmholtz\dotfill}}\hfill
\parbox[t]{1cm}{
\scriptsize\sf\hfill\pageref{quatre-axiomes-helmholtz}}}

\noindent
\parbox[t]{12.75cm}{
\parbox[t]{1.5cm}{\hfill{\footnotesize\bf }}
\parbox[t]{0.25cm}{}
\parbox[t]{10cm}{
\baselineskip=0.43cm
{\scriptsize\sf 2.5.~Linéarisation de l'isotropie\dotfill}}\hfill
\parbox[t]{1cm}{
\scriptsize\sf\hfill\pageref{linearisation-isotropie}}}

\noindent
\parbox[t]{12.75cm}{
\parbox[t]{1.5cm}{\hfill{\footnotesize\bf }}
\parbox[t]{0.25cm}{}
\parbox[t]{10cm}{
\baselineskip=0.43cm
{\scriptsize\sf 2.6.~Critique par Lie de l'erreur
principale de Helmholtz\dotfill}}\hfill
\parbox[t]{1cm}{
\scriptsize\sf\hfill\pageref{critique-Lie-Helmholtz}}}

\noindent
\parbox[t]{12.75cm}{
\parbox[t]{1.5cm}{\hfill{\footnotesize\bf }}
\parbox[t]{0.25cm}{}
\parbox[t]{10cm}{
\baselineskip=0.43cm
{\scriptsize\sf 2.7.~Calculs helmholtziens\dotfill}}\hfill
\parbox[t]{1cm}{
\scriptsize\sf\hfill\pageref{calculs-helmholtziens}}}

\noindent
\parbox[t]{12.75cm}{
\parbox[t]{1.5cm}{\hfill{\footnotesize\bf }}
\parbox[t]{0.25cm}{}
\parbox[t]{10cm}{
\baselineskip=0.43cm
{\scriptsize\sf 2.8.~Insuffisances et reprises\dotfill}}\hfill
\parbox[t]{1cm}{
\scriptsize\sf\hfill\pageref{insuffisances-et-reprises}}}

\noindent
\parbox[t]{12.75cm}{
\parbox[t]{1.5cm}{\hfill{\footnotesize\bf }}
\parbox[t]{0.25cm}{}
\parbox[t]{10cm}{
\baselineskip=0.43cm
{\scriptsize\sf 2.9.~L'approche infinitésimale systématique
de Engel et de Lie\dotfill}}\hfill
\parbox[t]{1cm}{
\scriptsize\sf\hfill\pageref{approche-infinitesimale-systematique}}}

%%%%%%%%%%%%%%%%%%%%%%%%%  CHAPITRE 3 %%%%%%%%%%%%%%%%%%%%%%%%%%%%%%%

%%%\medskip\noindent
%%%\parbox[t]{12.75cm}{
%%%\parbox[t]{1.5cm}{\hfill{\footnotesize\bf Chap.~3.}}
%%%\parbox[t]{0.25cm}{}
%%%\parbox[t]{10cm}{
%%%\baselineskip=0.43cm
%%%{\footnotesize\bf La méthode conceptuelle 
%%%synthétique de Lie\dotfill}}\hfill
%%%\parbox[t]{1cm}{\scriptsize\sf\hfill\pageref{Chapitre-3}}}

\bigskip
\centerline{\bf -----------------}
\bigskip

\noindent
{\large\bf Partie~II}\ \ \ \ \ \ \ \ \ \
{\large\bf\em Introduction mathématique à la théorie de Lie}

\medskip\noindent
\parbox[t]{12.75cm}{
\parbox[t]{1.5cm}{\hfill{\footnotesize\bf Prologue.}}
\parbox[t]{0.25cm}{}
\parbox[t]{10cm}{
\baselineskip=0.43cm
{\footnotesize\bf Trois principes de pensée qui 
gouvernent la théorie de Lie\dotfill}}\hfill
\parbox[t]{1cm}{\scriptsize\sf\hfill\pageref{Prologue}}}

%%%%%%%%%%%%%%%%%%%%%%%%%  CHAPITRE 4 %%%%%%%%%%%%%%%%%%%%%%%%%%%%%%%

\medskip\noindent
\parbox[t]{12.75cm}{
\parbox[t]{1.5cm}{\hfill{\footnotesize\bf Chap.~3.}}
\parbox[t]{0.25cm}{}
\parbox[t]{10cm}{
\baselineskip=0.43cm
{\footnotesize\bf Théorèmes fondamentaux sur les
groupes de transformations\dotfill}}\hfill
\parbox[t]{1cm}{\scriptsize\sf\hfill\pageref{Chapitre-3}}}

\noindent
\parbox[t]{12.75cm}{
\parbox[t]{1.5cm}{\hfill{\footnotesize\bf }}
\parbox[t]{0.25cm}{}
\parbox[t]{10cm}{
\baselineskip=0.43cm
{\scriptsize\sf 3.1.~Paramètres essentiels\dotfill}}\hfill
\parbox[t]{1cm}{\scriptsize\sf\hfill\pageref{parametres-essentiels}}}

\noindent
\parbox[t]{12.75cm}{
\parbox[t]{1.5cm}{\hfill{\footnotesize\bf }}
\parbox[t]{0.25cm}{}
\parbox[t]{10cm}{
\baselineskip=0.43cm
{\scriptsize\sf 3.2.~Concept de groupe de Lie local\dotfill}}\hfill
\parbox[t]{1cm}{\scriptsize\sf\hfill\pageref{concept-de-groupe-de-lie-local}}}

\noindent
\parbox[t]{12.75cm}{
\parbox[t]{1.5cm}{\hfill{\footnotesize\bf }}
\parbox[t]{0.25cm}{}
\parbox[t]{10cm}{
\baselineskip=0.43cm
{\scriptsize\sf 3.3.~Principe de raison suffisante
et axiome d'inverse\dotfill}}\hfill
\parbox[t]{1cm}{\scriptsize\sf\hfill\pageref{raison-suffisante-axiome-inverse}}}

\noindent
\parbox[t]{12.75cm}{
\parbox[t]{1.5cm}{\hfill{\footnotesize\bf }}
\parbox[t]{0.25cm}{}
\parbox[t]{10cm}{
\baselineskip=0.43cm
{\scriptsize\sf 3.4.~Introduction des 
transformations infinitésimales\dotfill}}\hfill
\parbox[t]{1cm}{\scriptsize\sf\hfill\pageref{introduction-transformations-infinitesimales}}}

\noindent
\parbox[t]{12.75cm}{
\parbox[t]{1.5cm}{\hfill{\footnotesize\bf }}
\parbox[t]{0.25cm}{}
\parbox[t]{10cm}{
\baselineskip=0.43cm
{\scriptsize\sf 3.5.~\'Equations différentielles
fondamentales\dotfill}}\hfill
\parbox[t]{1cm}{\scriptsize\sf\hfill\pageref{equations-differentielles-fondamentales}}}

\noindent
\parbox[t]{12.75cm}{
\parbox[t]{1.5cm}{\hfill{\footnotesize\bf }}
\parbox[t]{0.25cm}{}
\parbox[t]{10cm}{
\baselineskip=0.43cm
{\scriptsize\sf 3.6.~Champs de vecteurs
et groupes à un paramètre\dotfill}}\hfill
\parbox[t]{1cm}{\scriptsize\sf\hfill\pageref{champs-vecteurs-groupe-un-parametre}}}

\noindent
\parbox[t]{12.75cm}{
\parbox[t]{1.5cm}{\hfill{\footnotesize\bf }}
\parbox[t]{0.25cm}{}
\parbox[t]{10cm}{
\baselineskip=0.43cm
{\scriptsize\sf 3.7.~Le théorème de 
Clebsch--Lie-Frobenius\dotfill}}\hfill
\parbox[t]{1cm}{\scriptsize\sf\hfill\pageref{theoreme-Clebsch-Lie-Frobenius}}}

%\noindent
%\parbox[t]{12.75cm}{
%\parbox[t]{1.5cm}{\hfill{\footnotesize\bf }}
%\parbox[t]{0.25cm}{}
%\parbox[t]{10cm}{
%\baselineskip=0.43cm
%{\scriptsize\sf 3.8.~Systèmes de Pfaff et
%crochet de Lie\dotfill}}\hfill
%\parbox[t]{1cm}{\scriptsize\sf\hfill\pageref{pfaff-lie}}}

\noindent
\parbox[t]{12.75cm}{
\parbox[t]{1.5cm}{\hfill{\footnotesize\bf }}
\parbox[t]{0.25cm}{}
\parbox[t]{10cm}{
\baselineskip=0.43cm
{\scriptsize\sf 3.8.~Constantes de structure
et correspondance fondamentale\dotfill}}\hfill
\parbox[t]{1cm}{\scriptsize\sf\hfill\pageref{constantes-de-structure}}}

\noindent
\parbox[t]{12.75cm}{
\parbox[t]{1.5cm}{\hfill{\footnotesize\bf }}
\parbox[t]{0.25cm}{}
\parbox[t]{10cm}{
\baselineskip=0.43cm
{\scriptsize\sf 3.9.~Le problème de la classification
des groupes de transformations\dotfill}}\hfill
\parbox[t]{1cm}{\scriptsize\sf\hfill\pageref{probleme-classification}}}

%%%\noindent
%%%\parbox[t]{12.75cm}{
%%%\parbox[t]{1.5cm}{\hfill{\footnotesize\bf }}
%%%\parbox[t]{0.25cm}{}
%%%\parbox[t]{10cm}{
%%%\baselineskip=0.43cm
%%%{\scriptsize\sf 3.10.~Algébrisation de la genèse\dotfill}}\hfill
%%%\parbox[t]{1cm}{\scriptsize\sf\hfill\pageref{algebrisation-genese}}}

%%%\medskip\noindent
%%%\parbox[t]{12.75cm}{
%%%\parbox[t]{1.5cm}{\hfill{\footnotesize\bf Chap.~5.}}
%%%\parbox[t]{0.25cm}{}
%%%\parbox[t]{10cm}{
%%%\baselineskip=0.43cm
%%%{\footnotesize\bf Théorèmes de classification: 
%%%groupes et sous-groupes\dotfill}}\hfill
%%%\parbox[t]{1cm}{\scriptsize\sf\hfill\pageref{Chapitre-5}}}

%%%\noindent
%%%\parbox[t]{12.75cm}{
%%%\parbox[t]{1.5cm}{\hfill{\footnotesize\bf }}
%%%\parbox[t]{0.25cm}{}
%%%\parbox[t]{10cm}{
%%%\baselineskip=0.43cm
%%%{\scriptsize\sf 5.1.~Détermination des groupes de
%%%la variété une fois étendue\dotfill}}\hfill
%%%\parbox[t]{1cm}{\scriptsize\sf\hfill\pageref{une-fois-etendue}}}

\bigskip
\centerline{\bf -----------------}
\bigskip

\noindent
{\large\bf Partie~III}\ \ \ \ \ \ \ \ \ \ \ \ \
{\large\bf\em Traduction française commentée et annotée}

\smallskip\noindent
{\sf Sophus {\sc Lie}, unter Mitwirkung von Friedrich {\sc Engel}}

\smallskip\noindent
{\sf\em Theorie der Transformationsgruppen, Dritter und letzter
Abschnitt, Abtheilung V\,\footnotemark}

%%%%%%%%%%%%%%%%%%%%%%%%%  DIVISION V %%%%%%%%%%%%%%%%%%%%%%%%%%%%%%%

\medskip\noindent
\parbox[t]{12.75cm}{
\parbox[t]{1.5cm}{\hfill{\footnotesize\bf Divis.~V.}}
\parbox[t]{0.25cm}{}
\parbox[t]{10cm}{
\baselineskip=0.43cm
{\footnotesize\bf Recherches sur les fondements
de la G\'eom\'etrie\dotfill}}\hfill
\parbox[t]{1cm}{\scriptsize\sf\hfill\pageref{Division-5}}}

%%%%%%%%%%%%%%%%%%%%%%%%%  CHAPITRE 20 %%%%%%%%%%%%%%%%%%%%%%%%%%%%%%

\medskip\noindent
\parbox[t]{12.75cm}{
\parbox[t]{1.5cm}{\hfill{\footnotesize\bf Chap.~20.}}
\parbox[t]{0.25cm}{}
\parbox[t]{10cm}{
\baselineskip=0.43cm
{\footnotesize\bf D\'etermination des groupes de $R_3$ relativement
auxquels les paires de points poss\`edent un, et un seul invariant,
tandis que $s > 2$ points n'ont pas d'invariant essentiel\dotfill}}\hfill
\parbox[t]{1cm}{\scriptsize\sf\hfill\pageref{Chapitre-20}}}

\vspace{0.14cm}\noindent
\parbox[t]{12.75cm}{
\parbox[t]{1.5cm}{\hfill{\footnotesize\bf }}
\parbox[t]{0.25cm}{}
\parbox[t]{10cm}{
\baselineskip=0.43cm
{\scriptsize\sf \S\,\,85.~Propri\'et\'es 
caract\'eristiques des groupes recherch\'es\dotfill}}\hfill
\parbox[t]{1cm}{
\scriptsize\sf\hfill\pageref{S-85}}}

\noindent
\parbox[t]{12.75cm}{
\parbox[t]{1.5cm}{\hfill{\footnotesize\bf }}
\parbox[t]{0.25cm}{}
\parbox[t]{10cm}{
\baselineskip=0.43cm
{\scriptsize\sf \S\,\,86.~Groupes primitifs 
parmi les groupes recherch\'es\dotfill}}\hfill
\parbox[t]{1cm}{
\scriptsize\sf\hfill\pageref{S-86}}}

\noindent
\parbox[t]{12.75cm}{
\parbox[t]{1.5cm}{\hfill{\footnotesize\bf }}
\parbox[t]{0.25cm}{}
\parbox[t]{10cm}{
\baselineskip=0.43cm
{\scriptsize\sf \S\,\,87.~Groupes imprimitifs 
parmi les groupes recherch\'es\dotfill}}\hfill
\parbox[t]{1cm}{
\scriptsize\sf\hfill\pageref{S-87}}}

\noindent
\parbox[t]{12.75cm}{
\parbox[t]{1.5cm}{\hfill{\footnotesize\bf }}
\parbox[t]{0.25cm}{}
\parbox[t]{10cm}{
\baselineskip=0.43cm
{\scriptsize\sf \S\,\,88\dotfill}}\hfill
\parbox[t]{1cm}{
\scriptsize\sf\hfill\pageref{S-88}}}

\noindent
\parbox[t]{12.75cm}{
\parbox[t]{1.5cm}{\hfill{\footnotesize\bf }}
\parbox[t]{0.25cm}{}
\parbox[t]{10cm}{
\baselineskip=0.43cm
{\scriptsize\sf \S\,\,89.~R\'esolution du probl\`eme 
pour les groupes r\'eels\dotfill}}\hfill
\parbox[t]{1cm}{
\scriptsize\sf\hfill\pageref{S-89}}}

\noindent
\parbox[t]{12.75cm}{
\parbox[t]{1.5cm}{\hfill{\footnotesize\bf }}
\parbox[t]{0.25cm}{}
\parbox[t]{10cm}{
\baselineskip=0.43cm
{\scriptsize\sf \S\,\,90\dotfill}}\hfill
\parbox[t]{1cm}{
\scriptsize\sf\hfill\pageref{S-90}}}

%%%%%%%%%%%%%%%%%%%%%%%%%  CHAPITRE 21 %%%%%%%%%%%%%%%%%%%%%%%%%%%%%%

\medskip\noindent
\parbox[t]{12.75cm}{
\parbox[t]{1.5cm}{\hfill{\footnotesize\bf Chap.~21.}}
\parbox[t]{0.25cm}{}
\parbox[t]{10cm}{
\baselineskip=0.43cm
{\footnotesize\bf Critique des recherches helmholtziennes\dotfill}}\hfill
\parbox[t]{1cm}{\scriptsize\sf\hfill\pageref{Chapitre-21}}}

\noindent
\parbox[t]{12.75cm}{
\parbox[t]{1.5cm}{\hfill{\footnotesize\bf }}
\parbox[t]{0.25cm}{}
\parbox[t]{10cm}{
\baselineskip=0.43cm
{\scriptsize\sf \S\,\,91.~Les axiomes helmholtziens\dotfill}}\hfill
\parbox[t]{1cm}{
\scriptsize\sf\hfill\pageref{S-91}}}

\noindent
\parbox[t]{12.75cm}{
\parbox[t]{1.5cm}{\hfill{\footnotesize\bf }}
\parbox[t]{0.25cm}{}
\parbox[t]{10cm}{
\baselineskip=0.43cm
{\scriptsize\sf \S\,\,92.~Cons\'equences 
des axiomes helmholtziens\dotfill}}\hfill
\parbox[t]{1cm}{
\scriptsize\sf\hfill\pageref{S-92}}}

\noindent
\parbox[t]{12.75cm}{
\parbox[t]{1.5cm}{\hfill{\footnotesize\bf }}
\parbox[t]{0.25cm}{}
\parbox[t]{10cm}{
\baselineskip=0.43cm
{\scriptsize\sf \S\,\,93.~Formulation des axiomes helmholtziens
en termes de th\'eorie des groupes\dotfill}}\hfill
\parbox[t]{1cm}{
\scriptsize\sf\hfill\pageref{S-93}}}

\vspace{0.14cm}\noindent
\parbox[t]{12.75cm}{
\parbox[t]{1.5cm}{\hfill{\footnotesize\bf }}
\parbox[t]{0.25cm}{}
\parbox[t]{10cm}{
\baselineskip=0.43cm
{\scriptsize\sf \S\,\,94.~Critique des conclusions 
que Monsieur de Helmholtz
tire de ses axiomes\dotfill}}\hfill
\parbox[t]{1cm}{
\scriptsize\sf\hfill\pageref{S-94}}}

\vspace{0.21cm}\noindent
\parbox[t]{12.75cm}{
\parbox[t]{1.5cm}{\hfill{\footnotesize\bf }}
\parbox[t]{0.25cm}{}
\parbox[t]{10cm}{
\baselineskip=0.43cm
{\scriptsize\sf \S\,\,95.~Consid\'erations se rattachant 
aux calculs helmholtziens\dotfill}}\hfill
\parbox[t]{1cm}{
\scriptsize\sf\hfill\pageref{S-95}}}

\noindent
\parbox[t]{12.75cm}{
\parbox[t]{1.5cm}{\hfill{\footnotesize\bf }}
\parbox[t]{0.25cm}{}
\parbox[t]{10cm}{
\baselineskip=0.43cm
{\scriptsize\sf \S\,\,96.~Quelles conclusions peut-on tirer 
des axiomes helmholtziens?\dotfill}}\hfill
\parbox[t]{1cm}{
\scriptsize\sf\hfill\pageref{S-96}}}

%%%%%%%%%%%%%%%%%%%%%%%%%  CHAPITRE 22 %%%%%%%%%%%%%%%%%%%%%%%%%%%%%%

\medskip\noindent
\parbox[t]{12.75cm}{
\parbox[t]{1.5cm}{\hfill{\footnotesize\bf Chap.~22.}}
\parbox[t]{0.25cm}{}
\parbox[t]{10cm}{
\baselineskip=0.43cm
{\footnotesize\bf Premi\`ere solution 
du probl\`eme de Riemann-Helmholtz\dotfill}}\hfill
\parbox[t]{1cm}{\scriptsize\sf\hfill\pageref{Chapitre-22}}}

\noindent
\parbox[t]{12.75cm}{
\parbox[t]{1.5cm}{\hfill{\footnotesize\bf }}
\parbox[t]{0.25cm}{}
\parbox[t]{10cm}{
\baselineskip=0.43cm
{\scriptsize\sf \S\,\,97\dotfill}}\hfill
\parbox[t]{1cm}{
\scriptsize\sf\hfill\pageref{S-97}}}

\noindent
\parbox[t]{12.75cm}{
\parbox[t]{1.5cm}{\hfill{\footnotesize\bf }}
\parbox[t]{0.25cm}{}
\parbox[t]{10cm}{
\baselineskip=0.43cm
{\scriptsize\sf \S\,\,98\dotfill}}\hfill
\parbox[t]{1cm}{
\scriptsize\sf\hfill\pageref{S-98}}}

\noindent
\parbox[t]{12.75cm}{
\parbox[t]{1.5cm}{\hfill{\footnotesize\bf }}
\parbox[t]{0.25cm}{}
\parbox[t]{10cm}{
\baselineskip=0.43cm
{\scriptsize\sf \S\,\,99\dotfill}}\hfill
\parbox[t]{1cm}{
\scriptsize\sf\hfill\pageref{S-99}}}

\noindent
\parbox[t]{12.75cm}{
\parbox[t]{1.5cm}{\hfill{\footnotesize\bf }}
\parbox[t]{0.25cm}{}
\parbox[t]{10cm}{
\baselineskip=0.43cm
{\scriptsize\sf \S\,\,100.~Sur le discours 
d'habilitation de Riemann\dotfill}}\hfill
\parbox[t]{1cm}{
\scriptsize\sf\hfill\pageref{S-100}}}

%%%%%%%%%%%%%%%%%%%%%%%%%  CHAPITRE 23 %%%%%%%%%%%%%%%%%%%%%%%%%%%%%%

\medskip\noindent
\parbox[t]{12.75cm}{
\parbox[t]{1.5cm}{\hfill{\footnotesize\bf Chap.~23.}}
\parbox[t]{0.25cm}{}
\parbox[t]{10cm}{
\baselineskip=0.43cm
{\footnotesize\bf Deuxi\`eme solution du probl\`eme 
de Riemann-Helmholtz\dotfill}}\hfill
\parbox[t]{1cm}{\scriptsize\sf\hfill\pageref{Chapitre-23}}}

\noindent
\parbox[t]{12.75cm}{
\parbox[t]{1.5cm}{\hfill{\footnotesize\bf }}
\parbox[t]{0.25cm}{}
\parbox[t]{10cm}{
\baselineskip=0.43cm
{\scriptsize\sf \S\,\,101\dotfill}}\hfill
\parbox[t]{1cm}{
\scriptsize\sf\hfill\pageref{S-101}}}

\noindent
\parbox[t]{12.75cm}{
\parbox[t]{1.5cm}{\hfill{\footnotesize\bf }}
\parbox[t]{0.25cm}{}
\parbox[t]{10cm}{
\baselineskip=0.43cm
{\scriptsize\sf \S\,\,102\dotfill}}\hfill
\parbox[t]{1cm}{
\scriptsize\sf\hfill\pageref{S-102}}}

\noindent
\parbox[t]{12.75cm}{
\parbox[t]{1.5cm}{\hfill{\footnotesize\bf }}
\parbox[t]{0.25cm}{}
\parbox[t]{10cm}{
\baselineskip=0.43cm
{\scriptsize\sf \S\,\,103\dotfill}}\hfill
\parbox[t]{1cm}{
\scriptsize\sf\hfill\pageref{S-103}}}

\bigskip
\centerline{\bf -----------------}
\bigskip

\noindent
{\large\bf Bibliographie}

\medskip\noindent
\parbox[t]{12.75cm}{
\parbox[t]{1.5cm}{\hfill{\footnotesize\bf }}
\parbox[t]{0.25cm}{}
\parbox[t]{10cm}{
\baselineskip=0.43cm
{{\footnotesize\bf Bibliographie}\dotfill}}\hfill
\parbox[t]{1cm}{
\scriptsize\sf\hfill\pageref{-Bibliographie}}}

%%%\noindent
%%%\parbox[t]{12.75cm}{
%%%\parbox[t]{1.5cm}{\hfill{\footnotesize\bf }}
%%%\parbox[t]{0.25cm}{}
%%%\parbox[t]{10cm}{
%%%\baselineskip=0.43cm
%%%{{\footnotesize\bf Index}\dotfill}}\hfill
%%%\parbox[t]{1cm}{
%%%\scriptsize\sf\hfill\pageref{-Index}}

\vfill

\footnotetext{{\sc Lie}, S.: {\em Theorie der transformationsgruppen.
Dritter und Letzter Abschnitt. Unter Mitwirkung von Prof. Dr.
Friedrich Engel, bearbeitet von Sophus Lie}, B.G. Teubner, Leipzig,
1893. Reprinted by Chelsea Publishing Co. (New York, N.Y., 1970).}

%%%%%%%%%%%%%%%%%%%%%%%%%%%%%%%%%%%%%%%%%%%%%%%%%%%%%%%%%%%%%%%%%%%%%

\newpage

\thispagestyle{empty}
\setcounter{footnote}{0}

$\:$

\bigskip\bigskip\medskip

\thispagestyle{empty}
\setcounter{footnote}{0}

$\:$

\bigskip\bigskip\bigskip

\begin{center}
{\Large\bf
Préface
}\end{center}
\label{Preface}

\HEAD{Pr\'eface}{Jo\"el Merker}

\bigskip\medskip

La troisi\`eme et derni\`ere partie de cet ouvrage propose une
traduction française annot\'ee de la solution math\'ematique
compl\`ete qu'ont donn\'ee Friedrich Engel et Sophus Lie dans le
troisi\`eme volume de la monumentale {\em Theorie der
Transformationsgruppen} (\cite{ enlie1893}, pp.~393--523) 
au probl\`eme dit <<\,{\sl de Riemann-Helmholtz}\,>>:

\smallskip 

{\small\sf\em Par quels axiomes spatiaux peut-on caract\'eriser les
trois g\'eom\'etries \`a courbure constante, euclidienne,
lobatchevskienne ou riemannienne, \`a l'exclusion de toute autre
g\'eom\'etrie}?

\smallskip
Notre premier objectif est de faire revivre, au sein de la philosophie
contemporaine de la g\'eom\'etrie, l'incroyable puissance de pens\'ee
et de conception de Lie qui parvint, entre 1873 et 1893, \`a
\'eriger sur des
milliers de pages une th\'eorie r\'ealisant pleinement le fameux {\em
Programme d'Erlangen} de Klein (1872, \cite{ klei1974}). Lie partait
de l'id\'ee seule qu'il devait exister, dans le domaine des
transformations continues entre \'equations diff\'erentielles, un
analogue m\'etaphysiquement complet
\`a la th\'eorie des groupes de substitutions
des racines d'une \'equation alg\'ebrique que Serret et Jordan
venaient de d\'evelopper, quelques d\'ecennies apr\`es la parution
posthume du m\'emoire g\'enial~\cite{ ga1897} d'\'Evariste Galois.

Mais surtout, notre objectif principal est de montrer, en termes de
{\sl semailles} au sens de Grothendieck et dans la lign\'ee de
quelques
\'ecoles math\'ematiques contemporaines qui se sont d\'evelopp\'ees
aux \'Etats-Unis avec notamment Olver~\cite{ ol1986, ol1995},
Gardner~\cite{ gard1989} et Bryant~\cite{ bcggg1991}\,\,---\,\,\cf~aussi le beau
livre du su\'edois Stormark~\cite{ stk2000}\,\,---, que notre \'epoque
peut maintenant accueillir une sorte de r\'esurrection des
math\'ematiques qui se sont d\'evelopp\'ees \`a la charni\`ere du
19\textsuperscript{\`eme} et du 20\textsuperscript{\`eme} si\`ecle,
notamment dans ses techniques, dans ses mani\`eres de penser, et dans
ses probl\'ematiques mêmes, toujours riches, ouvertes et
in\'epuisables. \`A vrai dire, bien que le pr\'esent ouvrage ne
constitue en rien un travail d'histoire des math\'ematiques, c'est la
lecture des travaux d'Hawkins~\cite{ h2001, ha2005} qui nous a fait
prendre conscience, {\em philosophiquement}, de l'importance, pour
l'architecture des math\'ematiques, de classifier {\em a priori} et
abstraitement les groupes de transformations, qu'ils soient discrets
(Jordan) ou continus (Lie).

En effet, comme Lie lui-même l'a vite compris en profondeur d\`es
l'hiver 1873--74, la classification abstraite et {\em a priori} des
groupes poss\`ede une importance capitale dans l'\'edifice
math\'ematique, en tant que tronc commun de la connaissance \`a partir
duquel pourront être directement irrigu\'es puis r\'esolus tous les
probl\`emes sp\'ecialis\'es qui impliquent une structure de groupe
sous-jacente: transformations de contact, prolongements aux espaces de
jets, invariants diff\'erentiels, sym\'etries des \'equations aux
d\'eriv\'ees partielles, transformations projectives, transformations
conformes, structures symplectiques, th\'eorie des invariants
alg\'ebriques, \etc Toutefois, dans la litt\'erature math\'ematiques
actuelle, on ne trouve aucun article, aucun livre, aucune source qui
restitue la th\'eorie originale de Lie dans sa syst\'ematicit\'e
intrins\`eque, et donc pour prendre connaissance des d\'emonstrations
d\'etaill\'ees des th\'eor\`emes de classification, il faut se
reporter aux {\oe}uvres originales. Par ailleurs, le tournant de
globalisation et d'alg\'ebrisation de la th\'eorie des groupes
initi\'e par Hermann Weyl et \'Elie Cartan dans les ann\'ees 1930 et
continu\'e par Chevalley, Chern et Ehresmann dans les ann\'ees 1950 a
fait progressivement porter l'accent sur les th\'eor\`emes de
classification des alg\`ebres de Lie semi-simples complexes (Killing) et
r\'eelles (\'Elie Cartan), avec la combinatoire aff\'erente des
syst\`emes de racines (Weyl) et des diagrammes (Dynkin) qui est
aujourd'hui centrale dans la th\'eorie des repr\'esentations.

Mais puisque la th\'eorie initiale de Lie se ramifiait surtout en
direction des \'equations aux d\'eriv\'ees partielles et
s'architecturait afin que s'y inscrivent toutes les structures de
groupes continus possibles, ce premier ouvrage de traduction
comment\'ee sera suivi dans un avenir proche par d'autres travaux
proprement math\'ematiques (\cf~le texte en pr\'eparation~\cite{
merk2009b}) qui transcriront dans un langage math\'ematique moderne
d'autres r\'esultats de Lie, tout particuli\`erement ceux qui sont
consacr\'es, dans le Tome~III de la {\em Theorie der
Transformationsgruppen}, \`a la classification compl\`ete des actions
analytiques locales sur un espace \`a trois dimensions; celle-ci ne
fut en fait compl\'et\'ee qu'en 1902 par Amaldi (\cite{ a1901, a1902})
pour les actions doublement imprimitives qui stabilisent \`a la fois
un feuilletage local par des surfaces et un feuilletage subordonn\'e
par des courbes.

Voici pour terminer une br\`eve description du contenu de ce livre.
Deux chapitres constituent l'{\em Introduction philosophique
g\'en\'erale} (Partie~I) qui est suivie 
d'une {\em Introduction math\'ematique \`a la th\'eorie de Lie}
(Partie~II), utile pour aborder
notre traduction de la Division~V du Tome~III 
de la {\em Theorie der Transformationsgruppen} (Partie~III). 

Dans le premier chapitre que nous consacrons \`a un commentaire
d\'etaill\'e et cibl\'e de la c\'el\`ebre leçon orale d'habilitation
\deutsch{Probevorlesung} de Riemann (1854, \cite{
riem1892}, non publi\'ee de son vivant), nous insistons, en nous
servant notamment des \'etudes de Erhard 
Scholz~\cite{ scho1980,
scho1982b}, sur la capacit\'e si particuli\`ere que Riemann avait
d'{\em ouvrir sans les fermer} les questions de conceptualisation
math\'ematique.
\`A la fin de l'ann\'ee 1868, 
Helmholtz obtint une copie des notes manuscrites que Schering avait
prises pendant la soutenance de Riemann, et il publie rapidement une
\'etude (\cite{helm1868b}), maintenant centrale pour l'histoire des
fondements de la g\'eom\'etrie, dans laquelle il pr\'etend pouvoir
d\'emontrer rigoureusement que la postulation\,\,---\,\,physiquement
\'evidente sinon m\'etaphysiquement n\'ecessaire\,\,---\,\,de corps
rigides et n\'eanmoins maximalement mobiles dans un espace abstrait
quelconque implique l'existence d'une m\'etrique riemannienne \`a
courbure constante: inversion, donc, des points de vue de Riemann,
Christoffel et Lipschitz pour qui les m\'etriques quadratiques
possibles sont innombrables, puisqu'elles d\'ependent du comportement
des invariants diff\'erentiels d\'eriv\'es de la courbure.  Toutefois,
les arguments embryonnaires, riches d'id\'ees nouvelles, de Helmholtz,
sont en grande partie erron\'es ou incomplets, et il fallut attendre
les travaux de Lie~\cite{lie1886, lie1890a, lie1891, lie1892a,
lie1892b} pour qu'ils puissent être inscrits dans une vaste th\'eorie.
La traduction que nous proposons ici aura le m\'erite, nous
l'esp\'erons, de porter au moins \`a la connaissance des commentateurs
francophones de Helmholtz la discussion critique par Lie de la
pertinence des axiomes helmholtziens, {\em voir} notamment
l'introduction à la Division~5 au début de la Partie~III,
p.~\pageref{Division-5}, et aussi le Chapitre~21
p.~\pageref{Chapitre-21}.  Enfin, notre Partie~II est consacr\'ee \`a
une pr\'esentation des th\'eor\`emes fondamentaux de la th\'eorie de
Lie qui se base exclusivement sur le tr\`es syst\'ematique
Volume~I~\cite{ enlie1888} de la {\em Theorie der
Transformationsgruppen}, sans faire appel aux manuels modernes.

Non germaniste, l'auteur a b\'en\'efici\'e d'aides ponctuelles et
utiles concernant quelques difficult\'es de traduction, notamment de
la part de Maud Moreillon (Besançon), d'Egmont Porten (Berlin) et de
Jean Ruppenthal (Wuppertal). Françoise Panigeon a r\'eguli\`erement
contrôl\'e la correction lexicale et grammaticale du texte dans son
ensemble, qui lui doit donc beaucoup.  Toute notre reconnaissance
s'adresse aussi à Erhard Scholz, qui a accepté d'examiner l'ouvrage et
dont les travaux de stature internationale nous ont été d'une grande
utilité.  Le projet initial a \'et\'e f\'econd\'e grâce \`a certaines
s\'eances du <<\,{\sl S\'eminaire Riemann}\,>>, organis\'e en
collaboration avec Jean-Jacques Szczeciniarz et Ivahn Smadja \`a
l'{\sl \'Ecole Normale Sup\'erieure}.  Enfin, Jean-Jacques
Szczeciniarz, toujours disponible et ouvert, a constamment soutenu la
r\'ealisation de ce projet au cours d'\'echanges philosophiques
profonds et pr\'ecieux.

%%%%%%%%%%%%%%%%%%%%%%%%%%%%%%%%%%%%%%%%%%%%%%%%%%%%%%%%%%%%%%%%%%%%%

\newpage
\pagenumbering{arabic}

\thispagestyle{empty}
\setcounter{footnote}{0}

$\:$

\bigskip\bigskip\medskip

\begin{center}
{\Large\bf
Partie~I:
\\
\medskip
Introduction philosophique générale
}\end{center}
\label{Partie-I}

\bigskip\medskip

\begin{center}
\begin{minipage}[t]{8.5cm}
\baselineskip =0.35cm
{\scriptsize

\centerline{\footnotesize\bf Chapitres}

\medskip

\smallskip

{\bf 1.~L'ouverture riemannienne 
\dotfill \pageref{Chapitre-1}.}

{\bf 2.~La mobilité helmholtzienne de la rigidité
\dotfill \pageref{Chapitre-2}.}

%%%{\bf 3.~La méthode conceptuelle synthétique de Lie
%%%\dotfill \pageref{Chapitre-3}.}

}\end{minipage}
\end{center}

\bigskip

\begin{center}
{\large\bf Chapitre~1:

\smallskip
L'ouverture riemannienne}
\label{Chapitre-1}
\end{center}

\smallskip
\hfill
{\sl La math\'ematique est la seule bonne m\'etaphysique.}

\hfill
Lord {\sc Kelvin}

\HEAD{Chapitre~1.\,\,\,\,L'ouverture riemannienne}{
1.1.\,\,\,Circonstances historiques}

\medskip\noindent{\bf 1.1.~Circonstances historiques.}
\label{circonstances-historiques}
Le 10 juin 1854, à l'occasion de ses épreuves d'admission à
la c\'el\`ebre {\sl Université Georges Auguste} 
de Göttingen, Bernhard Riemann (1826--1866), alors
âgé de vingt-huit ans, défend oralement son {\em
Habilitationsvortrag}, qu'il a intitulée:

\smallskip

\centerline{\em Sur les hypothèses qui servent de fondement
à la géométrie\footnotemark.}

\smallskip

\footnotetext{\,
%%%%%%%%%%%%%%%%%%%%%%%-------DEBUT--------%%%%%%%%%%%%%%%%%%%%%%%%%%%
\deutsch{\"Uber die Hypothesen, welche der 
Geometrie zu Grunde liegen.} Paru en 1867 à titre posthume dans le
tome~XIII des {\em Mémoires de la Société Royale des Sciences de
Göttingen}, ce texte a immédiatement inspiré de nombreux travaux
mathématiques, notamment chez Dedekind, Gehring, Clifford, Helmholtz,
Christoffel, Lipschitz, Beltrami, et d'autres. En 1898, il a été
traduit en français par J.~{\sc Hoüel}, {\em voir}
\cite{ riem1898}, pp.~280--299.
} %%%%%%%%%%%%%%%%%%%%%%%%-----FIN-----%%%%%%%%%%%%%%%%%%%%%%%%%%%%%%%

En Allemagne au dix-neuvième siècle, le diplôme d'habilitation est
formellement requis pour être en mesure d'obtenir le statut de {\em
Privatdozent}, c'est-\`a-dire de charg\'e de cours. Non rétribuées
par l'université, ces positions sont néanmoins prisées; le salaire
afférent y dépend de la libéralité des étudiants qui assistent
r\'eguli\`erement aux leçons.

Au début des années 1850, on compte à Göttingen une dizaine de
professeurs permanents disposant d'une chaire. Le 16 d\'ecembre
1851, Riemann avait soutenu en latin sa th\`ese de doctorat
\deutsch{Inauguraldissertation} consacr\'ee aux fonctions d'une
variable 
complexe\footnote{\,
%%%%%%%%%%%%%%%%%%%%%%%-------DEBUT--------%%%%%%%%%%%%%%%%%%%%%%%%%%%
\deutschplain{Grundlagen für eine allgemeine Theorie der 
Functionen einer
veränderlichen complexen Grö{\ss}en}, Théorie générale des fonctions
d'une grandeur variable complexe, {\em voir}~\cite{ bell1939,
laug1999, remm1993} pour ce chapitre qui ne sera pas abordé dans cet
ouvrage.
}, %%%%%%%%%%%%%%%%%%%%%%%%-----FIN-----%%%%%%%%%%%%%%%%%%%%%%%%%%%%%%%
et le jury d'experts qui avait été consulté pour autoriser la
soutenance était composé de sept <<\,grands\,>> professeurs: Gauss en
mathématiques, Mitscherlich en rhétorique, Haussmann en minéralogie,
Ritter en philosophie, Hoeck et Hermann en philologie classique, Waitz
en histoire et Weber en physique (\cite{ remm1993}). Aussi n'est-il
pas étonnant que, dans un tel contexte d'universalité des compétences
et de proximité des savoirs, Riemann ait dû prononcer, deux ans et
demi plus tard pour l'habilitation, sa fameuse conférence d'épreuve
\deutsch{Probevorlesung} devant un auditoire majoritairement composé
de non-mathématiciens, dans le cadre de ce qu'on appelle aujourd'hui
un {\sl Colloquium} s'adressant en principe à tous les protagonistes
de l'université.

Frustration d'historien des mathématiques: les archives de la
Facult\'e philosophique de Göttingen n'ont conserv\'e aucune trace
(\cite{ remm1993}) concernant les rapports, les membres du jury, 
les discussions, les questions, 
ni pour la th\`ese, ni pour l'habilitation de Riemann. D'apr\`es
Dedekind (\cite{ ded1892}, pp.~547--548), les seuls \'el\'ements
incontestables concernant les circonstances de l'habilitation qui ne
rel\`event pas du mythe post-mortem sont deux lettres de Riemann, la
premi\`ere \'ecrite le 28 d\'ecembre 1853 \`a son fr\`ere Wilhelm au
moment de fixer tout ce
qui est inh\'erent \`a la
soutenance\footnote{\
%%%%%%%%%%%%%%%%%%%%%%%-------DEBUT--------%%%%%%%%%%%%%%%%%%%%%%%%%%%
<<\,Mon travail progresse raisonnablement: au
d\'ebut du mois de d\'ecembre, j'ai remis mon m\'emoire d'habilitation
et je devais proposer \`a ce moment-là trois sujets pour l'\'epreuve
orale, parmi lesquels la facult\'e devait en choisir un. J'avais
pr\'epar\'e les deux premiers et j'esp\'erais que l'un d'entre eux
serait s\'electionn\'e: malheureusement, Gauss opta pour le
troisi\`eme, et je suis
\`a pr\'esent un peu press\'e par le temps, car je dois encore
le pr\'eparer.\,>>
}, %%%%%%%%%%%%%%%%%%%%%%%%-----FIN-----%%%%%%%%%%%%%%%%%%%%%%%%%%%%%%%
et la seconde \`a ce même frère, le 26 juin 1854, deux semaines
apr\`es la soutenance\footnote{\, %%%%%%%%%%%%%%%%%%%%%%%%%%%%%%%%%%%%%
<<\,J'ai loué pour l'été une maison avec un jardin et, grâce à cela,
ma santé ne m'a plus tourmenté. Ayant terminé, deux semaines après
Pâques, une étude dont je ne pouvais pas venir à bout, je me suis
enfin mis à ma conférence d'épreuve et je l'ai terminée vers
Pentecôte.\,>>
}. %%%%%%%%%%%%%%%%%%%%%%%%-----FIN-----%%%%%%%%%%%%%%%%%%%%%%%%%%%%%%%
Maigre documentation, quand on pense à toutes les spéculations qui ont
été suscitées dans l'imagination des géomètres depuis plus d'un siècle
et demi.

Comme les règles universitaires l'exigent encore aujourd'hui en
Allemagne, le candidat \`a l'habilitation se doit de proposer trois
sujets distincts, afin de montrer au mieux l'extension th\'ematique de
ses travaux de recherche. Ainsi Riemann propose-t-il \`a la fin de
l'ann\'ee 1853 un sujet en Analyse, un sujet en Algèbre et un sujet en
Géométrie, à savoir:

\smallskip\noindent{\bf 1)}\, 
un mémoire abouti sur
les séries trigonométriques qu'il avait d\'ej\`a achevé pendant
l'automne\footnote{\,
%%%%%%%%%%%%%%%%%%%%%%%-------DEBUT--------%%%%%%%%%%%%%%%%%%%%%%%%%%%
Ce m\'emoire
ne fut publi\'e \`a titre posthume
qu'en 1868, par les soins de Dedekind.
}; %%%%%%%%%%%%%%%%%%%%%%%%-----FIN-----%%%%%%%%%%%%%%%%%%%%%%%%%%%%%%%

\smallskip\noindent{\bf 2)}\, 
un travail sur les intersections entre deux courbes planes 
du second degr\'e
\footnote{\,
%%%%%%%%%%%%%%%%%%%%%%%-------DEBUT--------%%%%%%%%%%%%%%%%%%%%%%%%%%%%
\deutsch{\"Uber die Auflösung zweier Gleichungen zweiten Grades mit
zwei unbekannten Grössen}, texte absent des
\deutschplain{Gesammelte Werke}. 
}; %%%%%%%%%%%%%%%%%%%%%%%%-----FIN-----%%%%%%%%%%%%%%%%%%%%%%%%%%%%%%%

\smallskip\noindent{\bf 3)}\, 
une réflexion générale sur les fondements de la g\'eom\'etrie.

\smallskip\noindent
Dans sa courte biographie~\cite{ ded1892}, Dedekind a \'ecrit que
Gauss aurait s\'electionn\'e le troisi\`eme sujet en
dérogeant à la convention académique habituelle
de choisir le premier, parce qu'il \'etait
curieux de voir comment un si jeune math\'ematicien pourrait traiter
une question qui demande tant de maturit\'e 
scientifique\footnote{\,
%%%%%%%%%%%%%%%%%%%%%%%-------DEBUT--------%%%%%%%%%%%%%%%%%%%%%%%%%%%
Remmert \cite{ remm1993} a reproduit le rapport de Gauss, écrit en
calligraphie pré-sütterline, sur la dissertation inaugurale de Riemann
de 1851; Gauss y appréciait déjà
l'<<\,indépendance productive et louable\,>>
\deutsch{rümliche productive Selbstthätigkeit} de Riemann. 
}. %%%%%%%%%%%%%%%%%%%%%%%%-----FIN-----%%%%%%%%%%%%%%%%%%%%%%%%%%%%%%%
Mais d'apr\`es Laugwitz~\cite{ laug1999}, le premier sujet propos\'e
\'etait th\'ematiquement trop proche de la th\`ese de Riemann; le
contenu du second est probablement apparu assez \'evident \`a Gauss;
de telle sorte que les professeurs concern\'es, se fiant à l'expertise
de Gauss, ne pouvaient qu'être conduits à s\'electionner le
troisi\`eme sujet sur la g\'eom\'etrie.

Il y a des raisons de penser que Riemann, en prenant le risque de
proposer un travail qui n'était que potentiellement en gestation,
avait l'intention de s'exposer \`a la contrainte, dans un cadre
institutionnel, de r\'ediger ses id\'ees nouvelles qu'il jugeait
fondamentales, bien qu'inachevées. En vérité, sur le moment, Riemann
ne semble pas avoir \'et\'e tellement préoccupé par cette tâche
suppl\'ementaire,
\'etant donn\'e
qu'après avoir programmé une soutenance vers la fin de
l'été\footnote{\,
%%%%%%%%%%%%%%%%%%%%%%%-------DEBUT--------%%%%%%%%%%%%%%%%%%%%%%%%%%%
Au printemps 1853, Gauss se plaignait dans une
lettre à Alexander von Humboldt de douleurs à la poitrine et au gosier,
d'essoufflements, de palpitations et d'insomnie. Un an plus tard, son
état s'était aggravé. Le vendredi 9 juin 1854, il apprend que Riemann
a officiellement déposé son texte, et il fixe la conférence au
lendemain (\cite{ bier1990}, pp.~201--202). 
}, %%%%%%%%%%%%%%%%%%%%%%%%-----FIN-----%%%%%%%%%%%%%%%%%%%%%%%%%%%%%%%
il ne s'y est consacré à plein temps qu'apr\`es Pâques 1854. En fait,
durant l'hiver 1854, il reprend ses recherches sur les relations entre
l'\'electricit\'e, le magn\'etisme, la lumi\`ere et la
gravitation\footnote{\,
%%%%%%%%%%%%%%%%%%%%%%%-------DEBUT--------%%%%%%%%%%%%%%%%%%%%%%%%%%%
Dans sa lettre du 5 février 1854 à son frère
Wilhelm (no.~65, \cite{ neue1981a} p.~109), Riemann exprime clairement
quelle est sa direction de recherche principale à cette époque-là:
<<\,\deutschplain{Ich hatte gleich nach Ablieferung meiner
Habilitationsschrift wieder meine Untersuchungen über den Zusammenhang
der Naturgesetze fortgesetzt, und mich so darin vertieft, da{\ss} ich
nicht davon loskommen konnte.}\,>>
}, %%%%%%%%%%%%%%%%%%%%%%%%-----FIN-----%%%%%%%%%%%%%%%%%%%%%%%%%%%%%%%
tout en travaillant comme assistant de Weber \`a l'institut de
physique math\'ematique. Malheureusement, le mauvais temps hivernal et
une crise de surmenage le conduisent \`a la maladie\,\,---\,\,il était
hypocondriaque et il souffrait régulièrement de la fragilité de ses
poumons\,\,---, ce qui le contraint à interrompre ses travaux pour se
reposer \`a la campagne. Ayant recouvr\'e la sant\'e, il r\'edigera
donc sa conf\'erence d'\'epreuve du 
lundi de Pâque au lundi de Pentecôte, en sept
semaines environ\footnote{\,
%%%%%%%%%%%%%%%%%%%%%%%-------DEBUT--------%%%%%%%%%%%%%%%%%%%%%%%%%%%
Le programme de travail de Riemann a certainement connu des
alternances complexes que la biographie précieuse~\cite{ ded1892} de
Dedekind (\cf aussi~\cite{ bell1939, tazz2002}) était naturellement
dans l'incapacité de reconstituer, puisque dans sa lettre du 5 février
1854 publiée par Neuenschwander (\cite{ neue1981a}, p.~110), Riemann
confie à son frère: <<\,\deutschplain{Seit acht Tagen geht es mir nun
wieder besser, die Probevorlesung, die ich beim Colloquium halten soll
ist halb ausgearbeitet, und Dein Brief und der Gedanke an Dich sollen
mir ein Sporn sein, mich durch nichts wieder von dieser Arbeit
abbringen zu lassen.}\,>>
}. %%%%%%%%%%%%%%%%%%%%%%%%-----FIN-----%%%%%%%%%%%%%%%%%%%%%%%%%%%%%%%

\HEAD{Chapitre~1.\,\,\,\,L'ouverture riemannienne}{
1.2.\,\,\,Appr\'eciations d'universalité}

\medskip\noindent{\bf 1.2.~Appr\'eciations d'universalité.}
\label{appreciations-d-universalite}
Newman (\cite{ new1956}) qualifie d'<<\,imp\'eris\-sable\,>> ce discours
d'habilitation, qui rayonne encore d'une puissance philosophique et
math\'ematique singulière. C'est aussi l'un des tr\`es rares exemples
d'accession, en math\'ematique, au statut de classique 
intemporel\footnote{\,
%%%%%%%%%%%%%%%%%%%%%%%-------DEBUT--------%%%%%%%%%%%%%%%%%%%%%%%%%%%
En litt\'erature et en philosophie, la fr\'equentation r\'eguli\`ere
et l'\'etude ex\'eg\'etique du corpus classique font partie
int\'egrante de la formation sp\'ecialis\'ee; tel n'est pas le cas en
math\'ematiques.
}. %%%%%%%%%%%%%%%%%%%%%%%%-----FIN-----%%%%%%%%%%%%%%%%%%%%%%%%%%%%%%%

Dans son essence même, l'\deutschplain{Habilitationsvortrag} de
Riemann est en effet 
un {\em chef-d'{\oe}uvre remarquable d'inach\`evement et
d'ouverture}; de par les cons\'equences multiples qu'elle recèle,
elle a eu en effet une influence d\'eterminante quant au destin de
branches math\'ematiques neuves qui devaient être développées
ult\'erieurement, telles que par exemple: les fondements de la
g\'eométrie, la topologie, la g\'eom\'etrie diff\'erentielle, la
g\'eom\'etrie riemannienne ou finsl\'erienne,
\etc 

Toutefois, même si les considérations de Riemann sont apparemment
tr\`es accessibles \`a la lecture et semblent avoir été dictées par
une langue philosophique universelle et intemporelle, elles renferment
nombre d'affirmations \'enigmatiques; et comme ces affirmations 
remarquables ne sont pas justifi\'ees par des d\'emonstrations
math\'ematiques, elles ont aiguis\'e la sagacité des géomètres pendant
des décennies. Sans concession, Sophus Lie commentera les zones de
pénombre\footnote{\,
%%%%%%%%%%%%%%%%%%%%%%%-------DEBUT--------%%%%%%%%%%%%%%%%%%%%%%%%%%%
Le \S100 du Chapitre~22 p.~\pageref{S-100} ci-dessous offre une
analyse du discours d'habilitation de Riemann, et tout
particulièrement du \S~III, 1 p.~265 de~\cite{ riem1892} qui prétend
avoir complètement résolu le problème auquel Helmholtz, et surtout
Lie, ont consacré leurs forces de genèse conceptuelle.
} %%%%%%%%%%%%%%%%%%%%%%%%-----FIN-----%%%%%%%%%%%%%%%%%%%%%%%%%%%%%%% 
qui touchent à cette théorie entièrement nouvelle des groupes continus
de transformations, th\'eorie que Riemann ne poss\'edait manifestement
pas, et dont Lie allait faire l'{\oe}uvre monumentale de sa vie. On
peut s'imaginer n\'eanmoins que Riemann a \'etay\'e par des recherches
analytiques rigoureuses la plupart des propositions qu'il \'enonce
seulement dans un langage conceptuel, eu égard au devoir qu'il avait
vis-\`a-vis de son auditoire d'user au minimum d'un appareil
technique\footnote{\,
%%%%%%%%%%%%%%%%%%%%%%%-------DEBUT--------%%%%%%%%%%%%%%%%%%%%%%%%%%%
Seules \label{commentatio}
les quatre courtes pages de la seconde et derni\`ere partie de
la {\em Commentatio mathematica, qua respondere tentatur quaestioni ab
III\textsuperscript{ma} Academia Parisiensi propositae} (\cite{
riem1892}, pp.~380--383) dévoilent quelques calculs elliptiques en
relation avec la définition de la courbure sectionnelle qui apparaît
dans l'habilitation. Ce manuscrit de 1861 ne fut pas publié car le
prix de l'Académie de Paris proposé en 1858 n'a finalement pas été
attribué à Riemann. Depuis les travaux de Lipschitz et de Christoffel,
les calculs visionnaires (et cryptiques) de Riemann peuvent être
interprétés comme associant aux quantités de courbure sectionnelle
introduites par Riemann une certaine forme bilinéaire symétrique sur
l'espace des 2-plans infinitésimaux qui est aujourd'hui appelée {\em
tenseur de Riemann-Christoffel} et dont la connaissance recouvre tous
les invariants locaux de la métrique (\cite{ spiv1970}).
}. %%%%%%%%%%%%%%%%%%%%%%%%-----FIN-----%%%%%%%%%%%%%%%%%%%%%%%%%%%%%%%

Grâce \`a sa p\'en\'etration conceptuelle, ce texte allait donc
devenir une source d'inspiration r\'ecurrente dans le dernier tiers du
19\textsuperscript{i\`eme} si\`ecle\,\,---\,\,et aussi \`a la
charni\`ere du 20\textsuperscript{i\`eme}\,\,---, au moment o\`u la
clarification et l'approfondissement des concepts fondamentaux
s'affirmaient comme l'une des tendances dominantes en
math\'ematiques. Ainsi, faudra-t-il attendre les travaux de Dedekind,
Gehring, Clifford, Helmholtz, Christoffel, Lipschitz, Beltrami,
Frobenius, Lie, Killing, Engel, Ricci-Curbastro, Levi-Civita,
Schouten, \'E.~Cartan et d'autres pour mesurer l'ampleur des
d\'eveloppements inattendus que ces idées à peine esquissées
contenaient en germe.

\HEAD{Chapitre~1.\,\,\,\,L'ouverture riemannienne}{
1.3.\,\,\,Assembler l'inachev\'e}

\medskip\noindent{\bf 1.3.~Assembler l'inachev\'e.}
\label{assembler-l-inacheve}
C'est certainement la citation que Riemann a choisi de mettre en
exergue à ses \deutschplain{Fragmente philosophischen Inhalts}\footnote{\,
%%%%%%%%%%%%%%%%%%%%%%%-------DEBUT--------%%%%%%%%%%%%%%%%%%%%%%%%%%%
---\,\,publiés à titre posthume en 1876 dans ses \deutschplain{Gesammelte
Mathematische Werke}\,\,---
} %%%%%%%%%%%%%%%%%%%%%%%%-----FIN-----%%%%%%%%%%%%%%%%%%%%%%%%%%%%%%%
qui caractérise le mieux sa propre position 
dans ses 
travaux 
scientifiques\footnote{\,
%%%%%%%%%%%%%%%%%%%%%%%-------DEBUT--------%%%%%%%%%%%%%%%%%%%%%%%%%%%
En 1840, Riemann quitte la maison familiale à Quickborn pour entrer au
lycée à Hanovre. C'est à ce moment-là que débute sa correspondance
régulière avec sa famille. Neueunschwander (\cite{ neue1981a},
p.~90) qui a transcrit des lettres inédites signale que Riemann avait
déjà en ce temps-là de la peine à mener ses compositions à
bonne fin, parce qu'il rejetait continuellement ce qu'il avait déjà
écrit.
}: %%%%%%%%%%%%%%%%%%%%%%%%-----FIN-----%%%%%%%%%%%%%%%%%%%%%%%%%%%%%%%

\CITATION{Ne rejetez pas avec mépris les présents que j'ai rassemblés
pour vous avec dévotion avant de les avoir compris.
\REFERENCE{Lucrèce,~De~Natura~Rerum}}

\noindent 
Plus de la moitié de l'{\oe}uvre fascinante de Riemann est en effet
constituée de travaux qu'il jugeait {\em inaboutis} et qu'il s'est
refusé, pour cette raison,
à publier\footnote{\,
%%%%%%%%%%%%%%%%%%%%%%%-------DEBUT--------%%%%%%%%%%%%%%%%%%%%%%%%%%%
Autre exemple dans le domaine de la cr\'eation litt\'eraire, analys\'e
par Roland Barthes (\cite{ bart2003}, p.~343): <<\,Flaubert (1871, 50
ans): `Comme si de rien n'\'etait, je prends des notes pour mon {\em
Saint Antoine} [ce sera la troisi\`eme version], que je suis bien
d\'ecid\'e \`a ne pas publier quand il sera fini, {\em ce qui fait que
je travaille en toute libert\'e d'esprit}' [Lettre \`a Ernest Feydeau,
8 août 1871]; probl\`eme bien \'enonc\'e [\dots]. `Ne pas publier',
sorte de figure mi-rh\'etorique, mi-magique, utilis\'ee par beaucoup
d'\'ecrivains.\,>>
}. %%%%%%%%%%%%%%%%%%%%%%%%-----FIN-----%%%%%%%%%%%%%%%%%%%%%%%%%%%%%%%
D'un point de vue philosophique, la règle {\em riemannienne} de
direction de l'esprit consiste donc à {\em rassembler des éléments qui
ne sont pas compris}, à {\em formuler des questions réflexives} les
concernant, à {\em renverser les interrogations spéculatives}, à {\em
désigner les questions non résolues}. Comme Socrate, Riemann exprime
qu'il ne sait pas; cet \'etat de fait qui est impersonnel et
universel, le math\'ematicien doit l'accepter, puisqu'il fait partie
intégrante de l'essence même des mathématiques.

Une telle posture générale s'apparente donc plus à une {\em volonté
d'ignorance} qu'au doute systématique de Descartes, à ceci près que la
{\em volonté d'ignorance}, en mathématiques, ne peut pas être une
aporétique de principe comme l'est la maïeutique socratique,
puisqu'elle doit déboucher \`a terme sur des propositions rigoureuses,
sur des théorèmes, sur des connaissances adéquates. Analyser la {\em
métaphysique des mathématiques} que nous a léguée Riemann, c'est d'une
certaine manière entrer dans une {\em topologie de l'ouverture
accept\'ee} de la pens\'ee.

\HEAD{Chapitre~1.\,\,\,\,L'ouverture riemannienne}{
1.4.\,\,\,Le myst\`ere des notions primitives de la g\'eom\'etrie}

\medskip\noindent{\bf 1.4.~Le myst\`ere des
notions primitives de la g\'eom\'etrie.} 
\label{mystere-notions-primitives}
D'apr\`es Clifford (\cite{
cli1873}, p.~565), <<\,C'est Riemann qui le premier a accompli la
tâche d'analyser toutes les hypoth\`eses de la g\'eom\'etrie, et de
montrer leur ind\'ependance mutuelle\,>>.

Au commencement de son discours, Riemann annonce en effet que la
Géométrie classique d'Euclide admet comme données préalables les
concepts de point, de droite, d'espace, ainsi que les notions
d'incidence, d'angle et d'intersection. De toutes ces notions, elle ne
donne que des définitions nominales, implicites sur le plan
logique\footnote{\,
%%%%%%%%%%%%%%%%%%%%%%%-------DEBUT--------%%%%%%%%%%%%%%%%%%%%%%%%%%%
D'après Euclide, un point est ce qui est sans partie; une ligne est
une longueur sans épaisseur; un angle est l'inclinaison l'une sur
l'autre de deux lignes. D'inspiration métaphysique et procédant par
prédication {\em négative}, les deux premières définitions occultent
en effet la question de savoir où et comment asseoir ou fonder ces
concepts qui doivent par ailleurs rester opératoirement évidents pour
tout géomètre.
}. %%%%%%%%%%%%%%%%%%%%%%%%-----FIN-----%%%%%%%%%%%%%%%%%%%%%%%%%%%%%%% 
En ce qui
concerne les 
représentations mentales, les axiomes qui donnent naissance aux
relations démonstratives fondamentales sont toujours sous-tendus, en
filigrane, par un réseau d'intuitions archaïques; et parce que ces
intuitions aisées s'enracinent dans l'expérience physique de la
pensée, elles occultent l'interrogation mathématique abstraite, pure
et {\em a priori}. Or d'après Riemann, la position des hypothèses
primitives doit par nature être indécise, elle doit par principe faire
question en tant que telle, et donc, elle constituer un thème de
recherche nouveau, abstrait, pur et {\em a priori}.

\CITATION{Les rapports mutuels
de ces données primitives restent enveloppés de mystère \deutsch{im
Dunkeln}; on n'aperçoit pas bien si elles sont nécessairement liées
entre elles, ni jusqu'à quel point elles le sont, ni même {\em a
priori} si elles peuvent l'être.
\REFERENCE{\cite{riem1898},~p.~280.}}

\noindent
Dans cette seule phrase qui semble énoncer une constatation inspirée,
s'exprime un des traits les plus caractéristiques de la philosophie
riemannienne des mathématiques: 

\smallskip

\centerline{l'{\sf ouverture},}

\smallskip\noindent
l'ouverture {\em dite}, 
l'ouverture {\em écrite}, 
l'ouverture {\em assumée}, 
l'ouverture {\em maintenue}, 
fussent-elles ombre, obscurité ou même {\em mystère}, comme Hoüel a si
bien choisi de le traduire.

\CITATION{Depuis 
Euclide jusqu'à Legendre, pour ne citer que le plus illustre
des réformateurs modernes de la Géométrie, personne, parmi les
mathématiciens ni parmi les philosophes, n'est parvenu à éclaircir ce
mystère \deutsch{Dunkelheit}.
\REFERENCE{\cite{riem1898},~p.~280.}}

\noindent
Riemann est le premier, dans l'histoire des mathématiques, à insister
explicitement dans ses écrits, sans en passer par la formulation de
problèmes ouverts précis, sur la présence constante de l'inachèvement.
C'est en cela qu'il est immortel. 

\HEAD{Chapitre~1.\,\,\,\,L'ouverture riemannienne}{
1.5.\,\,\,Fondements de la géométrie}

\medskip\noindent{\bf 1.5.~Fondements de la géométrie.}
\label{fondements-de-la-geometrie}
Riemann connaissait-il les travaux de géométrie non-euclidienne dus à
ses contemporains? Hormis les {\em \'Eléments de géométrie}~\cite{
lege1817} de Legendre, La bibliothèque de l'université de Göttingen
possédait une copie d'un travail de B\'olyai (\cf~\cite{ boly1903}),
et quelques publications de N.I.~Lobatchevski\u{\i}, notamment l'essai
sur la <<\,géométrie imaginaire\,>> paru au
\deutschplain{Journal für die reine und angewandte 
Mathematik}\footnote{\,
%%%%%%%%%%%%%%%%%%%%%%%-------DEBUT--------%%%%%%%%%%%%%%%%%%%%%%%%%%% 
\,\,---\,\,volume~17 (1837), 295--320\,\,---\,\,
} %%%%%%%%%%%%%%%%%%%%%%%%-----FIN-----%%%%%%%%%%%%%%%%%%%%%%%%%%%%%%%%
et son livre-fascicule {\em Recherches
géométriques sur la théorie des parallèles}
(\cf~\cite{ loba1914}) publié à Berlin en 1840.
Riemann pourrait avoir consulté ces sources, et d'après Neuenschwander
(\cf~aussi \cite{ scho1982a, scho1992}), il a effectivement emprunté
le volume~17 du Journal de Crelle le 15 février 1854.

Toutefois, l'étude du \deutschplain{Nachlass} conduite par Scholz
dans~\cite{ scho1982a} montre que Riemann n'était pas réellement
intéressé, comme avait pu l'être Gauss, par les fondements de la
géométrie élémentaire. Fait surprenant, les archives explorées ne
portent aucune indication du fait que Riemann ait pu se constituer une
connaissance circonstanciée des travaux de B\'olyai et de
Lobatchevski\u{\i}. En tout cas, dans son discours d'habilitation, il
mentionne seulement Legendre, qu'il considère comme <<\,le plus
illustre des réformateurs modernes de la géométrie\,>>. Fait tout
aussi étonnant, Riemann ne cite pas une seule fois 
l'axiome des parallèles,
même lorsqu'il traite des variétés de courbure constante au sein
desquelles l'existence et le comportement des parallèles peuvent être
caractérisés de manière extrêmement limpide en fonction du signe de la
courbure (constante). En vérité, l'analyse par Scholz~\cite{
scho1982a, scho1992} du feuillet no.~40 du dossier no.~16 du
\deutschplain{Nachlass} montre le peu d'intérêt que
Riemann aurait pu trouver à s'engager, dans des recherches de type
logique ou fondationnel, en partant des postulats
d'Euclide.

\CITATION{Même 
s'il est intéressant de saisir un tel mode de traitement de la
géométrie, une telle entreprise serait extrêmement infructueuse, 
car de la sorte, on ne trouverait pas de nouveaux
théorèmes, et ce qui apparaît simple et clair dans la présentation de
l'espace deviendrait de cette manière-là compliqué et difficile.
\REFERENCE{\cite{ scho1992},~pp.~28--29.}}

\noindent
Ainsi Riemann semble-t-il écarter la voix d'axiomatisation {\em a
posteriori} des géométries comme systèmes logiques clos et cohérents,
telle qu'elle devait naître dans les années 1880 avec les travaux de
Pasch\footnote{\,
%%%%%%%%%%%%%%%%%%%%%%%-------DEBUT--------%%%%%%%%%%%%%%%%%%%%%%%%%%%
Freudenthal (\cite{ freu1960b}, p.~617) souligne que Pasch dans ses
\deutsch{Vorlesungen über neuere Geometrie} (1882) anticipe très
largement le point de vue formaliste pur dont Hilbert s'était fait
l'ardent défenseur dans ses \deutschplain{Grundlagen der Geometrie}
(1899): <<\,\`A chaque fois\,>>, écrit en effet Pasch, <<\,que la
géométrie doit être réellement déductive, le procédé d'inférence doit
être indépendant aussi bien de la signification des notions
géométriques que des figures. Les seules choses qui comptent sont les
{\em relations} entre les notions géométriques, telles qu'elles sont
établies dans les théorèmes et utilisées dans les définitions\,>>. 
}, %%%%%%%%%%%%%%%%%%%%%%%%-----FIN-----%%%%%%%%%%%%%%%%%%%%%%%%%%%%%%%
Stolz, Schur,
et ultérieurement de Veronese, Killing, Enriques, Pieri, Padoa,
Russell, Hilbert, Poincaré. Mais sur la base d'une analyse de la
méthodologie philosophique de Herbart qui a influencé Riemann, on peut
néanmoins soutenir (\cf~ce qui va suivre) que Riemann {\em anticipe
l'organisation structurale et hiérarchique des concepts mathématiques
modernes}, et ce, bien avant que naisse la méthode axiomatique
hilbertienne proprement dite qui devait conférer un sens
métamathématique précis aux concepts de cohérence, de complétude,
d'indépendance, de suffisance, de nécessité, de spécialisation et de
catégoricité.
 
\HEAD{Chapitre~1.\,\,\,\,L'ouverture riemannienne}{
1.6.\,\,\,Le renversement riemannien}

\medskip\noindent{\bf 1.6.~Le renversement riemannien.}
\label{le-renversement-riemannien}
Riemann était principalement intéressé par des questions d'un type
nouveau, et qui font mystère en elles-mêmes et par elles-mêmes; en
acceptant d'explorer ces questions, il est possible de faire jouer aux
racines de la connaissance mathématique le rôle de branches nouvelles
en devenir. Renverser les questions d'essence et de conceptualisation,
c'est un acte métaphysique par excellence. Sur le plan de la théorie
de la connaissance, Riemann rejoint donc l'exigence critique de la
philosophie qui a marqué son temps.

En effet, Kant caractérisait de manière imagée sa <<\,solution
transcendantale\,>> au problème général de la raison pure comme un
{\sl renversement copernicien}\footnote{\,
%%%%%%%%%%%%%%%%%%%%%%%-------DEBUT--------%%%%%%%%%%%%%%%%%%%%%%%%%%%
L'astronome polonais Nicolas Copernic (1473--1543) proposa de
substituer au géocentrisme l'héliocentrisme (\voir~\cite{ szcz1998}
pour une étude philosophique). <<\,Il en est précisément ici\,>>,
écrit Kant (\cite{ kant1986}, p.~19) dans sa seconde Préface à la {\em
Critique de la Raison pure} afin de caractériser son hypothèse
fondamentale par une analogie imagée, <<\,comme de la première idée de
Copernic; voyant qu'il ne pouvait pas réussir à expliquer les
mouvements du ciel en admettant que toute l'armée des étoiles évoluait
autour du spectateur, il chercha s'il n'aurait pas plus de succès en
faisant tourner l'observateur lui-même autour des astres immobiles. Or
en Métaphysique, on peut faire un pareil essai, pour ce qui est de
l'intuition des objets, Si l'intuition devait se régler sur la nature
des objets, je ne vois pas comment on en pourrait connaître quelque
chose {\em a priori}; si l'objet, au contraire (en tant qu'objet
\deutsch{Object}
des sens), se règle sur la nature de notre pouvoir d'intuition, je
puis me représenter à merveille cette possibilité.\,>>
}: %%%%%%%%%%%%%%%%%%%%%%%%-----FIN-----%%%%%%%%%%%%%%%%%%%%%%%%%%%%%%%
le fait qu'il existe, en physique et en mathématiques, des
connaissances synthétiques {\em a priori}\,\,---\,\,fait incontestable
dont la philosophie (qui s'en étonne spontanément) se doit de rendre
compte\,\,---\,\,peut trouver toute son intelligibilité, d'apr\`es
Kant, si l'on fait l'hypothèse (révolutionnaire en philosophie) que
{\em les objets se règlent sur notre connaissance}, et que ces objets
sur lesquels se règle notre connaissance sont à distinguer
rigoureusement des choses telles qu'elles sont en elles-mêmes et
auxquelles notre entendement limité ne peut pas prétendre avoir
accès. Ainsi, autant une connaissance {\em synthétique} et {\em a
priori} des choses en elles-mêmes exposerait à d'invraisemblables
obscurités et à une anarchie sans fin de contradictions\footnote{\,
%%%%%%%%%%%%%%%%%%%%%%%-------DEBUT--------%%%%%%%%%%%%%%%%%%%%%%%%%%% 
<<\,Le terrain \deutsch{Kampfplatz} où se livrent ces combats sans fin
se nomme la {\sl Métaphysique}\,>> (\cite{ kant1986}, p.~5).
}, %%%%%%%%%%%%%%%%%%%%%%%%-----FIN-----%%%%%%%%%%%%%%%%%%%%%%%%%%%%%%%
autant une connaissance synthétique {\em a priori} des objets
d'expérience, en tant que ces objets trouvent leur source dans notre
sensibilité, dans notre intuition, et dans notre entendement, permet
d'établir un vrai rapport de structuration et de fonder ainsi sans
conteste la {\em part d'aprioricité} de la connaissance. 

Pour Kant, l'intuition ne se règle donc pas sur la nature de l'objet,
mais c'est l'objet, en tant qu'objet des sens, qui se {\em règle} sur
la nature de notre pouvoir d'intuition, ou plus exactement, qui se
dévoile dans et par ce que nous en explorons, dans et par ce que nous
en disons, et ainsi, avec le pouvoir de préformation qui nous
appartient en propre, il doit en aller de même et de manière
absolument générale quant aux {\em concepts} par lesquels
l'entendement élabore les déterminations de la
connaissance\footnote{\,
%%%%%%%%%%%%%%%%%%%%%%%-------DEBUT--------%%%%%%%%%%%%%%%%%%%%%%%%%%%
Postérité remarquable de cette th\`ese kantienne au sujet de notre
pratique de l'algèbre: bien que les id\'ealit\'es alg\'ebriques
(groupe de Galois; structure d'une alg\`ebre de Lie; diagramme de
Dynkin; groupe fini produit par g\'en\'erateurs et relations; syzygies
entre invariants alg\'ebriques) poss\`edent un r\'eseau {\em non
spatialis\'e} et {\em non temporalis\'e} de relations complexes, notre
compr\'ehension de ces relations et les d\'emonstrations que nous
\'elaborons pour les {\em parcourir} sont in\'evitablement {\em
lin\'eaires} et {\em successives}. L'entendement pr\'eformant
(et d\'eformant) n'approche le multiple alg\'ebrique que 
par actions discr\`etes.
} %%%%%%%%%%%%%%%%%%%%%%%%-----FIN-----%%%%%%%%%%%%%%%%%%%%%%%%%%%%%%%
(\cf~\cite{
kant1986}, p~19). <<\,Nous ne connaissons {\em a priori} des choses
que ce que nous y mettons nous-mêmes\,>>: c'est l'hypothèse
fondamentale (risquée) du système que Kant a érigé afin de délimiter
avec exactitude la portée des raisonnements métaphysiques qui
dépassent les limites de toute expérience.

Ainsi d'après Kant, l'entendement ne doit pas seulement être d\'efini
comme le pouvoir d'élaborer des règles empiriques, en raisonnant par
induction, en comparant, et en extrapolant, mais surtout, à un niveau
véritablement transcendantal, comme le {\em pouvoir structurel de
prescrire en quelque sorte ses propres règles à la nature}, au sens où
les objets d'expérience sont nécessairement conformes aux conditions
{\em a priori} dans lesquelles ils peuvent être perçus et
pensés\,\,---\,\,autrement dit: les formes de la sensibilité et les
contraintes logiques du raisonnement préforment en quelque sorte
toutes les connaissances qui sont synthétiques {\em a priori}.

En quoi consiste alors, chez Riemann, le renversement philosophique
par rapport à la tradition mathématique qui le précède? Certainement
pas dans l'élaboration d'une théorie générale de la connaissance qui
serait destinée à expliquer les antinomies de la raison, les
apparences transcendantales, et les paralogismes des preuves
métaphysiques. Les mathématiques ont rarement besoin qu'on leur
apprenne à corriger leurs raisonnements. Toutefois, dans la trajectoire
philosophique de Riemann, on trouve une analogie de fond avec la
solution kantienne au problème de la 
métaphysique\footnote{\, 
%%%%%%%%%%%%%%%%%%%%%%%-------DEBUT--------%%%%%%%%%%%%%%%%%%%%%%%%%%%
On ne méditera jamais assez les toutes premières lignes de la {\em
Critique de la raison pure} (\cite{ kant1986}, p.~5): <<\,La raison
humaine a cette destinée singulière, dans un genre de ses
connaissances, d'être accablée de questions qu'elle ne saurait éviter,
car elles lui sont imposées par sa nature même, mais auxquelles elle
ne peut répondre, parce qu'elles dépassent totalement le pouvoir de la
raison humaine\,>>. Ici, dans le parallélisme dialectique entre Kant
et Riemann, se dessine alors la question extrêmement délicate de la
{\em démarcation} entre les mathématiques et la métaphysique, en tant
que l'une et l'autre ont le même rapport à un champ de questions
spontanées (et <<\,accablantes\,>>), sans pour autant emprunter la
même voie quant aux réponses qu'elles sont susceptibles d'élaborer.
}. %%%%%%%%%%%%%%%%%%%%%%%%-----FIN-----%%%%%%%%%%%%%%%%%%%%%%%%%%%%%%%

Chez Riemann en effet, du point de vue de la théorie de la
connaissance en général, le renversement consiste principalement dans
la désignation\,\,---\,\,implicite, non théorisée et d'ordre
méta-mathématique\,\,---\,\,de {\em caractères universels du
questionnement mathématique}; presque {\em a priori}, ces caractères
s'intègrent aux structures fondamentales de la pensée mathématique,
mais ils ne la préforment pas, ils ne la prédéfinissent pas, et ils ne
la gouvernent pas. C'est que Riemann, dont l'humilité est légendaire,
fait spontanément preuve d'une extrême prudence quant à l'énonciation
de toute thèse métaphysique directrice. Ainsi est-il naturellement
préservé de tout engagement dogmatique ou de toute tentation de {\em
fermer} un ordre de questions. Jamais en effet il n'affirmerait avoir
résolu de manière entièrement satisfaisante un problème mathématique
donné. Son {\oe}uvre n'est pas seulement riche d'invention
conceptuelle, mais elle est aussi\,\,---\,\,et c'est en cela que
réside sa force philosophique majeure\,\,---\,\,confondante d'{\em
inachèvement volontaire}. Kant soutenait, lui, que grâce à
l'élaboration de cette <<\,science particulière\,>> qu'il appelait la
{\sl critique de la raison pure}, il avait pu résoudre les questions
en les dissolvant lorsqu'elles s'avéraient être non
résolubles\footnote{\,
%%%%%%%%%%%%%%%%%%%%%%%-------DEBUT--------%%%%%%%%%%%%%%%%%%%%%%%%%%%
<<\,Après avoir découvert le point du malentendu de la raison avec
elle-même, je les [ces questions] ai résolues à son entière
satisfaction. [\dots] Et je n'ose dire qu'il ne doit pas y avoir un
seul problème métaphysique qui ne soit ici résolu.\,>> (\cite{
kant1986}).
}. %%%%%%%%%%%%%%%%%%%%%%%%-----FIN-----%%%%%%%%%%%%%%%%%%%%%%%%%%%%%%%
Affirmer avoir résolu une question difficile, même au moyen d'un
grand sytème de pensée, c'est un risque que Riemann ne prend pas.

Conséquence limpide: la seule certitude qu'une <<\,méta\-physique
riemannienne\,>> proprement universelle pourrait affirmer, c'est qu'il
faut {\em maintenir ouvertes les questions inaugurales tout au long de
leur destin}. Pour un mathématicien comme Riemann, il faut donc
élaborer un suivi spéculatif constant des problèmes à travers leur
histoire, et donc par la même occasion, il faut m\'editer la
continuité des incertitudes qui s'expriment, se d\'eveloppent, et se
complexifient. Comme pourrait le soutenir Heidegger, il s'agit de ne
jamais {\em oublier} la question initiale dans toutes les questions
qui lui sont subordonnées. Une telle th\`ese est un universel
philosophique qui s'exprime aussi {\em dans} les math\'ematiques.

Il s'agit ainsi, pour la saine philosophie riemannienne qui s'exerce
toujours dans les mathématiques contemporaines, d'examiner
régulièrement dans quelle mesure les travaux qui paraissent d'année en
année apportent, ou n'apportent pas une réponse partielle, 
une réponse satisfaisante, une réponse complète, voire même
une réponse définitive\footnote{\,
%%%%%%%%%%%%%%%%%%%%%%%-------DEBUT--------%%%%%%%%%%%%%%%%%%%%%%%%%%%
En arrière-plan de ces considérations se profile donc un problème
particulièrement délicat pour la philosophie des mathématiques, à
savoir: quand, pourquoi, comment et à quelles conditions peut-on
déclarer qu'une question mathématique a été complètement résolue par
un théorème, ou par une théorie? \'Etant donné le renouveau des
problèmes d'effectivité en algèbre que suscite la maturité
grandissante des calculateurs électroniques, la théorie de Galois
fournit un exemple int\'eressant d'<<\,illusion d'aboutissement par
l'abstraction non calculatoire\,>>, puisque dans la pratique, la
détermination du groupe de Galois sur $\Q$ d'un polynôme à
coefficients entiers n'est en rien résolue par la correspondance
biunivoque (et presque tautologique) entre tous les sous-groupes du
groupe de Galois d'une extension finie normale, et tous les sous-corps
de cette extension qui contiennent $\Q$ (\cite{ esco2001}; \voir~le
m\'emoire de synth\`ese~\cite{ vali2008} et sa bibliographie pour les
r\'esultats les plus r\'ecents de la th\'eorie de Galois
algorithmique). Le quantificateur universel ``tout'' cache ici des
ordres de question nombreux, notamment le probl\`eme de {\em
classifier a priori} tous les groupes finis qui peuvent se r\'ealiser
comme groupes de Galois, qui est essentiellement le {\em probl\`eme
principal sur lequel d\'ebouche la correspondance de
Galois}. Signalons seulement que Jordan traitera de la classification
des sous-groupes du groupe des permutations de $n \geqslant 1$
lettres, et que Lie prendra aussi \`a bras-le-corps le probl\`eme de
{\em classifier} les groupes continus (finis ou infinis) de
transformations.
}. %%%%%%%%%%%%%%%%%%%%%%%%-----FIN-----%%%%%%%%%%%%%%%%%%%%%%%%%%%%%%%
Le traitement du {\em problème de Riemann} par Helmholtz et par Lie
avec la théorie des groupes continus de transformations (\cf~ce qui va
suivre), puis par Kolmogorov, Busemann
et par Freudenthal avec les concepts de
la topologie générale, montre qu'une question initiale peut se
métamorphoser, se ramifier et se diversifier en s'approfondissant.
Dans une telle optique, si dogmatisme doctrinaire il peut y avoir,
c'est un {\em dogmatisme de l'attente et de l'imprévisibilité}. La
position {\em philosophique} de Riemann, non exprimée par lui comme
système et intrinsèquement en adéquation avec le caractère
fondamentalement imprévisible et ouvert de l'{\sl
irréversible-synthétique} (\cf~\cite{ merk2008}) en mathématiques, ne
peut être entach\'ee par une d\'ecision ontologique unilat\'erale, que
ce soit un réalisme, un idéalisme, un transcendentalisme, un
empirisme, un relativisme, ou un scepticisme. Presque par pétition de
principe, certains problèmes acquièrent donc une ouverture dialectique
indéfinie qui les rend inépuisables.
\`A tout le moins, très peu de questions mathématiques
peuvent être considérées comme complètement résolues: là est toute la
teneur du {\em renversement riemannien}.

Voilà donc pour le niveau absolument général de la philosophie
riemannienne des mathématiques.

\HEAD{Chapitre~1.\,\,\,\,L'ouverture riemannienne}{
1.7.\,\,\,Décider l'ouverture problématique du conceptuel}

\medskip\noindent{\bf 1.7.~Décider 
l'ouverture problématique du conceptuel.}
\label{decider-l-ouverture}
Quant aux concepts fondamentaux de la géométrie, Riemann va opérer
dans son habilitation une {\em deuxième révolution de pensée} par
spécification de nouvelles arborescences mathématiques potentielles.

Aux structures prédéfinies de l'intuition et de l'entendement qui
exigeaient chez Kant l'élaboration d'une esthétique transcendantale et
d'une théorie de l'expérience synthétique, se substitue en effet chez
Riemann l'affirmation méthodologique (implicite) qu'il existe des
caractères {\em a priori}, constants et reproductibles de
l'interrogation mathématique, et que ces caractères recouvrent trois
grands domaines inépuisables d'investigation en agissant comme
principes de genèse pour les mathématiques tout entières:

\smallskip\noindent$\bullet$\,\,
la question de l'être et de la nature des objets
mathématiques; 

\smallskip\noindent$\bullet$\,\,
la question des liens de dépendance que ces êtres mathématiques
purement problématiques entretiennent entre eux;

\smallskip\noindent$\bullet$\,\,
la question, d'inspiration proprement kantienne, pour les intuitions
mathématiques comme pour les conceptions mathématiques, de leurs
conditions de possibilité.

\medskip\noindent
Absolument universelles, ces questions g\'en\'erales concernent toutes
les mathématiques. Ainsi donc, il ne s'agit pas pour Riemann, au début
de son discours d'habilitation, de mieux définir et de mieux fonder
les notions euclidiennes primitives de point, de droite, ou de plan,
mais il s'agit plutôt de s'interroger sur ce qu'est et sur ce que peut
être l'espace, et donc par voie de conséquence, sur les {\em
conditions mêmes de possibilité, {\em de diversité} et de connexion
pour les conceptualisations {\em éventuelles} de la notion
problématique d'<<\,espace\,>>}. Ceci, par excellence, c'est de la
philosophie, et une telle philosophie est appelée à se réaliser {\em
dans et par} les mathématiques.

\CITATION{La 
raison [de ce mystère] est que le concept général des grandeurs de
dimensions multiples, comprenant comme cas particulier les grandeurs
étendues, n'a jamais été l'objet d'aucune étude. En conséquence, je me
suis proposé d'abord le problème de construire, en partant du concept
général de grandeur \deutsch{aus
allgemeinen Grössenbegriffen}, le concept d'une grandeur de dimensions
multiples \deutsch{Begriff einer mehrfach ausgedehnten Grösse}.
\REFERENCE{\cite{riem1898},~pp.~280--281.}}

\HEAD{Chapitre~1.\,\,\,\,L'ouverture riemannienne}{
1.8.\,\,\,Nécessité, suffisance, bifurcation}

\medskip\noindent{\bf 1.8.~Nécessité, suffisance, bifurcation.}
\label{necessite-suffisance-bifurcation}
Renversement riemannien, donc, de la théorie des sciences
mathématiques: les concepts ne sont pas donnés dans l'expérience du
monde sensible, ils ne sont pas produits par une intuition {\em sui
generis}, et ils ne sont pas d'emblée présentés dans une évidence
définitionnelle\footnote{\,
%%%%%%%%%%%%%%%%%%%%%%%-------DEBUT--------%%%%%%%%%%%%%%%%%%%%%%%%%%%
Lire et méditer Riemann permet de se départir d'une certaine
<<\,naïveté r\'ealiste\,>> quant \`a l'<<\,\'evidence\,>> de l'acte de
pens\'ee que repr\'esente le fait de poser, ou de rappeler la
d\'efinition d'un objet math\'ematique, par exemple, dans un article
de recherche ou lors d'une conf\'erence.
}. %%%%%%%%%%%%%%%%%%%%%%%%-----FIN-----%%%%%%%%%%%%%%%%%%%%%%%%%%%%%%%
Au contraire, les concepts portent de manière permanente les marques
de leur propre problématicité; ils sont donc {\em à construire}, et
ils sont aussi à cause de cela {\em ouverts}.

Dans le \S~I où il introduit la notion, maintenant si fondamentale
en géométrie différentielle, de 
{\sl multiplicité}\footnote{\,
%%%%%%%%%%%%%%%%%%%%%%%-------DEBUT--------%%%%%%%%%%%%%%%%%%%%%%%%%%%
Vuillemin~\cite{ vui1962} et Chorlay~\cite{ chor2009} sugg\`erent de
ne pas traduire le terme <<\,\deutschplain{Mannigfaltigkeit}\,>>
introduit par Riemann par <<\,vari\'et\'e\,>>, mais par
<<\,multiplicit\'e\,>>.
} %%%%%%%%%%%%%%%%%%%%%%%%-----FIN-----%%%%%%%%%%%%%%%%%%%%%%%%%%%%%%%
\deutsch{Mannigfaltigkeit}\footnote{\, 
%%%%%%%%%%%%%%%%%%%%%%%-------DEBUT--------%%%%%%%%%%%%%%%%%%%%%%%%%%%
Pour une analyse historique du concept de <<\,variété\,>>, qui
débouchera ultérieurement avec les travaux de Weyl, Cartan, Veblen et
Whitehead sur la d\'efinition contemporaine en termes de cartes et
d'atlas maximal, nous renvoyons aux études~\cite{ scho1980, scho1982a,
flam1992, bfs1992, boi1995a, boi1995b, laug1999, chor2009}.
}, %%%%%%%%%%%%%%%%%%%%%%%%-----FIN-----%%%%%%%%%%%%%%%%%%%%%%%%%%%%%%% 
Riemann procède comme il l'a déjà fait dans sa {\em Dissertation
inaugurale} et dans son mémoire sur les séries trigonométriques: à
chaque étape d'une conceptualisation problématique, il recherche à la
fois les conditions qui sont {\em nécessaires}, et les conditions qui
sont {\em suffisantes} à la genèse du concept ou des concepts visés.
Les premières lignes du \S~1 en t\'emoignent: 

\CITATION{Les 
concepts de grandeur \deutsch{Grössenbegriffe} ne sont possibles
que là où il existe un concept général qui permette différents modes
de détermination
\deutsch{verschiedene Bestimmungsweisen}. 
\REFERENCE{\cite{riem1898},~p.~282.}}

\noindent
Condition nécessaire, donc: <<\,n'être possible que si\,\dots\,>>.
Avant de penser aux concepts de grandeur en question, il faut disposer
d'un concept préalable par lequel ces grandeurs pourraient être
déterminées, et il faut se donner les moyens de {\em comparer} ces
déterminations entre elles. Dans la donation, qui est nécessaire, l'Un
est pr\'ec\'ed\'e par le Multiple, qui est lui aussi n\'ecessaire.

Seconde constatation, frappante de lucidit\'e philosophique: comme
s'il respectait rigoureusement les antinomies entre le discret et le
continu qui ont pr\'eoccup\'e la philosophie grecque (\cf~notamment la
discusion par Aristote des paradoxes de Zénon dans la {\em Physique}
III), Riemann enracine d'emblée sa réflexion dans la {\em
bifurcation philosophique incontournable} entre le discret et le
continu.

\CITATION{Suivant 
qu'il est, ou non, possible de passer de l'un de ces modes de
détermination à un autre, d'une manière continue, ils forment une
multiplicité \deutsch{Mannigfaltigkeit}\footnotemark continue ou une
multiplicité discr\`ete; chacun en particulier de ces modes de
d\'etermination s'appelle, dans le premier cas, un point, dans le
second un \'el\'ement de cette multiplicité.
\REFERENCE{\cite{riem1898},~p.~282.}}

\noindent
Mais avant de poursuivre cette analyse philosophique ciblée du 
discours de Riemann (\voir~p.~\pageref{reprise-analyse} 
ci-dessous), un détour par une théorie générale
de la connaissance qui l'a marqué s'impose. 
\label{wait-zenon}

\footnotetext{\,
%%%%%%%%%%%%%%%%%%%%%%%-------DEBUT--------%%%%%%%%%%%%%%%%%%%%%%%%%%%
{\sl Varietas}, \cf~le mémoire de Gauss: {\em Theoria residuorum
biquadratorum}, et {\em Anzeige zu derselben}, Werke, Bd.~II, pp.~110,
116 u.~118. Dans ses leçons de 1850--51 sur la méthode des moindres
carrés, Gauss émettait des remarques occasionnelles sur l'extension de
sa théorie des surfaces à la dimension $n$ quelconque. Les notes
calligraphiées de ce cours furent prises par A.~Ritter qui soutint sa
thèse en 1853 sous la direction de Gauss. D'après Stäckel (\cite{
stac1873}, p.~55), A.~Ritter était ami proche de Riemann dans les
années 1850 à 1853, et c'est probablement par son intermédiaire que
Riemann est entré en contact avec les idées de Gauss sur la géométrie
(\cite{ krey1994}, p.~114).
} %%%%%%%%%%%%%%%%%%%%%%%%-----FIN-----%%%%%%%%%%%%%%%%%%%%%%%%%%%%%%%

\HEAD{Chapitre~1.\,\,\,\,L'ouverture riemannienne}{
1.9.\,\,\,L'influence épistémologique de Herbart}

\label{influence-de-Herbart}
\medskip\noindent{\bf 1.9.~L'influence épistémologique de 
Herbart\footnote{\,
%%%%%%%%%%%%%%%%%%%%%%%-------DEBUT--------%%%%%%%%%%%%%%%%%%%%%%%%%%%
Le lecteur est renvoy\'e \`a l'\'etude~\cite{ scho1982b} 
de Scholz pour de plus amples informations. 
} %%%%%%%%%%%%%%%%%%%%%%%%-----FIN-----%%%%%%%%%%%%%%%%%%%%%%%%%%%%%%%
sur Riemann.} En 1809, apr\`es des \'etudes en droit, philosophie,
litt\'erature et math\'ematiques \`a I\'ena, un doctorat puis une
habilitation \`a Göttingen dans lesquels il \'elabore les fondations
de son système, le philosophe allemand Herbart (1776--1841) succ\`ede
\`a Kant sur la chaire de philosophie de l'Universit\'e de
Königsberg. Il crée un s\'eminaire sur la p\'edagogie qu'il associe
\`a une \'ecole pratique exp\'erimentale, et
il participe aux r\'eformes de l'\'education en Prusse. En 1834, il
obtient la chaire de philosophie
\`a l'Universit\'e de Göttingen o\`u il enseigne la philosophie et la
p\'edagogie jusqu'\`a sa mort en 1841. Ses vues sont alors largement
diffus\'ees en Allemagne, et c'est par la lecture que Riemann entre en
contact avec sa métaphysique.

Le dossier no.~18 du \deutschplain{Nachlass} de Riemann, conservé à la
bibliothèque de l'Université de Göttingen, comporte 203 feuillets
manuscrits regroup\'es sous le titre de <<\,\deutschplain{Fragmente
naturphilosophischen Inhalts}\,>>. Ceux de ces fragments dont se
dégage une certaine cohérence de pensée ont été édités par Dedekind et
Weber à la fin des {\oe}uvres compl\`etes~\cite{ riem1892}\footnote{\,
%%%%%%%%%%%%%%%%%%%%%%%-------DEBUT--------%%%%%%%%%%%%%%%%%%%%%%%%%%%
{\em Voir}~\cite{ riefu2002} pour une traduction partielle en
français.
}. %%%%%%%%%%%%%%%%%%%%%%%%-----FIN-----%%%%%%%%%%%%%%%%%%%%%%%%%%%%%%%
D'apr\`es Scholz~\cite{ scho1982b} qui a sollicit\'e l'aide de
l'\deutschplain{Handschriftenabteilung der Göttinger
Universitätsbibliothek}, 12 feuillets contiennent des notes de lecture
consacr\'ees \`a Herbart, que Riemann enrichit parfois de ses
propres formulations. Les extraits sont tous tir\'es des travaux de
Herbart sur la m\'etaphysique et sur la psychologie.

\CITATION{L'auteur 
[Riemann lui-même] se consid\`ere comme herbartien en
psychologie et en \'epist\'emologie (m\'ethodologie et
\'eidolologie\footnotemark);
toutefois dans la plupart des cas, il ne peut s'accorder avec la
philosophie naturelle de Herbart, et avec les disciplines
m\'etaphysiques (ontologie et syn\'echologie ) qui s'y r\'ef\`erent.
\REFERENCE{\cite{riem1898},~p.~508.}}

\footnotetext{\,
%%%%%%%%%%%%%%%%%%%%%%%-------DEBUT--------%%%%%%%%%%%%%%%%%%%%%%%%%%%
D'apr\`es Herbart, la m\'etaphysique englobe quatre
disciplines: la m\'ethodologie, 
l'\'eidolologie \deutsch{Eidolologie}, l'ontologie
et la syn\'echologie \deutsch{Synechologie}. 
} %%%%%%%%%%%%%%%%%%%%%%%%-----FIN-----%%%%%%%%%%%%%%%%%%%%%%%%%%%%%%%

\noindent
Bien que la notion de {\sl présentations sérielles}
\deutsch{Vorstellungsreihen}\footnote{\,
%%%%%%%%%%%%%%%%%%%%%%%-------DEBUT--------%%%%%%%%%%%%%%%%%%%%%%%%%%%
Herbart considérait que les concepts géométriques tels qu'ils sont formés
et développés par les sciences trouvent leur origine dans les
représentations spatiales qui donnent accès à l'expérience du monde
sensible, et qui sont intrinsèquement {\em sérielles}. En effet,
toutes les formes de perception humaine de l'espace s'effectuent {\em
en séries}, c'est-à-dire de manière successive et temporalisée, dans
une continuité de discontinuités coprésentes, et elle s'exerce de
manière locale, car elle se restreint à des voisinages immédiats de
l'action individuelle. Aussi les {\em séries de présentation}
\deutsch{Vorstellungsreihen} s'ordonnent-elles et
se connectent-elles l'une à l'autre pour engendrer des représentations
spatiales. Une <<\,forme sérielle continue\,>> se présente
lorsqu'une classe spécifique de représentations est soumise à une
fusion graduelle et unifiée.
} %%%%%%%%%%%%%%%%%%%%%%%%-----FIN-----%%%%%%%%%%%%%%%%%%%%%%%%%%%%%%%
due à 
Herbart ait certainement pu inspirer Riemann dans son élaboration du
concept de <<\,multiplicit\'e\,>>, la {\sl synéchologie}\footnote{\,
%%%%%%%%%%%%%%%%%%%%%%%-------DEBUT--------%%%%%%%%%%%%%%%%%%%%%%%%%%%
\`A partir de la théorie psychologique 
des formes sérielles de représentations (\cf~la note précédente),
l'engendrement de {\em formes sérielles continues}
\deutsch{continuirliche Reihenformen} {\em de concepts}
constitue, selon une procédure théorique très élaborée, l'objectif de
la discipline métaphysique que Herbart appelle la {\sl synéchologie},
et qui désigne essentiellement sa théorie philosophique des fondements
du concept d'espace.
} %%%%%%%%%%%%%%%%%%%%%%%%-----FIN-----%%%%%%%%%%%%%%%%%%%%%%%%%%%%%%%
ne l'a en fait pas réellement convaincu, 
comme l'a clairement d\'emontr\'e Scholz~\cite{ scho1982b}. 
Il ressort de cette \'etude que Herbart a
beaucoup plus influencé Riemann quant à sa méthodologie de recherche
et à sa compréhension des aspects probl\'ematiques de la
conceptualisation math\'ematique, que pour ce qui concerne
l'élaboration d'une pensée spécifique de l'espace.

En premier lieu, Riemann \'etait vraisemblablement en accord avec
la réfutation par Herbart de la théorie kantienne\footnote{\,
%%%%%%%%%%%%%%%%%%%%%%%-------DEBUT--------%%%%%%%%%%%%%%%%%%%%%%%%%%%
D'après l'{\em Esthétique transcendantale} de Kant (\cite{ kant1986},
pp.~70--72), l'existence des jugements synthétiques {\em a
priori}\,\,---\,\,notamment en géométrie, mais aussi en arithmétique
et en analyse\,\,---\,\,interdit que l'on voie dans l'espace et le
temps des conditions de possibilité des choses mêmes. Quand il s'agit
par exemple
de tirer du concept de trois lignes droites la figure d'un triangle,
se donner des objets dans l'intuition est une condition nécessaire et
{\em sine qua non}. Autrement dit, pour que les démonstrations
synthétiques {\em a priori} de la géométrie soient possibles, il faut
disposer d'un pouvoir d'intuition {\em a priori}. Kant en déduit qu'il
est <<\,indubitablement certain\,>> <<\,que l'espace et le temps, en
tant que conditions nécessaires de toute expérience (extérieure ou
intérieure), ne sont que des conditions simplement subjectives de
notre intuition\,>>.
}, %%%%%%%%%%%%%%%%%%%%%%%%-----FIN-----%%%%%%%%%%%%%%%%%%%%%%%%%%%%%%%
de l'espace et du temps: l'idée de l'espace et du temps comme
<<\,récipients vides\,>> dans lesquels nos sens seraient censés
déverser leurs perceptions était critiquée par Herbart comme une
hypothèse <<\,complètement superficielle, vide de sens, et
inappropriée\,>> \deutsch{völlig gehaltlose, nichtssagende,
unpassende} (\cite{ scho1982b}, p.~422). Aussi Riemann 
(à la suite de Herbart)
considérait-il
peut-être que la question serait <<\,résolue\,>> au sein d'une plus
ample réflexion, grâce à des considérations nouvelles qui montreraient
l'ouverture d'une <<\,béance problématique\,>> beaucoup plus
considérable que ne l'avait soupçonné Kant.

\HEAD{Chapitre~1.\,\,\,\,L'ouverture riemannienne}{
1.10.\,\,\,Le r\'ealisme dialectique mod\'er\'e de Herbart}

\medskip\noindent{\bf 1.10.~Le r\'ealisme dialectique mod\'er\'e de 
Herbart.} 
\label{realisme-dialectique-Herbart}
Herbart distingue autant que possible <<\,le
donné\,>>\,\,---\,\,c'est-\`a-dire les phénomènes, les sensations et
les perceptions\,\,---\,\,du <<\,réel\,>> \deutsch{das Reale} qui
constitue l'objectif principal
de la connaissance philosophique et
scientifique. Le savoir doit être structur\'e en fonction des {\em
connexions conceptuelles possibles} entre les ph\'enom\`enes. Tandis
que l'expérience nous montre des propriétés et des faisceaux
\deutsch{Complexionen} de propriétés, le mouvement de la pens\'ee,
gouvern\'e par les principes de la m\'ethodologie conceptuelle (\cf~ce
qui va suivre), s'effectue par un arc d'embrassement, d'int\'egration
et de structuration de la r\'ealit\'e. Par contraste avec Platon,
Spinoza, ou encore Schelling, Herbart considérait donc qu'il existe
des relations d'homologie précises entre les phénomènes et le r\'eel,
et que le savoir doit proc\'eder de l'exp\'erience donn\'ee vers une
explication pens\'ee des ph\'enom\`enes, {\em via} des clarifications
conceptuelles ad\'equates.

Riemann a été marqué par 
la proposition de Herbart de d\'evelopper les math\'ematiques
<<\,philosophiquement\,>> comme une science pure des concepts;
une telle proposition
s'accordait en fait avec la tendance philosophique de ses
contemporains \`a
\'elaborer des th\'eories syst\'ematiques de la connaissance.
Form\'e dans la tradition de Kant et de Fichte, Herbart admettait
que l'objectif de la philosophie et de la science est de faire avancer
la connaissance en raffinant et en approfondissant le champ de
l'exp\'erience.

Pour cette raison, Herbart attribuait un rôle essentiellement
auxiliaire à la philosophie par rapport aux sciences, sans pour autant
assigner de démarcation nette entre ces deux domaines d'investigation.
L'un de ses buts principaux était d'\'etablir une {\em
méthodologie dialectique modérée} qui incorpore la contradiction comme
moyen transitoire pour l'élucidation de connexions entre concepts.
Toutefois, 
S'\'ecartant de Hegel, Herbart n'\'erigeait pas la contradiction
(et son relèvement, l'\deutschplain{Aufhebung} h\'eg\'elienne) au rang de
mouvement universel de la pensée dans sa totalité propre.

\smallskip

En r\'esum\'e, deux aspects de la philosophie de Herbart montrent une
affinit\'e avec les sciences physiques et math\'ematiques plus
marqu\'ee que chez les id\'ealistes allemands de son temps (Hegel,
Fichte, Schelling):

\smallskip\noindent$\bullet$\,\,
l'assignation d'un rôle essentiellement auxiliaire \`a la philosophie
par rapport \`a l'exploration scientifique du r\'eel;

\smallskip\noindent$\bullet$\,\,
l'\'elaboration d'un r\'ealisme dialectique mod\'er\'e, en opposition 
avec l'id\'ealisme strictement dialectique. 

\smallskip
En quoi Riemann alors a-t-il \'et\'e {\em influenc\'e} par 
Herbart\footnote{\,
%%%%%%%%%%%%%%%%%%%%%%%-------DEBUT--------%%%%%%%%%%%%%%%%%%%%%%%%%%%
Dedekind (\cite{ riem1892}, p.~545) décrivait comme suit la période où
Riemann commençait ses travaux de recherche à Göttingen en 1849:
<<\,Les premiers germes de ses idées en philosophie des sciences ont
dû se développer à cette époque, simultanément à ses occupations dans
les études philosophiques, largement orientées vers Herbart\,>>
\deutsch{In dieser Zeit müssen bei gleichzeitiger 
Beschäftigung mit philosophischen Studien, welche sich natmentlich auf
Herbart richteten, die ersten Keime seiner naturphilosophischen Ideen
sich entwicklelt haben}.
}? %%%%%%%%%%%%%%%%%%%%%%%%-----FIN-----%%%%%%%%%%%%%%%%%%%%%%%%%%%%%%%
Les conceptions philosophiques de Herbart sur les sciences ont \'et\'e
lues et assimil\'ees par Riemann \`a un niveau tr\`es g\'en\'eral,
mais en ce qui concerne la th\'eorie de la connaissance, des
divergences de vue existent, comme le montre une analyse fine des
\deutschplain{Fragmente philosophischen Inhalts}
(\cite{ scho1982b}). Par exemple, en
consid\'erant le d\'eveloppement historique de la connaissance
math\'ematique (toujours structur\'ee en syst\`emes coh\'erents
exempts de contradiction), Riemann ne pouvait pas s'accorder
compl\`etement avec Herbart, pour qui le changement historique des
concepts\,\,---\,\,et notamment dans l'ontologie qui fondait son
r\'ealisme mod\'er\'e\,\,---\,\,s'articule tout au long d'une chaîne
d'erreurs successivement corrigées\footnote{\,
%%%%%%%%%%%%%%%%%%%%%%%-------DEBUT--------%%%%%%%%%%%%%%%%%%%%%%%%%%%
Si ma théorie [de la substance, de l'ontologie] n'est
pas correcte, dit Herbart, alors elle confirme
mon affirmation présente, que les concepts
sont un travail encore inachevé. 
}. %%%%%%%%%%%%%%%%%%%%%%%%-----FIN-----%%%%%%%%%%%%%%%%%%%%%%%%%%%%%%% 
En contrepoint et parce qu'il \'etait math\'ematicien, Riemann
d\'eclarait que les relations d'un savoir nouvellement cr\'e\'e au
sujet d'un domaine de la r\'ealit\'e par rapport \`a un savoir ancien
dans le même domaine ne consistait pas n\'ecessairement en un lien de
correction et de falsification, mais en
un travail de {\em modification} et de
raffinement des structures conceptuelles\footnote{\,
%%%%%%%%%%%%%%%%%%%%%%%-------DEBUT--------%%%%%%%%%%%%%%%%%%%%%%%%%%%
<<\,Die Begriffssysteme, welche ihnen [Naturerkentnissen] jetzt
zu Grunde liegen, sind durch
allmählige Umwandlung älterer Begriffssysteme entstanden, und
die Gründe, welche zu neuen Erklärungensweisen trieben, lassen
sich stets auf Widersprüche oder Unwahrscheinlichkeiten, die
sich in den älteren Erklärungsweisen herausstellten, zurückführen.\,>>
(\cite{ riem1892}, p.~489). 
}. %%%%%%%%%%%%%%%%%%%%%%%%-----FIN-----%%%%%%%%%%%%%%%%%%%%%%%%%%%%%%%

Ainsi d'une mani\`ere g\'en\'erale, c'est principalement sur le plan de
la {\em m\'ethodologie de la genèse des concepts} que Riemann a
\'et\'e marqu\'e par les r\'eflexions de Herbart. 
Un fragment philosophique permet de se faire une première idée des
principes que Riemann a retenus. 

\CITATION{Si 
quelque chose a lieu qui n'était pas attendu selon nos suppositions
préalables, c'est-à-dire qui est impossible ou improbable selon ces
dernières, alors il faut effectuer un travail pour compléter [les
concepts] ou, si nécessaire, retravailler les axiomes afin que ce 
qui est perçu cesse d'être impossible ou improbable. Le complément ou
l'amélioration du système conceptuel forme l'<<\,explication\,>> de la
perception inattendue.
\REFERENCE{\cite{riefu2002},~pp.~20--21.}}

\noindent
Quant à la manière de conduire et d'organiser la recherche
philosophique ou scientifique, deux m\'ethodologies importantes ont
\'et\'e introduites par Herbart:

\smallskip
{\bf 1)} la m\'ethode des relations
\deutsch{Die Methode der Beziehungen}, et:

\smallskip
{\bf 2)} la m\'ethode de la sp\'eculation.

\smallskip\noindent
Riemann en extrait une troisième méthode, qu'il développera de manière
systématique:

\smallskip
{\bf 3)}\, la méthode des plus petits
changements conceptuels. 

\HEAD{Chapitre~1.\,\,\,\,L'ouverture riemannienne}{
1.11.\,\,\,La m\'ethode des relations}

\medskip\noindent{\bf 1.11.~La m\'ethode des relations.}
\label{methode-des-relations}
Aussi bien dans l'expérience donnée que dans l'analyse des concepts
déjà acquis, la contradiction constitue d'après Herbart
la force motrice principale
pour adapter, modifier, transformer,
définir et créer de nouveaux concepts. Par un travail
intellectuel adéquat, toute contradiction observée\,\,---\,\,tant dans
un système conceptuel donné qu'entre la théorie et
l'expérience\,\,---\,\,doit être résolue afin de faire progresser la
connaissance: c'est la {\sl méthode des relations} (dialectiques),
inspirée lointainement de la pensée hégélienne. Tout progrès
s'effectue alors par une {\em transition} 
à partir d'une raison
\deutsch{Grund}, d'une cause, d'un motif, d'un fond, vers
une certaine conclusion qui se situe dans une relation
d'opposition et d'éclaircissement avec
la raison <<\,contredisante\,>>\footnote{\,
%%%%%%%%%%%%%%%%%%%%%%%-------DEBUT--------%%%%%%%%%%%%%%%%%%%%%%%%%%%
<<\,Offenbar fordern wir von dem Grunde, da{\ss}, {\em indem er die
Folge erzeugt}, er selbst sich {\em ändert}. Seine Materie soll sich
verwandeln in die neue Materie der Folge. Hier kann nicht Wahrheit an
Wahrheit geknüpft werden, sondern, damit die Folge Wahrheit enthalte,
mu{\ss} der Grund das Gegentheil davon sein.
(\cite{ scho1982b}, p.~438.) 
}. %%%%%%%%%%%%%%%%%%%%%%%%-----FIN-----%%%%%%%%%%%%%%%%%%%%%%%%%%%%%%%
De fait, la résorption des contradictions débouche sur des
connaissances nouvelles, car la raison \deutsch{Grund}, en tant que
contradiction effective et actualisante, ouvre sur des systèmes
constitués de savoir potentiel auxquels elle n'appartient ni de
manière préalable, ni de manière explicite. Suivant l'aphorisme
(d'inspiration hégélienne) de Herbart: <<\,la raison est
contradiction\,>>, car pour la pensée, la causalité réelle
de l'essence des concepts se révèle dans, 
et seulement par le parcours des obstacles fondamentaux\footnote{\,
%%%%%%%%%%%%%%%%%%%%%%%-------DEBUT--------%%%%%%%%%%%%%%%%%%%%%%%%%%%
Darum sagen wir: {\em der Grund ist ein Widerspruch}. Die 
Schärfe dieser Behauptung abstumpfen, hei{\ss}t, dem
Grunde seine Kraft benehmen.\,>>
(\cite{ scho1982b}, p.~438.) 
}. %%%%%%%%%%%%%%%%%%%%%%%%-----FIN-----%%%%%%%%%%%%%%%%%%%%%%%%%%%%%%%

\`A un niveau
épistémologique supérieur, Herbart envisageait des situations
contextuelles complexes qui donnent naissance à des contradictions
{\em ouvertes}, et donc à une méthodologie ouverte pour l'exploration
de relations ouvertes. Pour que la connaissance progresse, la
procédure méthodologique générale consiste alors à résoudre les
contradictions en élargissant le système conceptuel initial, de
manière à mieux prendre en compte la situation contextuelle
considérée. La notion même de {\em relation}, dans la <<\,méthode des
relations\,>>, manifeste donc la problématicité réelle du {\em
rapport} aux contradictions que les systèmes rencontrent. La relation,
en tant que multiple exposé à l'inattendu, demande, pour faire face à
cet éclatement, la mise au point de méthodes d'investigation. Eu
égard à son intérêt pour l'éducation des esprits, l'une des questions
principales de Herbart était en effet:

\smallskip

\centerline{\em Comment procéder dans la recherche, en
philosophie ou en science?} 

\smallskip\noindent
Question aujourd'hui occultée, tue, bien qu'essentielle et incontournable. 

\HEAD{Chapitre~1.\,\,\,\,L'ouverture riemannienne}{
1.12.\,\,\,La m\'ethode de la sp\'eculation et les graphes de concepts}

\medskip\noindent{\bf 1.12.~La m\'ethode de la sp\'eculation
et les graphes de concepts.}
\label{methode-de-la-speculation}
En approfondissant la méthode des relations, Herbart a dégagé l'idée
générale d'{\sl élaboration de l'unité dans la diversité} et il a
proposé \`a cet effet que dans chaque domaine d'investigation
scientifique soit mis au point un {\sl concept central}
\deutsch{Hauptbegriff}\footnote{\,
%%%%%%%%%%%%%%%%%%%%%%%-------DEBUT--------%%%%%%%%%%%%%%%%%%%%%%%%%%%
Bien entendu, dans le discours d'habilitation de
Riemann, le concept central concerné n'est autre
que celui de {\em multiplicité} 
\deutsch{Mannigfaltigkeit}. 
} %%%%%%%%%%%%%%%%%%%%%%%%-----FIN-----%%%%%%%%%%%%%%%%%%%%%%%%%%%%%%%
 autour duquel s'organisent et
se reflètent les concepts dérivés, subordonn\'es, spécifiés, ou
catégorisés \footnote{\,
%%%%%%%%%%%%%%%%%%%%%%%-------DEBUT--------%%%%%%%%%%%%%%%%%%%%%%%%%%%
Alors que les sciences développent seulement des concepts centraux
reliés à leur domaine spécifique, la philosophie, lorsqu'elle est
entendue comme étude des sciences, se doit de former des concepts
unificateurs qui transcendent les contextes spécifiques.
}. %%%%%%%%%%%%%%%%%%%%%%%%-----FIN-----%%%%%%%%%%%%%%%%%%%%%%%%%%%%%%%

\smallskip

\hspace{-1cm}
\input concepts.pstex_t

\smallskip\noindent
La tâche de la philosophie est alors d'analyser les {\em relations}
entre ces concepts variés et d'examiner leurs caractéristiques
intrinsèques. De Herbart, Riemann a vraisemblablement adopté l'idée
que les liens potentiels entre une multiplicité virtuelle de concepts
satellitaires {\em font et doivent faire question}, d'où les points
d'interrogation dans notre diagramme ci-dessus, que nous appellerons
{\sl graphe problématique de concepts}\footnote{\,
%%%%%%%%%%%%%%%%%%%%%%%-------DEBUT--------%%%%%%%%%%%%%%%%%%%%%%%%%%%
Un tel graphe n'a aucune raison d'être hiérarchisé 
en branches linéaires, ni même d'être représentable
dans un espace à deux ou à trois dimensions. 
}. %%%%%%%%%%%%%%%%%%%%%%%%-----FIN-----%%%%%%%%%%%%%%%%%%%%%%%%%%%%%%%

Hiérarchies, dérivations, liens logiques, tout interroge dans un tel
graphe, tout est à explorer et donc, tout est à construire. Pour
Herbart, ces objectifs d'élucidation sont la matière même de la
philosophie comme savoir particulier
\deutsch{Philosophie als eigene Wissenschaft}. Certes, Herbart
considérait que la différence essentielle entre la philosophie et les
sciences consiste en ce que seules ces dernières ont directement
affaire au <<\,donné\,>>, mais d'après lui néanmoins, les
mathématiques ont une place privilégiée aux côtés de la philosophie,
en tant que toutes deux sont concernées par la création d'une science
de la quantité\footnote{\,
%%%%%%%%%%%%%%%%%%%%%%%-------DEBUT--------%%%%%%%%%%%%%%%%%%%%%%%%%%%
D'après Herbart, lorsqu'elles sont traitées de manière philosophique,
les mathématiques s'intègrent dans la philosophie en général, et au
début de son {\em \'Etude philosophique des sciences} (\cite{
herb1964}), il suggère que le mathématicien ressente le besoin
professionnel de dévoiler l'{\em esprit}
de ses formules pleines d'esprit
\deutsch{Der Mathematiker fühlt den Beruf, uns
den Geist seiner geistreichen Formeln zu enthüllen}. Scholz (\cite{
scho1982b}, p.~426) ajoute qu'il serait difficile d'imaginer une
meilleure caractérisation de la manière dont Riemann fait des
mathématiques.
} %%%%%%%%%%%%%%%%%%%%%%%%-----FIN-----%%%%%%%%%%%%%%%%%%%%%%%%%%%%%%%
\deutsch{Grössenlehre}.
 
En tout cas, pour ce qui concerne l'élaboration des systèmes, le but
ultime est la représentation du <<\,réel\,>>\footnote{\,
%%%%%%%%%%%%%%%%%%%%%%%-------DEBUT--------%%%%%%%%%%%%%%%%%%%%%%%%%%%
L'ontologie de Herbart postulait le <<\,réel\,>>
comme libre de contradictions et de changements. 
} %%%%%%%%%%%%%%%%%%%%%%%%-----FIN-----%%%%%%%%%%%%%%%%%%%%%%%%%%%%%%% 
en son architecture propre, et pour Herbart, la genèse des concepts
s'effectue spécialement grâce à ce qu'il désigne comme la {\sl méthode
de la spéculation}. Par là, il faut entendre <<\,toute entreprise
visant à tracer des {\em transitions} entre (des) concepts\,>>. Et
lorsque Riemann prend des notes au sujet du traité de Herbart {\em Sur
les études philosophiques}, il résume en quelques phrases-clé
(feuillet no.~177 du dossier~18)\footnote{\,
%%%%%%%%%%%%%%%%%%%%%%%-------DEBUT--------%%%%%%%%%%%%%%%%%%%%%%%%%%%
\deutschplain{Philosopie $=$ Untersuchung der Begriffe
\\
I.~Phil(osophische) Ansichten.
\\
II.~Speculation $=$ Streben zur Auflösung 
der Probleme.
\\
Nachweisung eines nothwendigen Zusammenhangs unter Begriffen
\\
Problem weiterer Specul(-ation) als Behmühung zwischen den Begriffen
die gehörigen \"Ubergänge zu bahnen.}
} %%%%%%%%%%%%%%%%%%%%%%%%-----FIN-----%%%%%%%%%%%%%%%%%%%%%%%%%%%%%%%
les idées principales qu'il retient au sujet de cette <<\,méthode de
la spéculation\,>>:

\CITATION{Philosophie 
$=$ recherche de concepts.
\\
I.~Vues philosophiques. 
\\
II.~Spéculation $=$ tendance \deutsch{Streben}
vers la résolution de problème.
\\
Démonstration d'une connexion nécessaire entre
concepts:
\\
problème de spéculation ultérieure, comme effort \deutsch{Bemühung}
pour tracer des transitions adéquates entre les concepts.
\REFERENCE{\cite{scho1982b},~p.~425.}}

\noindent
Ainsi la {\sl méthode de la spéculation} désigne-t-elle tout effort,
application ou travail visant à construire des {\em transitions
appropriées} entre les concepts, aussi bien entre ceux qui sont acquis
et connus, qu'entre ceux qui sont encore à inventer ou
inconnus. Chaque concept-nodal est impliqué dans une relation
problématique avec les concepts de son voisinage immédiat: c'est le
sens des points d'interrogation qui sont inscrits à titre symbolique
dans le diagramme abstrait ci-dessus. Le gain en adéquation augmente
lorsque les flèches du graphe de transition qui se rapportent à un
concept-nodal sont elles-mêmes construites en adéquation.

\smallskip

\centerline{\em la
recherche de l'adéquation porte aussi sur le relationnel.} 

\smallskip\noindent
Et l'exploration des <<\,arêtes relationnelles\,>> du graphe conceptuel
problématique fait parfois naître de nouveaux <<\,nodules
conceptuels\,>> au sein même des transitions qui soulèvent des
questions. La métaphore allusive
d'un <<\,réseau de neurones\,>> du conceptuel
pourrait éclaircir encore un peu plus cette mystérieuse {\em
méthode des relations}. 

Qu'il s'agisse ici de {\em transition} \deutsch{\"Ubergang} est
absolument crucial, car l'objectif est de créer les conditions de {\em
continuité et de complétude} tant pour l'{\em objet concret} que pour
la {\em pensée abstraite}, aussi bien en mathématiques que dans les
sciences de la nature. Et dans l'{\em Essai d'une théorie des
concepts fondamentaux des mathématiques et de la physique comme
fondement pour expliquer la nature}\footnote{\,
%%%%%%%%%%%%%%%%%%%%%%%-------DEBUT--------%%%%%%%%%%%%%%%%%%%%%%%%%%%
\deutschplain{Fragmente philosophischen Inhalts},
\cite{ riem1892}, pp.~477--488; 
\cite{ riefu2002}, pp.~15--20.
}, %%%%%%%%%%%%%%%%%%%%%%%%-----FIN-----%%%%%%%%%%%%%%%%%%%%%%%%%%%%%%%
Riemann exprime encore plus précisément son idée fixe 
au sujet de la {\em continuité}. 

\CITATION{Tout 
ce qui est observé est la transition \deutsch{\"Ubergang} d'une
chose d'un certain état vers un autre état ou, pour parler d'une
manière plus générale, d'un mode de détermination
\deutsch{Bestimmungsweise}
vers un autre, sans que soit perçu un saut au cours de la transition.
Pour compléter cette observation, nous pouvons soit supposer que la
transition se déroule selon un nombre très grand mais fini de sauts
imperceptibles pour nos sens, soit que la chose se déroule de manière
continue à travers {\em toutes} les étapes intermédiaires, en allant
d'un état vers l'autre.
\REFERENCE{\cite{riefu2002},~p.~21.}}

\noindent
Dans ce même passage, Riemann partait de l'hypothèse qu'on a déjà
formé le concept de chose existante en soi
\deutsch{der Begriff für sich bestehender Dinge}, 
et il demandait qu'on <<\,relève la tâche qui consiste à maintenir
aussi loin que possible\,>> ce concept déjà prouvé de chose existante
en soi, sachant que les processus de changement contredisent
manifestement un tel concept de chose stable et existante\footnote{\,
%%%%%%%%%%%%%%%%%%%%%%%-------DEBUT--------%%%%%%%%%%%%%%%%%%%%%%%%%%%
Héraclite {\em versus} Parménide.
}. %%%%%%%%%%%%%%%%%%%%%%%%-----FIN-----%%%%%%%%%%%%%%%%%%%%%%%%%%%%%%%
Autrement dit, dès qu'il s'agit d'assurer la cohérence et la cohésion
du concept d'existence en soi d'une chose soumise à l'étude, {\em
l'antinomie, dans les sciences, entre le discret et le continu devient
source de problèmes scientifiques}. De plus, la question très vaste
suggérée par Riemann débouche aussi sur une étude scientifique de la
psychologie de la perception et elle montre que l'{\em antinomie
problématique entre le discret et le continu, contagieuse, s'insinue
jusque dans les structures infimes de notre appareil perceptif
spatio-temporel pour accentuer encore plus son caractère
antinomique}.

Riemann suggère ainsi que le problème est encore plus complexe qu'il
n'y paraît, puisque nos structures perceptives sont elles-mêmes aussi
inscrites dans l'antinomie {\em physique} entre le discret et le
continu. Lorsque la question des moyens d'appropriation perceptive est
considérée comme question principale, aucune phénoménologie, fût-elle
<<\,transcendantale\,>>, ne peut se dispenser de psychologie
expérimentale. \'Eclatés, les problèmes se distribuent entre 
disciplines, et certains d'entre eux restent désespérément ouverts
malgré les progrès de la connaissance.

Ici à nouveau, Riemann, mathématicien pur, manifeste donc une {\em
conscience philosophique remarquable} des antinomies universelles
auxquelles l'entendement scientifique est exposé. De ces antinomies,
il faut déduire l'ouverture\,\,---\,\,tel est le message intemporel de
Riemann. Ces antinomies ne doivent en aucun cas déboucher sur la
fermeture des questions qu'elles ne parviennent pas à résoudre.

Par exemple, d'après un tel point de vue métaphysique, l'émergence des
théories physiques et astrophysiques va {\em ramifier et complexifier}
le canevas initial des deux premières thèses et antithèses que Kant a
formulées au sein des quatre conflits des idées
transcendantales\footnote{\,
%%%%%%%%%%%%%%%%%%%%%%%-------DEBUT--------%%%%%%%%%%%%%%%%%%%%%%%%%%%
Première antinomie: finitude/infinitude de l'univers dans le temps ou
dans l'espace; deuxième antinomie: discrétude/continuité ou
simplicité/composition de toutes les substances matérielles dans le
monde (\cite{ kant1986}, pp.~338--347). Pour
les questions cosmologiques, \voir~\cite{ szcz1998}. 
}. %%%%%%%%%%%%%%%%%%%%%%%%-----FIN-----%%%%%%%%%%%%%%%%%%%%%%%%%%%%%%%
On pourrait dire que la force philosophique intensive de Riemann
réside dans sa capacité à désigner des problèmes qui sont ouverts dans
leur plus grande généralité, dans un {\em amont} de l'amont, là où
l'inaugural est un initial pur, aux {\em racines mêmes} de l'arbre de
la {\em connaissance problématisante}.

Aussi Riemann ne se range-t-il pas à la <<\,solution\,>> que Kant a
voulu apporter aux antinomies de la raison 
pure\footnote{\,
%%%%%%%%%%%%%%%%%%%%%%%-------DEBUT--------%%%%%%%%%%%%%%%%%%%%%%%%%%%
Dans la septième section du chapitre~II: {\em Les antinomies de la
raison pure}, du Livre~II: {\em Les raisonnements dialectiques de la
raison pure} de la {\em Critique de la raison pure} qui s'intitule
{\em Décision critique du conflit cosmologique de la raison avec
elle-même} (et qui constitue l'un des passages-clés de l'ouvrage),
Kant conclut que toute l'antinomie de la raison pure repose sur la
fausse croyance que quand un phénomène conditionné est donné, la
synthèse qui constitue l'ensemble des conditions empiriques qui
contribuent à sa donation, puisse elle aussi être présentée en acte,
alors qu'en vérité, elle n'est pas présentable et s'éloigne en se
volatilisant dans la {\em régression inaccessible} des conditions
empiriques. Ainsi, l'erreur est-elle de croire que <<\,quand le
conditionné est donné [par exemple, tout objet des sens], la série
entière de toutes ses conditions est aussi donnée\,>> (\cite{
kant1986}, p.~376). <<\,Il ne nous reste pas d'autre moyen de
terminer définitivement la lutte, à la satisfaction des deux parties,
que de les convaincre qu'étant capables de se réfuter si bien
réciproquement, elles se disputent pour rien et qu'un mirage
transcendantal leur a fait voir une réalité là où il ne s'en trouve
pas\,>> (\cite{ kant1986}, p.~378).
}, %%%%%%%%%%%%%%%%%%%%%%%%-----FIN-----%%%%%%%%%%%%%%%%%%%%%%%%%%%%%%%
et qui consiste en quelque sorte à considérer que le {\em conflit
dialectique de la raison avec elle-même} provient d'une illusion
transcendantale. En s'inspirant à tort de l'expérience concrète, la
raison s'imagine en effet inconsidérément\,\,---\,\,d'après le
raisonnement de Kant\,\,---\,\,que le rapport du conditionné à la
série de ses conditions s'applique à des totalités de 
phénomènes\footnote{\,
%%%%%%%%%%%%%%%%%%%%%%%-------DEBUT--------%%%%%%%%%%%%%%%%%%%%%%%%%%%
Mais d'après Kant, ce rapport ne peut
exister en fait que dans et par les
représentations d'entendement. <<\,La majeure du syllogisme
cosmologique prend le conditionné dans le sens transcendantal d'une
catégorie pure et la mineure, dans le sens empirique d'un concept de
l'entendement appliqué à de simples phénomènes, et par conséquent, on
rencontre l'erreur dialectique qu'on nomme {\em sophisma figur{\ae}
dictionis}\,>> (\cite{ kant1986}, p.~377).
}, %%%%%%%%%%%%%%%%%%%%%%%%-----FIN-----%%%%%%%%%%%%%%%%%%%%%%%%%%%%%%%
et la raison en déduit sans plus de façons des conclusions métaphysiques
qui se donnent pour fermes et indubitables, mais qui s'entre-détruisent
parce qu'elles sont violemment antithétiques. Kant fait alors observer
que <<\,deux jugements opposés dialectiquement l'un à l'autre peuvent
être faux tous les deux\,>>, mais la troisième issue logique qu'il
théorise pour sortir de ces difficultés consiste à soutenir que
<<\,les phénomènes en général ne sont rien en dehors de nos
représentations\footnote{\,
%%%%%%%%%%%%%%%%%%%%%%%-------DEBUT--------%%%%%%%%%%%%%%%%%%%%%%%%%%%
Kant déduisait en effet de ces considérations sur les antinomies une
(seconde) <<\,preuve\,>> de l'<<\,idéalité transcendantale\,>> 
des phénomènes (\cite{ kant1986}, p.~381).
}\,>>. %%%%%%%%%%%%%%%%%%%%%%%%-----FIN-----%%%%%%%%%%%%%%%%%%%%%%%%%%%
Réhabilitation du tiers exclu à un niveau transcendantal, donc, mais à
une telle issue, Riemann aurait vraisemblablement pu formuler deux
objections inspirées de ses recherches en mathématiques et en
physique.

Premièrement, à cause de la limitation des moyens expérimentaux, la
régression du conditionné vers ses conditions, loin d'être indéfinie,
s'avère en fait toujours particulièrement limitée, ce qui compromet le
bien-fondé des raisonnements métaphysiques non seulement quant à leur
contenu conclusif (Kant), mais aussi quant à leur portée synthétique
organique\,\,---\,\,toujours trop rudimentaire et non légitimée par
l'existence d'une réalité dominatrice supérieure\footnote{\,
%%%%%%%%%%%%%%%%%%%%%%%-------DEBUT--------%%%%%%%%%%%%%%%%%%%%%%%%%%%
En mathématiques, la complexification incessante et imprévisible de
l'irréversible-synthétique garantit l'organicité du savoir, 
dans l'actuel et dans le potentiel.
}. %%%%%%%%%%%%%%%%%%%%%%%%-----FIN-----%%%%%%%%%%%%%%%%%%%%%%%%%%%%%%%
Pour cette raison, la portée des considérations de la {\em Critique de
la raison pure} est restreinte au champ synthétique ({\em a priori} ou
{\em a posteriori}) d'une époque, et il est impossible d'admettre que
la thétique et l'antithétique de la dialectique transcendantale aient
pu être définitivement établies par le discours kantien, alors que les
antinomies que dégage Kant continuent à se raffiner et à se
complexifier\,\,---\,\,{\em sans s'effacer}\,\,---\,\,grâce aux
progrès de la connaissance scientifique.

Deuxièmement, la représentation de régressions potentielles ne sert
qu'à {\em manifester la présence d'horizons de non-savoir} et à {\em
ouvrir} des domaines neufs pour l'investigation scientifique. Le
tiers-exclu riemannien, c'est l'inconnu devant soi, c'est
l'imprévisible, c'est la résolution inattendue d'un problème, ou mieux
encore, c'est l'ouverture constante de l'irréversible-synthétique.

\HEAD{Chapitre~1.\,\,\,\,L'ouverture riemannienne}{
1.13.\,\,\,La méthode des plus petits changements conceptuels}

\medskip\noindent{\bf 1.13.~La méthode des plus petits changements
conceptuels.}
\label{methode-plus-petits-changements}
Dans le feuillet manuscrit 141 du dossier~18 du
\deutschplain{Nachlass} que 
Scholz a découvert (\cite{ scho1982b}, p.~420 et p.~433), Riemann
s'interroge tout d'abord sur les étapes de la formation des concepts:
à partir du perçu; par induction; par abstraction, ou par synthèse
{\em a posteriori}; et ensuite surtout, il insiste sur le {\em
resserrement des transitions} entre les concepts:

\smallskip
\fbox{
\CITATION{Changement 
(ou complétion) des concepts qui soit
le plus petit possible [\dots]\footnotemark
\REFERENCE{\cite{scho1982b},~p.~433.}}}
\medskip

\noindent
Très certainement, cette pensée est directement inspirée d'une formule
aphoristique de Herbart que Riemann recopie en l'extrayant de son
contexte, et que nous encadrons aussi, vu son importance.

\footnotetext{\,
%%%%%%%%%%%%%%%%%%%%%%%-------DEBUT--------%%%%%%%%%%%%%%%%%%%%%%%%%%%
\deutschplain{Möglichst geringe Veränderung (oder Ergänzung) der
Begriffe.}
} %%%%%%%%%%%%%%%%%%%%%%%%-----FIN-----%%%%%%%%%%%%%%%%%%%%%%%%%%%%%%%

\smallskip

\centerline{\fbox{\em La m\'ethode des relations est la m\'ethode des plus
petits changements\footnotemark.}}

\footnotetext{\,
%%%%%%%%%%%%%%%%%%%%%%%-------DEBUT--------%%%%%%%%%%%%%%%%%%%%%%%%%%%
\deutschplain{Die Methode der Beziehungen ist die Methode der 
kleinsten Veränderungen.}
} %%%%%%%%%%%%%%%%%%%%%%%%-----FIN-----%%%%%%%%%%%%%%%%%%%%%%%%%%%%%%%

\smallskip\noindent
Resserrer les transitions conceptuelles, c'est tendre vers 
une plus grande {\em
continuité de la pensée}\footnote{\,
%%%%%%%%%%%%%%%%%%%%%%%-------DEBUT--------%%%%%%%%%%%%%%%%%%%%%%%%%%%
\`A la fin de la première partie des
\deutschplain{Fragmente philosophischen Inhalts}
intitulée \deutschplain{Zur Psychologie und Metaphysik}, Riemann
évoque la méthode des limites introduite par Newton pour fonder le
calcul infinitésimal, et il exprime spontanément une extension (par
analogie) de cette méthode à la recherche d'une {\em continuité} dans
la détermination\,\,---\,\,primitivement {\em
discontinue}\,\,---\,\,des concepts. La méthode de Newton, écrit
Riemann, <<\,consiste en ceci: au lieu de considérer une transition
continue d'une valeur d'une grandeur à une autre, d'une position à une
autre, ou plus généralement {\em d'un mode de détermination d'un
concept
\deutsch{von einer Bestimmungsweise eines Begriffes} à un autre}, on
considère d'abord une transition par un nombre fini d'étapes
intermédiaires, puis on permet au nombre de ces degrés intermédiaires
d'augmenter, de telle sorte que les distances entre deux degrés
intermédiaires consécutifs diminuent toutes à l'infini\,>> (\cite{
riem1892}, p.~487; nous soulignons). Cette <<\,continuification\,>>
des concepts qui se trouve\,\,---\,\,comme pour le calcul
infinitésimal\,\,---\,\,à la limite du représentable, n'est pas
directement accessible à notre réflexion, mais on peut s'imaginer
qu'il soit possibler de passer d'un système conceptuel à un autre par
une simple modification dans les grandeurs relatives, de telle sorte
que le système reste stable et inchangé dans la transition vers la
limite du continu.
}. %%%%%%%%%%%%%%%%%%%%%%%%-----FIN-----%%%%%%%%%%%%%%%%%%%%%%%%%%%%%%%
L' objectif affirmé est de {\em compléter par élimination} toutes les
différences, distorsions, sauts, fossés,
\etc Dans leur exposition
très systématique de la théorie des groupes de transformations, Engel
et Lie mettront à exécution ce principe fondamental que Riemann
n'avait pas eu le temps de développer.

\`A un niveau <<\,méta-mathématique\,>>, 
ce qui est désigné par Riemann dans la généralité la plus grande comme
{\em transition (continue) entre concepts} est susceptible de se
réaliser comme connexion, comme hiérarchisation, comme sériation,
comme interdépendance, comme nécessité, comme suffisance, ou enfin,
terme ultime du resserrement de la continuité conceptuelle, comme {\em
équivalence}, c'est-à-dire comme {\em nécessité et suffisance}. Tous
les liens dynamiques entre concepts, axiomes et hypothèses recouvrent
ainsi des distinctions en extension et en généralité qui les
inscrivent dans un graphe de transitions ouvertes sur un horizon de
continuité potentielles.

Ainsi s'exprime l'<<\,impératif catégorique\,>> riemannien qui fait
naître une {\em diversité éclatée de problèmes d'organisation
structurale} dont l'objecif est de créer des conditions pour la
continuité de la pensée spéculative:

\begin{itemize}

\smallskip\item[$\bullet$]
nature, conception, problématicité des notions primitives; 

\smallskip\item[$\bullet$]
liaisons, dépendances, indépendances entre concepts définis; 

\smallskip\item[$\bullet$]
nécessités, potentialités, {\em aprioricité};

\smallskip\item[$\bullet$]
principes de genèse conceptuelle; 

\smallskip\item[$\bullet$]
différenciations spécifiques. 

\end{itemize}\smallskip

\noindent
Immédiatement après avoir ouvert les questions de nature, Riemann ose
donc {\em ouvrir} le champ {\em encore inexistant} des relations
possibles entre concepts potentiels. On peut sans conteste affirmer
que la philosophie mathématique de Riemann anticipe la révolution
axiomatique\footnote{\,
%%%%%%%%%%%%%%%%%%%%%%%-------DEBUT--------%%%%%%%%%%%%%%%%%%%%%%%%%%%
Néanmoins, la pensée de Riemann ne se réduit pas à un procédé
entièrement logique, voire purement axiomatique (\cf~le
commentaire~\cite{ vui1962}, p.~408 de Vuillemin basé sur une lecture
de Russell), puisque la logique formelle ne questionne pas au sens
philosophique du terme, mais cherche à enfermer le relationnel du
conceptuel dans une syntaxe constituée {\em a posteriori} par rapport
aux recherches ouvertes. Il est beaucoup plus exact de dire que les
recherches de Riemann sur la définition de l'espace procèdent
<<\,comme il est naturel dans l'analyse des concepts\,>>, à savoir:
<<\,{\em per genus proximum et differentiam specificam}\,>> ({\em
ibidem}), sans toutefois voir dans les déterminations successives de
concepts des finalités prédéfinies et fermées, comme par exemple la
<<\,convergence\,>> vers un système d'axiomes caractérisant l'espace
euclidien. Très peu de commentateurs ont mis en lumière le souci
constant de problématisation et d'ouverture qui s'exprime de manière
très explicite dans tous les écrits de Riemann. }
%%%%%%%%%%%%%%%%%%%%%%%%-----FIN-----%%%%%%%%%%%%%%%%%%%%%%%%%%%%%%%
(Hilbert, Bourbaki), puisque le structural formalisé est préformé par
certaines structures dialectiques universelles du questionnement
mathématique. D'après un tel point de vue, l'axiomatisation ultérieure
des géométries est une tentative\,\,---\,\,parmi d'autres\,\,---\,\,de
systématiser certaines réponses {\em et certaines questions} que
l'exploration rencontre. Décider les ouvertures, accepter
l'ignorance, et maintenir intentionnellement le rapport à l'Inconnu
n'est peut-être pas la qualité la plus évidente des systèmes
formalisés, et pourtant, l'<<\,impératif catégorique riemannien\,>>
exigerait que lesdits systèmes expriment explicitement leur propre
inachèvement.

\HEAD{Chapitre~1.\,\,\,\,L'ouverture riemannienne}{
1.14.\,\,\,Modes amétriques de détermination}

\medskip\noindent{\bf 1.14.~Modes amétriques de détermination.}
\label{modes-ametriques}
\label{reprise-analyse}
Maintenant que nous avons élucidé les exigences méthodologiques
générales que Riemann s'est fixées au contact de la philosophie de
Herbart pour conduire ses recherches sur les fondements de la
géométrie, nous pouvons reprendre à présent l'analyse philosophique de
son discours, que nous avions interrompue p.~\pageref{wait-zenon}.

En partant, comme nous l'avons rappelé, de la bifurcation zénonienne
fondamentale entre le discret et le continu, Riemann insiste sur le
fait que les occasions de faire naître des concepts dont les modes de
détermination recourent à la continuité\footnote{\,
%%%%%%%%%%%%%%%%%%%%%%%-------DEBUT--------%%%%%%%%%%%%%%%%%%%%%%%%%%%
\`A nouveau, Riemann se révèle penseur du continu, 
en analyse, en géométrie différentielle et en physique. Weyl écrivait
(\cite{ riem1990}, p.~740) que la motivation de principe de Riemann
était de comprendre le monde par son comportement dans l'infiniment
petit \deutsch{die Welt aus ihrem Verhalten
im Unendlichkleinen zu verstehen}.
} %%%%%%%%%%%%%%%%%%%%%%%%-----FIN-----%%%%%%%%%%%%%%%%%%%%%%%%%%%%%%%
sont beaucoup moins fréquentes\footnote{\,
%%%%%%%%%%%%%%%%%%%%%%%-------DEBUT--------%%%%%%%%%%%%%%%%%%%%%%%%%%%
Riemann affirme même qu'elles sont plus rares \deutsch{selten}, or
cela n'est pas tout à fait exact, car les mouvements des corps dans
l'espace\,\,---\,\,cette réalité intuitive qui nous est omniprésente
et qui allait connaître un destin algébrique inattendu avec le
développement de la théorie des groupes de Lie\,\,---\,\,prouvent
manifestement le contraire. D'ailleurs, Herbart avait déjà montré
qu'il existe des espaces continus dont les modes d'existence sont
extrêmement variés, comme par exemple (à nouveau) les divers lieux que
peuvent occuper les objets sensibles, mais aussi la <<\,ligne du
son\,>> \deutsch{Tonlinie}, ou encore le triangle des couleurs, étudié
par Thomas Young et Clark Maxwell, avec le bleu, le rouge et le jaune
disposés à ses sommets, ces trois couleurs pouvant se fondre en
s'associant quantitativement pour produire toutes les couleurs
possibles dans le continu bidimensionnel de l'intérieur du
triangle. Au début de ce \S~I (\cite{ riem1898}, p.~282), Riemann a
probablement plutôt voulu suggérer que les modes de détermination
continue qui sont nécessaires pour penser les multiplicités continues
ne sont pas encore disponibles, car le besoin de les élaborer ne s'est
pas encore fait ressentir dans la vie ordinaire. }
%%%%%%%%%%%%%%%%%%%%%%%%-----FIN-----%%%%%%%%%%%%%%%%%%%%%%%%%%%%%%%
que lorsqu'il s'agit des concepts dont les modes de détermination
forment une multiplicité discrète. En effet, l'équivalence numérique
entre plusieurs collections d'objets concrets qui fonde, à un niveau
intuitif proprement archaïque la notion élémentaire de nombre entier,
ne pose quasiment aucun problème d'abstraction à la pensée, parce que
les bijections entre collections d'objets individués possèdent, pour
l'intuition, un sens physique extrêmement clair, du moins lorsqu'il
s'agit de petits nombres entiers. Riemann dévoile alors
sa motivation mathématique. 

\CITATION{De 
telles recherches sont devenues nécessaires dans plusieurs parties
des Mathématiques, notamment pour l'étude des fonctions analytiques à
plusieurs valeurs\footnotemark, et c'est surtout à cause de leur
imperfection que le célèbre théorème d'Abel, ainsi que les travaux de
Lagrange, de Pfaff, de Jacobi sur la théorie générale des équations
différentielles, sont restés si longtemps stériles.
\REFERENCE{\cite{riem1898},~p.~283.}}

\footnotetext{\,
%%%%%%%%%%%%%%%%%%%%%%%-------DEBUT--------%%%%%%%%%%%%%%%%%%%%%%%%%%%
Il s'agit bien sûr des surfaces étalées au-dessus de certaines régions
du plan complexe, ramifiées autour de certains points, et recousues le
long de certaines coupures que Riemann a introduites dans sa
dissertation inaugurale. Pour ne pas interrompre l'étude proprement
philosophique, nous laisserons de côté ces brèves allusions à
l'analyse complexe, au théorème d'Abel concernant les intégrales de
fonctions algébriques (\voir~\cite{ houz2002}), et au problème de
l'intégration des formes différentielles totales
(\voir~\cite{ ha2005}).
} %%%%%%%%%%%%%%%%%%%%%%%%-----FIN-----%%%%%%%%%%%%%%%%%%%%%%%%%%%%%%%

\noindent
Cette remarque réfère à une tendance naissante des mathématiques de la
première moitié du 19\textsuperscript{ème} siècle: essayer de
transférer le langage de la géométrie vers des systèmes analytiques ou
algébriques à {\em plusieurs variables}; 
une telle tendance était connue de
Riemann, au moins partiellement, {\em via} Gauss (\cite{ scho1980},
p.~15~sq., 53~sq.). La recherche de généralité exige de mettre sur
pied une théorie des grandeurs étendues qui soit indépendante des
déterminations métriques et dans laquelle <<\,on ne suppose {\em rien de
plus}\footnote{\,
%%%%%%%%%%%%%%%%%%%%%%%-------DEBUT--------%%%%%%%%%%%%%%%%%%%%%%%%%%%
L'exigence ainsi fondée par Riemann de se limiter
à des hypothèses minimales s'exprime aussi très explicitement
dans son \deutschplain{Habilitationsschrift} 
sur les séries trigonométriques. 
} %%%%%%%%%%%%%%%%%%%%%%%%-----FIN-----%%%%%%%%%%%%%%%%%%%%%%%%%%%%%%%
que ce qui est déjà renfermé dans le concept de ces
grandeurs\,>>.

\'Eliminer les hypothèses adventices, circonscrire des
hypothèses minimales, et s'en tenir rigoureusement à elles: impératif
méthodologique
riemannien; principe riemannien de genèse.

\CITATION{Dans 
cette branche générale de la théorie des grandeurs étendues, {\em
où l'on ne suppose rien de plus que ce qui est déjà renfermé dans le
concept de ces grandeurs}, il nous suffira, pour notre
objet actuel, de porter notre étude sur deux points, relatifs: le
premier, à la génération du concept d'une
multiplicité de plusieurs
dimensions; le second, au moyen de ramener les déterminations de lieu
dans une multiplicité donnée à des déterminations de quantité, et c'est ce
dernier point qui doit clairement faire ressortir le caractère
essentiel d'une étude de $n$ dimensions.
\REFERENCE{\cite{riem1898},~p.~283.}}

\noindent
Problème inaugural: il s'agit donc de penser, de
définir et d'engendrer une notion purement {\em amétrique}\footnote{\,
%%%%%%%%%%%%%%%%%%%%%%%-------DEBUT--------%%%%%%%%%%%%%%%%%%%%%%%%%%%
Si les moyens mathématiques ou axiomatiques manquent au géomètre pour
mesurer les grandeurs spatiales, et s'il est impossible de se déplacer
à la manière d'un arpenteur pour déposer chaînes, jalons et équerres
sur un hypothétique sol de la géométrie spatiale abstraite, alors la
question du plus grand et du plus petit perd son sens et son intérêt.
} %%%%%%%%%%%%%%%%%%%%%%%%-----FIN-----%%%%%%%%%%%%%%%%%%%%%%%%%%%%%%%
d'étendue, dont les <<\,quanta\,>> ou <<\,morceaux\,>> d'étendue
puissent être envisagés du point de vue de la comparaison inclusive,
sans pour autant l'être du point de vue de la mesure 
numérique\footnote{\,
%%%%%%%%%%%%%%%%%%%%%%%-------DEBUT--------%%%%%%%%%%%%%%%%%%%%%%%%%%%
Par son traitement détaillé des multiplicités continues qui précède
l'introduction de coordonnées numériques, Riemann affirmait clairement
que les nombres jouent un rôle auxilaire en géométrie, s'exprimant de
ce fait à contre-courant de la tendance générale à l'arithmétisation
de l'analyse et de la géométrie (\cite{ scho1980}, pp.~30--34).
}. %%%%%%%%%%%%%%%%%%%%%%%%-----FIN-----%%%%%%%%%%%%%%%%%%%%%%%%%%%%%%%
Il faut ainsi libérer la notion d'étendue de l'emprise des grilles et
des règles graduées, afin de penser l'étendue en termes de {\em
régions}\footnote{\,
%%%%%%%%%%%%%%%%%%%%%%%-------DEBUT--------%%%%%%%%%%%%%%%%%%%%%%%%%%%
Acte de saisie intuitive archaïque du spatial: comment le
caractériser? Les travaux de Sophus Lie déploieront une pensée du
spatial multidimensionel amétrique et mobile dans un langage
essentiellement rhétorique et algébrique qui cherche néanmoins à
transmettre chaque acte de pensée intuitive. }
%%%%%%%%%%%%%%%%%%%%%%%%-----FIN-----%%%%%%%%%%%%%%%%%%%%%%%%%%%%%%%
situées dans une certaine multiplicité \deutsch{als Gebiete in einer
Mannigfaltigkeit}.

Immédiatement et explicitement, Riemann divise donc le problème
concerné en deux sous-problèmes:
\label{deux-sous-problemes}

\smallskip\noindent
{\sf\em Problème~1:}
engendrer le concept de multiplicité par synthèse et par analyse, 
en ayant recours à des modes de détermination quantitifs non 
métriques. 

\smallskip\noindent
{\sf\em Problème~2:}
trouver des conditions pour ramener les déterminations de lieu
à des rapports métriques de type distance.

\smallskip\noindent
\'Evidemment, 
l'étude de ces deux questions générales devra toujours être encadrée
par les mêmes exigences <<\,méta-mathématiques universelles\,>> qui
animent constamment Riemann:

\smallskip$\bullet$\,\,
donner des conditions {\em suffisantes} pour la genèse;

\smallskip$\bullet$\,\,
donner des conditions {\em nécessaires} pour la genèse;

\smallskip$\bullet$\,\,
\'etablir la {\em nécessité} et la {\em suffisance} de 
conditions génétiques;

\smallskip$\bullet$\,\,
\'eliminer les caractères spécifiques étrangers au concept pur.

\smallskip\noindent
Pour accéder vraiment à une telle pensée génétique, il nous faut
maintenant mettre complètement entre parenthèses les actes de
postulation axiomatique par lesquels nous nous permettons aujourd'hui
d'accepter sans nous poser plus de questions aussi bien les
définitions formelles initiales structurées que les entames telles que
<<\,Soit $M$ une variété différentiable de dimension $n$\,>>, ou
<<\,Soient $(p_1, \dots, p_n, \, q_1, \dots, q_n)$ des coordonnées
locales sur l'espace des phases $\mathcal{ P}$ d'un système
hamiltonien quelconque\,>>. Riemann cherche en fait à {\em démontrer}
que l'on peut ramener la détermination locale des étendues
géométriques à des suites de quantités numériques:
%%%%%%%%%%%%%%%%%%%%%%%%-----FIN-----%%%%%%%%%%%%%%%%%%%%%%%%%%%%%%%
c'est un pur problème de genèse, qui
ne préoccupe plus notre
époque.

\HEAD{Chapitre~1.\,\,\,\,L'ouverture riemannienne}{
1.15.\,\,\,Modes d'engendrement du multidimensionnel}

\medskip\noindent{\bf 1.15.~Genèse du multidimensionnel.}
\label{genese-du-multidimensionnel}
Dans un tout premier moment, Riemann tente de caractériser
l'unidimensionnalité initiale d'une multiplicité non
discrète\,\,---\,\,dont il fera ensuite un principe fondamental de
genèse\,\,---\,\,par la propriété que ses modes continus de détermination
sont eux-mêmes déterminés\footnote{\,
%%%%%%%%%%%%%%%%%%%%%%%-------DEBUT--------%%%%%%%%%%%%%%%%%%%%%%%%%%%
Cette phrase longue et difficile (\cite{riem1898}, pp.~283--284)
semble être une tentative inaboutie pour trouver dans la variabilité
des modes de détermination d'une multiplicité donnée des conditions si
restrictives qu'elles en impliquent l'unidimensionnalité. Le caractère
inachevé de ce passage n'a pas échappé à Engel et à Lie, qui étaient
probablement déjà informés, en 1891--93, des travaux naissants de
Pasch, Stolz, Schur sur les fondements axiomatiques de la géométrie.
<<\,La véritable signification de la proposition d'après laquelle
l'espace est une variété numérique
\deutsch{Zählenmannigfaltigkeit} ne ressort pas du travail de
Riemann. Riemann cherche à démontrer cette proposition, mais sa
démonstration ne peut pas être prise au sérieux. Si l'on veut
véritablement démontrer que l'espace est une variété numérique, on
devra, à n'en pas douter, postuler auparavant un nombre non
négligeable d'axiomes, ce dont il semble que Riemann n'ait pas été
conscient\,>> (p.~\pageref{394} ci-dessous). Il est vrai en effet que
l'élaboration d'un raisonnement véritablement synthétique nécessite de
faire une différence marquée entre hypothèses et conclusion, 
et de s'interroger sur la nature des hypothèses qu'on prendra
comme axiomes. Comme
Riemann ne spécifie pas ce qu'il faut précisément entendre par l'idée
de lieu comme essence du géométral-local-continu, les raisonnements
logiques et les démonstrations qu'il cherche à conduire pour ramener
une telle notion de lieu à des grandeurs numériques sont donc encore
truffés de problèmes ouverts.
}. %%%%%%%%%%%%%%%%%%%%%%%%-----FIN-----%%%%%%%%%%%%%%%%%%%%%%%%%%%%%%%
D'après un raisonnement difficile à reconstituer, Riemann affirme
alors que de tels modes de détermination ne peuvent alors être
parcourus <<\,que dans un seul sens\,>>, c'est-à-dire <<\,en avant et
en arrière\,>>: c'est l'unidimensionnalité, appelée à devenir racine
et principe de genèse inductive pour la multidimensionnalité.

Ensuite, {\sl dyade}, 
{\sl soudure} et {\sl parcours} permettent d'engendrer le
bidimensionnel\footnote{\,
%%%%%%%%%%%%%%%%%%%%%%%-------DEBUT--------%%%%%%%%%%%%%%%%%%%%%%%%%%%
Pour de plus amples développements de
cette direction de pensée, le lecteur est renvoyé à la {\em
Science de la grandeur extensive} de Grassmann~\cite{ gras1994},
commentée par Flament~\cite{ flam1992} et par Gilles 
Châtelet~\cite{ chat1993}.
}. %%%%%%%%%%%%%%%%%%%%%%%%-----FIN-----%%%%%%%%%%%%%%%%%%%%%%%%%%%%%%%
{\sl Dyade}: dédoubler l'objet {\em un} qui est donné, c'est-à-dire se
donner deux multiplicités unidimensionnelles; {\sl soudure}:
transporter une multiplicité pour la <<\,riveter\,>> sur une autre en
un bipoint de coïncidence; {\sl parcours}: faire décrire à la seconde
multiplicité unidimensionnelle toute la multiplicité unidimensionnelle
de la première. Genèse: voir apparaître le bidimensionnel comme un
voile créé aux franges de l'unidimensionnel par développement
continu dans un éther extrinsèque.

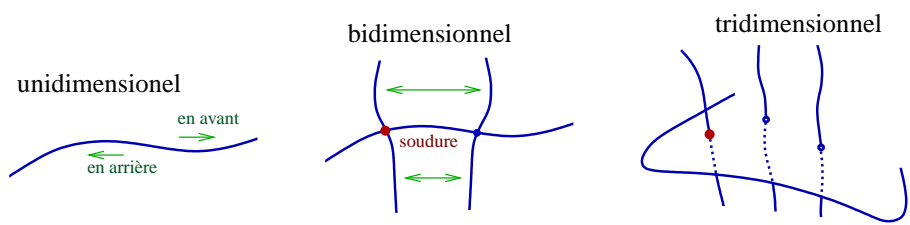
\begin{center}
\input genese-multiplicite.pstex_t
\end{center}

\noindent
Le tridimensionnel s'engendre alors de manière analogue par
déploiement extériorisé de l'unidimensionnel le long d'un bidimensionnel,
et ainsi de suite pour les dimensionnalités d'ordre supérieur. Autre
interprétation de ce procédé: la variation de la variabilité produit une
variabilité d'ordre supérieur.

\CITATION{Si, 
au lieu de considérer le concept comme déterminable, on considère
son objet comme variable, on pourra désigner cette construction comme
la composition d'une variabilité de $n+1$ dimensions, au moyen d'une
variabilité de $n$ dimensions et d'une variabilité d'une seule
dimension.
\REFERENCE{\cite{riem1898},~p.~284.}}

\noindent
Et maintenant, le resserrement des conditions: il faut à présent
s'interroger pour savoir si une {\em analyse inverse} est possible,
car l'existence de deux genèses réciproques montrerait que le concept
de multiplicité a été entièrement circonscrit.
\`A un procédé de composition par induction doit donc
succèder une analyse par décomposition dimensionnelle.

\CITATION{Je 
vais maintenant montrer réciproquement comment une variabilité,
dont le champ est donné, peut se décomposer en une variabilité d'une
dimension et une variabilité d'un nombre de dimensions moindre.
\REFERENCE{\cite{riem1898},~p.~284.}}

\noindent
\`A cet instant précis va se dévoiler pour
la première fois dans le texte de Riemann {\em l'interdépendance
fondamentale entre le fonctionnel et le géométral}. C'est le premier
moment où, à dessein, {\em l'analyse et l'algèbre commencent à
vouloir capturer la géométrie}. 

L'idée est simple: dans une portion de multiplicité, les éléments
divers de la multiplicité comptés à partir d'un point fixé à l'avance
doivent posséder un principe de différenciation spécifique par rapport
au point en question. Autrement dit, on peut s'imaginer qu'à
l'intérieur de ladite multiplicité, il existe au moins une certaine
{\em fonction du lieu} qui ne soit constante le long d'aucune
sous-portion de cette multiplicité. Ici, il ne s'agit pas de définir
axiomatiquement le géométral et le fonctionnel à partir d'atlas
maximaux de cartes locales à valeurs dans un ouvert de $\R^n$ (\cite{
spiv1970, scho1980, ol1995}). Il s'agit plutôt d'observer que le
fonctionnel naît d'emblée avec le géométral, par l'effet d'une dualité
ou d'une complémentarité qui se trouve à la racine des concepts.

\`A ce moment-là, bien que Riemann soit pertinemment
conscient\,\,---\,\,grâce à ses travaux sur les séries
trigonométriques\,\,---\,\,du fait que les fonctions diffèrent en
nature suivant qu'elles sont continues, différentiables, ou
analytiques, avec un ensemble éventuellement fini ou infini de
discontinuités, il semble vouloir ne préciser ici aucune hypothèse
technique au sujet de la régularité de la fonction en question. Par
conséquent, le fonctionnel est ici absolument ouvert à la généralité et à
la diversification des univers spatiaux. On évolue donc dans un monde
virtuel qui embrasse {\em a priori} les variétés topologiques, les
variétés de classe $\mathcal{ C}^1$, $\mathcal{ C}^2$, $\mathcal{
C}^k$ ou $\mathcal{ C}^\infty$, les variétés analytiques réelles,
complexes, ou quaternioniques, les espaces analytiques singuliers, les
espaces éventuellement fractals, non séparés, totalement discontinus,
ou encore absolument mixtes, c'est-à-dire qui incorporent
éventuellement, en des localités distinctes, chacun de ces aspects-là:
postérité stupéfiante de la généralité riemannienne. Les raisonnements
qui pourraient être considérés comme vagues et imprécis englobent
donc ici des pans entiers de la géométrie à plusieurs dimensions, que
Engel et Lie allaient être les premiers à développer d'une manière
vraiment systématique, dans une optique exclusivement locale et
générique.

Lorsqu'on fait varier la constante à laquelle on égale une telle
fonction, les lieux de points en lesquels sa valeur est fixe forment
alors une multiplicité continue\footnote{\,
%%%%%%%%%%%%%%%%%%%%%%%-------DEBUT--------%%%%%%%%%%%%%%%%%%%%%%%%%%%
En toute rigueur ici, il faut faire des hypothèses telles que par
exemple la différentiabilité au moins $\mathcal{ C}^1$ et la
non-annulation de la différentielle, puisque d'après un théorème de
Whitney, tout sous-ensemble fermé d'une variété, aussi pathologique
qu'il soit, peut être représenté comme lieu d'annulation d'une
certaine fonction continue (\cite{ malg1967}). Toutefois, dans ce
moment d'analyse, les raisonnements sont locaux et génériques: <<\,Les
cas d'exception, dont l'étude est importante [souci d'ouverture],
peuvent être ici laissés de côté\,>>. Rien n'empêche en tout cas de
pressentir que de tels raisonnements conservent un sens très précis
dans la catégorie des espaces analytiques complexes singuliers, où
l'ontologie parallèle du fonctionnel et du géométral, mieux contrôlée
par les séries entières convergentes, se prolonge toujours d'un point
vers un petit voisinage de ce point, grâce notamment au théorème de
préparation de Weierstrass, au théorème de paramétrisation locale de
Noether et au théorème de cohérence d'Oka.
} %%%%%%%%%%%%%%%%%%%%%%%%-----FIN-----%%%%%%%%%%%%%%%%%%%%%%%%%%%%%%%
d'un nombre de dimensions moindre que celui de la variété donnée.

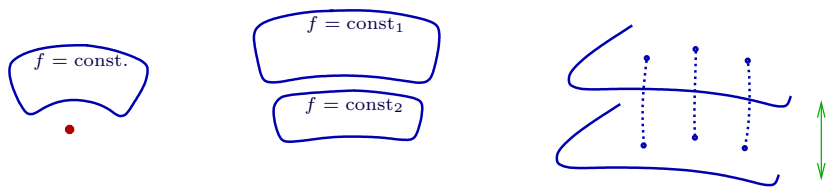
\begin{center}
\input f-cst.pstex_t
\end{center}

\CITATION{Ces
multiplicités, 
lorsqu'on fait varier la fonction, se transforment d'une
manière continue les unes dans les autres; on pourra donc admettre que
l'une d'entre elles engendre les autres, et cela pourra avoir lieu,
généralement parlant, de telle façon que chaque point de l'une se
transporte en un point déterminé de l'autre. 
\REFERENCE{\cite{riem1898},~pp.~284--285.}}

\noindent
Ainsi par décomposition, Riemann parvient-il en utilisant des
fonctions auxilaires, à montrer que l'analyse du concept de
multiplicité reconduit exactement au premier procédé de genèse:
amplification de la variabilité par variation inductive de la
variabilité.

En conclusion, que ce soit par l'analyse ou par la synthèse, <<\,la
détermination de lieu dans une multiplicité donnée, quand cela est
possible, se réduit à un nombre fini de déterminations de
quantité\footnote{\,
%%%%%%%%%%%%%%%%%%%%%%%-------DEBUT--------%%%%%%%%%%%%%%%%%%%%%%%%%%%
Engel et Lie appelleront <<\,variétés numériques\,>>
\deutsch{Zählenmannigfaltigkeit}
les multiplicités introduites par Riemann. 
}\,>>: %%%%%%%%%%%%%%%%%%%%%%%%-----FIN-----%%%%%%%%%%%%%%%%%%%%%%%%%%%%
ce qu'il fallait démontrer.

Contrairement aux définitions axiomatiques, Riemann a donc entrepris
d'engendrer le multidimensionnel à partir de l'unidimensionnel, et
surtout il a tenté d'articuler la genèse comme une {\em démonstration}
procédant par des conditions variées. Helmholtz quant à lui, puis
Russell (\cite{ russ1956}) et Vuillemin (\cite{ vui1962},
pp.~388--464), envisageront le problème sous l'angle de la postulation
d'axiomes, pourtant moins problématisant, et moins riche
d'imprévisible. Pour terminer sur ce chapitre, notons à nouveau que
Riemann ne peut se soustraire à un devoir constant de manifester le
souci abstrait d'ouverture.

\CITATION{Toutefois 
il y a aussi des multiplicités dans lesquelles la
détermination de lieu exige, non plus un nombre fini, mais soit une
série infinie, soit une multiplicité continue de déterminations de
grandeur. Telles sont, par exemple les multiplicités formées par les
déterminations possibles d'une fonction dans une région donnée, par
les formes possibles d'une figure de l'espace, \etc
\REFERENCE{\cite{riem1898},~p.~285.}}

\HEAD{Chapitre~1.\,\,\,\,L'ouverture riemannienne}{
1.16.\,\,\,Conditions pour la détermination des rapports métriques}

\medskip\noindent{\bf 1.16.~Conditions 
pour la détermination des rapports métriques.}
\label{conditions-suffisantes-metrique}
Après avoir libéré la notion archaïque d'étendue de toute saisie
métrique, il s'agit maintenant de traiter le deuxième sous-problème
(\cf~p.~\pageref{deux-sous-problemes}), que Riemann reformule comme
suit:

\smallskip

\centerline{\em De quels types de rapports métriques est
susceptible une multiplicité?}

\smallskip\noindent
Il s'agit de trouver des conditions qui caractérisent les
différentes manières possibles de munir une multiplicité donnée d'un
concept supplémentaire qui permette de parler de la {\em distance} qui
existe entre toutes les paires de points. Or le titre passablement
énigmatique de la sous-section concernée: 

\smallskip

{\em Rapports métriques dont est
susceptible une variété de $n$ dimensions, dans l'hypothèse où les
lignes possèdent une longueur, indépendamment de leur position, et où
toute ligne est ainsi mesurable par toute autre ligne} 

\smallskip\noindent 
a plongé dans la perplexité de nombreux mathématiciens, philosophes,
historiens et commentateurs. Tout d'abord, que signifie exactement
cette seconde hypothèse elliptique d'après laquelle toute ligne doit
être mesurable par toute autre ligne? Faut-il y voir un principe
d'arpentage: déplacement libre des règles (ou des lignes)? Doit-on en
faire toujours un axiome, en vertu d'une évidence physique imparable?
Par <<\,lignes\,>>, faut-il entendre n'importe quelle courbe tracée
dans la multiplicité? Un tel principe de métrisation implique-t-il des
comparaisons au niveau local fini, ou bien des comparaisons dans
l'infinitésimal?

\CITATION{Riemann 
n'est jamais facile \`a lire, mais sa c\'el\`ebre allocution
<<\,{\em Sur les hypoth\`eses qui servent de fondement \`a la
g\'eom\'etrie}\,>> pr\'esente des difficult\'es de compr\'ehension
d'un type tout \`a fait sp\'ecial.
\REFERENCE{\S~100,~p.~\pageref{Riemann-jamais-facile}~ci-dessous.}}

Autre énigme: que signifie la première hypothèse elliptique d'après
laquelle les lignes possèdent une longueur, indépendamment de leur
position? Faut-il y voir une allusion à la libre mobilité des courbes
et à l'invariance de leur longueur, lorsqu'on effectue une série de
transformations ponctuelles de l'espace? Mais une telle interprétation
entrerait en contradiction avec le fait que la métrisation est et doit
être absolument indépendante de la mobilité; une telle indépendance
vaut en effet déjà pour la théorie gaussienne des surfaces, et Riemann
a explicitement énoncé que cette théorie contient les fondements de la
question qu'il va traiter.
En tout cas, par rapport à la théorie gaussienne des
surfaces courbes qui héritent d'une métrique <<\,ondulée\,>> par
restriction de la métrique pythagoricienne <<\,plate\,>>
de
l'espace à trois dimensions, Riemann va {\em
renverser\footnote{\,
%%%%%%%%%%%%%%%%%%%%%%%-------DEBUT--------%%%%%%%%%%%%%%%%%%%%%%%%%%%
{\em Voir} le \S~1.6 p.~\pageref{le-renversement-riemannien}
} %%%%%%%%%%%%%%%%%%%%%%%%-----FIN-----%%%%%%%%%%%%%%%%%%%%%%%%%%%%%%%
complètement la réflexion}. 

En effet, sans s'autoriser à céder ni à la formulation aisée de
problèmes ouverts qui consiste simplement à augmenter le nombre de
variables, ni à la générativité symbolique des expressions formelles,
c'est-à-dire plus précisément, sans annoncer d'emblée à son auditoire:

\smallskip

<<\,{\em \'Etudions maintenant l'expression différentielle quadratique
$\sum_{ i, j= 1}^n \, g_{ ij} ( x_1, \dots, x_n) \, dx_i dx_j$ à $n$
variables qui généralise visiblement l'expression connue en
coordonnées paramétriques intrinsèques de la métrique $E ( u, v) \,
du^2 + 2\, F ( u, v) \, du dv + G(u, v) \, dv^2$ sur une surface
courbe, dont Gauss a montré, dans son {\em Theorema Egregium}, qu'elle
possède la mesure de courbure comme invariant à travers toute
transformation isométrique}\,>>,

\smallskip\noindent
Riemann va plutôt, en renversant le sens de son étude, {\em chercher à
trouver des principes de genèse a priori} qui montreront en quoi
l'expression {\em quadratique} gaussienne, ainsi que sa généralisation
à des dimensions supérieures, est en un certain sens naturelle,
nécessaire, ou tout du moins <<\,la plus simple possible\,>> qui
pourrait s'offrir à l'étude dans un {\em a priori} relatif,
reconstitué {\em a posteriori}, de la connaissance mathématique.
Régressive dans l'{\em a posteriori} par rapport à la théorie de
Gauss, l'étude riemannienne cherche à ouvrir une voie nouvelle vers
l'{\em a priori} génétique.

Par sa démarche, Riemann est donc un véritable {\em métaphysicien des
mathématiques}: nous devons nous interroger, dit-il en effet, sur
l'existence de causes profondes qui pourraient expliquer l'émergence
de telles formes symboliques, ou de telles structures
mathématiques.
Le point de vue riemannien se situe donc bien en amont de toute option
philosophique unilatérale sur l'essence des mathématiques et sait se
soustraire aux polémiques afférentes; idéalisme, platonisme, réalisme,
constructivisme, intuitionnisme, historicisme, essentialisme,
axiomatisme, formalisme: chacune de ces options philosophiques est
engagée dans une problématique d'essence tellement profonde que les
réponses possibles sont encore noyées dans l'ouverture et dans
l'indécision, et tout penseur d'inspiration riemannienne se voit dans
l'obligation philosophique d'accepter cet état de fait.

Ainsi la partie du discours de Riemann où il rappelle son principe
méthodologique général d'investigation n'est-elle plus maintenant
énigmatique pour notre analyse.

\CITATION{Nous 
arrivons au second des problèmes posés plus haut, savoir à
l'étude des rapports métriques dont une multiplicité est susceptible,
et des conditions suffisantes pour la détermination de ces rapports
métriques.
\REFERENCE{\cite{riem1898},~p.~285.}}

\noindent
Toute genèse doit en effet procéder par {\em conditions
démonstratives}, au moins suffisantes dans un premier moment, et si
possible ensuite, {\em nécessaires et suffisantes}, afin de resserrer
au mieux peut-être les {\em écarts synthétiques} entre concepts qui
pourraient cacher des pétitions de principe, des hypothèses
implicites, ou mieux encore, des concepts nouveaux et ouverts qui
pourraient connaître un destin inattendu dans l'histoire des
mathématiques\footnote{\,
%%%%%%%%%%%%%%%%%%%%%%%-------DEBUT--------%%%%%%%%%%%%%%%%%%%%%%%%%%%
Encore une fois, rappelons que la pensée mathématique structuraliste
contemporaine ne place jamais la question à un niveau aussi
problématisant. En effet, lorsqu'on s'autorise à commencer un article
ou un exposé par une phrase telle que <<\,Soit $M$ une variété
différentielle munie d'une métrique riemannienne $g$\,>>, ou telle que
<<\,Soit $\frac{ dz}{ 1 - \vert z \vert^2}$ la métrique de Poincaré
sur le disque unité $\Delta = \{ z \in \C : \, \vert z \vert < 1 \}$
dans $\C$\,>>, l'acte de position que désigne l'expression <<\,Soit
$X$ un objet mathématique défini\,>> réfère à un concept
considéré comme déjà donné dans une architecture paradigmatique
constituée. 
}. %%%%%%%%%%%%%%%%%%%%%%%%-----FIN-----%%%%%%%%%%%%%%%%%%%%%%%%%%%%%%%
Ainsi sur le chemin génétique qui conduit aux métriques
différentielles quadratiques maintenant dites <<\,riemanniennes\,>>,
Riemann va-t-il poser successivement plusieurs hypothèses ouvertes qui
pourraient conduire à d'autres types de rapports métriques possibles.

\HEAD{Chapitre~1.\,\,\,\,L'ouverture riemannienne}{
1.17.\,\,\,Genèse des métriques riemanniennes}

\medskip\noindent{\bf 1.17.~Genèse des métriques riemaniennes.}
\label{genese-des-metriques-riemanniennes}
Les déterminations de lieu étant ramenées aux déterminations
simultanées de $n$ grandeurs numériques $x_1, x_2, x_3, \dots, x_n$
(\S~1.15), le problème consiste maintenant à trouver une expression
mathématique pour la longueur des lignes courbes tracées dans la
multiplicité. Comme dans l'espace ordinaire, la donation d'une ligne
courbe revient à ce que les quantités $x_i$ dépendent paramétriquement
d'une seule variable auxiliaire.

\CITATION{Je 
ne traiterai ce problème que sous certaines restrictions, et je me
bornerai d'abord aux lignes dans lesquelles les rapports entre les
accroissements $dx$ des variables $x$ correspondantes varient d'une
manière continue.
\REFERENCE{\cite{riem1898},~p.~286.}}

\noindent {\em Leitmotiv} riemannien:
encore une annonce d'ouverture potentielle laissée de côté par le
choix d'une hypothèse déterminée. Poser une hypothèse, c'est
s'écarter éventuellement d'un (autre) univers mathématique, c'est
bifurquer vers une certaine branche de l'arbre mathématique, 
sans examiner d'autres branches, sans
explorer d'autres univers.

Par conséquent, Riemann décide d'{\em infinitésimaliser}\footnote{\,
%%%%%%%%%%%%%%%%%%%%%%%-------DEBUT--------%%%%%%%%%%%%%%%%%%%%%%%%%%%
En fait, dans l'en-tête (énigmatique) de ce paragraphe, Riemann aurait
été probablement mieux inspiré d'annoncer cette infinitésimalisation
comme l'une de ses hypothèses génétiques principales. }
%%%%%%%%%%%%%%%%%%%%%%%%-----FIN-----%%%%%%%%%%%%%%%%%%%%%%%%%%%%%%%
le problème: avec cette hypothèse de 
continuité\footnote{\,
%%%%%%%%%%%%%%%%%%%%%%%-------DEBUT--------%%%%%%%%%%%%%%%%%%%%%%%%%%%
---\,\,par laquelle
il faudrait entendre plus rigoureusement une hypothèse de
différentiabilité d'ordre au moins égal à $1$\,\,---
}, %%%%%%%%%%%%%%%%%%%%%%%%-----FIN-----%%%%%%%%%%%%%%%%%%%%%%%%%%%%%%%
les lignes
peuvent être décomposées en portions infinitésimales,
et si l'on s'autorise\footnote{\,
%%%%%%%%%%%%%%%%%%%%%%%-------DEBUT--------%%%%%%%%%%%%%%%%%%%%%%%%%%% 
\`A cause de
ce recours à l'intégration\,\,---\,\,concept d'analyse encore
problématique qui exige l'infini\,\,---, Engel et Lie objecteront
p.~\pageref{critique-integration} ci-dessous que les considérations de
Riemann fournissent peu d'éclaircissements quant aux fondements de la
géométrie purement élémentaire.
} %%%%%%%%%%%%%%%%%%%%%%%%-----FIN-----%%%%%%%%%%%%%%%%%%%%%%%%%%%%%%%
de plus à disposer de la théorie de l'intégration
afin de sommer toutes les longueurs infinitésimales des éléments de
ligne placés bout à bout, {\em il suffit alors}\footnote{\,
%%%%%%%%%%%%%%%%%%%%%%%-------DEBUT--------%%%%%%%%%%%%%%%%%%%%%%%%%%%
Riemann recherche en effet des conditions seulement {\em suffisantes}
pour la détermination des rapports métriques dans une multiplicité.
\`A chaque fois qu'une hypothèse simplificatrice ou
spécificatrice est admise, le choix d'une discontinuité conceptuelle
estompe une fraction de la {\em nécessité} qui doit être corrélative
de la {\em suffisance}.
} %%%%%%%%%%%%%%%%%%%%%%%%-----FIN-----%%%%%%%%%%%%%%%%%%%%%%%%%%%%%%%
de trouver une expression pour la longueur de tout élément linéaire
infinitésimal $(dx_1, dx_2, \dots, dx_n)$ qui est situé en un point
$(x_1, x_2, \dots, x_n)$.

Poursuite du raisonnement: dans l'infinitésimal et en un point $(x_1,
x_2, \dots, x_n)$ fixé, les rapports d'accroissement entre les
composantes $dx_i$ de l'élément infinitésimal en question le long
d'une ligne donnée peuvent être considérés comme {\em constants}.
Mais quand le point $(x_1, x_2, \dots, x_n)$ varie, ces rapports
cessent d'être constants, et puisque les lignes sont libres de
représenter, en un point donné quelconque, toutes les directions
possibles qui passent par ce point, il en découle que les rapports
d'accroissement infinitésimal le long d'une ligne doivent en fait
{\em dépendre du point}, et donc aussi: {\em seulement} du point. Par
conséquent, sous ces deux hypothèses fondamentales
d'infinitésimalisation première et d'intégration seconde, Riemann a
ramené la question géométrique de la genèse des rapports métriques à
une question d'Analyse, à savoir: déterminer une {\em fonction} du
lieu et de l'élément infinitésimal\footnote{\,
%%%%%%%%%%%%%%%%%%%%%%%-------DEBUT--------%%%%%%%%%%%%%%%%%%%%%%%%%%%
Dans le \S~100 p.~\pageref{Omega-metrique} sq. ci-dessous que le
lecteur est invité à lire en parallèle pour de plus amples
éclaircissements (\cf~aussi \cite{ vui1962}, pp.~409--412), Engel et
Lie réexpriment les raisonnements de Riemann en utilisant un langage
purement analytique. Par ailleurs, dans la recherche d'une expression
fonctionnelle pour la métrique, ils traitent simultanément du cas
local fini, inspiré de leur théorie des invariants, et du cas
infinitésimal (Riemann).
}: %%%%%%%%%%%%%%%%%%%%%%%%-----FIN-----%%%%%%%%%%%%%%%%%%%%%%%%%%%%%%%
\[
\Omega
=
\Omega(x_1,\dots,x_n;\,dx_1,\dots,dx_n)
\] 
la plus générale possible qui puisse fournir, en restriction sur les
courbes quelconques, la longueur de n'importe quelle ligne tracée dans
la multiplicité.

Autrement dit, la longueur en un point $(x_1,\dots, x_n)$ d'un
élément infinitésimal quelconque $(dx_1, \dots, dx_n)$ attaché en ce
point\,\,---\,\,que l'on note habituellement\footnote{\,
%%%%%%%%%%%%%%%%%%%%%%%-------DEBUT--------%%%%%%%%%%%%%%%%%%%%%%%%%%%
Cette notation utilisée depuis le 18\textsuperscript{ème} siècle et
reprise par Gauss, réfère à la longueur d'un élément d'arc
infinitésimal d'une courbe tracée dans le plan ou dans l'espace.
} %%%%%%%%%%%%%%%%%%%%%%%%-----FIN-----%%%%%%%%%%%%%%%%%%%%%%%%%%%%%%%
<<\,$ds$\,>>\,\,---\,\,est
égale à cette fonction pour l'instant inconnue: 
\[
\aligned
ds
&
=
\text{\rm longueur}_x(dx)
\\
&
=
\text{\rm distance}\big(x,\,x+dx\big)
\\
&
=
\Omega(x;\,dx).
\endaligned
\]

\CITATION{J'admettrai, 
en second lieu, que la longueur de l'élément linéaire,
abstraction faite des quantités du second ordre, reste invariable,
lorsque tous les points de cet élément subissent un même déplacement
infiniment petit, ce qui implique en même temps que, si toutes les
quantités $dx$ croissent dans un même rapport, l'élément linéaire
varie également dans ce même rapport.
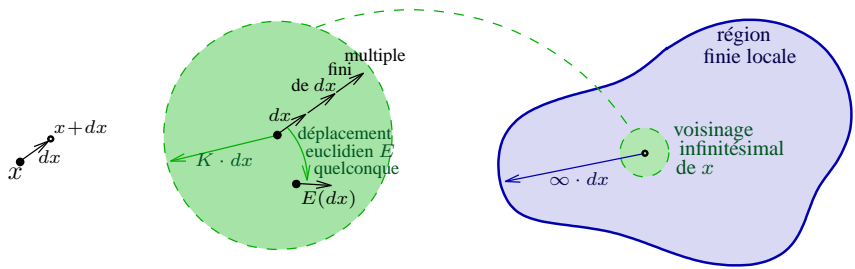
\REFERENCE{\cite{riem1898},~p.~286.}}

\begin{center}
\input euclidien-infinitesimal.pstex_t
\end{center}

\noindent
Ici seulement\,\,---\,\,mais exclusivement à niveau
infinitésimal\,\,---, Riemann utilise l'hypothèse énigmatique d'après
laquelle les lignes possèdent une longueur indépendamment de leur
position. Ainsi, la longueur de $dx$ doit être conservée lors de tout
déplacement euclidien $E$ qui est restreint à un voisinage
infinitésimal de $x$:
\[
\text{\rm longueur}(dx)
=
\text{\rm longueur}\big(E(dx)\big).
\]
Donc dans l'infiniment petit, la métrique recherchée semble être
supposée comme devant être euclidienne, ce qui équivaut à
dire\,\,---\,\,on peut le démontrer\,\,---\,\,que la métrique est
donnée par une forme quadratique différentielle positive. Mais ce
serait brûler les étapes et omettre de découvrir de nouveaux noyaux
conceptuels possibles dans le graphe problématique et
virtuel des concepts métriques.

En effet, Riemann ne se sert en fait de son hypothèse mystérieuse
(demandant que chaque ligne puisse être mesurée par toute autre ligne)
que pour en déduire que la longueur d'un multiple entier fini $k\, dx$
de $dx$ est égale à $k$ fois la longueur de $dx$\footnote{\,
%%%%%%%%%%%%%%%%%%%%%%%-------DEBUT--------%%%%%%%%%%%%%%%%%%%%%%%%%%%
---\,\,en effet, $k \, dx$ s'obtient en mettant bout à bout $k$ copies
de $dx$ dans la même direction que le $dx$ de départ (\voir~le
diagramme), et chacune de ces copies est tout simplement obtenue par
une translation parallèle à $dx$ 
dans le voisinage infinitésimal de $x$\,\,---
}, %%%%%%%%%%%%%%%%%%%%%%%%-----FIN-----%%%%%%%%%%%%%%%%%%%%%%%%%%%%%%%
et aussi en même temps, que la longueur de $-dx$\footnote{\,
%%%%%%%%%%%%%%%%%%%%%%%-------DEBUT--------%%%%%%%%%%%%%%%%%%%%%%%%%%%
---\,\,qui se déduit de $dx$ par une transformation
euclidienne standard, la symétrie orthogonale par rapport à
l'hyperplan orthogonal à $dx$,\,\,---
} %%%%%%%%%%%%%%%%%%%%%%%%-----FIN-----%%%%%%%%%%%%%%%%%%%%%%%%%%%%%%%
s'identifie à la longueur de $dx$. 

Ensuite, par un argument de continuité qui pourrait consister à faire
tendre $k$ vers l'infini tout en rapetissant $dx$ afin que $k \,
dx$ demeure une quantité infinitésimale, Riemann semble en déduire que
la fonction $\Omega ( x; \, dx)$ pourra être n'importe quelle fonction
homogène du premier degré en $dx$, à savoir qui 
satisfait\footnote{\,
%%%%%%%%%%%%%%%%%%%%%%%-------DEBUT--------%%%%%%%%%%%%%%%%%%%%%%%%%%%
Riemann sous-entend intuitivement qu'une telle fonction est en quelque
sorte semi-explicite, voire localement développable en série entière
par rapport à $dx$, puisqu'il se la représente comme une fonction
homogène du premier degré en les quantités $dx$ <<\,dans laquelle les
constantes arbitraires seront des fonctions continues de $x$\,>>.
}: %%%%%%%%%%%%%%%%%%%%%%%%-----FIN-----%%%%%%%%%%%%%%%%%%%%%%%%%%%%%%%
\[
\Omega(x;\,\lambda\,dx)
=
\vert\lambda\vert\,\Omega(x;\,dx),
\]
pour tout nombre réel fini $\lambda$. Or cette nouvelle conclusion
provisoire ne nécessite absolument pas que les rapports métriques
soient euclidiens dans l'infinitésimal. En effet, 
cette propriété demande
seulement que les longueurs se dilatent dans l'infinitésimal de
manière purement homothétique\,\,---\,\,exigence minimale qui laisse
encore disponible une très grande généralité. On peut donc dire qu'à
cet endroit-là (bien qu'il semble en avoir été clairement soucieux),
Riemann n'a pas réellement pris le temps de resserrer les hypothèses
minimales qui conduisent aux métriques dites de Finsler, 
nettement plus générales que les métriques 
riemanniennes ({\em voir}~\cite{ spiv1970, cher1996, bcs2000}). 

La poursuite du raisonnement marque alors un revirement inattendu
de la spéculation, puisqu'après fixation d'un point-origine $(x_1^0,
\dots, x_n^0)$ , Riemann cherche maintenant une fonction {\em non
infinitésimale} du lieu:
\[
\Omega\big(x_1,\dots,x_n;\,x_1^0,\dots,x_n^0\big)
\]
dont les ensembles de niveau $\{ x : \, \Omega ( x; \, x^0) = {\rm
const.} \}$ s'identifient aux lieux équidistants de l'origine. En
fait, ce retour au niveau macroscopique local fini va permettre de se
rapprocher en pensée des différentielles {\em quadratiques} qui
généralisent les métriques gaussiennes, car il est alors tout à fait
naturel que la {\em première différentielle par rapport à $x$} d'une
telle fonction $\Omega ( x; \, x^0)$ doive nécessairement 
s'annuler\footnote{\,
%%%%%%%%%%%%%%%%%%%%%%%-------DEBUT--------%%%%%%%%%%%%%%%%%%%%%%%%%%%
Si une fonction $\omega = \omega ( x_1, \dots, x_n)$ de $n$ variables
possède une dérivée partielle $\frac{ \partial \omega}{ \partial x_i}$
qui ne s'annule pas en un point $x^0$, alors $\omega$ croît ou décroît
{\em strictement} 
(selon le signe de cette dérivée partielle) le long d'un
petit segment affine parallèle à l'axe des $x_i$ qui passe par $x^0$.
} %%%%%%%%%%%%%%%%%%%%%%%%-----FIN-----%%%%%%%%%%%%%%%%%%%%%%%%%%%%%%%
pour $x = x^0$, puisque toute notion de distance que l'on peut
s'imaginer doit évidemment atteindre son minimum, égal à $0$, au point
de référence $x = x^0$. Le quadratique (ordre $2$) comme successeur du
linéaire (ordre $1$) s'introduit donc seulement à travers le principe
de stabilité existentielle des minima. On notera que presque
immédiatement après avoir voulu passer au macroscopique, Riemann {\em
réinfinitésimalise} le raisonnement en considérant la différentielle
$d_x \Omega ( x^0; \, x^0)$.

Puisque cette première différentielle $d_x \Omega ( x^0; \, x^0)$
s'annule, le comportement quantitatif de $\Omega ( x; \, x^0)$,
lorsque $x$ parcourt un voisinage infinitésimal de $x^0$, sera
entièrement représenté par sa {\em différentielle seconde}:
\[
d^2\Omega
=
\sum_{i=1}^n\,\sum_{j=1}^n\,
\frac{\partial^2\Omega}{\partial x_i\partial x_j}(x^0;\,x^0)\,
dx_i\,dx_j. 
\]
Par construction, cette différentielle seconde
donne donc une bonne approximation de type 
Taylor-Young pour la valeur: 
\[
\aligned
\Omega(x^0+dx;\,x^0)
-
\Omega(x^0;\,x^0)
&
=
\Omega(x^0+dx;\,x^0)
-
0
\\
&
=
\text{\rm longueur}_{x^0}(dx)
\\
&
=
ds\big\vert_{x^0}. 
\endaligned
\]
\`A présent, 
nouvelle bifurcation d'hypothèses: si tous les coefficients $\frac{
\partial^2 \Omega }{ \partial x_i \partial x_j }(x^0; \,x^0)$ de cette
différentielle seconde s'annulent, le développement devra
se poursuivre jusqu'aux termes d'ordre $3$. Mais comme tout produit
$dx_i dx_j dx_k$ de degré trois entre différentielles change de signe
quand on change $dx$ en $-dx$, et comme la fonction distance
recherchée doit forcément être positive, il en découle que dans ce
cas, la différentielle troisième $d^3 \Omega ( x^0; \, x^0)$ doit donc
nécessairement s'annuler. Ainsi, on doit alors tester si la
différentielle {\em quatrième} ne s'annule pas, et ainsi
de suite.

En toute généralité, ce raisonnement qui présuppose l'analyticité de
la fonction $\Omega$, montre que la métrique infinitésimale recherchée
doit s'identifier à la racine $2k$-ième d'une expression homogène de
degré $2k$ toujours positive dont les coefficients sont des fonctions
de $x$:
\[
ds
=
\sqrt[2\kappa]{
\sum_{i_1,\dots,i_{2\kappa}=1}^n\,
\omega_{i_1,\dots,i_{2\kappa}}(x_1,\dots,x_n)\,
dx_{i_1}\cdots dx_{i_{2\kappa}}}\,.
\]
Sans perte de généralité, on peut supposer que de tels coefficients
$\omega_{ i_1, \dots, i_{ 2\kappa}}$ sont complètement symétriques par
rapport à leurs indices inférieurs, à savoir:
\[
\omega_{i_1,\dots,i_{2\kappa}}(x)
\equiv
\omega_{i_{\sigma(1)},\dots,i_{\sigma(2\kappa)}}(x),
\]
pour toute permutation $\sigma$ de l'ensemble $\{ 1, \dots,
2\kappa\}$.

Le cas le plus simple est bien sûr celui des formes différentielles
{\em quadratiques} $(2 \kappa = 2)$, pour lequel le carré $ds^2$ de la
longueur d'un élément infinitésimal quelconque de coordonnées non
toutes nulles $(dx_1, \dots, dx_n)$ basé en un point $(x_1, \dots, x_n)$
 est donné par une expression du second ordre en les $dx_i$:
\[
ds^2
=
\sum_{i,j=1}^n\,
g_{i,j}(x_1,\dots,x_n)\,
dx_i\,dx_j
\]
à coefficients des fonctions arbitraires $g_{ i,j} \equiv g_{ j,i}$ de
$x$, de telle sorte que la somme ne prenne que des valeurs strictement
positives.

En conclusion, la genèse des métriques riemanniennes transcende tout
acte de postulation axiomatique {\em a posteriori}.  Problématisante,
la genèse riemannienne procède par spécification progressive
d'hypothèses qui sont hiérarchisées en ordre de généralité. Chaque
choix engage la pensée dans un nouvel irréversible-synthétique.

La coprésence de ces bifurcations spéculatives indique {\em
l'ouverture collatérale permanente de la pensée mathématique}. Au sein
même du concept final de différentielle quadratique infinitésimale
positive, le degré de liberté et d'arbitraire dans le choix des
fonctions $g_{ i,j} ( x)$ maintient l'ouverture intrinsèque du concept
et le prédispose à une plasticité remarquable, confirmée par sa
capacité à héberger des théories physiques aussi variées que la
cristallographie, la mécanique des milieux continus, ou encore la
théorie de la relativité généralisée.

\HEAD{Chapitre~1.\,\,\,\,L'ouverture riemannienne}{
1.18.\,\,\,Surfaces de courbure constante}

\medskip\noindent{\bf 1.18.~Surfaces de courbure constante.}
\label{surfaces-de-courbure-constante}
Le travail de Gauss sur la théorie intrinsèque des surface a trouvé
des continuateurs (\cite{ rei1973,rose1988, krey1994}) à l'Université
de Dorpat, maintenant Tart\`u, une ville
de langue germanique située dans une province
estonienne.  Senff a publié en 1831 les formules
qu'on attribue aujourd'hui à Frenet\footnote{\,
%%%%%%%%%%%%%%%%%%%%%%%-------DEBUT--------%%%%%%%%%%%%%%%%%%%%%%%%%%%
Pour une excellente présentation de la théorie des courbes et des
surfaces dans l'espace à trois dimensions, nous renvoyons aux
leçons de Do Carmo~\cite{ doca1976}. Les aspects philosophiques fins de
l'émergence de la théorie des surfaces de Gauss ne pourront pas être
abordés ici, et nous bornerons notre analyse à
l'examen résumé de la façon dont 
Riemann semble être parvenu au concept de
courbure, en nous basant sur Weyl
(\cite{ weyl1990, weyl1988}) et
sur Spivak (\cite{ spiv1970}). 
}; %%%%%%%%%%%%%%%%%%%%%%%%-----FIN-----%%%%%%%%%%%%%%%%%%%%%%%%%%%%%%%
Peterson a soutenu en 1853 une thèse sur les équations aujourd'hui
attribuées à Mainardi-Codazzi; et surtout Minding, figure la plus
influente, a travaillé sur le développement des lignes courbes à
l'intérieur de surfaces courbes, a introduit le concept de courbure
géodésique, et a étudié les surfaces dont la courbure gaussienne est
constante.

Pour ce qui nous intéresse, Minding a su exprimer la métrique
gaussienne d'une surface de courbure constante $\kappa > 0$ sous la
forme\footnote{\,
%%%%%%%%%%%%%%%%%%%%%%%-------DEBUT--------%%%%%%%%%%%%%%%%%%%%%%%%%%%
Lorsque le $ds^2$ est représentée en coordonnées polaires géodésiques
sous la forme normalisée $ds^2 = dp^2 + G ( p, q) \, dq^2$, sa courbure
de Gauss s'exprime alors par la formule relativement simple: $\kappa =
\kappa ( p) = -\frac{ 1}{ \sqrt{ G}}\,
\frac{ \partial^2 \sqrt{ G}}{ \partial p^2}$. 
Sachant que le signe de la dérivée seconde change suivant qu'on a
affaire au sinus (tout court): $\frac{ \partial^2}{ \partial p^2}
(\sin p \sqrt{ \kappa}) = -
\sqrt{ \kappa}^2 \, \sin p \sqrt{ \kappa}$,
ou au sinus hyperbolique: $\frac{ \partial^2}{ \partial p^2} (\sinh p
\sqrt{ \kappa}) = \sqrt{ \kappa}^2 \, \sinh p \sqrt{ \kappa}$, 
on retrouve effectivement $\kappa$ dans le premier cas, et $- \kappa$
dans le second cas.  
} %%%%%%%%%%%%%%%%%%%%%%%%-----FIN-----%%%%%%%%%%%%%%%%%%%%%%%%%%%%%%%
normalisée suivante: 
\[
ds^2
=
dp^2
+
\big(
{\textstyle{\frac{1}{\sqrt{\kappa}}}}\,\sin\,p\sqrt{\kappa}
\big)^2\,dq^2
\]
ou lorsque $-\kappa < 0$ sous la forme:
\[
ds^2
=
dp^2
+
\big(
{\textstyle{\frac{1}{\sqrt{\kappa}}}}\,\sinh\,p\sqrt{\kappa}
\big)^2\,dq^2
\]
Le cas de la courbure nulle s'obtient en prenant la limite, lorsque
$\kappa$ tend vers zéro, de chacune de ces deux formules\footnote{\,
%%%%%%%%%%%%%%%%%%%%%%%-------DEBUT--------%%%%%%%%%%%%%%%%%%%%%%%%%%%
Rappelons que $\sin t = t - \frac{ 1}{ 6}\, t^3 + \cdots$ et que
$\sinh t = t + \frac{ 1}{ 6} \, t^3 + \cdots$, d'où $\lim_{ \kappa \to
0}\, \frac{ 1}{ \sqrt{ \kappa}}\, \sin p \sqrt{ \kappa} = p$
et aussi $\lim_{ \kappa \to 0}\,  \frac{ 1}{ \sqrt{ \kappa}}\,
\sinh p \sqrt{ \kappa} = p$.  
}, %%%%%%%%%%%%%%%%%%%%%%%%-----FIN-----%%%%%%%%%%%%%%%%%%%%%%%%%%%%%%% 
et l'on retrouve ainsi l'expression de la métrique pythagoricienne en
coordonnées polaires (rayon $p$, angle $q$): $ds^2 = dp^2 + p^2 dq^2$.

Une autre expression normalisée des métriques gaussiennes de courbure
constante était très vraisemblablement connue de Riemann, et elle
possède l'avantage remarquable, par rapport aux formules
précédentes, de ne pas contraindre à distinguer plusieurs cas\footnote{\,
%%%%%%%%%%%%%%%%%%%%%%%-------DEBUT--------%%%%%%%%%%%%%%%%%%%%%%%%%%%
Lorsque le $ds^2$ est représenté en coordonnées isothermes sous la forme
normalisée: $ds^2 = \lambda^2 ( u, v) \big[ du^2 + dv^2 \big]$, la
courbure est donnée par une formule que Gauss possédait
déjà en 1822: 
\[
\kappa
=
\kappa(u,v)
=
-\frac{1}{\lambda^2}\,
\bigg(\frac{\partial^2\log\lambda}{\partial u^2}
+\frac{\partial^2\log\lambda}{\partial v^2}
\bigg),
\] 
et qui apparaissait dans son {\em Copenhagen Preisschrift} sur les
applications conformes qui lui a valu le prix de l'Académie de
Copenhague (\cite{ krey1994}).  L'application de cette formule
générale dans le cas où $\frac{ 1}{
\lambda^2} = \frac{ 1}{ 4} + \kappa ( x_1^2 + x_2^2)$
fournit effectivement la constante $\kappa$, quel que soit le nombre
réel $\kappa$ fixé à l'avance. 
}: %%%%%%%%%%%%%%%%%%%%%%%%-----FIN-----%%%%%%%%%%%%%%%%%%%%%%%%%%%%%%%
\[
ds^2
=
\frac{dx_1^2+dx_2^2}{
\big[1+\frac{\kappa}{4}(x_1^2+x_2^2)\big]^2}. 
\]
Dans son {\em Habilitationsvortrag}, la seule formule significative
que Riemann osera signaler à son auditoire d'universitaires issus de
tous les horizons sera la généralisation évidente
de cette expression à la dimension $n$
quelconque:
\[
ds^2 
= 
\frac{dx_1^2+\cdots+dx_n^2}{
\big[1+\frac{\kappa}{4}\,(x_1^2+\cdots+x_n^2)\big]^2}, 
\]
dans laquelle on peut immédiatement relire la
formule précédente en égalant à zéro $(n-2)$
variables $x_i$. 

\HEAD{Chapitre~1.\,\,\,\,L'ouverture riemannienne}{
1.19.\,\,\,Courbure sectionnelle de Riemann-Christoffel-Lipschitz}

\medskip\noindent{\bf 1.19.~Courbure sectionnelle de 
Riemann-Christoffel-Lipschitz.}
\label{courbure-Riemann-Christoffel-Lipschitz}
{\em Question}: 
en passant à la dimension quelconque $n \geqslant 2$ et pour des
métriques quelconques, pourquoi Riemann a-t-il envisagé de prolonger
la théorie de Gauss sous l'angle de la courbure dite {\sl
sectionnelle}, c'est-à-dire en {\em sectionnant} les multiplicité de
dimension $n$ par des {\em surfaces} de dimension $2$?  La {\em
Commentatio} (\cf~note p.~\pageref{commentatio}) fournit une réponse
qui témoigne clairement de l'enracinement de cette genèse
dans la matrice du tridimensionnel.

\CITATION{
L'expression $\sqrt{ \sum \, b_{ \iota, \iota'}\, ds_\iota \, ds_{
\iota'}}$ peut être envisagée comme l'élément linéaire dans un espace
généralisé de $n$ dimensions transcendant notre intuition. Si dans cet
espace on trace toutes les lignes les plus courtes issues du point
$(s_1, s_2, \dots, s_n)$, dans lesquelles les éléments initiaux de
variation des $s$ sont comme les rapports $\alpha ds_1 +
\beta \delta s_1 \colon \alpha ds_2 + \beta \delta s_2 \colon
\cdots \colon \alpha ds_n + \beta \delta s_n$, où 
$\alpha$ et $\beta$ désignent des quantités arbitraires, alors ces
lignes constituent une surface qui peut être développée dans l'espace
de notre intuition commune\footnotemark.
\REFERENCE{\cite{riem1892},~p.~382.}}

\footnotetext{\,
%%%%%%%%%%%%%%%%%%%%%%%-------DEBUT--------%%%%%%%%%%%%%%%%%%%%%%%%%%%
{\sl Expressio $\sqrt{ \sum \, b_{ \iota, \iota'}\, ds_\iota \, ds_{
\iota'}}$ spectari potest tanquam elementum lineare in spatio
generaliore $n$ dimensionum nostrum intuitum transcendente. Quodsi in
hoc spatio a puncto $(s_1, s_2, \dots, s_n)$ ducantur omnes lineae
brevissimae, in quarum elementis initialibus variationes ipsarum $s$
sunt ut $\alpha ds_1 +
\beta \delta s_1 \colon \alpha ds_2 + \beta \delta s_2 \colon
\cdots \colon \alpha ds_n + \beta \delta s_n$, 
denotantibus $\alpha$ et $\beta$ quantitates quaslibet, hae lineae
superficiem constituent, quam in spatium vulgare nostro intuitui
subjectum evolvere licet}.  Ici, les deux éléments infinitésimaux
$(ds_1, ds_2, \dots, ds_n)$ et $(\delta s_1, \delta s_2, \dots, \delta
s_n)$ basés en un point de coordonnées $(s_1, s_2, \dots, s_n)$ sont
supposés être linéairement indépendants, et la combinaison linéaire
générale:
\[
\big(
\alpha ds_1+\beta\delta s_1,\,
\alpha ds_2+\beta\delta s_2,\,
\dots,\,
\alpha ds_n+\beta\delta s_n
\big)
\]
comprend alors tous les éléments linéaires contenus dans le plan
qu'ils engendrent. Riemann considère donc la {\em surface} locale et
finie qui est obtenue en intégrant toutes les géodésiques issues du
point dans toutes ces directions et il s'imagine alors qu'une telle
surface, interne à la multiplicité (variété) initiale, pourrait en
être extraite afin de se réaliser visuellement dans un espace
tridimensionnel auxiliaire. \`A partir de cet extrait, 
on pourrait même s'imaginer que l'invention riemannienne de la courbure
par sectionnement obéissait à simple exigence d'appropriation
intuitive.
} %%%%%%%%%%%%%%%%%%%%%%%%-----FIN-----%%%%%%%%%%%%%%%%%%%%%%%%%%%%%%%

\noindent
Nécessité pour la {\em pensée}, pour la {\em conception}, et surtout
pour l'{\em intuition}: nécessité de ressaisir la courbure sous un angle
bidimensionnel et gaussien, car le pluridimensionnel {\em nous
transcende}. \'Etonnant coup de chance riemannien que confirmera
pleinement la réinterprétation tensorielle: tous les invariants
d'ordre deux d'une métrique quadratique infinitésimale peuvent être
obtenus en se restreignant à des {\em surfaces} qu'on inscrit dans la
variété et qu'on oriente à volonté dans des directions
arbitraires.

Toutefois, la réduction
s'arrête à la dimension $2$, car le tranchage par des
objets de dimension $1$ rend la courbure invisible.
En effet, chaque courbe 
(au moins de classe $\mathcal{ C}^1$)
est intrinsèquement
équivalente (isométrique\footnote{\,
%%%%%%%%%%%%%%%%%%%%%%%-------DEBUT--------%%%%%%%%%%%%%%%%%%%%%%%%%%%
---\,\,grâce à la paramétrisation par longueur d'arc, ou 
à la rectification d'un lacet par le
geste physique\,\,--- })
%%%%%%%%%%%%%%%%%%%%%%%%-----FIN-----%%%%%%%%%%%%%%%%%%%%%%%%%%%%%%% à
à un simple segment de droite: la {\em régression en dimension} du
caractère {\em intrinsèque} de la courbure doit s'arrêter net à la
dimension $2$.  Et inversement, le passage de la dimension $1$ à la
dimension 2 imposait
un saut qualitatif inattendu
que Gauss avait su découvrir: naissance de la courbure {\em
intrinsèque} des surfaces, courbure qui demeure ponctuellement
invariable dans toute transformation isométrique, alors que toutes les
lignes tracées dans une surface sont dépossédées de toute rigidité
intrinsèque.
 
Par ailleurs, la découverte de Gauss aurait pu faire croire qu'en
dimension $n \geqslant 3$, d'autres phénomènes spécifiques et d'autres
invariants inattendus nouveaux émergeraient, qui seraient eux aussi
{\em propres aux dimensions supérieures}. Peut-être même sans qu'il
s'en soit réellement douté, Riemann a-t-il été conduit à entrevoir la
bidimensionnalité pure de la courbure. Seuls les travaux de Lipschitz
et de Christoffel confirmeront cette intuition.  Impossible donc de se
faire une idée {\em a priori} de l'intrinsèque qui l'exempte de
l'imprévisibilité contingente des nécessités latérales inscrites dans
des modalités hypothétiques.

En toute dimension $n \geqslant 2$, Riemann va donc approcher le
concept de courbure par sections de surfaces, domaine où la théorie de
Gauss s'appliquera. C'est peut-être cette idée fondamentale qui a
guidé le mystérieux calcul que la {\em Commentatio} (\cf~note
p.~\pageref{commentatio}) nous transmet sans détails intermédiaires:
faire vivre le bidimensionnel gaussien dans le multidimensionnel.

Soient en effet $x = (x_1,
\dots, x_n)$ des coordonnées (numériques) locales dans lesquelles le
carré $ds^2$ de la longueur de l'élément linéaire $(dx_1, \dots,
dx_n)$ basé au point $x$ s'exprime par l'expression quadratique
$\sum_{ i, j = 1}^n \, g_{ i, j } ( x) \, dx_i dx_j$, le point central
étant l'origine $0 = ( 0, \dots, 0)$.  Riemann commence par effectuer
un simple développement de Taylor à l'ordre deux de tous les
coefficients métriques:
\[
\aligned
g_{i,j}(x)
=
g_{i,j}(0)
&
+
\sum_{k=1}^n\,
\frac{\partial g_{i,j}}{\partial x_k}(0)\,
x_k
+
{\textstyle{\frac{1}{2}}}\,
\sum_{k,l=1}^n\,
\frac{\partial^2g_{i,j}}{\partial x_k\partial x_l}(0)\,
x_kx_l
+
\cdots.
\endaligned
\]
Principe mathématique absolu: contracter, éliminer les termes
superflus, rendre visible l'être dans sa plus simple expression.  Tout
d'abord, la diagonalisation des formes quadratiques définies positives
à coefficients réels permet immédiatement, quitte à effectuer au
préalable un changement linéaire de coordonnées, de supposer qu'on a à
l'origine $g_{ i,j} = \delta_{ i,j}$, d'où:
\[
\sum_{i,j=1}^n\,g_{i,j}(0)\,dx_idx_j
=
dx_1^2+\cdots+dx_n^2. 
\]

\CITATION{
Si l'on introduit ces grandeurs, alors, pour des valeurs infiniment
petites des $x$, le carré de l'élément linéaire sera $= \sum dx^2$; le
terme de l'ordre suivant dans ce carré sera égal à une fonction
homogène du second degré des $n \frac{ n-1}{ 2}$ grandeurs $(x_1 dx_2
- x_2 dx_1)$, $(x_1 dx_3 - x_3 dx_1)$, \dots, c'est-à-dire qu'il sera
un infiniment petit du quatrième ordre; de telle sorte que l'on
obtient une grandeur finie en divisant ce terme par le carré du
triangle infiniment petit dont les sommets correspondent aux systèmes
de valeurs $(0, 0, 0, \dots)$, $(x_1, x_2, x_3, \dots)$, $(dx_1, dx_2,
dx_3, \dots)$ des variables.
\REFERENCE{\cite{riem1898},~p.~289.}}

\noindent
D'après Weyl~\cite{ weyl1990} et Spivak~\cite{ spiv1970}, il
semblerait que le principe de normalisation que Gauss avait élaboré en
termes de coordonnées isothermes ait été repris et généralisé par
Riemann. Aucun élément manuscrit ne nous est parvenu mais l'on peut
penser\footnote{\,
%%%%%%%%%%%%%%%%%%%%%%%-------DEBUT--------%%%%%%%%%%%%%%%%%%%%%%%%%%%
C'est là toute la limite de l'histoire des mathématiques lorsque, trop
pauvre en documents, mais consciente de la complexité des situations
et de la richesse éventuelle des échanges purement verbaux entre
acteurs, elle se trouve réduite à émettre une variété de conjectures
qui finissent à terme par circonscrire toutes les éventualités d'un
réel perdu.
} %%%%%%%%%%%%%%%%%%%%%%%%-----FIN-----%%%%%%%%%%%%%%%%%%%%%%%%%%%%%%%
que Riemann se soit proposé d'examiner, en dimension $n \geqslant 3$,
ce qui devait correspondre au lemme de Gauss sur l'orthogonalité
des géodésiques. En tout cas, notons $M$ la
variété riemannienne, fixons un point $p
\in M$ et considérons $n$ vecteurs $X_1 (p), \dots, X_n (p)$
dans l'espace tangent $T_p M$ à $M$ en $p$ qui forment une base
orthonormée de $T_pM$. L'application exponentielle locale
(\cite{ spiv1970, doca1992}):
\[
\exp
\colon 
T_pM\to M
\]
envoie tout vecteur $X (p)$ de norme riemannienne $g ( X ( p), X (
p))^{ 1/2}$ suffisamment petite sur le point qui est situé à la
distance $g ( X ( p), X ( p))^{ 1/2}$\,\,---\,\,égale 
à cette norme\,\,---\,\,sur l'unique
géodésique issue de $p$ et dirigée par $X (p)$.
Cette application définit un difféomorphisme local de $T_p M$ sur $M$
qui envoie l'origine $0 \in T_p M$ sur $p$. Si l'on note maintenant 
$\psi : T_p M \to \R^n$ l'isomorphisme de $T_p M$ avec $\R^n$ qui est
automatiquement fourni avec la base orthonormale:
\[
\diagram 
T_pM\rto^{\psi} 
\dto_{\exp} 
& \R^n 
\\
M 
\enddiagram
\ \ \ \ \ \ \ \ \
\psi\big(x_1X_1(p)+\cdots+x_nX_n(p)\big)
:=
(x_1,\dots,x_n)
\]
alors l'application $\psi \circ \exp^{ -1}$ fournit un système de
coordonnées locales $(x_1, \dots, x_n)$ sur $M$ qui sont appelées {\sl
coordonnées riemanniennes normales}.  Toute autre base orthonormée de
$T_p M$ fournirait un système essentiellement équivalent de
coordonnées locales. L'avantage principal de ces systèmes de
coordonnées est de donner l'accès le plus direct aux quantités
(tensorielles) de courbure, grâce à l'énoncé suivant, dont nous ne
reconstituerons pas la démonstration.

\smallskip\noindent{\bf Proposition.}
(\cite{ riem1898, weyl1990, spiv1970})
{\em 
Dans tout système de coordonnées riemanniennes normales $(x_1, \dots,
x_n)$, le développement de Taylor en $x = 0$ des coefficients $g_{
i,j} (x)$ de la métrique:
\[
g_{i,j}(x)
=
\delta_{i,j}
+
{\textstyle{\frac{1}{2}}}\,
\sum_{k,l=1}^n\,
\frac{\partial^2g_{i,j}}{\partial x_k\partial x_l}(0)\,
x_kx_l
+
\cdots,
\]
supprime tous les termes $\frac{ \partial g_{i,j }}{ \partial x_k
}(0)$ d'ordre $1$ et fait apparaître des termes d'ordre deux:
\[
{\sf R}_{i,j;\,k,l}
:=
{\textstyle{\frac{1}{2}}}\,
\frac{\partial^2g_{i,j}}{\partial x_k\partial x_l}(0)
\]
qui satisfont les relations de symétrie indicielle
évidentes\footnote{\,
%%%%%%%%%%%%%%%%%%%%%%%-------DEBUT--------%%%%%%%%%%%%%%%%%%%%%%%%%%%
---\,\,puisque $g_{ i,j} = g_{ j,i}$ et
que les deux dérivées partielles $\frac{ \partial}{ \partial x_k}$,
$\frac{ \partial }{ \partial x_l}$ commutent. 
}: %%%%%%%%%%%%%%%%%%%%%%%%-----FIN-----%%%%%%%%%%%%%%%%%%%%%%%%%%%%%%%
\[
{\sf R}_{i,j;\,k,l}
=
{\sf R}_{j,i;\,k,l}
=
{\sf R}_{i,j;\,l,k}
\]
ainsi que les symétries indicielles non triviales: 
\[
\aligned
\text{\bf (i)}
&
\ \ \ \ \ \ \ \ \ \
{\sf R}_{i,j;\,k,l}
=
{\sf R}_{k,l;\,i,j}
\\
\text{\bf (ii)}
&
\ \ \ \ \ \ \ \ \ \
{\sf R}_{i,j;\,k,l}
+
{\sf R}_{i,l;\,j,k}
+
{\sf R}_{i,k;\,l,j}
=
0.
\endaligned
\]
Alors dans ces conditions, le développement de Taylor à l'ordre deux
de la métrique quadratique infinitésimale:
\[
\aligned
\sum_{i,j=1}^n\,
g_{i,j}(x)\,
&
dx_idx_j
=
dx_1^2+\cdots+dx_n^2+
\\
&
\ \ \ \ \ 
+
{\textstyle{\frac{1}{3}}}\,
\sum_{i,j=1}^n\,\sum_{k,l=1}^n\,
{\sf R}_{i,j;\,k,l}\,
\big(x_idx_k-x_kdx_i\big)
\cdot
\big(
x_jdx_l-x_ldx_j
\big)
\endaligned
\]
peut être réécrit, à un facteur $\frac{ 1}{ 3}$ près, comme une forme
quadratique à coefficients ${\sf R}_{ i,j; \, k,l}$ sur les
coordonnées plückériennes:
\[
\left\vert
\begin{array}{cc}
x_{i_1} & dx_{i_1}
\\
x_{i_2} & dx_{i_2}
\end{array}
\right\vert
=
x_{i_1}dx_{i_2}
-
x_{i_2}dx_{i_1}
\ \ \ \ \ \ \ \ \ \ \ \ \
{\scriptstyle{(1\,\leqslant\,i_1\,<\,i_2\,\leqslant\,n)}}
\] 
du $2$-plan engendré par les deux éléments infinitésimaux $(x_1,
\dots, x_n)$ et $(dx_1, \dots, dx_n)$.
}\medskip

Précisions l'interprétation géométrique. Riemann s'imagine un triangle
infiniment petit variable (et surprenant) dans lequel non seulement
$(dx_1,
\dots, dx_n)$, mais encore $(x_1, \dots, x_n)$ sont des quantités
infinitésimales.  Si l'on note donc ce deuxième élément infinitésimal
$(\delta x_1, \dots, \delta x_n)$ avec un symbole $\delta$ pour plus
d'homogénéité conceptuelle, les termes d'ordre deux de la métrique
s'écrivent alors comme une certaine forme quadratique:
\[
{\textstyle{\frac{1}{3}}}\,
\sum_{i,j=1}^n\,\sum_{k,l=1}^n\,
{\sf R}_{i,j;\,k,l}\,
\big(\delta x_idx_k-\delta x_kdx_i\big)
\cdot
\big(
\delta x_jdx_l-\delta x_ldx_j
\big)
\]
en les coordonnées plückériennes $\delta x_{ i_1} dx_{ i_2} - \delta
x_{ i_2} d x_{ i_1}$ du $2$-plan infinitésimal engendré par $( dx_1,
\dots, dx_n)$ et $( \delta x_1, \dots, \delta x_n)$ 
basés au point de référence.  Pour obtenir la courbure de Gauss de la
surface formée des géodésiques dirigées par le $2$-plan $\alpha dx +
\beta \delta x$ (à un facteur constant près), il suffit de diviser
cette expression par le carré de l'aire infinitésimale du triangle
$0$, $dx$, $\delta x$.  Il est quasiment certain que Riemann s'est
inspiré de l'énoncé similaire en dimension $2$ connu par les
continuateurs de Gauss, et nous pouvons conclure que c'est {\em
l'exigence de représentation intuitive d'une multiplicité par des
tranches bidimensionnelles qui a conduit Riemann vers
la courbure sectionnelle}.

\HEAD{Chapitre~1.\,\,\,\,L'ouverture riemannienne}{
1.20.\,\,\,Caractérisation des variétés localement euclidiennes}

\medskip\noindent{\bf 1.20.~Caractérisation
des variétés localement euclidiennes.}
\label{caracterisation-par-courbure-nulle}
Le développement de Taylor dans des coordonnées géodésiques
riemanniennes normales offre l'accès le plus direct aux composantes de
courbure, mais cette approche présente l'inconvénient d'être confinée
à un seul point. Dans la deuxième et dernière partie de la {\em
Commentatio} (\cite{ riem1892}, pp.~380--383),
Riemann cherche à réduire l'équation
différentielle de la conduction de la chaleur: 
\[
\sum_{i}\,
\frac{\partial}{\partial s_i}
\bigg(
\sum_j\,b_{ij}\,\frac{\partial u}{\partial s_j}
\bigg)
=
h\,\frac{\partial u}{\partial t}
\]
à une forme la plus simple possible. Il ramène alors 
ce problème à la transformation d'une
métrique quadratique : 
\[
\sum_{\iota,\,\iota'}\,b_{\iota,\,\iota'}\,
ds_{\iota}ds_{\iota'}
\]
en une métrique plate euclidienne $\sum_{ \iota,
\iota'}\, a_{\iota, \, \iota'}\, 
ds_{ \iota} ds_{ \iota'}$ dont tous les $a_{\iota, \, \iota'}$ sont
constants, métrique qui est donc équivalente à $ds_1^2 + \cdots +
ds_n^2$. Par un calcul assez elliptique, Riemann
trouve la condition nécessaire que 
pour toute collection de quatre indices 
$\iota, \iota', \iota'', \iota''$, 
l'expression suivante: 
\def\theequation{I}\begin{equation}
\aligned
0
&
=
\frac{\partial^2b_{\iota,\,\iota''}}{
\partial s_{\iota'}\partial s_{\iota'''}}
+
\frac{\partial^2b_{\iota',\,\iota'''}}{
\partial s_{\iota}\partial s_{\iota''}}
-
\frac{\partial^2b_{\iota,\,\iota'''}}{
\partial s_{\iota'}\partial s_{\iota''}}
-
\frac{\partial^2b_{\iota',\,\iota''}}{
\partial s_{\iota}\partial s_{\iota'''}}
\\
&
\ \ \ \ \
+
{\textstyle{\frac{1}{2}}}\,
\sum_{\nu,\,\nu'}\,
\big(
p_{\nu,\,\iota',\,\iota'''}\,
p_{\nu',\,\iota,\,\iota''}
-
p_{\nu,\,\iota,\,\iota'''}\,
p_{\nu',\,\iota',\,\iota''}
\big)\,
\frac{\beta_{\nu,\,\nu'}}{B},
\endaligned
\end{equation}
doit s'annuler identiquement, où $B = \det \big( b_{ \iota, \, \iota'}
\big)$ et où $\big( \frac{ \beta_{ \nu, \, \nu'}}{ B} \big)$
est l'inverse de $\big( b_{ \iota, \, \iota'} \big)$. Aucun argument
ne vient supporter l'affirmation implicite que cette condition est
aussi nécessaire, bien que l'\deutschplain{Habilitationsvortrag} ait
énoncé une telle réciproque\footnote{\,
%%%%%%%%%%%%%%%%%%%%%%%-------DEBUT--------%%%%%%%%%%%%%%%%%%%%%%%%%%%
Spivak~\cite{ spiv1970} élabore cinq démonstrations distinctes
de cet énoncé fondamental.

\smallskip\noindent{\bf Théorème.}
{\em Si $(M, g)$ est une variété riemannienne de dimension $n
\geqslant 2$, les trois conditions suivantes sont équivalentes:

\smallskip$\bullet$\,\,
$(M, g)$ est localement isométrique 
à $\R^n$ muni de la métrique euclidienne; 

\smallskip$\bullet$\,\,
toutes les composantes ${A_{ijk}}^l$
du tenseur de courbure s'annulent identiquement; 

\smallskip$\bullet$\,\,
en tout point, la forme quadratique de Riemann-Christoffel 
s'annule identiquement;

\smallskip$\bullet$\,\,
en tout point, la courbure sectionnelle de Riemann-Gauss 
s'annule suivant au moins $\frac{ n ( n-1)}{ 2}$ 
directions superficielles indépendantes.
}\medskip
}. %%%%%%%%%%%%%%%%%%%%%%%%-----FIN-----%%%%%%%%%%%%%%%%%%%%%%%%%%%%%%%
Ensuite, Riemann désigne par la notation:
\[
\big(\iota\iota',\,\,\iota''\iota''')
\] 
le membre de droite de cette équation et 
en quelques lignes <<\,cryptiques\,>>, il prétend 
que la {\em variation seconde} de
la métrique: 
\[
\delta\delta\,\sum\,
b_{\iota,\,\iota'}\,
ds_{\iota}\,ds_{\iota'}
-
2d\delta\,\sum\,b_{\iota,\,\iota'}\,
ds_{\iota}\,ds_{\iota'}
+
dd\,\sum\,b_{\iota,\,\iota'}\,
\delta s_{\iota}\,\delta s_{\iota'}
\]
peut être réécrite de façon à faire apparaître ces expressions à
quatre indices:
\def\theequation{II}\begin{equation}
=
\sum\,
\big(\iota\iota',\,\,\iota''\iota'''\big)\,
\big(ds_\iota\delta s_{\iota'}-ds_{\iota'}\delta s_\iota\big)\,
\big(ds_{\iota''}\delta s_{\iota'''}-ds_{\iota'''}\delta s_{\iota''}\big).
\end{equation}
Plus encore, le quotient de cette expression 
par 

il est un invariant

\def\theequation{II}\begin{equation}
-{\textstyle{\frac{1}{2}}}\,
\frac{
\sum\,
\big(\iota\iota',\,\,\iota''\iota'''\big)\,
\big(ds_\iota\delta s_{\iota'}-ds_{\iota'}\delta s_\iota\big)\,
\big(ds_{\iota''}\delta s_{\iota'''}-ds_{\iota'''}\delta s_{\iota''}\big)
}{
\sum\,b_{\iota,\,\iota'}\,ds_{\iota}\,ds_{\iota'}\,
\sum\,b_{\iota,\,\iota'}\,\delta s_\iota\delta s_{\iota'}
-
\Big(
\sum\,b_{\iota,\,\iota'}\,ds_{\iota}\,\delta s_{\iota'}
\Big)^2}
\end{equation}
Ces affirmations sont énoncées de manière tellement elliptique et sans
aucune vérification que nous n'avons pas d'autre possibilité que de
nous imaginer que Riemann en contrôlait déjà parfaitement la justesse
en dimension deux, grâce aux travaux de Gauss, Minding, Peterson.

%%%%%%%%%%%%%%%%%%%%%%%%%%%%%%%%%%%%%%%%%%%%%%%%%%%%%%%%%%%%%%%%%%%%%

\newpage

\thispagestyle{empty}
\setcounter{footnote}{0}

$\:$

\bigskip\bigskip
\bigskip\bigskip\medskip

\begin{center}
{\large\bf Chapitre~2:

\smallskip
La mobilité helmholtzienne de la rigidité}
\label{Chapitre-2}
\end{center}

\HEAD{Chapitre~2.\,\,\,\,La mobilité helmholtzienne de la rigidité}{
2.1.\,\,\,Le problème de Riemann-Helmholtz}

\bigskip\bigskip

\medskip\noindent{\bf 2.1.~Le problème de Riemann-Helmholtz.}
\label{probleme-riemann-helmholtz}
Dans le plan de son \deutschplain{Habilitationsvortrag}, après avoir
exprimé son intention d'élaborer un concept général de multiplicité
continue, Riemann annonce qu'il se préoccupera d'appliquer ces
considérations abstraites à l'<<\,espace\,>>.  Chez les géomètres du
19\textsuperscript{ème} siècle, ce nom est strictement réservé à
l'<<\,espace physique réel\,>> tridimensionnel et euclidien.

Rappelons notamment qu'au début des années 1820, Gauss a été conduit à
mesurer le très grand triangle géodésique (et terrestre)
Brocken\,\,--\,\,Hohehagen\,\,--\,\,Inselberg dans la campagne du
royaume de Hannovre afin de <<\,tester\,>> la nature euclidienne de
l'espace physique (\cite{ scho2004, scho2009}).
\`A cette époque-là, Gauss venait d'être nommé
conseiller scientifique des gouvernements de Hanovre (dont dépendait
Göttingen) et du Danemark pour l'établissement géodésique du cadastre,
et il était déjà en pleine possession de sa théorie (extrinsèque) des
surfaces.  Aussi Riemann sait-il pertinemment que les propriétés par
lesquelles l'espace physique se distingue de toute autre multiplicité
abstraite à trois dimensions ne peuvent en effet être empruntées qu'à
l'{\em expérience}.

\CITATION{De 
là surgit le problème de rechercher les faits les plus simples au
moyen desquels puissent s'établir les rapports métriques de l'espace,
problème qui, par la nature même de l'objet, n'est pas complètement
déterminé; car on peut indiquer plusieurs systèmes de faits simples,
suffisants pour la détermination des rapports métriques de l'espace.
\REFERENCE{\cite{riem1898},~p.~281.}}

\noindent
Ainsi, pour décider quelle est <<\,la\,>> géométrie de
l'<<\,espace\,>>, l'expérience est susceptible de trancher, mais {\em
à l'avance}, par des considérations de géométrie pure, on peut {\em
rechercher} mathématiquement les hypothèses les plus naturelles et les
plus simples qui puissent caractériser complètement cette notion
d'<<\,espace physique\,>>, notion qui semble nous être donné dans
l'évidence de son euclidéanité présente. C'est donc dans le passage
cité à l'instant, extrait du plan de l'{\em Habilitationsvortrag},
qu'est formulé (implicitement) le célèbre:

\smallskip
\centerline{\fbox{\em Problème de Riemann-Helmholtz}\,.}
\smallskip

\noindent
<<\,{\em Problème}\,>> et
non <<\,théorème\,>>, parce que ni Riemann ni Helmholtz ne l'ont
réellement résolu, ce problème demande précisément par quelles
propriétés l'espace euclidien tridimensionnel muni de la métrique
pythagoricienne pourrait être caractérisé parmi {\em toutes} les
géométries possibles.  Mais au moment de soulever la question, Riemann
était certainement très loin de soupçonner l'incroyable diversité des
géométries que Sophus Lie allait découvrir par des procédés
algébriques uniformes, engendrant malgré lui, à la suite du fameux
{\em Programme d'Erlangen} (\cite{klei1974}), une prolifération de
groupes et de sous-groupes de transformations (Chapitre~5).

{\em Le problème surgit de lui-même}: nécessité sybilline, assertion
méta-mathématique, vision d'une question adéquate.  
Riemann, comme on le sait, suit tous les fils
d'Ariane de l'ouverture. \`A la donation d'un sens univoque dans
l'expérience physique du monde répond donc la surrection articulée des
questionnements mathématiques purs. Fait d'<<\,expérience
mathématique\,>>: les hypothèses abstraites possibles sont soumises à
la variation, à la diversité, à la multiplication, à la ramification,
et à l'éclatement.

\CITATION{
Ces faits, comme tous les faits possibles, ne sont pas nécessaires;
ils n'ont qu'une certitude empirique, ce sont des hypothèses.
\REFERENCE{\cite{riem1898},~p.~281.}}

\HEAD{Chapitre~2.\,\,\,\,La mobilité helmholtzienne de la rigidité}{
2.2.\,\,\,Incomplétudes riemanniennes}

\medskip\noindent{\bf 2.2.~Incomplétudes riemanniennes.}
\label{incompletudes-riemanniennes}
Dans la troisième et dernière partie de son
\deutschplain{Habilitationsvortrag}, Riemann annonce 
en quelques lignes qu'il a complètement résolu\footnote{\,
%%%%%%%%%%%%%%%%%%%%%%%-------DEBUT--------%%%%%%%%%%%%%%%%%%%%%%%%%%%
Deux raisons ont convaincu les mathématiciens des années 1868 à 1890
que le problème n'était en fait pas complètement résolu: 1) les
assertions de Riemann n'ont jamais été suivies de démonstrations
détaillées, même à titre posthume; 2) les raisonnements de Helmholtz
étaient entachés d'erreurs et d'imprécisions mathématiques. }
%%%%%%%%%%%%%%%%%%%%%%%%-----FIN-----%%%%%%%%%%%%%%%%%%%%%%%%%%%%%%%
le problème de caractériser l'espace euclidien standard parmi les
multiplicités à trois dimensions munies d'une métrique quadratique
infinitésimale définie positive.

Première solution, la plus simple et la 
plus économique sur le plan
démonstratif: demander que la mesure de courbure sectionnelle soit 
nulle\footnote{\,
%%%%%%%%%%%%%%%%%%%%%%%-------DEBUT--------%%%%%%%%%%%%%%%%%%%%%%%%%%%
{\em Cf.}~le \S~1.20 ci-dessus. 
} %%%%%%%%%%%%%%%%%%%%%%%%-----FIN-----%%%%%%%%%%%%%%%%%%%%%%%%%%%%%%%
en tout point suivant
trois directions de surface indépendantes\footnote{\,
%%%%%%%%%%%%%%%%%%%%%%%-------DEBUT--------%%%%%%%%%%%%%%%%%%%%%%%%%%%
\`A notre connaissance, l'examen du 
\deutschplain{Nachlass} n'a
fourni aucun document exploitable qui permette aux historiens des
mathématiques de se faire une idée des démonstrations de Riemann.
Contemporain de Klein, le bibliothécaire Distel de l'université de
Göttingen savait en revanche que le
\deutschplain{Nachlass} contenait des
notes intéressantes sur la théorie des nombres et sur la distribution
des zéros de la {\sl fonction zeta de Riemann}, définie par $\zeta (
s) = \sum_{ n \geqslant 1} \,
\frac{ 1}{ n^s}$ pour ${\rm Re}\, s > 1$ et
que Riemann avait introduite dans l'espoir de démontrer que le nombre
des nombres premiers inférieurs à un entier $n \geqslant 2$ se
comporte asymptotiquement comme $\frac{ \log n}{ n}$, ou mieux encore
(Gauss), comme le logarithme intégral $\int_2^n \frac{ 1}{ \log x}\,
dx$. En 1932, après des tentatives de Bessel-Hagen, Siegel a étudié
soigneusement ces notes manuscrites et il en a réorganisé le contenu
dans un mémoire de refonte générale~\cite{sieg1932}. Toutefois,
d'après Siegel: <<\,In Riemanns Aufzeichnungen zur Theorie der
Zetafunktion finden sich nirgendwo druckfertige Stellen; mitunter
stehen zusammenhanglose Formeln auf demselben Blatt; häufig ist von
Gleichungen nur eine Seite hingeschrieben; stets fehlen
Restabschätzungen und Konvergenzuntersuchungen, auch an wesentlichen
Punkten\,>>.
}. %%%%%%%%%%%%%%%%%%%%%%%%-----FIN-----%%%%%%%%%%%%%%%%%%%%%%%%%%%%%%%
Autrement dit, l'espace est localement euclidien {\em si et seulement
si} sa courbure riemannienne s'annule identiquement.
Dans ce cas, le théorème de Pythagore est alors satisfait pour tous
les triangles rectangles finis situés dans n'importe quelle surface
totalement géodésique, et aussi, la somme des angles de tout triangle
quelconque est égale à $\pi$.

Deuxième <<\,solution\,>>: Riemann affirme que l'existence d'un groupe
de mouvements isométriques qui permet de transférer un élément de
surface quelconque basé en un point vers un autre élément de surface
quelconque situé en un autre point arbitraire implique la constante de
la mesure de courbure sectionnelle.

\CITATION{Si
\label{riemann-rigide-mobile}
l'on suppose, en second lieu, comme Euclide, une existence
indépendante de la position, non seulement pour les lignes, mais
encore pour les corps, il s'ensuit que la mesure de courbure est
partout constante, et alors la somme des angles est déterminée dans
tous les triangles, lorsqu'elle l'est dans un seul.
\REFERENCE{\cite{riem1898},~p.~294.}}

\noindent
L'argument intuitif est simple: un corp rigide 
bidimensionnel
maximalement
mobile suffira pour propager partout dans la
multiplicité la courbure sectionnelle qu'il possède.
Toutes les surfaces géodésiques en tout point orientées
dans toute direction ont la même
courbure: la
courbure est {\em constante}. 
Alors dans cette circonstance et dans des coordonnées géodésiques
normales adéquates\footnote{\,
%%%%%%%%%%%%%%%%%%%%%%%-------DEBUT--------%%%%%%%%%%%%%%%%%%%%%%%%%%%
Dans~\cite{ weyl1990}, Weyl a reconstitué les calculs 
qui auraient pu conduire Riemann à cette expression. 
}, %%%%%%%%%%%%%%%%%%%%%%%%-----FIN-----%%%%%%%%%%%%%%%%%%%%%%%%%%%%%%%
Riemann représente la métrique par une formule (déjà mentionnée)
invariante symétrique par rapport à toutes les variables:
\[
ds^2 
= 
\frac{dx_1^2+\cdots+dx_n^2}{
\big[1+\frac{\kappa}{4}\,(x_1^2+\cdots+x_n^2)\big]^2}, 
\]
dans laquelle $\kappa$ est une {\em constante réelle}, que l'on peut
même supposer après dilatation être égale à $-1$ (géométrie
hyperbolique), à $0$ (géométrie euclidienne) ou à $1$ (géométrie
sphérique). Dans un autre passage, Riemann fournit quelques
explications.

\CITATION{
Le caract\`ere commun de ces vari\'et\'es dont la mesure de courbure
est constante, peut aussi \^etre exprim\'e en disant que les figures
peuvent se mouvoir\footnotemark sans
\'elargissement \deutsch{\sf Dehnung} en elles. Car il est \'evident
que les figures en elles ne pourraient pas coulisser
\deutsch{\sf verschiebar sein}
et pivoter \deutsch{\sf drehbar sein} librement, si la mesure de
courbure n'\'etait pas la m\^eme en chaque point et dans toutes les
directions. Mais d'autre part, les rapports m\'etriques de la
vari\'et\'e sont compl\`etement d\'etermin\'es par la mesure de
courbure; donc les rapports m\'etriques autour d'un point et dans
toutes les directions sont exactement les m\^emes qu'autour d'un autre
point, et par cons\'equent, \`a partir de ce premier point, les
m\^emes constructions peuvent \^etre transf\'er\'ees, d'o\`u il
s'ensuit que, dans les vari\'et\'es dont la mesure de courbure est
constante, on peut donner aux figures chaque position quelconque.
\REFERENCE{\cite{riem1898},~p.~281.}}

\footnotetext{\,
%%%%%%%%%%%%%%%%%%%%%%%-------DEBUT--------%%%%%%%%%%%%%%%%%%%%%%%%%%%
Note de Engel et Lie: <<\,Ici \`a vrai dire, Riemann
aurait m\^eme d\^u ajouter le mot `librement'\,>>.
} %%%%%%%%%%%%%%%%%%%%%%%%-----FIN-----%%%%%%%%%%%%%%%%%%%%%%%%%%%%%%%

\noindent
Toutefois, au-delà du niveau intuitif, la conceptualisation
mathématique rigoureuse de ces affirmations demeure essentiellement
problématique pour Engel et pour Lie. Riemann semble exprimer que
l'exigence d'apr\`es laquelle la mesure de courbure doit \^etre
partout constante possède la m\^eme signification que certaines
exigences concernant la mobilit\'e des figures. Dans un premier temps,
Engel et Lie cherchent à restituer l'encha\^{\i}nement des idées de
Riemann d'une mani\`ere quelque peu plus pr\'ecise, <<\,quoique non
absolument pr\'ecise\,>>, ajoutent-t-ils immédiatement.

\CITATION{
Riemann cherche, parmi les vari\'et\'es dont la longueur d'un
\'el\'ement courbe a la forme~\thetag{ 5}, toutes celles dans
lesquelles les figures peuvent occuper chaque position quelconque,
c'est-\`a-dire, dans lesquelles les figures peuvent coulisser et
tourner, sans subir d'\'elargissement. Il parvient \`a ce r\'esultat que
les vari\'et\'es dont la mesure de courbure est constante en tous
lieux sont les seules dans lesquelles les figures sont mobiles de
cette mani\`ere.
\REFERENCE{p.~\pageref{EL-RH}~ci-dessous.}}

\noindent
Plusieurs problèmes de conceptualisation se dessinent donc: qu'est-ce
que le <<\,mouvement\,>>?  Qu'entend-on par <<\,figures\,>>?  Que veut
dire <<\,coulisser\,>>?  Que veut dire <<\,tourner\,>>?  Que veut dire
<<\,sans subir d'élargissement\,>>? C'est le physicien allemand
Hermann von Helmholtz qui va tenter en 1868 de donner un sens
mathématique précis à ces notions intuitives.

\HEAD{Chapitre~2.\,\,\,\,La mobilité helmholtzienne de la rigidité}{
2.3.\,\,\,Rendre objectives les propriétés de la géométrie}

\medskip\noindent{\bf 2.3.~Rendre objectives les propositions
de la géométrie.}
\label{objectiver-la-geometrie}
Selon Lie\footnote{\,
%%%%%%%%%%%%%%%%%%%%%%%-------DEBUT--------%%%%%%%%%%%%%%%%%%%%%%%%%%%
En parallèle, le lecteur pourra découvrir les {\em Remarques
préliminaires} à la Division~V, traduites en français
p.~\pageref{Division-5} sq. ci-dessous. },
%%%%%%%%%%%%%%%%%%%%%%%%-----FIN-----%%%%%%%%%%%%%%%%%%%%%%%%%%%%%%%
le mérite principal de Helmholtz par rapport à Riemann est d'avoir
expressément construit la géométrie en stipulant un certain nombre
d'{\sl axiomes} (non problématisants) qui se rapportent à des
collections de points {\em finiment éloignés les uns des autres}, sans
utiliser la notion d'élément d'arc infinitésimal, et sans admettre de
pouvoir disposer de la théorie de l'intégration. Autre mérite de cette
méthode neuve: Helmholtz est le premier à avoir opéré avec la famille
des mouvements de l'espace (à trois dimensions), tout en
l'interprétant comme famille de transformations de la variété
numérique sous-jacente. On trouve en effet chez Helmholtz les tous
premiers concepts embryonnaires pour l'édification d'une théorie
mathématique des groupes de mouvements dans l'espace, et il est
très probable que son mémoire a été la pierre angulaire d'inspiration, de
motivation, et de problématisation lorsque, deux ou trois années plus
tard, Klein et Lie collaborèrent à l'élaboration du {\em Programme
d'Erlangen} (\cite{ klei1974})

Fortement stimulé par la publication posthume très récente de
l'\deutschplain{Habilitations\-vortrag} de Riemann\footnote{\,
%%%%%%%%%%%%%%%%%%%%%%%-------DEBUT--------%%%%%%%%%%%%%%%%%%%%%%%%%%%
En fait, il obtint une copie des notes manuscrites que Schering avait
prises sur le travail de Riemann. Schering s'est vu attribuer la chaire
de Riemann en 1866.
}, %%%%%%%%%%%%%%%%%%%%%%%%-----FIN-----%%%%%%%%%%%%%%%%%%%%%%%%%%%%%%%
Helmholtz déclare dans son mémoire~\cite{ helm1868b} de juin 1868
qu'il conduit en privé depuis de nombreuses années des réflexions
techniques sur la notion d'<<\,espace\,>>.  Il se félicite aussi de
trouver comme compagnon de pensée, dans ces investigations délicates
sur les fondements de la géométrie, un mathématicien de l'envergure de
Riemann, et il propose de présenter à la
\deutschplain{königlische Gesellschaft der Wissenschaften zu
Göttingen} un système alternatif d'axiomes afin de développer une
autre vision du problème mathématique de l'espace.

Pour Helmholtz qui a travaillé sur la physiologie et l'anatomie de
l'{\oe}il, la question de l'origine et de la nature essentielle de nos
intuitions générales de l'espace est primordiale.  Plus précisément,
il s'agit de déterminer dans quelle mesure les propositions de la
géométrie possèdent un sens qui est {\em objectivement} valide. Et \`a
l'opposé, il s'agit de savoir aussi quelle est la fraction des définitions
et des conséquences des définitions qui dépend seulement de
l'abstraction des descriptions mathématiques. Ces deux questions
helmholtziennes ne peuvent certainement pas recevoir de réponse
simple, car il y a un cercle mystérieux de correspondances intuitives
entre les représentations idéales des objets géométriques purs et leur
portrait physique approximatif, imparfait et défectueux.

En tout cas pour Helmholtz\,\,---\,\,que 
l'on a pu considérer comme {\em résolument empiriste} (\cite{
darr2003, darr2007})\,\,---, 
notre intuition archaïque, concrète, évidente et
profondément {\em non problématique} des objets dans le monde physique
{\em doit être incarnée fidèlement dans les sytèmes axiomatiques de la
géométrie pure}. Le point de vue de Helmholtz est donc celui d'un
penseur universel des sciences expérimentales qui admet deux faits
fondamentaux transmis par l'expérience:

\smallskip{\bf 1)}\,\,
la mesure de l'espace est et ne peut être basée que sur l'observation
de {\em congruences} entre des objets\footnote{\,
%%%%%%%%%%%%%%%%%%%%%%%-------DEBUT--------%%%%%%%%%%%%%%%%%%%%%%%%%%%
Par mise en congruence, il faut entendre tout procédé concret par
lequel les deux extrémités d'un objet-étalon sont placées en
coïncidence avec une paire de points repérables sur l'objet à mesurer;
ensuite si nécessaire, l'étalon invariable de mesure doit être reporté
jusqu'à épuisement de l'extension de l'objet, un certain nombre de
fois en fonction de la taille de l'objet. 
Lorsque la longueur cherchée n'est
pas un multiple entier de l'objet-règle utilisé,
graduations et subdivisions
complètent alors le défaut restant. L'horizon 
microscopique borne rapidement l'effort de raffinement dans la
mesure.  Helmholtz affirme donc clairement que l'ontologie
sous-jacente au mesurable repose sur l'équivalence, sur la
comparaison, sur la sommation finie et sur la (sub)division. Le
caractère archimédien de l'espace physique 
provient de ce qu'il héberge des
essences rigides {\em arithmétisantes}.
}; %%%%%%%%%%%%%%%%%%%%%%%%-----FIN-----%%%%%%%%%%%%%%%%%%%%%%%%%%%%%%%

\smallskip{\bf 2)}\,\,
toute forme de congruence demande pour son effectuation qu'il existe
des corps rigides, fixes en eux-mêmes, et qu'on puisse les déplacer
comme des règles sans qu'ils changent ni de forme, ni de
longueur\footnote{\,
%%%%%%%%%%%%%%%%%%%%%%%-------DEBUT--------%%%%%%%%%%%%%%%%%%%%%%%%%%%
Sinon, si des corps mobiles devaient nécessairement changer de
longueur lorsqu'ils sont déplacés d'un lieu vers un autre\,\,---\,\,la
modification de leur longueur pouvant d'ailleurs même {\em a priori}
dépendre du chemin qu'on leur fait emprunter\,\,---, il serait
absolument impossible de parler d'une longueur comme résultat
invariable d'une pluralité de mesures.  <<\,Les axiomes géométriques
ne concernent pas les relations spatiales {\em per se}, mais ils
impliquent aussi le comportement mécanique de nos corps {\em les plus
rigides} lorsqu'ils sont déplacés\,>> (cité p.~531 de~\cite{
darr2007}).  }. %%%%%%%%%%%%%%%%%%%%%%%%-----FIN-----%%%%%%%%%%%%%%%%%%%%%%%%%%%%%%% 

\smallskip\noindent
L'espace ubiquitaire est comme le milieu terrestre: baigné d'une
atmosphère invisible, il n'oppose en principe aucune résistance à la
libre mobilité des règles mesurantes. Corps mesurant et corps mesuré
sont structurellement homologues, voire interchangeables. Dans ce
monde homogène de corps rigides, l'essence fondamentale est la
distance; c'est aussi un principe d'équivalence ontologique.

On a donc affaire à une théorie essentiellement
active de la mesure qui semble s'imaginer comme absolument nécessaire
un sujet porteur muni de règles et de réglettes, forcé de parcourir
l'espace, depuis ses dimensions microscopiques jusqu'à ses dimensions
astronomiques, afin de s'assurer par des moyens rigoureusement
élémentaires que les distances relatives entre certains villages, que
la superficie d'un appartement, que les dimensions d'un meuble à la
vente, ou celles d'une boîte à bijoux, sont exactes au
mètre près, au centimètre près, au millimètre près\footnote{\,
%%%%%%%%%%%%%%%%%%%%%%%-------DEBUT--------%%%%%%%%%%%%%%%%%%%%%%%%%%%
Ce point vue qui s'imagine une universalisation possible de
l'empiricité archaïque est à distinguer nettement du point de vue de
Riemann, pour lequel la métrique est {\em intrinsèquement donnée dans
l'infinitésimal} et ne nécessite {\em aucun parcours}, ni même aucun
acte de mesure\,\,---\,\,si ce n'est pour estimer mathématiquement,
par intégration, la longueur d'une ligne courbe finie\,\,---\,\,
puisque {\em chaque différentielle infinitésimale $dx$ est d'emblée
accompagnée de sa longueur attitrée $\big( \sum_{ i, j = 1}^n \, g_{
ij} ( x ) \, dx_i dx_j \big)^{ 1/2}$}. Point crucial: {\em aucun
déplacement de corps ou de réglette n'est requis}. Plus précisément:
dans toute métrique riemannienne ou finslérienne, l'essence métrique
est déjà décalquée sur l'étendue comme le serait une infinité
d'étiquettes millimétriques disposées en tout point, et dans toutes
les directions issues d'un point. Mariage entre l'{\em intrinsèque
pur} et le champ physique possible: <<\,Il faut donc, ou que la
réalité sur laquelle est fondé l'espace forme une variété discrète, ou
que le fondement des rapports métriques soit cherché en dehors de lui,
dans les forces de liaison qui agissent en lui\,>> (\cite{ riem1898},
p.~297).
}. %%%%%%%%%%%%%%%%%%%%%%%%-----FIN-----%%%%%%%%%%%%%%%%%%%%%%%%%%%%%%%

Argument surprenant s'il en est: Helmholtz écarte la rectilinéarité
des rayons de lumière comme principe de référence pour la mesure,
parce que les aveugles à qui les vérités de la géométrie sont bien
connues n'ont accès qu'aux congruences.

\CITATION{Donc 
je suppose dès le début que la mesure de l'espace soit possible
par des congruences déterminées et je me propose comme but la tâche de
rechercher la forme analytique la plus générale d'une variété
multiplement étendue dans laquelle soient possibles les mouvements
ayant la constitution ainsi demandée.
\REFERENCE{\cite{helm1977},~p.~41.}}

\noindent
En fait, Helmholtz délimite pour ses recherches un objectif encore
plus ciblé, moins susceptible de l'exposer à l'inattendu mathématique:
il cherche en effet à formuler des hypothèses physiquement évidentes
(qu'on appellerait aujourd'hui {\sl axiomes}) afin de retrouver
seulement l'espace euclidien tridimensionnel comme concept compris
sous de tels axiomes. Approche axiomatique ou recherche d'un théorème
d'unicité? L'ambiguïté est réelle, d'autant plus que la mise au point
d'hypothèses spécifiques qui impliquent l'euclidéanité sera
complètement absorbée par la théorie des groupes continus, qui
montrera comment l'engendrement du divers déborde l'ontologie naïve
des unicités initialement espérées (\cf~le Chapitre~20
p.~\pageref{Chapitre-20} ci-dessous).

\HEAD{Chapitre~2.\,\,\,\,La mobilité helmholtzienne de la rigidité}{
2.4.\,\,\,Les quatre axiomes de Helmholtz}

\medskip\noindent{\bf 2.4.~Les quatre axiomes de Helmholtz.}
\label{quatre-axiomes-helmholtz}
Helmholtz reconnaît que son approche qui consiste à introduire
d'emblée la restriction de libre mobilité de corps rigides dans
l'espace embrasse beaucoup moins de concepts que l'analyse
problématisante de Riemann (\cf~le Chapitre~1).  Toutefois, il semble
être embarrassé par l'extrême généralité des considérations abstraites
de Riemann, qu'il préfère envisager comme une tentative de
caractériser l'espace euclidien tridimensionnel en recherchant les
meilleurs hypothèses suffisantes.  Ainsi a-t-il étudié de près la
restriction finale de mobilité introduite par Riemann pour distinguer
l'espace physique des autres multiplicités possibles, et il s'est
interrogé sur la possibilité de placer au fondement de la géométrie ce
qu'il considérait comme la conclusion principale de la leçon d'épreuve
de Riemann, à savoir: les propriétés de libre mobilité dont Riemann
prétendait qu'elles impliquent la constance de la courbure
sectionnelle.  En renversant complètement l'ordre de pensée
riemannien\footnote{\,
%%%%%%%%%%%%%%%%%%%%%%%-------DEBUT--------%%%%%%%%%%%%%%%%%%%%%%%%%%%
Chercher dans l'ouverture, c'est s'élever dans un arbre exponentiel
truffé de branchements imprévus et de bifurcations indécises. Prendre
connaissance de vérités par l'axiomatique {\em a posteriori}, c'est
descendre sans choix dans des branches qui canalisent la pensée.  },
%%%%%%%%%%%%%%%%%%%%%%%%-----FIN-----%%%%%%%%%%%%%%%%%%%%%%%%%%%%%%%
Helmholtz procède donc presque d'une manière purement axiomatique au
sens moderne du terme.  Car en effet, au début de son mémoire~\cite{
helm1977}, quatre <<\,axiomes\,>> concernant la notion
d'espace\footnote{\,
%%%%%%%%%%%%%%%%%%%%%%%-------DEBUT--------%%%%%%%%%%%%%%%%%%%%%%%%%%%
En parallèle, le lecteur pourra découvrir la traduction complète en
français de ces axiomes, reproduits par Engel et Lie dans le
\S~91 du Chapitre~21, p.~\pageref{S-91}
ci-dessous; ensuite, le \S~92 p.~\pageref{S-92} en offre une
formulation mathématique très précise.  },
%%%%%%%%%%%%%%%%%%%%%%%%-----FIN-----%%%%%%%%%%%%%%%%%%%%%%%%%%%%%%%
sont explicitement formulés et admis comme hypothèses dans les
démonstrations qui suivent.

\smallskip{\sf I: Axiome de multiplicité numérique et de
continuité des mouvements.}\, Les éléments individuels (points) d'une
multiplicité à $n$ dimensions peuvent être spécifiés par la {\em
mesure}\footnote{\,
%%%%%%%%%%%%%%%%%%%%%%%-------DEBUT--------%%%%%%%%%%%%%%%%%%%%%%%%%%%
Chez Helmholtz (\cite{ darr2007}), les coordonnées ont un sens {\em
métrique}: elles sont d'emblées
accompagnée d'un étalon d'unité de mesure
qui est représenté par la quantité $1$. }
%%%%%%%%%%%%%%%%%%%%%%%%-----FIN-----%%%%%%%%%%%%%%%%%%%%%%%%%%%%%%%
de exactement $n$ grandeurs numériques réelles $(x_1, \dots, x_n)$.
Le mouvement d'un point est représenté par une modification continue,
et même suffisamment différentiable si nécessaire, de ses coordonnées.
Tout est essentiellement local, implicitement plongé dans $\R^n$, sans
phénomène topologique, sans discontinuité, sans condition au bord,
sans contrainte spécifique.  L'étude des exceptions est
potentiellement réservée à des analyses ultérieures.

\smallskip{\sf II: Axiome de l'existence des corps
rigides mobiles.}  Cet axiome exprime l'idée nouvelle principale de
Helmhotz. Pour que l'espace soit physiquement mesurable, on présuppose
l'existence de systèmes de points rigides en eux-mêmes, et néanmoins
susceptibles d'être déplacés (axiome~III): 

\smallskip

\centerline{\em mobilité de la rigidité.}  

\smallskip\noindent
Mais quelle est ici la définition mathématique précise de
la <<\,rigidité\,>>?  Impossible de demander qu'il s'agisse d'une {\em
distance}, au sens euclidien ou métrique du terme, puisque c'est
justement l'euclidéanité locale que Helmholtz cherche à démontrer à
partir d'axiomes abstraits et naturels. 

Herméneutique {\em indécise} de la position
d'hypothèses, ou nécessité de réaliser une
{\em assomption en généralité} pour engendrer des essences
hypothétiques d'ordre supérieur: il faut donc
formuler un axiome qui ne demande presque rien d'explicite, et qui
s'inscrive {\em a priori} dans un très grand univers de potentialités
mathématiques.

L'espoir, {\em via} une démonstration synthétique éventuellement
longue et difficile, est de se persuader qu'une notion très archaïque
de <<\,stabilité-rigidité\,>> dans le mouvement redonnera celle qui
semble la plus naturellement transmise dans l'expérience physique, à
savoir la rigidité euclidienne. On est en pleine exploration
métaphysique des conditions d'aprioricité du savoir
physico-mathématique. On recherche en quelque sorte des signes qui
confirmeraient l'existence de vérités profondes (et vraisemblablement
cachées), lesquelles vérités profondes pourraient répondre à la
question très importante d'un <<\,pourquoi\,>>, à savoir: <<\,{\em
Pourquoi l'espace physique est-il tridimensionnel et euclidien}\,>>\,?
Avec son hypothèse de corps rigides, Helmholtz recherche donc un
<<\,parce que\,>> qui s'élève plus haut dans l'échelle métaphysique
qu'une simple confirmation expérimentale basée sur des mesures
microscopiques ou sur des mesures astronomiques.

La définition de Helmholtz est à la fois intuitive et très générale:
il doit exister une certaine fonction {\em à deux arguments} qui est
définie sur la totalité de toutes les paires de points appartenant au
corps rigide et qui satisfait la propriété d'invariance suivante: pour
toute paire de points fixée à l'avance, la valeur de cette fonction
sur les deux points en question devra rester {\em invariable} au cours
de tous les mouvements possibles (axiome~III) du corps en question.

Plus encore que dans la théorie gaussienne des formes quadratiques à
coefficients entiers, ou dans la théorie galoisienne des substitutions
de racines, c'est vraiment dans l'univers des {\em mouvements
continus} physico-mathématiques que l'Idée métaphysique d'{\sl
invariance} trouve son origine la plus profonde. Lie conceptualisera
cette Idée dominatrice dans de nombreuses branches de sa théorie des
groupes continus (équation invariante; système d'équations aux
dérivées partielles invariantes; algèbres classifiantes d'invariants
différentiels) en plaçant le concept à un très haut niveau
d'abstraction: l'invariant peut être une fonction arbitraire, réelle
ou complexe, souvent d'un très grand nombre d'arguments.  Et dans les
diverses solutions au problème de Riemann-Helmholtz qu'il expose avec
Engel, Lie formulera la notion d'invariant\footnote{\,
%%%%%%%%%%%%%%%%%%%%%%%-------DEBUT--------%%%%%%%%%%%%%%%%%%%%%%%%%%%
{\em Voir} l'équation~\thetag{ 3} au début du Chapitre~20 et aussi la
condition {\bf C)} p.~\pageref{452}). 
}, %%%%%%%%%%%%%%%%%%%%%%%%-----FIN-----%%%%%%%%%%%%%%%%%%%%%%%%%%%%%%%
{\em sans aucune référence implicite à la notion de distance}: pas de
valeurs positives, pas d'inégalité du triangle.

Plus encore, Lie extraira de cet axiome helmholtzien plusieurs
présuppositions implicites situées à l'interface entre l'axiome~II et
l'axiome~III, et que Helmholtz déduisait à tort de ses axiomes
incomplètement formulés: 1) l'invariant doit être non dégénéré; 2)
lorsqu'on fixe plusieurs points en position générale, les équations
d'invariance relativement à ces points doivent être mutuellement
indépendantes; 3) les invariants de toute collection de $s > 2$ points
se déduisent des invariants entre paires de points, 
c'est-à-dire: il n'y a pas
d'autres invariants que ceux qui existent entre les paires de points.

\smallskip{\sf III: Axiome de libre mobilité des corps
rigides.}
En s'inspirant de Riemann\footnote{\,
%%%%%%%%%%%%%%%%%%%%%%%-------DEBUT--------%%%%%%%%%%%%%%%%%%%%%%%%%%%
{\em Voir} la citation p.~\pageref{riemann-rigide-mobile}. 
}, %%%%%%%%%%%%%%%%%%%%%%%%-----FIN-----%%%%%%%%%%%%%%%%%%%%%%%%%%%%%%% 
Helmholtz demande que les corps rigides puissent être déplacés
continûment en tout lieu\footnote{\,
%%%%%%%%%%%%%%%%%%%%%%%-------DEBUT--------%%%%%%%%%%%%%%%%%%%%%%%%%%%
Cette condition est suffisante pour l'effectuation de mesures
physiques au niveau macroscopique. Toutefois, de nombreux
commentateurs (\cf~par exemple
\cite{ helm1977, darr2003, darr2007}) ont insisté sur le fait que
seule la mobilité de {\em règlettes} unidimensionnelles, et non la
mobilité de corps rigides arbitraires, peut être invoquée comme {\em
physiquement nécessaire} pour effectuer des mesures de distance dans
le monde. Ce point de vue laisse néanmoins complètement de côté la
question de savoir ce qu'il faut entendre {\em abstraitement} par un
déplacement de réglettes, et semble s'en remettre trop aisément à une
intuition approximative de l'arpentage. Seule la théorie abstraite des
groupes continus de transformations conceptualise les mouvements
possibles des figures, qu'elles soient règles ou corps, finies ou
infinitésimales.
}, %%%%%%%%%%%%%%%%%%%%%%%%-----FIN-----%%%%%%%%%%%%%%%%%%%%%%%%%%%%%%%
c'est-à-dire qu'un point quelconque spécifié à l'avance dans un corps
rigide peut être transféré vers tout autre point de l'espace au moyen
d'au moins un mouvement continu qui déplace le corps dans son
intégalité tout en respectant sa rigidité interne.  Si l'on fait
abstraction des conditions aux limites et des discontinuités
éventuelles, les seules contraintes qui pourraient faire obstacle au
mouvement ne peuvent provenir que des équations qui sont formées avec
la fonction invariante entre paires de points, par exemple si l'on
demande qu'un, ou deux, ou trois points, ou plus, du corps rigide
restent entièrement fixés au cours du mouvement.  Ensuite, Helmholtz
entreprend de <<\,démontrer\footnote{\,
%%%%%%%%%%%%%%%%%%%%%%%-------DEBUT--------%%%%%%%%%%%%%%%%%%%%%%%%%%%
{\em Voir} le \S93 p.~\pageref{S-93} ci-dessous pour une
critique de Lie et Engel.
}\,>>  %%%%%%%%%%%%%%%%%%%%%%%%-----FIN-----%%%%%%%%%%%%%%%%%%%%%%%%%%%
que les mouvements d'un corps rigide comprennent exactement $\frac{ n
( n+1)}{ 2}$ degrés de libertés, et donc notamment $6$ dans le cas de
l'espace physique euclidien (translations: $3$ paramètres; rotations:
$3$ paramètres).

Paradoxe: la notion de corps rigide continu comprenant une {\em
infinité} de points en cohésion n'est donc pas pleinement utilisée
dans ces raisonnements. Six points au maximum sont à distinguer.  En
tout cas, Engel et Lie démontreront rigoureusement cet énoncé sur le
nombre $\frac{ n ( n+1)}{ 2}$ de degrés de liberté, en analysant
finement les hypothèses qui étaient implicites dans les axiomes~II
et~III.  Par ailleurs, en supposant beaucoup moins que
Helmholtz\footnote{\,
%%%%%%%%%%%%%%%%%%%%%%%-------DEBUT--------%%%%%%%%%%%%%%%%%%%%%%%%%%%
<<\,Le procédé de Lie a lieu entièrement
{\em a priori}\,>>, écrit Jules Vuillemin
p.~420 de~\cite{ vui1962}. 
}, %%%%%%%%%%%%%%%%%%%%%%%%-----FIN-----%%%%%%%%%%%%%%%%%%%%%%%%%%%%%%%
le Chapitre~20 p.~\pageref{Chapitre-20} ci-dessous commencera par une
étude purement abstraite des groupes continus de transformations pour
lesquels deux points ont un, et un seul invariant, tandis que $s > 2$
points n'ont pas d'autre invariant que ceux qui se déduisent des
paires de points qui sont contenus en eux, et cela, d'abord dans le
domaine complexe, puis dans le domaine réel. Cette étude abstraite
extrêmement recherchée sur le plan mathématique 
repose sur un très grand nombre de résultats contenus
dans les trois volumes de la {\em Theorie der
Transformationsgruppen}.

\smallskip{\sf IV: Axiome de la monodromie ou de la
périodicité des rotations des corps rigides.}  C'est l'axiome le plus
controversé, car on 
verra\footnote{\,
%%%%%%%%%%%%%%%%%%%%%%%-------DEBUT--------%%%%%%%%%%%%%%%%%%%%%%%%%%%
En infinitésimalisant le problème, Helmholtz va se ramener à
considérer un système d'équations différentielles ordinaires $\frac{
dx_i}{ dt} = \sum_{ j=1}^3 \, a_{ ij}\, x_j$ d'ordre $1$ à
coefficients $a_{ ij}$ constants, et l'axiome de monodromie va
directement impliquer que la matrice $(a_{ ij})_{ 1 \leqslant i
\leqslant 3}^{ 1 \leqslant j \leqslant 3}$ possède
une valeur propre nulle, et deux valeurs propres imaginaires
conjuguées $\omega >0$ et $- \omega < 0$, ce qui fait que la matrice
$a_{ ij}$ représente tout
simplement une rotation dans l'espace {\em euclidien} à
trois dimensions. }
%%%%%%%%%%%%%%%%%%%%%%%%-----FIN-----%%%%%%%%%%%%%%%%%%%%%%%%%%%%%%%
qu'il exclut trop aisément un très grand nombre de géométries
possibles, au sens de Klein et Lie, c'est-à-dire de groupes de
transformations <<\,exotiques\,>>.  Lorsque $(n-1)$ points d'un corps
rigide sont fixés en position générale, Helmholtz prétend (comme
conséquence de l'Axiome~III) qu'il reste encore un, et un seul degré
de liberté autorisant un mouvement continu; ce mouvement dépend alors
d'un, et d'un seul paramètre, et nous venons de signaler que Engel et
Lie ont clarifié cette affirmation.  L'axiome de monodromie demande
alors que la <<\,{\em rotation}\footnote{\,
%%%%%%%%%%%%%%%%%%%%%%%-------DEBUT--------%%%%%%%%%%%%%%%%%%%%%%%%%%%
La terminologie utilisée par Helmholtz montre que la démonstration
mathématique {\em a priori} qu'il met en {\oe}uvre est sous-tendue,
dans l'intuition explorante, par l'idée que l'on a déjà affaire au
groupe orthogonal euclidien. Pour être rigoureux, il faudrait
qualifier ce mouvement non pas de <<\,rotation\,>>, mais de mouvement
continu à un paramètre.  }\,>>
%%%%%%%%%%%%%%%%%%%%%%%%-----FIN-----%%%%%%%%%%%%%%%%%%%%%%%%%%%%%%%
(sans retour en arrière) du corps rigide autour de $(n-1)$ de ses
points supposés fixés doit le reconduire, au bout d'un temps fini, à
sa position initiale: chaque point en lui revient coïncider avec la
position qu'il occupait au début. Cet axiome qui exclut donc tout
mouvement {\em en spirale} (contraction ou dilatation des
<<\,longueurs\,>> après un tour) et tout mouvement {\em en hélice}
(décalage le long d'un axe après un tour) semble parler un langage
évident pour l'intuition euclidienne.  Engel et Lie reformuleront cet
axiome de manière un peu plus précise\footnote{\,
%%%%%%%%%%%%%%%%%%%%%%%-------DEBUT--------%%%%%%%%%%%%%%%%%%%%%%%%%%%
{\em Voir} la condition {\bf E)} p.~\pageref{condition-E} ci-dessous, 
qui exprime tout simplement que le mouvement à un paramètre
encore possible est {\em périodique}.}
%%%%%%%%%%%%%%%%%%%%%%%%-----FIN-----%%%%%%%%%%%%%%%%%%%%%%%%%%%%%%%
en utilisant la notion de groupe à un paramètre, essentiellement
absente du mémoire de Helmholtz.

\HEAD{Chapitre~2.\,\,\,\,La mobilité helmholtzienne de la rigidité}{
2.5.\,\,\,Linéarisation de l'isotropie}

\medskip\noindent{\bf 2.5.~Linéarisation de l'isotropie.}
\label{linearisation-isotropie}
Les métriques quadratiques de Gauss $E(u,v) \, du^2 + 2
\, F (u,v) \, du dv + G (u, v)\, dv^2$ et de Riemann $\sum_{ i, j =
1}^n \, g_{ ij} (x) \, dx_i dx_j$ que Helmholtz cherche à ressaisir
s'expriment en termes {\em infinitésimaux}, et c'est certainement pour
cette raison que Helmholtz a considéré comme {\em naturel}
d'interpréter tous ses axiomes directement dans l'infiniment
petit\footnote{\,
%%%%%%%%%%%%%%%%%%%%%%%-------DEBUT--------%%%%%%%%%%%%%%%%%%%%%%%%%%%
Ici se trouve son erreur d'inadvertance principale,
\cf~le \S~2.6 ci-dessous. 
}. %%%%%%%%%%%%%%%%%%%%%%%%-----FIN-----%%%%%%%%%%%%%%%%%%%%%%%%%%%%%%%

\CITATION{J'utiliserai 
les hypothèses II, III et IV seulement pour 
des points dont les différences de cordonnées 
infiniment petites. Aussi la congruence indépendante
des limites sera supposée valide seulement
pour des éléments spatiaux
\REFERENCE{\cite{helm1977},~pp.~44--45.}}

\noindent
Examinons donc comment Helmholtz procède sur le plan mathématique.
Soient $(u, v, w)$ les coordonnées d'un point appartenant au corps
rigide dans une première situation de ce corps, et soient $(r, s, t)$
les coordonnées d'un même point dans une seconde situation du corps.
Alors ces coordonnées $(r, s, t)$ dépendent en toute généralité de
$(u, v, w)$ et de six constantes arbitraires qui 
expriment les degrés de libertés. Helmholtz sous-entend
que $(u, v, w)
\longmapsto (r, s, t)$ est un difféomorphisme\,\,---\,\,c'est
une conséquence de la rigidité\,\,---\,\,et il écrit la transformation
correspondante entre différentielles, donnée par une matrice
jacobienne. Ensuite, il fixe le point $(r, s, t)$ et il considère
seulement toutes les transformations encore possibles du corps rigide
qui envoient ce point $(r, s, t)$ sur un point de coordonnées $(\rho,
\sigma, \tau)$ satisfaisant $r = \rho$, $s = \sigma$, $t = \tau$. Dans
ce cas, la transformation entre les différentielles de
coordonnées\footnote{\,
%%%%%%%%%%%%%%%%%%%%%%%-------DEBUT--------%%%%%%%%%%%%%%%%%%%%%%%%%%%
En parallèle, on se reportera avantageusement à la discussion 
qu'en font Engel et Lie au début du \S~94 p.~\pageref{S-94}, 
notamment aux équations~\thetag{ 15}, \thetag{ 16} et~\thetag{ 16'}. 
}: %%%%%%%%%%%%%%%%%%%%%%%%-----FIN-----%%%%%%%%%%%%%%%%%%%%%%%%%%%%%%%:
\def\theequation{1}\begin{equation}
\left\{
\aligned
dr
&
=
A_0d\rho+B_0d\sigma+C_0d\tau
\\
ds
&
=
A_1d\rho+B_1d\sigma+C_1d\tau
\\
dt
&
=
A_2d\rho+B_2d\sigma+C_2d\tau
\endaligned\right.
\end{equation}
incorpore neuf fonctions $A_n$, $B_n$, $C_n$ ($n=1, 2, 3$) qui
dépendent encore de trois paramètres arbitraires\footnote{\,
%%%%%%%%%%%%%%%%%%%%%%%-------DEBUT--------%%%%%%%%%%%%%%%%%%%%%%%%%%%
Dans la théorie de Lie, cette donnée correspond essentiellement à
considérer le sous-groupe d'isotropie d'un point $(x_0, y_0, z_0)$
fixé.
} %%%%%%%%%%%%%%%%%%%%%%%%-----FIN-----%%%%%%%%%%%%%%%%%%%%%%%%%%%%%%%
que Helmholtz note $p'$, $p''$ et $p'''$. Ces équations
montrent comment sont transformés les éléments infinitésimaux basés
au point fixé. 

Dans ces trois équations~\thetag{ 1}, Helmholtz introduit
alors les notations:
\[
\aligned
dr
&
=
\varepsilon\,x
\ \ \ \ \ \ \ \ \ \ \ \ \ \ \ \ \ \ 
d\rho=\varepsilon\,\xi
\\
ds
&
=
\varepsilon\,y
\ \ \ \ \ \ \ \ \ \ \ \ \ \ \ \ \ \ 
d\sigma=\varepsilon\,\upsilon
\\
dt
&
=
\varepsilon\,z
\ \ \ \ \ \ \ \ \ \ \ \ \ \ \ \ \ \ \,
d\tau=\varepsilon\,\zeta,
\endaligned
\] 
où $\varepsilon$ est une quantité infinitésimale que l'on peut
visiblement supprimer en divisant de part et d'autre de chaque
équation. De cette manière, on envisage les transformations dans
l'infiniment petit comme des transformations {\em linéaires} agissant
dans le domaine {\em fini}\footnote{\,
%%%%%%%%%%%%%%%%%%%%%%%-------DEBUT--------%%%%%%%%%%%%%%%%%%%%%%%%%%%
La division par $\varepsilon$ agit comme un zoom infiniment puissant
au point considéré.  La matrice $3 \times 3$ ainsi obtenue de
fonctions de $p'$, $p''$, $p'''$ est localement inversible par
construction.
}: %%%%%%%%%%%%%%%%%%%%%%%%-----FIN-----%%%%%%%%%%%%%%%%%%%%%%%%%%%%%%%
\def\theequation{2}\begin{equation}
\left\{
\aligned
x
&
=
A_0\,\xi+B_0\,\upsilon+C_0\,\zeta
\\
y
&
=
A_1\,\xi+B_1\,\upsilon+C_1\,\zeta
\\
z
&
=
A_2\,\xi+B_2\,\upsilon+C_2\,\zeta.
\endaligned\right.
\end{equation}

D'après Engel et Lie (\S~94 p.~\pageref{S-94} ci-dessous),
lorsque le point fixé est rapporté à l'origine
d'un système de coordonnées $(x, y, z)$, le raisonnement de Helmholtz
revient à effectuer un développement en série entière par rapport aux
trois variables $x$, $y$ et $z$, mais en supprimant tous les termes
d'ordre $\geqslant 2$, ce qui donne un groupe linéaire
homogène qu'ils
écrivent quant à eux\footnote{\,
%%%%%%%%%%%%%%%%%%%%%%%-------DEBUT--------%%%%%%%%%%%%%%%%%%%%%%%%%%%
---\,\,sous-groupe à trois
paramètres de ${\rm GL}_3 ( \R)$\,\,---
}: %%%%%%%%%%%%%%%%%%%%%%%%-----FIN-----%%%%%%%%%%%%%%%%%%%%%%%%%%%%%%%
\[
\left\{
\aligned
x'
&
=
\lambda_1\,x+\lambda_2\,y+\lambda_3\,z
\\
y'
&
=
\mu_1\,x+\mu_2\,y+\mu_3\,z
\\
z'
&
=
\nu_1\,x+\nu_2\,y+\nu_3\,z. 
\endaligned\right.
\]
Ce sous-groupe du groupe linéaire homogène complet ${\rm GL}_3 ( \R)$
indique alors comment sont transformées les droites (infiniment
petites) passant par l'origine, lorsque le corps rigide pivote de
toutes les manières possibles autour de son point fixe.  

Ce groupe
linéaire agissant au niveau infinitésimal constitue une découverte
importante.  En effet, la considération de telles transformations
linéarisées (et projectivisées) en un point fixe constitue l'une des
plus idées les plus efficientes que Lie a mises au point afin de
classifier groupes et sous-groupes de transformations (\voir~le
Chapitre~5), et l'on peut s'imaginer que les considérations
embryonnaires et <<\,balbutiantes\,>> de Helmholtz ont pu constituer
pour Lie une source constante d'inspiration et de défi, un point
d'orgue à atteindre, ce qui expliquerait peut-être la raison pour
laquelle le chapitre sur le problème de Riemann-Helmholtz a été placé
à l'extrême fin de la {\em Theorie der
Transformationsgruppen}\footnote{\,
%%%%%%%%%%%%%%%%%%%%%%%-------DEBUT--------%%%%%%%%%%%%%%%%%%%%%%%%%%%
L'ultime Division~VI du Volume~III~\cite{ enlie1893} de la {\em
Theorie der Transformationsgruppen} comprend cinq chapitres consacrés
à des considérations générales, historiques et méthodologiques qui
reviennent en arrière sur l'ensemble de l'ouvrage.
}. %%%%%%%%%%%%%%%%%%%%%%%%-----FIN-----%%%%%%%%%%%%%%%%%%%%%%%%%%%%%%%

\smallskip

Ainsi, comme il l'a annoncé au tout début de son exposition
mathématique, Helmholtz applique <<\,seulement\,>> ses hypothèses à
des points infiniment proches les uns des autres, comme si cette
exigence semblait en demander {\em moins} que lorsqu'on les applique à
un domaine d'extension {\em fini}, puisque la portée des axiomes
serait de la sorte {\em restreinte} à un domaine d'extension encore
plus petit. Toutefois, la circulation imparfaite entre l'infinitésimal
et le fini va réserver des surprises à l'exploration 
rigoureuse en termes de contre-exemples. 

\HEAD{Chapitre~2.\,\,\,\,La mobilité helmholtzienne de la rigidité}{
2.6.\,\,\,Critique par Lie de l'erreur principale de Helhmoltz}

\medskip\noindent{\bf 2.6.~Critique par Lie de l'erreur
principale de Helhmoltz.}
\label{critique-Lie-Helmholtz}
Engel et Lie interprètent les exigences helmholtziennes en termes de
théorie des groupes continus de transformations\footnote{\,
%%%%%%%%%%%%%%%%%%%%%%%-------DEBUT--------%%%%%%%%%%%%%%%%%%%%%%%%%%%
Le lecteur est renvoyé au Chapitre~4 et au \S~85 p.~\pageref{S-85}
pour une présentation des éléments fondamentaux de la
théorie; sans cela, il peut aussi se reporter directement 
à la fin de ce paragraphe pour prendre connaissance 
de la teneur des contre-exemples de Lie. 
}. %%%%%%%%%%%%%%%%%%%%%%%%-----FIN-----%%%%%%%%%%%%%%%%%%%%%%%%%%%%%%%
En notant $p = \frac{ \partial }{ \partial x}$, $q := \frac{ \partial
}{ \partial y}$ et $r := \frac{ \partial }{ \partial z}$ pour abréger,
soient:
\[
X_k
=
\xi_k(x,y,z)\,p
+
\eta_k(x,y,z)\,q
+
\zeta_k(x,y,z)\,r 
\ \ \ \ \ \ \ 
{\scriptstyle{(k\,=\,1\,\cdots\,6)}}
\]
six transformations infinitésimales qui engendrent un groupe continu à
six paramètres de transformations locales de l'espace réel à trois
dimensions muni des coordonnées $(x, y, z)$, ce groupe étant supposé
satisfaire les axiomes II, III et IV de Helmholtz.  Le deuxième
axiome: parfaite mobilité des corps rigides, demande notamment que le
groupe soit transitif.  De manière équivalente, en tout point donné à
l'avance, les (champs de) vecteurs $X_1, \dots, X_6$
engendrent l'espace tangent en ce point.  Sans perte de
généralité, on peut supposer que le système des coordonnées est centré
au point considéré. 
Quitte à effectuer des combinaisons linéaires (pivot de
Gauss), on peut supposer aussi que les
trois premières transformations infinitésimales
s'écrivent:
\[
X_1
=
p+\cdots,
\ \ \ \ \ \ \ \ \ \ \ \
X_2
=
q+\cdots,
\ \ \ \ \ \ \ \ \ \ \ \
X_3
=
r+\cdots,
\]
où les termes <<\,$+\cdots$\,>> d'ordre supérieurs sont des restes de
la forme\footnote{\,
%%%%%%%%%%%%%%%%%%%%%%%-------DEBUT--------%%%%%%%%%%%%%%%%%%%%%%%%%%%
Ici, la notation ${\rm O} ( 1)$ désigne un reste analytique s'annulant
pour $x = y = z = 0$; de même, la notation ${\rm O} ( 2)$ désigne un
reste analytique s'annulant, ainsi que toutes ses dérivées partielles
d'ordre $1$, lorsque $x = y = z = 0$.  }
%%%%%%%%%%%%%%%%%%%%%%%%-----FIN-----%%%%%%%%%%%%%%%%%%%%%%%%%%%%%%%
${\rm O}(1)\,p+ {\rm O} ( 1) \, q + {\rm O} ( 1) \, r$, et
aussi en même temps que les
{\em trois} transformations infinitésimales linéairement indépendantes
restantes s'annulent à l'origine, de telle sorte que
l'on peut écrire:
\[
\aligned
X_k
:=
&
\quad 
(\alpha_{k1}\,x+\alpha_{k2}\,y+\alpha_{k3}\,z+\cdots)\,p+
\\
&
+
(\beta_{k1}\,x+\beta_{k2}\,y+\beta_{k3}\,z+\cdots)\,q+
\\
&
+
(\gamma_{k_1}\,x+\gamma_{k2}\,y+\gamma_{k3}\,z+\cdots)\,r
\\
&
\quad\quad\quad\quad\quad
{\scriptstyle{(k\,=\,4,\,5,\,6)}},
\endaligned
\]
où les termes <<\,$+\cdots$\,>> d'ordre supérieur sont de la forme
${\rm O}(2)\,p+ {\rm O} ( 2) \, q + {\rm O} ( 2) \, r$, et où les
$\alpha_{ kl}$, $\beta_{ kl}$, $\gamma_{kl}$ sont des constantes
réelles. Engel et Lie introduisent alors le {\sl groupe réduit}
associé au groupe $X_1, \dots, X_6$, qui est formé des transformations
infinitésimales $p$, $q$ et $r$ issues de $X_1$, $X_2$ et $X_3$ en
supprimant purement et simplement tous les termes d'ordre $\geqslant
1$, avec les trois transformations d'ordre exactement égal à $1$:
\[
\aligned
L_k
=
(\alpha_{k1}\,x+
&
\alpha_{k2}\,y+\alpha_{k3}\,z)\,p
+
(\beta_{k1}\,x+\beta_{k2}\,y+\beta_{k3}\,z)\,q
+
\\
&
+
(\gamma_{k1}\,x+\gamma_{k2}\,y+\gamma_{k3}\,z)\,r
\\
&
\quad\quad\quad
{\scriptstyle{(k\,=\,1,\,2,\,3)}},
\endaligned
\]
qui sont obtenues à partir de $X_4$, $X_5$ et $X_6$ en supprimant tous
les termes d'ordre $\geqslant 2$.  Alors l'hypothèse que les six
transformations infinitésimales initiales $X_1, \dots, X_6$ forment
une algèbre de Lie\footnote{\,
%%%%%%%%%%%%%%%%%%%%%%%-------DEBUT--------%%%%%%%%%%%%%%%%%%%%%%%%%%%
{\em Voir}~le \S~4.9 et la note
p.~\pageref{justification-groupe-reduit}.  }
%%%%%%%%%%%%%%%%%%%%%%%%-----FIN-----%%%%%%%%%%%%%%%%%%%%%%%%%%%%%%%
implique que les six transformations réduites $p$, $q$, $r$, $L_1$,
$L_2$ et $L_3$ forment elles aussi une algèbre de Lie.  Par
conséquent, d'après le troisième théorème fondamental de Lie (\S~4.9),
les six transformations infinitésimales réduites $p$, $q$, $r$, $L_1$,
$L_2$ et $L_3$ engendrent à nouveau un groupe de transformations, qui
constitue en fait un certain sous-groupe à six paramètres du groupe
affine complet de l'espace à trois dimensions.

Par une série d'exemples que Lie avait déjà fait paraître en 1892 dans
les {\em Comptes Rendus de l'Académie}, Engel et Lie montrent combien
il est périlleux d'extrapoler les axiomes supposés valides dans des
régions locales d'extension finie au sujet du comportement de points
qui sont infiniment proches les uns des autres.  Ainsi l'existence
d'un, et d'un seul invariant généralisé (distance abstraite) entre
paires de points ne se transfère-t-elle pas fidèlement d'un univers à
l'autre.  En effet, dans le \S~94 p.~\pageref{S-94} ci-dessous, Engel
et Lie décrivent deux exemples élémentaires:

\smallskip$\bullet$\,\,
un groupe $X_1, \dots, X_6$ pour lequel deux points ont {\em un et un
seul} invariant, alors que pour le groupe réduit $p$, $q$, $r$, $L_1$,
$L_2$, $L_3$, deux points ont {\em deux} invariants fonctionnellement
indépendants\footnote{\,
%%%%%%%%%%%%%%%%%%%%%%%-------DEBUT--------%%%%%%%%%%%%%%%%%%%%%%%%%%%
La raison formelle de ces décalages de structure est simple d'un point
de vue analytique: le passage du groupe~\thetag{ 21} de dimension $6$
au groupe réduit~\thetag{ 21'} du \S~94 p.~\pageref{S-94} supprime
tellement de termes, que le groupe réduit devient de dimension $4$!
}; %%%%%%%%%%%%%%%%%%%%%%%%-----FIN-----%%%%%%%%%%%%%%%%%%%%%%%%%%%%%%%

\smallskip$\bullet$\,\,
un autre groupe $X_1, \dots, X_6$ pour lequel deux points n'ont {\em
aucun} invariant, alors que pour le groupe réduit $p$, $q$, $r$,
$L_1$, $L_2$, $L_3$, deux points ont {\em un et un seul} invariant. 

\smallskip 

\label{erreur-helmholtz}
Voilà donc l'{\em erreur principale de Helmholtz}:
supprimer <<\,à la
physicienne\,>> tous les termes d'ordre supérieur à $1$, ce qui
détruit complètement l'harmonie euclidienne des corps rigides qui se
donnait pour évidente à l'intuition. Et l'axiome de monodromie recèle
un phénomène encore plus troublant. Dans le Théorème~37
p.~\pageref{Theorem-37-p-433} ci-dessous, sans utiliser les axiomes
III et IV de Helmholtz, Engel et Lie trouvent abstraitement {\em onze
groupes réels} pour lesquels deux points ont un et un seul invariant,
tandis qu'un nombre de point supérieur à deux n'a pas d'invariant
essentiel. Au cours de cette recherche, ils découvrent un groupe
particulier, le groupe~\thetag{ 24} p.~\pageref{457}, qui satisfait
l'axiome de monodromie, alors que son groupe réduit ne le satisfait
pas: même l'axiome le plus central dans les raisonnements
mathématiques de Helmholtz est remis en cause dans le passage à
l'infiniment petit!

\CITATION{Gr\^ace 
aux exemples pr\'ec\'edents, il a \'et\'e suffisamment
d\'emontr\'e que la supposition
que Monsieur de Helmholtz a introduite
tacitement [l'infinitésimalisation]
et qui a \'et\'e d\'ecrite plus pr\'ecis\'ement aux
pages~\pageref{455} sq. [dans la traduction
ci-dessous] est erron\'ee. Et maintenant,
comme ses consid\'erations ult\'erieures prennent enti\`erement cette
supposition comme point de d\'epart et n'ont force de preuve que sur
la base de cette supposition, nous parvenons donc au r\'esultat que
Monsieur de Helmholtz n'a pas d\'emontr\'e l'assertion qu'il \'enonce
\`a la fin de son travail, \`a savoir: il n'a pas d\'emontr\'e que
ses axiomes suffisent \`a caract\'eriser les mouvements euclidiens et
non-euclidiens.
\REFERENCE{p.~\pageref{non-complet-Helmholtz}~ci-dessous.}}

Cette observation conduit alors à contourner l'erreur principale de
Helmholtz en supposant directement que les axiomes II, III et IV sont
valides dans l'infiniment petit, c'est-à-dire plus précisément, que le
groupe réduit les satisfait. On se ramène ainsi à un problème plus
accessible: trouver, dans le groupe affine complet de l'espace à trois
dimensions, tous les sous-groupes transitifs à six paramètres pour
lesquels deux points ont un et un seul invariant et qui satisfont de
plus l'axiome de monodromie. En utilisant toute la force de la théorie
des groupes de transformations, Engel et Lie vérifieront que les
conclusions de Helmholtz peuvent être rendus rigoureuses et véridiques
par une voie différente (\voir~le
\S~2.8 ci-dessous). Et sans insister ici sur les erreurs
mathématiques de Helmholtz\footnote{\,
%%%%%%%%%%%%%%%%%%%%%%%-------DEBUT--------%%%%%%%%%%%%%%%%%%%%%%%%%%%
Le Chapitre~21 p.~\pageref{Chapitre-21} ci-dessous y consacre déjà une
analyse suffisamment minutieuse.  },
%%%%%%%%%%%%%%%%%%%%%%%%-----FIN-----%%%%%%%%%%%%%%%%%%%%%%%%%%%%%%%
eu égard à son but principal qui était de retrouver la métrique
pythagoricienne tridimensionnelle $ds^2 = dx_1^2 + dx_2^2 + dx_3^2$
telle qu'elle était représentée dans l'infiniment petit par Gauss et par
Riemann, l'approche qui consiste à infinitésimaliser d'emblée les
quatre axiomes était tout à fait
naturelle. 

\HEAD{Chapitre~2.\,\,\,\,La mobilité helmholtzienne de la rigidité}{
2.7.\,\,\,Calculs helmholtziens}

\medskip\noindent{\bf 2.7.~Calculs helmholtziens.}
\label{calculs-helmholtziens}
Pour l'instant, reprenons maintenant les raisonnements mathématiques
originaux de Helmholtz (\cite{ helm1977}) que
Lie ne cherchera pas à imiter\footnote{\,
%%%%%%%%%%%%%%%%%%%%%%%-------DEBUT--------%%%%%%%%%%%%%%%%%%%%%%%%%%%
Grande liberté de stratégie dans la 
théorie des groupes. 
}. %%%%%%%%%%%%%%%%%%%%%%%%-----FIN-----%%%%%%%%%%%%%%%%%%%%%%%%%%%%%%%
Dans les équations de transformations linéaires homogènes~\thetag{ 2},
les quantités $A_n$, $B_n$, $C_n$ dépendent de trois paramètres
indépendants $p'$, $p''$ et $p'''$. Supposons alors que $p'$, $p''$ et
$p'''$ dépendent {\em linéairement} d'une variable <<\,temporelle\,>>
auxiliaire notée $\eta$, introduisons les trois quotients
différentiels:
\[
\mathfrak{A}_n
:=
\frac{d}{d\eta}
\Big(
A_n
\big(
p'(\eta),p''(\eta),p'''(\eta)
\big)
\Big)
\ \ \ \ \ \ \ \ \ \ \ \ \
{\scriptstyle{(n\,=\,0,\,1,\,2)}},
\]
et de même, introduisons les quotients différentiels analogues
$\mathfrak{ B}_n$, $\mathfrak{ C}_n$ pour $n = 0, 1, 2$.  En
différentiant alors les équations~\thetag{ 2} par rapport à $\eta$, on
obtient des équations:
\[
\left\{
\aligned
\frac{dx}{d\eta}
&
=
\mathfrak{A}_0\,\xi+\mathfrak{B}_0\,\upsilon+\mathfrak{C}_0\,\zeta
\\
\frac{dy}{d\eta}
&
=
\mathfrak{A}_1\,\xi+\mathfrak{B}_1\,\upsilon+\mathfrak{C}_1\,\zeta
\\
\frac{dz}{d\eta}
&
=
\mathfrak{A}_2\,\xi+\mathfrak{B}_2\,\upsilon+\mathfrak{C}_2\,\zeta,
\endaligned\right.
\]
dans lesquelles on peut réexprimer $\xi$, $\upsilon$, $\zeta$ en
fonction de $x$, $y$, $z$ en inversant le système~\thetag{ 2}, ce qui
donne un système homogène d'équations différentielles ordinaires
d'ordre $1$:
\def\theequation{3}\begin{equation}
\left\{
\aligned
\frac{dx}{d\eta}
&
=
a_0x+b_0y+c_0z
\\
\frac{dy}{d\eta}
&
=
a_1x+b_1y+c_1z
\\
\frac{dz}{d\eta}
&
=
a_2x+b_2y+c_2z,
\endaligned\right.
\end{equation}
avec certaines fonctions $a_n$, $b_n$, $c_n$ ($n = 0, 1, 2$) du
paramètre individuel $\eta$.  Les paramètres initiaux $p'$, $p''$,
$p'''$ étant au nombre de trois, on peut en fait obtenir trois tels
systèmes indépendants en choisissant la dépendance (supposée linéaire)
par rapport à $\eta$ de ces paramètres le long de trois directions de
droite qui sont indépendantes
(\cf~ce qui va suivre). 

Dans un passage difficile à déchiffrer dont le contenu fut
ensuite très largement englobé par la théorie de Lie des
transformations infinitésimales associées à un groupe linéaire
homogène, Helmholtz montre que les fonctions $a_n$, $b_n$, $c_n$ sont
en fait {\em constantes}.  Ainsi peut-on appliquer au système~\thetag{
3} ci-dessus les théorèmes bien connus de la théorie des systèmes
d'ordre $1$ à coefficients constants.

Par ailleurs, en supposant l'existence d'un, et d'un seul invariant
pour toute paire de points relativement au groupe d'isotropie
linéarisé (validité de l'axiome~II dans l'infinitésimal), Helmholtz
affirme dans un court passage tout
aussi délicat\footnote{\,
%%%%%%%%%%%%%%%%%%%%%%%-------DEBUT--------%%%%%%%%%%%%%%%%%%%%%%%%%%%
Paul Hertz, qui a bénéficié d'échanges épistolaires avec le professeur
Engel, reconstitue l'argument (\cite{ helm1977}, pp.~67--68).  Pour
tout deuxième point $(x_0, y_0, z_0)$ distinct de l'origine, le
troisième axiome helmholtzien demande qu'un mouvement à un paramètre
soit encore possible. Puisque les quantités numériques $p', p'', p'''$
et $x_0, y_0, z_0$ sont toutes deux au nombre de trois, Helmholtz
semble admettre que tous les systèmes~\thetag{ 3} correspondent
biunivoquement avec les mouvements qui fixent deux points
distincts. 
Alors pour tout tel deuxième point
donné $(x_0, y_0, z_0)$ distinc
de l'origine, on doit
pouvoir choisir une certaine dépendance linéaire
de $(p', p'', p''')$ par rapport
à un paramètre $\eta$ tel que
le système~\thetag{ 3} décrive précisément tous les
mouvements fixant l'origine et le
deuxième point $(x_0, y_0, z_0)$. On en déduit donc
un système d'équations linéaires: 
\[
\aligned
0
&
=
{\textstyle{\frac{dx_0}{d\eta}}}
=
a_0x_0+b_0y_0+c_0z_0
\\
0
&
=
{\textstyle{\frac{dy_0}{d\eta}}}
=
a_1x_0+b_1y_0+c_1z_0
\\
0
&
=
{\textstyle{\frac{dz_0}{d\eta}}}
=
a_2x_0+b_2y_0+c_2z_0
\endaligned
\]
qui implique l'annulation du déterminant~\thetag{ 4}. 
} %%%%%%%%%%%%%%%%%%%%%%%%-----FIN-----%%%%%%%%%%%%%%%%%%%%%%%%%%%%%%%
que le déterminant du
système doit s'annuler:
\def\theequation{4}\begin{equation}
0
=
\left\vert
\begin{array}{ccc}
a_0 & b_0 & c_0 
\\
a_1 & b_1 & c_1
\\
a_2 & b_2 & c_2
\end{array}
\right\vert.
\end{equation}

Dans son approche du problème, Helmholtz semble connaître et maîtriser
la théorie qui permet de résoudre les systèmes linéaires homogènes
d'équations différentielles ordinaires d'ordre $1$, théorie que Lie
appliquera ensuite très fréquemment. Résumons brièvement les
résultats dans un langage contemporain, en admettant le théorème de
réduction des matrices à une forme normale de 
Jordan (\cf~\cite{ brec2006} pour un
historique).

L'annulation du déterminant montre qu'une valeur propre au moins de la
matrice vaut zéro.  Helmholtz démontre
avec clarté que les deux autres
valeurs propres $\lambda_2$ et $\lambda_3$ (éventuellement égales)
sont nécessairement imaginaires pures
conjuguées l'une de l'autre\footnote{\,
%%%%%%%%%%%%%%%%%%%%%%%-------DEBUT--------%%%%%%%%%%%%%%%%%%%%%%%%%%%
Ceci découle directement du fait que la matrice est réelle.
Comme ses valeurs propres sont alors distinctes, 
elle est diagonablisable sur $\C$. 
} %%%%%%%%%%%%%%%%%%%%%%%%-----FIN-----%%%%%%%%%%%%%%%%%%%%%%%%%%%%%%%
et non nulles. 
En effet, la solution au système~\thetag{ 3} de condition initiale le
point $(x_0, y_0, z_0)$ est donnée par l'exponentielle de matrice $e^{
M t}$ appliquée à ce point, vu comme vecteur colonne. Après un
changement de base éventuel, les trois éléments diagonaux de la
matrice $e^{ Mt}$ sont: $e^{ 0 \, t}$, $e^{ \lambda_2 t}$ et $e^{
\lambda_3 t}$. Dès que la partie réelle de $\lambda_2$ ou celle de
$\lambda_3$ est non nulle, il est impossible que le mouvement à un
paramètre représenté par le système~\thetag{ 3} soit périodique, comme
le demande l'axiome de monodromie, puisque $e^{ (a + ib)t}$ diverge
lorsque $t \to \pm \infty$, suivant le signe de $a$, lorsque $a \neq
0$. Enfin, le cas où {\em toutes} les valeurs propres sont nulles ne
se produit pas: ou bien la forme normale de Jordan de $M$ est
identiquement nulle, d'où tous les points restent au repos (cas
exclu par l'axiome
de mobilité); ou bien cette forme normale est nilpotente, et le mouvement
s'éloigne polynomialement vers l'infini\footnote{\,
%%%%%%%%%%%%%%%%%%%%%%%-------DEBUT--------%%%%%%%%%%%%%%%%%%%%%%%%%%%
Dans de nouvelles coordonnées normalisantes $(X, Y, Z)$, le système
s'écrit: $\frac{ d Z}{ dt} = 0$, $\frac{ dY}{ dt} = Z$ et $\frac{ dX}{
dt} = \varepsilon \, Y$, avec $\varepsilon = 0$ ou $= 1$ suivant que
le rang de la matrice (nilpotente) $M$ est égal à $1$ ou à $2$.
L'intégration donne un mouvement: $Z = Z_0$, $Y = Y_0 + Z_0 t$, $X =
\varepsilon \big(
X_0 + Y_0 t + \frac{ 1}{ 2}\, Z_0 \, t^2 \big)$ qui n'est pas périodique.
}. %%%%%%%%%%%%%%%%%%%%%%%%-----FIN-----%%%%%%%%%%%%%%%%%%%%%%%%%%%%%%%

Ainsi, dans de nouvelles coordonnées appropriées, le
système peut être écrit sous une forme
normale simplifiée: 
\[
\left\{
\aligned
\frac{dX}{d\eta}
&
=
0
\\
\frac{dY}{d\eta}
&
=
-\omega\,Z
\\
\frac{dZ}{d\eta}
&
=
\omega\,Y. 
\endaligned\right.
\]
Sans surprise, on retrouve les équations différentielles d'une
rotation euclidienne d'axe $\{ Y = Z = 0\}$.
\`A la très grande généralité initiale du propos
de Helmholtz succède donc la considération de mouvements qui étaient
déjà fort bien compris dans les mathématiques et dans la
mécanique de l'époque.  Lie au
contraire replacera l'ambition de généralité à un niveau largement
supérieur. L'axiome de monodromie qui semblait intuitivement si
évident impose donc une condition extrêmement forte qui permet
d'éliminer presque tous les systèmes~\thetag{ 3}, 
excepté ceux qui fournissent des rotations euclidiennes. 

Ensuite, en modifiant convenablement
la dépendance de $(p', p'', p''')$ par rapport à
$\eta$, Helmoltz prétend que l'on obtient deux autres rotations le
long des deux axes $\{ X = Z = 0 \}$ et $\{ X = Y = 0\}$ qu'il écrit
sous la forme générale suivante, en introduisant deux nouvelles
directions unidimensionnelles $\eta'$ et $\eta''$ dans l'espace des
paramètres $(p', p'', p''')$:
\[
\aligned
\left\{
\aligned
\frac{dX}{d\eta'}
&
=
\alpha_0\,X+0+\gamma_0\,Z
\\
\frac{dY}{d\eta'}
&
=
\alpha_1\,X+0+\gamma_1\,Z
\\
\frac{dZ}{d\eta'}
&
=
\alpha_2\,X+0+\gamma_2\,Z
\endaligned\right.
\endaligned
\ \ \ \ \ \ \ \ \ \ \ \
\text{\rm et}
\ \ \ \ \ \ \ \ \ \ \ \
\aligned
\left\{
\aligned
\frac{dX}{d\eta''}
&
=
\mathfrak{a}_0\,X+\mathfrak{b}_0\,Y+0
\\
\frac{dY}{d\eta''}
&
=
\mathfrak{a}_1\,X+\mathfrak{b}_1\,Y+0
\\
\frac{dZ}{d\eta''}
&
=
\mathfrak{a}_2\,X+\mathfrak{b}_2\,Y+0.
\endaligned\right.
\endaligned
\]
Pour chacun de ces deux systèmes, le déterminant
correspondant~\thetag{ 4} doit s'annuler, et en particulier, la trace
matricielle doit être nulle:
\[
0
=
\alpha_0+\gamma_2
\ \ \ \ \ \ \ \ \ \ \ \
\text{\rm et}
\ \ \ \ \ \ \ \ \ \ \ \
0
=
\mathfrak{a}_0+\mathfrak{b}_1. 
\]
Ensuite, Helmholtz utilise un argument indirect qu'il aurait pu
présenter de manière nettement plus limpide sur le plan technique: par
linéarité supposée de la dépendance de $(p', p'', p''')$ par rapport à
$\eta$, $\eta'$ et $\eta''$, le système obtenu par sommation de deux
quelconques systèmes parmi les trois systèmes considérés doit encore
constituer un système du même type, de telle sorte que le déterminant
correspondant~\thetag{ 4} doit à nouveau s'annuler.  Cet argument
combiné à un argument d'homogénéité des équations obtenues par rapport
aux constantes qui interviennent dans chaque système permet à
Helmholtz de déduire que\footnote{\,
%%%%%%%%%%%%%%%%%%%%%%%-------DEBUT--------%%%%%%%%%%%%%%%%%%%%%%%%%%%
Le raisonnement devrait
converger vers le fait que la matrice
combinaison linéaire générale des trois systèmes en question n'est
autre qu'une matrice 
$3 \times 3$ antisymétrique quelconque, quoiqu'un
tel énoncé parlant soit en fait absent du mémoire de Helmholtz.
}: %%%%%%%%%%%%%%%%%%%%%%%%-----FIN-----%%%%%%%%%%%%%%%%%%%%%%%%%%%%%%%
\[
0
=
\alpha_0=\alpha_1=\gamma_2
\ \ \ \ \ \ \ \ \ \ \ \
\text{\rm et}
\ \ \ \ \ \ \ \ \ \ \ \
0
=
\mathfrak{a}_0=\mathfrak{b}_1=\mathfrak{a}_2, 
\]
et que: 
\def\theequation{5}\begin{equation}
0
=
\gamma_0\mathfrak{a}_1
-
\mathfrak{b}_0\alpha_2.
\end{equation}
Enfin, en formant la somme des trois systèmes après annulation de ces
six constantes et en réappliquant l'argument d'annulation nécessaire
du déterminant~\thetag{ 4}, Helmholtz montre que\footnote{\,
%%%%%%%%%%%%%%%%%%%%%%%-------DEBUT--------%%%%%%%%%%%%%%%%%%%%%%%%%%%
L'étude de la matrice combinaison linéaire générale des trois systèmes
s'arrête là. Fait surprenant: Helmholtz qui semble être conscient de
retrouver les équations différentielles du groupe des rotations qui
fixent un point dans l'espace tridimensionnel ne poursuit néanmoins
pas le raisonnement de normalisation des constantes, et il oublie de
convoquer à nouveau le fait que les deux valeurs propres non nulles de
chaque système doivent être imaginaires conjuguées, ce qui lui aurait
donné les équations finales: $\mathfrak{ b}_0 = - \mathfrak{ a}_1$ et
$\gamma_0 = - \alpha_2$. \`A l'inverse chez Engel et Lie, le souci de
complétion absolue de la pensée technique ne laisse dans l'ombre
aucun calcul en vue des harmonies formelles
conclusives.  
}: %%%%%%%%%%%%%%%%%%%%%%%%-----FIN-----%%%%%%%%%%%%%%%%%%%%%%%%%%%%%%% 
\[
0
=
\gamma_1
\ \ \ \ \ \ \ \ \ \ \ \
\text{\rm et}
\ \ \ \ \ \ \ \ \ \ \ \
0
=
\mathfrak{b}_2.
\]
Pour terminer, en posant:
\[
\alpha_2
=
-\phi,
\ \ \ \ \ \ \ \ \ \ \ \
\gamma_0=\kappa\phi,
\ \ \ \ \ \ \ \ \ \ \ \
\mathfrak{a}_1=\psi,
\]
d'où il découle grâce à~\thetag{ 5} que:
\[
\mathfrak{b}_0
=
-\kappa\psi,
\]
Helmholtz obtient le système complet de toutes les transformations
possibles qui fixent l'origine:
\[
\left\{
\aligned
dX
&
=
0
+
\kappa\phi\,Z\,d\eta'
-
\kappa\psi\,Y\,d\eta''
\\
dY
&
=
-\omega\,Z\,d\eta
+
0
+
\psi\,X\,d\eta''
\\
dZ
&
=
\omega\,Y\,d\eta
-
\phi\,X\,d\eta'
+
0.
\endaligned\right.
\]
Ici, la quantité $\kappa$ doit être positive pour que les valeurs
propres non nulles soit imaginaires pures conjuguées l'une de l'autre.
Helmholtz observe alors qu'il découle de ce dernier système que la
quantité suivante s'annule:
\[
{\textstyle{\frac{1}{\kappa}}}\,XdX
+
YdY
+
ZdZ
=
0,
\]
c'est-à-dire après intégration: 
\[
X^2
+
\kappa\,Y^2
+
\kappa\,Z^2
=
{\rm const.}
\]
Le facteur positif $\kappa$ peut être supprimé en remplaçant $Y$, $Z$
par $\sqrt{ \kappa} \, Y$, $\sqrt{ \kappa} \, Z$, et l'on retrouve les
équations de la famille des sphères euclidiennes centrées en
l'origine.

\HEAD{Chapitre~2.\,\,\,\,La mobilité helmholtzienne de la rigidité}{
2.8.\,\,\,Insuffisances et reprises}

\medskip\noindent{\bf 2.8.~Insuffisances et reprises.}
\label{insuffisances-et-reprises}
Helmholtz semble être satisfait par cette conclusion, mais deux
erreurs éventuelles peuvent encore être commises.

Premièrement, par construction, les quantités $x$, $y$, $z$ puis $X$,
$Y$, $Z$ qui s'en déduisent par combinaison linéaire demeurent {\em
infinitésimales}. On ne retrouve donc la métrique pythagoricienne que
dans l'infiniment petit, et il est à nouveau absolument hors de
question d'en déduire que cette métrique vaut dans le domaine local
fini, voire dans le domaine global.  D'ailleurs, toute métrique
riemanienne quelconque est infinitésimalement équivalente à une
métrique euclidienne: la conclusion de Helmoltz ne prouve donc que le
caractère infinitésimalement 
riemannien de la métrique, ce dont Helmholtz semble être
néanmoins conscient dans le dernier paragraphe de son mémoire.

Mais deuxièmement, Helmholtz prétend qu'il a en fait exploré la
conséquence de la mobilité des corps rigides dans
l'extension locale finie en tenant compte de la courbure riemannienne,
et il annonce alors {\em sans aucune démonstration} le résultat final
de ses investigations. D'après lui, si les axiomes~I à~IV sont
satisfaits en toute généralité dans le domaine
fini, alors la seule géométrie
correspondante serait celle qui 
prévaut sur une sphère de rayon $R$
dans un espace à quatre dimensions muni des coordonnées réelles $(X,
Y, Z, S)$ d'équation:
\[
X^2+Y^2+Z^2+(S+R)^2
=
R^2,
\]
le cas $R = \infty$ correspondant à la géométrie euclidienne.  Et
puisque dans une telle famille, tous les cas où $R <
\infty$ fournissent une géométrie sur l'espace 
d'une sphère tridimensionnelle qui est {\em compact}, Helmholtz 
a cru pouvoir caractériser pleinement la 
géométrie euclidienne en ajoutant 
deux axiomes\footnote{\,
%%%%%%%%%%%%%%%%%%%%%%%-------DEBUT--------%%%%%%%%%%%%%%%%%%%%%%%%%%%
Implicitement, le premier de ces deux derniers axiomes avait déjà été
admis.  
}: %%%%%%%%%%%%%%%%%%%%%%%%-----FIN-----%%%%%%%%%%%%%%%%%%%%%%%%%%%%%%%

\smallskip$\bullet$\,\,
{\sf V:} L'espace possède trois dimensions. 

\smallskip$\bullet$\,\,
{\sf VI:} L'espace est infini en extension. 

\smallskip\noindent
Malheureusement, dans une lettre datant du 24 avril 1869, Eugenio
Beltrami
%\footnote{\,
%%%%%%%%%%%%%%%%%%%%%%%-------DEBUT--------%%%%%%%%%%%%%%%%%%%%%%%%%%%
%Beltrami, Minding, Lobatchevski\u{\i}, Crelle.
%\Fill
%} %%%%%%%%%%%%%%%%%%%%%%%%-----FIN-----%%%%%%%%%%%%%%%%%%%%%%%%%%%%%%%
informait Helmholtz de la consistance éclatante de la géométrie de
Lobatchevski\u{\i} grâce à sa réalisation sur une {\em pseudosphère
infinie tridimensionnelle}. Helmholtz était ainsi contredit. Pour
conclure, on retiendra l'énoncé suivant, qui laisse encore
essentiellement ouverte la question soulevée par Riemann.

\smallskip\noindent
{\bf Proposition de Helmholtz.}
\label{proposition-helmholtz}
{\em En supposant que la libre mobilité des corps rigides et l'axiome
de monodromie valent tous deux au niveau infinitésimal, le sous-groupe
d'isotropie linéarisé de tout point quelconque est isomorphe à ${\rm
SO}_3 ( \R)$. }

\HEAD{Chapitre~2.\,\,\,\,La mobilité helmholtzienne de la rigidité}{
2.9.\,\,\,L'approche infinitésimale systématique de Engel et de Lie}

\medskip\noindent{\bf 2.9.~L'approche infinitésimale systématique 
de Engel et de Lie.}
\label{approche-infinitesimale-systematique}
Engel et Lie donneront au moins trois solutions distinctes au problème
de Riemann-Helmholtz.  La première (\S~95 p.~\pageref{S-95}
ci-dessous) rebondit sur la pétition de principe helmholtzienne en
proposant tout simplement de formuler directement les axiomes au
niveau infinitésimal.  Mais une différence majeure s'introduit: on ne
suppose alors {\em nullement} l'existence d'une fonction de deux
points quelconques dont la valeur demeure invariante dans les
mouvements (rigidité).

Tout d'abord\footnote{\,
%%%%%%%%%%%%%%%%%%%%%%%-------DEBUT--------%%%%%%%%%%%%%%%%%%%%%%%%%%%
La lecture de ce dernier paragraphe nécessite une connaissance
préalable des fondements de la théorie de Lie qui seront exposées en
détail dans les Chapitres~4 et~5 ci-dessous.
}, %%%%%%%%%%%%%%%%%%%%%%%%-----FIN-----%%%%%%%%%%%%%%%%%%%%%%%%%%%%%%%
Engel et Lie partent de l'hypothèse que l'espace est une variété
numérique (locale) à trois dimensions et les mouvements de cet espace
forment un groupe continu de transformations (locales) qui est
transitif et à six paramètres. Soient $X_1, \dots, X_6$ six
transformations infinitésimales linéairement indépendantes et fermées
par crochet de Lie dans l'espace des $(x, y, z)$ qui engendrent un tel
groupe.  Comme le groupe est transitif, si l'on fixe un point réel
$(x_0, y_0, z_0)$ en position générale, il existe exactement $3$
transformations infinitésimales linéairement indépendantes dont les
combinaisons engendrent l'espace des $(x, y, z)$ en ce point, et il
existe aussi $3 = 6 - 3$ autres transformations infinitésimales
indépendantes s'annulant en ce point qui sont de la forme:
\[
\aligned
Y_k
&
=
\big(
\alpha_{k1}(x-x_0)+\alpha_{k2}(y-y_0)+\alpha_{k3}(z-z_0)+\cdots\big)\,p
+
\\
&
\ \ \ \ \ 
+
\big(
\beta_{k1}(x-x_0)+\beta_{k2}(y-y_0)+\beta_{k3}(z-z_0)+\cdots\big)\,q
+
\\
&
\ \ \ \ \ 
+
\big(
\gamma_{k1}(x-x_0)+\gamma_{k2}(y-y_0)+\gamma_{k3}(z-z_0)+\cdots\big)\,r
\\
&
\ \ \ \ \ \ \ \ \ \ \ \ \
\ \ \ \ \ \ \ \ \ \ \ \ \
\ \ \ \ \ \ \ \ \ \ \ \ \
{\scriptstyle{(k\,=\,1\,2\,3)}},
\endaligned
\]
où par convention, les termes <<\,$+\cdots$\,>> s'annulent à l'ordre
au moins deux en $(x_0, y_0, z_0)$.  En supprimant purement et
simplement tous ces restes, on obtient trois générateurs (qui peuvent
éventuellement devenir linéairement {\em dépendants}):
\[
\aligned
\overline{Y}_k
&
=
\big(
\alpha_{k1}(x-x_0)+\alpha_{k2}(y-y_0)+\alpha_{k3}(z-z_0)\big)\,p
+
\\
&
\ \ \ \ \ 
+
\big(
\beta_{k1}(x-x_0)+\beta_{k2}(y-y_0)+\beta_{k3}(z-z_0)\big)\,q
+
\\
&
\ \ \ \ \ 
+
\big(
\gamma_{k1}(x-x_0)+\gamma_{k2}(y-y_0)+\gamma_{k3}(z-z_0)\big)\,r
\\
&
\ \ \ \ \ \ \ \ \ \ \ \ \
\ \ \ \ \ \ \ \ \ \ \ \ \
\ \ \ \ \ \ \ \ \ \ \ \ \
{\scriptstyle{(k\,=\,1\,2\,3)}},
\endaligned
\]
d'un sous-groupe\footnote{\,
%%%%%%%%%%%%%%%%%%%%%%%-------DEBUT--------%%%%%%%%%%%%%%%%%%%%%%%%%%%
Les relations de crochets de Lie $\big[ Y_j , \, Y_k \big] =
\sum_{ s=1}^3 \, c_{ jks}\, Y_s$ 
où les $c_{ jks}$ sont des constantes réelles, sont en effet
directement héritées par les transformations réduites: $\big[
\overline{ Y}_j , \, \overline{ Y}_k \big] =
\sum_{ s=1}^3 \, c_{ jks}\, \overline{ Y}_s$ parce que
dans le groupe, il n'existe aucune transformation infinitésimale
s'annulant à l'ordre deux en $(x_0, y_0, z_0)$.
} %%%%%%%%%%%%%%%%%%%%%%%%-----FIN-----%%%%%%%%%%%%%%%%%%%%%%%%%%%%%%%
linéaire de ${\rm GL}_3 ( \R)$.  Ce sous-groupe (dit {\sl réduit}) 
détermine de quelle manière les éléments linéaires infinitésimaux
$(dx, dy, dz)$ sont transformés par les transformations
du groupe qui fixent le point $(x_0, y_0, z_0)$ considéré. 
 
Le troisième axiome p.~\pageref{III-infinitesimal} demande que ce
sous-groupe d'isotropie linéarisée demeure à trois paramètres.
Toutefois, aucune existence d'invariant n'est postulée.  Enfin, le
quatrième axiome p.~\pageref{IV-periodique} demande que {\em tout}
sous-groupe à {\em un paramètre} de ce groupe linéaire homogène 
à trois dimensions soit
constitué de mouvements qui 
agissent {\em périodiquement}
sur les éléments linéaires passant par l'origine
\footnote{\,
%%%%%%%%%%%%%%%%%%%%%%%-------DEBUT--------%%%%%%%%%%%%%%%%%%%%%%%%%%%
En première apparence, cette condition semble donc légèrement plus
exigente qu'un axiome de monodromie qui serait formulé au niveau
infinitésimal, en tant qu'un tel axiome ne devrait porter que sur les
sous-groupes à un paramètres qui fixent un deuxième point.
Par ailleurs, cette condition est
légèrement plus générale, en tant que l'on ne 
demande pas que l'élément linéaire revienne coïncider exactement
avec lui-même après un temps fini: on autorise
les dilatations éventuelles. 
}. %%%%%%%%%%%%%%%%%%%%%%%%-----FIN-----%%%%%%%%%%%%%%%%%%%%%%%%%%%%%%%
Plus précisément, pour tout sous-groupe à un paramètre
dudit groupe d'isotropie linéarisée, tout
élément linéaire subit une transformation 
qui le ramène à se retrouver dans
la même direction après un temps fini. 

Dans ces conditions, si l'on demandait de plus que
chaque élément linéaire revienne coïncider avec lui-même
après un temps fini, 
le raisonnement de Helmholtz s'appliquerait encore à
quelques modifications mineures près. En effet, soit $Y := \lambda_1
\overline{ Y}_1 + \lambda_2 \overline{ Y}_2 + \lambda_3 \overline{ Y}_3$ 
une combinaison linéaire quelconque des $\overline{ Y}_k$ à
coefficients réels $\lambda_1,
\lambda_2, \lambda_3$ non tous nuls.  Les valeurs propres de la
matrice $3 \times 3$ des coefficients (réels) de $\overline{ Y}$ par
rapport à $(x - x_0)$, $(y - y_0)$, $(z - z_0)$ sont ou bien réelles,
ou bien complexes conjuguées par paires, et jamais toutes nulles
puisque le groupe linéaire réduit $\overline{ Y}_1,
\overline{ Y}_2, \overline{ Y}_3$ est de dimension trois.
Lorsqu'elles ne sont pas toutes réelles, deux seulement peuvent être
complexes conjuguées\footnote{\,
%%%%%%%%%%%%%%%%%%%%%%%-------DEBUT--------%%%%%%%%%%%%%%%%%%%%%%%%%%%
---\,\,puisque $3$ est le nombre total de valeurs propres comptées
avec multiplicité\,\,--- },
%%%%%%%%%%%%%%%%%%%%%%%%-----FIN-----%%%%%%%%%%%%%%%%%%%%%%%%%%%%%%%
la dernière étant réelle. Par un raisonnement basée sur la forme
normale de Jordan de cette matrice (comme chez Helmholtz), on se
convainc alors aisément que la périodicité stricte du 
mouvement n'est possible que si
deux valeurs propres sont imaginaires {\em
pures} non nulles conjuguées l'une de l'autre, 
tandis que la troisième est
nécessairement {\em nulle}.  En appliquant cet énoncé à $\overline{ Y}
= \overline{ Y}_k$ pour $k = 1, 2, 3$, on trouve alors trois rotations
d'axes linéairement indépendants qui engendrent le groupe ${\rm SO}_3
( \R)$: c'est le groupe d'isotropie linéarisée en tout point, première
étape de la démonstration (la seconde étape est en fait
absente chez Helmholtz).

Toutefois, ce n'est pas du tout de cette manière-là que Engel et Lie
approchent le problème. Dans cette Division~V terminale de leur
ouvrage, ils peuvent se <<\,payer le luxe\,>> de convoquer pleinement
les théorèmes de classification qu'ils ont mis au point auparavant
afin de faire
voir que la caractérisation des groupes d'isotropie linéarisée tombe
comme un fruit mûr par un simple examen des listes
et des théorèmes qui ont été établis précédemment. 

En effet (\cf~le \S~95 p.~\pageref{S-95} ci-dessous), l'isotropie
linéarisée d'un point en position générale que l'on place au centre
(\ie à l'origine) d'un nouveau système de coordonnées $(x_1', x_2',
x_n')$ sera représentée par trois générateurs infinitésimaux
indépendants:
\[
\sum_{\mu\,\nu}^{1,2,3}\,
\alpha_{k\mu\nu}\,x_\mu'\,p_\nu'
\quad\quad\quad\quad
{\scriptstyle{(k\,=\,1,\,2,\,3)}},
\]
où l'on a noté $p_\mu ' := \frac{ \partial }{ \partial x_\mu'}$, qui
forment une algèbre de Lie que nous noterons $g'$.  Comme l'action est
linéaire, elle se transmet à l'espace projectif $P ( \R^3) = \P_2 (
\R)$ et l'on obtient une nouvelle algèbre de Lie $\mathfrak{ g}$ dont
on peut représenter les générateurs dans un système de coordonnées
inhomogènes $(\mathfrak{ x}, \mathfrak{ y})$. Seule la dilatation
$x_1' p_1' + x_2 ' p_2 ' + x_3 ' p_3'$ disparaît lorsqu'on
projectivise\footnote{\,
%%%%%%%%%%%%%%%%%%%%%%%-------DEBUT--------%%%%%%%%%%%%%%%%%%%%%%%%%%%
Elle seule donne la transformation projective infinitésimale
identiquement nulle.  },
%%%%%%%%%%%%%%%%%%%%%%%%-----FIN-----%%%%%%%%%%%%%%%%%%%%%%%%%%%%%%%
et donc la dimension de $\mathfrak{ g}$ vaut deux s'il existe une
combinaison linéaire des générateurs de $g'$ qui est égale à cette
dilatation (ce cas sera exclu dans un instant), et elle vaut trois
sinon.

Deux pièces maîtresses entrent alors en scène.  Nous avons signalé que
l'axiome le plus important par sa puissance de restriction demandait
que chaque mouvement à un paramètre soient périodique sur
les éléments linéaires, à 
un facteur dilatant près.  Quand
on projectivise, le facteur dilatant disparaît, 
chaque point subit un mouvement rigoureusement 
périodique\,\,---\,\,à moins qu'il ne reste fixe\,\,---,
et
cela impose notamment que tous les points décrivent
une courbe fermée. Mais une telle circonstance ne se produit que très
rarement.

Tout d'abord, en dimension $2$ et bien avant que n'apparaisse le
problème de Riemann-Helmholtz, Engel et Lie avaient déjà classifié
{\em toutes} les transformations infinitésimales projectives réelles à
changement de coordonnées projective près.
%%%\Fill 
Il y en a sept (\cf~\thetag{ 30} p.~\pageref{ch21-eq30} 
ci-dessous), et c'est la première pièce maîtresse:
\[
\left\{
\aligned
&
\quad\quad\quad
\mathfrak{p}+\mathfrak{y}\,\mathfrak{q};
\quad
\mathfrak{p}+\mathfrak{x}\,\mathfrak{q};
\quad
\mathfrak{y}\,\mathfrak{q};
\quad
\mathfrak{q};
\\
&
\mathfrak{x}\,\mathfrak{p}+\mathfrak{c}\,\mathfrak{y}\,\mathfrak{q}
\quad
{\scriptstyle{(\mathfrak{c}\,\neq\,0,\,1)}};
\quad
\mathfrak{y}\,\mathfrak{p}-\mathfrak{x}\,\mathfrak{q}
+
\mathfrak{c}\,(\mathfrak{x}\,\mathfrak{p}+\mathfrak{y}\mathfrak{q})
\quad
{\scriptstyle{(\mathfrak{c}\,\neq\,0)}};
\\
&
\quad\quad\quad\quad\quad
\mathfrak{y}\,\mathfrak{p}-\mathfrak{x}\,\mathfrak{q}.
\endaligned\right.
\]
Ici, $(\mathfrak{ x}, \mathfrak{ y})$ sont des coordonnées inhomogènes
sur $P_2 ( \R)$ et $\mathfrak{ p} = \frac{ \partial }{ \partial
\mathfrak{ x}}$, $\mathfrak{ q} = 
\frac{ \partial }{ \partial \mathfrak{ y}}$. 
On vérifie que les cinq premières transformations sont exclues, parce
que chacune d'entre elle laisse globalement invariante au moins une
droite projective sur laquelle la transformation infinitésimale se
ramène, si l'on note $x$ une coordonnée indépendante sur une telle
droite: 1) ou bien à $\frac{ \partial }{ \partial x}$; 2) ou bien à $x
\frac{ \partial}{ \partial x}$; et alors il est clair dans les deux cas
que les courbes intégrales correspondantes: $x (t) = x_0 + t$ et $x(t)
= x_0 e^t$ {\em ne} sont {\em pas} périodiques. La sixième
transformation:
\[
{\sf rotation}
+
\mathfrak{c}\cdot{\sf dilatation}
=
\mathfrak{y}\,\mathfrak{p}-\mathfrak{x}\,\mathfrak{q}
+
\mathfrak{c}\,(\mathfrak{x}\,\mathfrak{p}+\mathfrak{y}\mathfrak{q})
\quad
{\scriptstyle{(\mathfrak{c}\,\neq\,0)}}
\]
est elle aussi exclue, puisque tout point est entraîné autour de
l'origine par la présence du facteur de rotation, tout en étant attiré
(ou répulsé) par le facteur de dilatation: tout mouvement s'effectue
en spirale, sans périodiser.

Enfin, la septième et dernière transformation infinitésimale:
$\mathfrak{ y} \, \mathfrak{p} - \mathfrak{x}\, \mathfrak{q}$ engendre
le groupe des rotations, qui convient parfaitement.

L'utilisation de la deuxième pièce maîtresse est encore plus
magistrale.  Dans les chapitres précédents du Volume~III de la {\em
Theorie der Transformationsgruppen}, Engel et Lie avaient en effet
déjà classifié toutes les sous-algèbres de Lie de l'algèbre de Lie
$\mathfrak{ pgl}_2 ( \R)$ du groupe projectif réel en dimension deux.
%%%\Fill
Il leur suffit alors d'examiner, dans cette liste remarquable et
complète, seulement les algèbres de Lie à deux ou trois paramètres pour
exclure systématiquement toutes celles qui contiennent une
transformation infinitésimale de l'une des six formes précédentes.
Cet examen presque instantané montre qu'il existe {\em un seul} groupe
projectif réel à trois paramètres qui convient, à savoir:
\[
\mathfrak{p}
+
\mathfrak{x}\,(\mathfrak{x}\,\mathfrak{p}
+\mathfrak{y}\,\mathfrak{q}),
\quad\quad
\mathfrak{q}
+
\mathfrak{y}\,(\mathfrak{x}\,\mathfrak{p}
+\mathfrak{y}\,\mathfrak{q}),
\quad\quad
\mathfrak{y}\,\mathfrak{p}-\mathfrak{x}\,\mathfrak{q},
\]
et ce groupe n'est autre que le projectivisé du groupe des rotations
dans l'espace à trois dimensions, à une seule ambiguïté près: la
présence éventuelle d'un facteur de dilatation qui aurait disparu par
projectivisation. Autrement dit, si on revient aux coordonnées
homogènes initiales $(x_1', x_2', x_3')$, les trois générateurs de
l'isotropie linéarisée doivent forcément être de la forme:
\[
\aligned
&
x_\mu'\,p_\nu'-x_\nu'\,p_\mu'
+
\alpha_{\mu\nu}\,
(x_1'\,p_1'+x_2'\,p_2'+x_3'\,p_3')
\\
&
\quad\quad\quad
{\scriptstyle{(\mu,\,\nu\,=\,1,\,2,\,3;\,\,\,\,
\alpha_{\mu\nu}\,+\,\alpha_{\nu\mu}\,=\,0)}}.
\endaligned
\]
Mais alors un simple examen de la condition que ces trois
transformations infinitésimales doivent être fermées par crochet 
montre (\voir~la note p.~\pageref{verification-crochets-rotations}
ci-dessous) que toutes les constantes $\alpha_{ \mu \nu}$ doivent en
fait être nulles. En définitive, Engel et Lie ont 
essentiellement\footnote{\,
%%%%%%%%%%%%%%%%%%%%%%%-------DEBUT--------%%%%%%%%%%%%%%%%%%%%%%%%%%%
Rappelons toutefois que les hypothèses sont
légèrement différentes. 
} %%%%%%%%%%%%%%%%%%%%%%%%-----FIN-----%%%%%%%%%%%%%%%%%%%%%%%%%%%%%%%
redémontré la proposition de Helmholtz
p.~\pageref{proposition-helmholtz} avec de purs moyens de théorie
des groupes.

Maintenant que l'isotropie linéarisée est caractérisée et connue, il
reste à rechercher tous les groupes transitifs à six paramètres dont
l'isotropie linéarisée est constituée groupe des rotations euclidiennes
dans l'espace:
\[
x_1'p_2'-x_2'p_1',
\ \ \ \ \ \ \ \ \ \ \ \ \ 
x_1'p_3'-x_3'p_1',
\ \ \ \ \ \ \ \ \ \ \ \ \ 
x_2'p_3'-x_3'p_2'.
\]
Cette dernière étape de la démonstration, 
que Helmholtz n'a absolument pas traitée, comme nous
l'avons vu, Lie l'a complétée en détail, en 
appliquant les résultats de son traité. 

\smallskip\noindent{\bf Théorème de Lie.}
(\cf~p.~\pageref{Theorem-40})
{\em
Si un groupe réel continu de l'espace à trois dimensions est transitif
à six paramètres et si le groupe d'isotropie linéarisée de tout point
est formé des rotations euclidiennes, alors ce groupe est équivalent,
{\rm via} une transformation ponctuelle r\'eelle de cet espace, soit
au groupe des déplacements euclidiens de l'espace des $(x, y, z)$:
\[
p,\ \ \ \ \
q,\ \ \ \ \ 
r,\ \ \ \ \ 
xq-yp,\ \ \ \ \
xr-zp,\ \ \ \ \
yr-zq,
\]
soit à l'un des deux groupes à six paramètres de mouvements
non-eu\-cli\-diens:
\[
\aligned
&
x_2p_1-x_1p_2,\ \ \ \ \
x_3p_1-x_1p_3,\ \ \ \ \
x_3p_2-x_2p_3,
\\
&
x_1\big(x_1p_1+x_2p_2+x_3p_3\big)
\pm p_1,\ \ \ \ \ \ \ \ \
x_1\big(x_1p_1+x_2p_2+x_3p_3\big)
\pm p_2,
\\
&
x_1\big(x_1p_1+x_2p_2+x_3p_3\big)
\pm p_3,
\endaligned
\] 
par lequel la surface imaginaire: $x_1^2
+ x_2^2 + x_3^2 \pm 1 = 0$ reste 
invariante\footnote{\,
%%%%%%%%%%%%%%%%%%%%%%%-------DEBUT--------%%%%%%%%%%%%%%%%%%%%%%%%%%%
On vérifie en effet aisément que les trois premières transformations
infinitésimales annulent identiquement l'équation, tandis 
que les trois dernières reproduisent l'équation au 
facteur $2x_k$ près, $k = 1, 2, 3$. 
}. %%%%%%%%%%%%%%%%%%%%%%%%-----FIN-----%%%%%%%%%%%%%%%%%%%%%%%%%%%%%%%
}\medskip

%%%%%%%%%%%%%%%%%%%%%%%%%%%%%%%%%%%%%%%%%%%%%%%%%%%%%%%%%%%%%%%%%%%%%

\newpage

\HEAD{Présentation mathématique générale}{
Trois principes de pensée qui gouvernent la théorie de Lie}

\setcounter{footnote}{0}

\bigskip\bigskip\medskip

\begin{center}

{\Large\bf Partie~II:} 

\medskip
{\Large\bf Introduction mathématique} 

\medskip
{\Large\bf à la théorie de Lie}
\label{Partie-II}
\thispagestyle{empty}

\bigskip\medskip

\begin{center}
\begin{minipage}[t]{10.5cm}
\baselineskip =0.35cm
{\scriptsize

\centerline{\footnotesize\bf Chapitres}

\medskip

\smallskip

{\bf Prologue: Trois principes qui gouvernent la théorie de Lie 
\dotfill \pageref{Prologue}.}

{\bf Chapitre~3.~Théorèmes fondamentaux sur les groupes de transformations
\dotfill \pageref{Chapitre-3}.}

%%%{\bf 5.~Théorèmes de classification: groupes et sous-groupes
%%%\dotfill \pageref{Chapitre-5}.}

}\end{minipage}
\end{center}

\bigskip

{\large\bf Prologue: Trois principes de pensée}
\label{Prologue}

\smallskip
{\large\bf qui gouvernent la théorie de Lie}

\end{center}

\bigskip

\noindent{\bf Transformations ponctuelles param\'etr\'ees.}
Soient $x = (x_1, \dots, x_n )$ des coordonnées sur l'espace réel ou
complexe $\R^n$ ou $\C^n$ de dimension $n \geqslant 1$.  Si $(x_1',
\dots, x_n')$ sont des coordonn\'ees sur une 
seconde copie auxiliaire d'un même espace r\'eel ou complexe \`a $n$
dimensions, on peut consid\'erer les applications:
\[
x_1'
= 
f_1(x_1,\dots,x_n),\,\,
\cdots\cdots,\,\,
x_n'
=
f_n(x_1,\dots,x_n)
\]
d\'efinies dans le premier espace et
\`a valeurs dans le second espace, o\`u  les 
$f_i$ sont des fonctions quelconques de leurs $n$ arguments $(x_1,
\dots, x_n)$. L'application:
\[
x
\longmapsto 
f(x) 
=
x' 
\]
est alors appel\'ee une {\em transformation ponctuelle} lorsqu'elle
transforme de mani\`ere biunivoque et 
diff\'erentiable\footnote{\,
%%%%%%%%%%%%%%%%%%%%%%%-------DEBUT--------%%%%%%%%%%%%%%%%%%%%%%%%%%%
Au moins de classe $\mathcal{ C}^1$, c'est-\`a-dire que chaque
composante $f_i$ de l'application possède des d\'eriv\'ees
partielles $\frac{ \partial f_i}{ \partial x_k} ( x_1, \dots, x_n)$
dans toutes les directions de coordonn\'ees $x_k$, et que ces
d\'eriv\'ees partielles sont partout {\em continues}.
} %%%%%%%%%%%%%%%%%%%%%%%%-----FIN-----%%%%%%%%%%%%%%%%%%%%%%%%%%%%%%% 
la totalit\'e
des points $(x_1, \dots, x_n)$ de l'espace-source en des points
$(x_1', \dots, x_n')$ de l'espace d'arriv\'ee de telle sorte que
son application inverse:
\[
x'
\longmapsto 
f^{-1}(x')
=
x 
\]
soit elle aussi diff\'erentiable\footnote{\,
%%%%%%%%%%%%%%%%%%%%%%%-------DEBUT--------%%%%%%%%%%%%%%%%%%%%%%%%%%%
\`A nouveau: au moins de classe $\mathcal{ C}^1$.
} %%%%%%%%%%%%%%%%%%%%%%%%-----FIN-----%%%%%%%%%%%%%%%%%%%%%%%%%%%%%%% 
et biunivoque. Aujourd'hui, on dit que $x \mapsto f ( x) = x'$ est un
diff\'eomorphisme (de $\R^n$ ou $\C^n$ sur $\R^n$ ou $\C^n$,
respectivement). Plus généralement, on 
considère des difféomorphismes entre
sous-ensembles ouverts quelconques
$U$ et $U'$ (de $\R^n$ ou de $\C^n$). 

Dans la théorie de Lie des groupes continus de transformations,
l'objet-archétype est constitué de {\sl transformations
ponctuelles}:
\[
x_i'
=
f_i(x_1,\dots\,x_n;\,a_1,\dots,a_r)
\ \ \ \ \ \ \ \ \ \ 
{\scriptstyle{(i\,=\,1\,\cdots\,n)}},
\]
qui sont de plus
{\em paramétrées} par un nombre
fini $r$ de paramètres réels ou complexes
$(a_1, \dots, a_r )$, à savoir pour chaque $a$ fixé, 
chaque application: 
\[
x'= 
f(x;\,a)
=:
f_a(x)
\]
est supposée constituer un {\em difféomorphisme}
d'un certain domaine\footnote{\,
D'après une définition standard de topologie
générale,
un {\sl domaine} est un ouvert {\em connexe} non vide.
} %%%%%%%%%%%%%%%%%%%%%%%%%%%%%%%%%%%%%%%%%%%%%%%%%%%%%%%%%%% de
de l'espace-source sur un autre domaine situé dans un espace d'arrivée
ayant la même dimension $n$ et muni de coordonnées $(x_1', \dots,
x_n')$.  Ainsi, le déterminant jacobien:
\[
\small
\aligned
\det{\rm Jac}(f)
=
\left\vert
\begin{array}{cccc}
\frac{\partial f_1}{\partial x_1}
&\cdots\cdots&
\frac{\partial f_1}{\partial x_n}
\\
\vdots
&\ddots&
\vdots
\\
\frac{\partial f_n}{\partial x_1}
&\cdots\cdots&
\frac{\partial f_n}{\partial x_n}
\end{array}
\right\vert
=
\sum_{\sigma\in{\sf Perm}_n}\,
{\rm sgn}(\sigma)\,
\frac{\partial f_1}{\partial x_{\sigma(1)}}\,
\frac{\partial f_2}{\partial x_{\sigma(2)}}
\cdots
\frac{\partial f_n}{\partial x_{\sigma(n)}}
\endaligned
\]
ne s'annule en aucun point du domaine source.  Avant d'introduire les
axiomes de groupe ({\em voir} \S~4.2), la première question à résoudre
est: de combien de paramètres exactement les transformations $x_i ' =
f_i ( x;\, a)$ dépendent-elles réellement?  Certains paramètres
pourraient en effet être superflus, et devraient à ce titre être
supprimés à l'avance. Or à cette fin, il est crucial de formuler
explicitement dès le début et une fois pour toutes {\em trois
principes de pensée concernant l'admission des hypothèses
fondamentales qui gouvernent la théorie des groupes continus
développée par Lie}.

\medskip\noindent{\bf Hypothèse générale d'analyticité:}
Courbes, surfaces, variétés, sous-variétés, équations de
transformation, groupes, sous-groupes, \etcv
tous les objets
mathématiques de la théorie seront supposés {\em analytiques} (réels
ou complexes), c'est-à-dire que les fonctions qui les représentent
dans des systèmes de coordonnées locales seront systématiquement
supposées développables en série entière convergente dans un certain
domaine de $\R^k$ ou de $\C^k$, 
pour un certain $k \geqslant 1$. 

\medskip\noindent{\bf Principe de relocalisation libre 
au voisinage d'un point générique:} \label{relocalisation-libre}
Considérons un objet mathématique local représenté par des
fonctions qui sont analytiques dans un domaine $U_1$, et supposons
qu'un certain <<\,bon\,>> comportement <<\,générique\,>> se produise
dans $U_1 \backslash {\sf D}_1$ en dehors d'un certain sous-ensemble
analytique ${\sf D}_1 \subset U_1$; par exemple, l'inversion d'une
matrice carrée qui est constituée de fonctions analytiques est
possible seulement en dehors de l'ensemble ${\sf D}_1$ 
des zéros de son
déterminant\footnote{\,\label{definition-rang-generique}
%%%%%%%%%%%%%%%%%%%%%%%-------DEBUT--------%%%%%%%%%%%%%%%%%%%%%%%%%%%
Plus généralement, considérons une matrice rectangulaire $G(y) := (
g_i^j (y))_{ 1\leqslant i \leqslant n}^{ 1 \leqslant j \leqslant m}$
de taille $n \times m$ dont les éléments $g_i^j = g_i^j ( y)$ sont des
fonctions analytiques d'un certain nombre de variables $y = ( y_1,
\dots, y_q)$ qui sont définies dans un 
certain domaine $U$ de $\C^q$ (ou de $\R^q$). Pour tout entier $\rho$
tel que $1 \leqslant \rho \leqslant
\min ( m, n)$, on peut former la collection de tous les déterminants
de taille $\rho \times \rho$ (mineurs) qui sont extraits de cette
matrice:
\[
\Delta_{i_1,\dots,i_\rho}^{j_1,\dots,j_\rho}(y)
:=
\left\vert
\begin{array}{ccc}
g_{i_1}^{j_1}(y) & \cdots & g_{i_1}^{j_\rho}(y)
\\
\cdot\cdot & \cdots & \cdot\cdot
\\
g_{i_\rho}^{j_1}(y) & \cdots & g_{i_\rho}^{j_\rho}(y)
\end{array}
\right\vert. 
\]
En partant de $\rho := \min ( m, n)$, 
si tous ces déterminants sont identiquement nuls (en tant que
fonctions de la variable $y$), on passe alors
de la taille $\rho$ à la taille
juste inférieure $\rho-1$, on 
forme tous les mineurs, on teste leur annulation 
identique, 
et on recommence. 
Le {\em rang générique} $\rho^*$ de la matrice $G
(y)$ est alors 
le plus grand entier $\rho$ tel qu'il existe au moins un tel
mineur non identiquement nul, tous
les mineurs d'une taille
strictement supérieure étant
identiquement nuls. On a $\rho^* = 0$ si et seulement si
toutes les fonctions $g_i^j ( y)$ sont nulles (cas inintéressant), et
sinon, on a en toute généralité: $1 \leqslant \rho^* \leqslant
\min ( m, n)$. Enfin, si l'on introduit le lieu:
\[
{\sf D}^*
:=
\Big\{
y\in\C^q:\,
\Delta_{i_1,\dots,i_{\rho^*}}^{j_1,\dots,j_{\rho^*}}
(y)
=
0,\ 
\forall\,
i_1,\dots,i_{\rho^*},\
\forall\,
j_1,\dots,j_{\rho^*}
\Big\}
\]
des points $y$ en lesquels tous les mineurs de taille $\rho^* \times
\rho^*$ s'annulent, alors ce lieu ${\sf D}^*$ est un sous-ensemble
analytique {\em propre}\,\,---\,\,en particulier fermé et de
complémentaire $U \backslash {\sf D}^*$ ouvert et dense\,\,---\,\,qui
a par définition la propriété qu'en tout point $y \in U \backslash
{\sf D}^*$, au moins un mineur de taille $\rho^*
\times \rho^*$ ne s'annule pas, et comme tous les mineurs de taille
strictement supérieure s'annulent identiquement par construction, on
en déduit la propriété remarquable: {\em en tout point $y \in U
\backslash {\sf D}^*$, le rang de la matrice $G ( y)$ est maximal,
égal à son rang générique $\rho^*$}. Cas particulier: lorsqu'on a une
matrice carrée de déterminant non identiquement nul,
\ie $\rho^* = m = n$, cette matrice est inversible 
en tout point $y \in U \backslash {\sf D}^*$. 
}, %%%%%%%%%%%%%%%%%%%%%%%%-----FIN-----%%%%%%%%%%%%%%%%%%%%%%%%%%%%%%%
qui est une fonction analytique. Alors on s'autorise à
relocaliser les considérations dans tout sous-domaine $U_2 \subset U_1
\backslash {\sf D}_1$.

\begin{center}
\input delocalization.pstex_t
\end{center}

\noindent
Ensuite, dans $U_2$, des raisonnements ultérieurs peuvent demander
que l'on évite un autre sous-ensemble analytique propre ${\sf D}_2$,
et donc l'on doit à nouveau relocaliser les considérations dans un
sous-domaine $U_3 \subset U_2 \backslash {\sf D}_2$, et ainsi de
suite.  La plupart des démonstrations de la {\em Théorie des groupes de
transformations}, et tout particulièrement les théorèmes de
classification, autorisent un certain nombre de telles
relocalisations, souvent sans aucune mention de la part de Engel et
Lie, un tel {\em acte de pensée} étant considéré comme justifié 
\`a l'avance par
la nécessité d'{\em étudier en premier lieu les objets génériques},
eux-mêmes déjà extrêmement riches et diversifiés.

Toutefois, lorsque nous traduirons et ré-énoncerons 
au Chapitre~5
les principaux (et
magistraux) théorèmes de classification qui, dans le Tome~III,
précèdent les solutions complexes qu'Engel et Lie apportent au
problème de Riemann-Helmholtz, nous rappellerons
explicitement, précisément et rigoureusement toutes les hypothèses de
généricité qui sont requises pour la validité des résultats, afin
qu'il ne subsiste aucune ambiguïté d'interpr\'etation
pour le lecteur contemporain.

\medskip\noindent{\bf Ne nommer aucun domaine d'existence:}
Sans introduire de dénomination ou de notation spécifique, Engel et
Lie écrivent ordinairement {\em le} voisinage
\deutsch{\emphasis{der} Umgebung} d'un point
(ou au sens absolu), de la même manière que l'on parle {\em du}
voisinage d'une maison, ou {\em des} environs d'une ville, tandis que
la topologie contemporaine conceptualise {\em un} voisinage
(habituellement <<\,petit\,>>), parmi une infinité de voisinages
existants.  Contrairement à ce que la mythologie formaliste du
20\textsuperscript{ème} siècle aime à faire accroire depuis les
controverses entre les géomètres de Leipzig et l'école de Berlin,
Engel et Lie ont toujours mis en évidence, lorsque c'était nécessaire,
la nature locale des concepts de la théorie des groupes de
transformations. Nous illustrerons notamment cette constation en
étudiant les tentatives que Lie a entreprises pour économiser l'axiome
d'existence des éléments inverses dans sa théorie des groupes continus de
transformations (\S~4.3). Certes, il est vrai que les résultats de Lie,
gagneraient en pr\'ecision technique \`a être \'enonc\'es en
sp\'ecifiant tous les domaines d'existence, mais il est plausible que
Engel
et Lie ait rapidement réalisé que le fait de ne pas donner de nom aux
voisinages, et d'éviter de la sorte tout symbolisme essentiellement
superflu, était le moyen le plus efficace pour conduire à leur terme
des théorèmes de classification de grande envergure.

\smallskip

En conséquence, nous adopterons le style (économique) de pensée de
Engel et Lie, et nous présupposerons, sans pour autant effectuer des
rappels fréquents, que:

\begin{itemize}

\smallskip\item[$\bullet$]
tous les objets mathématiques sont analytiques; 

\smallskip\item[$\bullet$]
les relocalisations génériques sont librement autorisées;

\smallskip\item[$\bullet$]
les ensembles ouverts sont souvent petits, rarement nommés, 
non vides, et toujours
\emphasis{connexes}.

\end{itemize}

%%%%%%%%%%%%%%%%%%%%%%%%%%%%%%%%%%%%%%%%%%%%%%%%%%%%%%%%%%%%%%%%%%%%%

\newpage

\thispagestyle{empty}
\setcounter{footnote}{0}

$\:$

\bigskip\bigskip\bigskip

\begin{center}
{\large\bf Chapitre~3:

\smallskip
Théorèmes fondamentaux

\smallskip
sur les groupes de transformations}
\label{Chapitre-3}
\end{center}

\bigskip

\HEAD{Chapitre~3.\,\,\,\,Théorèmes fondamentaux sur les groupes
de transformations}{
3.1.\,\,\,Paramètres essentiels}

\medskip\noindent{\bf 3.1.~Paramètres essentiels.}
\label{parametres-essentiels}
Comme exemple d'application de ces trois principes généraux qui
gouvernent la pensée de Lie, montrons comment on peut s'assurer que
tous les paramètres dans une collection d'équations de transformation
sont effectivement présents; si tel n'est pas le cas, montrons comment
on peut ramener ces équations à d'autres équations de transformations
équivalentes qui comportent un nombre inférieur de paramètres. Ce
problème à résoudre obéit à une exigence incontournable
d'économie. Supprimer à l'avance tous les paramètres superflus
lorsqu'il en existe permettra certainement d'éviter d'avoir à traiter
des cas parasites dans l'énoncé des théorèmes principaux de la
théorie.

On supposera pour fixer les idées que les variables $x = (x_1, \dots,
x_n) \in \C^n$ sont complexes et que les paramètres $a = ( a_1, \dots,
a_r)$ sont eux aussi complexes, étant entendu que la théorie
fondamentale est inchangée\footnote{\,
%%%%%%%%%%%%%%%%%%%%%%%-------DEBUT--------%%%%%%%%%%%%%%%%%%%%%%%%%%%
Lorsqu'ils ne précisent pas le domaine de variation des quantités
numériques, Engel et Lie sous-entendent qu'elles sont complexes. Pour
les théorèmes de classification sur $\R$ qui prolongent et qui
utilisent les théorèmes de classification sur $\C$, les deux auteurs
précisent alors explicitement que variables et paramètres sont tous
deux {\em réels}.
} %%%%%%%%%%%%%%%%%%%%%%%%-----FIN-----%%%%%%%%%%%%%%%%%%%%%%%%%%%%%%%
lorsque $x \in \R^n$ et $a \in
\R^r$, ou même lorsque\footnote{\,
%%%%%%%%%%%%%%%%%%%%%%%-------DEBUT--------%%%%%%%%%%%%%%%%%%%%%%%%%%%
Ce cas intermédiaire n'est en général pas étudié par Engel et Lie,
mais il le sera par \'Elie Cartan.
} %%%%%%%%%%%%%%%%%%%%%%%%-----FIN-----%%%%%%%%%%%%%%%%%%%%%%%%%%%%%%%
$x \in \C^n$ et $a \in \R^r$. 

Dans des équations de transformations quelconques $x' = f ( x; \, a)$,
toutes les variables $(x_1, \dots, x_n)$ sont par hypothèse présentes,
puisque $x \mapsto f_a ( x)$ est un difféomorphisme local pour tout
$a$. Mais aucune supposition, excepté la finitude, n'a été faite sur
les paramètres $(a_1, \dots, a_r)$.

L'idée principale consiste à développer les fonctions $f_i$ des
équations de transformation:
\[
x_i '
=
f_i(x_1,\dots,x_n;\,a_1,\dots,a_r)
\ \ \ \ \ \ \ \ \ \ \ \ \
{\scriptstyle{(i\,=\,1\,\cdots\,n)}}
\] 
en série entière par rapport à $x - x_0$ dans un voisinage connexe
(non nommé) d'un point fixé $x_0 = (x_{ 10}, \dots, x_{ n0})$, 
ce qui nous donne:
\[
f_i(x;\,a)
=
\sum_{\alpha\in\N^n}\,\mathcal{F}_\alpha^i(a)\,
(x-x_0)^\alpha;
\]
ici, pour tout multiindice $\alpha = ( \alpha_1, \dots, \alpha_n) \in
\N^n$, nous avons noté de manière abrégée le monôme:
\[
(x-x_0)^\alpha
=
(x_1-x_{10})^{\alpha_1}
\cdots
(x_n-x_{n0})^{\alpha_n}.
\]
Comme coefficients des (multi-)puissances $(x-x_0)^\alpha$ dans une
telle série entière, il apparaît un nombre infini de certaines
fonctions analytiques $\mathcal{ F}_\alpha^i = \mathcal{ F}_\alpha^i
(a)$ des paramètres $(a_1, \dots, a_r)$ qui sont définies dans un
domaine fixe de $\C^r$. Alors nous prétendons que les paramètres $a_k$
qui sont éventuellement
superflus peuvent être détectés en étudiant le rang
générique de l'{\sl application infinie des coefficients}:
\[
{\sf F}_\infty:
\ \ \ \ \ \ \ \
\C^r\ni\ a
\longmapsto
\big(
\mathcal{F}_\alpha^i(a)
\big)^{\alpha\in\N^n,\,1\leqslant i\leqslant n}\
\in\C^\infty,
\]
c'est-à-dire, d'après la définition même du rang (générique) d'une
application, en étudiant le rang générique\footnote{\,
%%%%%%%%%%%%%%%%%%%%%%%-------DEBUT--------%%%%%%%%%%%%%%%%%%%%%%%%%%%
Pour une définition,
\voir~la note p.~\pageref{definition-rang-generique}, 
dont les considérations se généralisent sans modification au cas où il
y a une infinité de colonnes. Techniquement, $\rho_\infty$ est défini
comme le plus grand entier $\rho \leqslant \min ( r, \infty)$ tel
qu'il existe au moins un mineur de taille $\rho \times \rho$ dans
cette matrice jacobienne $\big( \frac{\partial \mathcal{U }_\alpha^i}{
\partial a_j} (a)\big)_{1 \leqslant j\leqslant r}^{
\alpha \in \N^n,\, 1\leqslant i\leqslant n
}$ (de taille $r \times \infty$) qui ne s'annule pas identiquement en
tant que fonction de $a$. Alors tous les mineurs de taille
$(\rho_\infty + 1)
\times ( \rho_\infty + 1)$ s'annulent
identiquement.
} %%%%%%%%%%%%%%%%%%%%%%%%-----FIN-----%%%%%%%%%%%%%%%%%%%%%%%%%%%%%%% 
$\rho_\infty$ de sa matrice jacobienne: 
\[
{\rm Jac}\,{\sf F}_\infty(a)
=
\Big(
\frac{\partial\mathcal{F}_\alpha^i}{\partial a_j}(a)
\Big)_{1\leqslant j\leqslant r}^{
\alpha\in\N^n,\,1\leqslant i\leqslant n},
\]
qui est considérée ici comme une matrice ayant $r$ lignes indexées par
l'entier $j$, et possédant une infinité de colonnes, indexées
simultanément par le multiindice $\alpha$ et par l'entier $i$.

\`A première vue, ces considérations techniques peuvent 
paraître complexes, mais il n'en est rien, puisqu'il est au contraire tout
à fait naturel que la manière dont les fonctions $f_i (x; \, a)$
dépendent réellement des paramètres $(a_1, \dots, a_r)$ doive être
déchiffrée sur l'ensemble des fonctions-coefficients $\mathcal{
F}_\alpha^i ( a_1, \dots, a_r)$.

Si par exemple il existe un paramètre, disons $a_1$, tel qu'aucune
fonction $f_i(x; \, a)$ n'en dépend, d'où il découle que tous les
coefficients $\mathcal{ F}_\alpha^i (a)$ sont indépendants de $a_1$,
alors cette application ${\sf F}_\infty$ a un rang qui est
trivialement $\leqslant r - 1$, parce que la première ligne de sa
matrice jacobienne ${\rm Jac}\, {\sf F}_\infty (a)$ est alors
identiquement nulle, et puisque le rang générique d'une matrice est en
toute circonstance inférieur au nombre de ses lignes.

En toute généralité, si $\rho_\infty$ désigne le rang générique de
${\rm Jac}\, {\sf F}_\infty$ et si on définit le sous-ensemble ${\sf
D}_\infty$ des zéros communs à tous ses mineurs
de taille $\rho_\infty \times \rho_\infty$:
\[
\left\vert
\begin{array}{ccc}
\frac{\partial\mathcal{F}_{\alpha_1}^{i_1}}{\partial a_{j_1}}(a)
& \cdots &
\frac{\partial\mathcal{F}_{\alpha_1}^{i_1}}{
\partial a_{j_{\rho_\infty}}}(a)
\\
\cdot\cdot & \cdots & \cdot\cdot
\\
\frac{\partial\mathcal{F}_{\alpha_{\rho_\infty}}^{i_{\rho_\infty}}}{
\partial a_{j_1}}(a)
& \cdots &
\frac{\partial\mathcal{F}_{\alpha_{\rho_\infty}}^{i_{\rho_\infty}}}{
\partial a_{j_{\rho_\infty}}}(a)
\end{array}
\right\vert,
\]
alors ${\sf D}_\infty$ est un sous-ensemble analytique propre de
l'espace des paramètres $a$, et pour tout paramètre $a_0$ qui
n'appartient pas à ce sous-ensemble de paramètres
<<\,exceptionnels\,>>, on peut tout
d'abord {\em relocaliser les considérations}
dans un petit voisinage convenable\footnote{\,
%%%%%%%%%%%%%%%%%%%%%%%-------DEBUT--------%%%%%%%%%%%%%%%%%%%%%%%%%%%
Ici et dans la ce qui va suivre, ce qui peut être dit au voisinage
des paramètres exceptionnels qui appartiennent à ${\sf D}_\infty$
requerrait des outils sophistiqués de la théorie des singularités qui
sont au-delà de la portée du présent travail. }
%%%%%%%%%%%%%%%%%%%%%%%%-----FIN-----%%%%%%%%%%%%%%%%%%%%%%%%%%%%%%%
de $a_0$, et ensuite, grâce à une application du théorème du rang
constant, on peut trouver un
difféomorphisme approprié de l'espace des paramètres $a
\mapsto \overline{ a} = \overline{ a} ( a)$ fixant $a_0$, ce
qui donne de nouvelles équations de transformations:
\[
\aligned
x_i '
&
=
f_i
\big(
x_1,\dots,x_n;\,a_1(\overline{a}),\dots,a_r(\overline{a})
\big)
\\
&
=:
\overline{f}_i
\big(
x_1,\dots,x_n;\,\overline{a}_1,\dots,\overline{a}_r
\big)
\ \ \ \ \ \ \ \ \ \ \ \ \
\ \ \ \ \ \ \ \ \ \ \ \ \
{\scriptstyle{(i\,=\,1\,\cdots\,n)}},
\endaligned
\]
de manière à ce que les nouveaux coefficients $\overline{
\mathcal{ F}}_\alpha^i \, \big( \overline{ a} \big)$
obtenus en développant les nouvelles fonctions $\overline{ f}_i$ en
série entière par rapport à $(x - x_0)$ deviennent {\em absolument
indépendants des $r - \rho_\infty$ derniers paramètres $\overline{
a}_{ \rho_\infty + 1}, \dots, \overline{ a}_r$},
qui sont ainsi 
devenus visiblement superflus. Sans démonstration
(\voir~\cite{ merk2009b}), formulons
l'énoncé rigoureux de réduction.

\smallskip\noindent{\bf Théorème.}
{\em Localement au voisinage de tout paramètre générique $a_0$ en
lequel l'application infinie des coefficients $a \mapsto {\sf
F}_\infty ( a)$ est de rang maximal, constant, et égal à son rang
générique $\rho_\infty$, il existe à la fois un changement local de
paramètres $a \mapsto
\big( b_1 ( a), \dots, b_{\rho_\infty} ( a) \big) =: (b_1,
\dots, b_{ \rho_\infty })$ 
qui fait baisser le nombre des paramètres de $r$ à $\rho_\infty$, et
des nouvelles équations de transformations:
\[
x_i'
=
g_i
\big(x;\,b_1,\dots,b_{\rho_\infty}\big)
\ \ \ \ \ \ \ \ \ \
{\scriptstyle{(i\,=\,1\,\cdots\,n)}}
\]
dépendant de seulement $\rho_\infty$ paramètres qui redonnent par
substitution les anciennes équations de transformations:}
\[
g_i\big(x;\,b(a)\big)
\equiv 
f_i(x;\, a)
\ \ \ \ \ \ \ \ \ \
{\scriptstyle{(i\,=\,1\,\cdots\,n)}}.
\]

\smallskip\noindent{\bf Définition.}
Les paramètres $(a_1, \dots, a_r)$ d'équations de transformations
données $x_i' = f_i ( x, a)$, $i = 1,
\dots, n$, sont dits {\sl essentiels}
si, après développement en série entière $f_i ( x, a) = \sum_{
\alpha \in \N^n}\, \mathcal{ F}_\alpha^i ( a) \, (x - x_0)^\alpha$
autour d'un point $x_0$, le rang générique $\rho_\infty$ de
l'application infinie des coefficients $\text{\sf F}_\infty: \, a
\longmapsto \big( \mathcal{ F}_\alpha^i ( a) \big)^{ \alpha \in
\N^n,\, 1 \leqslant i \leqslant n}$ est maximal possible, égal au
nombre $r$ des paramètres, à savoir: $\rho_\infty = r$.
\medskip

Dans ce cas, Engel et Lie
disent des équations de transformations qu'elles comportent
{\sl $r$-termes}\, \deutsch{$r$-gliedrig}. Grâce
à ce théorème, on peut toujours\,\,---\,\,pourvu que l'on admette le
principe de relocalisation libre en un point
générique\footnote{\,
%%%%%%%%%%%%%%%%%%%%%%%-------DEBUT--------%%%%%%%%%%%%%%%%%%%%%%%%%%%
Sans admettre ce principe, la théorie des singularités devrait être
convoquée. Néanmoins, lorsque les équations de $x_i' = f_i ( x; \,
a)$ constituent un groupe continu fini local de transformations au
sens qui est défini dans le \S~3.2 ci-dessous, 
on vérifie techniquement (\cf~\cite{ merk2009b}) qu'aucune
relocalisation n'est nécessaire pour éliminer les paramètres
superflus, la raison <<\,philosophique\,>> étant que tout groupe agit
transitivement sur lui-même, ce qui force les rangs (génériques) à
être {\em constants}.
}%%%%%%%%%%%%%%%%%%%%%%%%-----FIN-----%%%%%%%%%%%%%%%%%%%%%%%%%%%%%%%
\,\,---\,\,supposer sans aucune perte de généralité que les
paramètres de toutes les équations de transformations sont {\em
essentiels}. Du point de vue de la philosophie des mathématiques,
le paramétrique peut être et doit être ainsi définitivement
réduit à sa qualité essentielle. 
Par conséquent, 
dans ce qui va suivre, les paramètres seront toujours supposés
essentiels. 

Pour terminer sur l'essentialité des paramètres, 
voici encore un critère effectif qui
sera utilisé plus loin de manière incidente, mais que nous énoncerons
sans fournir de démonstration (\voir~\cite{ merk2009b}) afin de ne pas
retarder la présentation des concepts centraux.

\medskip\noindent{\bf Théorème.}
\label{theoreme-essentiel}
{\em Les trois conditions suivantes sont équivalentes:

\begin{itemize}

\smallskip\item[{\bf (i)}]
Dans les équations de transformations:
\[
x_i'
=
f_i(x_1,\dots,x_n;\,a_1,\dots,a_r)
=
\sum_{\alpha\in\N^n}\,\mathcal{F}_\alpha^i(a)\,
(x-x_0)^\alpha
\ \ \ \ \ \ \ \ \ \
{\scriptstyle{(i\,=\,1\,\cdots\,n)}},
\]
les paramètres $a_1, \dots, a_r$ {\em ne} sont {\em pas} 
essentiels. 

\smallskip\item[{\bf (ii)}]
(Par définition) Le rang
générique $\rho_\infty$ de la matrice
jacobienne infinie:
\[
{\rm Jac}\,{\sf F}_\infty(a)
=
\Big(
\frac{\partial\mathcal{F}_\alpha^i}{\partial a_j}(a)
\Big)_{1\leqslant j\leqslant r}^{\alpha\in\N^n,\,1\leqslant i\leqslant n}
\]
est strictement inférieur à $r$.

\smallskip\item[{\bf (iii)}]
Localement dans un voisinage de {\rm tout} $(x_0, a_0)$,
il existe un champ de vecteurs non identiquement nul
sur l'espace des paramètres: 
\[
\mathcal{T}
=
\sum_{k=1}^n\,\tau_k(a)\,\frac{\partial}{\partial a_k},
\]
dont les coefficients $\tau_k ( a)$ sont analytiques dans un voisinage
de $a_0$, qui annihile toutes les fonctions $f_i ( x; a)$:
\[
0
\equiv
\mathcal{T}\,f_i
=
\sum_{k=1}^n\,\tau_k\,\frac{\partial f_i}{\partial a_k}
=
\sum_{\alpha\in\N^n}\,\sum_{k=1}^r\,\tau_k(a)\,
\frac{\partial\mathcal{F}_\alpha^i}{\partial a_k}(a)\,
(x-x_0)^\alpha
\ \ \ \ \ \ \ \ \ \ \ \ \
{\scriptstyle{(i\,=\,1\,\cdots\,n)}}.
\]
\end{itemize}}\smallskip

\noindent
Lorsque l'on sait déjà (\cf~le premier théorème) que les $r -
\rho_\infty$ derniers paramètres $(a_{ \rho_\infty + 1}, \dots, a_r)$
sont purement absents de toutes les fonctions $f_i$, cette dernière
condition {\bf (iii)} exprime tout simplement 
le fait évident que
les fonctions $f_i$ sont annulées
identiquement par $\partial \big/ \partial a_{ \rho_\infty+1},
\dots, \partial \big/ \partial a_r$.

Plus généralement, si $\rho_\infty$ désigne
le rang générique de l'application infinie
des coefficients: 
\[
{\sf F}_\infty:\ \ \ \ \ \ \
a\ \ 
\longmapsto\ \
\big(
\mathcal{F}_\alpha^i(a)
\big)_{\alpha\in\N^n}^{1\leqslant i\leqslant n},
\]
et si $\rho_\infty \leqslant r-1$,
alors on démontre que localement dans un voisinage de {\em tout}
$(x_0, a_0)$, il existe exactement $r - \rho_\infty$
champs de vecteurs analytiques (mais pas plus):
\[
\mathcal{T}_\mu
=
\sum_{k=1}^n\,\tau_{\mu k}(a)\,
\frac{\partial}{\partial a_k}
\ \ \ \ \ \ \ \ \ \
{\scriptstyle{(\mu\,=\,1\,\cdots\,r-\rho_\infty)}},
\]
satisfaisant les deux propriétés suivantes:

\smallskip$\bullet$\,\,
la dimension de l'espace vectoriel qu'ils engendrent
\[
{\rm Vect}
\big(
\mathcal{T}_1\big\vert_a,
\dots,
\mathcal{T}_{r-\rho_\infty}\big\vert_a
\big)
\] 
est égale à $r - \rho_\infty$ en tout paramètre $a$ en lequel le rang
de ${\sf F}_\infty$ est maximal égal à $\rho_\infty$;

\smallskip$\bullet$\,\,
chacune de ces dérivations $\mathcal{ T}_\mu$
annihile identiquement toutes les fonctions $f_i ( x; a)$:
\[
0
\equiv
\mathcal{T}_\mu\,f_i
=
\sum_{k=1}^r\,\tau_{\mu k}(a)\,
\frac{\partial f_i}{\partial a_k}(x;a)
\ \ \ \ \ \ \ \ \ \
{\scriptstyle{(i\,=\,1\,\cdots\,n;\,\,\,\mu\,=\,1\,\cdots\,r-\rho_\infty)}}.
\]

\HEAD{Chapitre~3.\,\,\,\,Théorèmes fondamentaux sur les groupes
de transformations}{
3.2.\,\,\,Concept de groupe de Lie local}

\medskip\noindent 
{\bf 3.2.~Concept de groupe de Lie local.}
\label{concept-de-groupe-de-lie-local}
Nous restituons ici les définitions fondamentales et les théorèmes,
sans insister sur l'aspect purement technique qui exige, en toute
rigueur, de préciser quand et comment l'on doit rapetisser les ouverts
dans lesquels les objets sont définis. Dans les premières pages de la
{\em Theorie der Transformationsgruppen}, Engel et Lie 
exposent le principe général des 
voisinages emboîtés, avant de choisir
de se dispenser de les nommer afin d'alléger la présentation. 

En dimension $n \geqslant 1$ arbitraire, un {\sl groupe continu fini
de transformations} sur $\C^n$ est une famille de transformations
ponctuelles analytiques paramétrée par un nombre fini $r$ de
paramètres:
\[
x_i'
=
f_i(x_1,\dots,x_n;a_1,\dots,a_r)
\ \ \ \ \ \ \ \ \ \ 
{\scriptstyle{(i\,=\,1\,\cdots\,n)}}
\]
qui satisfait les trois propriétés suivantes.

\smallskip\noindent
{\small\sf\em Loi de composition:}
\`A chaque fois qu'elle est bien définie, la
succession $x ' = f ( x; a)$ et $x''
= f ( x'; b)$ de deux telles transformations, à savoir:
\[
x''
=
f\big(
f(x;a);b
\big)
=
f(x;c)
\]
s'identifie toujours à une transformation
de la {\em même} famille, pour
un certain nouveau paramètre: 
\[
c
= 
\text{\bf m}(a,b) 
\]
qui est défini de manière précise et unique par une certaine
application analytique locale $\text{\bf m} : \C^r \times \C^r \to
\C^r$, laquelle, de son
propre côté, hérite automatiquement de la propriété d'associativité
qu'ont les difféomorphismes par composition:
\[
\text{\bf m}
\big(\text{\bf m}(a,b),c\big) 
= 
\text{\bf m}\big(a,\text{\bf m}(b,c)\big).
\] 

\smallskip\noindent
{\small\sf\em Existence d'un élément identité:} Il existe un paramètre
spécial $e = (e_1, \dots, e_r)$ tel que $f ( x; e) \equiv x$
se réduit à l'application identique. 

\smallskip\noindent
{\small\sf\em Loi de multiplication de groupe sous-jacente:} L'application
analytique $( a, b)
\longmapsto \text{\bf m} ( a, b)$, 
qui peut aussi être écrite de manière alternative et plus brève en
utilisant un symbole de type multiplication: $( a, b) \longmapsto a
\cdot b$, est une {\sl loi de groupe} continue, au sens où:

\smallskip$\bullet$
Pour tout $a$, on doit avoir: $a \cdot e = e \cdot a = a$, une
propriété qui découle en fait de la loi de composition:
\[
f(x;a\cdot e)
=
f\big(f(x;a);e\big)
=
f( x; a) 
= 
f\big( f (x; e ); a \big)
=
f(x;e\cdot a),
\]
grâce à l'unicité postulée à l'instant de $c = a \cdot b$.

\smallskip$\bullet$
De même, l'associativité héritée $( a \cdot b) \cdot c = a \cdot (
b\cdot c)$ doit valoir.

\smallskip$\bullet$
{\small\sf\em Existence d'éléments inverses:} Enfin, comme dernier axiome
qui ne s'avère pas être conséquence 
élémentaire de la loi de multiplication de
groupe, il doit exister un difféomorphisme local analytique $\text{\bf
i} : \C^r \to \C^r$ défini dans un voisinage de $e$ avec $\text{\bf i}
( e) = e$ tel que:
\[
a\cdot\text{\bf i}(a) 
=
\text{\bf i}(a)\cdot 
a 
= 
e,
\]
à savoir: $\text{\bf i} ( a)$ représente l'{\sl inverse de $a$ pour la
structure de groupe}, et de plus, $a \mapsto
\text{\bf i} ( a)$ est une application analytique, nécessairement
un difféomorphisme local au voisinage de $e$. En particulier, si l'on
réécrit maintenant que la composition\footnote{\,
%%%%%%%%%%%%%%%%%%%%%%%-------DEBUT--------%%%%%%%%%%%%%%%%%%%%%%%%%%%
Dans la suite, pour tout $a$ fixé, le difféomorphisme local $x \mapsto
f( x; a)$ sera occasionnellement écrit $x \mapsto f_a ( x)$.
}: %%%%%%%%%%%%%%%%%%%%%%%%-----FIN-----%%%%%%%%%%%%%%%%%%%%%%%%%%%%%%%
\[
f\big(f(x;a);b\big) 
= 
f(x;a\cdot b)
\] 
est exécutée par la multiplication de groupe entre les paramètres, on
en déduit formellement que:
\[
f\big(f(x;a);
\text{\bf i}(a)\big)
\equiv
f(x;a\cdot\text{\bf i}(a))
\equiv
x
\equiv 
f(x;\text{\bf i}(a)\cdot a)
\equiv
f\big(f(x;\text{\bf i}(a));a\big), 
\]
ce qui montre que $x \mapsto f_{ \text{\bf i} ( a)} (x)$ est le
difféomorphisme inverse de $x \mapsto f_a ( x)$.

Il sera utile pour la suite de mémoriser le fait que dans la
terminologie de Lie, un <<\,{\em groupe continu fini de
transformations}\,>> signifie précisément une action analytique locale
effective\footnote{\,
%%%%%%%%%%%%%%%%%%%%%%%-------DEBUT--------%%%%%%%%%%%%%%%%%%%%%%%%%%%
Est effective (\cite{ ol1995}) une action locale 
de groupe de Lie local telle
que pour tout paramètre $a$ distinct du paramètre identité, le
difféomorphisme local $x \mapsto f_a ( x)$ ne se réduit pas à
l'identité. On démontre aisément qu'une action analytique de groupe de
Lie est effective si et seulement si les paramètres de la famille de
transformations qu'elle définit sont essentiels. }
%%%%%%%%%%%%%%%%%%%%%%%%-----FIN-----%%%%%%%%%%%%%%%%%%%%%%%%%%%%%%%
$x' = f ( x; \, a)$, sur une variété de dimension finie $n$, d'un
groupe de Lie analytique local de dimension $r$. Lie n'insiste pas en
général sur l'hypothèse sous-entendue d'analyticité, mais il utilise à
la place le mot <<\,{\em continu}\,>>, afin de marquer le contraste
entre sa propre théorie et la théorie de Galois discrète des groupes
de substitutions des racines d'une équation algébrique\footnote{\,
%%%%%%%%%%%%%%%%%%%%%%%-------DEBUT--------%%%%%%%%%%%%%%%%%%%%%%%%%%%
{\em Voir}~\cite{ h2001} pour une excellente étude historique.
}. %%%%%%%%%%%%%%%%%%%%%%%%-----FIN-----%%%%%%%%%%%%%%%%%%%%%%%%%%%%%%%

Ce que nous appelons aujourd'hui un {\sl groupe de Lie local}, à
savoir un espace $\C^r$ (ou $\R^r$) localisé autour d'un certain point
$e$ et muni d'une <<\,loi de multiplication de groupe\,>> analytique
locale $( a, b) \longmapsto \text{\bf m} ( a, b) = a \cdot b$ et d'une
application <<\,élément inverse\,>> $a \mapsto \text{\bf i} ( a)$, est
appelé par Engel et Lie <<\,{\sl groupe des paramètres\footnote{\,
%%%%%%%%%%%%%%%%%%%%%%%-------DEBUT--------%%%%%%%%%%%%%%%%%%%%%%%%%%%
Les pages 401--429 du volume~I de la \deutschplain{Theorie der
Transformationsgruppen} (\cite{ enlie1888}) sont consacrées à leur
étude générale. }
%%%%%%%%%%%%%%%%%%%%%%%%-----FIN-----%%%%%%%%%%%%%%%%%%%%%%%%%%%%%%%
d'un groupe de transformations}\,>>.

\medskip

{\footnotesize

\noindent{\bf Possibilité
de préciser rigoureusement les hypothèses de lieu.} Sans remettre en
cause les trois principes de pensée qui seront appliqués constamment
par la suite, signalons que l'on peut aisément formuler des hypothèses
précises et rigoureuses quant aux domaines d'existence d'un {\em
groupe de Lie local}.

\CITATION{
{\scriptsize
Ici, on doit observer que l'on a fixé le comportement des fonctions
$f_i ( x; \, a)$ seulement à l'intérieur des régions
\deutsch{Bereich} $(x)$ et $(a)$.

Par conséquent, nous avons la permission de substituer l'expression
$x_\nu ' = f_\nu ( x, a)$ dans les équations $x_i'' = f_i ( x', b)$
seulement lorsque le système de valeurs
\deutsch{Werthsystem} $x_1', \dots, x_n'$ se trouve dans la région
$(x)$. C'est pourquoi nous sommes forcé d'ajouter, aux hypothèses
remplies jusqu'à maintenant par les régions $(x)$ et $(a)$, la
supposition suivante: il doit être possible d'indiquer, à l'intérieur
des régions $(x)$ et $(a)$, des sous-régions respectives $(\!( x )\!)$
et $(\!( a )\!)$ d'une nature telle que les $x_i'$ restent toujours
dans la région $(x)$ lorsque les $x_i$ se meuvent \deutsch{laufen}
arbitrairement dans $(\!( x )\!)$ et quand les $a_k$ se meuvent
arbitrairement dans $(\!( a )\!)$; nous exprimons cela brièvement de
la manière suivante: la région $x' = f\big( (\!( x )\!) \, (\!( a )\!)
\big)$ doit tomber \deutsch{hineinfallen} entièrement dans la 
région $(x)$.

D'après ces conventions \deutsch{Festsetzungen}, si on choisit $x_1,
\dots, x_n$ dans la région $(\!( x )\!)$ et $a_1, \dots, a_r$ dans la
région $(\!( a )\!)$, on peut réellement effectuer la substitution
$x_k' = f_k ( x, a)$ dans l'expression $f_i ( x_1', \dots, x_n', \,
b_1, \dots, b_n)$; c'est-à-dire, lorsque que $x_1^0, \dots, x_n^0$
désigne un système de valeurs arbitraire dans la région $(\!( x)\!)$,
l'expression:
\[
f_i\big( f_1(x,a),\dots,f_n(x,a),\,b_1,\dots,b_n
\big)
\]
peut être développée, dans le voisinage du système de valeurs $x_k^0$,
comme une série entière ordinaire en $x_1 - x_1^0, \dots, x_n -
x_n^0$; les coefficients de cette série entière sont des fonctions de
$a_1, \dots, a_r$, $b_1, \dots, b_r$ et ils se comportent
régulièrement \deutsch{verhalten sich regulär}, lorsque les $a_k$ sont
arbitraires dans $(\!( a )\!)$ et les $b_k$ arbitraires dans $(a)$.
\REFERENCE{\cite{enlie1888},~pp.~15--16.}}}

Reformulons ces conditions avec une
optique purement locale. Si $( x_1,
\dots, x_n) \in \C^n$ sont les coordonnées complexes considérées, on
choisira la norme:
\[
\vert x\vert
:=
\max_{1\leqslant i\leqslant n}\,\vert x_i\vert,
\]
où $\vert w \vert = 
\sqrt{ w \overline{ w}}$ désigne le module d'un nombre complexe $w \in
\C$. Pour des <<\,rayons\,>> variés $\rho >0$, on considère alors les
ouverts spécifiques centrés en l'origine, appelés aujourd'hui {\sl
polydisques}:
\[
\Delta_\rho^n
:=
\big\{x\in\C^n:\,
\vert x\vert<\rho
\big\}.
\]
Par ailleurs, on munit aussi l'espace des paramètres $(a_1, \dots,
a_r) \in \C^r$ de la norme similaire:
\[
\vert
a
\vert
:=
\max_{1\leqslant k\leqslant r}\,
\vert a_k\vert, 
\]
et pour des réels strictement positifs $\sigma > 0$ variés, 
on introduit les ouverts: 
\[
\square_\sigma^r
:=
\big\{
a\in\C^r:\,\vert a\vert<\sigma
\big\}.
\]
Voici alors comment l'on peut formuler les axiomes de groupe de Lie en
localisant toutes les considérations autour de l'élément identité. Eu
égard au caractère purement local qui est visé, les fonctions $f_i (x;
a)$ seront supposées définies lorsque $\vert x \vert < \rho_1$ et
lorsque $\vert a \vert < \sigma_1$, pour un certain $\rho_1 > 0$ et
pour un certain $\sigma_1 > 0$, tous deux <<\,petits\,>>. L'élément
identité $e \in \C^r$ correspondra à l'origine $0 \in \C^r$, 
c'est-à-dire que l'on a:
\[
f_i(x_1,\dots,x_n:0,\dots,0)
\equiv
x_i
\ \ \ \ \ \ \ \ \ \
{\scriptstyle{(i\,=\,1\,\cdots\,n)}}.
\]
Afin que la composition de deux transformations successives $x' = f (
x; a)$ et $x'' = f ( x' ; b)$ soit bien définie, nous rapetissons
$\rho_1$ en un certain $\rho_2$ suffisamment petit avec $0 < \rho_2 <
\rho_1$ et nous rapetissons $\sigma_1$ en un certain $\sigma_2$
suffisamment petit avec $0 < \sigma_2 <
\sigma_1$ de telle sorte que l'on ait:
\[
\vert
f(x;a)
\vert
<
\rho_1
\ \ \ \ \ \ \ \
\text{\rm pour tout}
\ \ \ \ \
\vert x\vert<\rho_2\ \ \
\text{\rm et tout}\ \ \
\vert a\vert<\sigma_2,
\]
ce qui est possible par continuité grâce à la condition $f(x;0) \equiv 
x$.

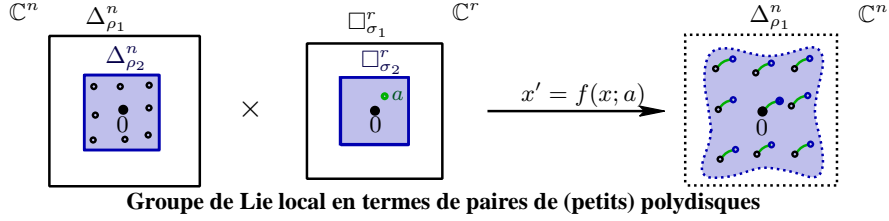
\medskip

\begin{center}
\input shrink-polydiscs.pstex_t
\end{center}

\smallskip\noindent
{\footnotesize\sf\em Loi de composition}:
Pour tout $x \in \Delta_{\rho_2}^n$, et tous 
$a, b \in \square_{\sigma_2}^r$, une composition arbitraire:
\def\theequation{1}\begin{equation}
x''
=
f\big(f(x;a);b\big)
=
f(x;c)
=
f(x;\text{\bf m}(a,b))
\end{equation}
s'identifie toujours à un élément $f ( x; c)$ de la même famille, pour
un paramètre unique $c = \text{\bf m} ( a,b)$ donné par une
application analytique locale (multiplication de groupe):
\[
\text{\bf m}:\ \ \ \
\square_{\sigma_1}^n\times\square_{\sigma_1}^n
\longrightarrow
\C^r
\]
qui satisfait $\text{\bf m} \big( \square_{\sigma_2}^r \times
\square_{ \sigma_2}^r \big) \subset \square_{ \sigma_1}^r$ et
$\text{\bf m} ( a, e) \equiv \text{\bf m} ( e, a) \equiv a$.

Si $a, b, c \in \square_{ \sigma_3}^n$ avec $0 < \sigma_3 < \sigma_2 <
\sigma_1$ suffisamment petit pour que trois compositions successives
soient bien définies, l'associativité de la composition entre
difféomorphismes garantit alors que:
\[
f\big(x;\text{\bf m}(\text{\bf m}(a,b),c)\big)
=
f\big(f(f(x;a);b);c\big)
=
f
\big(
x;
\text{\bf m}(a,\text{\bf m}(b,c))\big),
\]
d'où découle, grâce à l'unicité postulée de $c = \text{\bf m} ( a,
b)$, l'associativité de la loi de groupe: $\text{\bf m} \big(
\text{\bf m} ( a,b ), c\big) = \text{\bf m} \big( a, \text{\bf m} ( b,
c)\big)$ pour de telles valeurs (petites) de $a, b, c$.

\smallskip\noindent
{\footnotesize\sf\em Existence d'éléments inverses:}
Il existe une application analytique locale: 
\[
\text{\bf i}:\ \ \ 
\square_{\sigma_1}^r\longrightarrow\C^r
\]
satisfaisant $\text{\bf i} ( e) = e$, c'est-à-dire $\text{\bf i} ( 0)
= 0$, telle que pour tout $a \in \square_{ \sigma_2}^r$:
\[
\aligned
&
e
=
\text{\bf m}(a,\text{\bf i}(a)\big)
=
\text{\bf m}\big(\text{\bf i}(a),a\big)
\\
&
\text{\rm d'où de surcroît:}
\ \ \ \ \ 
x
=
f\big(f(x;a);\text{\bf i}(a)\big)
=
f\big(f(x;\text{\bf i}(a));a),
\endaligned
\]
pour tout $x \in \Delta_{\rho_2}^n$. Comme convenu, nous passerons
sous silence dans la suite la mention rigoureuse (qui serait
inélégante et lourde) des ouverts d'existence.

}

\HEAD{Chapitre~3.\,\,\,\,Théorèmes fondamentaux sur les groupes
de transformations}{ 3.3.\,\,\,Principe de raison suffisante et axiome
d'inverse}

\medskip\noindent{\bf 3.3.~Principe de raison
suffisante et axiome d'inverse.}
\label{raison-suffisante-axiome-inverse}
Un difféomorphisme analytique local quelconque d'une variété de
dimension $n$ peut être considéré comme effectuant une {\em
permutation} (différentiable) entre tous les points considérés, car en
particulier, un difféomorphisme, c'est une bijection. Ainsi, bien que
les difféomorphismes agissent sur un ensemble de cardinal infini (non
dénombrable), ils sont les {\em analogues continus} des {\em
permutations discontinues} d'un ensemble fini. En fait, dans les
années 1873 à 1880, l'idée fixe de Lie était d'ériger, dans le domaine
des {\em continua} $n$-dimensionnels, une théorie qui corresponde à la
théorie de Galois des substitutions des racines d'une équation
algébrique et qui lui soit en tout point analogue.

Comme dans les paragraphes qui précèdent, soit donc:
\[
x ' 
= 
f(x;\,a_1,\dots,a_r)
=: 
f_a(x)
\]
une famille de difféomorphismes locaux analytiques paramétrée par un
nombre fini $r$ de paramètres essentiels. Pour Lie, le seul axiome de
groupe vraiment significatif est celui qui demande qu'une telle
famille soit {\em fermée par composition}, à savoir que l'on a
toujours: 
\[
f_a\big(f_b(x)\big) 
\equiv 
f_c(x),
\]
pour un certain $c$ qui dépend de $a$ et de $b$, toujours avec la
restriction de localité assurant qu'une telle composition ait un
sens. En s'inspirant de sa connaissance du {\em Traité des
substitutions} de Jordan, Lie s'est demandé s'il était possible d'{\em
économiser} les deux autres axiomes standard de la structure de
groupe: l'axiome d'existence d'un élément identité, et l'axiome
d'existence d'un inverse pour tout élément du groupe.

\smallskip\noindent{\bf Assertion.}
(\cite{ jord1870})
{\em Soit $H$ un sous-ensemble quelconque d'un groupe abstrait $G$
dont le cardinal est {\em fini:} ${\rm Card}\, H < \infty$, et qui est
{\em fermé} par composition:
\[
h_1 h_2\in H 
\ \ \ \ \ \ \ \ \ \ \ \ 
\text{\rm toutes les fois que} 
\ \ \ \ \ \ \ 
h_1,\,h_2\in H. 
\]
Alors $H$ contient l'élément identité $e$ de $G$ et tout élément $h
\in H$ possède un inverse dans $H$, de telle sorte que $H$ lui-même
est un vrai sous-groupe de $G$.
}\medskip

{\em Preuve.}
En effet, soit $h \in H$ arbitraire. La suite infinie $h, h^2, h^3,
\dots, h^k, \dots$ d'éléments de l'ensemble fini $H$ doit
nécessairement devenir périodique: $h^a = h^{ a + n}$ pour un certain
$a \geqslant 1$ et pour un certain $n \geqslant 1$, d'où $e = h^n$,
donc $e \in H$ et $h^{ n - 1}$ est l'inverse de $h$.
\qed\smallskip

Dans ses travaux pionniers des années 1873 à 1880 et aussi dans la
{\em Theorie der Transformationsgruppen} qu'Engel a rédigée sous sa
direction entre 1884 et 1893, Lie est parvenu à transférer {\em tous}
les concepts de la théorie des groupes de substitutions
(Galois, Serret, Jordan) du discontinu
vers le continu: {\em loi de groupe, actions de groupe, sous-groupes,
sous-groupes normaux, groupe quotient, classification à conjugaison
près, groupe adjoint, représentation adjointe, formes normales,
(in)transitivité, (im)primitivité, prolongement holoédrique,
prolongement mériédrique, asystaticité},
\etc Le principe de raison suffisante
suggère alors naturellement qu'au fondement même de la théorie
générale, l'{\em élimination des deux axiomes} concernant l'élément
identité et l'existence d'inverses soit {\em aussi possible} dans
l'univers des groupes continus finis de transformations. Pendant plus
de dix années, Lie a en effet été convaincu qu'une assertion purement
similaire à celle énoncée ci-dessus devait être vraie dans le domaine
du continu, avec $G = {\sf Diff }_n$ le pseudo-groupe (infini,
continu) des difféomorphismes locaux de $\C^n$ (ou de $\R^n$) et en
prenant pour $H
\subset {\sf Diff }_n$ une collection finiment paramétrée d'équations 
de transformations $x' = f ( x; \, a)$ qui est stable par
composition\footnote{\,
%%%%%%%%%%%%%%%%%%%%%%%-------DEBUT--------%%%%%%%%%%%%%%%%%%%%%%%%%%%
\label{ferme-par-composition}
Sans présupposer l'existence d'un élément identité, la condition de
stabilité par composition qu'Engel et Lie ont dégagée (\cite{
enlie1888}, p.~14) s'énonce comme suit, lorsque $x \in \C^n$ et $a \in
\C^r$ sont à valeurs complexes. Il existe deux paires de domaines
$\mathcal{ X}^1 \subset \mathcal{ X}
\subset \C^n$ et $\mathcal{ A}^1
\subset \mathcal{ A} \subset \C^r$ 
dans l'espace des variables $x$ et dans l'espace des paramètres $a$
tels que pour tout $a
\in \mathcal{ A}$, l'application $x \mapsto f ( x; \, a)$ est un
difféomorphisme de $\mathcal{ X}$ sur son image, tels que de plus,
pour tout $a^1 \in \mathcal{ A}^1$ fixé, $f( \mathcal{ X}^1
\times\{ a^1 \} )
\subset \mathcal{ X}$, de telle sorte que pour tout 
$a^1 \in \mathcal{ A}^1$ et pour tout $b \in \mathcal{ A}$,
l'application composée $x \longmapsto f\big( f(x; \, a^1); \, b \big)$
est bien définie et établit un difféomorphe sur son image. De plus,
deux conditions significatives sont supposées.

$\bullet$\,\, Il existe une application analytique $\varphi = \varphi
( a, b)$ à valeurs dans $\C^m$ définie dans $\mathcal{ A} \times
\mathcal{ A}$ avec $\varphi ( \mathcal{ A}^1
\times \mathcal{ A}^1 ) \subset \mathcal{ A}$ telle que:

\centerline{$f \big( f(x; \,a); \,b \big)
\equiv f\big(x;\,\varphi (a,b) \big)
\ \ \ \ \ \ \ \ \ \ 
\text{\rm pour tous}\ \
x \in \mathcal{ X}^1,\ a\in \mathcal{ A}^1,\ b\in \mathcal{ A}^1$.}

$\bullet$\,\,
Il existe une application analytique $\chi = \chi ( a, c)$ à valeurs
dans $\C^m$ définie pour $a \in \mathcal{ A}$ et $c \in \mathcal{ A}$
avec $\chi \big( \mathcal{ A}^1 \times \mathcal{ A}^1 \big) 
\subset \mathcal{ A}$ qui résout $b$ en termes de $(a, c)$ 
dans l'équation $c = \varphi ( a, b)$, à savoir:

\centerline{$c \equiv \varphi \big( a,\,\chi (a, c) \big) \ \ \ \ \ \
\ \ \ \ \text{\rm pour tous}\ \ a \in \mathcal{ A}^1,\ c \in \mathcal{
A}^1$.}

\noindent
C'est avec ces hypothèses précises que Lie pensait déduire:
1) la présence
d'un élément identité $e \in \mathcal{ A}$; et: 2) l'existence d'un
difféomorphisme analytique local $a \mapsto \iota ( a)$ défini près de
$e$ et fixant $e$, tel que $f ( x; \, \iota ( a))$ est la
transformation inverse de $f ( x; \, a)$.
}. %%%%%%%%%%%%%%%%%%%%%%%%-----FIN-----%%%%%%%%%%%%%%%%%%%%%%%%%%%%%%%

Citons un extrait du premier mémoire systématique de Lie, paru en 1880
aux {\em Mathematische Annalen}. 

\CITATION{
Comme on le sait, on montre dans la théorie des substitutions que les
permutations d'un groupe de substitutions peuvent être ordonnées en
paires de permutations inverses l'une de l'autre. Mais comme la
différence entre un groupe de substitutions et un groupe de
transformations réside seulement dans le fait que le premier contient
un ensemble fini, et le second un ensemble infini d'opérations, il est
naturel de conjecturer que les transformations d'un groupe de
transformations puissent aussi être ordonnées par paires de
transformations inverses l'une de l'autre. Dans des travaux
antérieurs, je suis parvenu à la conclusion que tel devrait être le
cas. Mais comme au cours de mes investigations en question, certaines
hypothèses {\em implicites} se sont introduites au sujet des fonctions
qui apparaissent, je pense alors qu'il est nécessaire d'{\em ajouter
expressément l'exigence que les transformations du groupe puissent
être ordonnées par paires de transformations inverses l'une de
l'autre}. En tout cas, je conjecture que cette exigence est une
condition nécessaire de ma définition originale du concept
\deutsch{Begriff} de groupe de transformations. Toutefois, il m'a été
impossible de démontrer cela en général.
\REFERENCE{\cite{lie1880},~p.~444--445.}}

\smallskip
Cependant, en 1884, dans sa toute première année de travail de
rédaction en collaboration, Engel proposa un contre-exemple à cette
conjecture de Lie. Considérons en effet la famille d'équations de
transformations:
\[
x'
=
\zeta\,x,
\]
où $x, \, x' \in \C$ sont les coordonnées de l'espace-source et de
l'espace-image, et où le paramètre $\zeta \in \C$ est restreint à
$\vert \zeta \vert < 1$. \'Evidemment, cette famille est fermée par
composition, à savoir: lorsque $x' = \zeta_1 \, x$ et lorsque $x'' =
\zeta_2 \, x'$, la composition donne $x'' = \zeta_2 \, x' =
\zeta_2 \zeta_1 x$, et
elle appartient en effet à la famille en question, puisque $\vert
\zeta_2 \,\zeta_1 \vert < 1$ découle de $\vert \zeta_1 \vert, \, \vert
\zeta_2 \vert < 1$. Par contre, ni l'élément identité, ni l'inverse
de toute transformation n'appartiennent à la famille ainsi définie.

Comme on l'aura noté, cette proposition de contre-exemple n'est en
fait pas réellement convaincante. En effet, la condition $\vert \zeta
\vert < 1$ est ici visiblement 
artificielle, puisque la famille se prolonge en fait trivialement
comme groupe complet des dilatations $\big( x' =
\zeta\, x \big)_{ \zeta \in \C}$ de la droite complexe. L'idée de
Engel était d'en appeler à une application holomorphe univalente
$\omega : \Delta \to \C$ 
du disque unité $\Delta = \{ \zeta \in \C \colon
\vert \zeta \vert < 1\}$
à valeurs dans $\C$ qui possède le cercle unité $\{
\vert \zeta \vert = 1 \}$ comme coupure, à savoir: 
$\omega$ ne peut être prolongée holomorphiquement au-delà d'aucun
point $\zeta_0 \in \partial \Delta = \{
\vert \zeta \vert = 1 \}$ du bord
du disque unité. L'application $\omega$ utilisée (sans référence
bibliographique) par Engel est la suivante (\cite{ enlie1888},
p.~164). Si $k\geqslant 1$ est un entier quelconque, soit ${\sf
od}_k$ le nombre de diviseurs impairs de $k$, y compris $1$ et $k$
(lorsque $k$ est impair).

\smallskip\noindent{\bf Assertion.}
{\em La série infinie:
\[
\omega(\zeta)
:=
\sum_{\nu\geqslant 1}\,
\frac{\zeta^\nu}{1-\zeta^{2\nu}}
=
\sum_{\nu\geqslant 1}\,
\big(
\zeta^\nu+\zeta^{3\nu}+\zeta^{5\nu}+\zeta^{7\nu}+\cdots
\big)
=
\sum_{k\geqslant 1}\,
{\sf od}_k\,\zeta^k
\]
converge absolument dans tout disque $\Delta_\rho = \{ \zeta
\in \C: \,
\, \vert \zeta \vert < \rho \}$ de rayon $\rho < 1$ et y
définit une application holomorphe univalente $\Delta \to \C$ du
disque unité $\Delta = \{ \vert \zeta \vert < 1 \}$ à valeurs dans $\C$
qui {\em ne} se prolonge holomorphiquement au-delà d'{\em aucun} point
$\zeta_0 \in \partial \Delta := \{ \vert \zeta \vert = 1\}$. 
}\medskip

\label{engel-counterexample}
En fait, à la place de cette série explicite, on pourrait utiliser aussi
n'importe quelle autre application de type uniformisante de
Riemann\footnote{\,
%%%%%%%%%%%%%%%%%%%%%%%-------DEBUT--------%%%%%%%%%%%%%%%%%%%%%%%%%%%
Le théorème de Riemann énonce que pour tout ouvert connexe et
simplement connexe $\Lambda$ de $\C$ dont le complémentaire $\C
\backslash
\Lambda$ contient au moins un point, il existe une
application biholomorphe\,\,---\,\,c'est-à-dire holomorphe, bijective
et d'inverse holomorphe\,\,---\,\,$\chi :
\Lambda \to \Delta$ de $\Lambda$ sur le
disque unité $\Delta \subset \C$, appelée {\sl uniformisante de
Riemann}. Le théorème réflexion de Schwarz stipule qu'en tout point
$\lambda_0 \in \partial \Lambda$ du bord de $\Lambda$ autour duquel
$\partial \Lambda$ est un petit morceau de courbe analytique réelle,
l'uniformisante $\chi$ se prolonge holomorphiquement dans un voisinage
de $\lambda_0$. Réciproquement, si $\chi$ se prolonge
holomorphiquement au-delà d'un point du bord de $\Lambda$, ce bord est
alors nécessairement analytique réel au voisinage de ce point.
} %%%%%%%%%%%%%%%%%%%%%%%%-----FIN-----%%%%%%%%%%%%%%%%%%%%%%%%%%%%%%%
$\zeta \longmapsto \omega ( \zeta) =: \lambda$ qui établit un
biholomorphisme du disque unité $\Delta$ sur un domaine simplement
connexe $\Lambda := \omega ( \Delta)$ dont le bord est une courbe de
Jordan nulle part localement analytique réelle\footnote{\,
%%%%%%%%%%%%%%%%%%%%%%%-------DEBUT--------%%%%%%%%%%%%%%%%%%%%%%%%%%%
On peut même considérer un domaine dont le bord est une courbe de
Jordan continue nulle part différentiable, voire dont le bord est un
ensemble fractal, comme par exemple le flocon de Von Koch; on trouvera
une présentation concise et lumineuse de la théorie de Carathéodoy sur
le prolongement au bord des uniformisantes de Riemann dans le
Chapitre~17 de~\cite{ mi2006}.
}. %%%%%%%%%%%%%%%%%%%%%%%%-----FIN-----%%%%%%%%%%%%%%%%%%%%%%%%%%%%%%% 
Désignons alors par $\lambda \longmapsto \chi ( \lambda) =:
\zeta$ l'inverse d'une telle
application, et considérons la famille d'équations de transformations:
\[
\big(
x'
=
\chi(\lambda)\,x
\big)_{\lambda\in\Lambda}.
\]
Par construction, $\vert \chi ( \lambda) \vert < 1$ pour tout $\lambda
\in \Lambda$. Toute composition de $x' = \chi ( \lambda_1) \, x$ et de
$x'' = \chi ( \lambda_2) \, x'$ est de la forme $x'' = \chi ( \lambda)
\, x$, avec le paramètre défini
de manière unique $\lambda := \omega \big(
\chi ( \lambda_1)\, \chi ( \lambda_2) \big)$, donc
l'axiome de composition de groupe est satisfait. Cependant, il n'y a à
nouveau pas d'élément identité, et à nouveau, aucune transformation ne
possède un inverse. Et de plus crucialement (et pour terminer), il
n'existe aucun prolongement de cette famille à un domaine plus grand
$\widetilde{ \Lambda} \supset \Lambda$ qui soit accompagné d'un
prolongement holomorphe $\widetilde{
\chi}$ de $\chi$ à $\widetilde{ \Lambda}$ 
de telle sorte que $\widetilde{ \chi}
\big( \widetilde{ \Lambda} \big)$
contienne un voisinage de $\{ 1 \}$ (afin d'atteindre l'élément
identité), et même \emphasis{a fortiori}, il n'existe aucun
prolongement holomorphe tel de $\chi$ à un voisinage de $\overline{
\Delta}$ (afin d'obtenir les inverses des transformations $x' = \chi(
\lambda) \, x$ avec $\lambda \in \Lambda$ proche de $\partial
\Lambda$).

Dans le volume~I de la \emphasis{Theorie der Transformationsgruppen},
cet exemple apparaît seulement au Chapitre~9, pp.~163--165, et il est
écrit en petits caractères. Constation frappante: pour des raisons de
pureté et de systématicité dans la pensée, la structure des neuf
premiers chapitres est organisée afin de faire
autant que possible l'économie de l'existence d'un élément identité et
de transformations inverses l'une de l'autre par paires. Ce choix
théorique complexifie considérablement la présentation de la théorie
fondamentale, pourtant censée être relativement facile d'accès
afin de toucher un public assez large de mathématiciens. Mais il est
bien connu que la meilleure qualité dans l'exposition des premiers
éléments d'un ouvrage et que les meilleurs choix stratégiques quant à
son organisation thématique, ne peuvent être atteints, paradoxalement,
qu'à la fin du processus de mise en forme dans son
ensemble. Impossible, donc, de demander au jeune et novice Friedrich
Engel (alors âgé de 24 ans) de faire prévaloir un point de vue
d'accessibilité et de simplicité dans la présentation. Au contraire,
l'objectif affirmé de son maître Sophus Lie était d'atteindre la plus
grande généralité, de s'élever le plus haut possible dans un tel
traité.

Et malgré le contre-exemple du jeune Engel que nous venons de
détailler, Lie était quand même persuadé que l'analogie profonde de sa
théorie des groupes continus avec la théorie des groupes finis de
substitutions n'était pas, dans sa racine métaphysique profonde,
réellement remise en cause. Ainsi, toutes les fois que cela est
possible du point de vue abstrait de Lie, les énoncés des neuf
premiers chapitres de la {\em Theorie der Transformationsgruppen}
n'utilisent ni l'existence d'un élément identité, ni l'existence de
transformations inverses l'une de l'autre par paires. Engel et Lie
étudient seulement les familles continues finies d'équations de
transformations $x_i' = f_i ( x; \, a_1, \dots, a_r)$, $i = 1, \dots,
n$ qui sont fermées par composition au sens de la note
p.~\pageref{ferme-par-composition}, sans hypothèse supplémentaire:
grand degré d'abstraction et de généralité. En fait, à partir de cette
seule condition de fermeture par composition, Engel et Lie déduisent
que les équations finies $x_i' = f_i ( x; \, a)$ satisfont certaines
{\em équations différentielles fondamentales}, stipulées dans le
Théorème p.~\pageref{Theoreme-3} ci-dessous (Théorème~3 p.~33
dans~\cite{ enlie1888}). C'est alors l'existence de telles équations
différentielles qu'ils prennent systématiquement comme hypothèse
principale, à la place de la fermeture par composition, 
toujours sans élément identité et sans transformations
inverses. 

Le Théorème~26
est énoncé à la page~163 de~\cite{ enlie1888}
({\em voir} p.~\pageref{Theoreme-26}), 
et il est
démontré tout juste avant que n'apparaisse le contre-exemple de Engel.
Ce théorème raffiné et subtil confirme la croyance métaphysique de
Lie, montre sa persévérance intellectuelle, et prouve à nouveau sa
remarquable force de conceptualisation. En résumé, il s'énonce comme
suit
\footnote{\,
%%%%%%%%%%%%%%%%%%%%%%%-------DEBUT--------%%%%%%%%%%%%%%%%%%%%%%%%%%%
En première approche, ce passage doit être lu en admettant quelques
notions qui ne seront présentées que dans les paragraphes qui suivent.
}. %%%%%%%%%%%%%%%%%%%%%%%%-----FIN-----%%%%%%%%%%%%%%%%%%%%%%%%%%%%%%%

Soit $x_i' = f_i ( x; \, a_1, \dots, a_r)$, $i = 1, \dots, n$ une
collection de transformations fermée par compositions locales.
D'après le théorème énoncé p.~\pageref{Theoreme-3} ci-dessous, il
existe un système d'équations différentielles fondamentales de la
forme:
\[
\frac{\partial x_i'}{\partial a_k}
=
\sum_{j=1}^r\,\psi_{kj}(a)\cdot\xi_{ji}(x')
\ \ \ \ \ \ \ \ \ \ \ \ \
{\scriptstyle{(i\,=\,1\,\cdots\,n\,;\,k\,=\,1\,\cdots\,r)}}.
\]
qui est satisfait identiquement par les fonctions $f_i ( x, a)$, 
où les $\psi_{ kj}$ sont
certaines fonctions analytiques des paramètres
$(a_1, \dots, a_r)$. Si
l'on introduit alors les $r$ transformations infinitésimales (\voir~le
\S~3.4 ci-dessous) qui sont définies par:
\[
\sum_{i=1}^n\,\xi_{ki}(x)\,\frac{\partial f}{\partial x_i}
=:
X_k(f)
\ \ \ \ \ \ \ \ \ \ \ \ \
{\scriptstyle{(k\,=\,1\,\cdots\,r)}},
\]
et si l'on forme les équations finies: 
\[
\aligned
x_i'
&
=
\exp\big(\lambda_1X_1+\cdots+\lambda_kX_k\big)(x)
\\
&
=:
g_i(x;\,\lambda_1,\dots,\lambda_r)
\ \ \ \ \ \ \ \ \ \ \ \ \
{\scriptstyle{(i\,=\,1\,\cdots\,n)}}
\endaligned
\]
du groupe à $r$ paramètres qui est engendré par ces $r$
transformations infinitésimales, alors ce groupe contient l'élément
identité $g( x; \, 0)$ et ses transformations sont ordonnées par
paires inverses l'une de l'autre: $g(x; \, - \lambda) = g (x; \,
\lambda)^{ -1}$. Enfin, 
le Théorème~26 
en question
(p.~\pageref{Theoreme-26} ci-dessous) 
énonce que dans ces équations finies $x_i'
= g_i ( x; \, \lambda)$, il est possible d'introduire de {\em
nouveaux} paramètres locaux $\overline{ a}_1,
\dots, \overline{ a}_r$ à la place de $\lambda_1, \dots, \lambda_r$ 
de telle sorte que les équations de transformations qui en résultent:
\[
\aligned
x_i'
&
=
g_i\big(x;\,\lambda_1(\overline{a}),\dots,\lambda_r(\overline{a})\big)
\\
&
=:
\overline{f}_i
(x_1,\dots,x_n,\,\overline{a}_1,\dots,\overline{a}_r)
\ \ \ \ \ \ \ \ \ \ \ \ \
{\scriptstyle{(i\,=\,1\,\cdots\,n)}}
\endaligned
\]
représentent une famille de $\infty^r$ transformations qui 
embrasse, {\em après un prolongement analytique
éventuel}, toutes les $\infty^r$ transformations initiales:
\[
x_i'
=
f_i(x_1,\dots,x_n,\,a_1,\dots,a_r).
\]
Ainsi Lie réalise-t-il et confirme-t-il son idée d'épuration
axiomatique. En répondant à Engel que l'on doit s'autoriser à changer
éventuellement de paramètres\footnote{\,
%%%%%%%%%%%%%%%%%%%%%%%-------DEBUT--------%%%%%%%%%%%%%%%%%%%%%%%%%%%
En changeant de paramètre dans la famille $x' = \chi ( \lambda) \, x$,
on retrouve évidemment le groupe complet des dilatations $x' = \zeta
\, x$.
}, %%%%%%%%%%%%%%%%%%%%%%%%-----FIN-----%%%%%%%%%%%%%%%%%%%%%%%%%%%%%%%
on peut toujours élargir le domaine initial d'existence pour capturer
l'identité et les transformations inverses. \`A la
contre-exemplification contrariante répond donc la dialectique
complexifiante de vérités supérieures qui doivent
être quêtées sans relâche.

Afin de poursuivre et de rendre plus compréhensible cet énoncé, il
nous faut à présent exposer: 1) les transformations infinitésimales
(\S~3.4); 2) les équations différentielles fondamentales (\S~3.5); 3)
la reconstitution des équations finies du groupe à partir d'une
collection de transformations infinitésimales indépendantes (\S~3.6);
4) le théorème de Clebsch-Lie-Frobenius sur les systèmes complets et
indépendants d'équations aux dérivées partielles (\S~3.7).

\HEAD{Chapitre~3.\,\,\,\,Théorèmes fondamentaux sur les groupes
de transformations}{ 3.4.\,\,\,Introduction des transformations
infinitésimales}

\medskip\noindent{\bf 3.4.~Introduction 
des transformations infinitésimales.}
\label{introduction-transformations-infinitesimales}
Soit $\varepsilon$ une quantité infiniment petite au sens de Leibniz,
ou une quantité arbitrairement petite soumise au formalisme rigoureux
de Weierstrass. Pour chaque $k \in \{ 1, 2, \dots, r
\}$ fixé à l'avance, considérons le point:
\[
\label{introduce-infinitesimal-transformations}
\aligned
x_i'
&
=
f_i
\big(x;\,e_1,\dots,e_k+\varepsilon,\dots,e_n\big)
\\
&
=
x_i
+
\frac{\partial f_i}{\partial a_k}(x;e)\cdot
\varepsilon
+\cdots
\ \ \ \ \ \ \ \ \ \
{\scriptstyle{(i\,=\,1\,\cdots\,n)}}
\endaligned
\]
qui est déplacé infinitésimalement en partant du point initial $x = f
( x; e)$ en ajoutant l'incrément infime $\varepsilon$ seulement à la
$k$-ième coordonnée $e_k$ du paramètre identité $e$. Grâce à un
développement de Taylor à l'ordre $1$, on peut interpréter ce
mouvement spatial infinitésimal en introduisant,
pour tout $k \in \{ 1, 2, \dots, r\}$, le champ de vecteurs
(ainsi qu'une notation raccourcie appropriée pour 
désigner ses coefficients):
\[
\aligned
X_k^e
:=
\sum_{i=1}^n\,
\frac{\partial f_i}{\partial a_k}(x;e)\,
\frac{\partial}{\partial x_i}
&
:=
\sum_{i=1}^n\,
\xi_{ki}(x)\,
\frac{\partial}{\partial x_i},
\endaligned
\]
écrit ici sous la forme d'un opérateur de dérivation d'ordre un. On
peut aussi le considérer comme une famille de vecteurs colonnes:
\[
{}^\tau \, \big( 
{\textstyle{ \frac{ \partial f_1 }{\partial a_k}}}, \cdots,
{\textstyle{ 
\frac{ \partial f_n }{\partial a_k}}} \big) \Big\vert_x = 
{}^\tau \big( \xi_{ k1}, \dots, \xi_{ kn} \big)
\Big\vert_x 
\] 
basés aux points $x$, où ${}^\tau ( \cdot)$ désigne l'opérateur de
transposition des matrices, qui transforme bien sûr les lignes
horizontales en lignes verticales. Alors on peut réécrire ce
mouvement infime sous la forme $x ' = x + \varepsilon\, X_k^e +
\cdots$, ou bien encore, de manière équivalente:
\[
x_i'
=
x_i
+
\varepsilon\,\xi_{ki}(x)
+\cdots
\ \ \ \ \ \ \ \ \ \ \ \ \ \
{\scriptstyle{(i\,=\,1\,\cdots\,n)}},
\]
où les termes supprimés ``$+ \cdots$'' sont bien entendu des ${\rm O}
( \varepsilon^2)$ uniformes par rapport à $x$, de telle sorte que d'un
point de vue géométrique, $x '$ est <<\,poussé\,>> infinitésimalement
à partir de $x$ d'une longueur $\varepsilon$ le long du vecteur $X_k^e
\big\vert_x$, comme l'illustre la partie
gauche de notre figure (avec $k = 1$).

\begin{center}
\input push-points.pstex_t
\end{center}
\label{push-points}

Plus généralement, toujours à partir du paramètre
identité $e$, on peut ajouter à $e$ un 
incrément infinitésimal arbitraire: 
\[
\big(
e_1+\varepsilon\,\lambda_1,\dots,
e_k+\varepsilon\,\lambda_k,\dots,
e_r+\varepsilon\,\lambda_r
\big),
\]
où ${}^\tau(\lambda_1, \dots, \lambda_r) \big\vert_e$ est un vecteur
tangent constant basé en $e$ dans l'espace des paramètres. Alors grâce
à la linéarité de l'application tangente, c'est-à-dire grâce à la
règle de dérivation composée en coordonnées, il en découle que:
\[
\aligned
f_i(x;\,e+\varepsilon\,\lambda)
&
=
x_i
+
\sum_{k=1}^n\,\varepsilon\,\lambda_k\cdot
\frac{\partial f_i}{\partial a_k}(x;e)
+\cdots
\\
&
=
x_i
+
\varepsilon\,
\sum_{k=1}^n\,\lambda_k\cdot\xi_{ki}(x)
+\cdots,
\endaligned
\]
de telle sorte que tous les points $x ' = x + \varepsilon \, X +
\cdots$ sont simultanément et infinitésimalement déplacés le long du
champ de vecteurs:
\[
X
:=
\lambda_1\,X_1^e
+\cdots+
\lambda_r\,X_r^e
\]
qui est la combinaison linéaire générale des $r$ précédents champs de
vecteurs basiques $X_k^e$, $k=1, \dots, r$.

Occasionnellement, Engel et Lie écrivent qu'un tel champ de vecteurs
$X$ {\sl appartient au groupe} $x ' = f ( x; a)$, pour
signifier que $X$ vient accompagné du mouvement infinitésimal $x ' = x
+ \varepsilon \, X$ qu'il est supposé exécuter (les points de
suspension sont censés être supprimés dans l'intuition). En
accord avec cet acte de pensée synthétique, Lie appelle $X$\, une {\sl
transformation infinitésimale}, considérant en effet que $x ' = x +
\varepsilon \, X$ est juste un cas particulier de 
$x ' = f ( x, a)$. Une autre raison fondamentale et très profonde,
pour laquelle Lie dit que $X$ {\sl appartient} au groupe $x ' = f (
x;\, a)$ est qu'il a démontré que les groupes continus finis de
transformations locales sont en correspondance biunivoque avec les
espaces vectoriels purement linéaires:
\[
\text{\rm Vect}_\C
\big(
X_1,X_2,\dots,X_r
\big),
\]
de transformations infinitésimales, qui héritent en fait aussi,
directement à partir de la loi de multiplication de groupe, d'une
structure {\em algébrique} additionnelle, comme nous allons maintenant
le rappeler. Quelques préparatifs (\S\S~3.5, 3.6 et 3.7) sont
nécessaires.

Tout d'abord (\S~3.5), il faut {\em
infinitésimaliser}\,\,---\,\,c'est-à-dire différentier\,\,---\,\,la
loi de composition de groupe, réorganiser les équations obtenues, et
interpréter géométriquement leur signification.

Ensuite (\S~3.6), il faut traduire dans la théorie des
groupes les théorèmes classiques sur l'intégration des systèmes
d'équations différentielles ordinaires d'ordre un. Ces systèmes
correspondent à l'intégration d'un, et d'un seul
champ de vecteurs. Une dernière
pièce préliminaire 
est donc nécessaire (\S~3.7): le théorème de Clebsch-Frobenius, 
qui correspond à l'intégration de {\em plusieurs}
champs de vecteurs soumis à une certaine
condition de compatibilité. 

\HEAD{Chapitre~3.\,\,\,\,Théorèmes fondamentaux sur les groupes
de transformations}{ 3.5.\,\,\,\'Equations différentielles fondamentales}

\medskip\noindent{\bf 3.5.~\'Equations différentielles fondamentales.}
\label{equations-differentielles-fondamentales}
Partons de la loi de composition de groupe, que
nous réécrivons comme suit: 
\[
\label{beginning-computation}
x''
=
f\big(f(x;a);b\big)
=
f(x;a\cdot b)
=:
f(x;c).
\]
Ici, $c := a\cdot b$ dépend de $a$ et $b$, mais à la place de $a$ et
$b$, nous allons considérer $a$ et $c$ comme paramètres indépendants.
Ainsi, en remplaçant $b = a^{ - 1} \cdot c =: b(a,c)$, les équations:
\[
f_i\big(f(x;a);\,b(a,c)\big)
\equiv
f_i(x;c)
\ \ \ \ \ \ \ \ \ \
{\scriptstyle{(i\,=\,1\,\cdots\,n)}}
\]
sont satisfaites identiquement pour tout $x$, tout $a$ et tout
$c$. Ensuite, pour $k \in \{ 1, 2, \dots, r\}$ fixé, différentions ces
identités par rapport à $a_k$, en notant de manière abrégée 
$f_i' \equiv f_i( x'; b)$ et $x_j' \equiv f_j(x;a)$, ce
qui nous donne:
\[
\frac{\partial f_i'}{\partial x_1'}\,
\frac{\partial x_1'}{\partial a_k}
+\cdots+
\frac{\partial f_i'}{\partial x_n'}\,
\frac{\partial x_n'}{\partial a_k}
+
\frac{\partial f_i'}{\partial b_1}\,
\frac{\partial b_1}{\partial a_k}
+\cdots+
\frac{\partial f_i'}{\partial b_r}\,
\frac{\partial b_r}{\partial a_k}
\equiv
0
\ \ \ \ \ \ \ \ \ \
{\scriptstyle{(i\,=\,1\,\cdots\,n)}}.
\]
Ici bien sûr à nouveau, l'argument de $f_i'$ est $\big( f( x, a); b (
a,c) \big)$, l'argument des $x_l'$ est $( x; a)$ et l'argument 
des $b_j$ est
$(a, c)$. Grâce à $x'' ( x'; e) \equiv x'$, la matrice $\frac{\partial
f_i'}{\partial x_k'} \big( f(x;e); b(e, e) \big)$ est l'identité $I_{
n \times n}$. Par conséquent, en appliquant la règle de
Cramer\footnote{\,
%%%%%%%%%%%%%%%%%%%%%%%-------DEBUT--------%%%%%%%%%%%%%%%%%%%%%%%%%%%
---\,\,et en rapetissant aussi si 
nécessaire les domaines
d'existence, mais sans introduire de
dénomination spéciale\,\,---
}, %%%%%%%%%%%%%%%%%%%%%%%%-----FIN-----%%%%%%%%%%%%%%%%%%%%%%%%%%%%%%%
pour tout $k$ fixé, 
nous pouvons résoudre les $n$ 
équations linéaires précédentes par rapport
aux $n$ inconnues $\frac{ \partial x_1' }{ \partial a_k}, \dots, \frac{
\partial x_n' }{ \partial a_k}$, 
et nous obtenons des expressions de la forme:
\def\theequation{2}\begin{equation}
\aligned
\label{2-computation}
\frac{\partial x_\nu'}{\partial a_k}
(x;\,a)
&
=
\Xi_{1\nu}(x',b)\,\frac{\partial b_1}{\partial a_k}(a,c)
+\cdots+
\Xi_{r\nu}(x',b)\,\frac{\partial b_r}{\partial a_k}(a,c)
\\
&
\ \ \ \ \ \ \ \ \ \ \ \ \ \ \ \
{\scriptstyle{(\nu\,=\,1\,\cdots\,n\,;\,k\,=\,1\,\cdots\,r)}},
\endaligned
\end{equation}
pour certaines fonctions analytiques $\Xi_{ j\nu} ( x', b)$ qui sont
indépendantes de $k$.

D'autre part, afin de pouvoir substituer la dérivée partielle
$\frac{ \partial b_j}{ \partial a_k}$ par
une expression adéquate, différentions par rapport à $a_k$ les
équations suivantes, qui sont satisfaites identiquement:
\[
c_\mu
\equiv
\text{\bf m}_\mu\big(a,b(a,c)\big)
\ \ \ \ \ \ \ \ \ \
{\scriptstyle{(\mu\,=\,1\,\cdots\,r)}}.
\]
De cette manière-là, nous obtenons: 
\[
0
\equiv
\frac{\partial\text{\bf m}_\mu}{\partial a_k}
+
\sum_{\pi=1}^r\,
\frac{\partial\text{\bf m}_\mu}{\partial b_\pi}\,
\frac{\partial b_\pi}{\partial a_k}
\ \ \ \ \ \ \ \ \ \
{\scriptstyle{(\mu\,=\,1\,\cdots\,r)}}.
\]
Mais puisque la matrice $\frac{ \partial \text{\bf m}_\mu}{ \partial
b_\pi}$ se réduit à la matrice
identité $I_{ r \times r}$ pour\footnote{\,
%%%%%%%%%%%%%%%%%%%%%%%-------DEBUT--------%%%%%%%%%%%%%%%%%%%%%%%%%%%
---\,\,seulement
parce que $\text{\bf m} (e, b) \equiv b$\,\,---
} %%%%%%%%%%%%%%%%%%%%%%%%-----FIN-----%%%%%%%%%%%%%%%%%%%%%%%%%%%%%%%
$\big( a, b (a,
c) \big) \big\vert_{ (a, c) = (e,e)} = ( e, e)$, la règle de Cramer,
à nouveau, nous permet de résoudre ce système par rapport aux $r$
inconnues $\frac{ \partial b_\pi}{ \partial a_k}$, et nous obtenons
ainsi des expressions de la forme:
\[
\frac{\partial b_\pi}{\partial a_k}
(a,c)
=
\Psi_{k\pi}(a,c),
\]
avec certaines fonctions $\Psi_{ k\pi}$,
qui sont définies dans un domaine éventuellement
plus petit. En posant $(a, c) = ( e, e)$
dans le même système avant de le résoudre, nous
obtenons en fait (noter le signe <<\,moins\,>>): 
\[
\Big(
\frac{\partial b_\pi}{\partial a_k}(e,e)
\Big)_{1\leqslant k\leqslant r}^{1\leqslant\pi\leqslant r}
=
-I_{r\times r},
\]
d'où $\Psi_{ k \pi} ( e, e) = - \delta_k^\pi$, où le symbole
$\delta_k^\pi$ vaut $1$ si $k = \pi$ et $0$ sinon.

Nous pouvons donc maintenant insérer dans~\thetag{ 2} la valeur
obtenue $\Psi_{ k \pi}$ pour $\frac{ \partial b_\pi}{ \partial a_k}$,
ce qui nous donne des équations aux dérivées partielles qui sont
absolument cruciales dans toute la théorie de Lie:

\CITATION{
\def\theequation{2'}\begin{equation}
\frac{\partial x_\nu'}{\partial a_k}(x;a)
=
\sum_{\pi=1}^r\,\Psi_{k \pi}(a,b)\cdot\Xi_{\pi\nu}(x',b)
\ \ \ \ \ \ \ \ \ \
{\scriptstyle{(\nu\,=\,1\,\cdots\,n,\,\,\,k\,=\,1\,\cdots\,r)}}
\end{equation}
\emphasis{Ces équations sont
de la plus haute importance \deutsch{\"ausserst
wichtig}, comme nous allons
le voir plus tard.}
\REFERENCE{\cite{enlie1888},~p.~29.}}

Ici, nous avons remplacé $c$ par $c = c ( a, b) = a \cdot b$, d'où $b
\big( a, c (a, b) \big) \equiv b$,
et nous avons reconsidéré $(a, b)$ commes variables indépendantes.

\smallskip

Maintenant, en posant $b:=e$ dans ces équations~\thetag{ 2'}, les
dérivées partielles par rapport aux paramètres $a_k$ des
transformations $x_i' = x_i' ( x; a)$ s'expriment par des équations
différentielles:
\def\theequation{2''}\begin{equation}
\label{fundamental-differential-equations}
\boxed{
\boxed{\fboxrule=0.25pt
\frac{\partial x_i'}{\partial a_k}(x;a)
=
\sum_{j=1}^r\,\psi_{kj}(a)\cdot
\xi_{ji}\big(x'(x;a)\big)
}}
\ \ \ \ \ \ \ \ \ \ 
{\scriptstyle{(i\,=\,1\,\cdots\,n,\,\,\,k\,=\,1\,\cdots\,r)}},
\end{equation}
qui les font apparaître comme {\em combinaisons linéaires}, avec
certains coefficients $\psi_{ kj} (a) := \Psi_{ kj} ( a, e)$ {\em qui
dépendent seulement de $a$} des quantités $\xi_{ ji} ( x') := \Xi_{ji}
( x', b) \big\vert_{ b = e}$. Mais nous connaissons en fait déjà ces
fonctions $\xi_{ ji} (x')$. 

\smallskip

En effet, si nous posons $a = e$ dans ces équations, alors grâce à
$\psi_{ kj} ( e) = - \delta_k^j$, nous obtenons immédiatement:
\[
\xi_{ki}(x)
=
\xi_{ki}\big(x'(x;e)\big)
=
-
\frac{\partial x_i'}{\partial a_k}(x;e),
\]
d'où les $\xi_{ ki} ( x)$ coïncident, modulo un signe <<\,moins\,>>
uniforme, avec les coefficients des $r$ transformations infinitésimales
que nous avons déjà introduites
p.~\pageref{introduce-infinitesimal-transformations}, et que nous
avons considérées comme des opérateurs de dérivations
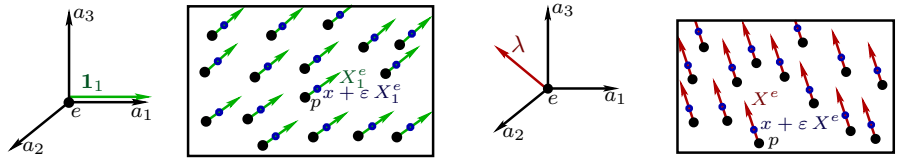
donnant la direction d'un mouvement infinitésimal de tous
les points de l'espace:
\def\theequation{3}\begin{equation}
\aligned
X_k^e\big\vert_x
&
=
\frac{\partial f_1}{\partial a_k}(x;e)\,
\frac{\partial}{\partial x_1}
+\cdots+
\frac{\partial f_n}{\partial a_k}(x;e)\,
\frac{\partial}{\partial x_n}
\\
&
=:
-\xi_{k1}(x)\,\frac{\partial}{\partial x_1}
-\cdots-
\xi_{kn}(x)\,\frac{\partial}{\partial x_n}
\ \ \ \ \ \ \ \ \ \ \ \ \ \ \ \ \ \ \ \
{\scriptstyle{(k\,=\,1\,\cdots\,r)}}.
\endaligned
\end{equation}
{\em Maintenant enfin}, après avoir reproduit à l'aveugle ces calculs
arides d'élimination différentielle dûs à 
Lie, nous pouvons dévoiler encore plus
précisément l'interprétation géométrique adéquate dont sont riches les
équations différentielles~\thetag{ 2''} doublement encadrées.

\begin{center}
\input geom-fund-diff-eqs.pstex_t
\end{center}

\noindent
Plutôt que de différentier par rapport à $a_k$ au point spécial $(x;
e)$, nous pouvons, par souci de généralité, différentier en un point
quelconque $(x;a)$, ce qui nous donne alors les champs de vecteurs:
\[
X_k^a\big\vert_{f_a(x)}
:=
\frac{\partial f_1}{\partial a_k}(x;a)\,
\frac{\partial}{\partial x_1}
+\cdots+
\frac{\partial f_n}{\partial a_k}(x;a)\,
\frac{\partial}{\partial x_n}
\ \ \ \ \ \ \ \ \ \ \ \ \ \ \ \ \ \ \ \
{\scriptstyle{(k\,=\,1\,\cdots\,r)}},
\]
et alors les équations différentielles fondamentales expriment que ces
$r$ nouvelles transformations infinitésimales:
\[
\boxed{
X_k^a\big\vert_{f_a(x)}
=
-\psi_{k1}(a)\,X_1^e\big\vert_{f_a(x)}
-\cdots-
\psi_{kr}(a)\,X_r^e\big\vert_{f_a(x)}}
\]
sont des combinaisons linéaires {\em à coefficients qui ne dépendent
que des paramètres de groupe $(a_1,
\dots, a_r)$}, des $r$ transformations infinitésimales
$X_1^e, \dots, X_r^e$ calculées au paramètre-identité $e$, {\em et
reconsidérées au point $f_a ( x)$ issu de $x$ qui est <<\,poussé\,>>
par la transformation $f_a ( \cdot)$}. Le diagramme illustre
géométriquement cette interprétation.

Enfin, puisque la matrice $\psi (a)$ possède un inverse $\widetilde{
\psi} ( a)$ qui est analytique dans un voisinage de $e$, les équations
différentielles fondamentales peuvent aussi être écrites sous la forme
réciproque, parfois utile:
\def\theequation{4}\begin{equation}
\xi_{ji}\big(x'(x;a)\big)
=
\sum_{k=1}^r\,\widetilde{\psi}_{jk}(a)\cdot
\frac{\partial x_i'}{\partial a_k}(x;a).
\ \ \ \ \ \ \ \ \ \ \ \ \
{\scriptstyle{(i\,=\,1\,\cdots\,n\,;\,\,\,j\,=\,1\,\cdots\,r)}}.
\end{equation}
Nous pouvons donc maintenant énoncer le tout premier théorème
fondamental de la théorie de Lie.

\smallskip\noindent{\bf Théorème.} 
\label{Theoreme-3}
(\cite{enlie1888}, pp.~33--34)
{\em Si les $n$ équations de transformations:
\[
x_i'
=
f_i(x_1,\dots,x_n,\,a_1,\dots,a_n)
\ \ \ \ \ \ \ \ \ \
{\scriptstyle{(i\,=\,1\,\cdots\,n)}}
\]
dont les paramètres $a_1, \dots, a_r$ sont tous essentiels
représentent un groupe continu fini local de
transformations\footnote{\,
%%%%%%%%%%%%%%%%%%%%%%%-------DEBUT--------%%%%%%%%%%%%%%%%%%%%%%%%%%%
\label{sans-inverse}
Ici, nous avons admis l'existence d'un élément identité et nous avons
supposé que les transformations peuvent être ordonnées par paires de
transformations inverses l'une de l'autre. Mais Engel et Lie déduisent
l'existence d'équations différentielles du type~\thetag{ 5} sous la
seule hypothèse que la famille $x' = f ( x; \, a)$ soit fermée par
composition, au sens explicité dans la note
p.~\pageref{ferme-par-composition}, {\em voir} aussi le \S~3.9
ci-dessous. },
%%%%%%%%%%%%%%%%%%%%%%%%-----FIN-----%%%%%%%%%%%%%%%%%%%%%%%%%%%%%%%
alors $x_1', \dots, x_n'$, considérées comme fonctions de 
$a_1, \dots, a_r, \, x_1, \dots, x_n$ satisfont certaines 
équations différentielles de la forme spécifique:
\def\theequation{5}\begin{equation}
\frac{\partial x_i'}{\partial a_k}
=
\sum_{j=1}^r\,\psi_{kj}(a_1,\dots,a_r)\cdot
\xi_{ji}(x_1',\dots,x_n')
\ \ \ \ \ \ \ \ \
{\scriptstyle{(i\,=\,1\,\cdots\,n;\,\,\,k\,=\,1\,\cdots\,r)}},
\end{equation}
où la matrice $\psi_{ kj} (a)$ de fonctions analytiques définie au
voisinage de l'identité $e$ satisfait $\psi_{ kj} ( e) = -
\delta_k^j$, et où les fonctions 
$\xi_{ ki} ( x)$, égales à $- \frac{ \partial x_i'}{ \partial a_k} (
x; \, e)$, coïncident, modulo un signe <<\,moins\,>> uniforme, avec
les coefficients des $r$ transformations infinitésimales basiques:
\[
X_k^e
\big\vert_{x}
=
\sum_{i=1}^n\,
\frac{\partial f_i}{\partial a_k}(x;\,e)\,
\frac{\partial}{\partial x_i}
\ \ \ \ \ \ \ \ \ \ \ \ \
{\scriptstyle{(k\,=\,1\,\cdots\,r)}}. 
\]
De plus, il est impossible de trouver $r$ quantités constantes $\lambda_1,
\dots, \lambda_r$ non toutes nulles telles que les $n$
expressions:
\[
\lambda_1\,\xi_{1i}(x')
+\cdots+
\lambda_r\,\xi_{ri}(x')
\ \ \ \ \ \ \ \ \ \
{\scriptstyle{(i\,=\,1\,\cdots\,n)}}
\]
s'annulent simultanément. 
}\medskip

Cette dernière propriété signifie évidemment que les $r$
transformations infinitésimales $X_1, \dots, X_r$ sont {\em
linéairement indépendantes}\footnote{\,
%%%%%%%%%%%%%%%%%%%%%%%-------DEBUT--------%%%%%%%%%%%%%%%%%%%%%%%%%%%
---\,\,sur $\C$ ou sur $\R$, 
suivant que $a \in \C^r$
ou $a \in \R^r$\,\,---
} %%%%%%%%%%%%%%%%%%%%%%%%-----FIN-----%%%%%%%%%%%%%%%%%%%%%%%%%%%%%%%
et elle découle du
fait que les paramètres sont tous essentiels\footnote{\,
%%%%%%%%%%%%%%%%%%%%%%%-------DEBUT--------%%%%%%%%%%%%%%%%%%%%%%%%%%%
Si les paramètres $(a_1, \dots, a_r)$ ne sont pas essentiels, d'après
le théorème p.~\pageref{theoreme-essentiel}, le rang générique
$\rho_\infty$ de l'application infinie ${\sf F}_\infty$ des
coefficients est $\leqslant r - 1$. Mais lorsque les équations de
transformations constituent un {\em groupe} continu fini, on peut
établir (\voir~\cite{ merk2009b}) qu'en fait, le {\em rang} (tout
court) de ${\sf F}_\infty$ est {\em égal} à $\rho_\infty$ au paramètre
identité $e$, et donc aussi par voie de conséquence, égal à
$\rho_\infty$ dans un voisinage de $e$. Il existe donc au moins un
champ de vecteurs $\mathcal{ T} =
\sum_{ k=1}^n \, \tau_k ( a) \, \frac{ \partial }{ \partial 
a_k}$ {\em non nul en $e$} à coefficients analytiques tel que: $0
\equiv \mathcal{T}\,f_i(x;a) =
\sum_{k=1}^r\,\tau_k(a)\, \frac{\partial x_i'}{\partial a_k}(x;a)$ 
pour tout $i = 1, \dots, n$. En remplaçant maintenant chacune des
dérivées partielles $\frac{
\partial x_i'}{\partial a_k}$ par sa valeur déduite des équations
différentielles fondamentales~\thetag{ 5}, on obtient:
\[
0
\equiv
\sum_{j=1}^r\,
\sum_{k=1}^r\,\tau_k(a)\,\psi_{kj}(a)\,\xi_{ji}
\big(x'(x;a)\big)
\ \ \ \ \ \ \ \ \ \ \ \ \
{\scriptstyle{(i\,=\,1\,\cdots\,n)}}.
\]
Enfin, en posant $a = e$ dans ces équations, en rappelant $\psi_{
kj} ( e) = - \delta_k^j$, et en introduisant les
constantes $\lambda_j := {\textstyle{ \sum_{k=1}^r}}\,
\tau_k(e)\, \psi_{kj}(e) = -\tau_j(e)$ qui ne sont pas toutes 
nulles par hypothèse, on en déduit des équations:
\[
0
\equiv
\lambda_1\,\xi_{1i}(x)
+\cdots+
\lambda_r\,\xi_{ri}(x)
\ \ \ \ \ \ \ \ \ \ \ \ \
{\scriptstyle{(i\,=\,1\,\cdots\,n)}}
\]
qui montrent que les champs de vecteurs $X_k$ {\em ne}
sont {\em pas} linéairement indépendants. 
}. %%%%%%%%%%%%%%%%%%%%%%%%-----FIN-----%%%%%%%%%%%%%%%%%%%%%%%%%%%%%%%
Après avoir formé les transformations infinitésimales $X_k :=
\frac{ \partial f}{ \partial a_k} ( x; \, e)$, 
on peut donc tester de manière directe, effective et finitaire
l'essentialité des paramètres, en examinant seulement\footnote{\,
%%%%%%%%%%%%%%%%%%%%%%%-------DEBUT--------%%%%%%%%%%%%%%%%%%%%%%%%%%%
C'est
donc un exemple d'irréversible-synthétique: pour un groupe de
transformations, on contourne ainsi l'infinité potentiellement
imparcourable de toutes les vérifications qui seraient {\em a priori}
nécessaires en toute généralité afin de connaître le rang générique
{\em exact} de la matrice jacobienne ${\rm Jac}\, {\sf
F}_\infty$, dont les colonnes sont en nombre {\em infini}. La
généralité initiale du concept d'essentialité des paramètres 
se transforme en un critère concret, calculable et effectif. 
} %%%%%%%%%%%%%%%%%%%%%%%%-----FIN-----%%%%%%%%%%%%%%%%%%%%%%%%%%%%%%%
s'il existe des
combinaisons linéaires à coefficients {\em constants} entre les
colonnes de la matrice $\frac{ \partial f_i}{
\partial a_k} ( x; \, e)$.

\smallskip\noindent{\bf Définition.}
On dira que $r$ champs de vecteurs quelconques (qui ne proviennent pas
forcément d'un groupe de transformations) à coefficients analytiques
$X_k = \sum_{ i=1}^n \, \xi_{ ki} ( x) \, \frac{ \partial }{ \partial
x_i}$, $k = 1, \dots, r$, sont {\sl indépendants} s'ils sont {\em
linéairement indépendants}, à savoir, s'il n'existe pas de constantes
$\lambda_1, \dots, \lambda_r$ non toutes nulles telles que $\lambda_1
X_1 + \cdots + \lambda_r X_r
\equiv 0$.
\medskip

Ainsi un groupe dont les paramètres sont essentiels donne-t-il
naissance à une collection de transformations infinitésimales
{\em indépendantes}. 

\smallskip\noindent
{\bf Question.}
Réciproquement, peut-on reconstituer le groupe
de transformations $x' = f( x;
\, a_1, \dots, a_r)$ à partir des $r$ transformations infinitésimales
$X_1, \dots, X_r$ qui donnent toutes les directions de mouvement
infinitésimal possible, lorsqu'on fait subir un incrément
infinitésimal $(e_1, \dots, e_k + \varepsilon, \dots, e_n)$ à la
$k$-ième coordonnée du paramètre-identité?
\medskip

Principe de retour en arrière de la genèse: on sait bien que la
théorie de l'intégration permet de reconstituer le fini à partir de
l'infinitésimal. Cette question en soulève d'autres. Comment intégrer
un système de champs de vecteurs $X_1, \dots, X_r$? Retrouve-t-on
toujours une famille de transformations fermée par composition?

Lorsque $r \geqslant 2$, la <<\,non-canonicité\,>> des $r$
transformations infinitésimales $X_1, \dots, X_r$ se révèle
immédiatement. En effet, dès qu'on effectue un changement de
paramètres $a \mapsto b = b ( a)$ fixant 
l'identité\footnote{\,
%%%%%%%%%%%%%%%%%%%%%%%-------DEBUT--------%%%%%%%%%%%%%%%%%%%%%%%%%%%
On note alors $b \mapsto a = b ( a)$ le changement
inverse de paramètres. 
}, %%%%%%%%%%%%%%%%%%%%%%%%-----FIN-----%%%%%%%%%%%%%%%%%%%%%%%%%%%%%%%
ce qui transforme $f ( x; \, a)$ en $g ( x; \, b ) := f ( x; \, a (
b))$, les $r$ transformations infinitésimales $Y_k := \frac{ \partial
g}{ \partial b_k} ( x; \, e)$ calculées de la même manière dans le nouvel
espace des paramètres deviennent certaines combinaisons\footnote{\,
%%%%%%%%%%%%%%%%%%%%%%%-------DEBUT--------%%%%%%%%%%%%%%%%%%%%%%%%%%%
La règle de dérivation des fonctions composées donne
en effet immédiatement:
$\frac{ \partial g_i}{ \partial b_k} (x; \, e) = \sum_{ l=1}^r \,
\frac{ \partial f_i}{ \partial a_l} ( x; \, e) \, 
\frac{ \partial a_l}{ \partial b_k} (e)$. 
} %%%%%%%%%%%%%%%%%%%%%%%%-----FIN-----%%%%%%%%%%%%%%%%%%%%%%%%%%%%%%%
linéaires à coefficients constants des précédentes transformations
infinitésimales $X_1, \dots, X_r$. 
De telles combinaisons linéaires interviennent aussi nécessairement
lorsqu'on calcule les transformations infinitésimales $X_k$ non pas en
l'identité $e$, mais en un paramètre $a$ quelconque. Ainsi ces deux
observations montrent clairement que la structure linéaire est
centrale.

\HEAD{Chapitre~3.\,\,\,\,Théorèmes fondamentaux sur les groupes
de transformations}{ 3.6.\,\,\,Champs
de vecteurs et groupes à un paramètre}

\medskip\noindent{\bf 3.6.~Champs
de vecteurs et groupes à un paramètre.}
\label{champs-vecteurs-groupe-un-parametre}
Pour commencer, étudions tout d'abord la question posée à l'instant
dans le cas d'une seule transformation, c'est-à-dire lorsque $r = 1$.
Dans la théorie de Lie s'exprime une affirmation d'équivalence
ontologique qui a de nombreuses implications et qui demande des
explications:

\smallskip
\centerline{{\em transformation infinitésimale}
$\equiv$ {\em champ de vecteurs quelconque}}

\smallskip\noindent
En effet et pour commencer, rappelons que tout champ de vecteurs:
\[
X 
=
\sum_{i=1}^n\,\xi_i(x)\,
\frac{\partial}{\partial x_i} 
\]
n'est pas seulement un opérateur de dérivation, il vient aussi
accompagné de son {\sl flot}\,\,---\,\,dont nous admettrons
l'existence, \voir~\cite{ kono1963, spiv1970, lang1983,
ol1995}\,\,---\,\,noté:
\[
(x;\,t)\longmapsto
\exp(tX)(x).
\]
Ce flot est défini géométriquement en
suivant la courbe intégrale de $X$ issue d'un point $x$ jusqu'au
<<\,temps\,>> $t$, ou bien, de manière équivalente et d'un point de vue
analytique, en intégrant le système d'équations différentielles
ordinaires:
\[
\frac{dx_i}{dt}
=
\xi_i\big(x_1(t),\dots,x_n(t)\big)
\ \ \ \ \ \ \ \ \ \ \ \ \
{\scriptstyle{(i\,=\,1\,\cdots\,n)}}
\]
avec la condition initiale $x_i ( 0) = x_i$, la valeur au temps $t$ de
la solution étant justement
$\exp ( tX) ( x)$. Or le flot est un groupe à
un paramètre de difféomorphismes locaux:
\[
\exp(t_2X)\big(\exp(t_1X)(x)\big)
=
\exp\big((t_1+t_2)(X)\big)(x),
\]
l'identité correspondant au paramètre $t = 0$ et l'inverse de $x \mapsto
\exp ( tX) ( x)$ étant tout simplement
$x \mapsto \exp( - t X) ( x)$.

Soit donc maintenant
$x_i' = f_i ( x; \, a)$ un groupe continu de transformations 
(au sens de Lie) à {\em un
seul} paramètre $a \in \C$ (ou bien $a \in \R$)
pour lequel la transformation
identique correspond au paramètre $a = 0$. Les équations
différentielles fondamentales~\thetag{ 5} ci-dessus:
\[
\frac{dx_i'}{da}
=
\psi(a)\cdot\xi_i(x_1',\dots,x_n')
\ \ \ \ \ \ \ \ \ \ \ \ \
{\scriptstyle{(i\,=\,1\,\cdots\,n)}}
\]
consistent alors en un système d'équations différentielles
ordinaires du premier ordre paramétrées par $a$. Mais en introduisant
le nouveau <<\,paramètre-temps\,>> défini par:
\[
t
=
t(a)
:=
\int_0^a\,
\psi(\tilde{a})\,d\tilde{a},
\]
on peut transférer immédiatement ces équations différentielles
fondamentales vers un système de $n$ équations différentielles
ordinaires dont tous les seconds membres deviennent
{\em indépendants du temps}:
\def\theequation{1}\begin{equation}
\frac{dx_i'}{dt}
=
\xi_i(x_1',\dots,x_n')
\ \ \ \ \ \ \ \ \ \ \ \ \
{\scriptstyle{(i\,=\,1\,\cdots\,n)}}.
\end{equation}
L'intégration revient par conséquent à calculer le {\sl flot} du champ
de vecteurs:
\[
X'
:=
\sum_{i=1}^n\,\xi_i(x')\,
\frac{\partial}{\partial x_i'}, 
\]
vu dans l'espace des $(x_1', \dots, x_n')$. Avec les mêmes lettres
$f_i$, nous écrirons $f_i ( x; \, t)$ à la place de $f_i ( x; \,
a(t))$ en supposant que ce changement de paramètre a déjà été
exécuté. Alors en fait, l'unique solution du système~\thetag{ 1} avec
la condition initiale $x_i' ( x; \, 0) = x_i$ n'est autre que $x_i' =
f_i ( x; \, t)$: le flot était donné et connu depuis le début. De
plus, l'unicité du flot et le fait que les coefficients $\xi_i$ sont
indépendants de $t$ impliquent tous deux que la loi de composition de
groupe correspond alors juste à l'addition des paramètres
<<\,temporels\,>>:
\def\theequation{2}\begin{equation}
f_i\big(f(x;\,t_1);\,t_2)
\equiv
f_i(x;\,t_1+t_2)
\ \ \ \ \ \ \ \ \ \ \ \ \
{\scriptstyle{(i\,=\,1\,\cdots\,n)}}.
\end{equation}
En fait, on peut redresser localement\footnote{\,
%%%%%%%%%%%%%%%%%%%%%%%-------DEBUT--------%%%%%%%%%%%%%%%%%%%%%%%%%%%
Rappelons que le principe de relocalisation sous-entend que l'on se
replace en un point {\em générique} où le champ $X$ ne s'annule pas
(sinon il est identiquement nul, cas inintéressant).
} %%%%%%%%%%%%%%%%%%%%%%%%-----FIN-----%%%%%%%%%%%%%%%%%%%%%%%%%%%%%%%
le champ $X'$ en le champ très simple $\frac{
\partial }{\partial y_1'}$ dans un certain nouveau 
système de coordonnées $(y_1', y_2', \dots, y_n')$, et alors
le fait que les paramètres temporels s'additionnent
devient évident. 

\smallskip\noindent{\bf Théorème.}
\label{redressement-champ}
{\em Tout groupe continu à un paramètre:
\[
x_i'
=
f_i(x_1,\dots,x_n;\,t)
\ \ \ \ \ \ \ \ \ \ \ \ \
{\scriptstyle{(i\,=\,1\,\cdots\,n)}}
\]
est localement équivalent, à travers un
changement approprié de variables $x\mapsto y =
y ( x)$, à un groupe de translations: }
\[
y_1'
=
y_1+t,
\ \ \ \ \ \ \ \ \
y_2'
=
y_2,\ \ \
\dots\dots,\ \ \ 
y_n'
=
y_n.
\]

\begin{center}
\input flow-straightening.pstex_t
\end{center}

{\em Preuve.}
On peut supposer depuis le début que les coordonnées $x_1, \dots, x_n$
ont été choisies et adaptées de telle sorte que $\xi_1 ( 0) = 1$ et
que $\xi_2 ( 0 ) = \cdots = \xi_n ( 0) = 0$. Dans un espace auxiliaire
$y_1, \dots, y_n$ représenté sur la gauche de la figure, considérons
tous les points $\widehat{ y} := (0, y_2,
\dots, y_n)$ proches de l'origine qui se trouvent sur l'hyperplan de
coordonnées complémentaire de l'axe des $y_1$, et introduisons le
difféomorphisme:
\[
y
\longmapsto
x
=
x(y)
:=
f\big(0,\widehat{y};\,y_1\big)
=:
\Phi(y),
\]
défini en suivant le flot $f ( y; \, t)$ jusqu'à un temps $t = y_1$ en
partant de tous ces points possibles $\widehat{ y}$ dans l'hyperplan
en question. Cette application est effectivement un difféomorphisme
local fixant l'origine, grâce à $\frac{ \partial \Phi_1}{ \partial y_1
} ( 0) = \frac{
\partial f_1 }{ \partial t} ( 0) = \xi_1 ( 0) = 1$, grâce à $\frac{
\partial \Phi_k }{ \partial y_1} ( 0) = \frac{ \partial f_k}{
\partial t} ( 0) = \xi_k(0) = 0$ et grâce à $\Phi_k (0, \widehat{ y })
\equiv \widehat{ y }_k$ for $k = 2, \dots, n$. 
Nous affirmons maintenant que le flot (courbé) qui est représenté dans
la partie droite de la figure a ainsi été redressé à gauche pour
devenir juste une translation uniforme dirigée par l'axe des $y_1$.
En effet\footnote{\,
%%%%%%%%%%%%%%%%%%%%%%%-------DEBUT--------%%%%%%%%%%%%%%%%%%%%%%%%%%%
\`A travers le difféomorphisme local
$y \mapsto x = \Phi (y)$, 
le flot $x' = f ( x; \, t)$
se transforme naturellement par substitution en le flot 
bien déterminé
$\Phi (y') = f \big( \Phi ( y); \, t\big)$, 
ou bien, de manière équivalente en le flot: 
$y' = \Phi^{ -1} \big( f \big( \Phi ( y); \, t\big) \big)$. 
}, %%%%%%%%%%%%%%%%%%%%%%%%-----FIN-----%%%%%%%%%%%%%%%%%%%%%%%%%%%%%%%
si l'on suppose à l'avance que le nouveau flot est
effectivement donné par $\widehat{ y}' = \widehat{ y}$ et par $y_1' =
y_1 + t$, et si on effectue ensuite des substitutions:
\[
x'
=
f\big(
0,\widehat{y}';\,y_1'\big)
=
f\big(
0,\widehat{y};\,y_1+t\big)
\overset{(2)}{=}
f\big(f\big(
0,\widehat{y};\,y_1
\big);\,t\big)
=
f(x;\,t),
\]
nous retrouvons de cette manière-là le flot $x' = f ( x; \, t)$ qui
était défini de manière unique.
\qed\smallskip

En général, dans les traités contemporains de géométrie
différentielle, les flot sont étudiés sous l'hypothèse que les champs
de vecteurs soient différentiables, de classe au moins $\mathcal{
C}^1$, voire lipschitziens. Mais dans le cas où les coefficients des
champs de vecteurs sont tous analytiques (hypothèse générale admise
par Engel et Lie), il existe une série entière explicite et simple,
appelée {\sl série de Lie} (et étudiée entre autres par Gröbner
\cite{
gr1960}), qui redonne le flot {\em sans aucune intégration},
directement à partir des coefficients du champ de vecteurs $X' =
\sum_{ i = 1}^n \, \xi_i ( x') \, \frac{ \partial }{ \partial x_i'}$.

Admettons donc le résultat d'après lequel le flot (local ou global)
est analytique, si les $\xi_i ( x')$ le sont\footnote{\,
%%%%%%%%%%%%%%%%%%%%%%%-------DEBUT--------%%%%%%%%%%%%%%%%%%%%%%%%%%%
{\em Voir}~\cite{ lang1983}; nous
pourrions en fait ré-établir un tel résultat 
{\em via} des techniques de
séries majorantes en partant de la formule exponentielle de Lie
que nous allons donner à l'instant. 
}. %%%%%%%%%%%%%%%%%%%%%%%%-----FIN-----%%%%%%%%%%%%%%%%%%%%%%%%%%%%%%%
La solution formelle $x' (
x; \, t)$ du système d'équations différentielles ordinaires~\thetag{
1} satisfaisant $x' ( x; \, 0) = x$ peut être recherchée en développant
l'inconnue $x'$ en série formelle par rapport à la variable temporelle
$t$, les deux premiers termes étant évidents:
\[
x_i'(x;\,t)
=
\sum_{k\geqslant 0}\,\Xi_{ik}(x)\,t^k
=
x_i
+
t\,\xi_i(x)
+\cdots
\ \ \ \ \ \ \ \ \ \ \ \ \
{\scriptstyle{(i\,=\,1\,\cdots\,n)}}.
\]
Ainsi donc, quelles sont les fonctions-coefficients $\Xi_{ ik} ( x)$?
En différentiant~\thetag{ 1} une première fois par rapport à $t$, puis
en redérivant encore une fois le résultat obtenu, tout en
resubstituant comme il se doit, on obtient par exemple:
\[
\small
\left[
\aligned
\frac{d^2x_i'}{dt^2}
&
=
\sum_{k=1}^n\,\frac{\partial\xi_i}{\partial x_k'}\,
\frac{dx_k'}{dt}
=
X'(\xi_i)
\ \ \ \ \ \ \ \ \ \ \ \ \
{\scriptstyle{(i\,=\,1\,\cdots\,n)}}
\\
\frac{d^3x_i'}{dt^3}
&
=
\sum_{k=1}^n\,\frac{\partial X'(\xi_i)}{\partial x_k'}\,
\frac{dx_k'}{dt}
=
X'\big(X'(\xi_i)\big),
\endaligned\right.
\]
où $X' ( \xi_i)$ et $X '( X '( \xi_i))$ sont
des fonctions de $(x_1', \dots, x_n')$, 
et donc généralement, par une récurrence évidente:
\[
\small
\aligned
\frac{d^kx_i'}{dt^k}
=
\underbrace{X'\big(\cdots\big(
X'}_{k-1\,\,\text{\rm fois}}(\xi_i)\big)\cdots\big)
=
\underbrace{X'\big(\cdots\big(
X'\big(X'}_{k\,\,\text{\rm fois}}(x_i')\big)\big)\cdots\big),
\endaligned
\]
pour tout entier $k \geqslant 1$, sachant que $X' ( x_i' ) = \xi_i =
\xi_i ( x')$;
il sera utile de convenir que ${X'}^0 x_i' = x_i '$ lorsque $k = 0$.

En posant maintenant $t = 0$ dans ces équations $\frac{ d^k x_i'}{
dt^k} = {X'}^k ( x_i')$, nous obtenons donc l'expression recherchée
des fonctions $\Xi_{ ik} ( x)$:
\[
\small
\aligned
k!\,\,\Xi_{ik}(x)
\equiv
{X'}^k(x_i')
\big\vert_{t=0}
=
X^k(x_i),
\endaligned
\]
où $X := \sum_{ i=1}^n \, \xi_i ( x) \,
\frac{ \partial }{ \partial x_i}$ est le
même champ de vecteurs que $X'$, vu dans les coordonnées $(x_1
,\dots, x_n)$. Ainsi de manière surprenante, l'{\em intégration} d'un
flot dans le cas analytique revient à la sommation d'un nombre infini de termes {\em
différentiés}.

\smallskip\noindent{\bf Proposition.}
{\em L'unique solution $x' ( x; \, t)$ d'un système d'équations
différentielles ordinaires:
\[
\left[
\aligned
&
\frac{dx_i'}{dt}
=
\xi_i(x_1',\dots,x_n')
\\ 
&
x_i'(x;\,0) 
= 
x_i
\ \ \ \ \ \ \ \ \ \ \ \ \ \ \ \ \ \ \ \ \
{\scriptstyle{(i\,=\,1\,\cdots\,n)}}
\endaligned\right.
\]
qui est associé ou bien à un groupe de transformations à un seul
paramètre {\em via} les équations différentielles fondamentales
de Lie, ou
bien à un champ de vecteurs quelconque $X' = \sum_{ i=1}^n \,
\xi_i ( x') \, \frac{ \partial}{ \partial x_i'}$, 
est fournie par le développement en série entière convergent:
\def\theequation{3}\begin{equation}
x_i'(x;\,t)
=
x_i
+
t\,X(x_i)
+\cdots+
{\textstyle{\frac{t^k}{k!}}}\,
X^k(x_i)
+\cdots
\ \ \ \ \ \ \ \ \ \ \ \ \
{\scriptstyle{(i\,=\,1\,\cdots\,n)}},
\end{equation}
où $X = \sum_{ i=1}^n \, \xi_i ( x) \, \frac{ \partial }{\partial x_i}$ 
est le même champ que $X'$ vu dans les coordonnées $(x_1, \dots, x_n)$,
qui agit sur $x_i$ comme dérivation $X^k$ d'ordre $k$ arbitraire. De
plus, ce développement peut aussi être réécrit de manière appropriée
au moyen d'une simple notation exponentielle: }
\def\theequation{3'}\begin{equation}
\label{flow-exponential-formula}
\aligned
x_i'
=
\exp\big(t\,X\big)(x_i)
=
\sum_{k\geqslant 0}\,\frac{(t\,X)^k}{k!}(x_i)
\ \ \ \ \ \ \ \ \ \ \ \ \
{\scriptstyle{(i\,=\,1\,\cdots\,n)}}.
\endaligned
\end{equation}

Ainsi cette proposition établit-elle l'équivalence ontologique
fondamentale: 

\smallskip
\centerline{\fbox{\bf \ groupe local à 
un paramètre $\,\equiv\,$ transformation infinitésimale }\,,}
\label{equivalence-ontologique}

\smallskip\noindent
qui s'insère plus généralement dans l'équivalence fonctionnelle entre
le différentiel infinitésimal et le local fini, tout en développant
les premiers éléments d'une théorie géométrique du
mouvement. Toutefois, cette équivalence ne saurait s'affirmer comme
principe d'égalité absolue entre deux êtres initialement distincts qui
deviendraient par là-même strictement interchangeables. Comme dans
toute autre équivalence mathématique, l'ontologie est interrogation
{\em en devenir} sur la structure et sur la constitution d'un être
mathématique problématique.
\`A travers les équivalences que
l'être mathématique découvre à son 
propre sujet par l'analyse ou par la
synthèse, l'être vise en effet à se déployer dans des espaces neufs
qui soient propices à révéler sa nature intrinsèque, en supprimant
toutes les formes d'arbitraire et de non-compréhension dont est
entachée sa donation initiale.

\`A strictement parler, aucun énoncé 
mathématique d'équivalence entre deux concepts ou entre deux
conditions spécifiques n'est réellement transparent dans une double
circulation du sens. L'équivalence, en mathématiques, transcende tout
concept logique ou méta-mathématique de termes formels syntaxiquement
substituables.

En vérité, dans l'équivalence, il doit se manifester un
différentiel-synthétique du potentiel interrogatif, comme par l'effet
d'une révélation progressive qui autoriserait à {\em oublier} presque
définitivement le membre initial de l'<<\,équivalence\,>> pour ne
retenir que le membre final, plus rapproché, bien que peut-être encore
fort éloigné, de l'essence de la chose à comprendre. Quant au choix de
brisure de symétrie dans l'équivalence, c'est-à-dire quant au choix de
l'initial et du final dans l'équivalence en question, c'est bien sûr
le concept de transformation infinitésimale qui doit en quelque sorte
{\em éliminer} le concept de groupe local à un paramètre, pour s'y
substituer comme objet d'étude principal. Et c'est effectivement ce
que Lie affirmera systématiquement dans ses travaux: grâce à
l'équivalence en question, on peut {\em mettre entre parenthèses}
l'intervention de l'analyse comme procédé d'intégration pour se
concentrer seulement sur la classification des transformations
infinitésimales. L'infinitésimal se substitue au fini, car il est plus
simple: il est {\em linéaire}.

\`A vrai dire, 
si l'on accepte comme Lie les règles de la généricité comme hypothèse
générale d'étude (\cf~p.~\pageref{relocalisation-libre}), le théorème
de redressement montre que le concept de transformation infinitésimale
individuelle est entièrement circonscrit. 
En effet, toute transformation
infinitésimale se réduit au champ de vecteurs le plus simple possible:
$\frac{ \partial }{ \partial x_1}$, somme réduite à un seul terme,
dérivation le long d'un seul axe de coordonnées, avec un coefficient
constant égal à $1$. L'infinitésimalisation des transformations
transfère la théorie des groupes vers l'algèbre, et tout d'abord, vers
l'algèbre {\em linéaire}. On peut donc maintenant se demander quelle
doit être l'équivalence entre groupe à $r$ paramètres et système de
$r$ transformations infinitésimales mutuellement indépendantes $X_1,
\dots, X_r$. Ici, l'Un va provoquer l'imprévu du Multiple, 
c'est-à-dire forcer la genèse de concepts inattendus. 
Mais pour l'instant, poursuivons l'étude 
des groupes à un seul paramètre. 

\smallskip

Engel et Lie écrivent les champs de vecteurs $X$\,\,---\,\,qu'ils
appellent systématiquement {\sl transformations
infinitésimales}\,\,---\,\,toujours sous la forme <<\,$Xf$\,>>, non
pas donc comme opérateur abstrait de dérivation, mais comme
action effective sur une fonction-test
toujours notée $f = f ( x_1, \dots,
x_n)$. L'action concrète
de $X$ sur $f$ consiste bien entendu à le faire agir
comme dérivation:
\[
Xf
= 
\sum_{i=1}^n\,\xi_i(x)\,
\frac{\partial f}{\partial x_i}. 
\]
Eu égard à l'équivalence ontologique sus-mentionnée, il est naturel de se
demander maintenant (scholie) quelle est l'action d'un groupe fini à un
paramètre sur les fonctions. Si donc $f = f ( x_1, \dots, x_n)$ est
une fonction analytique arbitraire, si l'on compose $f$ avec le flot
du groupe à un paramètre:
\[
f'
:=
f(x_1',\dots,x_n')
=
f\big(x_1'(x;\,t),\dots,x_n'(x;\,t)\big),
\]
et si l'on développe en série entière cette composition par rapport à
$t$:
\[
f'
=
\big(f'\big)_{t=0}
+
{\textstyle{\frac{t}{1!}}}\,
\big(
{\textstyle{\frac{df'}{dt}}}
\big)_{t=0}
+
{\textstyle{\frac{t^2}{2!}}}\,
\big(
{\textstyle{\frac{d^2f'}{dt^2}}}
\big)_{t=0}
+
\cdots,
\]
on doit calculer les quotients différentiels $\frac{ df'}{ dt}$,
$\frac{ d^2 f'}{ dt^2}$, \dots, ce qui ne pose à vrai dire aucune
difficulté:
\[
\small
\left[
\aligned
\frac{df'}{dt}
&
=
\sum_{i=1}^n\,
\frac{dx_i'}{dt}\,
\frac{\partial f'}{\partial x_i'}
=
\sum_{i=1}^n\,\xi_i'\,
\frac{\partial f'}{\partial x_i'}
=
X'(f'),
\\
\frac{d^2f'}{dt^2}
&
=
X'\Big(
\frac{df'}{dt}
\Big)
=
X'\big(
X'(f')\big),
\endaligned\right.
\]
et ainsi de suite. En posant $t = 0$, chaque $x_i'$ devient $x_i$, la
fonction $f'$ devient $f$ et $X' ( f')$ devient $X (f)$, et ainsi de
suite, et grâce à ces observations, on obtient le développement
recherché:
\def\theequation{4}\begin{equation}
\aligned
f(x_1',\dots,x_n')
&
=
f(x_1,\dots,x_n)
+
{\textstyle{\frac{t}{1!}}}\,X(f)
+
\cdots
+
{\textstyle{\frac{t^k}{k!}}}\,X^k(f)
+\cdots
\\
&
=
{\textstyle{\sum_{k\geqslant 0}}}\,
{\textstyle{\frac{(tX)^k}{k!}}}\,(f)
=
\exp(tX)(f),
\endaligned
\end{equation}
qui fait bien sûr réapparaître la symbolique
exponentielle. \'Evidemment, par différentiation, on doit retrouver
l'action infinitésimale:
\[
Xf 
= 
\frac{d}{dt}
\big(\exp(tX)(f)\big)_{t=0}.
\] 

Maintenant, soit $X = \sum_{ i=1 }^n\, \xi_i ( x) \, \frac{ \partial
}{\partial x_i}$ une transformation infinitésimale qui engendre le
groupe à un paramètre $x' = \exp ( t X) ( x)$. Qu'advient-il lorsque
les variables $x_1, \dots, x_n$ sont soumises à un changement de
coordonnées analytiques locales? Cette question générale est tout à
fait cruciale quant à l'individuation effective. En effet, toute
donation dans un système de coordonnées est entachée de quelconque et
d'arbitraire. Principe métaphysique: user au maximum de la liberté que
l'on a de changer le système de coordonnées afin de spécifier et de
simplifier au mieux les individus géométriques.

En premier lieu, on doit donc comprendre {\em a priori} (et en général)
de quelle façon les transformations infinitésimales arbitraires sont
modifiées, simplifiées ou complexifiées lorsqu'on effectue un
changement arbitraire de variables.

Soit donc $x \mapsto \overline{ x} = \overline{ x} ( x)$ un
changement de coordonnées analytiques locales, c'est-à-dire en faisant
apparaître les indices:
\[
(x_1,\dots,x_n)
\longmapsto 
(\overline{x}_1,\dots,\overline{x}_n)
= 
\big(
\overline{x}_1(x),\dots,
\overline{x}_n(x)
\big). 
\] 
Ici, en respectant la pensée de Lie, aucun symbole fonctionnel n'est
introduit, mais pour les besoins occasionnels de l'exégèse moderne, on
pourrait convenir d'appeler $x \mapsto \Phi (x) =
\overline{ x}$ ce difféomorphisme local, avec une lettre auxiliaire
$\Phi$. Néanmoins, le formalisme de Lie présente un avantage
considérable par rapport au formalisme contemporain, et ce, pour au
moins deux raisons.

\smallskip$\bullet$\,\,D'un point
de vue symbolique, le difféomorphisme inverse: 
\[
(\overline{x}_1,\dots,\overline{x}_n)
\longmapsto
(x_1,\dots,x_n)
=
\big(
x_1(\overline{x}),\dots,x_n(\overline{x})
\big)
\]
s'exprime exactement comme le difféomorphisme
initial, en intervertissant
seulement\footnote{\,
%%%%%%%%%%%%%%%%%%%%%%%-------DEBUT--------%%%%%%%%%%%%%%%%%%%%%%%%%%%
On se dispense ainsi d'avoir à résoudre la micro-question de
présentation symbolique formelle: <<\,est-il plus adapté de choisir la
lettre $\Phi$ ou bien la lettre $\Phi^{ -1}$ pour désigner le
difféomorphisme initial\,?\,>>.
} %%%%%%%%%%%%%%%%%%%%%%%%-----FIN-----%%%%%%%%%%%%%%%%%%%%%%%%%%%%%%%
les rôles de $x$ et de $\overline{ x}$. 

\smallskip$\bullet$\,\,Aussi bien dans les applications concrètes
que dans les démonstrations des théorèmes de classification, les
fonctions coordonnées-images sont toujours données {\em explicitement}
ou {\em spécifiquement} en fonction des coordonnées-sources $(x_1,
\dots, x_n)$, si bien que l'introduction de symboles fonctionnels
$\Phi_i$ serait superflue ou redondante. En fait, tout changement
de coordonnées sous-entend une {\em procédure de remplacement
automatique d'anciennes variables par de nouvelles variables}, ce que
la notation $x = x (\overline{ x})$ montre 
de manière adéquate.

\smallskip

Ainsi, une fonction arbitraire $f = f (x_1, \dots, x_n)$ définie dans
l'espace-source donne alors naissance, à travers le <<\,miroir\,>> du
changement de coordonnées, à la fonction correspondante:
\def\theequation{5}\begin{equation}
\overline{f}
(\overline{x}_1,\dots,\overline{x}_n)
:=
f\big(
x_1(\overline{x}),\dots,x_n(\overline{x})
\big)
\end{equation}
qui est définie dans l'espace-image. 
Cette relation s'écrit aussi de
manière symétrique:
\[
\overline{f}(\overline{x})
=
f(x).
\]
Le transfert d'une dérivation
arbitraire $X = \sum_{ i=1}^n \, \xi_i ( x) \, \frac{ \partial }{ 
\partial x_i}$ à travers ce miroir doit alors être tel que la 
relation purement symétrique:
\[
X(f)
=
\overline{X}(\overline{f}), 
\]
soit satisfaite pour toute fonction $f$, c'est-à-dire que, après
remplacement à droite de $\overline{ x}$ par $\overline{ x} (x)$ et
en utilisant la
différentiation de~\thetag{ 5} par rapport à $\overline{ x}_k$,
le transfert doit satisfaire:
\[
\aligned
X(f)
=
\sum_{i=1}^n\,\xi_i(x)\,
\frac{\partial f}{\partial x_i}
&
=
\sum_{k=1}^n\,\overline{\xi}_k(\overline{x})\,
\frac{\partial \overline{f}}{\partial\overline{x}_k}
\Big\vert_{\overline{x}=\overline{x}(x)}
\\
&
=
\sum_{k=1}^n\,\overline{\xi}_k\big(\overline{x}(x)\big)\,
\sum_{i=1}^n\,
\frac{\partial f}{\partial x_i}(x)\,
\frac{\partial x_i}{\partial\overline{x}_k}\big(\overline{x}(x)\big)
\\
&
=
\sum_{i=1}^n
\bigg(
\sum_{k=1}^n\,
\overline{\xi}_k\big(\overline{x}(x)\big)\,
\frac{\partial x_i}{\partial\overline{x}_k}\big(\overline{x}(x)\big)
\bigg)
\frac{\partial f}{\partial x_i}(x).
\endaligned
\]
Autrement dit, de manière abrégée et en supprimant le symbole de
fonction-test:
\[
X
=
\sum_{i=1}^n\,\overline{X}(x_i)\,
\frac{\partial}{\partial x_i}. 
\]
En intervertissant les rôles de $X$ et de $\overline{ X}$, nous avons
donc démontré l'énoncé suivant.

\smallskip\noindent{\bf Lemme.}
\label{changement-coordonnees-champ}
{\em
\`A travers un changement arbitraire de coordonnées locales:
\[
(x_1,\dots,x_n)
\longmapsto 
(\overline{x}_1,\dots,\overline{x}_n)
= 
\big(
\overline{x}_1(x),\dots,
\overline{x}_n(x)
\big),
\]
un champ de vecteurs quelconque $X = \sum_{ i=1}^n \, \xi_i ( x) \, 
\frac{ \partial }{ \partial x_i}$ est transféré vers:}
\[
\overline{X}
:=
\sum_{i=1}^n\,X(\overline{x}_i)\,
\frac{\partial}{\partial\overline{x}_i}. 
\]

De la même manière que l'on pouvait écrire l'équation $\overline{ f} =
f$, on peut maintenant écrire $\overline{ X} = X$, étant entendu que
ces deux écritures n'ont de sens que si l'on remplace partout
$\overline{ x}$ par $x$ (ou partout $x$ par $\overline{ x}$).

\`A présent, comment sont transformés les groupes
à un paramètre? L'énoncé suivant (\cf~\cite{ enlie1888, kono1963,
spiv1970}), que nous ne redémontrerons pas, découle aisément des
considérations précédentes.

\smallskip\noindent{\bf Proposition.}
{\em Le nouveau groupe à un paramètre $\overline{ x}' = \exp ( t
\overline{ X}) ( \overline{ x})$ associé
au nouveau champ de vecteurs $\overline{ X} =
\sum_{i=1}^n\,X(\overline{x}_i)\,
\frac{\partial}{\partial\overline{x}_i}$ peut être retrouvé à
partir de l'ancien $x' = \exp ( t X) ( x)$ grâce à la formule
canonique:
\[
\exp\big(t\overline{X}\big)(\overline{x})
=
\exp(tX)(x)
\big\vert_{x=x(\overline{x})}.
\]
En d'autres termes, l'ancien groupe à un paramètre:
\[
x_i'
=
x_i
+
{\textstyle{\frac{t}{1}}}\,X(x_i)
+
{\textstyle{\frac{t^2}{1\cdot 2}}}\,
X^2(x_i)
+\cdots
\ \ \ \ \ \ \ \ \ \ \ \ \
{\scriptstyle{(i\,=\,1\,\cdots\,n)}}
\]
est transféré vers le nouveau groupe à un paramètre:
\[
\overline{x}_i'
=
\overline{x}_i
+
{\textstyle{\frac{t}{1}}}\,\overline{X}(\overline{x}_i)
+
{\textstyle{\frac{t^2}{1\cdot 2}}}\,
\overline{X}^2(\overline{x}_i)
+\cdots
\ \ \ \ \ \ \ \ \ \ \ \ \
{\scriptstyle{(i\,=\,1\,\cdots\,n)}},
\]
où la variable $t$ reste la même dans les deux collections de $n$
équations.
}\medskip

\smallskip

Après ces préparatifs, nous pouvons maintenant revenir à la question
soulevée à la fin du \S~3.5: comment reconstituer les équations $x_i'
= f_i ( x; \, a)$ d'un groupe continu fini de transformations à partir
de ses transformations infinitésimales $X_1, \dots, X_r$\,? En
prenant une combinaison linéaire arbitraire $X := \lambda_1 X_1 +
\cdots + \lambda_r X_r$ de ces transformations, après
fixation des constantes $\lambda_k$, on peut considérer le groupe à un
paramètre $\exp (t X) (x)$ qui est engendré par $X$. Autrement dit, en
considérant après-coup que les constantes $\lambda_k$ peuvent 
aussi redevenir
variables, on obtient de nouvelles
équations de transformations:
\[
\aligned
x_i'
&
=
\exp\big(
t\,\lambda_1\,X_1
+\cdots+
t\,\lambda_r\,X_r\big)(x_i)
\\
&
=:
h_i\big(x;\,t,\lambda_1,\dots,\lambda_r\big)
\ \ \ \ \ \ \ \ \ \ \ \ \ \ \ \ \ \ \ \ \ \ \
{\scriptstyle{(i\,=\,1\,\dots\,n)}}
\endaligned
\]
paramétrées non seulement par le <<\,temps\,>> $t$, mais aussi par ces
constantes arbitraires $\lambda_k$. La réponse positive que l'on a
déjà devinée énoncera que ces équations de transformations, appelées
{\em équations finies canoniques} du groupe par Engel et Lie,
recouvrent complètement les équations originales $x_i' = f_i ( x; \,
a)$.

Avant de poursuivre, résumons le parcours spéculatif et soulevons des
questions nouvelles. Le Multiple $X_1,
\dots, X_r$ comprend l'Un-divers $X = \lambda_1 X_1 + \cdots +
\lambda_r X_r$ de transformations infinitésimales 
individuelles possibles. Chaque tel
Un-divers jouit pleinement de
l'équivalence ontologique avec le groupe à un paramètre qu'il
engendre. Ainsi le groupe à un paramètre de l'Un-divers fournit-il
gratuitement des équations de transformations à plusieurs
paramètres. La question de savoir comment ces équations de
transformations $x_i' = h_i ( x; \, t, \lambda_1, \dots, \lambda_r)$
sont reliées aux anciennes équations se divise alors en deux questions.

\smallskip$\bullet$\,\,
\'Etant donné $r$ transformations infinitésimales linéairement
indépendantes $X_1, \dots, X_r$ qui proviennent d'un groupe continu à
$r$ paramètres essentiels $x' = f ( x; \, a_1, \dots, a_r)$, comment
choisir\footnote{\,
%%%%%%%%%%%%%%%%%%%%%%%-------DEBUT--------%%%%%%%%%%%%%%%%%%%%%%%%%%%
Nous venons d'annoncer que l'on retrouve les équations d'origine, mais
il aurait très bien pu se produire que les équations $x ' = \exp \big(
t\lambda_1 X_1 + \cdots + t\lambda_r X_r) (x)$ leur soient purement
étrangères. }
%%%%%%%%%%%%%%%%%%%%%%%%-----FIN-----%%%%%%%%%%%%%%%%%%%%%%%%%%%%%%%
les paramètres $t, \lambda_1, \dots, \lambda_r$ dans les équations de
transformations $x ' = \exp \big( t\lambda_1 X_1 + \cdots + t\lambda_r
X_r) (x)$ pour retrouver $x' = f ( x; \, a_1, \dots, a_r)$?

\smallskip$\bullet$\,\,
\'Etant donné $r$ transformations infinitésimales linéairement
indépendantes quelconques $X_1, \dots, X_r$ qui ne proviennent pas
forcément d'un groupe continu fini, les équations de transformations $x' 
= \exp \big( t\lambda_1 X_1 + \cdots + t\lambda_r X_r) (x)$
constituent-elles un groupe continu fini? Et si tel n'est pas le cas,
sous quelles conditions, nécessaires, suffisantes, ou mieux encore:
{\em nécessaires et suffisantes}\footnote{\,
%%%%%%%%%%%%%%%%%%%%%%%-------DEBUT--------%%%%%%%%%%%%%%%%%%%%%%%%%%%
Nécessité et suffisance en toute circonstance: exigence riemannienne
universelle.
} %%%%%%%%%%%%%%%%%%%%%%%%-----FIN-----%%%%%%%%%%%%%%%%%%%%%%%%%%%%%%%
cette conclusion est-elle satisfaite? 

\smallskip

La seconde question exige la notion cruciale de crochet de Lie
(\S\S~3.8 et 3.9 ci-dessous), tandis que la réponse à la première peut
s'en dispenser: le passage à un niveau supérieur d'abstraction dans
les conditions de donation implique un renversement
complet du champ synthétique.

\smallskip

Tout d'abord, au sujet des fonctions $h_i ( x; \, t, \lambda_1, \dots,
\lambda_r)$ introduites à l'instant, une question se pose
immédiatement: les $r+1$ paramètres $t$ et $\lambda_1, \dots,
\lambda_r$ y sont-ils essentiels? Certainement pas: ces $r+1$
paramètres se réduisent en fait à $r$ paramètres (au maximum), puisque
chaque $\lambda_i$ apparaît multiplié par $t$ dans $\exp ( t\lambda_1
X_1 + \cdots + t \lambda_r X_r)(x)$. Toutefois, en posant $t = 1$,
Engel et Lie établissent l'essentialité des $r$ paramètres restants
$\lambda_1, \dots, \lambda_r$ par une démonstration sophistiquée et
ingénieuse que nous restituons, 
en la modernisant et avec de plus amples détails, à la fin
de ce paragraphe. 

\smallskip\noindent{\bf Théorème 8.}
\label{Theoreme-8}
(\cite{ enlie1888}, p.~65)
{\em Si $r$ transformations infinitésimales:
\[
X_k(f)
=
\sum_{i=1}^n\,\xi_{ki}(x_1,\dots,x_n)\,
\frac{\partial f}{\partial x_i}
\ \ \ \ \ \ \ \ \ \ \ \ \
{\scriptstyle{(k\,=\,1\,\cdots\,r)}}
\]
sont mutuellement [linéairement] indépendantes et si, de plus,
$\lambda_1, \dots, \lambda_r$ sont des paramètres arbitraires, alors
la totalité \deutsch{Inbegriff} des groupes\footnote{\,
%%%%%%%%%%%%%%%%%%%%%%%-------DEBUT--------%%%%%%%%%%%%%%%%%%%%%%%%%%%
Le principe d'équivalence ontologique énoncé
p.~\pageref{equivalence-ontologique} s'exerce déjà implicitement ici:
tout groupe à un paramètre est {\em identifié} par Lie à son
générateur infinitésimal.
} %%%%%%%%%%%%%%%%%%%%%%%%-----FIN-----%%%%%%%%%%%%%%%%%%%%%%%%%%%%%%%
à un paramètre $\lambda_1 \, X_1 ( f) + \cdots +
\lambda_r \, X_r (f)$ forme une famille de transformations:
\[
x_i'
=
x_i
+
\sum_{k=1}^r\,\lambda_k\cdot\xi_{ki}
+
\sum_{k,\,j}^{1\dots r}\,
\frac{\lambda_k\,\lambda_j}{1\cdot 2}\cdot
X_k(\xi_{ji})
+\cdots
\ \ \ \ \ \ \ \ \ \ \ \ \
{\scriptstyle{(i\,=\,1\,\cdots\,n)}},
\]
dans lesquelles les $r$ paramètres $\lambda_1, \dots, \lambda_r$ sont
tous \emphasis{essentiels}, donc une famille de $\infty^r$
transformations différentes\footnote{\,
%%%%%%%%%%%%%%%%%%%%%%%-------DEBUT--------%%%%%%%%%%%%%%%%%%%%%%%%%%%
La notation <<\,$\infty^r$\,>> désigne le nombre de paramètres {\em
continus} (d'où le symbole d'infinité $\infty$) dont dépend un objet
analytique ou géométrique.
}. %%%%%%%%%%%%%%%%%%%%%%%%-----FIN-----%%%%%%%%%%%%%%%%%%%%%%%%%%%%%%%
}\medskip

Fait remarquable: toute question spéculative ou rhétorique qui
apparaît naturellement est traitée par Engel et Lie {\em au moment
approprié} dans le continu temporel et mémoriel du déploiement
irréversible de la théorie. S'il existe une systématique du
questionnement, c'est dans les mathématiques d'inspiration
riemannienne qui se sont développées pendant la deuxième moitié du
19\textsuperscript{ième} siècle qu'il faut en trouver les racines,
bien avant que l'axiomatique formelle du 20\textsuperscript{ième}
siècle ne l'enfouisse sous des strates de reconstitution {\em a
posteriori} et {\em non ouvertement problématisante}.

Maintenant, examinons la première question. Pour $\lambda_1, \dots,
\lambda_r$ fixés, la transformation infinitésimale combinaison
linéaire $X := \lambda_1 \, X_1 + \cdots + \lambda_r \, X_r$ a pour
coefficients:
\[
\xi_i(x)
:=
\sum_{j=1}^n\,\lambda_j\,\xi_{ji}(x)
\ \ \ \ \ \ \ \ \ \ \ \ \
{\scriptstyle{(i\,=\,1\,\dots\,n)}}. 
\]
Par définition du flot $x' = \exp ( tX) ( x)$, les fonctions:
\[
h_i(x;\,t,\lambda_1,\dots,\lambda_r)
=
\exp\big(t\lambda_1X_1+\cdots+t\lambda_rX_r)(x_i)
\]
satisfont le système d'équations différentielles ordinaires:
\[
\frac{dh_i}{dt}
= 
\xi_i(h_1,\dots,h_n)
\ \ \ \ \ \ \ \ \ \ \ \ \
{\scriptstyle{(i\,=\,1\,\dots\,n)}}, 
\]
ou bien, de manière équivalente: 
\def\theequation{6}\begin{equation}
\frac{dh_i}{dt}
=
\sum_{j=1}^r\,\lambda_j\,\xi_{ji}(h_1,\dots,h_n)
\ \ \ \ \ \ \ \ \ \ \ \ \
{\scriptstyle{(i\,=\,1\,\dots\,n)}},
\end{equation}
avec bien sûr la condition initiale $h ( x; 0, \lambda) = x$ lorsque
$t = 0$. D'un autre côté, d'après le
paragraphe qui précède le Théorème p.~\pageref{Theoreme-3},
les fonctions $f_i$ des équations de transformations
d'origine $x_i' =
f_i ( x; \, a)$ satisfont les équations différentielles 
fondamentales:
\def\theequation{4}\begin{equation}
\sum_{k=1}^r\,
\frac{\partial f_i}{\partial a_k}\,
\widetilde{\psi}_{jk}(a)
=
\xi_{ji}(f_1,\dots,f_n)
\ \ \ \ \ \ \ \ \ \ \ \ \ \
{\scriptstyle{(i\,=\,1\,\dots\,n;\,\,\,j\,=\,1\,\dots\,r)}}.
\end{equation}
La correspondance recherchée entre les équations de transformations
initiales $x' = f ( x; \, a)$ du groupe et ses équations finies
canoniques $x' = h (x; \, t, \lambda)$ est maintenant
fournie par la proposition
suivante, qui résout la première question. 

\smallskip\noindent{\bf Proposition.}
\label{p-69}
{\em Si les paramètres $a_1, \dots, a_r$ sont les uniques solutions
$a_k ( t, \lambda)$ du système d'équations différentielles
ordinaires du premier ordre:
\[
\frac{da_k}{dt}
=
\sum_{j=1}^r\,\lambda_j\,
\widetilde{\psi}_{jk}(a)
\ \ \ \ \ \ \ \ \ \ \ \ \
{\scriptstyle{(k\,=\,1\,\dots\,r)}}
\]
satisfaisant la condition initiale 
d'après laquelle $a ( 0, \lambda) = e$ est
l'élément identité, alors les équations suivantes sont identiquement
satisfaites:
\[
\aligned
f_i\big(x;\,a(t,\lambda)\big)
&
\equiv
\exp
\big(
t\,\lambda_1\,X_1
+\cdots+
t\,\lambda_r\,X_r
\big)
(x_i)
=
h_i\big(x;\,t,\lambda_1,\dots,\lambda_r\big)
\\
&
\ \ \ \ \ \ \ \ \ \ \ \ \ \ \ \ \ \ \ \ \ \ \ \ \ \
{\scriptstyle{(i\,=\,1\,\dots\,n)}}
\endaligned
\]
et elles montrent comment les fonctions $h_i$
se déduisent des fonctions $f_i$. 
}\medskip

{\em Preuve.}
En effet, multiplions~\thetag{ 4} ci-dessus par $\lambda_j$ et sommons
par rapport à $j$ pour $j$ allant
de $1$ jusqu'à $r$, ce qui nous donne:
\[
\sum_{k=1}^r\,
\frac{\partial f_i}{\partial a_k}\,
\sum_{j=1}^r\,\lambda_j\,\widetilde{\psi}_{jk}(a)
=
\sum_{j=1}^r\,\lambda_j\,\xi_{ji}(f_1,\dots,f_n)
\ \ \ \ \ \ \ \ \ \ \ \ \
{\scriptstyle{(i\,=\,1\,\dots\,n)}}.
\]
Grâce à l'hypothèse principale concernant les $a_k$, nous pouvons
remplacer la seconde somme du membre de gauche par $\frac{ da_k}{
dt}$, et nous obtenons ainsi des identités:
\[
\sum_{k=1}^r\,
\frac{\partial f_i}{\partial a_k}\,
\frac{da_k}{dt}
\equiv
\sum_{j=1}^r\,\lambda_j\,\xi_{ji}(f_1,\dots,f_n)
\ \ \ \ \ \ \ \ \ \ \ \ \
{\scriptstyle{(i\,=\,1\,\dots\,n)}}
\]
dans la partie gauche desquelles nous reconnaissons une simple
dérivation par rapport à $t$:
\[
\frac{df_i}{dt}
=
\frac{d}{dt}
\big[
f_i\big(x;\,a(t,\lambda)\big)
\big]
\equiv
\sum_{j=1}^r\,\lambda_j\,\xi_{ji}(f_1,\dots,f_n)
\ \ \ \ \ \ \ \ \ \ \ \ \
{\scriptstyle{(i\,=\,1\,\dots\,n)}}.
\]
Mais puisque $f \big( x; \, a ( 0, \lambda) \big) = f ( x; \, e) = x$
est soumis à la même condition initiale que la solution $h \big( x; \,
t, \lambda\big)$ du système~\thetag{ 6}, lorsque $t = 0$, l'unicité des
solutions à un système d'équations différentielles ordinaires du
premier ordre implique immédiatement la coïncidence annoncée: $f \big(
x; \, a ( t, \lambda) \big) \equiv h\big( x; \, t, \lambda\big)$.
\qed

\bigskip

{\small

{\em Démonstration du Théorème~8 p.~\pageref{Theoreme-8}.} Ici
exceptionnellement, nous avons observé une légère incorrection
technique (la seule que nous ayons pu découvrir!)
dans la preuve écrite par Engel et par Lie (\cite{
enlie1888}, pp.~62--65) au sujet du lien entre le rang générique de
$X_1 \big\vert_x, \dots, X_r \big\vert_x$ et une borne inférieure pour
le nombre des paramètres essentiels\footnote{\,
%%%%%%%%%%%%%%%%%%%%%%%-------DEBUT--------%%%%%%%%%%%%%%%%%%%%%%%%%%%
\`A la page~63, il est dit que si le
nombre $r$ de transformations infinitésimales indépendantes $X_k$ est
$\leqslant n$, alors la matrice $\big( \xi_{ ki} ( x) \big)_{
1\leqslant k \leqslant r}^{ 1 \leqslant i \leqslant n}$ (de taille $r
\times n$) de leurs coefficients est de rang générique égal à $r$,
mais cette assertion est contredite pour $n = r = 2$ par les deux
champs de vecteurs $x \, \frac{ \partial }{\partial x} + y \,
\frac{ \partial }{\partial y}$ et $xx \, \frac{ \partial }{\partial
x} + xy \, \frac{ \partial }{\partial y}$. Toutefois les idées et les
arguments de la preuve présentée (qui ne nécessite en fait pas de
telle assertion) sont parfaitement corrects.
}. %%%%%%%%%%%%%%%%%%%%%%%%-----FIN-----%%%%%%%%%%%%%%%%%%%%%%%%%%%%%%%

L'idée géométrique principale de Lie est astucieuse et pertinente: elle
consiste à introduire exactement $r$\,\,---\,\,le nombre des
$\lambda_k$\,\,---\,\,copies du même espace $x_1, \dots, x_n$ dont les
coordonnées seront notées $x_1^{ ( \mu)}, \dots, x_n^{ ( \mu)}$ pour
$\mu = 1, \dots, r$ et à considérer la famille des équations de
transformations qui sont induites par les {\em mêmes équations de
transformations}:
\[
{x_i^{(\mu)}}' 
=
\exp(C)\big(x_i^{(\mu)}\big)
= 
h_i\big(x^{(\mu)};\,
\lambda_1,\dots,\lambda_r\big)
\ \ \ \ \ \ \ \ \ \ \ \ \
{\scriptstyle{(i\,=\,1\,\cdots\,n;\,\,\,\mu\,=\,1\,\cdots\,r)}}
\]
dans chaque copie de l'espace, avec le paramètre
<<\,temporel\,>> $t =
1$, où l'on 
pose pour abréger $h ( x; \, \lambda) := 
h ( x; \, 1, \lambda )$. 
Géométriquement, on voit ainsi de quelle manière les équations de
transformations initiales $x_i' = h_i ( x; \,
\lambda_1, \dots, \lambda_r)$ agissent 
\emphasis{simultanément} sur les $r$-uples de points. 
Si on les développe en série entière
par rapport aux paramètres 
$\lambda_k$, ces transformations s'écrivent:
\def\theequation{5}\begin{equation}
\aligned
{x_i^{(\mu)}}'
=
x_i^{(\mu)}
+
&
\sum_{k=1}^r\,\lambda_k\cdot\xi_{ki}^{(\mu)}
+
\sum_{k,\,j}^{1\dots r}\,
\frac{\lambda_k\,\lambda_j}{1\cdot 2}\cdot
X_k^{(\mu)}\big(\xi_{ji}^{(\mu)}\big)
+
\cdots
\\
&
\ \ \ \ \ \ \ \ \ \ \ \ \
{\scriptstyle{(i\,=\,1\,\cdots\,n;\,\,\,\mu\,=\,1\,\cdots\,r)}},
\endaligned
\end{equation}
où nous avons bien sûr posé: $\xi_{ ki}^{ ( \mu)} := \xi_{ ki} ( x^{ (
\mu)})$ et $X_k^{ ( \mu)} := \sum_{ i=1}^n\, \xi_{ ki} ( x^{ ( \mu)})
\, \frac{ \partial }{\partial x_i^{ (\mu)}}$.
Une telle idée se révèlera fructueuse dans d'autres contextes, \cf~la
démonstration du Théorème~24 p.~\pageref{Theoreme-24} ci-dessous.

D'après le théorème p.~\pageref{theoreme-essentiel}, afin 
d'établir que les paramètres $\lambda_1, \dots, \lambda_r$ sont
essentiels, on doit seulement développer 
$x'$ en série entière par rapport à $x$ à l'origine:
\[
x_i'
=
\sum_{\alpha\in\N^n}\,
\mathcal{F}_\alpha^i(\lambda)\,
x^\alpha
\ \ \ \ \ \ \ \ \ \ \ \ \
{\scriptstyle{(i\,=\,1\,\dots\,n)}},
\]
et montrer que le rang générique de la matrice infinie des
coefficients $\lambda \longmapsto \big( \mathcal{ F}_\alpha^i (
\lambda) \big)_{
\alpha \in \N^n}^{ 1 \leqslant i \leqslant n}$ 
est le maximal possible, égal à $r$. De même et immédiatement, on
obtient le développement correspondant des équations de
transformations dans les espaces-copies:
\def\theequation{6}\begin{equation}
{x_i^{(\mu)}}'
=
\sum_{\alpha\in\N^n}\,
\mathcal{F}_\alpha^{i,(\mu)}(\lambda)\,(x^{(\mu)})^\alpha
\ \ \ \ \ \ \ \ \ \ \ \ \
{\scriptstyle{(i\,=\,1\,\dots\,n;\,\,\,\mu\,=\,1\,\dots\,r)}},
\end{equation}
avec, pour tout $\mu = 1, \dots, r$, les
\emphasis{mêmes} fonctions coefficients: 
\[
\mathcal{F}_\alpha^{i,(\mu)}(\lambda)
\equiv 
\mathcal{F}_\alpha^i
(\lambda)
\ \ \ \ \ \ \ \ \ \ \ \ \
{\scriptstyle{(i\,=\,1\,\dots\,n;\,\,\,\alpha\,\in\,\N^n;
\,\,\,\mu\,=\,1\,\dots\,r)}}. 
\]
Il en découle que le rang générique de la matrice infinie des
coefficients correspondante, qui n'est autre qu'une copie de $r$ fois
la même application $\lambda
\longmapsto \big( \mathcal{ F}_\alpha^i ( \lambda)
\big)_{ \alpha \in \N^n}^{ 1 \leqslant i \leqslant n}$, 
ni ne croît, ni ne décroit.

\smallskip
{\em 
Ainsi, les paramètres $\lambda_1, \dots, \lambda_r$ pour les équations
de transformations $x' = h ( x; \, \lambda)$ sont essentiels si et
seulement si ils sont essentiels pour les équations de transformations
diagonales ${ x^{ ( \mu )}} ' = h \big( x^{ ( \mu)}; \, \lambda\big)$,
$\mu = 1, \dots, r$, induites sur le produit de $r$ copies de l'espace
des $x_1,
\dots, x_n$. }

\smallskip
Donc il nous faut démontrer que le rang générique de la copie de $r$
matrices infinies de coefficients $\lambda
\longmapsto \big( \mathcal{ F}_\alpha^{ i, (\mu)} ( \lambda) \big)_{
\alpha \in \N^n}^{ 1 \leqslant i \leqslant n, \, \, 1 \leqslant \mu
\leqslant r}$ est égal à $r$. Nous
allons en fait établir que le rang en $\lambda = 0$ de cette même
application est déjà égal à $r$, ou, de manière équivalente, que la
matrice infinie {\em constante}:
\[
\Big(
\frac{\partial\mathcal{F}_\alpha^{i,(\mu)}}{\partial\lambda_k}(0)
\Big)_{1\leqslant k\leqslant r}^{
1\leqslant i\leqslant n,\,\alpha\in\N^n,\,
1\leqslant\mu\leqslant r}
\]
dont les $r$ lignes sont indexées par les dérivées partielles, est de
rang égal $r$.

Afin de préparer cette matrice infinie, si nous différentions les
développements~\thetag{ 5}\,\,---\,\,qui s'identifient à~\thetag{
6}\,\,---\,\,par rapport à $\lambda_k$ en $\lambda = 0$, et si nous
développons les coefficients de nos transformations infinitésimales:
\[
\xi_{ki}(x^{(\mu)})
=
\sum_{\alpha\in\N^n}\,
\xi_{ki\alpha}\cdot(x^{(\mu)})^\alpha
\ \ \ \ \ \ \ \ \ \ \ \ \
{\scriptstyle{(i\,=\,1\,\dots\,n;\,\,\,k\,=\,1\,\dots\,r;
\,\,\,\mu\,=\,1\,\dots\,r)}}
\]
par rapport aux puissances de $x_1, \dots, x_n$, 
nous obtenons une expression appropriée de cette
matrice: 
\[
\small
\aligned
\Big(
\frac{\partial\mathcal{F}_\alpha^{i,(\mu)}}{\partial\lambda_k}(0)
\Big)_{1\leqslant k\leqslant r}^{
1\leqslant i\leqslant n,\,\alpha\in\N^n,\,
1\leqslant\mu\leqslant r}
&
\equiv
\Big(
\big(
\xi_{ki\alpha}
\big)_{1\leqslant k\leqslant r}^{1\leqslant i\leqslant n,\,\alpha\in\N^n}
\,\,\cdots\,\,
\big(
\xi_{ki\alpha}
\big)_{1\leqslant k\leqslant r}^{1\leqslant i\leqslant n,\,\alpha\in\N^n}
\Big)
\\
&
=:
\Big(
{\sf T}^\infty\Xi(0)
\,\,\cdots\,\,
{\sf T}^\infty\Xi(0)
\Big).
\endaligned
\]
Comme nous l'avons dit, il suffit donc de démontrer que cette matrice
est de rang $r$. Visiblement, cette matrice s'identifie à $r$ copies
de la même matrice infinie ${\sf T}^\infty \Xi(0)$ de
tous les coefficients
de Taylor en $0$ de la matrice:
\[
\Xi(x)
:=
\left(
\begin{array}{ccc}
\xi_{11}(x) & \cdots & \xi_{1n}(x)
\\
\cdots & \cdots & \cdots
\\
\xi_{r1}(x) & \cdots & \xi_{rn}(x)
\end{array}
\right)
\]
des coefficients des transformations infinitésimales $X_k$. 
\`A présent, nous pouvons formuler un lemme
auxiliaire qui va nous permettre de conclure. 

\smallskip\noindent{\bf Lemme.}
{\em Soit $n\geqslant 1$, $q \geqslant 1$, $m \geqslant 1$ des entiers,
soit ${\sf x} \in \C^n$ et soit:
\[
{\sf A}({\sf x})
=
\left(
\begin{array}{ccc}
a_{11}({\sf x}) & \cdots & a_{1m} ({\sf x})
\\
\cdots & \cdots & \cdots
\\
a_{q1}({\sf x}) & \cdots & a_{qm}({\sf x})
\end{array}
\right)
\]
une matrice arbitraire $q \times m$ de fonctions analytiques:
\[
a_{ij}({\sf x})
=
\sum_{\alpha\in\N^n}\,a_{ij\alpha}\,{\sf x}^\alpha
\ \ \ \ \ \ \ \ \ \ \ \ \
{\scriptstyle{(i\,=\,1\,\dots\,q;\,\,\,j\,=\,1\,\dots\,m)}}
\]
qui sont toutes définies dans un voisinage fixé de l'origine dans
$\C^n$, et soit la matrice constante $q \times \infty$ de tous ses
coefficients de Taylor à l'origine:
\[
{\sf T}^\infty{\sf A}(0)
:=
\big(
a_{ij\alpha}
\big)_{1\leqslant i\leqslant q}^{1\leqslant j\leqslant m,\,\alpha\in\N^n}
\]
dont les $q$ lignes sont étiquetées par l'indice $i$. Alors
l'inégalité suivante entre rangs génériques est satisfaite:}
\[
\text{\small\sf rang}\,{\sf T}^\infty{\sf A}(0)
\geqslant
\text{\small\sf rang-générique}\,{\sf A}({\sf x}).
\]

\smallskip{\em Preuve.}
Ici, notre matrice infinie ${\sf T}^\infty {\sf A} ( 0)$ sera
considérée comme agissant par multiplication {\em à gauche} sur des
vecteurs {\em horizontaux} $u = (u_1, \dots, u_q)$, de telle sorte que
$u \cdot {\sf T}^\infty {\sf A} ( 0)$ est une matrice $\infty \times
1$, c'est-à-dire un vecteur horizontal infini. De manière similaire,
${\sf A} ( {\sf x})$ agira sur des vecteurs horizontaux de fonctions
analytiques $(u_1 ( {\sf x}), \dots, u_r ( {\sf x}))$.

Supposons que $u = (u_1, \dots, u_q) \in \C^q$ est un vecteur
quelconque dans le noyau de ${\sf T}^\infty {\sf A} ( 0)$, à savoir:
$0 = u \cdot {\sf T}^\infty {\sf A} ( 0)$, 
c'est-à-dire avec des indices:
\[
0
=
u_1\cdot a_{1j\alpha}
+\cdots+
u_q\cdot a_{qj\alpha}
\ \ \ \ \ \ \ \ \ \ \ \ \
{\scriptstyle{(j\,=\,1\,\dots\,m;\,\,\,\alpha\,\in\,\N^n)}};
\]
nous déduisons alors immédiatement, après avoir multiplié chaque telle
équation par ${\sf x}^\alpha$ et après avoir sommé sur tous les
$\alpha \in \N^n$ que:
\[
0
\equiv
u_1\cdot a_{1j}({\sf x})
+\cdots+
u_q\cdot a_{qj}({\sf x})
\ \ \ \ \ \ \ \ \ \ \ \ \
{\scriptstyle{(j\,=\,1\,\dots\,m)}},
\]
de telle sorte que le même vecteur constant $u = (u_1, \dots, u_q)$
satisfait aussi $0 \equiv u \cdot A ({\sf x})$. Il en découle que la
dimension du noyau de ${\sf T}^\infty {\sf A} ( 0)$ est inférieure ou
égale à la dimension du noyau de $A ( {\sf x})$ (en un point 
générique ${\sf x}$): 
ceci coïncide clairement avec l'inégalité entre rangs (génériques) écrite
ci-dessus. 
\qed\smallskip

Maintenant, pour tout $q = 1, 2, \dots, r$, nous voulons appliquer le
lemme avec la matrice ${\sf A} ({\sf x})$ égale à:
\[
\big(\Xi(x^{(1)})\,
\big(\Xi(x^{(2)})
\,\cdots\, 
\Xi(x^{(q)})
\big),
\] 
c'est-à-dire égale à: 
\[
\label{q-fold-extended-coefficient-matrix}
\Xi_q
\big(
\widetilde{\sf x}_q
\big)
:=
\left(
\begin{array}{cccccccccc}
\xi_{11}^{(1)} & \cdots & \xi_{1n}^{(1)} &
\xi_{11}^{(2)} & \cdots & \xi_{1n}^{(2)} &
\cdots\cdots &
\xi_{11}^{(q)} & \cdots & \xi_{1n}^{(q)}
\\
\cdots & \cdots & \cdots &
\cdots & \cdots & \cdots &
\cdots\cdots &
\cdots & \cdots & \cdots
\\ 
\xi_{r1}^{(1)} & \cdots & \xi_{rn}^{(1)} &
\xi_{r1}^{(2)} & \cdots & \xi_{rn}^{(2)} &
\cdots\cdots &
\xi_{r1}^{(q)} & \cdots & \xi_{rn}^{(q)}
\end{array}
\right),
\]
où nous avons abrégé: 
\[
\widetilde{\sf x}_q 
:= 
\big(x^{(1)},\dots, 
x^{(q)}\big).
\]

\smallskip\noindent{\bf Assertion.}
{\em C'est une conséquence du fait que $X_1, X_2, \dots, X_r$ sont
linéairement indépendants que pour tout $q = 1, 2 ,\dots, r$, on a:}
\[
\text{\small\sf rang-générique}\,
\Big(
\Xi\big(x^{(1)}\big)\ \ \
\Xi\big(x^{(2)}\big)\ \ \
\cdots\ \ \
\Xi\big(x^{(q)}\big)
\Big)
\geqslant q.
\]

{\em Preuve.}
En effet, pour $q = 1$, il est en premier lieu clair que:
\[
\text{\small\sf rang-générique}\,
\big(\Theta(x^{(1)})\big) 
\geqslant 1, 
\]
puisque l'un au moins des $\xi_{ ki} (x)$ ne s'annule pas identiquement.
\'Etablissons ensuite par récurrence que, aussi 
longtemps qu'ils restent $< r$, les rangs génériques augmentent au
moins d'une unité à chaque pas:
\[
\text{\small\sf rang-générique}
\Big(
\Xi_{q+1}
\big(\widetilde{\sf x}_{q+1}\big)
\Big)
\geqslant 
1
+
\text{\small\sf rang-générique}
\Big(
\Xi_q\big(\widetilde{\sf x}_q\big)
\Big),
\]
une propriété qui concluera visiblement la preuve de l'Assertion.

En effet, si au contraire les rangs génériques se stabilisaient
lorsqu'on passe de $q$ à $q+1$, tout en restant demeuraient $< r$,
alors localement au voisinage d'un point générique, fixé $\widetilde{
\sf x}_{ q+1}^0$, les deux matrices $\Xi_{ q+1}$ et $\Xi_q$ auraient
le même rang, localement constant. Par conséquent, les solutions
$\big(
\vartheta_1 ( \widetilde{ \sf x}_q) \, \, \cdots \, \, \vartheta_r (
\widetilde{ \sf x}_q ) \big)$
au système d'équations linéaires écrit sous forme matricielle:
\[
0
\equiv
\big(
\vartheta_1(\widetilde{\sf x}_q)
\,\,\cdots\,\,
\vartheta_r(\widetilde{\sf x}_q)
\big)
\cdot
\Xi_q\big(\widetilde{\sf x}_q\big),
\]
lesquelles
sont analytiques près de $\widetilde{ \sf
x}_q^0$\,\,---\,\,grâce à une application de la règle de Cramer et
grâce à la constance du rang\,\,---\,\, seraient automatiquement
solutions du système étendu:
\[
0
\equiv
\big(
\vartheta_1(\widetilde{\sf x}_q)
\,\,\cdots\,\,
\vartheta_r(\widetilde{\sf x}_q)
\big)
\cdot
\big(
\Xi_q(\widetilde{\sf x}_q)\,\,\,\,
\Xi(x^{(q+1)})
\big),
\]
et donc il existerait des solutions \emphasis{non nulles} 
$(\vartheta_1, \dots, \vartheta_r)$ aux équations de dépendance 
linéaire:
\[
0
=
\big(
\vartheta_1\,\,\cdots\,\,\vartheta_r
\big)\cdot
\Xi
\big(x^{(q+1)}\big)
\]
qui seraient \emphasis{constantes} par rapport à la variable $x^{ (
q+1)}$, puisqu'elles dépendent seulement de $\widetilde{ \sf
x}_q$. Ceci contredirait précisément l'hypothèse que $X_1^{ (q+1)},
\dots, X_r^{ (q+1)}$ sont linéairement indépendants.
\qed\smallskip

Pour terminer, nous pouvons enchaîner une série d'inégalités
qui sont maintenant des conséquences évidentes du Lemme
et de l'Assertion: 
\[
\aligned
\text{\small\sf rang}\,
\Big(
{\sf T}^\infty\Xi(0)
\,\,\cdots\,\,
{\sf T}^\infty\Xi(0)
\Big)
&
=
\text{\small\sf rang}\,
{\sf T}^\infty\Xi_r(0)
\\
&
\geqslant
\text{\small\sf rang-générique}\,
\Xi_r\big(\widetilde{\sf x}_r\big)
\geqslant
r,
\endaligned
\]
et puisque tous ces rangs (génériques)
sont en tout cas $\leqslant r$, nous
obtenons l'estimation promise:
\[
r
=
{\rm rank}\,\Big(
{\sf T}^\infty\Xi(0)
\,\,\cdots\,\,
{\sf T}^\infty\Xi(0)
\Big),
\]
ce qui achève finalement la démonstration du théorème. 
\qed\smallskip

Afin de mémoriser le prolongement des transformations au produit de
$r$ copies de l'espace des $x_1, \dots, x_n$, formulons une
proposition qui sera utilisée dans la démonstration du Théorème~24
p.~\pageref{Theoreme-24}.

\smallskip\noindent{\bf Proposition.}
\label{Proposition-p-66}
{\em Si les $r$ transformations infinitésimales:
\[
X_k(f)
=
\sum_{i=1}^n\,\xi_{ki}(x_1,\dots,x_n)\,
\frac{\partial f}{\partial x_i}
\ \ \ \ \ \ \ \ \ \ \ \ \
{\scriptstyle{(k\,=\,1\,\dots\,r)}}
\]
sont linéairement indépendantes et si de plus:
\[
x_1^{(\mu)},\dots,x_n^{(\mu)}
\ \ \ \ \ \ \ \ \ \ \ \ \
{\scriptstyle{(\mu\,=\,1\,\dots\,r)}}
\]
sont $r$ systèmes distincts de $n$ variables, et si enfin on pose pour
abréger:
\[
X_k^{(\mu)}(f)
=
\sum_{i=1}^n\,
\xi_{ki}
\big(
x_1^{(\mu)},\dots,x_n^{(\mu)}
\big)\,
\frac{\partial f}{\partial x_i^{(\mu)}}
\ \ \ \ \ \ \ \ \ \ \ \ \
{\scriptstyle{(k,\,\mu\,=\,1\,\dots\,r)}},
\]
alors les $r$ transformations infinitésimales:
\[
W_k(f)
=
\sum_{\mu=1}^r\,X_k^{(\mu)}(f)
\ \ \ \ \ \ \ \ \ \ \ \ \
{\scriptstyle{(k\,=\,1\,\dots\,r)}}
\]
en les $nr$ variables $x_i^{ ( \mu)}$ ne satisfont
\emphasis{aucune} relation de la forme:}
\[
\sum_{k=1}^n\,
\chi_k
\big(
x_1^{(1)},\dots,x_n^{(1)},\dots\dots,
x_1^{(r)},\dots,x_n^{(r)}
\big)\cdot
W_k(f)
\equiv
0.
\]

}

\HEAD{Chapitre~3.\,\,\,\,Théorèmes fondamentaux sur les groupes
de transformations}{ 3.7.\,\,\,Le théorème de Clebsch-Frobenius}

\noindent{\bf 3.7.~Le théorème de Clebsch-Lie-Frobenius.}
\label{theoreme-Clebsch-Lie-Frobenius}
{\em Question préliminaire:} Quelles sont les solutions générales
$\omega$ à une équation scalaire aux dérivées partielles du premier
ordre $X \omega = 0$ associée à un champ de vecteurs $X = \sum_{
i=1}^n\, \xi_i ( x) \, \frac{ \partial }{\partial x_i}$ à 
coefficients analytiques\,?

Les relocalisations au voisinage d'un point générique étant
autorisées, nous pouvons supposer, après une renumérotation éventuelle
des variables, que $\xi_n$ ne s'annule pas en un point que nous
choisissons comme origine d'un système de coordonnées $(x_1, \dots,
x_n)$. En divisant par $\xi_n (x)$, il est équivalent de rechercher
les fonctions $\omega$ qui sont annihilées par l'opérateur
différentiel:
\[
X
=
\sum_{i=1}^{n-1}\,
\frac{\xi_i(x)}{\xi_n(x)}\,
\frac{\partial}{\partial x_i}
+
\frac{\partial}{\partial x_n},
\]
toujours noté $X$ et qui satisfait maintenant $X ( x_n) \equiv 1$.
Rappelons que le système d'équations différentielles ordinaires qui
définit les courbes intégrales de ce champ $X$, à savoir le système:
\[
\frac{d{\sf x}_1}{dt}
=
\frac{\xi_1\big({\sf x}(t)\big)}{\xi_n\big({\sf x}(t)\big)},\ \
\dots\dots,\ \ \ \
\frac{d{\sf x}_{n-1}}{dt}
=
\frac{\xi_{n-1}\big({\sf x}(t)\big)}{\xi_n\big({\sf x}(t)\big)},\ \ \ \ 
\frac{d{\sf x}_n(t)}{dt}
=
1,
\]
avec la condition initiale pour $t = 0$: 
\[
{\sf x}_1(0)
=
x_1,\dots\dots,
{\sf x}_{n-1}(0)
=
x_{n-1},\ \ \
{\sf x}_n(0)
=
0,
\]
est {\em résoluble}; plus précisément, il possède une unique solution
$( {\sf x}_1 (t), \dots, {\sf x}_{ n-1} (t), {\sf x}_n (t))$ qui est
analytique dans un voisinage de l'origine. En vérité, nous connaissons
déjà la technique de résolution. Tout d'abord, par une intégration
évidente, on a ${\sf x}_n (t) = t$; ensuite, les $(n-1)$ autres
fonctions ${\sf x}_k (t)$ sont données par la merveilleuse formule
exponentielle (Proposition p.~\pageref{flow-exponential-formula}):
\[
{\sf x}_k(t)
=
\exp(tX)(x_k)
=
\sum_{l\geqslant 0}\,
\frac{t^l}{l!}\,
X^l(x_k)
\ \ \ \ \ \ \ \ \ \ \ \ \
{\scriptstyle{(k\,=\,1\,\cdots\,n-1)}}.
\]
Posons alors $t := -\, x_n$ dans cette formule (le signe <<\,moins\,>> va
être crucial), et définissons les $(n-1)$ fonctions qui seront
intéressantes:
\def\theequation{1}\begin{equation}
\aligned
\omega_k(x_1,\dots,x_n)
:=
{\sf x}_k(-x_n)
&
=
\exp\big(-x_nX\big)(x_k)
\\
&
=
\sum_{l\geqslant 0}\,
(-1)^l\,
\frac{(x_n)^l}{l!}\,
X^l(x_k)
\ \ \ \ \ \ \ \ \ \ \ \ \
{\scriptstyle{(k\,=\,1\,\cdots\,n\,-\,1)}}.
\endaligned
\end{equation}

\noindent{\bf Proposition.}
\label{equation-ordre-1}
{\em Les $(n-1)$ fonctions ainsi définies $\omega_1, \dots, \omega_{
n-1}$ sont des solutions de l'équation aux dérivées partielles $X
\omega = 0$. 
De plus, ces solutions sont {\em fonctionnellement indépendantes}, au
sens où le rang de leur matrice jacobienne $\big(
\frac{ \partial
\omega_k}{ \partial x_i} \big)_{ 1 \leqslant i \leqslant n}^{ 1
\leqslant k \leqslant n-1}$ 
est égal à $n-1$ à l'origine. Enfin, pour toute autre solution
$f$ de $X f = 0$, il existe une fonction analytique locale
$F = F \big( \omega_1, \dots, \omega_{ n-1} \big)$
définie au voisinage de l'origine dans $\C^{ n-1}$ telle que:}
\[
f(x)
\equiv
F\big(\omega_1(x),\dots,\omega_{n-1}(x)\big).
\]

{\em Preuve.}
En effet, lorsqu'on applique $X$ à la série ci-dessus qui définit les
fonctions $\omega_k$, on observe que tous les termes s'annihilent
grâce à une simple application élémentaire de la formule de Leibniz,
développée sous la forme:
\[
X\big[
(x_n)^l\,X^l(x_k)
\big]
=
l\,(x_n)^{l-1}\,
X^l(x_k)
+
(x_n)^l\,X^{l+1}(x_k).
\]
Ensuite, l'assertion d'après laquelle l'application:
\[
x \longmapsto
\big(\omega_1(x),\dots,\omega_{n-1}(x)\big)
\]
est de rang $n-1$ en $x = 0$ provient tout simplement du fait que l'on
a par construction:
\[
\omega_k(x_1,\dots,x_{n-1},0) 
\equiv 
x_k. 
\]

Enfin, après redressement ({\em voir} le Théorème
p.~\pageref{redressement-champ}) de $X$ en $X' :=
\frac{ \partial }{ \partial x_n'}$ 
dans certaines nouvelles coordonnées appropriées
$(x_1', \dots, x_n')$, la solution générale $f ' ( x')$ de
l'équation aux dérivées partielles (redressée):
\[
X'(f') 
= 
{\textstyle{\frac{\partial f'}{\partial x_n'}}}
=
0
\]
s'avère alors trivialement être une fonction arbitraire: 
\[
F'\big(x_1',\dots,x_{n-1}'\big)
\]
des $(n-1)$ premières variables $x_1', \dots, x_{ n-1}'$, lesquelles
s'identifient, dans ce système de coordonnées, aux fonctions
$\omega_1' \equiv x_1', \dots, \omega_{ n-1}' \equiv x_{ n-1}'$ 
définies par la formule~\thetag{ 1}. Ceci montre que dans l'ancien
système de coordonnées $(x_1, \dots, x_n)$, la solution générale de
$Xf = 0$ s'écrit comme annoncé: $f =
F ( \omega_1, \dots, \omega_{ n-1})$.
\qed\smallskip

{\em Question:} Qu'advient-il en présence de {\em plusieurs} équations
aux dérivées partielles? La manière dont Engel et Lie présentent la
résolution de cette question dans le Chapitre~5 de~\cite{ enlie1888}
est symptomatique sur le plan du {\em contrôle par la pensée des
principes de genèse}, et nous nous proposons à présent d'en restituer
la teneur.

Premier principe: examiner la dyade comme germe du Multiple pur. Par
exemple, si une fonction $f( x_1, \dots, x_n)$ satisfait {\em deux}
équations aux dérivées partielles du premier
ordre:
\[
X_1(f)
=
0,
\ \ \ \ \ \ \ \ \
X_2(f)
=
0,
\]
alors elle satisfait naturellement aussi les deux équations
différentielles du second ordre:
\[
X_1\big(X_2(f)\big)
=
0,
\ \ \ \ \ \ \ \ \ 
X_2\big(X_1(f)\big)
=
0,
\]
et par conséquent, aussi l'équation: 
\[
X_1\big(X_2(f)\big)
-
X_2\big(X_1(f)\big)
=
0,
\]
qui est obtenue en soustrayant ces deux équations. Or si l'on
introduit les expressions développées de ces deux opérateurs d'ordre
$1$:
\[
X_k
=
\sum_{i=1}^n\,
\xi_{ki}(x_1,\dots,x_n)\,
\frac{\partial}{\partial x_i}
\ \ \ \ \ \ \ \ \ \ \ \ \
{\scriptstyle{(k\,=\,1,\,2)}},
\]
alors cette dernière équation ne dérive la fonction $f$ qu'à l'ordre
$1$:
\[
X_1\big(X_2(f)\big)
-
X_2\big(X_1(f)\big)
=
\sum_{i=1}^n\,
\big[
X_1(\xi_{2i})
-
X_2(\xi_{1i})
\big]\,
\frac{\partial f}{\partial x_i},
\]
puisque tous les termes qui incorporent des dérivées partielles du
second ordre s'annihilent dans la soustraction. Grâce à ce procédé,
on obtient alors un nouvel opérateur qui a la même structure que les
deux opérateurs initiaux: {\em homologie de l'ontologie}. 
On notera alors:
\[
\big[X_1,\,X_2\big]
=
-
\big[X_2,\,X_1\big]
\]
cet opérateur que l'on appellera 
{\sl crochet}\footnote{\,
%%%%%%%%%%%%%%%%%%%%%%%-------DEBUT--------%%%%%%%%%%%%%%%%%%%%%%%%%%%
Hawkins~\cite{ h2001} relate la réinterprétation de ce concept dans
les premiers travaux de Lie sur l'intégration des systèmes d'équations
aux dérivées partielles. Dans la {\em Theorie der
Transformationsgruppen}, c'est peu fréquemment que Engel et Lie
nomment le fameux <<\,{\sl crochet}\,>> souvent considéré comme une
étape de calcul; en certains endroits, ils l'appellent <<\,{\sl
combinaison}\,>>
\deutsch{Combination} (entre deux transformations
infinitésimales), ou simplement <<\,{\sl équation}\,>>
\deutsch{Gleichung}, et
à la fin du traité, ils utilisent en général la terminologie
\deutsch{Klammeroperation}, <<\,{\sl opération de
crochet}\,>>, ou <<\,de parenthésage\,>>. Il est noté $( X_1 X_2
\big)$, toujours avec des parenthèses simples, sans le symbole de
fonction $f$, et sans virgule, sauf quand les deux éléments insérés
sont des champs de vecteurs en coordonnées.
} %%%%%%%%%%%%%%%%%%%%%%%%-----FIN-----%%%%%%%%%%%%%%%%%%%%%%%%%%%%%%%
(de Jacobi ou de Lie) entre $X_1$ et $X_2$ et qui est bien s\^ur
antisymétrique par rapport à ses deux arguments. Toute solution $f$
de $X_1f = X_2 f = 0$ est alors solution de $\big[ X_1, \, X_2 \big] f
= 0$. L'opérateur $\big[ X_1, \, X_2 \big]$ s'ajoute alors en quelque
sorte gratuitement aux équations initiales. C'est donc un deuxième
principe: {\em engendrement du tiers par antisymétrisation de la dyade}.
En résumé, on a une: 

\smallskip\noindent{\bf Observation fondamentale.}
{\em Si une fonction $\psi ( x_1, \dots, x_n)$ 
satisfait les deux équations aux dérivées partielles
du premier ordre:
\[
X_k(f)
=
\sum_{i=1}^n\,\xi_{ki}(x_1,\dots,x_n)\,
\frac{\partial f}{\partial x_i}
=
0
\ \ \ \ \ \ \ \ \ \ \ \ \
{\scriptstyle{(k\,=\,1,\,2)}},
\]
alors elle satisfait aussi l'équation
différentielle du premier ordre:}
\[
X_1\big(X_2(f)\big)
-
X_2\big(X_1(f)\big)
=
\sum_{i=1}^n\,
\big[
X_1(\xi_{2i})
-
X_2(\xi_{1i})
\big]\,
\frac{\partial f}{\partial x_i}
=
0.
\]

{\em Question suscitée [dans l'arbre de scindage de
l'irréversible-syn\-thé\-tique]}: l'opérateur $\big[ X_1, \, X_2
\big]$ apporte-t-il, ou n'apporte-t-il pas d'information nouvelle? En
admettant les principes de pensée que nous avons formulés
p.~\pageref{Prologue} sq., la réponse est simple. S'il existe deux
fonctions (analytiques) $\chi_1 ( x)$ et $\chi_2 ( x)$ telles qu'on
peut réécrire (après relocalisation éventuelle):
\[
\big[X_1,\,X_2\big]
=
\chi_1\,X_1+\chi_2\,X_2,
\]
alors le fait que toute solution $f$ de $X_1 f = X_2 f = 0$ est aussi
solution de $\big[ X_1, \, X_2 \big] f = 0$ sera une conséquence
triviale des hypothèse et n'apportera aucune connaissance nouvelle,
puisque de $X_1 f = X_2 f = 0$ on déduit immédiatement par des
opérations algébriques {\em non différentielles}:
\[
\chi_1X_1f
=
\chi_2X_2f
=
0
\ \ \ \ \ \ \ \
\text{\rm puis}
\ \ \ \ \ \ \ \
\chi_1X_1f+\chi_2X_2f=0. 
\]
Au contraire, lorsque le crochet $\big[ X_1, \, X_2 \big]$ ne peut pas
être exprimé comme combinaison linéaire de $X_1$ et de $X_2$ (même
après relocalisation en un point générique), l'équation $\big[ X_1, \,
X_2 \big] f = 0$ doit être considérée comme {\em nouvelle et
nécessaire}. Ce cas se produit par exemple avec $X_1 = \frac{
\partial}{ \partial x_1}$ et $X_2 = \frac{ \partial}{ \partial x_2} +
x_1 \, \frac{ \partial}{ \partial x_3}$ d'où $\big[ X_1, \, X_2\big] =
\frac{ \partial }{ \partial x_3}$. 

Ainsi, en partant de deux opérateurs distincts, la considération d'un
troisième opérateur peut s'avérer incontournable. Autrement dit, le
principe de différenciation par examen de la dyade force à envisager
le multiple général.

Considérons donc maintenant {\em d'emblée} un nombre quelconque $q
\geqslant 2$ d'opérateurs d'ordre $1$ à coefficients analytiques:
\[
X_k
=
\sum_{i=1}^n\,
\xi_{ki}(x_1,\dots,x_n)\,
\frac{\partial}{\partial x_i}
\ \ \ \ \ \ \ \ \ \ \ \ \
{\scriptstyle{(k\,=\,1\,\cdots\,q)}}. 
\]
Tout d'abord, il peut se produire qu'il existe des relations
de dépendance de la forme:
\[
\sum_{k=1}^q\,\chi_k(x)\cdot
X_k
\equiv
0,
\]
où les $\chi_k ( x)$ sont des fonctions analytiques non toutes nulles.
Après relocalisation et renumérotation éventuelle, on peut résoudre
une telle équation sous la forme: $X_q = \tau_1 X_1 + \cdots + \tau_{
q-1} X_{ q-1}$. Si une telle équation résolue non triviale existe,
alors parmi les $q$ équations $X_1 f = \cdots = X_{ q-1} f = X_q f =
0$, la dernière sera manifestement conséquence des $(q-1)$ premières,
et elle pourra donc d'ores et déjà être laissée de côté. Aussi est-il
parfaitement légitime, lorsqu'on veut résoudre les équations $X_k f =
0$, de supposer qu'il n'existe pas de telle relation de
dépendance. Après relocalisation éventuelle et renumérotation
éventuelle des variables, cela revient, d'après un théorème d'algèbre
linéaire, à admettre que les opérateurs $X_1, \dots, X_q$ sont
résolubles par rapport aux $q$ quotients différentiels $\frac{
\partial f }{ \partial x_i}$, $i = 1, \dots, q$, autrement dit,
qu'il existe une matrice $q \times q$ de fonctions analytiques
$\varpi_{ j k} ( x)$, inversible dans un ouvert relocalisé, telle que
les nouveaux 
opérateurs $Y_j := \sum_{ k=1 }^q \, \varpi_{j k}(x)\, X_k$ sont
de la forme normalisée:
\[
Y_j
=
\frac{\partial}{\partial x_j}
+
\sum_{q+1\leqslant i\leqslant n}\,
\theta_{ji}(x)\,
\frac{\partial}{\partial x_i}
\ \ \ \ \ \ \ \ \ \ \ \ \
{\scriptstyle{(j\,=\,1\,\cdots\,q)}}. 
\]
Bien entendu, l'étude du sytème $X_1 f = \cdots = X_q f = 0$ se ramène
à celle du système $Y_1 f = \cdots = Y_q f = 0$, puisque
le déterminant de la matrice $\varpi_{ j k} ( x)$ ne s'annule
en aucun point de l'ouvert relocalisé. 

D'après l'observation fondamentale, les solutions communes possibles
aux $q$ équations $X_1f = \cdots = X_q f = 0$ satisfont aussi toutes
les équations par paires de la forme:
\[
X_i\big(X_k(f)\big)
-
X_k\big(X_i(f)\big)
=
0
\ \ \ \ \ \ \ \ \ \ \ \ \
{\scriptstyle{(i,\,k\,=\,1\,\cdots\,q)}}.
\]
Et maintenant, deux circonstances distinctes peuvent se produire.

Premièrement, chacune des 
$\frac{ q ( q - 1)}{ 2}$ équations ainsi obtenues
peut s'avérer être conséquence des précédentes. Tel est le cas si et
seulement si, pour tout $i \leqslant q$ et pour tout $k \leqslant q$,
une relation de dépendance de la forme:
\[
X_i\big(X_k(f)\big)
-
X_k\big(X_i(f)\big)
=
\chi_{ik1}(x)\,X_1(f)
+\cdots+
\chi_{ikq}(x)\,X_q(f)
\]
est satisfaite. {\sc Clebsch} dit alors les $q$ équations
indépendantes $X_1 (f) = \cdots = X_q ( f) = 0$ forment {\sl
système complet}.

Mais en général, c'est le second cas, plus délicat, qui se produit.
Parmi les nouvelles équations:
\[
X_i\big(X_k(f)\big)
-
X_k\big(X_i(f)\big)
=
0,
\]
un certain nombre, disons $s \geqslant 1$, seront indépendantes des
$q$ équations $X_1 f = \cdots = X_q f = 0$. Ajoutons alors ces
nouvelles équations aux $q$ équations initiales, notons-les: 
\[
X_{q+1}(f)=0,
\dots\dots,
X_{q+s}(f)
=
0,
\]
et traitons maintenant le système obtenu de ces $q+s$ équations
exactement comme nous avons traité précédemment les $q$ équations de
départ. Bien entendu, la relocalisation autour d'un point générique
est toujours admise. Comme on ne peut pas obtenir plus de $n$
équations $X_i (f) = 0$ qui sont indépendantes l'une de l'autre en un
point générique, on doit aboutir, au bout d'un nombre fini d'étapes, à
un système complet qui consiste en un nombre $q \leqslant n$
d'équations indépendantes.

\smallskip\noindent{\bf Proposition.}
(\cite{ enlie1888}, p.~86)
{\em 
La détermination des solutions communes de $q$ équations linéaires aux
dérivées partielles du premier ordre $X_1 ( f) = \cdots = X_q ( f)$
peut toujours être ramenée, par différentiation et élimination
algébrique linéaire, à l'intégration d'un système complet
d'équations indépendantes. 
}\medskip

On peut donc supposer maintenant sans perte de généralité que le
système à étudier $X_1 f = \cdots = X_q f = 0$ est complet et qu'il
est formé de $q$ équations indépendantes.

\smallskip\noindent{\bf Proposition.}
(\cite{ enlie1888}, pp.~86--87)
{\em 
Si l'on résout un système complet et indépendant de $q$ équations:
\[
X_1(f)=0,\dots\dots,
X_q(f)=0
\]
par rapport à $q$ quotients différentiels, disons $\frac{ \partial}{
\partial x_{ n-q+1}}, \dots, \frac{ \partial }{ \partial x_n}$ après
renumérotation éventuelle des variables, alors les $q$ équations
obtenues:
\def\theequation{4}\begin{equation}
Y_k(f)
=
\frac{\partial f}{\partial x_{n-q+k}}
+
\sum_{i=1}^{n-q}\,
\eta_{ki}\,
\frac{\partial f}{\partial x_i}
=
0
\ \ \ \ \ \ \ \ \ \ \ \ \
{\scriptstyle{(k\,=\,1\,\cdots\,q)}}
\end{equation}
satisfont les relations de commutation par paires:}
\def\theequation{5}\begin{equation}
Y_j\big(Y_k(f)\big)
-
Y_k\big(Y_j(f)\big)
=
0
\ \ \ \ \ \ \ \ \ \ \ \ \
{\scriptstyle{(j,\,k\,=\,1\,\cdots\,q)}}.
\end{equation}

\smallskip{\em Preuve.}
En effet, observons que l'expression semi-développée:
\[
Y_j\big(Y_k(f)\big)
-
Y_k\big(Y_j(f)\big)
=
\sum_{i=1}^{n-q}\,
\big[
Y_j(\eta_{ki})
-
Y_k(\eta_{ji})
\big]
\frac{\partial f}{\partial x_i}
\ \ \ \ \ \ \ \ \ \ \ \ \
{\scriptstyle{(j,\,k\,=\,1\,\cdots\,q)}}
\]
ne fait intervenir aucun des quotients différentiels $\frac{
\partial}{ \partial x_{ n-q+1}}, \dots, \frac{ \partial 
}{ \partial x_n}$, puisque la dérivation $Y_j ( 1)$ du coefficient
constant égal à $1$ de la première dérivation $\frac{\partial
f}{\partial x_{n-q+k}}$ de $Y_k$ s'annule trivialement. Mais alors une
relation de dépendance linéaire de la forme:
\[
Y_j\big(Y_k(f)\big)
-
Y_k\big(Y_j(f)\big)
=
\sum\,{\sf coeff}\cdot
\bigg(
\frac{\partial f}{\partial x_{n-q+k}}
+
\sum_{i=1}^{n-q}\,
\eta_{ki}\,
\frac{\partial f}{\partial x_i}
\bigg)
\]
ne peut manifestement être possible que si
tous les coefficients présents s'annulent. 
\qed\smallskip

La généralisation conjointe du Théorème
p.~\pageref{redressement-champ} et de la Proposition
p.~\pageref{equation-ordre-1} à un système complet de $q$ équations
indépendantes s'énonce alors comme suit\footnote{\,
%%%%%%%%%%%%%%%%%%%%%%%-------DEBUT--------%%%%%%%%%%%%%%%%%%%%%%%%%%%
Pour la démonstration détaillée de ce théorème standard de calcul
différentiel aujourd'hui dit <<\,{\sl de Frobenius}\,>>, outre~\cite{
enlie1888}, on pourra consulter~\cite{ kono1963, spiv1970, sha1997,
stk2000, ha2005}, ou bien compléter les arguments en s'inspirant des
raisonnements qui précèdent.
}. %%%%%%%%%%%%%%%%%%%%%%%%-----FIN-----%%%%%%%%%%%%%%%%%%%%%%%%%%%%%%% 

\smallskip\noindent{\bf Théorème.}
\label{12}
\label{Theoreme-12}
{\sc (Clebsch-Lie-Frobenius)}
{\em Tout système complet formé de $q$ équations résolues:
\[
\frac{\partial f}{\partial x_{n-q+k}}
+
\sum_{i=1}^{n-q}\,
\eta_{ki}(x_1,\dots,x_n)\,
\frac{\partial f}{\partial x_i}
=
0
\ \ \ \ \ \ \ \ \ \ \ \ \
{\scriptstyle{(k\,=\,1\,\cdots\,q)}}
\]
dont les coefficients $\eta_{ ki}$ sont analytiques au voisinage d'un
point $(x_1^0, \dots, x_n^0)$ possède $n - q$ solutions indépendantes
$\omega_1 (x), \dots, \omega_{ n-q} ( x)$ qui sont analytiques dans un
certain voisinage $(x_1^0, \dots, x_n^0)$ et qui, de plus, se
réduisent à $x_1, \dots, x_{ n-q}$ lorsqu'on pose $x_{ n-q + 1} = x_{
n-q + 1}^0$, \dots, $x_n = x_n^0$.

De plus, pour toute autre solution $f$ de $X_1 f = \cdots = X_q f =
0$, il existe une fonction analytique locale $F = F ( \omega_1, \dots,
\omega_{ n-q})$ définie au voisinage de $(x_1^0, \dots, x_{ n-q}^0)$
dans $\C^{ n-q}$ telle que:}
\[
f(x)
\equiv
F
\big(\omega_1(x),\dots,\omega_{n-q}(x)\big). 
\]

D'après Lie (\cite{ enlie1888}, p.~91), ce théorème central de la
théorie des systèmes complet n'avait pas été énoncé d'une manière
aussi précise, ni par Jacobi, ni par Clebsch, bien qu'il soit
implicitement contenu dans des mémoires de Cauchy, de Weierstrass, de
Briot et Bouquet, de Kowalevsky et de Darboux consacrés à l'existence
des solutions de systèmes d'équations aux dérivées partielles.
Frobenius (\cite{ frob1877}) en fera la synthèse 
finale\footnote{\,
%%%%%%%%%%%%%%%%%%%%%%%-------DEBUT--------%%%%%%%%%%%%%%%%%%%%%%%%%%%
{\em
Voir}~\cite{ ha2005} pour une excellente mise en perspective
historique et philosophique.
}. %%%%%%%%%%%%%%%%%%%%%%%%-----FIN-----%%%%%%%%%%%%%%%%%%%%%%%%%%%%%%%

\HEAD{Chapitre~3.\,\,\,\,Théorèmes fondamentaux sur les groupes
de transformations}{
3.8.\,\,\,\'Equations de structure}

\medskip\noindent{\bf 3.8.~Constantes de structure et correspondance
fondamentale.}
\label{constantes-de-structure} 
Expliquons maintenant dans le détail comment, au Chapitre~9 de~\cite{
enlie1888}, Engel et Lie organisent\footnote{\,
%%%%%%%%%%%%%%%%%%%%%%%-------DEBUT--------%%%%%%%%%%%%%%%%%%%%%%%%%%%
Dans ce paragraphe~3.8, nous traduisons en l'adaptant librement la
majeure partie du Chapitre~9 de~\cite{ enlie1888}. Les théorèmes
fondamentaux de la théorie y apparaissent déployés dans une pensée
beaucoup plus systématique que ce que la postérité en a retenu.  En
particulier, la classification en trois {\em Théorèmes fondamentaux}
telle qu'établie après la rédaction du premier volume, à savoir: 1)
l'existence d'équations différentielles fondamentales; 2) l'existence
de constantes dans les crochets entre générateurs infinitésimaux; 3)
la reconstition d'un groupe de Lie local à partir d'un système de
constantes satisfaisant les identités algébriques naturelles
correspondant à l'antisymétrie et à l'identité de Jacobi
(\cf~Schur, Scheffers, Cartan, Campbell, Bianchi et d'autres) 
simplifie l'exposition d'origine, toute entière orientée 
vers l'exploration spéculative des axiomes fondamentaux. 
}
%%%%%%%%%%%%%%%%%%%%%%%%-----FIN-----%%%%%%%%%%%%%%%%%%%%%%%%%%%%%%%
la présentation de ces théorèmes fondamentaux, eu égard à cette idée
fixe de la théorie commençante: économiser à la fois l'axiome de
l'élément identité et l'axiome des éléments 
inverses. Filigrane de
tension métaphysique: il s'agit d'engendrer la théorie du continu en
complète analogie harmonique avec la théorie discrète des groupes de
substitutions.

Soit donc une famille d'équations de transformations à $r$ paramètres
essentiels:
\def\theequation{1}\begin{equation}
x_i'
=
f_i\big(x_1,\dots,x_n;\,a_1,\dots,a_r\big)
\ \ \ \ \ \ \ \ \ \ \ \ \
{\scriptstyle{(i\,=\,1\,\cdots\,n)}}. 
\end{equation}
Ici, les variables $(x_1, \dots, x_n)$ varient dans un 
domaine\footnote{\,
%%%%%%%%%%%%%%%%%%%%%%%-------DEBUT--------%%%%%%%%%%%%%%%%%%%%%%%%%%%
Nous considérerons en effet ici explicitement les domaines d'existence. 
} %%%%%%%%%%%%%%%%%%%%%%%%-----FIN-----%%%%%%%%%%%%%%%%%%%%%%%%%%%%%%%
$\mathcal{ X} \subset \C^n$, les paramètres $(a_1, \dots, a_r)$
varient dans un domaine $\mathcal{ A} \subset \C^r$, et l'application
$x \mapsto f_a ( x) = f ( x; \, a)$ est, pour tout $a$, un
difféomorphisme de $\mathcal{ X}$ sur son image $f_a ( \mathcal{ X})$.
Comme cela a déjà été incidemment signalé dans les deux notes
p.~\pageref{ferme-par-composition} et p.~\pageref{sans-inverse}, il se
trouve que l'existence d'équations différentielles fondamentales
(\cf~le Théorème p.~\pageref{Theoreme-3}) peut être dérivée de la seule
condition que les transformations du groupe sont stables par
composition au sens technique suivant:
\label{debut-38}
\[
\left[
\aligned
&
f\big(f(x;\,a);\,b\big)
\equiv
f\big(x;\,\varphi(a,b)\big)
\ \ \ \ \ \ \ \ \
\text{\rm pour tous}\ \
x\in\mathcal{X}^1,\ \
a\in\mathcal{A}^1,\ \
b\in\mathcal{A}^1,
\\
&
c
\equiv
\varphi\big(a,\,\chi(a,c)\big)
\ \ \ \ \ \ \ \ \ \ \ \ \ \ \ \ \ \ \ \ \ \ \ \ \ \ \ \ \ \,
\text{\rm pour tous}\ \
a\in\mathcal{A}^1,\ \
c\in\mathcal{A}^1,
\endaligned\right.
\]
sans supposer ni l'existence d'un élément identité, ni l'existence de
transformations inverses l'une de l'autre par paires. Plus
précisément, sous les hypothèses spécifiques de la note
p.~\pageref{ferme-par-composition}, on démontre la proposition
suivante en s'inspirant des raisonnements du \S~3.5.

\smallskip\noindent{\bf Proposition.}
(\cite{ enlie1888}, pp.~33--34)
{\em 
Il existe une matrice $( \psi_{ kj} ( a))_{ 1 \leqslant k \leqslant
r}^{ 1 \leqslant j \leqslant r}$ de taille $r \times r$ de fonctions
qui sont holomorphes et inversibles dans $\mathcal{ A}^1$, et il
existe des fonctions $\xi_{ j i} ( x)$ holomorphes dans $\mathcal{ X}$
telles que les équations différentielles:
\def\theequation{2}\begin{equation}
\frac{\partial x_i'}{\partial a_k}(x;\,a)
=
\sum_{j=1}^r\,\psi_{kj}(a)\cdot\xi_{j i}(x')
\ \ \ \ \ \ \ \ \ \ \ \ \
{\scriptstyle{(i\,=\,1\,\cdots\,n\,;\,\,\,k\,=\,1\,\cdots\,r)}}
\end{equation}
sont identiquement satisfaites pour tout $x \in \mathcal{ X}^1$ et
tout $a \in \mathcal{ A}^1$, après remplacement de $x'$ par $f( x; \,
a)$.
}\medskip

\CITATION{
Pour le moment, nous voulons nous retenir de supposer que les
équations $x_i' = f_i ( x_1, \dots, x_n, \, a_1, \dots, a_r)$ doivent
représenter un groupe à $r$ termes. En ce qui concerne les
équations~\thetag{ 1}, nous voulons plutôt supposer:
\emphasis{premièrement}, qu'elles représentent une famille de
$\infty^r$ transformations différentes, donc que les $r$ paramètres
$a_1, \dots, a_r$ sont tous essentiels, et \emphasis{deuxièmement},
qu'elles satisfont des équations différentielles de la forme
spécifique~\thetag{ 2}.
\REFERENCE{\cite{enlie1888},~pp.~67--68.}}

\noindent
L'essentialité des paramètres $(a_1, \dots, a_r)$ assure alors (\cite{
enlie1888}, p.~68) que:

\begin{itemize}

\smallskip\item[$\bullet$]
{\em le déterminant des $\psi_{ kj} (a)$ ne s'annule pas
identiquement}; et que:

\smallskip\item[$\bullet$]
{\em les $r$ transformations infinitésimales:}
\[
X_k'
=
\sum_{i=1}^n\,\xi_{ki}(x')\,
\frac{\partial}{\partial x_i'}
\ \ \ \ \ \ \ \ \ \ \ \ \
{\scriptstyle{(k\,=\,1\,\cdots\,r)}}
\]
{\em sont linéairement indépendantes.}

\end{itemize}\smallskip

\smallskip\noindent
{\small\sf\em Nouvelles hypothèses économiques.}
On supposera premièrement que les équations de
transformations~\thetag{ 1}, définies pour $x \in \mathcal{ X}$ et $a
\in \mathcal{ A}$, constituent une famille de difféomorphismes $x
\mapsto f_a ( x) = x'$ du domaine $\mathcal{ X} \subset
\C^n$ sur son image $f_a ( \mathcal{ X}) \subset \C^n$ 
dont les paramètres $(a_1, \dots, a_r)$ sont essentiels, et
deuxièmement que cette famille satisfait des équations différentielles
du type~\thetag{ 2} ci-dessus, {\em en ajoutant expressément l'hypothèse
que le déterminant des $\psi_{ kj} ( a)$ ne s'annule en aucun
paramètre $a \in \mathcal{ A}^1$}.

\medskip

Trois moments théoriques majeurs entrent alors en scène au Chapitre~9
(Vol.~I) de la {\em Theorie der Transformationsgruppen}.

\begin{itemize}

\smallskip\item[$\bullet$]
{\small\sf Premier moment: Constantes de structure.}

\smallskip\item[$\bullet$]
{\small\sf Deuxième moment: Réciproque intermédiaire.}

\smallskip\item[$\bullet$]
{\small\sf Troisième moment: \'Elimination des transformations
auxilaires}. 

\end{itemize}\smallskip

\smallskip\noindent
Le but principal est d'établir que si $r$ transformations
infinitésimales $X_1, \dots, X_r$ satisfont les relations par paires
$\big[ X_k, \, X_j \big] =
\sum_{ j=1}^r \, c_{ kjs}\, X_s$, où les
$c_{ kjs}$ sont des constantes, alors la totalité des transformations
$x' = \exp \big( \lambda_1 X_1 + \cdots +
\lambda_r X_r \big) ( x)$ constitue un groupe
continu local de transformations à $r$ paramètres essentiels. C'est le
Théorème~24, obtenu à l'issu du troisième moment, qui va
aboutir à cette conclusion.

\smallskip\noindent
{\small\sf\em \'Enoncé technique auxiliaire.} Mais tout d'abord, le
théorème suivant sera utilisé d'une manière essentielle par Lie pour
établir, au cours du second moment, la fermeture par composition d'une
famille d'équations de transformation construites en intégrant un
système d'équations aux dérivées partielles construites dans l'espace
produit des $x$ et des $a$.
\`A noter: on n'utilise ici
que la connaissance des groupes à un paramètres.

\smallskip\noindent{\bf Théorème~9.}
(\cite{enlie1888}, p.~72)
\label{Theoreme-9} 
{\em Si, dans les équations de transformations définies pour $(x, a)
\in \mathcal{ X} \times \mathcal{ A}${\rm :}
\def\theequation{1}\begin{equation}
x_i'
=
f_i(x_1,\dots,x_n;\,a_1,\dots,a_r)
\ \ \ \ \ \ \ \ \ \ \ \ \
{\scriptstyle{(i\,=\,1\,\cdots\,n)}},
\end{equation}
les $r$ paramètres $a_1, \dots, a_r$ sont tous essentiels et si, de
plus, certaines équations différentielles de la forme:
\def\theequation{2}\begin{equation}
\frac{\partial x_i'}{\partial a_k}
=
\sum_{j=1}^r\,\psi_{kj}(a_1,\dots,a_r)\cdot
\xi_{ji}(x_1',\dots,x_n')
\ \ \ \ \ \ \ \ \ \ \ \ \
{\scriptstyle{(i\,=\,1\,\cdots\,n\,;\,\,\,k\,=\,1\,\cdots\,r)}}
\end{equation}
sont identiquement satisfaites par $x_1' = f_1 (x; \, a), \dots, x_n'
= f_n ( x; \, a)$, où la matrice $\psi_{ kj} (a)$ est holomorphe et
inversible dans un sous-domaine non vide $\mathcal{ A}^1 \subset
\mathcal{ A}$, et où les fonctions $\xi_{ ji} ( x')$ 
sont holomorphes dans $\mathcal{ X}$, alors en introduisant les $r$
transformations infinitésimales:
\[
X_k
:=
\sum_{i=1}^n\,\xi_{ki}(x)\,
\frac{\partial}{\partial x_i},
\]
l'assertion suivante est vérifiée: toute transformation $x_i' = f_i (
x; \, a)$ dont les paramètres $a_1, \dots, a_r$ se trouvent dans un
petit voisinage d'un paramètre quelconque fixé $a^0 \in \mathcal{
A}^1$ peut être obtenue en exécutant en premier lieu la
transformation:
\[
\overline{x}_i
=
f_i(x_1,\dots,x_n;\,a_1^0,\dots,a_r^0)
\ \ \ \ \ \ \ \ \ \ \ \ \
{\scriptstyle{(i\,=\,1\,\cdots\,n)}},
\]
et en second lieu, en exécutant une certaine transformation:
\[
x_i'
=
\exp\big(t\lambda_1X_1+\cdots+t\lambda_rX_r\big)(\overline{x}_i)
\ \ \ \ \ \ \ \ \ \ \ \ \
{\scriptstyle{(i\,=\,1\,\cdots\,n)}}
\]
du groupe à un paramètre qui est engendré par une combinaison linéaire
appropriée des $X_k$, où $t$ et les $\lambda_1, \dots, \lambda_r$ sont
des nombres complexes <<\,petits\,>>.
}\medskip

Spécialement, cet énoncé technique sera utilisé par Engel et Lie pour
établir que toutes les fois que $r$ transformations infinitésimales
$X_1, \dots, X_r$ constituent une algèbre de Lie ({\em voir}
ci-dessous), la composition de deux transformations de la forme $x' =
\exp \big( t \lambda_1 X_1 + \cdots +
t \lambda_r X_r \big)$ est {\em à nouveau de la même forme}, de telle
sorte que la totalité de toutes ces transformations constitue un {\em
groupe}.

\smallskip
{\em Démonstration.}
Les arguments s'inspirent de ceux qui ont été développés
p.~\pageref{p-69} dans un contexte local; ici, c'est $a^0 \in
\mathcal{ A}^1$ qui remplace le paramètre identité.

D'un premier côté, fixons donc $a^0 \in \mathcal{ A}^1$ et introduisons 
les solutions $a_k = a_k ( t, \lambda_1, \dots, \lambda_r)$ du 
système suivant d'équations différentielles ordinaires:
\[
\frac{da_k}{dt}
=
\sum_{j=1}^r\,\lambda_j\,\widetilde{\psi}_{jk}(a)
\ \ \ \ \ \ \ \ \ \ \ \ \
{\scriptstyle{(k\,=\,1\,\cdots\,r)}},
\]
avec la condition initiale $a_k ( 0, \lambda_1, \dots, \lambda_r) =
a_k^0$, où $\lambda_1, \dots, \lambda_r$ sont des paramètres complexes
(petits) et où, comme précédemment, l'inverse $\widetilde{ \psi}_{ jk}
( a)$ de la matrice $\psi_{ jk} ( a)$ est holomorphe dans $\mathcal{
A}_1$.

D'un deuxième côté, introduisons le flot local:
\[
\exp\big(t\lambda_1X_1+\cdots+t\lambda_rX_r\big)(\overline{x})
=:
h\big(\overline{x};\,t,\lambda\big)
\]
d'une combinaison linéaire générale $\lambda_1 X_1 + \cdots +
\lambda_r X_r$ des $r$ transformations infinitésimales $X_k = \sum_{ i
= 1}^n \, \xi_{ ki} ( x) \, \frac{ \partial}{ \partial x_i}$, où
$\overline{ x}$ est supposé être dans $\mathcal{ A}^1$. Ainsi par sa
définition même, ce flot intègre les équations différentielles
ordinaires:
\[
\frac{dh_i}{dt}
=
\sum_{j=1}^r\,\lambda_j\,\xi_{ji}(h_1,\dots,h_n)
\ \ \ \ \ \ \ \ \ \ \ \ \
{\scriptstyle{(i\,=\,1\,\cdots\,n)}}
\]
avec la condition initiale $h \big( \overline{ x}; \, 0, \lambda \big) 
= \overline{ x}$.

D'un troisième côté, rappelons que l'on peut résoudre les $\xi_{ ji}$
dans les équations différentielles fondamentales~\thetag{ 2} en
inversant la matrice $\widetilde{ \psi}$:
\[
\xi_{ji}(f_1,\dots,f_n)
=
\sum_{k=1}^r\,\widetilde{\psi}_{jk}(a)\,
\frac{\partial f_i}{\partial a_k}
\ \ \ \ \ \ \ \ \ \ \ \ \
{\scriptstyle{(i\,=\,1\,\cdots\,n\,;\,\,\,j\,=\,1\,\cdots\,r)}}.
\]
\`A $i$ et à 
$j$ fixés, multiplions ensuite par $\lambda_j$ les deux 
membres de cette dernière
équation, sommons pour $j$ allant de $1$ jusqu'à $r$ et reconnaissons
$\frac{ da_k}{ dt}$, que nous pouvons donc faire apparaître:
\[
\aligned
\sum_{j=1}^r\,\lambda_j\,\xi_{ji}(f_1,\dots,f_n)
&
=
\sum_{k=1}^r\,
\frac{\partial f_i}{\partial a_k}\,
\sum_{j=1}^r\,\lambda_j\,\widetilde{\psi}_{jk}(a)
\\
&
=
\sum_{k=1}^r\,
\frac{\partial f_i}{\partial a_k}\,
\frac{da_k}{dt}
\\
&
=
\frac{d}{dt}
\big[
f_i\big(x;\,a(t,\lambda)\big)
\big]
\ \ \ \ \ \ \ \ \ \ \ \ \
{\scriptstyle{(i\,=\,1\,\cdots\,n)}}.
\endaligned
\]
Donc les $f_i \big( x; \, a ( t, \lambda)\big)$ satisfont les
\emphasis{mêmes} équations différentielles que les $h_i \big(
\overline{ x}; \, t, \lambda
\big)$, et en outre, si nous assignons à
$\overline{ x}$ la valeur $f ( x; \, a^0)$, les deux collections de
solutions auront de surcroît la \emphasis{même} condition initiale
pour $t = 0$, à savoir: $f(x; \, a^0)$. En conclusion, cette
observation que les $f_i$ et les $h_i$ satisfont les mêmes équations
jointe à la propriété fondamentale d'unicité pour les solutions d'un
système d'équations différentielles ordinaires fournit l'identité de
coïncidence:
\[
\boxed{
f\big(x;\,a(t,\lambda)\big)
\equiv
\exp\big(t_1\lambda_1 X_1
+\cdots+
t\lambda_rX_r\big)\big(f(x;\,a^0)\big)
}
\]
qui exprime que chaque transformation $x' = f( x; \, a)$ pour $a$ dans
un voisinage de $a^0$ s'avère être la composition de la transformation
fixée $\overline{ x} = f ( x; \, a^0)$, suivi de la transformation du
groupe à un paramètre $\exp \big( t \lambda_1 X_1 +
\cdots + t \lambda_r X_r \big) ( \overline{ x})$: c'est ce 
que l'on voulait démontrer. 
\qed\smallskip

\smallskip\noindent
{\small\sf\em Premier moment: Constantes de structure.}
En produisant par différentiation un système d'équations soumis à la
condition de Clebsch-Lie-Frobenius, la démonstration remarquable qui
suit et qui semble ne plus apparaître dans les traités contemporains
va dévoiler l'existence de constantes fondamentales $c_{ kjs}$ en
faisant jouer un rôle quasiment symétrique à l'espace des variables et
à celui des paramètres. Dans la dipolarité du $x$ et du $a$, c'est en
ne se désolidarisant pas de l'action spatiale que la structure
(implicite) de groupe abstrait (dans l'espace des paramètres $a$) va
articuler le premier aspect de l'autonomisation algébrique de sa
genèse ultérieure. L'argument-clé et purement archaïque sera de
fractionner tout crochet entre deux sommes de deux transformations
infinitésimales strictement attachées à l'espace des $x$ et à
l'espaces des $a$ (respectivement):
\[
\big[X+A,\,X'+A'\big]
=
\big[X,\,X']
+
\big[A,\,A'\big]
\]
---\,\,identité garantie par les relations triviales: $0 = \big[ X, \,
A' \big] = \big[ A,\, X'\big]$ qui proviennent de l'hétérogénéité l'un à
l'autre des deux espaces: $0 = \partial_{ x_i} ( a_k) = \partial_{
a_k} ( x_i)$. La constance des constantes $c_{ kjs}$ proviendra alors
d'une réalité de fait évidente et tout aussi archaïque, sachant que la
condition fondamentale de Clebsch:
\[
\small
\aligned
\big[X_k+A_k,\,X_j+A_j\big]
&
=
\text{\small\sf combinaison linéaire}_{a,x}\big(X_s+A_s\big)
\\
&
=
\big[X_k,\,X_j\big]
+
\big[A_k,\,A_j\big]
\\
&
=
\text{\small\sf combinaison linéaire}_x\big(X_t\big)
+
\text{\small\sf combinaison linéaire}_a\big(A_u\big)
\endaligned
\]
force toutes les fonctions coefficients de ces combinaisons
linéaires: 
\[
{\sf coeff}_{a,x}
=
{\sf coeff}_x
=
{\sf coeff}_a
\]
qui dépendent {\em a priori} en toute généralité conjointement de $x$
et de $a$, à ne dépendre en fait que de $x$, et aussi que de $a$, d'où
leur constance {\sl absolue}.

\smallskip\noindent{\bf Théorème~21.}
\label{Theoreme-21}
(\cite{ enlie1888}, pp.~149--150)
{\em 
Si une famille de $\infty^r$ transformations:
\[
x_i'
=
f_i(x_1,\dots,x_n,\,a_1,\dots,a_r)
\ \ \ \ \ \ \ \ \ \ \ \ \
{\scriptstyle{(i\,=\,1\,\cdots\,n)}}
\]
satisfait certaines équations différentielles de la forme spécifique:
\[
\frac{\partial x_i'}{\partial a_k}
=
\sum_{j=1}^r\,\psi_{kj}(a_1,\dots,a_r)\cdot
\xi_{ji}(x_1',\dots,x_n')
\ \ \ \ \ \ \ \ \ \ \ \ \
{\scriptstyle{(i\,=\,1\,\cdots\,n;\,\,\,k\,=\,1\,\cdots\,r)}},
\]
et si on écrit ces équations, ce qui est toujours possible, sous la
forme:
\[
\xi_{ji}(x_1',\dots,x_n')
=
\sum_{k=1}^r\,\widetilde{\psi}_{jk}(a_1,\dots,a_r)\,
\frac{\partial x_i'}{\partial a_k}
\ \ \ \ \ \ \ \ \ \ \ \ \
{\scriptstyle{(i\,=\,1\,\cdots\,n;\,\,\,j\,=\,1\,\cdots\,r)}},
\]
alors il existe entre les $2r$ transformations infinitésimales:
\[
\aligned
X_j'(F)
&
=
\sum_{i=1}^n\,\xi_{ji}(x_1',\dots,x_n')\,
\frac{\partial F}{\partial x_i}
\ \ \ \ \ \ \ \ \ \ \ \ \
{\scriptstyle{(j\,=\,1\,\cdots\,r)}}
\\
A_j(F)
&
=
\sum_{\mu=1}^r\,\widetilde{\psi}_{j\mu}(a_1,\dots,a_r)\,
\frac{\partial F}{\partial a_\mu}
\ \ \ \ \ \ \ \ \ \ \ \ \
{\scriptstyle{(j\,=\,1\,\cdots\,r)}}
\endaligned
\]
des relations de la forme: 
\def\theequation{3}\begin{equation}
\left\{
\aligned
X_k'\big(X_j'(F)\big)
-
X_j'\big(X_k'(F)\big)
&
=
\sum_{s=1}^r\,c_{kjs}\,X_s'(F)
\ \ \ \ \ \ \ \ \ \ \ \ \
{\scriptstyle{(k,\,j\,=\,1\,\cdots\,r)}},
\\
A_k\big(A_j(F)\big)
-
A_j\big(A_k(F)\big)
&
=
\sum_{s=1}^r\,c_{kjs}\,A_s(F)
\ \ \ \ \ \ \ \ \ \ \ \ \
{\scriptstyle{(k,\,j\,=\,1\,\cdots\,r)}},
\endaligned\right.
\end{equation}
où les $c_{ kjs}$ désignent des constantes numériques. En conséquence
de cela, les $r$ équations:
\[
X_j'(F)
+
A_j(F)
=
0
\ \ \ \ \ \ \ \ \ \ \ \ \
{\scriptstyle{(k\,=\,1\,\cdots\,r)}},
\]
qui sont résolubles par rapport à $\frac{ \partial F}{ \partial a_1}$,
\dots, $\frac{ \partial F}{ \partial a_r}$, constituent un système
complet à $r$ termes en les $n+r$ variables $x_1', \dots, x_n'$, $a_1,
\dots, a_r$; si l'on résout les $n$ équations $x_i' = f_i ( x, \, a)$
par rapport à $x_1, \dots, x_n$:
\[
x_i
=
F_i(x_1',\dots,x_n',\,a_1,\dots,a_r)
\ \ \ \ \ \ \ \ \ \ \ \ \
{\scriptstyle{(i\,=\,1\,\cdots\,n)}},
\]
alors $F_1 (x', \, a)$, \dots, $F_n( x', \, a)$ sont des solutions
indépendantes de ce système complet. 
}\medskip

{\em Démonstration.}
Pour commencer, résolvons donc les équations: 
\[
x_i'
=
f_i(x_1,\dots,x_n,\,a_1,\dots,a_r)
\ \ \ \ \ \ \ \ \ \ \ \ \
{\scriptstyle{(i\,=\,1\,\cdots\,n)}}
\]
par rapport à $x_1, \dots, x_n$, ce qui donne:
\[
x_i
=
F_i(x_1',\dots,x_n',\,a_1,\dots,a_r)
\ \ \ \ \ \ \ \ \ \ \ \ \
{\scriptstyle{(i\,=\,1\,\cdots\,n)}}. 
\]
Alors on peut aisément déduire certaines équations différentielles qui
sont satisfaites par $F_1, \dots, F_n$. Si en effet nous différentions
simplement les identités:
\[
F_i\big(f_1(x,\,a),\dots,f_n(x,\,a),\,
a_1,\dots,a_r\big)
\equiv
x_i
\]
par rapport à $a_k$, nous obtenons pour $i$ fixé l'identité: 
\[
\sum_{\nu=1}^n\,
\frac{\partial F_i(x',a)}{\partial x_\nu'}\,
\frac{\partial f_\nu(x,a)}{\partial a_k}
+
\frac{\partial F_i(x',a)}{\partial a_k}
\equiv 
0,
\]
pourvu que l'on pose partout $x_\nu ' = f_\nu (x, \, a)$. Multiplions
maintenant cette identité par $\widetilde{ \psi}_{ jk} (a)$ et sommons
le résultat obtenu pour $k$ allant de $1$ jusqu'à $r$; alors en tenant
compte de l'hypothèse:
\[
\sum_{k=1}^r\,\widetilde{\psi}_{jk}(a)\,
\frac{\partial f_\nu(x,a)}{\partial a_k}
\equiv
\xi_{j\nu}(f_1,\dots,f_n),
\]
nous obtenons les équations suivantes:
\[
\aligned
\sum_{\nu=1}^n\,\xi_{j\nu}(x_1',\dots,x_n')\,
&
\frac{\partial F_i}{\partial x_\nu'}
+
\sum_{k=1}^r\,\widetilde{\psi}_{jk}(a_1,\dots,a_r)\,
\frac{\partial F_i}{\partial a_k}
=
0
\\
&
{\scriptstyle{(i\,=\,1\,\cdots\,n;\,\,\,j\,=\,1\,\cdots\,r)}}.
\endaligned
\]

D'après la manière dont ces équations ont été dérivées, elles sont
satisfaites identiquement lorsqu'on y fait la substitution $x_\nu ' =
f_\nu (x, \, a)$. Mais puisqu'elles ne contiennent pas du tout $x_1,
\dots, x_n$, elle doivent en fait être satisfaites identiquement,
c'est-à-dire: les fonctions $F_1, \dots, F_n$ sont toutes solutions
des équations linéaires aux dérivées partielles suivantes:
\def\theequation{4}\begin{equation}
\aligned
\Omega_j(F)
=
\sum_{\nu=1}^n\,\xi_{j\nu}(x')\,
&
\frac{\partial F}{\partial x_\nu'}
+
\sum_{\mu=1}^r\,\widetilde{\psi}_{j\mu}(a)\,
\frac{\partial F}{\partial a_\mu}
=
0
\\
&
{\scriptstyle{(j\,=\,1\,\cdots\,r)}}.
\endaligned
\end{equation}
Ces $r$ équations contiennent $n + r$ variables, à savoir $x_1',
\dots, x_n'$ et $a_1, \dots, a_r$; de plus, elles sont indépendantes
les unes des autres, puisque le déterminant des $\widetilde{ \psi}_{
j\mu} (a)$ ne s'annule pas, et par conséquent, une résolution par
rapport aux $r$ quotients différentiels $\frac{ \partial F}{ \partial
a_1}, \dots, \frac{ \partial F}{ \partial a_r}$ est possible. Mais
par ailleurs, les équations~\thetag{ 4} possèdent $n$ solutions
indépendantes en commun: les fonctions $F_1 (x', \, a), \dots, F_n(
x', \, a)$ dont le déterminant fonctionnel par rapport aux $x'$:
\[
\sum\,\pm\,
\frac{\partial F_1}{\partial x_1'}\,\dots\,
\frac{\partial F_n}{\partial x_n'}
=
\frac{1}{
\sum\,\pm\,
\frac{\partial f_1}{\partial x_1}\,\dots\,
\frac{\partial f_n}{\partial x_n}
}
\]
ne s'annule pas identiquement, parce que par hypothèse, les équations
$x_i' = f_i( x, \, a)$ représentent des transformations difféomorphes.
Ainsi en définitive, les hypothèses du théorème de
Clebsch-Lie-Frobenius sont satisfaites par les équations~\thetag{ 4},
c'est-à-dire: ces équations constituent un système complet à $r$
termes.

Si nous posons maintenant: 
\[
\sum_{k=1}^r\,\widetilde{\psi}_{jk}(a)\,
\frac{\partial F}{\partial a_k}
=
A_j(F)
\]
et si nous posons aussi: 
\[
\sum_{\nu=1}^n\,\xi_{j\nu}(x')\,
\frac{\partial F}{\partial x_\nu'}
=
X_j'(F),
\]
en accord avec les notations employées dans le \S~3.6, alors les
équations~\thetag{ 4} reçoivent la forme brève:
\[
\Omega_j(F)
=
X_j'(F)
+
A_j(F)
=
0
\ \ \ \ \ \ \ \ \ \ \ \ \
{\scriptstyle{(j\,=\,1\,\cdots\,r)}}.
\]
Comme nous le savons, le fait que ces équations constituent un système
complet revient à ce que des équations de dépendance de la
forme:
\[
\aligned
\Omega_k\big(\Omega_j(F)\big)
-
\Omega_j\big(\Omega_k(F)\big)
&
=
\sum_{s=1}^r\,
\vartheta_{kjs}(x_1',\dots,x_n',a_1,\dots,a_r)\cdot
\Omega_s(F)
\\
&
{\scriptstyle{(k,\,j\,=\,1\,\cdots\,r)}}
\endaligned
\] 
doivent être satisfaites identiquement, quelles que peuvent être les
$F$ comme fonctions de $x_1', \dots, x_n', \, a_1, \dots, a_r$: ce
sont en effet des identités entre champs de vecteurs. Mais puisque
ces identités peuvent aussi être écrites de manière développée comme:
\[
\aligned
X_k'\big(X_j'(F)\big)
&
-
X_j'\big(X_k'(F)\big)
+
A_k\big(A_j(F)\big)
-
A_j\big(A_k(F)\big)
=
\\
&
=
\sum_{s=1}^r\,\vartheta_{kjs}\,X_s'(F)
+
\sum_{s=1}^r\,\vartheta_{kjs}\,A_s(F),
\endaligned
\]
nous pouvons immédiatement les diviser en deux collections 
d'identités:
\def\theequation{5}\begin{equation}
\left\{
\aligned
X_k'\big(X_j'(F)\big)
-
X_j'\big(X_k'(F)\big)
&
=
\sum_{s=1}^r\,\vartheta_{kjs}\,X_s'(F)
\\
A_k\big(A_j(F)\big)
-
A_j\big(A_k(F)\big)
&
=
\sum_{s=1}^r\,\vartheta_{kjs}\,A_s(F),
\endaligned\right.
\end{equation}
et ici, la seconde série d'équations peut encorer être à nouveau
décomposée en:
\[
A_k\big(\widetilde{\psi}_{j\mu}\big)
-
A_j\big(\widetilde{\psi}_{k\mu}\big)
=
\sum_{s=1}^r\,\vartheta_{kjs}\,\widetilde{\psi}_{s\mu}
\ \ \ \ \ \ \ \ \ \ \ \ \
{\scriptstyle{(k,\,j,\,\mu\,=\,1\,\cdots\,r)}}.
\]
Maintenant, comme le déterminant des $\widetilde{ \psi}_{ s\mu}$ ne
s'annule pas identiquement, les fonctions 
$\vartheta_{ kjs}$ sont complètement
déterminées par ces dernière conditions, et il en découle que les
$\vartheta_{ kjs}$ peuvent seulement dépendre de $a_1, \dots, a_r$,
c'est-à-dire qu'elles sont en tout cas libres de $x_1', \dots, x_n'$.
Mais on se convainc aussi aisément que les $\vartheta_{ kjs}$ sont
aussi libres de $a_1, \dots, a_r$: en effet, si dans la première série
d'identités~\thetag{ 5}, on considère $F$ comme une fonction
arbitraire des seules variables $x_1', \dots, x_n'$, alors on obtient
par différentiation par rapport à $a_\mu$ les identités suivantes:
\[
0
\equiv
\sum_{s=1}^r\,\frac{\partial\vartheta_{kjs}}{\partial a_\mu}\,
X_s'(F)
\ \ \ \ \ \ \ \ \ \ \ \ \
{\scriptstyle{(k,\,j,\,\mu\,=\,1\,\cdots\,r)}}.
\]
Mais puisque $X_1' (F), \dots, X_r' ( F)$ sont des transformations
infinitésimales indépendantes, et puisque de plus, les $\frac{
\partial \vartheta_{ kjs}}{\partial a_\mu}$ 
ne dépendent pas de $x_1', \dots, x_n'$, toutes les dérivées
partielles $\frac{ \partial \vartheta_{ kjs }}{ \partial a_\mu}$
s'annulent identiquement; autrement dit, les $\vartheta_{ kjs}$ sont
aussi libres de $a_1, \dots, a_r$: ce sont des constantes numériques, 
comme annoncé.
\qed\smallskip

En particulier, 
ce théorème peut maintenant être immédiatement appliqué à tous les
groupes à $r$ paramètres qui contiennent la transformation identité.

\smallskip\noindent{\bf Théorème~22.}
\label{Theoreme-22}
(\cite{ enlie1888}, p.~150)
{\em Entre les $r$ transformations infinitésimales:
\[
X_k
:=
\frac{\partial f_1}{\partial a_k}(x;\,e)\,
\frac{\partial}{\partial x_1}
+\cdots+
\frac{\partial f_n}{\partial a_k}(x;\,e)\,
\frac{\partial}{\partial x_n}
\ \ \ \ \ \ \ \ \ \ \ \ \
{\scriptstyle{(k\,=\,1\,\cdots\,r)}}
\]
d'un groupe de transformations ponctuelles locales $x' = f ( x; \ a)$
qui contient l'élément identité $e$, il existe des relations par paires de
la forme:
\[
X_k\big(X_j(f)\big)
-
X_j\big(X_k(f)\big)
=
\sum_{s=1}^r\,c_{kjs}\,X_s(f),
\]
où les $c_{ kjs} \in \C$ sont des constantes numériques.}\medskip

En particulier, si un groupe
continu de transformations contient
deux transformations infinitésimales: 
\[
X(f)
=
\sum_{i=1}^n\,\xi_i(x_1,\dots,x_n)\,
\frac{\partial f}{\partial x_i},
\ \ \ \ \ \ 
Y(f)
=
\sum_{i=1}^n\,\eta_i(x_1,\dots,x_n)\,
\frac{\partial f}{\partial x_i},
\]
alors il contient aussi la transformation
infinitésimale:
\[
X\big(Y(f)\big)
-
Y\big(X(f)\big).
\]

\smallskip\noindent
{\small\sf\em Deuxième moment: Réciproque intermédiaire.} L'objectif
est, réciproquement, de démontrer qu'une paire d'algèbres de Lie de
champs de vecteurs $X_1', \dots, X_r'$ et $A_1', \dots, A_r'$ permet
de reconstituer les équations de transformation d'un certain groupe
fini. Spécialement, il va être établi au cours de la démonstration
({\em voir} {\em infra}) que la famille exponentielle $x' = \exp \big(
\lambda_1 X_1 + \cdots +
\lambda_r X_r \big) ( x)$ est stable par
composition, ce que l'on pourrait écrire informellement comme:
\[
\exp\circ\exp 
\equiv 
\exp. 
\]
C'est la première
fois qu'apparaît cette propriété de stabilité lorsque le nombre $r$ de
paramètres est $\geqslant 2$, et l'on peut se convaincre que seule la
condition de fermeture par crochet est à même de garantir cette
propriété.

\smallskip\noindent{\bf Théorème~23.}
\label{Theoreme-23}
(\cite{ enlie1888}, pp.~154--155)
{\em Si $r$ transformations infinitésimales indépendantes:
\[
X_k'(f)
=
\sum_{i=1}^n\,\xi_{ki}(x_1',\dots,x_n')\,
\frac{\partial f}{\partial x_i'}
\ \ \ \ \ \ \ \ \ \ \ \ \
{\scriptstyle{(k\,=\,1\,\cdots\,r)}}
\]
dans les variables
$x_1', \dots, x_n'$ satisfont des conditions par paires de 
la forme: 
\[
X_k'\big(X_j'(f)\big)
-
X_j'\big(X_k'(f)\big)
=
\sum_{s=1}^r\,c_{kjs}\,X_s'(f),
\]
si de plus $r$ transformations infinitésimales: 
\[
A_k(f)
=
\sum_{\mu=1}^r\,
\widetilde{\psi}_{k\mu}(a_1,\dots,a_r)\,
\frac{\partial f}{\partial a_\mu}
\ \ \ \ \ \ \ \ \ \ \ \ \
{\scriptstyle{(k\,=\,1\,\cdots\,r)}}
\]
dans un espace auxiliaire 
$a_1, \dots, a_r$ satisfont les conditions analogues:
\[
A_k\big(A_j(f)\big)
-
A_j\big(A_k(f)\big)
=
\sum_{s=1}^r\,c_{kjs}\,A_s(f)
\]
avec les mêmes constantes $c_{ kjs}$, et si enfin, le déterminant
$\sum \, \pm \, \widetilde{ \psi}_{ 11} (a) 
\cdots \widetilde{ \psi}_{ rr} (
a)$ ne s'annule pas, alors on obtient comme suit les équations d'un
groupe à $r$ paramètres essentiels: on forme le sytème complet à $r$
termes:
\[
X_k'(f)+A_k(f)
=
0
\ \ \ \ \ \ \ \ \ \ \ \ \
{\scriptstyle{(k\,=\,1\,\cdots\,r)}}
\]
et l'on détermine ses solutions générales relativement à un système
approprié de valeurs $a_k = a_k^0$. Si $x_i = F_i ( x_1', \dots, x_n',
\, a_1, \dots, a_r)$ sont ces solutions
générales, alors les équations résolues $x_i' = f_i ( x_1, \dots, x_n,
\, a_1, \dots, a_r)$ représentent un groupe continu de transformations
à $r$ paramètres essentiels. Ce groupe contient la transformation
identique et pour chacune de ses transformations, il contient aussi la
transformation inverse; il est engendré par les $\infty^{ r - 1}$
tranformations infinitésimales:
\[
\lambda_1\,X_1'(f)
+\cdots+
\lambda_r\,X_r'(f),
\]
où $\lambda_1, \dots, \lambda_r$ désignent des constantes
arbitraires. En introduisant des nouveaux paramètres à la place des
$a_k$, les équations du groupe peuvent donc être rapportées à la
forme:}
\[
x_i'
=
x_i
+
\sum_{k=1}^r\,\lambda_k\,\xi_{ki}(x)
+
\sum_{k,\,j}^{1\dots r}\,
\frac{\lambda_k\,\lambda_j}{1\cdot 2}\,
X_j(\xi_{ki})
+\cdots
\ \ \ \ \ \ \ \ \ \ \ \ \
{\scriptstyle{(i\,=\,1\,\cdots\,n)}}.
\]

{\em Démonstration.}
Il est clair que le système d'équations aux dérivées
partielles d'ordre un: 
\[
\Omega_j(F)
=
X_j'(F)
+
A_j(F)
=
0
\ \ \ \ \ \ \ \ \ \ \ \ \
{\scriptstyle{(j\,=\,1\,\cdots\,r)}},
\]
est complet, puisque les hypothèses garantissent que l'on a des
relations par paires de la forme:
\[
\Omega_k\big(\Omega_j(F)\big)
-
\Omega_j\big(\Omega_k(F)\big)
=
\sum_{s=1}^r\,c_{kjs}\,\Omega_s(F),
\] 
et ce système est indépendant, puisque l'hypothèse d'invertibilité de
la matrice $\widetilde{ \psi}_{ k\mu}$ garantit que les équations
$\Omega_1 ( F) = 0$, \dots, $\Omega_r ( F) = 0$ sont résolubles 
par rapport à $\frac{ \partial F}{\partial a_1}$, \dots, 
$\frac{ \partial F}{ \partial a_r}$.

Maintenant, soit $a_1^0, \dots, a_r^0$ un système de valeurs des $a$
dans un voisinage duquel les fonctions $\widetilde{ \psi}_{ jk} (a)$
se comportent régulièrement. D'après le Théorème
p.~\pageref{Theoreme-12} de Clebsch-Lie-Frobenius, le système complet
$\Omega_j ( F) = 0$ possède $n$ solutions $F_1 ( x',\, a), \dots, F_n
( x', \, a)$ qui se réduisent à $x_1', \dots, x_n'$ (respectivement)
pour $a_k = a_k^0$. Imaginons maintenant
que ces solutions générales sont données, formons les $n$ équations:
\[
x_i
=
F_i(x_1',\dots,x_n',\,a_1,\dots,a_r)
\ \ \ \ \ \ \ \ \ \ \ \ \
{\scriptstyle{(i\,=\,1\,\cdots\,n)}},
\]
et résolvons-les par rapport à $x_1', \dots, x_n'$, ce qui est toujours
possible, puisque $F_1, \dots, F_n$ sont évidemment indépendantes
l'une de l'autre, pour autant que seules les variables $x_1', \dots,
x_n'$ sont concernées. Les équations obtenues de cette manière:
\[
x_i'
=
f_i(x_1,\dots,x_n,\,a_1,\dots,a_r)
\ \ \ \ \ \ \ \ \ \ \ \ \
{\scriptstyle{(i\,=\,1\,\cdots\,n)}}
\]
représentent alors, comme nous allons maintenant le démontrer, un
groupe à $r$ paramètres, et en fait naturellement, un groupe contenant
la transformation identique, car pour $a_k = a_k^0$, on a $x_i' =
x_i$.

Tout d'abord, on a identiquement: 
\def\theequation{6}\begin{equation}
\aligned
&
\sum_{\nu=1}^n\,\xi_{j\nu}(x')\,\frac{\partial F_i}{\partial x_\nu'}
+
\sum_{\mu=1}^r\,\widetilde{\psi}_{j\mu}(a)\,
\frac{\partial F_i}{\partial a_\mu}
=
0
\\
&
\ \ \ \ \ \ \ \ \ \ \ \ \
{\scriptstyle{(i\,=\,1\,\cdots\,n;\,\,\,j\,=\,1\,\cdots\,r)}}.
\endaligned
\end{equation}
D'un autre côté, en différentiant $x_i = F_i (x', \, a)$
par rapport à $a_\mu$, on obtient l'équation:
\[
0
=
\sum_{\nu=1}^n\,\frac{\partial F_i}{\partial x_\nu'}\,
\frac{\partial x_\nu'}{\partial a_\mu}
+
\frac{\partial F_i}{\partial a_\mu},
\]
qui est satisfaite identiquement après la substitution $x_\nu ' =
f_\nu ( x,
\, a)$. Multiplions cette équation par $\widetilde{ \psi}_{ j\mu} ( a)$
et sommons pour $\mu$ allant de $1$ jusqu'à $r$, ce qui nous donne une
équation qui devient, en tenant compte de~\thetag{ 6}:
\[
\aligned
&
\sum_{\nu=1}^n\,
\frac{\partial F_i}{\partial x_\nu'}\,
\Big(
\sum_{\mu=1}^r\,\widetilde{\psi}_{j\mu}(a)\,
\frac{\partial x_\nu'}{\partial a_\mu}
-
\xi_{j\nu}(x')
\Big)
=
0
\\
&
\ \ \ \ \ \ \ \ \ \ \ \ \ \ \ \ \ \ \ \ \ \ \
{\scriptstyle{(i\,=\,1\,\cdots\,n;\,\,\,\mu\,=\,1\,\cdots\,r)}}.
\endaligned
\]
Mais comme le déterminant $\sum \, \pm \frac{ \partial F_1}{\partial
x_1'} \, \dots \, \frac{ \partial F_n}{\partial x_n'}$
ne s'annule pas, on peut donc en déduire les équations:
\[
\sum_{\mu=1}^r\,\widetilde{\psi}_{j\mu}(a)\,
\frac{\partial x_\nu'}{\partial a_\mu}
=
\xi_{j\nu}(x')
\ \ \ \ \ \ \ \ \ \ \ \ \
{\scriptstyle{(j\,=\,1\,\cdots\,r\,;\,\,\nu\,=\,1,\,\dots,\,n)}},
\]
qui forment alors un système que nous pouvons à nouveau résoudre par
rapport aux $\frac{ \partial x_\nu'}{ \partial a_\mu}$, 
puisque par hypothèse, le déterminant de $\widetilde{ \psi}_{ j\mu}
(a)$ ne s'annule pas. Ainsi, nous obtenons finalement des équations
différentielles de la forme\footnote{\,
%%%%%%%%%%%%%%%%%%%%%%%-------DEBUT--------%%%%%%%%%%%%%%%%%%%%%%%%%%%
Moment crucial: des équations différentielles fondamentales du type
déjà rencontré sont à présent reconstituées, et c'est grâce à elles
que la structure de groupe fermé par composition va pouvoir renaître
plus bas. 
}: %%%%%%%%%%%%%%%%%%%%%%%%-----FIN-----%%%%%%%%%%%%%%%%%%%%%%%%%%%%%%%
\def\theequation{7}\begin{equation}
\aligned
&
\frac{\partial x_\nu'}{\partial a_\mu}
=
\sum_{j=1}^r\,\psi_{\mu j}(a_1,\dots,a_r)\cdot
\xi_{j\nu}(x_1',\dots,x_n')
\\
&
\ \ \ \ \ \ \ \ \ \ \ \ \ \ \ \ \ \ \ \ \ \ \
{\scriptstyle{(\nu\,=\,1\,\cdots\,n\,;\,\,\,\mu\,=\,1\,\cdots\,r)}}
\endaligned
\end{equation}
qui se réduisent naturellement à des identités après la substitution
$x_\nu ' = f_\nu ( x, \, a)$.

\`A ce point, la démonstration que les
équations $x_i' = f_i (x, \, a)$ représentent un groupe à $r$
paramètres essentiels ne présente pas de difficulté.

En effet, il est facile tout d'abord de voir que les équations $x_i' =
f_i( x, \, a)$ représenent $\infty^r$ transformations distinctes, donc
que les paramètres sont $a_1, \dots, a_r$ sont tous essentiels. Sinon
(\cf~le Théorème p.~\pageref{theoreme-essentiel}), toues les fonctions
$f_1 ( x, \, a)$, \dots, $f_n( x, \, a)$ devraient satisfaire une
équation linéaire aux dérivées partielles de la forme:
\[
\sum_{k=1}^r\,\chi_k(a_1,\dots,a_r)\,
\frac{\partial f}{\partial a_k}
=
0,
\]
dans laquelle les $\chi_k$ seraient libres de $x_1, \dots, x_n$.
En vertu de~\thetag{ 7}, on obtiendrait alors:
\[
\sum_{k,\,j}^{1\dots r}\,
\chi_k(a)\cdot \psi_{kj}(a)\cdot\xi_{j\nu}(f_1,\dots,f_n)
\equiv 
0
\ \ \ \ \ \ \ \ \ \ \ \ \
{\scriptstyle{(\nu\,=\,1\,\cdots\,n)}},
\]
d'où, puisque $X_1' ( F), \dots, X_r' ( F)$ sont
des transformations infinitésimales indépendantes: 
\[
\sum_{k=1}^r\,\chi_k(a)\cdot\psi_{kj}(a)
=
0
\ \ \ \ \ \ \ \ \ \ \ \ \
{\scriptstyle{(j\,=\,1\,\cdots\,r)}};
\]
mais d'après cela, il vient immédiatement: $\chi_1 ( a) = 0$, \dots,
$\chi_r ( a) = 0$, parce que le déterminant des $\psi_{ kj} ( a)$ ne
s'annule pas.

Ainsi, les équations $x_i' = f_i (x, \, a)$ représentent effectivement
une famille de $\infty^r$ transformations différentes. Mais
maintenant, puisque cette famille satisfait certaines équations
différentielles de la forme spécifique~\thetag{ 7}, nous pouvons
appliquer immédiatement le Théorème~9 p.~\pageref{Theoreme-9}.
D'après ce théorème, si $\overline{ a}_1, \dots,
\overline{ a}_r$ sont des paramètres fixés, toute
transformation $x_i' = f_i( x, \, a)$ dont les paramètres $a_1, \dots,
a_r$ se trouvent dans un certain voisinage de $\overline{ a}_1, \dots,
\overline{ a}_r$ peut être obtenue en exécutant
d'abord la transformation:
\[
\overline{x}_i
=
f_i(x_1,\dots,x_n,\,\overline{a}_1,\dots,\overline{a}_r)
\ \ \ \ \ \ \ \ \ \ \ \ \
{\scriptstyle{(i\,=\,1\,\cdots\,n)}},
\]
et ensuite une transformation: 
\[
x_i'
=
\overline{x}_i
+
\sum_{k=1}^r\,\lambda_k\,\xi_{ki}(\overline{x})
+\cdots
\ \ \ \ \ \ \ \ \ \ \ \ \
{\scriptstyle{(i\,=\,1\,\cdots\,n)}}
\]
d'un groupe à un paramètre $\lambda_1 \, X_1 (f) + \cdots + \lambda_r
\, X_r ( f)$, où il est entendu que $\lambda_1, \dots, \lambda_r$ sont
des constantes appropriées. Si nous posons en particulier $\overline{
a}_k = a_k^0$, nous obtenons $\overline{ x}_i = x_i$, donc nous voyons
que la famille des $\infty^r$ transformations $x_i' = f_i (x, \, a)$
coïncide, dans un certain voisinage de $a_1^0, \dots, a_r^0$, avec la
famille des transformations:
\def\theequation{8}\begin{equation}
\aligned
x_i'
=
x_i
&
+
\sum_{k=1}^r\,\lambda_k\,\xi_{ki}(x)
+\cdots
\\
&
{\scriptstyle{(i\,=\,1\,\cdots\,n)}}.
\endaligned
\end{equation}

Une fois ce point atteint, la trame archaïque de l'argumentation
spéculative se résume à opérer une transsubstantiation, contagieuse et
homogénéisante, des types de transformations. En effet, le Théorème~9
p.~\pageref{Theoreme-9} montrait non pas que $f \circ f = f$ ou
que $\exp \circ \exp = \exp$, mais seulement qu'il y a une stabilité
par composition {\em entre transformations d'un type hétérogène}:
\[
f\circ\exp\equiv f, 
\]
ou encore, avec de plus
amples détails, que l'on a:
\def\theequation{$*$}\begin{equation}
\left(
\aligned
&
\overline{x}=f\big(x;\,\overline{a}\big)
\\
&\ \ 
\overline{a}\,\,\text{\rm near}\,\,a^0
\endaligned
\right)
\circ
\left(
\aligned
&
x'=\exp\big(\lambda\,X\big)(\overline{x})
\\
&\ \ \ \ \ \ \
\lambda\,\,\text{\rm near}\,\,0
\endaligned
\right)
\equiv
\left(
\aligned
&
x'=f(x;\,a)
\\
&\ \ 
a\,\,\text{\rm near}\,\,\overline{a}
\endaligned
\right).
\end{equation}
Mais si on applique maintenant cet énoncé au paramètre $\overline{ a}
:= a^0$ considéré comme initial lors de la résolution du système
complet $\Omega_1 ( F) = \cdots = \Omega_r ( f) = 0$, alors puisque
par construction ce paramètre $a^0$ produit la transformation
identique: $\overline{ x} = f ( x; \, a^0 ) = x$, il en découle que
l'on obtient\,\,---\,\,si l'on pose donc $\overline{ a} = a^0$
dans~\thetag{ $*$}\,\,---\,\,l'identité:
\[
\left(
\aligned
&
x'=\exp\big(\lambda\,X\big)(\overline{x})
\\
&\ \ \ \ \ \ \
\lambda\,\,\text{\rm near}\,\,0
\endaligned
\right)
\equiv
\left(
\aligned
&
x'=f(x;\,a)
\\
&\ \ 
a\,\,\text{\rm near}\,\,a^0
\endaligned
\right), 
\]
une coïncidence d'essence que l'on peut réexprimer
dans le langage archaïque comme: 
\[
\exp
\equiv 
f.
\]
On peut donc alors remplacer, dans l'identité de composition
hétérogène $f \circ \exp \equiv f$, non seulement $\exp$ par $f$ pour
obtenir une identité de composition {\em homogène}: $f \circ f = f$,
mais aussi $f$ par $\exp$ pour obtenir une deuxième identité de
composition équivalente: $\exp \circ \exp$, elle aussi {\em homogène}.
Le ré-engendrement de la stabilité par composition 
repose, à la fin de la démonstration, sur le passage
à une communauté de types. 

\smallskip
Voici maintenant les arguments tels qu'écrits dans la langue de Engel.
Si nous choisissons $\overline{ a}_1, \dots, \overline{
a}_r$ arbitrairement dans un certain voisinage de $a_1^0, \dots,
a_r^0$, alors la transformation $\overline{ x}_i = f_i ( x, \,
\overline{ a})$ appartient 
toujours à la famille~\thetag{ 8}. Mais si nous exécutons tout d'abord
la transformation $\overline{ x}_i = f_i(x, \, \overline{ a})$ et
ensuite une transformation appropriée:
\[
x_i'
=
\overline{x}_i
+
\sum_{k=1}^r\,\lambda_k\,\xi_{ki}(\overline{x})
+\cdots
\]
de la famille~\thetag{ 8}, alors d'après ce qui a été dit plus haut,
nous obtenons une transformation $x_i' = f_i ( x, \, a)$ dans laquelle
$a_1, \dots, a_r$ peut prendre toutes les valeurs dans un certain
voisinage de $\overline{ a}_1, \dots,
\overline{ a}_r$. En particulier, si 
nous choisissons $a_1, \dots, a_r$ dans le voisinage de $a_1^0, \dots,
a_r^0$ mentionné plus haut, ce qui est toujours possible, alors à
nouveau la transformation $x_i' = f_i ( x, \, a)$ appartient aussi à
la famille~\thetag{ 8}. Par conséquent, nous voyons que deux
transformations de la famille~\thetag{ 8}, lorsqu'elles sont exécutées
l'une après l'autre, donnent à nouveau une transformation de cette
famille. En définitive, cette famille\,\,---\,\,et naturellement aussi
la famille $x_i' = f_i (x, \, a)$ qui s'identifie à
elle\,\,---\,\,forme un groupe continu de transformations à $r$
paramètres qui contient la transformation identique et des
transformations inverses l'une de l'autre par paires.
\qed\smallskip

\smallskip\noindent
{\small\sf\em Troisième moment: \'Elimination des transformations
auxiliaires.} Retour en arrière et examen du gain synthétique obtenu:
à présent, il faut rebondir et s'interroger sur la possibilité d'une
réciproque plus forte qui ferait l'économie d'hypothèses 
secondaires\,\,---\,\,au prix d'un plus grand effort de pensée. 

\CITATION{
Les hypothèses de l'important Théorème~23
peuvent être simplifiées d'une manière essentielle. 

Le théorème exprime que les $2r$ transformations infinitésimales $X_k
(f)$ et $A_k (f)$ déterminent un certain groupe à $r$ paramètres dans
l'espace des $x$; mais au même moment, il y a une représentation de ce
groupe qui est absolument indépendante des $A_k ( f)$; en effet,
d'après le théorème cité, le groupe en question s'identifie à la
famille des $\infty^{ r-1}$ groupes à un paramètre $\lambda_1 \, X_1
(f) + \cdots + \lambda_r \, X_r (f)$, et cette famille est déjà
déterminée par les seuls $X_k (f)$. Cette circonstance nous conduit à
conjecturer
\deutsch{Dieser Umstand führt uns auf die Vermuthung}
que la famille des $\infty^{ r-1}$ groupes à un paramètre $\lambda_1
\, X_1 (f) + \cdots + \lambda_r \, X_r( f)$ forme toujours un groupe à
$r$ si et seulement si les transformations infinitésimales
indépendantes satisfont des relations par paires de la forme:
\[
X_k\big(X_j(f)\big)
-
X_j\big(X_k(f)\big)
=
\left[
X_k,\,X_j
\right]
=
\sum_{s=1}^r\,c_{kjs}\,X_s(f).
\]
D'après le Théorème~22 p.~\pageref{Theoreme-22}, cette condition est
nécessaire pour que les $\infty^{ r-1}$ groupes à un paramètres $\sum
\, \lambda_k \, X_k (f)$ forment un groupe à $r$ paramètres. Donc
notre conjecture \deutsch{unsere Vermuthung}
revient à suspecter que cette condition nécessaire soit aussi
suffisante.
\REFERENCE{\cite{enlie1888},~pp.~155--156.}}

Cette présomption serait changée en certitude si l'on pouvait, pour
tout système de telles transformations infinitésimales $X_1, \dots,
X_r$, réussir à produire $r$ autres transformations infinitésimales
auxiliaires
\def\theequation{9}\begin{equation}
A_k(f)
=
\sum_{\mu=1}^r\,\widetilde{\psi}_{k\mu}(a_1,\dots,a_r)\,
\frac{\partial f}{\partial a_\mu}
\ \ \ \ \ \ \ \ \ \ \ \ \
{\scriptstyle{(k\,=\,1\,\cdots\,r)}}
\end{equation}
en des variables auxiliaires $(a_1, \dots, a_r)$ destinées
à jouer le rôle de paramètres, 
et qui satisfassent en outre bien sûr les relations correspondantes: 
\[
A_k\big(A_j(f)\big)
-
A_j\big(A_k(f)\big)
=
\sum_{s=1}^r\,c_{kjs}\,A_s(f),
\]
sans que le déterminant $\sum \, \pm\, \widetilde{\psi }_{ 11}
\cdots \widetilde{ \psi}_{ rr }$ ne s'annule.
Ce troisième moment est lui aussi riche d'une métaphysique génétique
où Lie se révèle surprenant d'inventivité.

Comme il n'est question que de l'espace des variables $x_1, \dots,
x_n$ dans les transformations infinitésimales données:
\[
X_k(f)
=
\sum_{i=1}^n\,\xi_{ki}(x)\,\frac{\partial f}{\partial x_i}
\ \ \ \ \ \ \ \ \ \ \ \ \
{\scriptstyle{(k\,=\,1\,\cdots\,r)}}, 
\]
on ne voit pas bien comment engendrer, presque {\em ex nihilo}, une
essence paramétrique. Choisir tout simplement $A_1 := X_1$, \dots,
$A_r := X_r$ en y remplaçant la variable $x$ par la variable $a$ ne
marche {\em certainement pas}, puisqu'il n'y a aucune raison que le
nombre $r$ de paramètres soit égal à la dimension $n$. \`A vrai dire,
dès que $r > n$, les $r$ transformations infinitésimales sont
nécessairement linéairement {\em dépendantes} en tout point $x^0$
fixé, tandis que les $r$ transformations auxilaires du type~\thetag{
9} recherchées pour pouvoir appliquer le Théorème~23 doivent être
linéairement indépendantes en tout $a^0$ fixé, puisque la matrice des
$\widetilde{ \psi}_{ k\mu} ( a)$ doit être inversible. Cette
dépendance linéaire, inévitable lorsque $r > n$, est d'ailleurs le
défaut le plus gênant de $X_1, \dots, X_r$. 

\smallskip

L'idée (remarquable) de Lie consiste à considérer l'action du groupe,
non pas sur les points $x$ pris un à un dans l'espace initial, mais
sur les collections d'un certain nombre $k$ de points {\em en
simultané}, ou, ce qui revient au même, sur les points pris un à un
dans le produit de $k$ {\em copies} de cet espace. Introduisons donc
à cet effet un nombre, ici exactement égal à $r$ pour les besoins
indiqués, de copies de l'espace initial, ces copies étant chacunes
munies de coordonnées notées:
\[
\big(x_1^{(\mu)},\dots,x_n^{(\mu)}\big)
\ \ \ \ \ \ \ \ \ \ \ \ \
{\scriptstyle{(\mu\,=\,1\,\cdots\,r)}}.
\]
Introduisons aussi les transformations infinitésimales qui sont les
reflets des $X_k$ dans chaque espace:
\[
\aligned
X_k^{(\mu)}(f)
=
\sum_{i=1}^n\,\xi_{ki}\big(x_1^{(\mu)},\dots,x_n^{(\mu)}\big)\,
\frac{\partial f}{\partial x_i^{(\mu)}},
\endaligned
\]
et considérons les $r$ transformations infinitésimales:
\[
W_k(f)
=
\sum_{\mu=1}^r\,X_k^{(\mu)}(f).
\]
D'après la Proposition p.~\pageref{Proposition-p-66} ci-dessus, 
ces transformations infinitésimales sont telles
qu'aucune relation de la forme: 
\[
\sum_{k=1}^r\,\psi_k\big(
x_1',\dots,x_n',x_1'',\dots,x_n'',\cdots\cdots,
x_1^{(r)},\dots,x_n^{(r)}
\big)\cdot W_k(f)
=
0
\]
n'est possible, c'est-à-dire: ces transformations infinitésimales sont
indépendantes. Maintenant, puisqu'on a de plus:
\[
W_k\big(W_j(f)\big)
-
W_j\big(W_k(f)\big)
=
\sum_{s=1}^r\,c_{kjs}\,W_s(f),
\]
les $r$ équations indépedantes l'une de l'autre: 
\[
W_1(f)=0,\cdots\cdots,W_r(f)=0
\]
forment un système complet à $r$ termes en les $rn$ variables $x_1',
\dots, x_n', \cdots, x_1^{(r)}, \dots, x_n^{ (r)}$. Ce système
complet possède $r ( n-1)$ solutions indépendantes, que l'on peut
appeler $u_1, u_2, \dots, u_{ rn-r}$. Donc si l'on sélectionne $r$
fonctions $y_1, \dots, y_r$ des $rn$ quantités $x_i^{ ( \mu)}$ 
qui sont indépendantes l'une de l'autre et indépendantes de $u_1,
\dots, u_{ rn - r}$, on peut alors introduire les $y$ et les $u$ comme
nouvelles variables indépendantes à la place des $x_i^{ ( \mu)}$. En
effectuant cela, on obtient:
\[
W_k(f)
=
\sum_{\pi=1}^r\,
W_k(y_\pi)\,
\frac{\partial f}{\partial y_\pi}
+
\sum_{\tau=1}^{rn-r}\,W_k(u_\tau)\,\frac{\partial f}{\partial
u_\tau},
\]
ou encore, puisque tous les $W_k ( u_\tau)$ s'annulent
identiquement:
\[
W_k(f)
=
\sum_{\pi=1}^r\,\omega_{k\pi}
(y_1,\dots,y_r,u_1,\dots,u_{rn-r})\,
\frac{\partial f}{\partial y_\pi},
\]
et ici, les transformations infinitésimales $W_1 ( f), \dots, W_r (
f)$ ne sont reliés par aucune relations de la forme:
\[
\sum_{k=1}^r\,\varphi_k(y_1,\dots,y_r,\,
u_1,\dots,u_{rn-r})\cdot
W_k(f)
=
0.
\]
Naturellement, cette propriété des $W_k( f)$ reste aussi vraie lorsque
l'on confère à $u_\tau$ des valeurs appropriées $u_\tau^0$. Si on pose
alors $\omega_{ k\pi} ( y, u^0) = \omega_{ k\pi}^0 ( y)$, les $r$
transformations infinitésimales indépendantes en les variables
indépendantes $y_1, \dots, y_r$:
\[
V_k(f)
=
\sum_{\pi=1}^r\,\omega_{k\pi}^0(y_1,\dots,y_r)\,
\frac{\partial f}{\partial y_\pi}
\]
satisfont les relations par paires: 
\[
V_k\big(V_j(f)\big)
-
V_j\big(V_k(f)\big)
=
\sum_{s=1}^r\,c_{kjs}\,V_s(f)
\]
et de plus, elles ne sont liées par aucune relation 
de la forme: 
\[
\sum_{k=1}^r\,\varphi_k(y_1,\dots,y_r)\cdot V_k(f)
=
0.
\]
Par conséquent, les $V_k ( f)$ sont des transformations
infinitésimales ayant la constitution recherchée. Ains, nous pouvons
immédiatement appliquer le Théorème~23 p.~\pageref{Theoreme-23} aux
$2r$ transformations infinitésimales $X_1 ( f), \dots, X_r ( f)$, $V_1
( f), \dots, V_r ( f)$ et nous avons donc démontré que les $\infty^{ r
- 1}$ groupes à un paramètre $\sum \, \lambda_k \, X_k ( f)$
constituent un groupe à $r$ paramètres (essentiels). En définitive, la
réciproque <<\,conjecturée\,>> ci-dessus est vraie. 

\smallskip\noindent
{\bf Théorème~24.}
\label{Theoreme-24}
(\cite{enlie1888}, p.~158)
{\em Si $r$ transformations
infinitésimales indépendantes:
\[
X_k(f)
=
\sum_{i=1}^n\,\xi_{ki}(x_1,\dots,x_n)\,
\frac{\partial f}{\partial x_i}
\ \ \ \ \ \ \ \ \ \ \ \ \
{\scriptstyle{(k\,=\,1\,\cdots\,r)}}
\]
satisfont les relations par paires: 
\[
X_k\big(X_j(f)\big)
-
X_j\big(X_k(f)\big)
=
\left[
X_k,\,X_j
\right]
=
\sum_{s=1}^r\,c_{kjs}\,X_s(f),
\]
où les $c_{ kjs}$ sont des constantes, alors la totalité des $\infty^{
r - 1}$ groupes à un paramètre:
\[
\lambda_1\,X_1(f)
+\cdots+
\lambda_r\,X_r(f)
\]
forme un groupe continu à $r$ paramètres, qui contient la
transformation identité, et dont les transformations peuvent être
ordonnées par paires inverses l'une de l'autre.}\medskip

Si les équations $x_i' = f_i ( x_1, \dots, x_n, \, a_1, \dots, a_r)$
représentent une famille de $\infty^r$ transformations et si en outre,
elles satisfont des équations différentielles de la forme spécifique:
\[
\frac{\partial x_i'}{\partial a_k}
=
\sum_{j=1}^r\,\psi_{kj}(a_1,\dots,a_r)\cdot
\xi_{ji}(x_1',\dots,x_n')
\ \ \ \ \ \ \ \ \ \ \ \ \
{\scriptstyle{(i\,=\,1\,\cdots\,n;\,\,\,k\,=\,1\,\cdots\,r)}},
\]
alors comme nous le savons, les $r$ transformations infinitésimales:
\[
X_k(f)
=
\sum_{i=1}^n\,\xi_{ki}(x_1,\dots,x_n)\,
\frac{\partial f}{\partial x_i}
\ \ \ \ \ \ \ \ \ \ \ \ \
{\scriptstyle{(k\,=\,1\,\cdots\,r)}}
\]
sont linéairement indépendantes et de plus, d'après le Théorème~21
p.~\pageref{Theoreme-21}, elles sont reliées entre elles par des
relations de la forme: 
\[
X_k\big(X_j(f)\big)
-
X_j\big(X_k(f)\big)
=
\left[X_k,\,X_j\right]
=
\sum_{s=1}^r\,c_{kjs}\cdot X_s(f).
\]
Ainsi, la famille des $\infty^{ r-1}$ groupes à un paramètre:
\[
\lambda_1\,X_1(f)
+\cdots+
\lambda_r\,X_r(f)
\]
forme un groupe à $r$ paramètres contenant la transformation identité.
Par conséquent, nous pouvons ré-énoncer comme suit le Théorème~9
p.~\pageref{Theoreme-9}.

\smallskip\noindent{\bf Théorème~25.}
(\cite{ enlie1888}, p.~160)
{\em Si une famille de $\infty^r$ transformations:
\[
x_i'
=
f_i(x_1,\dots,x_n,\,a_1,\dots,a_r)
\ \ \ \ \ \ \ \ \ \ \ \ \
{\scriptstyle{(i\,=\,1\,\cdots\,n)}}
\]
satisfait certaines équations différentielles de la forme:
\[
\frac{\partial x_i'}{\partial a_k}
=
\sum_{j=1}^r\,\psi_{kj}(a_1,\dots,a_r)\cdot
\xi_{ji}(x_1',\dots,x_n')
\ \ \ \ \ \ \ \ \ \ \ \ \
{\scriptstyle{(i\,=\,1\,\cdots\,n\,;\,\,\,k\,=\,1\,\cdots\,r)}},
\]
et si le déterminant:
\[
\sum\,\pm\,\psi_{11}(a)\cdots\psi_{rr}(a)
\]
ne s'annule pas, alors toute transformation $x_i' = f_i ( x, \, a)$
dont les paramètres $a_1, \dots, a_r$ se trouvent dans un petit
voisinage d'un paramètre quelconque fixé $a_1^0, \dots, a_r^0$ peut
être obtenue en exécutant en premier lieu la transformation
$\overline{ x}_i = f_i ( x, \, a^0)$ et en second lieu, une
transformation complètement déterminée:
\[
x_i'
=
\overline{x}_i
+
\sum_{k=1}^r\,\lambda_k\,\xi_{ki}(\overline{x})
+\cdots
\ \ \ \ \ \ \ \ \ \ \ \ \
{\scriptstyle{(i\,=\,1\,\cdots\,n)}}
\]
du {\small\sf groupe} à $r$ paramètres qui, sous les hypothèses
supposées, est engendré par les $r$ transformations infinitésimales
indépendantes:}
\[
X_k(f)
=
\sum_{i=1}^n\,\xi_{ki}(x_1,\dots,x_n)\,
\frac{\partial f}{\partial x_i}
\ \ \ \ \ \ \ \ \ \ \ \ \
{\scriptstyle{(k\,=\,1\,\cdots\,r)}}.
\]

{\small 

Pour terminer ce paragraphe, considérons le cas où les équations $x_i'
= f_i ( x, \, a)$, $i=1, \dots, n$ ne contiennent pas la transformation
identité dans le domaine $\mathcal{ A}^1 \subset
\mathcal{ A}$, en supposant bien entendu
comme à la p.~\pageref{debut-38} que deux transformations:
\[
\aligned
x_i'
&
=
f_i(x_1,\dots,x_n,\,a_1,\dots,a_r)
\\
x_i''
&
=
f_i(x_1',\dots,x_n',\,b_1,\dots,b_r)
\endaligned
\] 
exécutées l'une après l'autre, avec $a \in \mathcal{ A}^1$ et $b \in
\mathcal{ A}^1$, produisent la transformation:
\[
x_i''
=
f_i(x_1,\dots,x_n,\,c_1,\dots,c_r)
=
f_i\big(x_1,\dots,x_n,\,\varphi_1(a,b),\dots,\varphi_r(a,b)\big).
\]
Il en découle que ces transformations satisfont des équations
différentielles fondamentales de la forme: 
\[
\aligned
\frac{\partial x_i'}{\partial a_k}
&
=
\sum_{j=1}^r\,\psi_{kj}(a_1,\dots,a_r)\cdot
\xi_{ji}(x_1',\dots,x_n')
\\
&
\ \ \ \ \ \ \ \ \ \ \ \ \
{\scriptstyle{(i\,=\,1\,\cdots\,n\,;\,k\,=\,1\,\cdots\,r)}},
\endaligned
\]
et l'on supposera comme auparavant que le déterminant des $\psi_{ kj}
( a)$ ne s'annule en aucun $a \in \mathcal{ A}^1$.

Dans ce qui va suivre, $a_1^0, \dots, a_r^0$ et pareillement $b_1^0,
\dots, b_r^0$ 
désigneront des points déterminés (fixes) du domaine $\mathcal{ A}^1$,
et $\varphi_k ( a^0, b^0)$ sera noté $c_k^0$. En revanche,
$\overline{ a}_1, \dots, \overline{ a}_r$ 
désignera un point arbitraire du domain $\mathcal{ A}$, de telle sorte
que les équations:
\[
\overline{x}_i
=
f_i(x_1,\dots,x_n,\,\overline{a}_1,\dots,\overline{a}_r)
\]
représenteront une transformation quelconque du groupe.

Toute transformation de la forme $\overline{ x}_i = f_i ( x,
\overline{ a})$ peut être obtenue en exécutant 
d'abord la transformation $x_i' = f_i ( x, a^0)$ et ensuite une
certaine seconde transformation définie comme suit. Résolvons les
équations $x_i' = f_i (x, a^0)$ par rapport à $x_1, \dots, x_n$, 
ce qui donne:
\[
x_i
=
F_i(x_1',\dots,x_n',\,a_1^0,\dots,a_r^0),
\]
et introduisons ces valeurs des $x_i$ dans $\overline{ x}_i = f_i (
x, \overline{ a})$. De cette manière, nous
obtenons, pour les équations cherchées, une
expression de la forme: 
\def\theequation{10}\begin{equation}
\overline{x}_i
=
\Phi_i
\big(x_1',\dots,x_n',\,\overline{a}_1,\dots,\overline{a}_r)
\ \ \ \ \ \ \ \ \ \ \ \ \
{\scriptstyle{(i\,=\,1\,\cdots\,n)}},
\end{equation}
dans laquelle nous n'écrivons pas les
$a_k^0$, puisque nous voulons les considérer 
comme des constantes numériques. 

La transformation~\thetag{ 10}
et bien définie pour tous les systèmes de valeurs 
$\overline{ a}_k$ dans le domaine $\mathcal{ A}$
et son expression peut être prolongée analytiquement
(au sens de Weierstra\ss) à ce domaine 
entier $\mathcal{ A}$; 
cela découle en effet des hypothèses que nous avons
effectuées au sujet de la nature des
fonctions $f_i$ et $F_i$. 

Nous affirmons maintenant que pour certaines
valeurs des paramètres $\overline{
a}_k$, les transformations de la famille $\overline{ x}_i = \Phi_i (
x', \, \overline{ a})$ 
appartiennent au groupe initialement donné $x_i' = f_i
( x, a)$, alors
que par contraste, pour certaines autres valeurs
des 
$\overline{ a}_k$, elles appartiennent au groupe 
$X_1 f, \dots, X_r f$ qui contient la transformation
identité.

\'Etablissons la première partie de cette
assertion. Nous savons que les deux transformations:
\[
x_i'
=
f_i(x_1,\dots,x_n,\,a_1^0,\dots,a_r^0),
\ \ \ \ \ \ \ \ \ \ \
\overline{x}_i
=
f_i(x_1',\dots,x_n',\,b_1,\dots,b_r)
\]
exécutées l'une après l'autre produisent la transformation $\overline{
x}_i = f_i ( x, c)$, où $c_k = \varphi_k ( a^0, b)$; ici, nous pouvons
donner à $b_1, \dots, b_r$ tout système de valeurs dans le domaine
$\mathcal{ A}^1$, tandis que le système de valeurs $c_1, \dots, c_r$
appartient au domaine $\mathcal{ A}$, dans un certain voisinage de
$c_1^0, \dots, c_r^0$. Mais d'après ce qui a été dit précédemment, la
transformation $\overline{ x}_i = f_i ( x, c)$ est aussi obtenue en
exécutant les deux transformations:
\[
x_i'
=
f_i(x_1,\dots,x_n,\,a_1^0,\dots,a_r^0),
\ \ \ \ \ \ \ \ \ \ 
\overline{x}_i
=
\Phi_i(x_1',\dots,x_n',\,\overline{a}_1,\dots,\overline{a}_r)
\]
l'une après l'autre, et en choisissant $\overline{ a}_k = c_k$. Par
conséquent, après la substitution $\overline{ a}_k =
\varphi_k ( a^0, b)$, la transformation $\overline{ x}_i =
\Phi_i ( x', \overline{ a})$ 
est identique à la transformation $\overline{ x}_i = f_i ( x', b)$,
c'est-à-dire: \emphasis{toutes les transformations $\overline{ x}_i =
\Phi_i ( x', \overline{ a})$ dont les paramètres $\overline{ a}_k$ se
trouvent dans un certain voisinage de $c_1^0, \dots, c_r^0$ qui est
défini {\em via} l'équation $\overline{ a}_k =
\varphi_k ( a^0, b)$, appartiennent au 
groupe donné $x_i' = f_i( x, a)$}.

Afin d'établir la seconde partie de notre assertion, rappelons le
Théorème~25. Si $\overline{ a}_1, \dots, \overline{ a}_r$ se trouvent
dans un certain voisinage de $a_1^0, \dots, a_r^0$, alors en vertu de
ce théorème, la transformation $\overline{ x}_i = f_i ( x, \overline{
a})$ peut être obtenue en exécutant d'abord la transformation:
\[
x_i'
=
f_i(x_1,\dots,x_n,\,a_1^0,\dots,a_r^0)
\]
et ensuite une transformation complètement déterminée:
\def\theequation{11}\begin{equation}
\overline{x}_i
=
x_i'
+
\sum_{k=1}^r\,\lambda_k\,\xi_{ki}(x')
+\cdots
\end{equation}
du goupe à $r$ paramètres qui est engendré par
les $r$ transformations infinitésimales indépendantes: 
\[
X_k(f)
=
\sum_{i=1}^n\,\xi_{ki}(x_1,\dots,x_n)\,
\frac{\partial f}{\partial x_i}
\ \ \ \ \ \ \ \ \ \ \ \ \
{\scriptstyle{(k\,=\,1\,\cdots\,r)}}.
\]
D'après la preuve du Théorème~9 p.~\pageref{Theoreme-9}, nous savons
de plus que l'on trouve la transformation~\thetag{ 11} en question en
choisissant d'une manière appropriée $\overline{ a}_1, \dots,
\overline{ a}_r$ comme fonctions indépendantes de
$\lambda_1, \dots, \lambda_r$ et en déterminant inversement
$\lambda_1, \dots, \lambda_r$ comme fonction de $\overline{ a}_1,
\dots, \overline{ a}_r$. Mais d'un autre côté, nous obtenons aussi la
transformation $\overline{ x}_i = f_i ( x, \overline{ a} )$ en
exécutant d'abord la transformation $x_i' = f_i ( x, a^0)$, puis la
transformation $\overline{ x}_i = \Phi_i ( x',
\overline{ a})$. Par conséquent, 
la transformation $\overline{ x}_i =
\Phi_i ( x', \overline{ a})$ appartient au groupe engendré par
$X_1f, \dots, X_r f$ dès que le système de valeurs $\overline{ a}_1,
\dots, \overline{ a}_r$ se trouve dans un certain voisinage de $a_1^0,
\dots, a_r^0$. Pour l'exprimer différemment: les équations $\overline{
x}_i = \Phi_i ( x', \overline{ a})$
sont transformées en les équations~\thetag{ 11} lorsque $\overline{
a}_1, \dots, \overline{ a}_r$ est remplacé par les fonctions
mentionnées de $\lambda_1, \dots,
\lambda_r$. Ceci prouve la deuxième
partie de notre assertion.

Ainsi, les équations de transformation $\overline{ x}_i =
\Phi_i ( x', \overline{ a })$ possèdent la propriété
importante suivante: si à la place de $\overline{ a}_k$, on introduit
les nouveaux paramètres $b_1, \dots, b_r$ au moyen des équations
$\overline{ a}_k = \varphi_k ( a^0, b)$, alors pour un certain domaine
des variables, les équations $\overline{ x}_i = \Phi_i ( x',
\overline{ a})$ prennent la forme $\overline{ x}_i = f_i ( x', b)$;
d'un autre côté, si on introduit les nouveaux paramètres $\lambda_1,
\dots, \lambda_r$ à la place de $\overline{ a}_k$, alors pour un
certain domaine, les équations $\overline{ x}_i = \Phi_i ( x',
\overline{ a})$
se convertissent en:
\[
x_i
=
x_i'
+
\sum_{k=1}^r\,\lambda_k\,\xi_{ki}(x')
+\cdots
\ \ \ \ \ \ \ \ \ \ \ \ \
{\scriptstyle{(i\,=\,1\,\cdots\,n)}}
\]

Voilà donc en définitive une caractéristique importante du groupe
initialement donné $x_i' = f_i ( x, a)$: quand on introduit dans les
équations $x_i' = f_i ( x, a)$ les nouveaux paramètres $\overline{
a}_1, \dots, \overline{ a}_r$ à la place des $a_k$ au moyen de
$\overline{ a}_k = \varphi_k ( a^0, a)$, alors on obtient un système
d'équations de transformations $x_i' = \Phi_i ( x,
\overline{ a})$ qui représentent, 
lorsqu'on les prolonge analytiquement, une famille de transformations
à laquelle appartiennent toutes les transformations d'un certain
groupe à $r$ paramètres contenant la transformations identité.

Nous pouvons exprimer cela comme suit. 

\smallskip\noindent{\bf Théorème~26.}
\label{Theoreme-26}
(\cite{ enlie1888}, p.~163)
{\em 
Tout groupe $x_i' = f_i ( x_1, \dots, x_n, \, a_1, \dots, a_r)$ à $r$
paramètres qui n'est pas engendré par $r$ transformations
infinitésimales indépendantes dérive d'un groupe contenant $r$
transformations infinitésimales indépendantes de la manière suivante:
former tout d'abord les équations différentielles: 
\[
\frac{\partial x_i'}{\partial a_k}
=
\sum_{j=1}^r\,\psi_{kj}(a)\cdot\xi_{ji}(x')
\ \ \ \ \ \ \ \ \ \ \ \ \
{\scriptstyle{(i\,=\,1\,\cdots\,n\,;\,k\,=\,1\,\cdots\,r)}},
\]
qui sont satisfaites par les équations $x_i' = 
f_i ( x, a)$, puis poser:
\[
\sum_{i=1}^n\,\xi_{ki}(x)\,\frac{\partial f}{\partial x_i}
=
X_k(f)
\ \ \ \ \ \ \ \ \ \ \ \ \
{\scriptstyle{(k\,=\,1\,\cdots\,r)}}
\]
et former ensuite les équations finies:
\[
x_i'
=
x_i
+
\sum_{k=1}^r\,\lambda_k\,\xi_{ki}(x)
+\cdots
\ \ \ \ \ \ \ \ \ \ \ \ \
{\scriptstyle{(i\,=\,1\,\cdots\,n)}}
\]
du groupe à $r$ paramètres contenant la transformation identité qui
est engendré par les $r$ transformations infinitésimales indépendantes
$X_1 f, \dots, X_r f$. Avec ces données, il est alors possible, dans
ces équations finies, d'introduire des nouveaux paramètres $\overline{
a}_1, \dots, \overline{ a}_r$ à la place de $\lambda_1, \dots,
\lambda_r$ de manière à ce que les équations de transformations qui en
résultent:
\[
x_i'
=
\Phi_i(x_1,\dots,x_n,\,\overline{a}_1,\dots,\overline{a}_r)
\ \ \ \ \ \ \ \ \ \ \ \ \
{\scriptstyle{(i\,=\,1\,\cdots\,n)}}
\]
représentent une famille de $\infty^r$ transformations qui embrasse,
après prolongement analytique, toutes les $\infty^r$ transformations:
\[
x_i'
=
f_i(x_1,\dots,x_n,\,a_1,\dots,a_r)
\ \ \ \ \ \ \ \ \ \ \ \ \
{\scriptstyle{(i\,=\,1\,\cdots\,n)}}
\]
du groupe.}\medskip

}

\`A l'issue de ce premier trajet fondamental que clôt la fin 
du Chapitre~9 de la {\em Theorie der Transformationsgruppen}, les
Théorème~22 et~24 vont permettre dans la suite à Engel et à Lie d'{\em
identifier systématiquement} tout groupe continu de transformations:
\[
x_i'
=
f_i(x_1,\dots,x_n;\,a_1,\dots,a_r)
\ \ \ \ \ \ \ \ \ \ \ \ \
{\scriptstyle{(i\,=\,1\,\cdots\,n)}}
\]
à une collection de transformations infinitésimales linéairement
indépendantes:
\[
X_k(f)
=
\sum_{i=1}^n\,\xi_{ki}(x_1,\dots,x_n)\,
\frac{\partial f}{\partial x_i}
\ \ \ \ \ \ \ \ \ \ \ \ \
{\scriptstyle{(i\,=\,1\,\cdots\,n)}}
\]
dont les coefficients $\xi_{ ki} ( x)$ sont analytiques (réels ou
complexes) et qui est linéairement fermée par crochets:
\[
\big[X_j,\,X_k\big]
=
\sum_{s=1}^r\,c_{j,k}^s\,
X_s
\ \ \ \ \ \ \ \ \ \ \ \ \
{\scriptstyle{(j,\,k\,=\,1\,\cdots\,r)}}, 
\]
où les $c_{ j,k}^s$ sont des constantes. Une telle 
identification\footnote{\,
%%%%%%%%%%%%%%%%%%%%%%%-------DEBUT--------%%%%%%%%%%%%%%%%%%%%%%%%%%%
Tous les énoncés précédents
sont clairement et explicitement locaux, et ce 
serait se méprendre sur la portée
rigoureuse de la théorie de Lie que de
} %%%%%%%%%%%%%%%%%%%%%%%%-----FIN-----%%%%%%%%%%%%%%%%%%%%%%%%%%%%%%%
d'un
groupe à des générateurs infinitésimaux n'a pas seulement un caractère
d'interchangeabilité ontologique, elle {\em métamorphose} aussi l'être
d'un groupe local arbitraire en linéarisant ses caractéristiques
fondamentales. Avec ces théorèmes fondamentaux, non seulement le
différentiel supplante le fini, mais encore: on s'apprête à algébriser
définitivement la genèse.  Un point de seuil, en effet, est atteint,
c'est un point de non retour: la connaissance mathématique
initialement indécise et problématisante peut à présent se décider à
réenvisager l'objet <<\,groupe continu\,>>\,\,---\,\,maintenant moins
opaque\,\,---\,\,sous un angle absolument neuf: celui des
transformations infinitésimales, plus riches de virtualités et de
manipulations possibles.  La genèse synthétique du concept de groupe
continu de transformations a donc ceci d'<<\,irréversible\,>> que ses
caractéristiques initiales sont destinées à s'effacer devant
l'approfondissement incessant de leur compréhension. 

Sur le plan algébrique, les deux seules contraintes qui s'exercent sur
une collection de transformations infinitésimales sont l'antisymétrie
du crochet:
\[
\sum_{s=1}^r\,c_{j,k}^s\,X_s
=
\big[X_j,\,X_k\big]
=
-\big[X_k,\,X_j\big]
=
\sum_{s=1}^r\,-\,c_{k,j}^s\,X_s
\]
qui se lit simplement: 
\[
0
=
c_{j,k}^s+c_{k,j}^s,
\]
et les identités de type Jacobi; ces dernières sont satisfaites
automatiquement entre triplets de champs de vecteurs et elles
s'écrivent ici pour tous $j, k, l = 1,
\dots, r$:
\[
\aligned
0
&
=
\big[X_j,\,[X_k,\,X_l]\big]
+
\big[X_l,\,[X_j,\,X_k]\big]
+
\big[X_k,\,[X_l,\,X_j]\big]
\\
&
=
\big[X_j,\,
{\textstyle{\sum_{s=1}^r}}\,c_{k,l}^s\,X_s\big]
+
\big[X_l,\,
{\textstyle{\sum_{s=1}^r}}\,c_{j,k}^s\,X_s\big]
+
\big[X_k,\,
{\textstyle{\sum_{s=1}^r}}\,c_{l,j}^s\,X_s\big]
\\
&
=
\sum_{t=1}^r\,X_t
\sum_{s=1}^r\,
\big[
c_{k,l}^s\,c_{j,s}^t
+
c_{j,k}^s\,c_{l,s}^t
+
c_{l,j}^s\,c_{k,s}^t
\big], 
\endaligned
\]
ce qui équivaut aux relations quadratiques:
\[
0
=
\sum_{s=1}^r\,
\big[
c_{k,l}^s\,c_{j,s}^t
+
c_{j,k}^s\,c_{l,s}^t
+
c_{l,j}^s\,c_{k,s}^t
\big].
\]

\HEAD{Chapitre~3.\,\,\,\,Théorèmes fondamentaux sur les groupes
de transformations}{
3.9.\,\,\,Le problème de la classification des groupes}

\noindent{\bf 3.9.~Le problème de la classification
des groupes de transformations.}
Pour Lie, la question dominante dans la théorie qu'il a érigée était
de {\em classifier, à équivalence près, tous les groupes de
transformations possibles}, localement, génériquement:
\label{probleme-classification}
\[
\boxed{
\rule[-3pt]{0pt}{22pt}
\aligned
\text{\bf
Classification des groupes continus finis locaux}\
\\
\text{\bf
de transformations analytiques\,}
\ \ \ \ \ \ \ \ \ \ \ \ \,
\endaligned}
\]
et surtout pour l'espace réel à trois dimension, le seul qui possède
un sens <<\,physique\,>>. Du point de vue des équations de
transformations finies, il est bien entendu naturel de déclarer que
deux groupes de transformations: 
\[
x'
=
f(x;\,a)
\ \ \ \ \ \ 
\text{\rm et}
\ \ \ \ \ \
y' 
=
g(y;b) 
\]
qui agissent sur des espaces
respectifs $x_1, \dots, x_n$ et $y_1, \dots, y_n$ de
la même dimension $n \geqslant 1$ avec le même nombre $r \geqslant 1$
de paramètres essentiels $a_1, \dots, a_r$ et $b_1, \dots, b_r$ sont
{\sl équivalents}\footnote{\,
%%%%%%%%%%%%%%%%%%%%%%%-------DEBUT--------%%%%%%%%%%%%%%%%%%%%%%%%%%%
L'adjectif <<\,semblable\,>> appartenant trop au langage
non-conceptuel,
\deutsch{ähnlich} sera traduit par <<\,équivalent\,>>, 
en référence à la {\sl méthode d'équivalence} que \'Elie Cartan a
développée à la suite de Lie pour une pluralité de structures
géométriques qu'il a interprétées en termes de systèmes différentiels
extérieurs (\cite{ ca1937, koba1972, ster1989, gard1989, ol1995}). }
%%%%%%%%%%%%%%%%%%%%%%%%-----FIN-----%%%%%%%%%%%%%%%%%%%%%%%%%%%%%%%
\deutsch{ähnlich} s'il existe à la fois un
changement de paramètres $b = \beta ( a)$ et un changement de
coordonnées $y = \tau ( x)$ dans l'espace-source qui s'effectue
simultanément aussi dans l'espace-image: $y ' = \tau ( x ')$, de
telle sorte que, après substitutions adéquates, on a la relation:
\[
x'
=
\tau^{-1}(y')
=
\tau^{-1}\big(g(y;b)\big)
=
\tau^{-1}
\big(
g(\tau(x);\,\beta(a))
\big)
\equiv
f(x;\,a),
\] 
la dernière égalité étant
identiquement satisfaite pour tout $x$ et tout $a$.

Mais grâce aux théorèmes fondamentaux, ce problème de
classification revient en fait à classifier les algèbres de Lie locales
finies de champs de vecteurs et leurs sous-algèbres, à changement de
coordonnées près:
\[
\boxed{
\rule[-3pt]{0pt}{22pt}
\aligned
\text{\bf
Classification des algèbres de Lie}
\ \ \ \ \ \ \ \ \ \ \ \
\ \ \ \ \ \ \ \ \
\\
\text{\bf
de champs de vecteurs analytiques
locaux\,}
\ \ \ \ \ \ \ \ \ \ \ \ \ \
\\
\text{\bf
\,en dimensions 1, 2 et 3 et au voisinage de points 
génériques\,}
\endaligned}
\]
Au niveau infinitésimal, il est bien entendu naturel de déclarer que
deux algèbres de Lie de champs de vecteurs analytiques locaux $X_1,
\dots, X_r$ et $Y_1, \dots, Y_r$ de la même dimension $r$ sur deux
espaces de coordonnées $x_1, \dots, x_n$ et $y_1, \dots, y_n$ de la
même dimension $n$ sont (localement) {\sl équivalentes} s'il existe un
difféomorphisme $x \mapsto y = y (x)$ qui envoie chaque $X_k$ sur une
combinaison linéaire $\lambda_{ k1} \, Y_1 +
\cdots + \lambda_{ kr}\, Y_r$ des $Y_l$ à coefficients constants
$\lambda_{ kl}$.

\smallskip

Ainsi, le problème de classification revient-il à trouver des formes
normales les plus simples possibles pour les algèbres de Lie de
transformations infinitésimales. Plus précisément, il s'agit 
d'entreprendre l'étude suivante. 

\smallskip$\bullet$
Déterminer
\deutsch{bestimmen} 
toutes les algèbres de Lie de dimension finie $r$ de champs de
vecteurs analytiques complexes locaux:
\[
X_k
=
\sum_{i=1}^n\,\xi_{ki}(x)\,
\frac{\partial}{\partial x_i}
\ \ \ \ \ \ \ \ \ \
{\scriptstyle{(k\,=\,1\,\cdots\,r)}}
\]
définies dans un certain ouvert initial $U \subset \C^n$; 
il est permis de relocaliser les considérations
un nombre fini de fois à un sous-domaine plus petit dès qu'une
opération mathématique nécessite qu'un certain objet soit nondégénéré,
ou qu'une certaine fonction soit non nulle\footnote{\,
%%%%%%%%%%%%%%%%%%%%%%%-------DEBUT--------%%%%%%%%%%%%%%%%%%%%%%%%%%%
En admettant les relocalisations libres, on évite notamment de se
confronter au difficile problème de trouver des formes normales pour
un {\em unique} champ de vecteur analytique $X$ au voisinage d'un
point où tous ses coefficients s'annulent, une question toujours non
résolue et probablement non résoluble en toute généralité, même en
dimension $n = 2$.  Au contraire, d'après le théorème de redressement
local p.~\pageref{redressement-champ} toute transformation
infinitésimale non identiquement nulle relocalisée en un point
générique est localement équivalent à une unique forme normale:
$\frac{ \partial }{ \partial x_1}$.
}. %%%%%%%%%%%%%%%%%%%%%%%%-----FIN-----%%%%%%%%%%%%%%%%%%%%%%%%%%%%%%%

\smallskip$\bullet$
Rapporter chaque tel système $X_1, \dots, X_r$ à une forme normale la
plus simple possible, \eg assurer que la plupart des coefficients
sont nuls, monomiaux, égaux à des fonctions élémentaires, ou qu'ils
dépendent d'un nombre de variables qui est strictement inférieur à
$n$.

\smallskip$\bullet$
Distinguer précisément tous les systèmes possibles de champs de
vecteurs en introduisant des {\em concepts} géométriques ou
algébriques indépendants des coordonnées afin de ranger tous les
groupes dans des catégories et dans des sous-catégories qui soient
précisément et quasi-instantanément {\em discernables par la pensée}.

\medskip

{\footnotesize\sf

Les problèmes de classification: un groupe de taille importante et
d'une complexité invisible agit de manière quasiment incontrôlable sur
une catégorie d'objets. L'objet quelconque flotte alors, transporté
passivement par les ambiguïtés de sa donation initiale. La saisie
vraie en tant que telle ne peut qu'exprimer fondamentalement la
mobilité unique qui est consubstantielle à l'objet qu'elle vise: le
transport possible d'un être mathématique par une transformation $f$
quelconque:
\[
f_*\big(\text{\small\sf être mathématique})
=
\text{\small\sf le même être vu autrement},
\]
transport qui préserve la nature abstraite générale de l'être en
spécifiant seulement les modes généraux par lesquels il se voit
déterminé, re-déterminé, et déterminé à nouveau. Identité, symétrie et
transitivité: par un acte de pensée réduit à sa plus simple
expression, l'homogénéité ontologique du symbole de transformation:

\[
f_*\big(f_*(\cdot)\big)
=
(ff)_*(\cdot)
=
f_*(\cdot),
\] 
garantit la permanence de cet être, de tous ces êtres.

Mais le subjectif pour soi de la pensée perceptive accentue les
exigences de la saisie. Le subjectif biologique du perceptif
biologique structure en effet fortement la nature désirée de toute
saisie par la pensée. Le champ mathématique doit se restructurer en
accord avec une demande impérieuse d'immédiateté. Il s'agit de
capturer, en un même moment d'intuition globale, tous les objets
possibles dans leur individualité propre, clairement et distinctement,
avec la même netteté que l'{\oe}il embrasse dans son champ visuel
toutes les frontières qui coexistent entre centaines de régions de
couleur.  Reproduire le perceptif dans le champ mathématique, c'est
désirer reproduire l'immédiateté immanente à soi du champ biologique
de la conscience. Le perceptif cherche à s'incrire dans la pensée, et
la pensée cherche à s'inscrire dans le perceptif: double circulation
étrangère à soi du soi de l'être qui pense et qui perçoit, au contact
d'une mathématique qui ne pense pas par elle-même et ne perçoit rien
ni en elle-même, ni par elle-même.

}

%%%\HEAD{Chapitre~3.\,\,\,\,Théorèmes fondamentaux sur les groupes
%%%de transformations}{
%%%3.10.\,\,\,Algébrisation de la genèse}

%%%\medskip\noindent{\bf 3.10.~Algébrisation de la genèse.} 
%%%\label{algebrisation-genese}

%%%%%%%%%%%%%%%%%%%%%%%%%%%%%%%%%%%%%%%%%%%%%%%%%%%%%%%%%%%%%%%%%%%%%

\newpage

% 6   :   393--398

\setcounter{footnote}{0}

$\:$
\bigskip\bigskip\medskip

\begin{center}
{\Large\bf
Partie~III:
\\
\medskip
Traduction française commentée et annotée}

\medskip
{\sf Sophus {\sc Lie}, unter Mitwirkung von Friedrich {\sc Engel}}

{\sf\em Theorie der Transformationsgruppen}

{\sf Dritter und letzter Abschnitt, Abtheilung V}

\end{center}
\label{Partie-III}

\medskip

\begin{center}
\begin{minipage}[t]{10.25cm}
\baselineskip =0.35cm
{\scriptsize

\centerline{\footnotesize\bf Chapitres}

\medskip

\smallskip

{\bf Division~5: Recherche sur les fondements de la géométrie
\dotfill \pageref{Division-5}.}

{\bf Chap.~20: D\'etermination des groupes de $R_3$ relativement
auxquels les paires de points poss\`edent un, et un seul invariant,
tandis que $s > 2$ points n'ont pas d'invariant essentiel
\dotfill \pageref{Chapitre-20}.}

{\bf Chap.~21: Critique des recherches helmholtziennes
\dotfill \pageref{Chapitre-21}.}

{\bf Chap.~22: Premi\`ere solution 
du probl\`eme de Riemann-Helmholtz
\dotfill \pageref{Chapitre-22}.}

{\bf Chap.~23: Deuxi\`eme solution du probl\`eme 
de Riemann-Helmholtz 
\dotfill \pageref{Chapitre-23}.}

}\end{minipage}
\end{center}

\bigskip

\centerline{\Large\sf Division~5}
\label{Division-5}
\thispagestyle{empty}

\bigskip

\HEAD{Remarques préliminaires}{
Division\,\,V.}

\centerline{\large\bf Recherches sur les fondements de la G\'eom\'etrie}

\bigskip\medskip
Euclide a d\'evelopp\'e la G\'eom\'etrie <<\,purement
g\'eom\'etrique\,>> en partant d'un certain nombre d'axiomes simples
et de notions fondamentales \'el\'ementaires, sans utiliser aucun
outil analytique. Aussi admirable que soit son syst\`eme d\'eductif,
celui-ci laisse encore cependant quelque peu \`a d\'esirer, lorsque
l'on consid\`ere la fa\c con dont sont employ\'es les \'el\'ements
fondamentaux.

Premi\`erement, on ne saisit pas si le syst\`eme euclidien d'axiomes
et de notions fondamentales est r\'eellement complet. Il est en effet
maintenant g\'en\'eralement admis qu'au cours de ses d\'eveloppements,
Euclide a introduit des hypoth\`eses tacites qu'il aurait d\^u
formuler comme axiomes. Par exemple, l'introduction du concept
d'espace \`a deux dimensions\footnote{\, 
Pr\'ecis\'ement~: <<\,espace de type surface\,>>, 
\deutschplain{Fl\"achenraum}.
} %%%%%%%%%%%%%%%%%%%%%%%%%%%%%%%%%%%%%%%%%%%%%%%%%%%%%%%%%%%%%%%%%
repose chez Euclide sur un
v\'eritable axiome qu'il n'a pas explicit\'e.

Deuxi\`emement, on peut s'imaginer que certains des axiomes euclidiens
sont superflus, c'est-\`a-dire qu'ils pourraient \^etre d\'emontr\'es
\`a partir des d\'efinitions et axiomes pr\'ec\'edents.

Mais au fond, il est beaucoup plus important de s'assurer avant tout
qu'on poss\`ede un syst\`eme complet et suffisant d'axiomes et de
notions fondamentales, plut\^ot que de se demander si certains axiomes
sont \'eventuellement superflus. Tout de m\^eme, la question de savoir
jusqu'\`a quel point le syst\`eme euclidien d'axiomes doit \^etre
enrichi ou compl\'et\'e reste en attente, sans m\^eme parler de
r\'ealisation d\'efinitive. Cependant, on s'est pr\'eoccup\'e avec
d'autant plus de z\`ele de la deuxi\`eme question, que les recherches
r\'ecentes sur les fondements de la G\'eom\'etrie ont seulement
\'et\'e incit\'ees \`a \'etudier la question de la d\'emontrabilit\'e
ou de la non-d\'emontrabilit\'e du onzi\`eme axiome d'Euclide~: le
postulat des parall\`eles.

\bigskip

Apr\`es que plusieurs math\'ematiciens, notamment Legendre, ont
effectu\'e de nombreuses tentatives infructueuses pour d\'emontrer
l'axiome des parall\`eles, Lobatchevski\u{\i}
(1829) tout d'abord et peu
apr\`es B\'olyai (1832) ont r\'eussi \`a exposer indirectement la
non-d\'emontrabilit\'e de l'axiome des parall\`eles, \`a savoir, en
laissant une g\'eom\'etrie se construire, dans laquelle l'axiome des
parall\`eles n'est pas du tout utilis\'e. D'apr\`es des lettres
anciennes de Gauss, qui ne furent \`a vrai dire publi\'ees que tr\`es
tardivement, il ressort que Gauss \'etait d\'ej\`a parvenu depuis
longtemps \`a des r\'esultats similaires.

Aujourd'hui, on doit v\'eritablement s'\'etonner que les
math\'ematiciens aient d\^u s'\'eclairer sur la n\'ecessit\'e de
l'axiome des parall\`eles seulement par un tel d\'etour, alors qu'un
seul coup d'{\oe}il sur la surface d'une sph\`ere aurait pu leur
montrer qu'une g\'eom\'etrie libre de contradiction est aussi possible
sans l'axiome des parall\`eles, une g\'eom\'etrie qui satisfait en
tout cas les axiomes introduits pr\'ec\'edemment par Euclide, \`a
l'int\'erieur d'un domaine choisi convenablement.

Lobatchevski\u{\i} et B\'olyai ont enti\`erement d\'evelopp\'e leur
g\'eom\'etrie \`a la mani\`ere d'Euclide, \ie de fa\c con
purement g\'eom\'etrique. Riemann fut le premier \`a employer des
instruments analytiques, afin d'en tirer des \'eclaircissements sur
les fondements de la G\'eom\'etrie. Malheureusement, nous ne
poss\'edons de sa main aucune pr\'esentation d\'etaill\'ee de ses
recherches, et nous devons nous en remettre aux explications fort
succinctes et souvent difficiles \`a comprendre que contient sa
soutenance d'habilitation de 1854.

Riemann place au tout d\'ebut de ses recherches la proposition
d'apr\`es laquelle l'espace est une vari\'et\'e num\'erique,
partant que les points de l'espace peuvent \^etre rep\'er\'es par des
coordonn\'ees. Ensuite, il demande quelles propri\'et\'es doivent
\^etre attribu\'ees \`a cette vari\'et\'e num\'erique pour qu'elle
repr\'esente la g\'eom\'etrie euclidienne, ou bien une autre
g\'eom\'etrie semblable. La r\'eponse \`a cette question est
manifestement un probl\`eme purement analytique, qui peut
\^etre r\'esolu en tant que tel.

Cependant, la v\'eritable signification de la proposition d'apr\`es
laquelle l'espace est une vari\'et\'e num\'erique ne ressort pas du
travail de Riemann. Riemann cherche
\`a d\'emontrer cette proposition, \label{394}
mais sa d\'emonstration ne peut pas \^etre prise au s\'erieux. Si l'on
veut v\'eritablement d\'emontrer que l'espace est une vari\'et\'e
num\'erique, on devra, \`a n'en pas douter, postuler auparavant un
nombre non n\'egligeable d'axiomes, ce dont il semble que Riemann
n'ait pas \'et\'e conscient. Il faut cependant prendre en
consid\'eration le fait que par l'introduction des vari\'et\'es
num\'eriques, il importait \`a Riemann de donner une version purement
analytique du probl\`eme, et il faut faire observer que les
hypoth\`eses justifiant la possibilit\'e d'introduire au commencement
les vari\'et\'es num\'eriques n'\'etaient pour lui qu'accessoires.

Ajoutons \`a cela que son \'etude constituait une soutenance
orale \deutsch{Probevorlesung}, et
qu'elle n'\'etait pas destin\'ee \`a l'impression; s'il avait voulu
lui-m\^eme la publier, il l'aurait s\^urement \'ecrite d'une tout
autre mani\`ere.

Du reste, l'introduction du concept de vari\'et\'e num\'erique dans
les recherches sur les fondements de la G\'eom\'etrie ne poss\`ede en
rien un caract\`ere arbitraire, mais s'inscrit objectivement dans
l'essence des choses. En effet, en ce qui concerne l'\'edification de
la G\'eom\'etrie, il faut distinguer plusieurs niveaux. L'un d'entre
eux est ind\'ependant aussi bien de l'axiome des parall\`eles, que du
concept intrins\`eque de surface \deutsch{Fl\"acheninhalt} et de celui
des nombres irrationnels. Mais il y a aussi des niveaux plus \'elev\'es,
entre autres, un niveau o\`u l'on ne peut pas \'eluder le concept de
nombre irrationnel et o\`u on doit en tout cas introduire \`a titre
d'axiome le fait que la droite soit une vari\'et\'e
num\'erique. M. G. Cantor est le premier \`a avoir
\'et\'e confront\'e \`a la
n\'ecessit\'e d'introduire un axiome de cette sorte, si l'on souhaite
pousser l'\'edification de la G\'eom\'etrie jusqu'\`a son
ach\`evement.

Riemann fonde la g\'eom\'etrie des vari\'et\'es num\'eriques sur le
concept de longueur d'un \'el\'ement courbe 
\deutsch{Bogenelement}, dont
d\'ecoule par int\'egration le concept de longueur d'une ligne finie.
Il demande que le carr\'e de la longueur d'un \'el\'ement courbe soit
une fonction compl\`etement homog\`ene du second degr\'e par rapport
aux diff\'erentielles des coordonn\'ees; en ajoutant, entre autres,
l'exigence que chaque ligne arbitraire puisse \^etre d\'eplac\'ee sans
modifier sa longueur, Riemann est alors parvenu \`a ce r\'esultat
qu'en dehors de la g\'eom\'etrie euclidienne, seulement deux autres
g\'eom\'etries sont possibles. Parmi ces deux derni\`eres, l'une
s'identifie \`a celle qui a \'et\'e r\'ealis\'ee par 
Lobatchevski\u{\i}, et
l'autre correspond \`a la g\'eom\'etrie inscrite sur la surface d'une
sph\`ere.

Avant toute chose, du point de vue de l'Analyse, les recherches de
Riemann sur la longueur d'un \'el\'ement courbe\,\,---\,\,qu'il n'a
toutefois qu'esquiss\'ees\,\,---\,\,sont du plus haut int\'er\^et;
entre autres, elles ont vraisemblablement donn\'e l'impulsion aux
d\'eveloppements que Messieurs Lipschitz et Christoffel ont entrepris
ult\'erieurement (depuis 1870) sur les expressions diff\'erentielles
du second degr\'e, Monsieur Lipschitz ayant en tout cas poursuivi au
cours du d\'eveloppement de sa th\'eorie l'objectif de trancher en
m\^eme temps quant \`a la justesse des consid\'erations de
Riemann. Mais on ne peut pas nier que les assertions de
Riemann ne fournissent que peu d'\'eclaircissements sur l'objet
v\'eritable de la recherche, \ie sur les fondements de la
G\'eom\'etrie. Les axiomes de Riemann se r\'ef\`erent en effet tous
\`a la longueur d'un \'el\'ement courbe, donc seulement aux
propri\'et\'es de l'espace dans l'infinit\'esimal; si l'on veut en
d\'eduire quelque chose quant \`a la constitution de l'espace \`a
l'int\'erieur d'une r\'egion d'extension finie, on devra d'abord
effectuer une int\'egration. 
\label{critique-integration}
Mais puisque les axiomes qui doivent
servir \`a l'\'edification d'une G\'eom\'etrie purement
g\'eom\'etrique doivent n\'ecessairement \^etre \'el\'ementaires, et
puisque ni le concept de longueur d'un \'el\'ement courbe ni celui
d'int\'egration ne sont \'el\'ementaires, il est clair que les axiomes
de Riemann sont inutilisables pour un tel objectif.

\renewcommand{\thefootnote}{\fnsymbol{footnote}}

Dans son travail c\'el\`ebre\footnote[1]{\, <<\,{\em Sur les faits qui
se trouvent au fondement de la G\'eom\'etrie}\,>>,
G\"ott. Nachr. 1868, pp.~193--221; voir aussi ses {\oe}uvres
scientifiques compl\`etes, Tome~II, pp.~618--639. 
} %%%%%%%%%%%%%%%%%%%%%%%%%%%%%%%%%%%%%%%%%%%%%%%%%%%%%%%%%%%%%%%%% 
de l'ann\'ee 1868,\label{396} M. de Helmholtz s'est soustrait, quoique
de mani\`ere inconsciente, aux insuffisances des axiomes riemanniens
dont nous venons justement de discuter, notamment lorsqu'il postula
certains axiomes se r\'ef\'erant \`a un nombre fini de points
\'eloign\'es les uns des autres, et lorsqu'il tenta d'en d\'eduire
l'axiome de Riemann sur la longueur d'un \'el\'ement courbe.

\renewcommand{\thefootnote}{\arabic{footnote}}

M. de Helmholtz place express\'ement en premi\`ere position l'axiome
d'apr\`es lequel l'espace est une vari\'et\'e num\'erique; ceci
constitue un progr\`es par rapport \`a Riemann, bien que, \`a vrai
dire, M. de Helmholtz s'enqui\`ere lui aussi de la port\'ee de cet
axiome lorsque, \`a la fin de son travail (\`a la page~221), il
affirme que l'axiome pos\'e par lui <<\,demande d'en accepter moins
que ce que l'on pr\'esuppose dans les d\'emonstrations
g\'eom\'etriques que l'on conduit ordinairement\,>>. Mais cette
affirmation est totalement fausse, d'autant plus que les autres
axiomes helmholtziens renferment des pr\'esuppositions superflues.

Un progr\`es suppl\'ementaire en comparaison avec Riemann consiste en
ce que M. de Helmholtz op\`ere directement avec la famille des
mouvements de l'espace, tout en l'interpr\'etant comme famille des
transformations de la vari\'et\'e num\'erique concern\'ee; il
ex\'ecute m\^eme {\em une fois} deux tels mouvements l'un apr\`es
l'autre et il met \`a profit le fait que ces deux mouvements pris
ensemble peuvent \^etre substitu\'es \`a un troisi\`eme mouvement; il
fait donc en quelque sorte usage des caract\'erististiques de type
<<\,groupe\,>> que poss\`edent les mouvements, sans cependant
conna\^{\i}tre le concept g\'en\'eral de groupe.

M\^eme si le travail de Helmholtz, quant \`a ses hypoth\`eses, marque
un certain progr\`es par rapport aux recherches plus anciennes de
Riemann, on ne doit cependant pas ignorer que le travail de Riemann ne
lui c\`ede en rien quant \`a la valeur math\'ematique. En effet,
tandis qu'il s'est av\'er\'e que la m\'ethode esquiss\'ee par Riemann
pouvait r\'eellement se d\'eduire des propositions \'enonc\'ees par
lui, nous verrons ult\'erieurement que les ressources analytiques dont
Monsieur de Helmholtz s'est servi ne sont pas suffisantes pour
atteindre l'objectif, et qu'au cours de ses recherches, Monsieur de
Helmholtz introduit {\em toute une s\'erie d'hypoth\`eses
incorrectes}.

\renewcommand{\thefootnote}{\fnsymbol{footnote}} 

Si nous embrassons les pens\'ees communes qui se trouvent au fondement
des consid\'erations de Riemann et de M. de Helmholtz, nous pouvons
dire que les deux chercheurs ont pos\'e un probl\`eme de type nouveau,
quoiqu'ils ne l'aient fait que de mani\`ere implicite, probl\`eme que
nous souhaitons appeler <<\,{\small\sf Probl\`eme de
Riemann-Helmholtz}\footnote[1]{\, Au fond, le probl\`eme provient
plut\^ot v\'eritablement de Riemann lui-m\^eme. Nous croyons cependant
que notre appellation du probl\`eme est tout \`a fait l\'egitime, car
M. de Helmholtz est le premier \`a avoir r\'ev\'el\'e que le
probl\`eme n'avait pas du tout \'et\'e r\'esolu par la th\'eorie de
Riemann. Une formulation v\'eritable du probl\`eme n'a \'et\'e fournie
par aucun des deux auteurs.
}\,>>, %%%%%%%%%%%%%%%%%%%%%%%%%%%%%%%%%%%%%%%%%%%%%%%%%%%%%%%%%%%
et qui peut s'\'enoncer bri\`evement comme suit~: \label{397} {\em
Trouver des propri\'et\'es qui permettent de distinguer non seulement
la famille des mouvements euclidiens, mais aussi les deux familles de
mouvements non-euclidiens, et gr\^ace auxquelles ces trois familles
appara{\^i}tront alors comme remarquables par rapport \`a toutes les
autres familles de mouvements d'une vari\'et\'e num\'erique.}

\renewcommand{\thefootnote}{\arabic{footnote}} 

%%%\Fill 
%%%[[Traduire ou citer le passage de la fin du m\'emoire de 1880.
%%%$+$ Klein.]]

Lorsqu'en 1869, Lie a communiqu\'e \`a son ami F.~Klein ses
premi\`eres recherches sur les groupes continus, Klein a attir\'e
tr\`es t\^ot l'attention de Lie sur les recherches de Riemann et de
Helmholtz, en soulignant que le concept de groupe continu y jouait un
r\^ole implicite (\cf le m\'emoire de Lie dans le tome 16
des Math. Ann., p.~527). Mais c'est seulement au cours de l'ann\'ee
1884 qu'\`a la suite des invitations renouvel\'ees de Klein, Lie a
entrepris de mettre strictement \`a l'\'epreuve les r\'ealisations
helmholtziennes et d'apporter un traitement approfondi au probl\`eme
de Riemann-Helmholtz dans l'espace ordinaire \`a trois dimensions, au
moyen de la th\'eorie des groupes. Lie n'a rencontr\'e \`a cet \'egard
aucune difficult\'e, sachant qu'il avait d\'ej\`a d\'etermin\'e depuis
longtemps tous les groupes continus finis de cet espace. Lie a tout
d'abord fait conna\^{\i}tre les r\'esultats de cette recherche dans le
{\em Leipziger Berichten} en 1886, sans que leur d\'eduction soit
accompagn\'ee des calculs n\'ecessaires. Ensuite, dans deux gros
m\'emoires de l'ann\'ee 1890, parus \'egalement dans le {\em Leipziger
Berichten}, non seulement il formula une critique d\'etaill\'ee des
d\'eveloppements helmholtziens, mais encore, ce qui constituait le but
principal de tout ce travail, il apporta plusieurs nouvelles solutions
au probl\`eme de Riemann-Helmholtz.

Les chapitres qui suivent constituent un remaniement et en 
partie des compl\'ements aux travaux sus-mentionn\'es, que Lie a
publi\'es \`a ce sujet.

\bigskip

Le Chapitre~20 contient certaines th\'eories qui seront expos\'ees au
mieux pour elles-m\^emes, avant que nous rendions compte des
d\'eveloppements helmholtziens. Une partie des axiomes qui ont \'et\'e
pos\'es par Helmholtz peut en effet s'exprimer, comme nous le verrons,
de la fa\c con suivante~: la totalit\'e des mouvements constitue un
groupe, et relativement \`a ce groupe, deux points poss\`edent un et
un seul invariant, tandis que tous les invariants d'un nombre de
points sup\'erieur \`a deux peuvent s'exprimer comme fonctions des
invariants de paires de points. Afin de mesurer la port\'ee des
axiomes helmholtziens, il faut conna\^{\i}tre \`a l'avance tous les
groupes qui satisfont pr\'ecis\'ement \`a cette exigence, et puisque
cette t\^ache, d\'ej\`a importante en elle m\^eme, qui consiste \`a
d\'eterminer tous ces groupes, ne peut pas \^etre accomplie sans
quelques d\'epenses de calcul, il semble recommand\'e d'entreprendre
sp\'ecialement une telle d\'etermination.

Dans le Chapitre~21, nous formulerons une critique d\'etaill\'ee des
axiomes et des r\'esultats expos\'es par Helmholtz. Ensuite, dans les
Chapitres~22 et~23, nous exposerons diff\'erentes solutions du
probl\`eme de Riemann-Helmholtz, dans l'espace ordinaire et dans
l'espace $n$ fois \'etendu. Enfin, dans le Chapitre~24, nous
discuterons et nous critiquerons quelques recherches r\'ecentes
sur les fondements de la G\'eom\'etrie.

Nous ne pouvons pas conclure ces remarques pr\'eliminaires sans
souligner express\'ement que les recherches qui vont suivre n'ont pas
la pr\'etention de constituer des sp\'eculations philosophiques sur
les fondements de la G\'eom\'etrie; elles ont seulement pour objet
d'apporter un traitement soign\'e de type <<\,th\'eorie des
groupes\,>> \`a ce probl\`eme de la th\'eorie des groupes\footnote{\,
Certainement intentionnelle, la redondance est explicite~: {\em Eine
sorgf\"altige gruppentheoretische Behandlung des gruppentheoretischen
Problems, das wir als das Riemann-Helmholtzsche Problem bezeichet
haben}.
} %%%%%%%%%%%%%%%%%%%%%%%%%%%%%%%%%%%%%%%%%%%%%%%%%%%%%%%%%%%
que nous avons appel\'e {\sl Probl\`eme de
Riemann-Helmholtz}. \`A la fin de la pr\'esente Division~V, nous
\'evoquerons le b\'en\'efice que la r\'esolution de ce probl\`eme peut
apporter \`a l'\'edification d'un syst\`eme de G\'eom\'etrie.

M\^eme si nous n'entreprenons pas ici une telle tentative, nous
souhaitons toutefois exprimer comme \'etant notre conviction profonde,
l'opinion d'apr\`es laquelle il n'est nullement impossible de mettre
sur pied un syst\`eme d'axiomes g\'eom\'etriques qui soit suffisant et
ne renferme pas de condition superflue. Malheureusement, on doit
indubitablement reconna\^{\i}tre qu'il n'existe que tr\`es peu de
recherches qui ont v\'eritablement d\'evelopp\'e les fondements de la
G\'eom\'etrie.

%%%%%%%%%%%%%%%%%%%%%%%%%%%%%%%%%%%%%%%%%%%%%%%%%%%%%%%%%%%%%%%%%%%%%

\newpage

% 39   :   399--437

\setcounter{footnote}{0}

$\:$
\bigskip\bigskip\bigskip

\centerline{\Large Chapitre~20.}
\label{Chapitre-20}
\thispagestyle{empty}

\bigskip

\noindent
\begin{center}
{\large\bf 
D\'etermination des groupes de $R_3$ relativement auxquels les paires
de points poss\`edent un, et un seul invariant, tandis que $s > 2$
points n'ont pas d'invariant essentiel.}
\end{center}

\bigskip\medskip
Les d\'eveloppements du pr\'esent chapitre se rattachent au \S59 du
Tome~I (p.~218 sq.). \`A ce moment-l\`a, nous\footnote{\, Ambiguïté
éventuelle, ici, sur le <<\,nous\,>>, qui pourrait aussi bien renvoyer à
un impersonnel mathématique qu'au duo de rédacteurs formé par Engel et
Lie.
} 
avons introduit le concept d'{\em invariant de plusieurs points} et
maintenant, nous voulons d\'eterminer tous les groupes continus finis
de $R_3$ qui satisfont certaines exigences quant aux invariants d'un
nombre quelconque fini de points. Notamment, deux points doivent avoir
un et un seul invariant relativement \`a chaque tel groupe, tandis que
tous les invariants de $s > 2$ points doivent toujours pouvoir
s'exprimer comme fonctions des invariants des paires de points qui
sont comprises dans ces $s$ points. Si, en concordance avec le Tome~I,
p.~219, nous disons qu'un invariant de $s$ points est {\em essentiel}
s'il ne peut pas \^etre exprim\'e au moyen des invariants d'un
sous-syt\`eme de $s-1$ (voire moins) points, alors nous pouvons
\'enoncer notre probl\`eme (qui se trouve d\'ej\`a exprim\'e dans le
titre du chapitre) de la mani\`ere suivante:
 
\smallskip

{\em D\'eterminer~\label{399} tous les groupes continus finis de $R_3$
relativement auxquels deux points poss\`edent un et un seul invariant,
tandis qu'un nombre de points sup\'erieur \`a deux n'a jamais
d'invariant essentiel}.

\smallskip

Nous allons d'abord r\'esoudre ce probl\`eme sans tenir compte de la
condition de r\'ealit\'e. Par cons\'equent, nous allons chercher en
premier lieu tous les groupes de transformations complexes qui
poss\`edent la propri\'et\'e indiqu\'ee, et plus tard seulement, nous
r\'esoudrons aussi le probl\`eme pour les groupes de
transformations\renewcommand{\thefootnote}{\fnsymbol{footnote}}
r\'eels\footnote[1]{\ En 1886, Lie a d\'ej\`a esquiss\'e dans le {\em
Leipziger Berichten}, p.~337~sq., les r\'esultats qui alimentent le
pr\'esent chapitre; d\'ebut 1890, il en a donn\'e des
d\'emonstrations d\'etaill\'ees (ib., pp.~355--418). Ce qui suit est
un remaniement de cette derni\`ere \'etude.
}. 
\renewcommand{\thefootnote}{\arabic{footnote}}

Nous nous sommes d\'ej\`a expliqu\'es p.~397~sq. sur les raisons qui
nous ont pouss\'es \`a placer ce chapitre en position pr\'eliminaire,
avant de passer aux recherches proprement dites sur les fondements de
la G\'eom\'etrie.

\HEAD{Groupes de $R_3$ pour lesquels deux points ont un seul
invariant.}{Division\,\,V.\,\,\,Chapitre\,\,20.\,\,\,\S\,\,85.}

\sectiondritterV{\S\,\,\,85.
\\ 
Propri\'et\'es caract\'eristiques des groupes recherch\'es.}
\label{S-85}
\setcounter{footnote}{0}

Soit:
\[
X_kf
=
\xi_k(x,y,z)\,p
+
\eta_k(x,y,z)\,q
+
\zeta_k(x,y,z)\,r 
\ \ \ \ \ \ \ 
{\scriptstyle{(k\,=\,1\,\cdots\,m)}}
\]
un groupe\footnote{\, Tout groupe de transformations fini et continu
est systématiquement identifié par Lie à un système de générateurs 
infinitésimaux clos par crochets.  
} %%%%%%%%%%%%%%%%%%%%%%%%%%%%%%%%%%%%%%%%%%%%%%%%%%%%%%%%%%%
\`a $m$ param\`etres\footnote{\, Le texte allemand imprim\'e
comporte ici l'une des tr\`es rares coquilles de tout le trait\'e
(Tomes~I, II et~III), Engel et Lie ayant not\'e par r\'eflexe $r$ au
lieu de $m$ ce nombre de param\`etres (nous rectifions), comme ils ont
l'habitude de le faire lorsqu'ils envisagent un groupe continu
g\'en\'eral, mais dans ce chapitre, la lettre $r$ n'est plus
autoris\'ee, car elle entrerait en confusion avec la troisi\`eme
transformation infinit\'esimale de l'espace, traditionnellement
not\'ee $r = \frac{\partial f}{\partial z}$.
} %%%%%%%%%%%%%%%%%%%%%%%%%%%%%%%%%%%%%%%%%%%%%%%%%%%%%%%%%%% 
poss\'edant la
constitution exig\'ee. Si ensuite:
\def\theequation{1}\begin{equation}
x_1,y_1,z_1;\ \ \ \ \
x_2,y_2,z_2;\ \ \ \ \
\dots\dots;\ \ \ \ \
x_s,y_s,z_s
\end{equation}
sont $s$ points arbitraires\footnote{\, Ici, $R_3$ désigne l'espace à
trois dimensions, en fait complexe: c'est $\C^3$, et non $\R^3$. Engel
et Lie travaillent toujours d'abord sur $\C$ (sans le préciser), puis
sur $\R$ (en le précisant), 
mais ils notent à nouveau $R_3$ (ou $R_n$)
l'espace réel. 
} %%%%%%%%%%%%%%%%%%%%%%%%%%%%%%%%%%%%%%%%%%%%%%%%%%%%%%%%%%% 
de $R_3$ et si nous posons:
\[
\xi_k(x_\nu,y_\nu,z_\nu)\,p_\nu
+
\eta_k(x_\nu,y_\nu,z_\nu)\,q_\nu
+
\zeta_k(x_\nu,y_\nu,z_\nu)\,r_\nu
=
X_k^{(\nu)}f,
\]
alors premi\`erement, les $m$ \'equations lin\'eaires aux d\'eriv\'ees
partielles:
\def\theequation{2}\begin{equation}
\label{eq-2-p-400}
X_k^{(1)}f
+
X_k^{(2)}f
=
0
\ \ \ \ \ \ \
{\scriptstyle{(k\,=\,1\,\cdots\,m)}}
\end{equation}
en les six variables: $x_1, y_1, z_1, x_2, y_2, z_2$ doivent
poss\'eder en commun une et une seule solution\footnote{\, Plus
précisément, la solution générale du système~\thetag{ 2} est 
nécessairement une
fonction arbitraire de l'unique invariant $J$,
de la forme $\Phi \big( J ( x_1, y_1, z_1; x_2, y_2,
z_2)\big)$, où 
$\Phi$ est une fonction arbitraire. 
En permutant les variables $(x_1, y_1, z_1)$ et $(x_2,
y_2, z_2)$, les équations~\thetag{ 2} restent inchangées, donc $J
(x_2, y_2, z_2; \, x_1, y_1, z_1)$ est aussi une solution 
de~\thetag{ 2}, mais {\em non
essentielle}, puisqu'elle est de la forme $\Phi \big( J
( x_1, y_1, z_1; x_2, y_2, z_2)\big)$. 
}: %%%%%%%%%%%%%%%%%%%%%%%%%%%%%%%%%%%%%%%%%%%%%%%%%%%%%%%%%%%
\def\theequation{3}\begin{equation}
J
\big(
x_1,y_1,z_1;\,
x_2,y_2,z_2
\big);
\end{equation}
et deuxi\`emement, lorsque $s > 2$, toutes les solutions communes des
$m$ \'equations:
\def\theequation{4}\begin{equation}
X_k^{(1)}f
+
X_k^{(2)}f
+
\cdots
+
X_k^{(s)}f
=
0
\ \ \ \ \ \ \
{\scriptstyle{(k\,=\,1\,\cdots\,m)}}
\end{equation}
en les $3s$ variables $x_k, y_k, z_k$ doivent se laisser exprimer au
moyen des $\frac{ s ( s-1)}{ 1 \, \cdot \, 2}$ fonctions\footnote{\,
Ici, $1 \leqslant \lambda < \mu
\leqslant s$, car les autres fonctions sont des
invariants non essentiels. 
}: %%%%%%%%%%%%%%%%%%%%%%%%%%%%%%%%%%%%%%%%%%%%%%%%%%%%%%%%%%%
\def\theequation{5}\begin{equation}
J
\big(
x_\lambda,y_\lambda,z_\lambda;\
x_\mu,y_\mu,z_\mu
\big)
\ \ \ \ \ \ \
{\scriptstyle{(\lambda\,=\,1\,\cdots\,s-1;\,\,
\mu\,=\,\lambda\,+\,1\,\cdots\,s)}}.
\end{equation}

Avant toute chose, il faut remarquer que chaque groupe poss\'edant la
constitution exig\'ee\footnote{\, Pour l'instant,  
les groupes sont inconnus. Le Théorèmes~37
(sur $\R$) p.~\pageref{Theorem-37-p-433} 
fournira la liste des onze groupes réels 
répondant au problème posé.  
} %%%%%%%%%%%%%%%%%%%%%%%%%%%%%%%%%%%%%%%%%%%%%%%%%%%%%%%%%%%
doit \^etre transitif. \label{400} En effet, si
le groupe: $X_1 f \dots X_m f$ \'etait intransitif, les
\'equations:
\[
X_kf=0
\ \ \ \ \ \ \
{\scriptstyle{(k\,=\,1\,\cdots\,m)}}
\]
auraient d\'ej\`a en tout cas une solution commune, et les
\'equations~\thetag{ 2} poss\`ederaient
alors au minimum deux solutions
communes ind\'ependantes\footnote{\, Si $\varphi ( x, y, z)$ 
est une solution non constante du système 
$X_1 f = \cdots = X_m f = 0$, alors
$\varphi ( x_1, y_1, z_1)$ et $\varphi ( x_2, y_2, z_2)$ sont deux
solutions indépendantes de~\thetag{ 2}. 
}, %%%%%%%%%%%%%%%%%%%%%%%%%%%%%%%%%%%%%%%%%%%%%%%%%%%%%%%%%%%
en contradiction avec notre exigence.

Afin de trouver encore d'autres propri\'et\'es des groupes
recherch\'es, rappelons que les\footnote{\, Seconde inadvertance du
texte imprim\'e, que nous rectifions \`a nouveau, mais qui ne se
reproduira plus.
} %%%%%%%%%%%%%%%%%%%%%%%%%%%%%%%%%%%%%%%%%%%%%%%%%%%%%%%%%%%
$m$ transformations infinit\'esimales:
\def\theequation{6}\begin{equation}
X_k^{(1)}f
+
X_k^{(2)}f
+\cdots+
X_k^{(s)}f
\ \ \ \ \ \ 
{\scriptstyle{(k\,=\,1\,\cdots\,m)}}
\end{equation}
produisent un groupe \`a $m$ param\`etres en les $3s$
variables~\thetag{ 1}, lequel indique comment le syst\`eme des $s$
points~\thetag{ 1} de $R_3$ est transform\'e par le groupe: $X_1 f
\dots X_m f$. En outre, nous avons l'intention de supposer que les $s$
points~\thetag{ 1} soient non seulement chacun en position
g\'en\'erale vis-\`a-vis du groupe: $X_1 f \dots X_mf$, mais encore
soient {\em mutuellement}\, en position g\'en\'erale\footnote{\,
\^Etre mutuellement en position générale, pour un système de $s$
points $p_1, \dots, p_s$ de $R_3$ signifie simplement que le point
produit $(p_1, \dots, p_s)$ dans $R_{ 3s}$ soit en position
générale. 
}, %%%%%%%%%%%%%%%%%%%%%%%%%%%%%%%%%%%%%%%%%%%%%%%%%%%%%%%%%%%
de telle sorte que le syst\`eme des valeurs~\thetag{ 1}
soit en position g\'en\'erale vis-\`a-vis du groupe~\thetag{ 6}.

Sous ces hypoth\`eses, on peut facilement embrasser d'un coup
d'{\oe}il [{\sc übersehen}] quelles sont toutes les positions
nouvelles que le syst\`eme des $s$ points~\thetag{ 1} peut prendre \`a
travers les transformations du groupe~\thetag{ 6}; en effet,
d'apr\`es le Tome~I, p.~216, la mobilit\'e g\'en\'erale de ces
syst\`emes de points n'est limit\'ee par aucune autre
condition\footnote{\, \`A l'endroit cité, il a été établi qu'à un
groupe quelconque de transformations $X_1f, \dots, X_m f$ à $m$
paramètres agissant sur un espace $(x_1, \dots, x_n)$
de dimension $n$ quelconque est toujours attaché, localement au
voisinage d'un point générique, un nombre déterminé $n - q$
(éventuellement nul) de fonctions indépendantes $\Omega_1 ( x), \dots,
\Omega_{ n-q} ( x)$ découpant l'espace en le feuilletage:
\[
\Omega_1(x)
=
{\rm const}_1,\ \ 
\dots,\ \ 
\Omega_{n-q}(x)
=
{\rm const}_{n-q}
\]
qui reste stable par l'action du groupe. Alors l'action du groupe est
localement transitive en famille sur les feuilles $\Omega_k ( x) =
{\rm const}_k$, $k = 1, \dots, n-q$. De plus, l'exigence
helmholtzienne reformulée abstraitement par Lie 
demande que les positions
$(x_1',
\dots, x_n')$ que peut prendre $(x_1, \dots, x_n)$ par l'action du
groupe soient entièrement déterminées par les fonctions invariantes
$J$, lesquelles
satisfont $J ( x_1', \dots, x_n') = J ( x_1, \dots,
x_n)$ pour toute transformation finie $x_i' = f_i( x_1,
\dots, x_n;  a_1, \dots, a_m)$ du groupe recherché. 
} %%%%%%%%%%%%%%%%%%%%%%%%%%%%%%%%%%%%%%%%%%%%%%%%%%%%%%%%%%%
que la
condition que chaque invariant du groupe doit conserver la m\^eme
valeur num\'erique pour toutes les positions du syst\`eme de
points. Si nous notons alors $x_k', y_k', z_k'$ la position que prend
le point $x_k, y_k, z_k$ apr\`es une transformation quelconque du
groupe: $X_1 f \dots X_m f$ et si nous nous rappelons que tous les
invariants du groupe~\thetag{ 6} doivent pouvoir \^etre exprim\'es au
moyen des invariants~\thetag{ 5}, nous reconnaissons par cons\'equent
que les $s$ points~\thetag{ 1} peuvent se transformer\footnote{\, Les
$s$ points $p_1, \dots, p_s$ mutuellement en position générale de
coordonnées $(x_\lambda, y_\lambda, z_\lambda)$, $\lambda = 1, \dots,
s$, sont transformés par le groupe (inconnu, que l'on recherche) en
des points $(x_\lambda', y_\lambda', z_\lambda')$ de telle sorte
qu'une <<\,pseudodistance\,>> soit conservée: $J ( p_\lambda', p_\mu'
) = J ( p_\lambda, p_\mu)$, et c'est cette seule hypothèse
qui va <<\,faire naître\,>> les 11 groupes possibles 
collectés dans le Théorème~37, p.~\pageref{Theorem-37-p-433}. 
} %%%%%%%%%%%%%%%%%%%%%%%%%%%%%%%%%%%%%%%%%%%%%%%%%%%%%%%%%%%
en tous les
points\footnote{\, Les points sont représentés par leurs coordonnées
qui sont des systèmes de valeurs. 
}: %%%%%%%%%%%%%%%%%%%%%%%%%%%%%%%%%%%%%%%%%%%%%%%%%%%%%%%%%%%
$x_k', y_k', z_k'$ ($k=1 \dots
s$) qui satisfont les
\'equations:
\def\theequation{7}\begin{equation}
\label{eq-7-p-10}
\left\{
\aligned
J
\big(
x_\lambda',y_\lambda',z_\lambda';\
&
x_\mu',y_\mu',z_\mu'
\big)
=
J
\big(
x_\lambda,y_\lambda,z_\lambda;\
x_\mu,y_\mu,z_\mu
\big)
\\
&\
{\scriptstyle{(\lambda\,=\,1\,\cdots\,s\,-\,1\,;\,\,
\mu\,=\,\lambda\,+\,1\,\cdots\,s)}},
\endaligned\right.
\end{equation}
et qui sont procur\'es de telle sorte que le syst\`eme de valeurs:
\def\theequation{1'}\begin{equation}
x_1',y_1',z_1':\
x_2',y_2',z_2';\
\dots;\
x_s',y_s',z_s'
\end{equation}
reste dans un certain voisinage du syst\`eme de valeurs~\thetag{ 1}.

Si nous fixons\footnote{\, Fixer un ou
plusieurs points signifie considérer le
sous-groupe consistant seulement en les
transformations ponctuelles qui fixent ce ou ces points. 
} %%%%%%%%%%%%%%%%%%%%%%%%%%%%%%%%%%%%%%%%%%%%%%%%%%%%%%%%%%% 
maintenant les $s-1$
premiers points~\thetag{ 1}:
\[
x_k'
=
x_k,\ \
y_k'
=
y_k,\ \
z_k'
=
z_k
\ \ \ \ \ \ \
{\scriptstyle{(k\,=\,1\,\cdots\,s\,-\,1)}},
\]
alors les conditions~\thetag{ 7} pour lesquelles $\lambda$ et $\mu$
sont tous deux $< s$ sont satisfaites automatiquement, et il ne
reste plus que les $s-1$ conditions:
\def\theequation{8}\begin{equation}
\left\{
\aligned
J
\big(
x_\lambda,y_\lambda,z_\lambda;\
&
x_s',y_s',z_s'
\big)
=
J
\big(
x_\lambda,y_\lambda,z_\lambda;\
x_s,y_s,z_s
\big)
\\
&
{\scriptstyle{(\lambda\,=\,1,\,2\,\cdots\,s\,-\,1)}},
\endaligned\right.
\end{equation}
c'est-\`a-dire: apr\`es fixation de ces $s-1$ points-l\`a, le point
$x_s, y_s, z_s$ peut encore se transformer en tous les points $x_s',
y_s', z_s'$ qui satisfont les \'equations~\thetag{ 8} et qui se
trouvent dans un voisinage donn\'e de $x_s, y_s, z_s$. Mais puisque le
groupe: $X_1 f \dots X_mf$ est fini, lorsqu'on fixe un nombre
suffisant de points qui sont mutuellement en position g\'en\'erale, il
doit finalement se produire le cas que tous les points restent
g\'en\'eralement au repos; par cons\'equent, les \'equations~\thetag{
8} doivent \^etre constitu\'ees de telle sorte que, lorsque $s$ est
suffisamment grand, plus aucune famille continue de syst\`emes de
valeurs $x_s', y_s', z_s'$ ne les satisfait, et qu'on peut en tirer
les \'equations:
\def\theequation{9}\begin{equation}
x_s'
=
x_s,\ \
y_s'
=
y_s,\ \
z_s'
=
z_s.
\end{equation}
D'apr\`es le Tome~I, p.~490, ce cas se produit\footnote{\, Voici
l'argument.  Soit $G$ un groupe de transformations fini et continu 
comportant $m$
paramètres (essentiels) qui agit sur $R_n$.  Il existe au moins un
point $p_1$ en position générale qui n'est pas fixé par $G$.  Le
sous-groupe (d'isotropie) $G_1 \subset G$ des transformations qui
laissent $p_1$ au repos a donc $m_1 \leqslant m - 1$ paramètres. Si
$m_1 \geqslant 1$, il existe au moins un point $p_2$ en position
générale qui n'est pas fixé par $G_1$. Le sous-groupe $G_2 \subset
G_1$ des transformations fixant $p_2$ possède $m_2 \leqslant m_1 - 1
\leqslant m - 2$ paramètres, {\em etc.} 
} %%%%%%%%%%%%%%%%%%%%%%%%%%%%%%%%%%%%%%%%%%%%%%%%%%%%%%%%%%%
au plus tard lorsque $s =
m+1$, mais nous verrons cependant que sous les hypoth\`eses pos\'ees
ici, ce cas se produit d\'ej\`a pour un nombre inf\'erieur\footnote{\,
Les développements ultérieurs montreront que cela se produit en fait déjà
pour $s < 4$, meilleure inégalité que $s < 7 = m +1$ 
lorsque $m = 6$. 
} %%%%%%%%%%%%%%%%%%%%%%%%%%%%%%%%%%%%%%%%%%%%%%%%%%%%%%%%%%%
\`a $m+1$.

Si nous fixons seulement un point en position g\'en\'erale, par
exemple $x_1, y_1, z_1$, alors chaque autre point $x_2, y_2, z_2$ en
position g\'en\'erale peut encore
\^etre transform\'e en les $\infty^2$ points\footnote{\,
Le symbole courant
<<\,$\infty^k$\,>> dénote le nombre $k$ de paramètres
(réels ou complexes) dont un objet géométrique
(ou analytique) dépend effectivement. Chaque paramètre, 
susceptible de parcourir une infinité de valeurs continues, 
compte pour une et une seule puissance
du symbole <<\,$\infty$\,>>.  
} %%%%%%%%%%%%%%%%%%%%%%%%%%%%%%%%%%%%%%%%%%%%%%%%%%%%%%%%%%%
$x_2', y_2', z_2'$ qui satisfont l'\'equation:
\[
J
\big(
x_1,y_1,z_1;\
x_2',y_2',z_2'
\big)
=
J
\big(
x_1,y_1,z_1;\
x_2,y_2,z_2
\big).
\]
Et maintenant, comme notre groupe: $X_1 f \dots X_m f$ est
transitif, apr\`es fixation d'un point en position g\'en\'erale, les
param\`etres du groupe sont soumis exactement \`a trois conditions;
par cons\'equent le point $x_2, y_2, z_2$ ne peut \'evidemment occuper
encore $\infty^2$ positions que s'il reste encore au minimum deux
param\`etres arbitraires; donc notre groupe poss\`ede au moins cinq
param\`etres.

Il est clair qu'apr\`es fixation du point: $x_1, y_1, z_1$, les
$\infty^3$ points de l'espace peuvent se disposer sur les $\infty^1$
surfaces invariantes:
\def\theequation{10}\begin{equation}
J
\big(
x_1,y_1,z_1;\
x,y,z
\big)
=
\text{\rm const.}
\end{equation}
Relativement au groupe \`a six param\`etres des mouvements euclidiens,
qui satisfait \'evidemment toutes les conditions pos\'ees \`a la
page~\pageref{399}, ces $\infty^1$ surfaces-l\`a ne sont autres que
les $\infty^1$ sph\`eres~\thetag{ 10} centr\'ees au point: $x_1, y_1,
z_1$. Afin de pouvoir nous exprimer bri\`evement, nous voulons donc
appeller simplement {\sl pseudosph\`eres \label{402-0} centr\'ees au
point $x_1, y_1, z_1$ relatives au groupe: $X_1f \dots X_m f$} les
$\infty^1$ surfaces~\thetag{ 10}. Ensuite, nous pouvons aussi dire:
si un point $x_1, y_1, z_1$ qui est en position g\'en\'erale
vis-\`a-vis du groupe: $X_1 f \dots X_m f$ est fix\'e, alors chaque
autre point en position g\'en\'erale peut se mouvoir sur la
pseudosph\`ere centr\'ee en $x_1, y_1, z_1$ passant par lui.
La totalité 
\deutsch{Inbegriff} de toutes les pseudosph\`eres existantes forme
naturellement une famille de surfaces invariantes~\label{402} par le
groupe\footnote{\, En effet, chaque transformation
finie du groupe conserve la pseudodistance $J$, donc envoie
chaque pseudosphère centrée en un point sur la 
pseudosphère de même pseudorayon qui est centrée au point image. 
Par conséquent, l'ensemble de toutes les pseudosphères 
se transforme en lui-même. 
}: %%%%%%%%%%%%%%%%%%%%%%%%%%%%%%%%%%%%%%%%%%%%%%%%%%%%%%%%%%%
$X_1 f \dots X_m f$. On peut encore mentionner aussi que les
$\infty^1$ pseudosph\`eres~\thetag{ 10} peuvent \^etre d\'efinies par
une \'equation de Pfaff int\'egrable de la forme:
\def\theequation{11}\begin{equation}
\alpha
(x_1,y_1,z_1;\
x,y,z)\,dx
+
\beta\,dy
+
\gamma\,dz
=
0,
\end{equation}
qui est obtenue par diff\'erentiation\footnote{\, 
---\,\,par rapport aux trois variables
$x, y, z$\,\,---
} %%%%%%%%%%%%%%%%%%%%%%%%%%%%%%%%%%%%%%%%%%%%%%%%%%%%%%%%%%%
de~\thetag{ 10}, et dans
laquelle $x_1, y_1, z_1$ jouent le r\^ole de constantes.

\`A pr\'esent, imaginons-nous que deux points: $x_1, y_1, z_1$ et
$x_2, y_2, z_2$ sont fix\'es; alors, comme nous le savons\footnote{\, 
---\,\,d'après les équations~\thetag{ 7} et~\thetag{ 1'}\,\,---
}, %%%%%%%%%%%%%%%%%%%%%%%%%%%%%%%%%%%%%%%%%%%%%%%%%%%%%%%%%%%
tout autre
troisi\`eme point: $x_3, y_3, z_3$ en position g\'en\'erale peut
encore occuper toutes les positions: $x_3', y_3', z_3'$ qui satisfont
les deux \'equations:
\def\theequation{12}\begin{equation}
\left\{
\aligned
J
\big(
x_1,y_1,z_1;\
x_3',y_3',z_3'
\big)
&
=
J
\big(
x_1,y_1,z_1;\
x_3,y_3,z_3
\big)
\\
J
\big(
x_2,y_2,z_2;\
x_3',y_3',z_3'
\big)
&
=
J
\big(
x_2,y_2,z_2;\
x_3,y_3,z_3
\big),
\endaligned\right.
\end{equation}
et qui se trouvent dans un certain voisinage de $x_3, y_3, z_3$, donc
ce point peut encore occuper au moins $\infty^1$ positions.
Maintenant, comme par fixation des deux points, les param\`etres du
groupe sont soumis \`a cinq conditions, il en d\'ecoule que le nombre
de param\`etres du groupe ne peut en tout cas pas \^etre inf\'erieur
\`a six\footnote{\, 
Fixer un point de coordonnées $(x_1, y_1, z_1)$ absorbe au moins trois
paramètres du groupe, puisqu'il est supposé transitif.  Le second
point $(x_2, y_2, z_2)$ se meut alors, par l'action du sous-groupe
d'isotropie fixant le premier point, sur une pseudosphère de centre
$(x_1, y_1, z_1)$, avec deux degrés de liberté, et de manière
localement transitive. Fixer ensuite ce deuxième point absorbe au
moins deux paramètres supplémentaires du groupe.  Enfin, 
une mobilité comportant au moins un dernier
paramètre est encore possible sur l'intersection (généralement
non vide) des deux familles de pseudosphères
centrées en $(x_1, y_1, z_1)$ et en 
$(x_2, y_2, z_2)$. 
}. %%%%%%%%%%%%%%%%%%%%%%%%%%%%%%%%%%%%%%%%%%%%%%%%%%%%%%%%%%% 

Si les deux \'equations~\thetag{ 12} n'\'etaient pas ind\'ependantes
l'une de l'autre relativement \`a $x_3', y_3', z_3'$, la seconde
\'equation serait alors une cons\'equence de la premi\`ere, et il en
d\'ecoulerait \'evidemment que toutes les $s-1$
\'equations~\thetag{ 8} seraient aussi des cons\'equences de la premi\`ere
d'entre elles, aussi grand que l'on choisisse $s$; on ne pourrait
donc pas choisir $s$ assez grand pour que les \'equations~\thetag{ 9}
soient tir\'ees des \'equations~\thetag{ 8}, ou, ce qui revient au m\^eme,
il en d\'ecoulerait, m\^eme
si $s$ \'etait assez grand, que tous les
invariants de $s$ points ne pourraient pas se laisser exprimer au
moyen des invariants des paires de points, en contradiction avec ce qui
a \'et\'e dit ci-dessus. Par cons\'equent, nous pouvons conclure que les
\'equations~\thetag{ 12}, relativement \`a $x_3', y_3', z_3'$, sont
ind\'ependantes l'une par rapport \`a l'autre et que le point $x_3,
y_3, z_3$ peut seulement occuper encore $\infty^1$ positions apr\`es
fixation de $x_1, y_1, z_1$ et de $x_2, y_2, z_2$. En d'autres
termes, les $\infty^1$ pseudosph\`eres de centre: $x_2, y_2, z_2$
doivent couper les $\infty^1$ pseudosph\`eres de centre: $x_1, y_1,
z_1$ en les $\infty^2$ courbes qui sont d\'etermin\'ees par les deux
\'equations:
\def\theequation{13}\begin{equation}
\left\{
\aligned
J
\big(
x_1,y_1,z_1;\
x,y,z
\big)
&
=
\text{\rm const.}
\\
J
\big(
x_2,y_2,z_2;\
x,y,z
\big)
&
=
\text{\rm const.},
\endaligned\right.
\end{equation}
ou par le syst\`eme simultan\'e:
\def\theequation{14}\begin{equation}
\left\{
\aligned
\alpha(x_1,y_1,z_1;\
x,y,z)\,dx
+
\beta\,dy
+
\gamma\,dz
&
=
0
\\
\alpha(x_2,y_2,z_2;\
x,y,z)\,dx
+
\beta\,dy
+
\gamma\,dz
&
=
0.
\endaligned\right.
\end{equation}
En particulier, on obtient que les $\infty^1$ pseudosph\`eres de
centre: $x_1, y_1, z_1$ ne peuvent pas \^etre ind\'ependantes de leur
centre\footnote{\, 
---\,\,sinon l'intersection des deux familles de surfaces~\thetag{ 14}
ne se réduirait pas à des courbes\,\,---
}, %%%%%%%%%%%%%%%%%%%%%%%%%%%%%%%%%%%%%%%%%%%%%%%%%%%%%%%%%%%.
et donc qu'il y a au minimum $\infty^2$ pseudosph\`eres
diff\'erentes\footnote{\,
Le pseudorayon, à savoir la constante dans
l'équation $J ( p_1; p) = {\rm const.}$, constitue 
un premier paramètre (évident) pour la famille des pseudosphères. 
Les équations des pseudosphères ne pouvant pas être
indépendantes de leur centre $p_1$, elles dépendent
d'au moins un paramètre supplémentaire. Au total, 
il y a au moins deux paramètres. 
} %%%%%%%%%%%%%%%%%%%%%%%%%%%%%%%%%%%%%%%%%%%%%%%%%%%%%%%%%%%.

Enfin, si nous nous imaginons que trois points distincts: $x_k, y_k,
z_k$ ($k=1, 2, 3$) sont fix\'es, alors chaque autre quatri\`eme point
$x_4, y_4, z_4$ peut encore occuper toutes les positions $x_4', y_4',
z_4'$ dans un certain voisinage de $x_4, y_4, z_4$ qui satisfont les
trois \'equations:
\def\theequation{15}\begin{equation}
\left\{
\aligned
J
\big(
x_1,y_1,z_1;\
x_4',y_4',z_4'
\big)
&
=
J
\big(
x_1,y_1,z_1;\
x_4,y_4,z_4
\big)
\\
J
\big(
x_2,y_2,z_2;\
x_4',y_4',z_4'
\big)
&
=
J
\big(
x_2,y_2,z_2;\
x_4,y_4,z_4
\big)
\\
J
\big(
x_3,y_3,z_3;\
x_4',y_4',z_4'
\big)
&
=
J
\big(
x_3,y_3,z_3;\
x_4,y_4,z_4
\big).
\endaligned\right.
\end{equation}
Si maintenant, 
parmi ces trois \'equations, il n'y avait que deux d'entre elles qui
\'etaient ind\'ependantes relativement \`a $x_4', y_4', z_4'$, la
troisi\`eme d\'ecoulant par exemple des deux premi\`eres, alors parmi
les $s-1$ \'equations~\thetag{ 8}, les $s-3$ derni\`eres
d\'ecouleraient aussi toujours des deux premi\`eres, donc les
\'equations~\thetag{ 9} ne pourraient jamais \^etre tir\'ees des
\'equations~\thetag{ 8}, si grand que soit 
choisi $s$. Par cons\'equent, les
\'equations suivantes:
\[
x_4'
=
x_4,\ \ \
y_4'
=
y_4,\ \ \
z_4'
=
z_4
\]
doivent d\'ej\`a se tirer des \'equations~\thetag{ 15}, ce qui veut
dire qu'apr\`es fixation de trois points qui sont mutuellement en
position g\'en\'erale, tous les points de l'espace doivent
rester généralement
\renewcommand{\thefootnote}{\fnsymbol{footnote}}
au repos\footnote[1]{\, 
Entre autres choses, il d\'ecoule
encore de l\`a que les \'equations:
\def\theequation{A}\begin{equation}
J
\big(
x_k',y_k',z_k';\
x',y',z'
\big)
=
J
\big(
x_k,y_k,z_k;\
x,y,z
\big)
\ \ \ \ \ \ \
{\scriptstyle{(k\,=\,1,\,2,\,3)}}
\end{equation}
sont r\'esolubles par rapport \`a $x', y', z'$. Si on soumet les
quantit\'es: $x_k', y_k', z_k', x_k, y_k, z_k$ aux \'equations:
\[
\aligned
J
\big(
x_k',y_k',z_k';\
x_j',
&
y_j',z_j'
\big)
=
J
\big(
x_k,y_k,z_k;\
x_j,y_j,z_j
\big)
\\
&
{\scriptstyle{(k,\,j\,=\,1,\,2,\,3,\ \ \
k\,<\,j)}},
\endaligned
\]
et si on les interpr\`ete comme des param\`etres, la r\'esolution des
\'equations~\thetag{ A} par rapport \`a $x', y', z'$ repr\'esentera,
comme on s'en convaincra facilement, le groupe le plus g\'en\'eral de
l'espace pour lequel les deux points: $x_1, y_1, z_1; \ x_2, y_2,
z_2$ poss\`edent l'invariant: $J \big( x_1, y_1, z_1; \ x_2, y_2, z_2
\big)$. Ce groupe comprend \'evidemment le groupe: $X_1 f \dots X_m
f$, mais il n'est pas n\'ecessairement continu; n\'eanmoins, il peut
se d\'ecomposer en plusieurs familles s\'epar\'ees de transformations
continues ({\em voir}\, le Chap.~18 du Tome~I).

Nous voulons encore mentionner que l'on obtient les \'equations finies
du groupe: $X_1 f \dots X_m f$ \`a partir des \'equations~\thetag{
A}, lorsqu'on choisit la r\'esolution en $x', y', z'$ qui pour:
$x_k' = x_k$, $y_k' = y_k$, $z_k' = z_k$ ($k=1, 2, 3$) se r\'eduit aux
\'equations: $x' = x$, $y' = y$, $z' = z$.
}.%%%%%%%%%%%%%%%%%%%%%%%%%%%%%%%%%%%%%%%%%%%%%%%%%%%%%%%%%%%
\renewcommand{\thefootnote}{\arabic{footnote}}\footnote{\, Ici, dans
les applications, pour un groupe concret dont l'invariant $J$ est
connu explicitement en termes de fonctions algébriques ou
transcendantes, il est {\em a priori}\, possible, d'après ce
qui vient d'être vu, de résoudre les équations~\thetag{ A} par rapport
à $(x', y', z')$ comme fonctions de $(x, y, z)$ et des 18 paramètres
$x_k, y_k, z_k, x_k', y_k', z_k'$, $k = 1, 2, 3$.  En général, cela
peut donner une représentation des équations {\em finies}\, du groupe,
à ceci près que $18 - 6 = 9$ paramètres sont nécessairement superflus,
et que les équations finies ainsi obtenues peuvent parfois représenter
une famille qui ne contient pas l'identité, par exemple si les deux
<<\,repères\,>> formés des trois points $(x_k, y_k, z_k)$ et des trois
autres points $(x_k', y_k', z_k')$ sont disposés d'une façon spéciale.
Effectivement, un groupe de Lie $G$ peut posséder plusieurs
composantes connexes, par exemple le groupe des déplacements
euclidiens $E_3 ( \R) = E^+ ( \R) \cup E^- ( \R)$, suivant que
l'orientation est préservée ou inversée. Dans le Chap.~18 du Tome~I,
Engel et Lie élaborent une théorie générale des groupes de
transformation qui se décomposent en un nombre fini de groupes
connexes; seule la composante de l'identité est engendrée, via
l'exponentielle, par des transformations infinitésimales
$X_1, \dots, X_m$.  Cette stratégie éventuelle: trouver les invariants
possibles et en déduire les groupes possibles par le théorème des
fonctions implicites, ne sera pas adoptée.
En fait, Lie évite généralement d'écrire les
équations finies d'un groupe concret. 
} %%%%%%%%%%%%%%%%%%%%%%%%%%%%%%%%%%%%%%%%%%%%%%%%%%%%%%%%%%%
\label{404}

De ceci, il r\'esulte en premier lieu que $m$ doit \^etre
pr\'ecis\'ement \'egal \`a six. Et il suit en second lieu que les
$\infty^2$ courbes~\thetag{ 13} ne peuvent pas toutes \^etre contenues
dans les $\infty^1$ pseudosph\`eres de centre $x_3, y_3, z_3$, mais
que chacune de ces courbes n'a en g\'en\'eral qu'un point en commun
avec chacune de ces pseudosph\`eres\footnote{\, L'indépendance des
deux équations~\thetag{ 13} lorsque leurs centres sont mutuellement en
position générale équivaut à ce que l'intersection des pseudosphères
se réduise (en général) à une courbe. Si l'on considère une troisième
pseudosphère centrée en un autre point, l'intersection triple devient,
grâce à~\thetag{ 15} et au paragraphe qui suit, un simple point, ce
qui revient à dire que l'intersection de cette pseudosphère avec les
$\infty^2$ courbes d'intersection entre les deux premières familles de
pseudosphères, se réduit à un point.
}, %%%%%%%%%%%%%%%%%%%%%%%%%%%%%%%%%%%%%%%%%%%%%%%%%%%%%%%%%%%
et donc que le syst\`eme
simultan\'e~\thetag{ 14} ne peut pas \^etre ind\'ependant de $x_1,
y_1, z_1, x_2, y_2, z_2$. Ainsi, on a d\'emontr\'e que dans $R_3$, il
n'existe pas de famille de $\infty^2$ courbes qui est constitu\'ee de
telle sorte que dans toute pseudosph\`ere se trouvent $\infty^1$
courbes de cette famille. Et maintenant,
comme chaque famille de $\infty^2$
courbes qui remplit compl\`etement $R_3$, est 
repr\'esentable\footnote{\, 
Localement, 
toute famille régulière de courbes se redresse en
la famille 
$\{ x = {\rm const.}, \, 
y = {\rm const.} \}$ des droites parallèles à l'axe des
$z$, et ces droites sont les courbes
intégrales de tout champ de vecteurs de la forme
$\nu (x, y, z) \, \frac{ \partial }{\partial z}$ avec $\nu \neq 0$. 
}, %%%%%%%%%%%%%%%%%%%%%%%%%%%%%%%%%%%%%%%%%%%%%%%%%%%%%%%%%%%
comme
famille des courbes caract\'eristiques
de Monge d'une \'equation lin\'eaire aux
d\'eriv\'ees partielles:
\[
\lambda(x,y,z)\,
\frac{\partial f}{\partial x}
+
\mu(x,y,z)\,
\frac{\partial f}{\partial y}
+
\nu(x,y,z)\,
\frac{\partial f}{\partial z}
=
0,
\]
nous pouvons aussi dire ceci: les pseudosph\`eres ne sont pas 
enti\`erement surfaces int\'egrales d'une m\^eme \'equation
lin\'eaire aux d\'eriv\'ees partielles.

Cette propri\'et\'e des pseudosph\`eres qui vient d'\^etre \'etablie
permet encore de d\'eduire une propri\'et\'e importante des groupes
recherch\'es, \`a savoir, elle permet de d\'emontrer que dans aucun des
groupes recherch\'es ne peuvent se trouver deux transformations
infinit\'esimales qui poss\`edent les m\^emes courbes int\'egrales.

En effet, s'il y avait deux transformations infinit\'esimales du
groupe: $X_1 f \dots X_6 f$, par exemple: $X_1f$ et $X_2 f$, qui
ont exactement les m\^emes courbes int\'egrales, il existerait une
identit\'e de la forme:
\[
X_2f
\equiv
\varphi(x,y,z)\,
X_1f,
\]
o\`u $\varphi$ n'est pas une simple constante\footnote{\, Sinon, si
$X_2 = \lambda\, X_1$ avec $\lambda$ constant, l'algèbre de Lie
(espace vectoriel) $X_1, X_2, \dots, X_6$ ne serait pas de dimension
six. Donc la différentielle $d \varphi = \varphi_x dx + \varphi_y dy +
\varphi_z dz$ n'est pas identiquement nulle.
}. %%%%%%%%%%%%%%%%%%%%%%%%%%%%%%%%%%%%%%%%%%%%%%%%%%%%%%%%%%% 
Si nous formons maintenant les six
\'equations lin\'eaires aux d\'eriv\'ees partielles:
\def\theequation{16}\begin{equation}
X_kf
+
X_k^{(1)}f
=
0
\ \ \ \ \ \ \
{\scriptstyle{(k\,=\,1\,\cdots\,6)}}
\end{equation}
dont une solution commune:
\[
J
\big(
x,y,z;\
x_1,y_1,z_1
\big)
\]
d\'etermine imm\'ediatement\,\,---\,\,\'etant donn\'e une constante
arbitraire\,\,---\,\,les $\infty^1$ pseudosph\`eres de centre $x_1,
y_1, z_1$, alors nous remarquons que dans ces six \'equations sont
comprises les deux suivantes:
\[
\aligned
X_1f
&
+
X_1^{(1)}f
=
0
\\
\varphi(x,y,z)\,X_1f
&
+
\varphi(x_1,y_1,z_1)\,X_1^{(1)}f
=
0,
\endaligned
\]
d'o\`u suivent imm\'ediatement\footnote{\, 
Le déterminant $2 \times 2$ du système, 
égal à $\varphi ( x_1, y_1, z_1) - 
\varphi ( x, y, z)$, ne s'annule pas identiquement, 
puisque sa différentielle par
rapport aux variables $( x, y, z)$ est non nulle. 
} %%%%%%%%%%%%%%%%%%%%%%%%%%%%%%%%%%%%%%%%%%%%%%%%%%%%%%%%%%%
les \'equations:
\[
X_1f
=
0,
\ \ \ \ \ \ \
X_1^{(1)}f
=
0.
\]
Et maintenant, comme la premi\`ere de ces deux \'equations est
ind\'ependante de $x_1, y_1, z_1$, il se produit que sous les
hypoth\`eses pos\'ees, toutes les pseudosph\`eres de l'espace
satisfont une m\^eme \'equation lin\'eaire aux d\'eriv\'ees partielles
du premier ordre, ce qui, comme nous venons \`a l'instant de le
d\'emontrer, ne peut pas \^etre le cas\footnote{\,
Par exemple, une application de cette propriété sera
utilisée pour traiter, au Chap.~21, le contre-exemple~\thetag{ 23'}
à une assertion de Helmholtz. 
}. %%%%%%%%%%%%%%%%%%%%%%%%%%%%%%%%%%%%%%%%%%%%%%%%%%%%%%%%%%%

Par suite, nous pouvons maintenant \'enoncer 
la proposition suivante:

\medskip
{\bf Proposition~1.} \label{Satz-1}
{\em Si un groupe de transformations ponctuelles de $R_3$ continu et
fini est constitu\'e de telle sorte que, relativement \`a son action,
deux points poss\`edent un et un seul invariant, tandis que $s > 2$
points n'ont pas d'invariant essentiel, alors ce groupe est transitif
\`a six param\`etres et en outre, il ne comprend jamais deux
transformations infinit\'esimales dont les courbes int\'egrales
co\"{\i}ncident; de plus, la famille des pseudosph\`eres relatives au
groupe se compose au minimum de $\infty^2$ surfaces, et il n'y a pas
dans $R_3$ de famille doublement infinie de courbes qui produisent
toutes les pseudosph\`eres existantes.

}\medskip

Avant d'aller plus loin, nous voulons encore mettre sous une forme
appropri\'ee la condition d'apr\`es laquelle deux points ont un et un
seul invariant, relativement \`a un groupe transitif: $X_1 f \dots
X_6 f$ \`a six param\`etres.

Tout d'abord, il sort de la condition susnomm\'ee que les six
\'equations lin\'eaires aux d\'eriv\'ees partielles~\thetag{ 16}
poss\'edent une et une seule solution commune, et donc que parmi ces
\'equations, il doit s'en trouver exactement cinq\footnote{\,
Si l'on note $X_k = \sum_{ i=1}^3 \, \xi_{ ki} ( x) \,
\frac{ \partial}{ \partial x_i}$, pour $k = 1, \dots, 6$, 
les six équations~\thetag{ 16} d'inconnues des fonctions $J = J ( x,
y, z; \, x_1, y_1, z_1)$ s'écrivent
\[
{\textstyle{\sum_{i=1}^3\,\xi_{ki}(x,y,z)
\frac{\partial J}{\partial x_i}}}
+
{\textstyle{\sum_{i=1}^3\,\xi_{ki}(x_1,y_1,z_1)
\frac{\partial J}{\partial x_{1i}}}}
=
0
\ \ \ \ \ \ \ \ \ \ \ \
{\scriptstyle{(k\,=\,1\,\cdots\,6)}}.
\]
L'hypothèse qu'elles possèdent une seule solution 
équivaut à ce que le rang générique de la matrice $6 
\times 6$: 
\[
\left(
\begin{array}{cccccc}
\xi_{11}(x,y,z) & \cdots & \xi_{13}(x,y,z) & \xi_{11}(x_1,y_1,z_1)
& \cdots & \xi_{13}(x_1,y_1,z_1)
\\
\cdot\cdot & \cdots & \cdot\cdot & \cdot\cdot & \cdots & \cdot\cdot
\\
\xi_{61}(x,y,z) & \cdots & \xi_{63}(x,y,z) & \xi_{61}(x_1,y_1,z_1)
& \cdots & \xi_{63}(x_1,y_1,z_1)
\end{array}
\right)
\]
soit égal à cinq. La suite des raisonnements explore cette
hypothèse. 
} %%%%%%%%%%%%%%%%%%%%%%%%%%%%%%%%%%%%%%%%%%%%%%%%%%%%%%%%%%% 
qui sont ind\'ependantes les unes des autres. Par ailleurs, comme le
groupe: $X_1 f \dots X_6f$ est transitif, on aura par exemple que:
$X_1 f, X_2 f, X_3 f$ ne sont li\'es par aucune relation lin\'eaire
homog\`ene, tandis que: $X_4 f, X_5 f, X_6 f$ se laissent exprimer de
mani\`ere lin\'eaire\footnote{\, ---\,\,à coefficients {\em non}
constants, le terme <<\,linéaire\,>> ne possédant pas ici le sens
moderne\,\,---
} %%%%%%%%%%%%%%%%%%%%%%%%%%%%%%%%%%%%%%%%%%%%%%%%%%%%%%%%%%%
et homog\`ene au moyen de: $X_1 f, X_2 f,
X_3f$, à savoir:
\def\theequation{17}\begin{equation}
X_{3+k}f 
\equiv
\sum_{1\leqslant j\leqslant 3}\,
\varphi_{kj}(x,y,z)\,X_jf
\ \ \ \ \ \ \
{\scriptstyle{(k\,=\,1,\,2,\,3)}}.
\end{equation}
En tenant compte de ces identit\'es, nous pouvons clairement 
remplacer
les six \'equations~\thetag{ 16} par les suivantes:
\def\theequation{18}\begin{equation}
\left\{
\aligned
X_kf
&
+
X_k^{(1)}f
=
0
\\
\sum_{1\leqslant j\leqslant 3}\,
\big\{
\varphi_{kj}(x,y,z)\,X_jf
&
+
\varphi_{kj}(x_1,y_1,z_1)\,X_j^{(1)}f
\big\}
=
0
\\
&
{\scriptstyle{(k\,=\,1,\,2,\,3)}},
\endaligned\right.
\end{equation}
ou encore, par les suivantes:
\def\theequation{19}\begin{equation}
\left\{
\aligned
X_kf
&
+
X_k^{(1)}f
=
0
\\
\sum_{1\leqslant j\leqslant 3}\,
\big\{
\varphi_{kj}(x,y,z)
&
-
\varphi_{kj}(x_1,y_1,z_1)
\big\}
X_jf
=
0
\\
&
{\scriptstyle{(k\,=\,1,\,2,\,3)}}.
\endaligned\right.
\end{equation}

Sous les hypoth\`eses pos\'ees, parmi les six \'equations~\thetag{
19}, les trois premi\`eres sont r\'esolubles\footnote{\, Le groupe
étant transitif, la matrice $\big( \xi_{ ki} ( x_1, y_1, z_1) \big)_{
1\leqslant k \leqslant 3}^{ 1 \leqslant i \leqslant 3}$ formée des
coefficients des trois transformations infinitésimales $X_k^{ (1)}$,
$k= 1, 2, 3$, par rapport aux champs basiques $p_1, q_1, r_1$, est de
déterminant non nul.
} %%%%%%%%%%%%%%%%%%%%%%%%%%%%%%%%%%%%%%%%%%%%%%%%%%%%%%%%%%%
par rapport \`a $p_1, q_1, r_1$; en revanche, les trois derni\`eres
sont compl\`etement libres de $p_1, q_1, r_1$; comme parmi les
\'equations~\thetag{ 19}, il doit y en avoir exactement
cinq\footnote{\, ---\,\,puisque, rappelons-le, 
les $m = 6$ équations~\thetag{ 2}
doivent posséder une et une seule solution commune
qui constitue le seul invariant entre deux points\,\,---
} %%%%%%%%%%%%%%%%%%%%%%%%%%%%%%%%%%%%%%%%%%%%%%%%%%%%%%%%%%%
qui sont ind\'ependantes les unes
des autres, il est n\'ecessaire et suffisant que les trois derni\`eres
d'entre elles, que nous pouvons aussi \'ecrire sous la forme:
\def\theequation{20}\begin{equation}
X_{3+k}f
-
\sum_{1\leqslant j\leqslant 3}\,
\varphi_{kj}(x_1,y_1,z_1)\,X_jf
=
0
\ \ \ \ \ \ \
{\scriptstyle{(k\,=\,1,\,2,\,3)}},
\end{equation}
se r\'eduisent pr\'ecis\'ement \`a deux \'equations ind\'ependantes.
De l\`a, si nous prenons en consid\'eration le fait que les trois
transformations infinit\'esimales:
\[
X_{3+k}f
-
\sum_{1\leqslant j\leqslant 3}\,
\varphi_{kj}(x_1,y_1,z_1)\,X_jf,
\ \ \ \ \ \ \
{\scriptstyle{(k\,=\,1,\,2,\,3)}}
\]
laissent compl\`etement invariant\footnote{\, En effet,
d'après~\thetag{ 17}, au point $(x_1, y_1, z_1)$, ces trois
transformations infinitésimales s'annulent, donc les trois groupes à
un paramètre qu'elles engendrent via l'exponentielle laissent fixe le
point $(x_1, y_1, z_1)$.
} %%%%%%%%%%%%%%%%%%%%%%%%%%%%%%%%%%%%%%%%%%%%%%%%%%%%%%%%%%%
le point $x_1, y_1, z_1$ et qu'elles peuvent s'exprimer (quand ce
point est un point en position g\'en\'erale) comme combinaisons
lin\'eaires des transformations infinit\'esimales du groupe: $X_1 f
\dots X_6 f$ fixant le point en question ({\em voir}\footnote{\,
D'après cette proposition générale qui repose seulement sur des
considérations d'algèbre linéaire, si $X_k = \sum_{ i = 1}^n \,
\xi_{ ki} ( x_1, \dots, x_n) \, 
\frac{ \partial }{ \partial x_i}$, $k = 1, \dots, m$ sont
$m$ transformations infinitésimales indépendantes d'un groupe fini
continu quelconque à $m$ paramètres, et si $q$ désigne le rang
générique de la matrice $\big( \xi_{ ki} ( x) \big)_{ 1\leqslant k
\leqslant m}^{ 1 \leqslant i \leqslant n}$ formée de leurs
coefficients, alors localement au voisinage d'un point
en position générale et après
renumérotation éventuelle des $X_k$: 

\begin{itemize}

\smallskip\item[$\bullet$]
les $q$ premières transformations infinitésimales
$X_1, \dots, X_q$ ne sont liées par aucune relation de la
forme:
\[
\chi_1(x_1,\dots,x_n)\cdot X_1
+\cdots+
\chi_q(x_1,\dots,x_n)\cdot X_q
\equiv
0\,; 
\]

\smallskip\item[$\bullet$]
les $m - q$ transformations infinitésimales restantes $X_{ q+1},
\dots, X_m$ s'expriment comme combinaisons linéaires à coefficients
fonctionnels de $X_1, \dots, X_q$:
\[
X_{q+j}
\equiv
{\textstyle{\sum_{k=1}^q}}\,\varphi_{jk}(x_1,\dots,x_n)\cdot X_k
\ \ \ \ \ \ \ \ \ \ \ \ \
{\scriptstyle{(j\,=\,1\,\cdots\,m\,-\,q)}}\,;
\]

\smallskip\item[$\bullet$]
la sous-algèbre de Lie constituée des transformations infinitésimales
qui s'annulent en un point fixé $(x_1^0,
\dots, x_n^0)$ est précisément de dimension $m
- q$ et elle est engendrée par les $m - q$ transformations explicites:
\[
\label{isotropie-general}
X_{q+j}
-
{\textstyle{\sum_{k=1}^q}}\,\varphi_{jk}(x_1^0,\dots,x_n^0)\cdot X_k
\ \ \ \ \ \ \ \ \ \ \ \ \
{\scriptstyle{(j\,=\,1\,\cdots\,m\,-\,q)}}.
\]

\end{itemize}\smallskip

Appliquée ici avec $n = 3$, $m = 6$ et $q = 3$, cette proposition
donne les $m - q = 3$ transformations infinitésimales $X_{ 3+k} -
\sum_{ 1 \leqslant j \leqslant 3}\,
\varphi_{ kj} ( x_1, y_1, z_1) \, 
X_j$, $k=1, 2, 3$ indépendantes s'annulant au point $(x_1, y_1, z_1)$
et dont la matrice $3 \times 3$ des coefficients possède un rang générique
égal à $2$.  Cette propriété se transmet donc à tout système $Y_1,
Y_2, Y_3$ de trois transformations infinitésimales indépendantes du
groupe qui s'annulent en $(x_1, y_1, z_1)$: c'est la Proposition~2.
} %%%%%%%%%%%%%%%%%%%%%%%%%%%%%%%%%%%%%%%%%%%%%%%%%%%%%%%%%%%
Tome~I, p.~203, Proposition~7), nous obtenons
la:

\medskip
{\bf Proposition~2.}\label{Satz-2}
{\em Soit~{\rm :} $X_1 f \dots X_6 f$ un groupe \`a six param\`etres
de $R_3$ et soit}:
\[
\aligned
Y_kf
&
=
\alpha_k(x,y,z)\,p
+
\beta_k(x,y,z)\,q
+
\gamma_k(x,y,z)\,r
\\
&
\ \ \ \ \ \ \ \ \ \ \ \ \ \ \ \ \ \ \ \ \ \ \ \ \ \
{\scriptstyle{(k\,=\,1,\,2,\,3)}}
\endaligned
\]
{\em trois transformations infinit\'esimales quelconques du groupe qui
laissent invariant un point en position g\'en\'erale. Alors,
relativement au groupe: $X_1 f \dots X_6 f$, deux points poss\`edent
toujours un et un seul invariant lorsque, et seulement lorsque, le
d\'eterminant}:
\def\theequation{21}\begin{equation}
\left\vert
\begin{array}{ccc}
\alpha_1 & \beta_1 & \gamma_1 \\
\alpha_2 & \beta_2 & \gamma_2 \\
\alpha_3 & \beta_3 & \gamma_3 \\
\end{array}
\right\vert
\end{equation}
{\em s'annule identiquement, sans que tous ses sous-d\'eterminants
d'ordre deux s'annulent\footnote{\, Intuitivement, le sous-groupe à
trois paramètres $Y_1, Y_2, Y_3$ qui fixe un point donné déplace
encore transitivement tous les autres points sur la pseudosphère à
laquelle ils appartiennent.  Ces pseudosphères représentent donc les
seules familles de sous-variétés invariantes par l'action du
sous-groupe $Y_1, Y_2, Y_3$, et par conséquent, puisqu'elles sont de
codimension $1$, d'après une propriété générale, le rang générique de
la matrice $3 \times 3$ des coefficients de $Y_1, Y_2, Y_3$ doit être
égal à $3 - 1 = 2$.
}. %%%%%%%%%%%%%%%%%%%%%%%%%%%%%%%%%%%%%%%%%%%%%%%%%%%%%%%%%%%

}\medskip

En utilisant le crit\`ere qui est contenu dans cette proposition, nous
pouvons maintenant d\'eduire encore une autre propri\'et\'e importante
des groupes recherch\'es; toutefois, cette propri\'et\'e n'est
poss\'ed\'ee que par ceux de ces groupes qui sont imprimitifs.

\'Etant donn\'e un groupe imprimitif \`a six param\`etres ayant la
constitution demand\'ee, il peut donc arriver que ce groupe laisse
invariante une famille de $\infty^1$ surfaces: $\omega ( x, y, z) =
\text{\rm const.}$ Si ce cas se produit, nous pouvons
toujours nous imaginer que les variables $x, y, z$ sont choisies de
telle sorte que la famille de surfaces invariantes est repr\'esent\'ee
par l'\'equation: $x = \text{\rm const.}$; le groupe en question
est ensuite de la 
forme\footnote{\, 
%%%%%%%%%%%%%%%%%%%%%%%-------DEBUT--------%%%%%%%%%%%%%%%%%%%%%%%%%%%
Le groupe stabilise les hyperplans $\{ x = {\rm const.} \}$ 
si et seulement si le coefficient
$\xi_k$ de $p$ dans chaque $X_k$ ne dépend que de $x$. 
}: %%%%%%%%%%%%%%%%%%%%%%%%-----FIN-----%%%%%%%%%%%%%%%%%%%%%%%%%%%%%%%
\[
X_kf
=
\xi_k(x)\,p
+
\eta_k(x,y,z)\,q
+
\zeta_k(x,y,z)\,r
\ \ \ \ \ \ \
{\scriptstyle{(k\,=\,1\,\cdots\,6)}},
\]
o\`u $\xi_1 \dots \xi_6$ ne s'annulent
de toute fa\c con pas tous, parce
que sinon, le groupe serait intransitif\footnote{\,
Si tous les $\xi_k (x)$ étaient identiquement nuls, 
la direction $p$ manquerait.
}. %%%%%%%%%%%%%%%%%%%%%%%%%%%%%%%%%%%%%%%%%%%%%%%%%%%%%%%%%%%

En principe \deutsch{an und für sich}, trois cas seulement
sont {\em a priori}\, imaginables: en
effet, d'apr\`es\footnote{\,
Ce théorème, le premier du Tome~III, énonce
que tout groupe de transformations continu fini de la droite
des $x$ est de dimension 
$\leqslant 3$, et est équivalent soit à 
$\partial_x$ (groupe
des translations), soit à $\partial_x$, $x\, \partial _x$
(groupe affine), soit
à $\partial_x$, $x\, \partial_x$, $x^2\, \partial_x$
(groupe projectif). 
} %%%%%%%%%%%%%%%%%%%%%%%%%%%%%%%%%%%%%%%%%%%%%%%%%%%%%%%%%%%
le Th\'eor\`eme~1, p.~6, le groupe: $X_1 f \dots X_6
f$ peut transformer les surfaces: $x = \text{\rm const.}$ de une,
deux ou trois mani\`eres diff\'erentes
\deutsch{dreigliedrig transformiren}; nous allons cependant
d\'emontrer ici que seul le troisi\`eme cas se produit.

Si les surfaces: $x = \text{\rm const.}$ se transformaient seulement
d'une mani\`ere, on pourrait, d'apr\`es le
th\'eor\`eme cit\'e \`a l'instant
choisir la variable $x$ de telle mani\`ere que chaque $\xi_k ( x) \,
p$ prenne la forme: $a_k \, p$. Si l'on avait 
maintenant par exemple $a_1 \neq
0$, on pourrait introduire comme nouveau $y$ et nouveau $z$ les deux
solutions ind\'ependantes de l'\'equation\footnote{\,
Le changement de coordonnées est alors de la forme: 
\[
\overline{x}=x,
\ \ \ \ \ \ \ \ \
\overline{y}=\overline{y}(x,y,z),
\ \ \ \ \ \ \ \ \
\overline{z}=\overline{z}(x,y,z),
\]
et il induit la transformations suivante sur les
transformations infinitésimales basiques:
\[
p
=
\overline{p}
+
\overline{y}_x\overline{q}
+
\overline{y}_x\overline{r},
\ \ \ \ \ \ \ \ \
q
=
\overline{y}_y\overline{q}
+
\overline{z}_y\overline{r},
\ \ \ \ \ \ \ \ \
r
=
\overline{y}_z\overline{q}
+
\overline{z}_z\overline{r},
\] 
ce qui ne change rien au fait que les autres transformations
infinitésimales $X_2, \dots, X_6$ n'incorporent pas $p$: les nouvelles
transformations $\overline{ X}_2, \dots, \overline{ X}_6$
n'incorporent pas $\overline{ p}$.
}: %%%%%%%%%%%%%%%%%%%%%%%%%%%%%%%%%%%%%%%%%%%%%%%%%%%%%%%%%%%
$X_1 f = 0$, et rapporter
par cette op\'eration le groupe \`a la forme:
\[
\aligned
X_1f
=
p,
\ \ \ \ \ \ \
X_kf
&
=
\eta_k(x,y,z)\,q
+
\zeta_k(x,y,z)\,r
\\
&
{\scriptstyle{(k\,=\,2\,\cdots\,6)}}.
\endaligned
\]
Mais alors l'invariant des deux points: $x, y, z$ et $x_1, y_1, z_1$
serait n\'ecessairement de la forme: $x- x_1$, donc les
pseudosph\`eres de centre $x_1, y_1, z_1$ seraient repr\'esent\'ees
par: $x - x_1 = \text{\rm const.}$, ou plus simplement par: $x =
\text{\rm const.}$, c'est-\`a-dire qu'il y aurait en g\'en\'eral
seulement $\infty^1$ pseudosph\`eres distinctes, ce qui est exclu
d'apr\`es ce qui pr\'ec\`ede.

Si d'autre part les surfaces: $x = \text{\rm const.}$ se
transformaient de deux mani\`eres diff\'erentes, nous pourrions,
d'apr\`es le Th\'eor\`eme~1, p.~6, par un choix appropri\'e de $x$,
rapporter le groupe: $X_1 f \dots X_6 f$ \`a la forme:
\[
\aligned
X_1f
&
=
\ \
p
+
\eta_1(x,y,z)\,q
+
\zeta_1(x,y,z)\,r
\\
X_2f
&
=
xp
+
\eta_2(x,y,z)\,q
+
\zeta_2(x,y,z)\,r
\\
X_kf
&
=
\ \ \ \ \ \ \ \ \,
\eta_k(x,y,z)\,q
+
\zeta_k(x,y,z)\,r
\\
&
\ \ \ \ \ \ \ \ \ \ \ \ \ \ \ \ \
{\scriptstyle{(k\,=\,3\,\cdots\,6)}},
\endaligned
\]
et par cons\'equent, un point $x_1, y_1, z_1$ fix\'e en position
g\'en\'erale admettrait trois transformations infinit\'esimales de la
forme
\footnote{\,
La transformation $Y_1 := X_2 - x_1 X_1$ s'annule au 
point $p_1$. Les quatre transformations restantes 
$X_3, X_4, X_5, X_6$ n'incorporant que 
$q$ et $r$, il existe deux combinaisons
linéaires indépendantes 
à coefficients constants $Y_2$ et $Y_3$ qui s'annulent
aussi en $p_1$. 
}: %%%%%%%%%%%%%%%%%%%%%%%%%%%%%%%%%%%%%%%%%%%%%%%%%%%%%%%%%%%
\[
\aligned
Y_1f
&
=
(x-x_1)\,p
+
\overline{\eta}_1\,q
+
\overline{\zeta}_1\,r
\\
Y_2f
&
=
\ \ \ \ \ \ \ \ \ \ \ \ \ \ \ \ \ \ \ \ \
\overline{\eta}_2\,q
+
\overline{\zeta}_2\,r
\\
Y_3f
&
=
\ \ \ \ \ \ \ \ \ \ \ \ \ \ \ \ \ \ \ \ \
\overline{\eta}_3\,q
+
\overline{\zeta}_3\,r,
\endaligned
\]
et le d\'eterminant correspondant:
\[
\left\vert
\begin{array}{ccc}
x-x_1 & \overline{\eta}_1 & \overline{\zeta}_1
\\
0 & \overline{\eta}_2 & \overline{\zeta}_2
\\
0 & \overline{\eta}_3 & \overline{\zeta}_3
\end{array}
\right\vert
=
(x-x_1)
\left\vert
\begin{array}{cc}
\overline{\eta}_2 & \overline{\zeta}_2
\\ 
\overline{\eta}_3 & \overline{\zeta}_3
\end{array}
\right\vert
\]
devrait, d'apr\`es la Proposition~2, s'annuler identiquement. Mais
alors les deux transformations infinit\'esimales $Y_2 f$ et $Y_3 f$
seraient li\'ees par une identit\'e de la forme:
\[
Y_3f
\equiv
\omega(x,y,z)\,Y_2f
\]
et poss\`ederaient par suite des courbes int\'egrales en commun, ce
qui est exclu d'apr\`es la Proposition~1.

Ainsi, on a d\'emontr\'e que les surfaces invariantes $x = \text{\rm
const.}$ doivent \^etre transform\'ees par l'action d'un groupe \`a
trois param\`etres:

\medskip
{\bf Proposition~3.} \label{Satz-3}
{\em Si, relativement \`a un groupe \`a six param\`etres de $R_3$,
deux points poss\`edent un et un seul invariant, tandis que $s > 2$
points n'ont pas d'invariant essentiel, alors le groupe transforme de
trois mani\`eres diff\'erentes chaque famille de $\infty^1$ surfaces
qu'il laisse \'eventuellement invariantes.}

\medskip

Nous allons maintenant d\'eterminer tous les groupes transitifs \`a
six param\`etres de $R_3$ qui poss\`edent certaines propri\'et\'es
\'enonc\'ees dans les Propositions~1 et~3, donc pour nous exprimer de
mani\`ere plus pr\'ecise: nous cherchons maintenant tous les groupes
transitifs: $X_1 f \dots X_6 f$ de $R_3$ \`a six param\`etres tels que:

\smallskip

{\small\sf premi\`erement}~\label{408} deux points ont un et un seul
invariant; tels que {\small\sf deuxi\`emement} deux transformations
infinit\'esimales ind\'ependantes n'ont jamais les m\^emes courbes
int\'egrales et tels que {\small\sf troisi\`emement} chaque famille
invariante \'eventuelle de $\infty^1$ surfaces est transform\'ee de
trois mani\`eres diff\'erentes.

\smallskip

Afin de r\'esoudre compl\`etement le probl\`eme pos\'e auparavant,
p.~\pageref{399}, nous devrons encore d\'ecider au final, quant
\`a chacun des groupes qui se pr\'esentera, si $s > 2$ points ont un
invariant essentiel ou non, car alors, seuls les groupes pour lesquels
$s >2$ points n'ont pas d'invariant essentiel seront int\'eressants
pour nous. La d\'ecision \`a ce sujet sera essentiellement facilit\'ee
par la proposition suivante.

\medskip
{\bf Proposition~4.} \label{Satz-4}
{\em Si, relativement \`a un groupe transitif \`a six param\`etres:
$X_1 f \dots X_6 f$ de $R_3$, deux points $x, y, z$ et $x_1, y_1, z_1$
poss\`edent un et un seul invariant: $J \big( x, y, z; \ x_1, y_1,
z_1 \big)$, alors $s > 2$ points n'ont aucun invariant essentiel
\text{\rm si et seulement si:} premi\`erement, au minimum $\infty^2$
pseudosph\`eres diff\'erentes appartiennent au groupe; et
deuxi\`emement, il n'existe pas de famille de $\infty^2$ courbes qui
sont invariantes\footnote{\,  
Dans la Proposition~1, les familles de $\infty^2$ 
courbes exclues n'étaient
pas forcément invariantes par l'action du groupe. 
Mais les deux énoncés sont en fait équivalents, parce que les
pseudosphères elles-mêmes et les familles
de leurs courbes d'intersection deux à deux sont
évidemment invariantes par le groupe.  
} %%%%%%%%%%%%%%%%%%%%%%%%%%%%%%%%%%%%%%%%%%%%%%%%%%%%%%%%%%%
par le groupe et au moyen desquelles sont produites
toutes les pseudosph\`eres existantes.}

\medskip

On se convainc sans difficult\'e de la justesse de cette
proposition. En effet, la Proposition~1 p.~\pageref{Satz-1} montre
imm\'ediatement que dans la Proposition~4, les conditions indiqu\'ees
sont n\'ecessaires; qu'elles soient aussi suffisantes, on peut s'en
rendre compte comme suit.\label{409}

Comme deux points doivent poss\'eder un et un seul invariant essentiel
et comme il doit y avoir au minimum $\infty^2$ pseudosph\`eres
diff\'erentes, il est clair que deux pseudosph\`eres, dont les
centres: $x_1, y_1, z_1$ et $x_2, y_2, z_2$ sont mutuellement en
position g\'en\'erale, n'ont en g\'en\'eral qu'une seule courbe en
commun, et donc que les deux fonctions:
\[
J
\big(
x_1,y_1,z_1;\
x_3,y_3,z_3
\big),
\ \ \ \ \ \ \
J
\big(
x_2,y_2,z_2;\
x_3,y_3,z_3
\big),
\]
que nous voulons \'ecrire en abr\'eg\'e $J_{ 1, 3}$ et $J_{ 2, 3}$, sont
ind\'ependantes l'une de l'autre relativement \`a $x_3, y_3, z_3$. De
l\`a il suit manifestement que les trois fonctions: $J_{ 1, 2}$, $J_{
1, 3}$, $J_{ 2, 3}$ sont des solutions ind\'ependantes des
\'equations:
\[
X_k^{(1)}f
+
X_k^{(2)}f
+
X_k^{(3)}f
=
0
\ \ \ \ \ \ \
{\scriptstyle{(k\,=\,1\,\cdots\,6)}};
\]
mais comme, sous les hypoth\`eses pos\'ees, ces \'equations sont
tr\`es certainement r\'esolubles par rapport \`a six des neuf
quotients diff\'erentiels $p_\nu$, $q_\nu$, $r_\nu$ ($\nu = 1, 2, 3$),
en cons\'equence de quoi elles n'ont que trois solutions communes
ind\'ependantes, on a donc d\'emontr\'e que, sous les hypoth\`eses
pos\'ees, tous les invariants de trois points se laissent exprimer au
moyen des invariants des paires de points, et donc que trois points
n'ont pas d'invariant essentiel.

Par ailleurs, si trois pseudosph\`eres dont les centres: $x_\nu,
y_\nu, z_\nu$ ($\nu = 1, 2, 3$) sont mutuellement en position
g\'en\'erale ne s'intersectaient pas toujours en un seul point, mais
poss\'edaient toujours une courbe entière en commun, alors les
pseudosph\`eres centr\'ees en deux points quelconques parmi ces points
d\'etermineraient une famille de $\infty^2$ courbes qui serait
constitu\'ee de telle sorte que sur chaque pseudosph\`ere seraient
dispos\'ees $\infty^1$ courbes de la famille; et maintenant, comme
cette famille de courbes serait \'evidemment ind\'ependante de la
position de ces trois points-l\`a, elle resterait \'egalement
invariante par le groupe: $X_1 f \dots X_6 f$, puisque la famille de
toutes les pseudosph\`eres est elle-m\^eme invariante 
(\cf~p.~\pageref{402}), alors que dans notre proposition, l'existence
d'une telle famille de courbes invariantes est cependant
express\'ement exclue. Ainsi, on a d\'emontr\'e que trois
pseudosph\`eres ayant trois centres distincts n'ont en g\'en\'eral
qu'un seul point en commun. Si l'on entend maintenant par $x_4, y_4,
z_4$ un quatri\`eme point quelconque en position g\'en\'erale, les
trois fonctions: $J_{ 1, 4}$, $J_{ 2, 4}$, $J_{ 3, 4}$ seront tr\`es
certainement ind\'ependantes l'une de l'autre relativement \`a $x_4,
y_4, z_4$. Si de surcro\^{\i}t, on tient compte du fait que quatre
points ont pr\'ecis\'ement six invariants ind\'ependants relativement
au groupe: $X_1 f \dots X_6 f$, et aussi du fait que: $J_{ 1, 2},
J_{ 1, 3}, J_{ 2, 3}, J_{ 1, 4}, J_{ 2, 4}, J_{ 3, 4}$ constituent
d\'ej\`a six invariants tels, on obtient que, sous les hypoth\`eses
pos\'ees, quatre points n'ont pas non plus d'invariant essentiel. Enfin,
il est d\'esormais tout aussi imm\'ediatement clair que $s
> 4$ points sont d\'epourvus d'invariant essentiel; notre proposition
est donc compl\`etement d\'emontr\'ee.

\bigskip

Maintenant, nous allons tout d'abord d\'eterminer tous les groupes
primitifs qui poss\`edent les propri\'et\'es \'enonc\'ees
p.~\pageref{408} et ensuite, nous d\'eterminerons ceux qui sont
imprimitifs.

\HEAD{Groupes de $R_3$ pour lesquels deux points ont un seul
invariant.}{Division\,\,V.\,\,\,Chapitre\,\,20.\,\,\,\S\,\,86.}

\sectiondritterV{\sf\S\,\,86.
\\
Groupes primitifs parmi les groupes recherch\'es.}
\label{S-86}
\setcounter{footnote}{0}

Parmi les groupes primitifs \`a six param\`etres de $R_3$, il n'y a
apr\`es tout que deux types\footnote{\,
%%%%%%%%%%%%%%%%%%%%%%%-------DEBUT--------%%%%%%%%%%%%%%%%%%%%%%%%%%%
La référence concerne le Théorème~9 du Volume~III
qui classifie tous les groupes continus finis de transformations
holomorphes locaux sur un espace à trois
dimensions qui sont primitifs. 
} %%%%%%%%%%%%%%%%%%%%%%%%-----FIN-----%%%%%%%%%%%%%%%%%%%%%%%%%%%%%%%
({\em voir} Th\'eor\`eme~9, p.~139),
%%%\Fill
lesquels sont repr\'esent\'es premi\`erement par le groupe:
\def\theequation{22}\begin{equation}
\boxed{
\ \,
p,\ \
q,\ \
r,\ \
xq-yp,\ \
yr-zq,\ \ 
zp-xr\ \
}
\end{equation}
des mouvements euclidiens, et deuxi\`emement par le groupe projectif
\`a six param\`etres:
\def\theequation{23}\begin{equation}
\boxed{
\aligned
\
p+xU,\ \
q+yU,\ \
r+zU\
\\
xq-yp,\ \
yr-zq,\ \
zp-xr\
\endaligned
}
\end{equation}
de la surface du second degr\'e: $x^2 + y^2 + z^2 + 1 = 0$, o\`u,
comme \`a l'accoutum\'ee, on \'ecrit $U$ pour $xp + yq + zr$.

Comme on le sait, relativement au groupe~\thetag{ 22}, deux points
poss\`edent un et un seul invariant, \`a savoir notre distance au sens
usuel:
\[
(x_1-x_2)^2
+
(y_1-y_2)^2
+
(z_1-z_2)^2.
\]
Il en va de m\^eme pour le groupe~\thetag{ 23}, 
et pour pr\'eciser, l'invariant en question:
\begin{footnotesize}\label{410}
\[
\frac{
(x_1\!\!-\!x_2)^2
\!\!+\!
(y_1\!\!-\!y_2)^2
\!\!+\!
(z_1\!\!-\!z_2)^2
\!\!+\!
[x_1y_2\!\!-\!y_1x_2]^2
\!\!+\!
[y_1z_2\!\!-\!z_1y_2]^2
\!\!+\!
[z_1x_2\!\!-\!x_1z_2]^2
}
{
\{1
\!+\!
x_1x_2
\!+\!
y_1y_2
\!+\!
z_1z_2
\}^2}
\]
\end{footnotesize}

\noindent
constitue, comme on le sait, le birapport entre $x_1, y_1, z_1; \ x_2,
y_2, z_2$ et les deux points que la droite de liaison
\deutsch{Verbindungslinie} entre $x_1, y_1, z_1$ et 
$x_2, y_2, z_2$ d\'ecoupent sur la surface invariante du second degr\'e.
%%%\Mathematiques

En outre, il est d\'ej\`a connu depuis longtemps que pour ces deux
groupes, les invariants de plusieurs points quelconques peuvent
s'exprimer au moyen des invariants des paires de points, et donc que
$s > 2$ points n'ont pas d'invariant essentiel; si on ne connaissait
pas encore cette propri\'et\'e, on pourrait la conclure
imm\'ediatement de la Proposition~4, p.~\pageref{Satz-4}, puisque,
premi\`erement: pour chacun de ces deux groupes, deux points
poss\`edent un et un seul invariant et deuxi\`emement: \'etant
primitifs, ces deux groupes ne laissent invariante absolument aucune
famille de $\infty^2$ courbes.

Le cas des groupes primitifs est ainsi achev\'e; nous voulons
seulement mentionner encore que pour chacun des deux groupes
trouv\'es, une \'equation du second degr\'e:
\[
\alpha_{11}\,dx^2
+
\alpha_{22}\,dy^2
+
\alpha_{33}\,dz^2
+
2\alpha_{12}\,dxdy
+
2\alpha_{23}\,dydz
+
2\alpha_{31}\,dzdx,
\]
dont le d\'eterminant ne s'annule pas identiquement,
reste invariante, et que les $\infty^2$ directions passant
par un point fix\'e en position g\'en\'erale sont transform\'ees
par une action \`a trois 
param\`etres\footnote{\,
%%%%%%%%%%%%%%%%%%%%%%%-------DEBUT--------%%%%%%%%%%%%%%%%%%%%%%%%%%%
Pour chacun des deux groupes~\thetag{ 22}
et~\thetag{ 23}, l'origine $0$ est un point de position 
générale, puisqu'ils sont tous deux transitifs. 
Les trois transformations infinitésimales
$xq - yp$, $yr - zq$ et 
$zp - xr$ qu'ils possèdent en commun 
s'annulent en $0$ et elles engendrent 
le groupe complet des rotations (complexes ou réelles)
autour de $0$, qui est évidemment transitif sur le 
plan projectif des $\infty^2$ éléments linéaires
passant par $0$.   
}. %%%%%%%%%%%%%%%%%%%%%%%%-----FIN-----%%%%%%%%%%%%%%%%%%%%%%%%%%%%%%%

\HEAD{Groupes de $R_3$ pour lesquels deux points ont un seul
invariant.}{Division\,\,V.\,\,\,Chapitre\,\,20.\,\,\,\S\,\,87.}

\sectiondritterV{\sf\S\,\,87.
\\
Groupes imprimitifs parmi les groupes recherch\'es.}
\label{S-87}

Tous les groupes concern\'es sont transitifs \`a six
param\`etres. Maintenant, si l'on fixe un point en position
g\'en\'erale par rapport \`a l'action d'un groupe transitif \`a six
param\`etres, la vari\'et\'e \`a deux dimensions des $\infty^2$
\'el\'ements lin\'eaires passant par ce point sera transform\'ee par
un groupe projectif qui comporte au plus trois param\`etres. Par
ailleurs, nous savons d'apr\`es l'\'enum\'eration de tous les groupes
projectifs du plan (p.~106 sq.) que chaque groupe projectif du plan
qui ne comporte pas plus de trois param\`etres, soit laisse un point
invariant, soit se compose des $\infty^3$ transformations projectives
qui laissent invariante une quadrique non d\'eg\'en\'er\'ee. Nous
pouvons conclure de l\`a que pour les groupes transitifs de $R_3$ \`a
six param\`etres, seulement deux cas peuvent se produire: ou bien,
apr\`es fixation d'un point en position g\'en\'erale, il existe un
\'el\'ement lin\'eaire traversant ce point qui reste en m\^eme temps
fix\'e, ou bien il existe une quadrique non d\'eg\'en\'er\'ee qui est
invariante et les $\infty^2$ \'el\'ements lin\'eaires traversant ce
point se transforment de trois mani\`eres diff\'erentes.
%%%\Mathematiques

Ce dernier cas ne peut pas se produire ici, 
%%%\Mathematiques
parce que les groupes \`a
six param\`etres en question seraient alors primitifs\footnote{\
D\'ej\`a \'etudi\'es et traités il y a un instant 
au \S~86.} (\cf p.~123 sq.),
et par cons\'equent, chaque groupe poss\'edant la constitution ici
exig\'ee doit laisser invariant
un syst\`eme simultan\'e:
\[
\frac{dx}{\alpha(x,y,z)}
=
\frac{dy}{\beta(x,y,z)}
=
\frac{dz}{\gamma(x,y,z)}
\]
ou, ce qui revient au m\^eme, une famille de $\infty^2$ courbes:
\[
\varphi(x,y,z)
=
\text{\rm const.},
\ \ \ \ \ \ \
\psi(x,y,z)
=
\text{\rm const.}
\]
Si nous choisissons alors 
les coordonn\'ees $x, y, z$ de telle sorte que la
famille invariante de courbes re\c coive
la forme: $x = \text{\rm const.}$, $y =
\text{\rm const.}$, alors tous les groupes
recherch\'es sont de la forme:
\[
\aligned
X_kf
=
\xi_k(x,y)\,p
&
+
\eta_k(x,y)\,q
+
\zeta_k(x,y,z)\,r
\\
&
{\scriptstyle{(k\,=\,1\,\cdots\,6)}}.
\endaligned
\]

D'apr\`es le Tome~1, p.~307, Proposition~4, les six transformations
infinit\'esimales r\'eduites:
\[\label{412a}
\overline{X}_kf
=
\xi_k(x,y)\,p
+
\eta_k(x,y)\,q
\ \ \ \ \ \ \
{\scriptstyle{(k\,=\,1\,\cdots\,6)}}
\]
produisent un groupe dans l'espace des deux variables $x, y$ qui est
isomorphe-holo\'edrique, ou 
isomorphe-m\'ero\'edrique\footnote{\,
Un espace vectoriel $F$ de transformations infinitésimales est dit
{\sl isomorphe-holoédrique} \deutsch{holoedrisch isomorph} à un autre
espace vectoriel $E$ de transformations infinitésimales lorsqu'il
existe une application linéaire respectant les crochets (et naturelle
en un certain sens) {\em bijective} de $E$ dans $F$. Si l'application
(naturelle) de $E$ dans $F$ est seulement surjective, mais pas
bijective, $F$ est dit {\sl isomorphe-méroédrique}
\deutsch{meroedrisch isomorph} à $E$, ou {\sl isomorphe}
\deutsch{isomorph} (tout court) à $E$.  Ici, l'application naturelle
de <<\,projection\,>> qui, à
chaque transformation infinitésimale $X_k$, associe sa réduite
$\overline{ X}_k$ respecte les crochets, car $\big[ \overline{ X}_j,
\, \overline{ X}_k \big] = \sum_{ l=1}^6 \, c_{jkl}\, \overline{ X}_l$
découle en effet de de $\big[ X_j, \, X_k \big] = \sum_{ l=1}^6 \,
c_{jkl}\, X_l$ (avec les mêmes constantes de structure), mais elle
n'est pas forcément injective.  }
%%%%%%%%%%%%%%%%%%%%%%%%%%%%%%%%%%%%%%%%%%%%%%%%%%%%%%%%%%% au groupe:
$X_1 f \dots X_6 f$. Nous voulons tout d'abord consid\'erer ce groupe
r\'eduit.

Si le groupe r\'eduit comportait moins de cinq param\`etres, le groupe
$X_1 f \dots X_6 f$ comprendrait deux transformations
infinit\'esimales de la forme:
\[
\omega_1(x,y,z)\,r,
\ \ \ \ \ \ \
\omega_2(x,y,z)\,r,
\]
donc il comprendrait deux transformations infinit\'esimales ayant les
m\^emes courbes int\'egrales, ce qui, d'apr\`es
la Proposition~1, p.~\pageref{Satz-1}, 
n'est pas possible, et par cons\'equent, le groupe
r\'eduit: $\overline{ X}_1 f \dots \overline{ X}_6 f$ doit
comporter soit cinq, soit six param\`etres. Ainsi, il faut seulement
se demander encore quelles sont les formes diff\'erentes que le groupe
r\'eduit peut prendre.

Si le groupe: $\overline{ X}_1 f \dots \overline{ X}_6 f$ dans
l'espace des deux variables $x, y$ est primitif, alors, par un choix
appropri\'e de $x$ et de $y$, il peut toujours, d'apr\`es la page~71
%%%\Fill
recevoir l'une des deux formes
\footnote{\,
La référence concerne le Théorème~6
du Volume~III qui classifie tous
les groupes continus finis de transformations holomorphes
locaux sur un espace à deux dimensions. 
Le premier est le groupe affine; le second est le groupe
spécial affine. 
}: %%%%%%%%%%%%%%%%%%%%%%%%%%%%%%%%%%%%%%%%%%%%%%%%%%%%%%%%%%%
\def\theequation{24}\begin{equation}\label{412b}
\left\{
\aligned
&
p,\ \
q, \ \
xq,\ \
xp-yq,\ \
yp,\ \
xp+yq\ \
\\
&
p,\ \
q,\ \
xq,\ \
xp-yq,\ \
yp.
\endaligned\right.
\end{equation} 
D'un autre c\^ot\'e, si le groupe: $\overline{ X}_1 f \dots
\overline{ X}_6 f$ est imprimitif, il laisse alors invariante,
lorsqu'on l'interpr\`ete comme groupe du plan $x, y$, au moins une
famille de $\infty^1$ courbes: $\varphi ( x, y) = \text{\rm const.}$,
mais ensuite alors, dans l'espace des $x, y, z$, l'\'equation $\varphi
( x, y) = \text{\rm const.}$ repr\'esente \'evidemment une famille de
$\infty^1$ surfaces qui reste invariante\footnote{\, Analytiquement,
les surfaces $\varphi (x, y) = {\rm const.}$ sont invariantes par le
groupe $\overline{ X}_k$, $k = 1, \dots, 6$ si et seulement si
$\overline{ X}_k ( \varphi)$ s'exprime comme une certaine fonction
$\omega_k ( \varphi)$ de la {\em même} fonction $\varphi$, pour tout
$k = 1, \dots, 6$. Puisque $\varphi$ est indépendant de $z$, il est
alors évident qu'on a de même $X_k ( \varphi) = \overline{ X}_k (
\varphi) =
\omega_k ( \varphi)$, $k = 1, \dots, 6$. Géométriquement, la
famille des surfaces réglées formées au-dessus des courbes $\varphi (
x, y) = {\rm const.}$ par ajout en chaque point d'une droite parallèle
à l'axe des $z$ reste invariante par les $X_k$, puisqu'ils se
déduisent des $\overline{ X}_k$ par simple ajout d'une
transformation $\zeta \, \partial_z$ parallèle
à l'axe des $z$. 

} %%%%%%%%%%%%%%%%%%%%%%%%%%%%%%%%%%%%%%%%%%%%%%%%%%%%%%%%%%%
par le groupe: $X_1 f \dots
X_6 f$. Maintenant, comme d'apr\`es la Proposition~3,
p.~\pageref{Satz-3} le groupe: $X_1 f \dots X_6 f$ doit transformer
de trois mani\`eres diff\'erentes 
\deutsch{dreigliedrig transformiren}
toute famille de $\infty^1$ surfaces
qu'il laisse invariante, il s'ensuit que le groupe r\'eduit:
$\overline{ X}_1 f \dots \overline{ X}_6 f$, interpr\'et\'e comme
groupe du plan des $x, y$, transforme de trois mani\`eres
diff\'erentes toute famille de $\infty^1$ courbes qu'il laisse
invariante. Mais parmi les groupes imprimitifs du plan \`a cinq ou \`a
six param\`etres, il n'y en a proportionnellement qu'un petit nombre
qui satisfont \`a cette exigence, comme le montre le Tableau p.~71
sq.\footnote{\,
%%%%%%%%%%%%%%%%%%%%%%%-------DEBUT--------%%%%%%%%%%%%%%%%%%%%%%%%%%%
\`A nouveau, il s'agit du Théorème~6. 
}, %%%%%%%%%%%%%%%%%%%%%%%%-----FIN-----%%%%%%%%%%%%%%%%%%%%%%%%%%%%%%%
%%%\Fill
et pour pr\'eciser, par un choix appropri\'e de $x$ et de $y$,
ceux qui ont six param\`etres peuvent \^etre rapport\'es \`a l'une des
trois formes:
\def\theequation{25}\begin{equation}
\left\{
\aligned
&
q,\ \
xq,\ \
yq,\ \
p,\ \
xp,\ \
x^2p+xyq
\\
&
q,\ \
xq,\ \
x^2q,\ \
p,\ \
xp+yq,\ \
x^2p+2xyp
\\
&
q,\ \
yq,\ \
y^2q,\ \
p,\ \
xp,\ \
x^2p,
\endaligned\right.
\end{equation}
et ceux qui ont cinq param\`etres \`a la forme:
\def\theequation{26}\begin{equation}
q,\ \
xq,\ \
p,\ \
2xp+yq,\ \
x^2p+xyq.
\end{equation}

Ainsi, nous avons trouv\'e toutes les formes possibles du groupe
r\'eduit: $\overline{ X }_1 f \dots \overline{ X }_6 f$
et nous avons seulement \`a 
d\'eterminer encore, pour
chacun de ces groupes r\'eduits, 
quels sont les groupes \`a six param\`etres: $X_1 f
\dots X_6 f$ en les trois variables $x, y, z$ qui produisent de tels
groupes r\'eduits. Avec cela, nous devrons naturellement tenir compte
du fait que le groupe: $X_1 f \dots X_6 f$ ne doit pas comporter
deux transformations infinit\'esimales ayant les m\^emes courbes
int\'egrales. Quand nous aurons d\'etermin\'e les diff\'erentes formes
du groupe: $X_1 f \dots X_6 f$ qui satisfait cette derni\`ere
condition, nous devrons finalement examiner 
encore, pour chacun des groupes trouv\'es,
si deux points poss\`edent
un invariant, et si $s > 2$ points n'ont aucun invariant essentiel.

Nous supposerons en premier lieu que le groupe r\'eduit: $\overline{
X}_1 f \dots \overline{ X}_6 f$ a six\footnote{\,
Il s'agit des quatre
groupes réduits~$(24)_1$ et $(25)_1$, $(25)_2$, $(25)_3$. 
} %%%%%%%%%%%%%%%%%%%%%%%%%%%%%%%%%%%%%%%%%%%%%%%%%%%%%%%%%%%
param\`etres, et ensuite, 
nous supposerons qu'il
en a cinq\footnote{\,
Il s'agit des deux
groupes réduits~$(24)_2$ et $(26)$. 
}. %%%%%%%%%%%%%%%%%%%%%%%%%%%%%%%%%%%%%%%%%%%%%%%%%%%%%%%%%%%

\bigskip

\begin{center}
{\large
A)
D\'etermination de tous les groupes 
\\
dont le groupe r\'eduit poss\`ede six param\`etres
}
\end{center} 

\bigskip

Comme nous venons de le voir, si le groupe r\'eduit: $\overline{ X}_1
f \dots \overline{ X}_6 f$ a six param\`etres, il peut \^etre
rapport\'e \`a l'une des quatre formes suivantes
\footnote{\,
Pour étudier le groupe~\thetag{\bf II}, qui correspond exactement à
$(25)_1$, on remplace les deux générateurs infinitésimaux $xp$ et $yq$
par $xp - yq$ et $xp + yq$. Ces quatre générateurs du groupe réduit
$\overline{ X}_1, \dots, \overline{ X}_6$, lorsqu'il possède six
paramètres, sont réécrits dans un autre ordre, adaptés à l'avance aux
calculs de réduction et de normalisation du groupe $X_1, \dots, X_6$
qui vont suivre.
}: %%%%%%%%%%%%%%%%%%%%%%%%%%%%%%%%%%%%%%%%%%%%%%%%%%%%%%%%%%%

\[\label{413}
\aligned
{\bf (I)}\ \ \
&
p,\ \
q,\ \
xq,\ \
xp+yq,\ \
xp-yq,\ \
yp
\\
{\bf (II)}\ \ \
&
p,\ \
q,\ \
xq,\ \
xp+yq,\ \
xp-yq,\ \
x^2p+xyp
\\
{\bf (III)}\ \ \
&
p,\ \
q,\ \
xq,\ \
xp+yq,\ \
x^2q,\ \
x^2p+2xyq
\\
{\bf (IV)}\ \ \
&
p,\ \
q,\ \
xp,\ \
yq,\ \
x^2p,\ \
y^2q.
\endaligned
\]
Nous devons donc examiner ces quatre cas l'un apr\`es
l'autre.

\bigskip

\centerline{\sf
Premier cas}\label{414}

\medskip

En premier lieu, nous recherchons tous les groupes \`a six
param\`etres dans l'espace des $x, y, z$ qui satisfont nos exigences
formul\'ees \`a la page~\pageref{408} et dont les groupes 
r\'eduits sont de la
forme~\thetag{ I}. Chacun de ces groupes est de la forme:
\def\theequation{27}\begin{equation}
\left\{
\aligned
p+\varphi_1r,
\ \ \ \ \
&
q+\varphi_2r,
\ \ \ \ \
xq+\varphi_3r,
\ \ \ \ \
xp+yq+\varphi_4r
\\
&
xp-yq+\varphi_5r,
\ \ \ \ \
yp+\varphi_6r,
\endaligned\right.
\end{equation}
o\`u l'on entend
par $\varphi_1 \dots \varphi_6$ des fonctions de $x, y, z$.
Ainsi, nous devons maintenant d\'eterminer les fonctions $\varphi_1
\dots \varphi_6$ de la mani\`ere la plus g\'en\'erale possible pour
que les transformations infinit\'esimales~\thetag{ 27} constituent un
groupe \`a six param\`etres qui satisfont nos exigences formul\'ees
page~\pageref{408}, et nous pouvons encore choisir la variable $z$
d'une mani\`ere telle que la forme des fonctions $\varphi_1 \dots
\varphi_6$ soit la plus simple possible.

On peut supposer depuis le d\'ebut qu'\`a la place de $z$, on a
introduit une solution de l'\'equation diff\'erentielle:
\[
\frac{\partial f}{\partial x}
+
\varphi_1(x,y,z)\,
\frac{\partial f}{\partial z}
=
0
\]
qui n'est pas ind\'ependante de $z$; gr\^ace \`a ce choix, la
transformation infinit\'esimale $p + \varphi_1 r$ prend la forme
simple: $p$, tandis que la forme des autres transformations n'est pas
essentiellement modifi\'ee\footnote{\, En fait, tout changement de
coordonnées locales de la forme $\overline{ x} = x$, $\overline{ y} =
y$ et $\overline{ z} = \overline{ z} ( x, y, z)$ où la fonction
$\overline{ z}$, quelconque, est seulement supposée satisfaire la
condition $\overline{ z}_z \neq 0$ assurant qu'il s'agit d'un
difféomorphisme local, a la propriété de laisser inchangés les parties
en $p$ et en $q$ de chacune des six transformations infinitésimales du
groupe. Par exemple, $xp + yq + \varphi_4 r$ se transforme en
$\overline{ x} \overline{ p} +
\overline{ y} \overline{ q} +
\big( \varphi_4 \overline{ z}_z + x 
\overline{ z}_x + y \overline{ z}_y \big) \overline{ r}$.
Ce phénomène est bien sûr général et sera utilisé
dans d'autres contextes.
}. %%%%%%%%%%%%%%%%%%%%%%%%%%%%%%%%%%%%%%%%%%%%%%%%%%%%%%%%%%%
Nous pouvons donc supposer d\'esormais que
$\varphi_1$ est simplement nul.

\`A pr\'esent, il s'ensuit\footnote{\,
%%%%%%%%%%%%%%%%%%%%%%%-------DEBUT--------%%%%%%%%%%%%%%%%%%%%%%%%%%%
La plupart du temps, le <<\,crochet\,>> n'est pas nommé par
Engel et Lie et il est toujours noté au moyen de parenthèses.  
}: %%%%%%%%%%%%%%%%%%%%%%%%-----FIN-----%%%%%%%%%%%%%%%%%%%%%%%%%%%%%%%
\[
\big[
p,\
q+\varphi_2r
\big]
=
\frac{\partial\varphi_2}{\partial x}\,
r.
\]
Mais comme notre groupe ne renferme pas de transformation
infinit\'esimale de la forme\footnote{\,
---\,\,puisque le groupe réduit $\overline{ X}_1, 
\dots, \overline{ X}_6$ est de dimension six, ou
parce qu'on s'en convainc
aisément par un examen des six générateurs~\thetag{ 27}\,\,---
}: %%%%%%%%%%%%%%%%%%%%%%%%%%%%%%%%%%%%%%%%%%%%%%%%%%%%%%%%%%%
$\psi ( x, y, z)\, r$,
le coefficient $\frac{ \partial \varphi_2 }{ \partial x}$ en facteur
devant $r$ doit s'annuler identiquement, c'est-\`a-dire que
$\varphi_2$ doit \^etre libre de
$x$. Si de plus nous introduisons: $z_1 = \omega ( y, z)$
comme nouvelle variable, $p$ conserve sa forme\footnote{\,
---\,\,tandis que les quatre transformations infinitésimales
restantes conservent elles aussi essentiellement leur forme\,\,---
} %%%%%%%%%%%%%%%%%%%%%%%%%%%%%%%%%%%%%%%%%%%%%%%%%%%%%%%%%%%
tandis que la seconde
transformation infinit\'esimale devient:
\[
q
+
\varphi_2(y,z)\,r
=
q
+
\Big(
\frac{\partial\omega}{\partial y}
+
\varphi_2\,
\frac{\partial\omega}{\partial z}
\Big)\cdot
r_1.
\]
Si nous
choisissons $\omega$ de telle sorte que l'\'equation: $\omega_y
+ \varphi_2 \, \omega_z = 0$ est satisfaite, alors les
deux premi\`eres transformations infinit\'esimales de notre groupe
re\c coivent la forme:
\[
p,\ \
q.
\]

Afin de d\'eterminer $\varphi_3$, calculons les crochets:
\[
\aligned
\big[
p,\
xq+\varphi_3r
\big]
&
=
q
+
\frac{\partial\varphi_3}{\partial x}\,r
\\
\big[
q,\
xq+\varphi_3r
\big]
&
=
\ \ \ \ \ \ \,
\frac{\partial\varphi_3}{\partial y}\,r,
\endaligned
\]
qui montrent que $\varphi_3$ d\'epend seulement de $z$; mais par
ailleurs, $\varphi_3$ ne peut pas non plus
s'annuler, parce que sinon, notre groupe
comprendrait deux transformations infinit\'esimales: $q$ et $xq$,
dont les courbes int\'egrales co\"{\i}ncident, or ce cas est exclu,
d'apr\`es la page~\pageref{408}. Nous pouvons donc introduire comme
nouvelle coordonn\'ee $z$ une fonction de $z$ qui rend $\varphi_3 =
1$, et comme cela n'a pas d'influence sur $p$ et sur $q$, il en d\'ecoule
que les trois premi\`eres transformations infinit\'esimales de notre
groupe se pr\'esentent maintenant sous la forme:
\[
p,\ \
q,\ \
xq+r.
\]
En calculant les crochets de: $xp + yq + \varphi_4 r$ avec $p$ et
$q$, on v\'erifie, exactement comme pour $\varphi_3$, que
$\varphi_4$ d\'epend seulement de $z$; mais comme par ailleurs
on obtient aussi\footnote{\, Il sera maintenant fréquemment
sous-entendu que les crochets $\big[ X_j, \, X_k \big]$ calculés par
la suite doivent, d'après l'hypothèse que $X_1, \dots, X_6$ forment
une algèbre de Lie, s'exprimer comme certaines combinaisons linéaires
$\lambda_1 X_1 + \cdots +
\lambda_6 X_6$
à coefficients constants de $X_1, \dots, X_6$ eux-mêmes. Dans la
plupart des cas, en inspectant la partie en $p$ et en $q$ de chaque
$X_k$, on voit d'un seul coup d'{\oe}il les seuls constantes
$\lambda_k$ possibles qui peuvent apparaître, sans avoir à résoudre un
système linéaire.  Une fois que ces constantes sont déterminées, on en
déduit une condition intéressante sur les fonctions inconnues
$\varphi_1, \dots, \varphi_6$.  Par exemple, lorsque le crochet
calculé est de la forme $\psi \, r$, nécessairement toutes les
constantes $\lambda_k$ sont nulles, donc $\psi \equiv 0$,
et cela donne une équation différentielle sur 
les inconnues $\varphi_k$. Engel et Lie organisent les
calculs de telle sorte que la détermination et la
normalisation  
de ces fonction s'effectuent de manière progressive et graduelle. 
}: %%%%%%%%%%%%%%%%%%%%%%%%%%%%%%%%%%%%%%%%%%%%%%%%%%%%%%%%%%%
\[
\big[
xq+r,\ \,
xp+yq+\varphi_4(z)\,r
\big]
=
\varphi_4'(z)\,r,
\]
il s'ensuit que $\varphi_4 ' ( z) = 0$, d'o\`u $\varphi_4 ( z) = c$,
de telle sorte que nous obtenons la transformation infinit\'esimale:
\[
xp+yq+cr.
\]

Exactement comme pour $\varphi_3$ et $\varphi_4$, on v\'erifie
aussi
que $\varphi_5$ ne d\'epend que de $z$; de plus, nous obtenons:
\[
\big[
xq+r,\ \,
xp-yq
+
\varphi_5(z)\,r
\big]
=
-2xq
+
\varphi_5'(z)\,r,
\]
et donc\footnote{\,
Ce crochet doit être combinaison linéaire
des six transformations infinitésimales~\thetag{ 27}, dont
les quatre premières: $p$, $q$, $xq + r$ et $xp + yq + cr$
ont déjà été normalisées. Puisque le
membre de droite du crochet calculé contient seulement le
terme $-2xq$ en $p$ et en $q$, la combinaison
linéaire en question ne peut être que: 
$-2 \, ( xq+ r)$, d'où $\varphi_5'(z) \, r = -2 \, r$. 
}: %%%%%%%%%%%%%%%%%%%%%%%%%%%%%%%%%%%%%%%%%%%%%%%%%%%%%%%%%%%
$\varphi_5 ' ( z) = - 2$, c'est-\`a-dire: $\varphi_5 ( z) =
- 2 z + \text{ \rm const.}$ Mais en introduisant une nouvelle
coordonn\'ee $z$, on peut annuler la constante d'int\'egration, sans
modifier la forme des pr\'ec\'edentes\footnote{\,
Les quatre premières transformations infinitésimales
déjà normalisées sont de la forme $\xi \, p + \zeta\, q$, sans
terme en $\zeta \, r$, et elles
restent invariables par tout changement de
coordonnées de la forme 
$\overline{ x} = x$, $\overline{ y} = y$, 
$\overline{ z} = \overline{ z} ( x, y, z)$. 
} %%%%%%%%%%%%%%%%%%%%%%%%%%%%%%%%%%%%%%%%%%%%%%%%%%%%%%%%%%%
transformations
infinit\'esimales. Ensuite, en calculant le crochet:
\[
\big[
xp+yq+cr,\ \,
xp-yq-2zr
\big]
=
-2cr
\]
on obtient que la constante $c$ doit aussi s'annuler.

Finalement, pour d\'eterminer encore $\varphi_6$\,\,---\,\,qui ne d\'epend
naturellement que de $z$\,\,---, formons les crochets:
\[
\aligned
\big[
xq+r,\ \,
yp+\varphi_6(z)\,r
\big]
&
=
xp-yq+\varphi_6'(z)\,r
\\
\big[
xp-yq-2zr,\ \,
yp+\varphi_6(z)\,r
\big]
&
=
-2yp
-
2\{
z\varphi_6'(z)-\varphi_6(z)\}\,r,
\endaligned
\]
d'o\`u il d\'ecoule:
\[
\varphi_6'(z)
=
-2z,
\ \ \ \ \
z\,\varphi_6'(z)
=
2\,\varphi_6(z)
=
-2z^2.
\]
Par cons\'equent, nous trouvons que les transformations
infinit\'esimales de notre groupe peuvent recevoir la forme simple:
\def\theequation{28}\begin{equation}\label{415}
\left\{
\aligned
&
p,\ \
q,\ \
xq+r,\ \
yq+zr,\ \
xp-zr
\\
&
yp-z^2r.
\endaligned\right.
\end{equation}
Accessoirement, on peut remarquer que ce groupe appara\^{\i}t aussi
par
prolongement\footnote{\,
Un difféomorphisme local $(x, y) \mapsto (\overline{ x},
\overline{ y}) = 
\big( \overline{ x} ( x, y), \, \overline{ y} ( x, y) \big)$ du 
plan qui est proche de l'identité transforme un graphe $\{ y = y ( x)
\}$ en un autre graphe $\{ \overline{ y} =
\overline{ y} ( 
\overline{ x}) \big\}$
et transforme aussi les dérivées $y_x$ du graphe en:
\[
\overline{y}_{\overline{x}}
:=
\frac{d\overline{y}}{d\overline{x}}
= 
\frac{dx\cdot\partial\overline{y}(x,y(x))/\partial x}{
dx\cdot\partial\overline{x}(x,y(x))/\partial x}
=
\frac{\overline{x}_x+y_x\overline{y}_y}{\overline{x}_x+y_x
\overline{x}_y}.
\]
En appliquant cette
formule au groupe local à 
un paramètre: 
\[
(x,y)
\longmapsto 
\exp(tX)(x,y) 
=: 
\big(
\overline{x}(x,y,t),\,\overline{y}(x,y,t)
\big)
\]
engendré par une transformation infinitésimale
quelconque $X = \xi ( x, y) \, p
+ \zeta ( x, y) \, q$,
et en la différentiant par rapport à $t$ en $t = 0$, on obtient par le
calcul (\cite{ol1986, bk1989, ol1995, merk2009}):
\[
{\textstyle{\frac{d}{dt}}}
\overline{y}_{\overline{x}}
\big\vert_{t=0}
=
\eta_x
+
\left[
\eta_y
-
\xi_x
\right]
y_x
+
\left[
-
\xi_y
\right]
(y_x)^2.  
\]
Autrement dit, si l'on introduit le
{\sl prolongement $X^{ (1)}$ de $X$ aux jets d'ordre $1$}
qui est la transformation infinitésimale
en les trois variables $(x, y, y_x)$ définie par:
\[
X^{(1)}
:=
\xi(x,y)\,
{\textstyle{\frac{\partial}{\partial x}}}
+
\eta(x,y)\,
{\textstyle{\frac{\partial}{\partial y}}}
+
\big(
\eta_x
+
[\eta_y
-
\xi_x]
y_x
+
[-\xi_y]
(y_x)^2\big)
{\textstyle{\frac{\partial}{\partial y_x}}},
\]
alors on réobtient: 
\[
\exp\big(tX^{(1)}\big)(x,y,y_x)
=
\big(\overline{x}(x,y,t),\,\overline{y}(x,y,t),\,
\overline{y}_{\overline{x}}(x,y,y_x,t)\big)
\]
Les six transformations infinitésimales~\thetag{ 28} sont donc bien
des prolongements, à l'espace $(x, y, z \equiv y_x)$ des jets d'ordre
$1$, des six transformations $p$, $q$, $xq$, $yq$, $xp$ et $yp$.
} %%%%%%%%%%%%%%%%%%%%%%%%%%%%%%%%%%%%%%%%%%%%%%%%%%%%%%%%%%%
\deutsch{durch Hinzunahme} 
(\cf p.~165) du groupe lin\'eaire g\'en\'eral du
plan:
\[
p,\ \
q,\ \
xq,\ \
yq,\ \
xp,\ \
yp
\]
au moyen du quotient diff\'erentiel:
\[
\frac{dy}{dx}
=
z.
\]

Il reste maintenant encore \`a examiner si notre groupe~\thetag{ 28}
satisfait aux exigences de la page~\pageref{399}, donc avant tout si
deux points de l'espace des $x, y, z$ du groupe ont un et un seul
invariant relativement au groupe.

Notre groupe renferme trois transformations infinit\'esimales qui sont
d'ordre z\'ero par rapport \`a $x, y, z$ et dont on ne peut tirer, par
combinaison lin\'eaire, aucune transformation d'ordre un, ou d'un
ordre plus \'elev\'e, \`a savoir:
\[
p,\ \
q,\ \
xq+r.
\]
Nous en d\'eduisons que l'origine des coordonn\'ees: $x = y = z = 0$
est un point en position g\'en\'erale\footnote{\,
Dès que le groupe est (localement) transifif, 
les sous-algèbres d'isotropie de tout points sont
isomorphes deux à deux et tout
point est en position générale. 
}. %%%%%%%%%%%%%%%%%%%%%%%%%%%%%%%%%%%%%%%%%%%%%%%%%%%%%%%%%%%
Ce point reste au repos par les
trois transformations infinit\'esimales:
\[
yp-z^2r,
\ \ \ \ \
xp-yq-2zr,
\ \ \ \ \
xp+yq
\]
de notre groupe. Par cons\'equent, d'apr\`es la Proposition~2,
p.~\pageref{Satz-2}, l'existence d'un invariant pour les paires de
points d\'epend seulement du fait que le d\'eterminant:
\[
\left\vert
\begin{array}{ccc}
y & 0 & - z^2
\\
x & - y & - 2z
\\
x & y & 0
\end{array}
\right\vert
=
2y^2z
-
2xyz^2
\]
s'annule, ou ne s'annule pas, identiquement. Mais comme ce
d\'eterminant ne s'annule pas identiquement, deux points de $R_3$
n'ont s\^urement aucun invariant relativement au groupe~\thetag{ 28}.
Ainsi, le groupe~\thetag{ 28} ne fait pas partie des groupes que nous
recherchons \deutsch{ist also für uns unbrauchbar}.

\bigskip

\centerline{\sf
Deuxi\`eme cas}

\medskip

Le groupe r\'eduit: $\overline{ X}_1 f \dots \overline{ X}_6 f$ a
maintenant la forme~\thetag{ {\bf II}} explicit\'ee
page~\pageref{413}, par suite de quoi le groupe: $X_1 f \dots X_6 f$
est de la forme:
\def\theequation{29}\begin{equation}
\left\{
\aligned
&
p+\varphi_1r,\ \
q+\varphi_2r,\ \
xq+\varphi_3r,\ \
xp+yq+\varphi_4r
\\
&
xp-yq+\varphi_5r,
\ \
x^2p+xyq+\varphi_6r,
\endaligned\right.
\end{equation}
o\`u $\varphi_1 \dots \varphi_6$ sont des fonctions de $x, y, z$.

Les calculs effectu\'es pour le premier cas montrent directement que
les cinq premi\`eres transformations infinit\'esimales~\thetag{ 29}
peuvent \^etre mises sous la forme:
\[
p,\ \
q,\ \
xq+r,\ \
xp+yq,\ \
xp-yq-2zr.
\]
D'un autre c\^ot\'e, comme les crochets:
\[
\aligned
\big[
p,\ \
x^2p+xyq+\varphi_6r
\big]
&
=
2xp+yq
+
\frac{\partial\varphi_6}{\partial x}\,r
\\
\big[
q,\ \
x^2p+xyq+\varphi_6r
\big]
&
=
xq
+
\frac{\partial\varphi_6}{\partial y}\,r
\\
\big[
xp+yq,\ \
x^2p+xyq+\varphi_6r
\big]
&
=
x^2p+xyq+
\Big\{
x\,\frac{\partial\varphi_6}{\partial x}
+
y\,\frac{\partial\varphi_6}{\partial y}
\Big\}\,r
\endaligned
\]
donnent imm\'ediatement\footnote{\,
Rappelons que les membres de droite
de ces trois crochets doivent être combinaisons
linéaires à coefficients constants des six
transformations infinitésimales~\thetag{ 29}\,\,---\,\,dont 
les cinq premières sont déjà normalisées\,\,---,
}: %%%%%%%%%%%%%%%%%%%%%%%%%%%%%%%%%%%%%%%%%%%%%%%%%%%%%%%%%%%
\[
\aligned
\frac{\partial\varphi_6}{\partial x}
=
-z,
\ \ \ \ \ \ \
\frac{\partial\varphi_6}{\partial y}
=
1,
\\
x\,\frac{\partial\varphi_6}{\partial x}
+
y\,\frac{\partial\varphi_6}{\partial y}
=
\varphi_6
=
y-xz,
\endaligned
\]
nous obtenons les six transformations infinit\'esimales:
\def\theequation{30}\begin{equation}
\left\{
\aligned
&
p,\ \
q,\ \
xq+r,\ \
xp+yq,\ \
xp-yq-2zr
\\
&
x^2p+xyq
+
(y-zx)\,r.
\endaligned\right.
\end{equation}
Celles-ci constituent manifestement un groupe \`a six 
param\`etres\footnote{\,
En effet, on peut se convaincre sans aucun calcul de la stabilité par
crochets grâce à la propriété générale (\cite{ol1986, bk1989, ol1995,
merk2009}) que les crochets entre prolongements $\big[ X_j^{ (1)}, \,
X_k^{ (1)} \big] = \sum_{ l=1}^6 \, c_{ jkl} \, X_l^{ (1)}$ se
comportent exactement de la même manière que les crochets $\big[ X_k,
\, X_l \big] = \sum_{ l=1}^6 \, c_{ jkl}\, X_l$ entre les
transformations infinitésimales initiales.
} %%%%%%%%%%%%%%%%%%%%%%%%%%%%%%%%%%%%%%%%%%%%%%%%%%%%%%%%%%%
qui provient aussi par extension
\deutsch{durch Erweiterung} du groupe projectif \`a six
param\`etres du plan:
\[
p,\ \
q,\ \
xq,\ \
xp+yq,\ \
xp-yq,\ \
x^2p+xyq,
\]
de la m\^eme mani\`ere que le groupe~\thetag{ 28} provient par
extension du groupe lin\'eaire g\'en\'eral du plan.

Cependant, le groupe~\thetag{ 30} ne fait pas non plus partie des
groupes que nous recherchons. En effet, l'origine: $x = y = z = 0$
qui est \`a nouveau un point en position g\'en\'erale, reste invariante
par les trois transformations infinit\'esimales:
\[
xp+yq,
\ \
xp-yq-2zr,\ \
x^2p+xyq+(y-xz)\,r
\]
du groupe~\thetag{ 30}, mais le d\'eterminant correspondant:
\[
\left\vert
\begin{array}{ccc}
x & y & 0
\\
x & -y & -2z
\\
x^2 & xy & y-xz
\end{array}
\right\vert
=
-2xy(y-xz)
\]
ne s'annule pas identiquement, donc deux points de $R_3$ n'ont \`a nouveau
aucun invariant relativement au groupe~\thetag{ 30}.

\bigskip

\centerline{\sf
Troisi\`eme cas}

\medskip

Ici, le groupe r\'eduit: $\overline{ X}_1 f \dots \overline{ X}_6
f$ a la forme~\thetag{ {\bf III}}, explicit\'ee page~\pageref{413},
d'o\`u le groupe: $X_1 f \dots X_6 f$ est de la forme:
\def\theequation{31}\begin{equation}
\left\{
\aligned
&
p+\varphi_1r,\ \
q+\varphi_2r,\ \
xq+\varphi_3r,\ \
xp+yq+\varphi_4r
\\
&
x^2q+\varphi_5r,\ \
x^2p+2xyq+\varphi_6r,
\endaligned\right.
\end{equation}
o\`u l'on entend par $\varphi_1 \dots \varphi_6$ des fonctions de
$x, y, z$.

Exactement comme pour le premier cas, on v\'erifie que les quatre
premi\`eres transformations infinit\'esimales~\thetag{ 31} peuvent
\^etre mises sous la forme:
\[
p,\ \
q,\ \
xq+r,\ \
xp+yq+cr.
\]

Afin de d\'eterminer $\varphi_5$, formons les crochets:
\[
\aligned
\big[
p,\ \,
x^2q+\varphi_5\,r
\big]
&
=
2xq+
\frac{\partial\varphi_5}{\partial x}\,r
\\
\big[
q,\ \,
x^2q+\varphi_5\,r
\big]
&
=
\ \ \ \ \ \ \
\frac{\partial\varphi_5}{\partial y}\,r
\\
\big[
xq+r,\ \,
x^2q+\varphi_5\,r
\big]
&
=
\Big(
x\,\frac{\partial\varphi_5}{\partial y}
+
\frac{\partial\varphi_5}{\partial z}
\Big)\,r
\\
\big[
xp+yq+cr,\ \,
x^2q+\varphi_5\,r
\big]
&
=
x^2q+
\Big\{
x\,\frac{\partial\varphi_5}{\partial x}
+
y\,\frac{\partial\varphi_5}{\partial y}
+
c\,\frac{\partial\varphi_5}{\partial z}
\Big\}\,r,
\endaligned
\]
\`a partir desquels nous obtenons les \'equations:
\[
\aligned
\frac{\partial\varphi_5}{\partial x}
=2,
\ \ \ \ \ \ \
\frac{\partial\varphi_5}{\partial y}
=
\frac{\partial\varphi_5}{\partial z}
=
0,
\\
x\,\frac{\partial\varphi_5}{\partial x}
+
y\,\frac{\partial\varphi_5}{\partial y}
+
c\,\frac{\partial\varphi_5}{\partial z}
=
\varphi_5
=
2x.
\endaligned
\]
De la m\^eme mani\`ere, en calculant les crochets:
\[
\aligned
\big[
p,\ \,
x^2p+2xyq+\varphi_6r
\big]
&
=
2xp+2yq+
\frac{\partial\varphi_6}{\partial x}\,r
\\
\big[
q,\ \,
x^2p+2xyq+\varphi_6r
\big]
&
=
2xq+
\frac{\partial\varphi_6}{\partial y}\,r
\\
\big[
xq+r,\ \,
x^2p+2xyq+\varphi_6r
\big]
&
=
x^2q+
\Big(
x\,\frac{\partial\varphi_6}{\partial y}
+
\frac{\partial\varphi_6}{\partial z}
\Big)\,r
\\
\big[
xp+yq+cr,\ \,
x^2p+2xyq+\varphi_6r
\big]
&
=
x^2p+2xyq+
\\
&
+
\Big\{
x\,\frac{\partial\varphi_6}{\partial x}
+
y\,\frac{\partial\varphi_6}{\partial y}
+
c\,\frac{\partial\varphi_6}{\partial z}
\Big\}\,r,
\endaligned
\]
nous pouvons d\'eterminer la derni\`ere fonction 
inconnue $\varphi_6$ comme suit:
\[
\frac{\partial\varphi_6}{\partial x}
=
2c,
\ \ \ \ \
\frac{\partial\varphi_6}{\partial y}
=
2,
\ \ \ \ \
\frac{\partial\varphi_6}{\partial z}
=
0,
\]
\[
x\,\frac{\partial\varphi_6}{\partial x}
+
y\,\frac{\partial\varphi_6}{\partial y}
+
c\,\frac{\partial\varphi_6}{\partial z}
=
\varphi_6
=
2\,(y+cx).
\]
Par cons\'equent, notre groupe~\thetag{ 31} re\c coit la forme
\footnote{\,
Le fait que le groupe obtenu soit encadré laisse entendre à l'avance
qu'il sera retenu dans la liste de tous les groupes 
qui satisfont les conditions formulées à la page~\pageref{408}.
Mais la démonstration et les examens ne sont pas terminés.  
}: %%%%%%%%%%%%%%%%%%%%%%%%%%%%%%%%%%%%%%%%%%%%%%%%%%%%%%%%%%%
\def\theequation{32}\begin{equation}\label{418}
\boxed{
\aligned
p,\ \
q,\ \
xq+r,\ \
xp+yq+cr
\ \ \ \ \ \ \
\\
x^2q+2xr,\ \
x^2p+2xyq+2\,(y+cx)\,r
\endaligned
}\ 
.
\end{equation}

Ici aussi: $x = y = z = 0$ est un 
point en position g\'en\'erale. Ce point
reste invariant par les trois transformations
infinit\'esimales:
\def\theequation{33}\begin{equation}
\left\{
\aligned
&
xp+(y-cx)\,q,\ \
x^2q+2xr
\\
&
x^2p+2xyq+2\,(y+cx)\,r
\endaligned\right.
\end{equation}
et le d\'eterminant correspondant:
\def\theequation{34}\begin{equation}
\left\vert
\begin{array}{ccc}
x & y-cx & 0
\\
0 & x^2 & 2x
\\
x^2 & 2xy & 2(y+cx)
\end{array}
\right\vert
=
2x^3\,(y+cx)
-
4x^3y
+
2x^3(y-cx)
\end{equation}
s'annule identiquement, sans que tous ses sous-d\'eterminants
d'ordre deux ne s'annulent. Par cons\'equent ({\em
voir}\, Proposition~2,
p.~\pageref{Satz-2}), deux points de $R_3$: $x_1, y_1, z_1$ et $x_2,
y_2, z_2$ poss\`edent un et un seul invariant 
\label{419-0}
relativement au
groupe~\thetag{ 32}, et pour pr\'eciser, 
on trouve par le calcul\footnote{\,
Soit 
on forme et on résout le système~\thetag{ 2} avec
les six transformations infinitésimales
$X_k$, $k = 1, \dots, 6$ du groupe~\thetag{ 32},
soit on vérifie directement que l'expression 
de cet invariant est annulée identiquement par
$X_k^{ (1)} + X_k^{ (2)}$ pour
$k = 1, \dots, 6$.  
} %%%%%%%%%%%%%%%%%%%%%%%%%%%%%%%%%%%%%%%%%%%%%%%%%%%%%%%%%%%
l'expression
suivante pour cet invariant:
\def\theequation{35}\begin{equation}
z_1+z_2
-
c\,\log(x_2-x_1)^2
-
2\,\frac{y_2-y_1}{x_2-x_1}.
\end{equation}
Il s'agit alors de savoir aussi quelle valeur la constante $c$ peut
prendre; mais comme cette constante ne peut visiblement pas \^etre
supprim\'ee en introduisant de nouvelles variables $x, y, z$, elle
constitue alors un param\`etre essentiel, ce qui veut dire que le
groupe~\thetag{ 32} repr\'esente en fait $\infty^1$ types diff\'erents
de groupes.

Maintenant, il reste encore \`a d\'eterminer si le groupe~\thetag{ 32}
appartient r\'eellement aux groupes que nous recherchons, donc si les
invariants d'un nombre quelconque de points peuvent s'exprimer au
moyen des invariants de paires de points. Mais comme nous avons
d\'ej\`a vu que deux points poss\`edent un et un seul invariant, et
comme par ailleurs, la famille des pseudosph\`eres (\cf
p.~\pageref{402}):
\def\theequation{36}\begin{equation}\label{419}
z+z_0
-
c\,\log(x-x_0)^2
-
2\,
\frac{y-y_0}{x-x_0}
=
\text{\rm const.}
\end{equation}
n'est pas seulement constitu\'ee de $\infty^1$, mais 
m\^eme\footnote{\,
Il y a quatre paramètres: $z_0$, $x_0$, $y_0$ et
${\rm const.}$, mais
la constante absorbe $z_0$ par passage
du membre de gauche au membre de droite.
} %%%%%%%%%%%%%%%%%%%%%%%%%%%%%%%%%%%%%%%%%%%%%%%%%%%%%%%%%%%
de
$\infty^3$ surfaces diff\'erentes, il ne nous reste plus, d'apr\`es la
Proposition~4, p.~\pageref{Satz-4}, qu'\`a savoir s'il
existe une famille de $\infty^2$ courbes qui engendre la totalit\'e
des pseudosph\`eres.

Comme le point: $x = y = z = 0$ est un point en position
g\'en\'erale, alors, quand il est fix\'e, pour chaque famille de
$\infty^2$ courbes qui est invariante, la courbe qui passe par ce
point doit demeurer au repos. Or le point: $x = y = z = 0$ admet les
trois transformations infinit\'esimales~\thetag{ 33}, qui produisent
naturellement un sous-groupe \`a trois param\`etres de~\thetag{
32}. En annulant alors tous les sous-d\'eterminants d'ordre 
deux\footnote{\,
En restriction à une courbe quelconque passant
par l'origine qui est invariante
par le sous-groupe d'isotropie fixant
l'origine, le rang  
de l'espace vectoriel engendré par
les trois transformations infinitésimales~\thetag{ 33}, donc
de la matrice~\thetag{ 34}, 
doit être $\leqslant 1$.
} %%%%%%%%%%%%%%%%%%%%%%%%%%%%%%%%%%%%%%%%%%%%%%%%%%%%%%%%%%%
du
d\'eterminant~\thetag{ 34} relatif \`a ce sous-groupe, on trouve que
$x = y = 0$ est la seule courbe passant par le point: $x = y = z = 0$
qui reste invariante en m\^eme temps que ce point. Mais comme, par
l'action du groupe~\thetag{ 32}, la courbe: $x = y = 0$ prend toutes
les $\infty^2$ positions: $x = \text{\rm const.}$, $y = \text{\rm
const.}$, on a donc \'etabli que: $x = \text{\rm const.}$, $y =
\text{\rm const.}$ est l'unique famille de $\infty^2$ courbes qui est
invariante par le groupe~\thetag{ 32}. Enfin, puisque les
pseudosph\`eres~\thetag{ 36} ne sont \'evidemment pas compos\'ees des
courbes de la famille: $x = \text{\rm const.}$, $y = \text{\rm
const.}$, il appara\^{\i}t que tout nombre de points sup\'erieur \`a
deux n'a pas d'invariant essentiel relativement au groupe~\thetag{
32}.

Mentionnons encore que: $x = 0$ est l'unique surface passant par le
point: $x = y = z = 0$ qui reste invariante par le groupe~\thetag{
33}, et que par cons\'equent, mis \`a part la famille: $x = \text{\rm
const.}$, il n'existe pas d'autre famille de $\infty^1$ surfaces qui
admet le groupe~\thetag{ 32}.

\bigskip

\centerline{\sf
Quatri\`eme cas}

\medskip

Maintenant, nous cherchons tous les groupes \`a six param\`etres de
l'espace des $x, y, z$ qui satisfont aux exigences de la
page~\pageref{408} et dont le groupe r\'eduit poss\`ede la forme:
\def\theequation{IV}\begin{equation}
p,\ \
q,\ \
xp,\ \
yq,\ \
x^2p,\ \
y^2q,\ \
\end{equation}
Les transformations infinit\'esimales de ce groupe 
s'\'ecrivent:
\def\theequation{37}\begin{equation}
\left\{
\aligned
&
p+\varphi_1\,r,\ \
q+\varphi_2\,r,\ \
xp+\varphi_3\,r,\ \
yq+\varphi_4\,r,
\\
&
x^2p+\varphi_5\,r,\ \
y^2q+\varphi_6\,r,
\endaligned\right.
\end{equation}
o\`u $\varphi_1 \dots \varphi_6$ sont des
fonctions de $x, y, z$.

Comme dans le premier cas, on v\'erifie que 
[les fonctions]
$\varphi_1$ et $\varphi_2$
peuvent \^etre suppos\'ees nulles sans perte de g\'en\'eralit\'e. En
calculant le crochet de $xp + \varphi_3\, r$ avec $p$ et avec $q$, on
obtient que $\varphi_3$ ne d\'epend ni de $x$ ni de $y$; 
de plus, comme $\varphi_3$ ne peut \'evidemment pas \^etre 
nul\footnote{\,
Si la constante $c$ était nulle, le
groupe possèderait les deux transformations
$p$ et $xp$ dont les courbes
intégrales coïncideraient, ce qui serait exclut par 
la Proposition~1 p.~\pageref{Satz-1}. 
}, %%%%%%%%%%%%%%%%%%%%%%%%%%%%%%%%%%%%%%%%%%%%%%%%%%%%%%%%%%%
nous pouvons le rendre \'egal \`a $1$ en introduisant une nouvelle
variable $z$, et nous avons \`a pr\'esent les trois transformations
infinit\'esimales:
\[
p,\ \
q,\ \
xp+r.
\]

En calculant maintenant les crochets de ces transformations
simplifi\'ees avec: $yq + \varphi_4 \, r$, nous trouvons que
$\varphi_4$ ne d\'epend pas de $x, y, z$, donc est simplement une
constante; nous pouvons donc poser $\varphi_4 = c$, o\`u cependant la
constante $c$ ne doit pas s'annuler\footnote{\, Si la fonction
$\varphi_4$ était identiquement nulle, le groupe possèderait les deux
transformations $q$ et $yp$ dont les courbes intégrales coïncideraient.
}, %%%%%%%%%%%%%%%%%%%%%%%%%%%%%%%%%%%%%%%%%%%%%%%%%%%%%%%%%%%.

Afin de d\'eterminer aussi $\varphi_5$, formons les crochets:
\[
\aligned
\big[
p,\ \,
x^2p+\varphi_5\,r
\big]
&
=
2xp
+
\frac{\partial\varphi_5}{\partial x}\,r
\\
\big[
q,\ \,
x^2p
+
\varphi_5\,r
\big]
&
=
\ \ \ \ \ \ \ \ \ \
\frac{\partial\varphi_5}{\partial y}\,r
\\
\big[
yq+cr,\ \,
x^2p+\varphi_5\,r
\big]
&
=
\Big(
y\,\frac{\partial\varphi_5}{\partial y}
+
c\,\frac{\partial\varphi_5}{\partial z}
\Big)\,r
\\
\big[
xp+r,\ \,
x^2p+\varphi_5\,r
\big]
&
=
x^2p
+
\Big(
x\,\frac{\partial\varphi_5}{\partial x}
+
\frac{\partial\varphi_5}{\partial z}
\Big)\,r.
\endaligned
\]
Les trois premi\`eres de ces \'equations montrent que:
\[
\frac{\partial\varphi_5}{\partial x}
=
2,
\ \ \ \ \ \ \
\frac{\partial\varphi_5}{\partial y}
=
\frac{\partial\varphi_5}{\partial z}
=
0,
\]
et la derni\`ere donne alors: $\varphi_5 = 2x$.

Par un calcul compl\`etement similaire, on trouve: $\varphi_6 = 2\,
cy$. En d\'efinitive, nous parvenons donc au groupe:
\def\theequation{38}\begin{equation}
\boxed{
\aligned
p,\ \,
q,\ \,
xp+r,\ \,
&
yq+cr,\ \,
x^2p+2xr,\ \,
y^2q+2\,cyr
\\
&
\ \ \ \ \ \ \ \ 
{\scriptstyle{(c\,\neq\,0)}}
\endaligned
}\,.
\end{equation}
Ici, le param\`etre $c$ est manifestement essentiel et ne peut 
pas \^etre
supprim\'e.

La pr\'esence des trois transformations infinit\'esimales: $p$, $q$,
$xp + r$ montre que l'origine des coordonn\'ees est un
point en position g\'en\'erale. Ce point reste au repos
lorsqu'agissent les trois transformations infinit\'esimales:
\def\theequation{39}\begin{equation}
yq
-
cxp,\ \,
x^2p+2xr,\ \,
y^2q+2\,cyr
\ \ \ \ \ \ \
{\scriptstyle{(c\,\neq\,0)}}
\end{equation}
du groupe~\thetag{ 38}. Maintenant comme le d\'eterminant correspondant:
\def\theequation{40}\begin{equation}
\left\vert
\begin{array}{ccc}
-cx & y & 0
\\
x^2 & 0 & 2x
\\
0 & y^2 & 2cy
\end{array}
\right\vert
=
2\,cx^2y^2
-
2\,cx^2y^2
\end{equation}
s'annule identiquement, tandis que ses sous-d\'eterminants d'ordre
deux ne sont pas tous nuls, nous pouvons en d\'eduire que deux points
quelconques de l'espace des $x, y, z$ ont un et un seul invariant
relativement au groupe~\thetag{ 38}. En fait, on se persuade par un
calcul ais\'e que l'expression:
\def\theequation{41}\begin{equation}
z_1+z_2
-
\log(x_2-x_1)^2
-
c\cdot 
\log(y_2-y_1)^2
\end{equation}
constitue un invariant, et \`a vrai dire, l'unique invariant des deux
points: $x_1, y_1, z_1$ et $x_2, y_2, z_2$.

La famille des pseudosph\`eres de notre groupe~\thetag{ 38} est de la
forme:
\def\theequation{42}\begin{equation}
z+z_0
-
\log(x-x_0)^2
-
c\cdot 
\log(y-y_0)^2
=
\text{\rm const.}
\end{equation}
Ainsi, elle est constitu\'ee d'un nombre de surfaces qui est
sup\'erieur \`a $\infty^1$. Par ailleurs, on v\'erifie, comme dans le
cas pr\'ec\'edent, en annulant les sous-d\'eterminants d'ordre deux
de~\thetag{ 40}, que $x = y = 0$ est la seule courbe passant par le
point: $x = y = z = 0$ qui reste au repos en m\^eme temps que ce point
et donc aussi que la famille des $\infty^2$ courbes: $x = \text{\rm
const.}$, $y = \text{\rm const.}$ est la seule qui reste invariante
par le groupe~\thetag{ 38}. Mais comme la pseudosph\`ere~\thetag{ 42}
n'est manifestement pas constitu\'ee par les courbes de cette famille,
on en d\'eduit avec certitude, gr\^ace
\`a la Proposition~4, p.~\pageref{Satz-4}, qu'un nombre
de points sup\'erieur \`a deux ne poss\`ede aucun invariant essentiel
relativement au groupe~\thetag{ 38}.

On peut encore mentionner qu'en dehors des deux familles de surfaces:
$x = \text{\rm const.}$ et $y = \text{\rm const.}$, il n'existe pas
d'autre famille de $\infty^1$ surfaces qui soit invariante par le
groupe~\thetag{ 38}.

\bigskip

\begin{center}
{\large
B)
D\'etermination de tous les groupes 
\\
dont le groupe r\'eduit poss\`ede cinq param\`etres
}
\end{center} 

\bigskip

En cette circonstance, d'apr\`es les pages~\pageref{412a} sq., le
groupe r\'eduit: $\overline{ X}_1 f \dots \overline{ X}_6 f$ peut
\^etre rapport\'e aux deux formes suivantes:
\[\label{422}
\aligned
{\bf (V)}\ \ \
&
p,\ \
q,\ \
xq,\ \
xp-yq,\ \
yp
\\
{\bf (VI)}\ \ \
&
p,\ \
q,\ \
xq,\ \
2xp+yq,\ \
x^2p+xyq.
\endaligned
\]
Par cons\'equent, nous devons encore traiter maintenant ces deux cas.

\bigskip

\centerline{\sf
Cinqui\`eme cas}

\medskip

Si le groupe r\'eduit est de la forme~\thetag{ V}, le groupe: $X_1 f
\dots X_6 f$ peut \^etre mis sous la forme:
\def\theequation{43}\begin{equation}
\left\{
\aligned
&
\varphi\,r,\ \
p+\varphi_1\,r,\ \
q+\varphi_2\,r,\ \
xp-yq+\varphi_3\,r,\ \
\\
&
xq+\varphi_4\,r,\ \
yp+\varphi_5\,r,
\endaligned\right.
\end{equation}
o\`u l'on entend par $\varphi, \varphi_1, \dots, \varphi_5$ des
fonctions de $x, y, z$.

Ce cas a d\'ej\`a \'et\'e compl\`etement \'etudi\'e aux pages~157 sq.,
o\`u l'on a d\'emontr\'e\footnote{\, Par souci de complétude, voici
les arguments. En introduisant $\int\, \frac{ 1}{ \varphi}\, dz$ comme
nouveau $z$, on remplace $\varphi \, r$ par $r$ et on obtient le
groupe:
\[
r,\ \ \
p+\varphi_1\,r,\ \ \ \ \
q+\varphi_2\,r,\ \ \ \ \
xp-yq+\varphi_3\,r,\ \ \ \ \
xq+\varphi_4\,r,\ \ \ \ \
yp+\varphi_5\,r. 
\]
Le crochet avec $r$ des cinq derniers générateurs montre que chaque
$\varphi_k$ a la forme $\varphi_k = a_k z + \psi_k (x, y)$, $k = 1,
\dots, 5$.  Mais comme le crochet:
\[
\big[p+(a_1z+\psi_1)\,r,\,\,
xp-yq+(a_4z+\psi_4)\,r\big]
=
p
+
\chi(x,y)\,r,
\]
où $\chi$ est une fonction de $x$ et de $y$ qu'il est inutile de
calculer précisément, ne peut qu'être égal à $p + ( a_1 z + \psi_1)\,
r$, on obtient que $a_1 = 0$, et d'une manière analogue aussi, que
$a_2 = \cdots = a_5 = 0$.  En introduisant maintenant $z - \int \,
\psi_1 dx$ comme nouveau $z$, on redresse la transformation
infinitésimale $p + \psi_1 r$ en $p$.  Le crochet $\big[ p, \, q +
\psi_2 r \big] = \frac{ \partial \psi_2}{ \partial x}\, r$
montre alors que $\varphi_2 = C x + \chi ( y)$, et on peut annuler
ensuite la fonction $\chi ( y)$ en introduisant $z - \int \, \chi\,
dy$ comme nouveau $z$. Le groupe considéré est ainsi réduit à la
forme:
\[
r,\ \ \ \ \
p,\ \ \ \ \
q+Cx\,r,\ \ \ \ \
xp-yq+\psi_3\,r,\ \ \ \ \
xq+\psi_4\,r,\ \ \ \ \ 
yp+\psi_5\,r,
\]
où $\psi_3$, $\psi_4$, $\psi_5$ sont des fonctions de $(x, y)$.  Le
crochet $\big[ p, \, xp - yq + \psi_3 \, r \big] = p + \frac{ \partial
\psi_3 }{ \partial x}\, r$ montre que $\psi_3 = a x + \tau ( y)$. D'un
autre côté, on a:
\[
\big[
q+Cx\,r,\,\,
xp-yq+(ax+\tau)\,r
\big]
=
-q
+
\big(\tau'(y)-Cx)\,r,
\]
d'où $\tau' ( y)$ est constant $=b$
et $\tau ( y) = by$, si l'on supprime la constante
d'intégration grâce à la présence de la (première)
transformation infinitésimale $r$. 
En introduisant encore $z - ax + by$ comme
nouveau $z$, on obtient simplement
la transformation infinitésimale $xp - yq$, tandis
que $r$, $p$ et $q + Cxr$ ne changent pas essentiellement de
forme. 
De plus, on a: 
\[
\big[
p,\,\,xq+\psi_4\,r
\big]
=
q
+
{\textstyle{\frac{\partial\psi_4}{\partial x}}}\,r,
\ \ \ \ \
\big[
q+Cxr,\,\,xq+\psi_4\,r
\big]
=
{\textstyle{\frac{\partial\psi_4}{\partial y}}}\,r,
\]
d'où en faisant abstraction de la constante d'intégration
(superflue): 
$\psi_4 = \frac{ 1}{ 2}\, C x^2 + ky$. 
Mais le crochet: 
\[
\big[xp-yq,\,\,
xq+({\textstyle{\frac{1}{2}}}\,Cx^2+ky)\,r\big]
=
2xq+(Cx^2-ky)\,r
\]
montre que la constante $k$ s'annule.  Enfin, on montre d'une manière
analogue que $\psi_5 = \frac{ 1}{ 2}\, Cy^2$.  Si $C = 0$, on obtient
le premier groupe; si $C \neq 0$, en remplaçant $z$ par $\frac{ 1}{
C}\, z$, on obtient le second groupe~\thetag{ 44}.
} %%%%%%%%%%%%%%%%%%%%%%%%%%%%%%%%%%%%%%%%%%%%%%%%%%%%%%%%%%%
que le groupe~\thetag{ 43} peut \^etre, par
un choix appropri\'e de la variable $z$, rapport\'e soit \`a la
forme:
\[
r,\ \
p,\ \
q,\ \
xq,\ \
xp-yq,\ \
yp,
\]
soit \`a la forme:
\def\theequation{44}\begin{equation}
\left\{
\aligned
&
r,\ \
p,\ \
q+xr,\ \
xp-yq,\ \
\\
&
yp+{\textstyle{\frac{1}{2}}\,}
y^2r,\ \
xq+{\textstyle{\frac{1}{2}}\,}
x^2r.
\endaligned\right.
\end{equation}
Le premier de ces groupes n'entre pas en ligne de compte, puisque $q$
et $xq$ poss\`edent les m\^emes courbes int\'egrales. Le
groupe~\thetag{ 44} nous est d\'ej\`a connu, d'apr\`es le Tome~II
({\em voir} le Th\'eor\`eme~66, p.~421); de plus, nous savons qu'en
introduisant les expressions $x$, ${\scriptstyle{ \frac{ 1}{ 2}}} \,
y$, $z - {\scriptstyle{ \frac{ 1}{ 2}}} xy$ comme nouvelles variables
\`a la place
de $x$, $y$, $z$, ce groupe
se transforme\footnote{\,
On vérifie en effet par le calcul que~\thetag{ 44}
devient~\thetag{ 45}, lequel 
est effectivement un sous-groupe
du groupe projectif
engendré par les 8 générateurs infinitésimaux
$p$, $q$, $xp$, $yp$, $xq$, $yq$, 
$xxp + xy q$, $xyp + yyq$. 
} %%%%%%%%%%%%%%%%%%%%%%%%%%%%%%%%%%%%%%%%%%%%%%%%%%%%%%%%%%% 
en le groupe
projectif suivant ({\em loc. cit.}, pp.~445 sq.):
\def\theequation{45}\begin{equation}
\boxed{
p-yr,\ \
q+xr,\ \
r,\ \
xq,\ \
xp-yq,\ \
yp
}\,.
\end{equation}
Ainsi, pour notre recherche, nous souhaitons consid\'erer ce groupe
fondamentalement sous cette forme projective.

La pr\'esentation des trois transformations infinit\'esimales: $p -
yr$, $q + xr $, $r$ montre que $x = y = z = 0$ est un point en
position g\'en\'erale. Mais comme ce point reste au repos par
l'action des trois transformations infinit\'esimales:
\def\theequation{46}\begin{equation}
xq,\ \
xp-yq,\ \
yp
\end{equation}
du groupe~\thetag{ 45}, et comme le d\'eterminant correspondant:
\def\theequation{47}\begin{equation}
\left\vert
\begin{array}{ccc}
0 & x & 0
\\
x & -y & 0
\\
y & 0 & 0
\end{array}
\right\vert
\end{equation}
s'annule identiquement, sans que tous ses sous-d\'eterminants d'ordre
deux en fassent de m\^eme, il d\'ecoule
tr\`es certainement de la Proposition~2
p.~\pageref{Satz-2} que deux points: $x_1, y_1, z_1$ et $x_2, y_2,
z_2$ poss\`edent un et un seul invariant relativement au
groupe~\thetag{ 45}. Pour cet invariant, on trouve tr\`es facilement
l'expression: 
\def\theequation{48}\begin{equation}
z_2-z_1
+
x_1y_2-x_2y_1.
\end{equation}

La famille des pseudosph\`eres:
\def\theequation{49}\begin{equation}
z-z_0
+
x_0y-y_0x
=
\text{\rm const.}
\end{equation}
appartenant \`a ce groupe consiste en $\infty^3$ surfaces
diff\'erentes et pour pr\'eciser, 
elle consiste en tous les plans de l'espace des
$x, y, z$. Comme de
plus, d'apr\`es le Th\'eor\`eme~73 p.~445 du Tome~II, le
groupe~\thetag{ 44} ne laisse invariante aucune autre famille de
$\infty^2$ courbes que la famille: $x = \text{\rm const.}$, $y =
\text{\rm const.}$, il en va \'evidemment de m\^eme pour le
groupe~\thetag{ 45}, et comme les $\infty^3$ pseudosph\`eres~\thetag{
49} ne sont pas enti\`erement constitu\'ees des courbes de la
famille: $x = \text{\rm const.}$, $y = \text{\rm const.}$, 
on d\'eduit donc imm\'ediatement de la Proposition~4,
p.~\pageref{Satz-4}, qu'un nombre de points sup\'erieur \`a deux ne
poss\`ede aucun invariant essentiel relativement au
groupe~\thetag{ 44}.

On doit encore mentionner que le groupe~\thetag{ 45} ne laisse en
g\'en\'eral invariante aucune famille de $\infty^1$ surfaces ({\em
voir} Tome~II, {\em loc. cit.}).

\bigskip

\centerline{\sf
Sixi\`eme cas}

\medskip

Lorsque le groupe r\'eduit: $\overline{ X}_1 f \dots \overline{
X}_6 f$ a la forme~\thetag{VI}, p.~\pageref{422}, le groupe: $X_1 f
\dots X_6 f$ est n\'ecessairement de la forme:
\def\theequation{50}\begin{equation}
\left\{
\aligned
&
\varphi\,r,\ \
p+\varphi_1\,r,\ \
q+\varphi_2\,r,\ \
2xp+yq+\varphi_3\,r,
\\
&
xq+\varphi_4\,r,\ \
x^2p+xyq+\varphi_5\,r,
\endaligned\right.
\end{equation}
o\`u $\varphi, \varphi_1, \dots, \varphi_5$ sont des fonctions de $x,
y, z$.

Exactement comme aux pages~157 sq.\footnote{\,
Voir la note de rappel au début de l'étude du 
cinquième cas. 
}, %%%%%%%%%%%%%%%%%%%%%%%%%%%%%%%%%%%%%%%%%%%%%%%%%%%%%%%%%%%
on peut parvenir \`a
ce que $\varphi = 1$,
$\varphi_1 = 0$, $\varphi_2 = Cx$; alors
en m\^eme temps, les fonctions $\varphi_3$,
$\varphi_4$ et $\varphi_5$ sont certainement toutes ind\'ependantes de
$z$.

Pour la d\'etermination de $\varphi_3 (x, y)$, formons les crochets:
\[
\aligned
\big[
p,\ \,
2xp+yq+\varphi_3\,r
\big]
&
=
2\,p
+
\frac{\partial\varphi_3}{\partial x}\,r
\\
\big[
q+Cxr,\ \,
2xp+yq+\varphi_3\,r
\big]
&
=
q
+
\Big(
\frac{\partial\varphi_3}{\partial y}
-
2Cx
\Big)\,r,
\endaligned
\]
qui donnent:
\[
\frac{\partial\varphi_3}{\partial x}
=
D,
\ \ \ \ \ \ \ \ \
\frac{\partial\varphi_3}{\partial y}
=
3\,Cx+E,
\]
d'o\`u: $C = 0$ et $\varphi_3 = D\, x + E\, y + H$, o\`u la constante
$H$ peut \^etre simplement supprim\'ee\footnote{\, \,\,---\,\,en
soustrayant $H r$ à la transformation infinitésimale 
$2xp + yq + \varphi_3 \,
r$.
}, %%%%%%%%%%%%%%%%%%%%%%%%%%%%%%%%%%%%%%%%%%%%%%%%%%%%%%%%%%% 
tandis que les constantes $D$ et $E$ sont ramen\'ees \`a
z\'ero lorsqu'on introduit: $z - {\scriptstyle \frac{ 1}{ 2}}\, Dx - E
y$ comme nouvelle variable $z$. Nous avons
donc maintenant ramen\'e les quatres
premi\`eres transformations~\thetag{ 50} \`a la forme simplifi\'ee:
\[
p,\ \
q,\ \
r,\ \
2\,xp+yq.
\]

En calculant les crochets de $x q + \varphi_4 (x, y) \, r$ avec $p$ et
avec
$q$, on v\'erifie que $\varphi_4$ a la forme: $\varphi_4 = L x + My$,
o\`u la constante d'int\'egration (superflue) a d\'ej\`a \'et\'e
supprim\'ee\footnote{\, \,\,---\,\,\`a nouveau, en soustrayant un multiple
de la transformation $r$.
}. %%%%%%%%%%%%%%%%%%%%%%%%%%%%%%%%%%%%%%%%%%%%%%%%%%%%%%%%%%%  
En outre, on obtient:
\[
\big[
2xp+yq,\, \
xq+(Lx+My)\,r
\big]
=
xq
+
(2Lx+My)\,r,
\]
de telle sorte que 
$L$ doit \^etre nul, tandis que $M$ ne peut naturellement pas
s'annuler, sinon $q$ et $x q$ auraient les m\^emes courbes
int\'egrales. En introduisant ${\scriptstyle \frac{ 1}{ M}} \, z$
comme nouvelle variable $z$, on peut faire que la constante $M$ soit
simplement \'egale \`a $1$.

Enfin, formons les crochets:
\[
\aligned
\big[
p,\ \,
x^2p+xyq+\varphi_5\,r
\big]
&
=
2xp+
yq
+
\frac{\partial\varphi_5}{\partial x}\,r
\\
\big[
q,\ \,
x^2p+xyq+\varphi_5\,r
\big]
&
=
\ \ \ \ \ \ \ \ 
xq
+
\frac{\partial\varphi_5}{\partial y}\,r
\\
\big[
xq+yr,\ \,
x^2p+xyq+\varphi_5\,r
\big]
&
=
\Big(
x\,\frac{\partial\varphi_5}{\partial y}
-
xy
\Big)\,r
\\
\big[
2xp+yq,\ \,
x^2p+xyq+\varphi_5\,r
\big]
&
=
2(x^2p+xyq)
+
\Big(
2x\,\frac{\partial\varphi_5}{\partial x}
+
y\,\frac{\partial\varphi_5}{\partial y}
\Big)\,r,
\endaligned
\]
\`a partir desquels nous trouvons:
\[
\aligned
\frac{\partial\varphi_5}{\partial x}
=
c,
\ \ \ \ \ \ \
\frac{\partial\varphi_5}{\partial y}
=
y
\\
2x\,\frac{\partial\varphi_5}{\partial x}
+
y\,\frac{\partial\varphi_5}{\partial y}
=
2cx+y^2
=
2\,\varphi_5
+
K,
\endaligned
\]
o\`u cependant la constante $K$ peut \^etre simplement r\'eduite \`a
z\'ero. Ainsi, nous sommes parvenus au groupe:
\def\theequation{51}\begin{equation}
\left\{
\aligned
p,\ \
q,\ \
r,\ \
2xp+yq,\ \
xq+yr
\\
x^2p+xyq+({\textstyle{\frac{1}{2}}}y^2+cx)\,r.
\endaligned\right.
\end{equation}

Puisque les transformations infinit\'esimales $p$, $q$, $r$ sont
pr\'esentes, le point $x = y = z = 0$ est \`a nouveau un point en
position g\'en\'erale. Les transformations infinit\'esimales du
groupe~\thetag{ 51} qui laissent invariant ce point sont:
\[
2xp+yq,\ \ \
xq+yr,\ \ \
x^2p+xyq+({\textstyle{\frac{1}{2}}}y^2+cx)\,r,
\]
et le d\'eterminant correspondant: 
\[
\left\vert
\begin{array}{ccc}
2x & y & 0
\\
0 & x & y
\\
x^2 & xy & {\textstyle{\frac{1}{2}}}y^2+cx
\end{array}
\right\vert
\]
a la valeur: 
\[
x^2y^2
+
2cx^3
-
2x^2y^2
+
x^2y^2
=
2cx^3,
\]
qui ne s'annule donc identiquement que lorsque la constante $c$ a la
valeur z\'ero, tandis que les sous-d\'eterminants d'ordre deux ne
s'annulent pas tous, y compris pour $c = 0$. Par cons\'equent, parmi
tous les groupes de la forme~\thetag{ 51}, le groupe:
\def\theequation{52}\begin{equation}
\boxed{
\aligned
&
\ \
p,\ \
q,\ \
r,\ \
2xp+yq,\ \
xq+yr\ \
\\
&
\ \ \ \ \ \ \ \ \ \ \
x^2p+xyq+{\textstyle{\frac{1}{2}}}y^2r
\endaligned
}
\end{equation}
est le seul relativement auquel deux points poss\`edent un et un seul
invariant. Par le calcul, on trouve que cet invariant poss\`ede la forme:
\def\theequation{53}\begin{equation}
z_2-z_1
-
\frac{(y_2-y_1)^2}{2\,(x_2-x_1)}.
\end{equation}

La famille des pseudosph\`eres:
\def\theequation{54}\begin{equation}
z-z_0
-
\frac{(y-y_0)^2}{2\,(x-x_0)}
=
\text{\rm const.}
\end{equation}
qui appartient \`a notre groupe~\thetag{ 52} consiste manifestement en
$\infty^3$ surfaces diff\'erentes. Par ailleurs, on se convainc
facilement que les transformations infinit\'esimales du
groupe~\thetag{ 52} qui fixent le point en position g\'en\'erale: $x
= y = z = 0$, ne laissent au repos nulle autre courbe passant par ce
point que la droite: $x = y = 0$; par cons\'equent: $x = \text{\rm
const.}$, $y = \text{\rm const.}$ est la seule famille de $\infty^2$
courbes qui reste invariante par le groupe~\thetag{ 52} et la
Proposition~4, p.~\pageref{Satz-4} montre \`a nouveau qu'un nombre
de points sup\'erieur \`a deux n'a pas d'invariant essentiel
relativement au groupe~\thetag{ 52}.

Pour terminer, remarquons encore que le groupe~\thetag{ 52} ne laisse
invariante qu'une seule famille de $\infty^1$ surfaces, \`a savoir la
famille: $x = \text{\rm const.}$

\bigskip

{\em Ainsi, nous avons aussi trouv\'e tous les groupes imprimitifs de
$R_3$ qui satisfont l'exigence que nous avons formul\'ee
page~\pageref{399}.}\label{425}

\HEAD{Groupes de $R_3$ pour lesquels deux points ont un seul
invariant.}{Division\,\,V.\,\,\,Chapitre\,\,20.\,\,\,\S\,\,88.}

\sectiondritterV{\sf\S\,\,88.}
\label{S-88}
\setcounter{footnote}{0}

D'apr\`es les r\'esultats des deux pr\'ec\'edents paragraphes, le
probl\`eme \'enonc\'e dans l'ouverture du pr\'esent chapitre possède
la solution suivante\footnote{\, Le résultat vaut pour un groupe de
transformations holomorphes locales agissant sur $\C^3$, le cas des
groupes réels étant traité au \S\,\,89 qui suit.
}. %%%%%%%%%%%%%%%%%%%%%%%%%%%%%%%%%%%%%%%%%%%%%%%%%%%%%%%%%%%

\medskip
{\bf Th\'eor\`eme~36.}\label{Theorem-36} 
{\em Si un groupe continu fini de l'espace des $x, y, z$ est
constitu\'e de telle sorte que relativement
\`a son action, deux points ont un et un seul invariant tandis que
tous les invariants d'un nombre de points sup\'erieur \`a deux peuvent
s'exprimer au moyen des invariants des paires de points, alors le
groupe est transitif \`a six param\`etres et pour pr\'eciser, il est
semblable\footnote{\, Le langage contemporain appelle {\sl
équivalentes} deux structures géométriques qui ne diffèrent l'une de
l'autre que par un changement de coordonnées ponctuelles.
} %%%%%%%%%%%%%%%%%%%%%%%%%%%%%%%%%%%%%%%%%%%%%%%%%%%%%%%%%%%
\deutsch{ähnlich}, 
{\em via}\, une transformation ponctuelle de l'espace des $x, y, z$,
ou bien au groupe des mouvements euclidiens, ou bien au groupe
projectif \`a six param\`etres qui stabilise une surface du second
degr\'e non-d\'eg\'en\'er\'ee, ou bien encore \`a l'un des quatre
groupes suivants}~{\rm :}
\[
\label{426}
\text{\bf [1]}
\ \ \ \ \ \ \ \ \
\boxed{
\aligned
\ \
p,\ \
q,\ \
xp+r,\ \
&
yq+cr,\ \
x^2p+2xr,\ \
y^2q+2cyr\ \
\\
&
\ \ \ \ \ \ \ \ \
c\neq 0
\endaligned
}
\ \ \ \ \ \ \ \ \
\]
\[
\text{\bf [2]}
\ \ \ \ \ \ \ \ \ \ \
\boxed{
\aligned
\ \
p,\ \
q,\ \
x&q+r,\ \
x^2q+2xr,\ \
xp+yq+cr\ \
\\
&
x^2p+2xyq+2(y+cx)\,r
\endaligned
}
\ \ \ \ \ \ \ \ \ \ \
\]
\[
\text{\bf [3]}
\ \ \ \ \ \ \ \ \
\boxed{
\ \
p-yr,\ \
q+xr,\ \
r,\ \
xq,\ \
xp-yq,\ \
yp\ \
}
\ \ \ \ \ \ \ \ \
\]
\[
\text{\bf [4]}
\ \ \ \ \ \ \ \ \ \ \
\boxed{
\aligned
\ \
p,\ \
q,&\ \
r,\ \
xq+yr,\ \
2xp+yq\ \
\\
&
x^2p+xyq+
{\textstyle{\frac{1}{2}}}y^2r
\endaligned
}
\ \ \ \ \ \ \ \ \ \ \
\]
{\em Le param\`etre $c$ est ici essentiel dans les deux cas {\rm [1]} et
{\rm [2]} et ne peut pas \^etre supprim\'e.}

\medskip

Comme nous l'avons montr\'e dans les pr\'ec\'edents paragraphes,
chacun des quatre groupes [1] \dots \ [4] laisse invariante une unique
famille de $\infty^2$ courbes, \`a savoir la famille: $x = \text{\rm
const.}$, $y= \text{ \rm const.}$ Maintenant, nous voulons encore
ajouter que ces quatre groupes sont enti\`erement 
{\sl systatiques}\footnote{\, Soit $X_1, \dots, X_m$ un groupe fini
continu de transformations locales d'une variété à $n$ dimension munie
des coordonnées $(x_1, \dots, x_n)$.  Un point $(x_1^0, \dots, x_n^0)$
en position générale admet un certain nombre $m - q$ de
transformations infinitésimales $Y_1, \dots, Y_{ m-q}$ qui le laissent
au repos, {\em cf.}
les rappels de la note p.~\pageref{isotropie-general}.  Dans certaines
circonstances, il peut exister une sous-variété connexe non vide
$\Lambda$ (avec des singularités éventuelles) passant par ce point
dont {\em chaque point individuel} est laissé fixe par toutes les
transformations du sous-groupe engendré par $Y_1, \dots, Y_{ m-q}$,
qui n'est autre que le sous-groupe d'isotropie de $(x_1^0, \dots,
x_n^0)$.  Dans ce cas, le groupe est dit {\sl systatique}, 
concept introduit par Lie en 1884--85. D'après le
Chap.~24 du Tome~I, en reprenant les notations de la note
p.~\pageref{isotropie-general}, les équations du plus grand tel
ensemble $\Lambda$ sont: $\varphi_{ jk} ( x_1, \dots, x_n) =
\varphi_{ jk} ( x_1^0, \dots, x_n^0)$ pour
$j = 1, \dots, m-q$ et $k = 1, \dots, q$.  Si $\Lambda$ est l'ensemble
vide, le groupe est dit {\sl asystatique}.

}. %%%%%%%%%%%%%%%%%%%%%%%%%%%%%%%%%%%%%%%%%%%%%%%%%%%%%%%%%%%

En effet, les transformations infinit\'esimales des quatre groupes
sus-nomm\'es commutent toutes avec la transformation
infinit\'esimale $r$, ce qui montre, gr\^ace \`a la Proposition~2
p.~510 du Tome~I, que ces groupes sont systatiques; on peut aussi
s'en convaincre directement, puisque pour chacun des groupes [1]
\dots \ [4], les transformations infinit\'esimales qui fixent un point
$x_0, y_0, z_0$ en position g\'en\'erale laissent en m\^eme temps au
repos la droite: $x = x_0$, $y = y_0$.

Enfin, mentionnons encore que par l'action des groupes [1], [2], [3], les
$\infty^2$ \'el\'ements lin\'eaires qui passent par chaque point
fix\'e en
position g\'en\'erale se transforment de trois mani\`eres
diff\'erentes, mais seulement de deux mani\`eres par l'action du groupe
[4]; c'est pr\'ecis\'ement parce que le groupe [4] renferme, dans le
voisinage de chaque point: $x_0, y_0, z_0$ en position g\'en\'erale,
une transformation infinit\'esimale du second ordre en: $x-x_0$, $y -
y_0$, $z- z_0$.

On pourrait se demander si les groupes [1], [2], [4], qui ne sont pas
pr\'esent\'es sous forme projective, sont transformables en des
groupes projectifs par un choix appropri\'e de coordonn\'ees $x, y,
z$. Pour le groupe [4], la r\'eponse \`a cette question est tr\`es
certainement non, car s'il devait se transformer en un groupe
projectif, la transformation infinit\'esimale du second ordre qu'il
renfermerait dans le voisinage du point: $x = y = z = 0$ en position
g\'en\'erale devrait \^etre de la forme\footnote{\,
En effet, le groupe projectif contient toutes les transformations
infinitésimales d'ordre zéro et un, et en 
dimension trois, ses transformations
du second ordre, au nombre de trois, sont:
$x(x\,p + y\, q + z\, r)$, 
$y(x\,p + y\, q + z\, r)$
et
$z(x\,p + y\, q + z\, r)$.
}: %%%%%%%%%%%%%%%%%%%%%%%%%%%%%%%%%%%%%%%%%%%%%%%%%%%%%%%%%%%
\[
\big(
\lambda\,x
+
\mu\,y
+
\nu\,z
\big)
\big(
xp+yq+zr
\big)
+
\cdots,
\]
o\`u les termes supprim\'es sont d'ordre trois et d'ordre
sup\'erieur en $x, y, z$; 
mais comme cette condition n\'ecessaire n'est pas
remplie, le groupe [4] ne peut s\^urement pas \^etre semblable \`a un
groupe projectif.

\HEAD{Groupes de $R_3$ pour lesquels deux points ont un seul
invariant.}{Division\,\,V.\,\,\,Chapitre\,\,20.\,\,\,\S\,\,89.}

\sectiondritterV{\sf\S\,\,89.
\\
R\'esolution du probl\`eme pour les groupes r\'eels.}
\label{S-89}
\setcounter{footnote}{0}

Jusqu'\`a maintenant, nous avons consid\'er\'e les variables $x, y, z$
comme des quantit\'es complexes. \`A pr\'esent, nous voulons les
limiter au domaine des nombres r\'eels et nous voulons r\'esoudre le
probl\`eme pos\'e dans l'introduction de ce chapitre, en supposant que
les groupes sont r\'eels et que seules les transformations ponctuelles
qui sont r\'eelles sont admissibles.

{\em Nous recherchons donc maintenant tous les groupes {\small\sf continus
finis r\'eels} de $R_3$ relativement auxquels deux points r\'eels ont
un et un seul invariant, tandis que trois points r\'eels ou plus ne
poss\`edent aucun invariant essentiel.} Avec cela, nous nous
restreignons naturellement aux groupes dont les \'equations finies
autorisent un certain nombre de diff\'erentiations par rapport aux
variables et aux param\`etres\footnote{\, On part d'un groupe de
transformations réel fini continu dont les équations $x_i' = f_i (
x_1, \dots, x_n; \, a_1, \dots, a_m)$, $i = 1,
\dots, n$, ne sont pas forcément analytiques.  En supposant seulement
que les $f_i$ possèdent des dérivées continues d'ordre un et deux par
rapport à toutes leurs variables, Lie et Engel énoncent et démontrent
(Théorème~33, p.~366 du Volume~III) que lorsque le groupe est {\em
transitif}, il est toujours possible d'introduire localement un
changement simultané de paramètres $a \mapsto b = b ( a)$ et de
coordonnées $x \mapsto y = y ( x)$ qui transfère le groupe en un
groupe équivalent $y_i ' = g_i ( y_1, \dots, y_n; \, b_1, \dots,
b_r)$, $i = 1, \dots, n$, dont toutes les équations sont analytiques.
L'énoncé n'est plus valable lorsque le groupe n'est pas transitif,
comme le montre l'exemple $q$, $xq$, $f(x) q$ sur le plan des $(x,
y)$, où $f (x)$ est une fonction différentiable non analytique (\cf \cite{
enlie1893}, p.~368). Le cinquième problème de Hilbert
(1900) demande s'il est
possible d'éliminer l'hypothèse de différentiabilité, et
de supposer seulement que les fonctions $f_i$ sont
continues. En 1952, 
Gleason, et indépendamment, 
Montgomery-Zippin ont établi que tout groupe
topologique localement euclidien, \ie qui forme aussi une
variété topologique, admet une structure
de variété analytique réelle pour laquelle il
constitue un groupe de Lie analytique. 
En 1955, Montgomery-Zippin ont établi le
théorème
suivant: {\em Si un groupe topologique
localement compact $G$ agit transitivement et
effectivement sur un espace topologique $X$
compact et localement connexe, alors $G$
possède une structure de groupe de Lie, et
$X$ possède une structure de variété analytique
réelle de telle sorte que l'action est analytique}
(\cite{ mz1955}; \cite{ gov1997}, pp.~87--88). 
} %%%%%%%%%%%%%%%%%%%%%%%%%%%%%%%%%%%%%%%%%%%%%%%%%%%%%%%%%%%
(\cf p.~365).

\bigskip

D\`es le d\'ebut, il est clair que les groupes recherch\'es sont
enti\`erement transitifs (\cf p.~\pageref{400}). Ainsi,
d'apr\`es la page~366, nous pouvons nous imaginer dans chaque cas
individuel qu'au moyen d'une transformation ponctuelle r\'eelle, de
nouvelles variables $x, y, z$ sont introduites de telle sorte que le
groupe concern\'e est repr\'esent\'e par des \'equations r\'eelles,
dans lesquelles les fonctions qui apparaissent sont des fonctions {\em
analytiques}\, des variables et des param\`etres, au sens de
Weierstra{\ss}. Ainsi, nous devons seulement consid\'erer les groupes
\`a $m$ param\`etres tels que, dans leurs transformations
infinit\'esimales:
\[
\aligned
X_kf
=
\xi_k(x,y,z)\,p
&
+
\eta_k(x,y,z)\,q
+
\zeta_k(x,y,z)\,r
\\
&
{\scriptstyle{(k\,=\,1,\,2\,\cdots\,m)}},
\endaligned
\]
les coefficients $\xi_k$, $\eta_k$, $\zeta_k$ sont des s\'eries de
puissances ordinaires par rapport \`a $x - x_0$, $y - y_0$, $z - z_0$
dans un voisinage d'un point r\'eel quelconque en position
g\'en\'erale, et bien entendu, des s\'eries de puissances avec
seulement des coefficients r\'eels.

Soit: $X_1 f \dots X_m f$, ou bri\`evement $G_m$, un tel groupe
r\'eel \`a $m$ param\`etres qui satisfait notre exigence concernant
les invariants de deux points et plus. Ensuite, l'expression:
\[
e_1\,X_1f
+\cdots+
e_m\,X_mf,
\]
repr\'esente la transformation infinit\'esimale g\'en\'erale d'un
certain groupe \`a $m$ param\`etres $G_m'$ de transformations
complexes, lorsqu'on y consid\`ere $x, y, z$ comme des variables
complexes et $e_1 \dots e_m$ comme des param\`etres
complexes\footnote{\, Techniquement, cette opération consiste tout
simplement à remplacer par des variables complexes les variables
réelles qui apparaissent dans les séries entières par lesquelles on
peut développer tous les coefficient des transformations
infinitésimales $X_k$, $k = 1, \dots, m$. Cela a un sens, parce que
les séries en question demeurent convergentes lorsque les modules des
variables complexes satisfont les mêmes inégalités que celles que
devaient satisfaire les variables réelles initiales pour garantir la
convergence.  Ensuite, pour reconstituer les équations finies du
groupe {\em via} l'application exponentielle, on peut faire prendre
aussi des valeurs complexes aux variables <<\,temporelles\,>> $t$
qui apparaissent dans les
groupes à un paramètre $\exp ( tX_k) ( x)$ engendrés par les $X_k$.
} %%%%%%%%%%%%%%%%%%%%%%%%%%%%%%%%%%%%%%%%%%%%%%%%%%%%%%%%%%%
(\cf p.~362). Il
est maintenant facile de voir que ce groupe $G_m'$ satisfait
\'egalement notre exigence concernant les invariants de deux points ou
plus, c'est-\`a-dire: relativement \`a $G_m'$, deux points ont un et
un seul invariant, tandis qu'un nombre de points sup\'erieur \`a deux
n'a pas d'invariant essentiel. En effet, aussi bien pour $G_m$ que
pour $G_m'$, les invariants de $s$ points sont simplement les
solutions du syt\`eme complet:
\[
X_k^{(1)}f
+
X_k^{(2)}f
+\cdots+
X_k^{(s)}f
=
0
\ \ \ \ \ \ \ \
{\scriptstyle{(k\,=\,1\,\cdots\,m)}},
\]
et le fait que l'on consid\`ere $x, y, z$ comme des variables
r\'eelles ou comme des variables complexes n'a aucune influence sur le
nombre des solutions de ce syst\`eme complet qui sont ind\'ependantes
l'une de l'autre. Donc si deux points ont un et un seul invariant et
si tous les invariants d'un nombre de points sup\'erieur \`a deux
peuvent s'exprimer au moyen des invariants de paires de points,
lorsque $x, y, z$ sont r\'eels, cela vaut aussi lorsque $x, y, z$ sont
complexes et vice-versa: on peut m\^eme toujours assurer que tous les
invariants de $s$ points sont des fonctions r\'eelles des $3s$
variables: $x_1, y_1, z_1, \dots, x_s, y_s, z_s$.

De cette fa\c con, l'\'enonc\'e suivant est d\'emontr\'e: {\em Si un
groupe r\'eel \`a $m$ param\`etres: $X_1f \dots X_m f$ satisfait
notre exigence concernant les invariants de plusieurs points, alors,
d\`es qu'on interpr\`ete $x, y, z$ comme des variables complexes, il
fournit un groupe qui poss\`ede toutes les propri\'et\'es d\'egag\'ees
au \S85; le groupe: $X_1 f \dots X_m f$ est alors n\'ecessairement
\`a six param\`etres et il est semblable \`a l'un des groupes du
Th\'eor\`eme~36, p.~\pageref{Theorem-36}, gr\^ace \`a une
transformation r\'eelle ou complexe de l'espace des $x, y, z$}.
\label{428}

Les groupes \'enum\'er\'es par le Th\'eor\`eme~36 se rangent
dans deux classes. Chaque groupe de la premi\`ere classe\footnote{\,
---\,\,à savoir, les deux groupes imprimitifs~\thetag{ 22} 
et~\thetag{ 23}\,\,---
} %%%%%%%%%%%%%%%%%%%%%%%%%%%%%%%%%%%%%%%%%%%%%%%%%%%%%%%%%%%
laisse invariante une \'equation du second degr\'e:
\def\theequation{55}\begin{equation}
\alpha_{11}dx^2
+
\alpha_{22}dy^2
+
\alpha_{33}dz^2
+
2\alpha_{12}dxdy
+
2\alpha_{23}dydz
+
2\alpha_{31}dzdx
=
0,
\end{equation}
o\`u les $\alpha$ sont des fonctions de $x, y, z$ et o\`u le
d\'eterminant correspondant ne s'annule pas identiquement; en outre
ces groupes sont constitu\'es de telle sorte que les $\infty^2$
\'el\'ements lin\'eaires passant par chaque point fix\'e en position
g\'en\'erale sont transform\'es par une action
\`a trois param\`etres. Les groupes de
la deuxi\`eme classe, \ie les groupes [1] \dots \ [4], se distinguent
en ce que chacun d'entre eux laisse invariante la famille des
$\infty^2$ droites: $x = \text{\rm const.}$, $y = \text{\rm const.}$,
sans stabiliser aucune autre famille de $\infty^2$ courbes.

Si un groupe r\'eel \`a six param\`etres est semblable \`a l'un des
groupes qui appartient \`a la premi\`ere de nos deux classes, alors il
doit \'evidemment laisser invariante une \'equation de la
forme~\thetag{ 55} dans laquelle les $\alpha$ sont des fonctions
r\'eelles de $x, y, z$ et dont le d\'eterminant ne s'annule pas
identiquement; et de plus, les $\infty^2$ \'el\'ements lin\'eaires
passant par chaque point r\'eel fix\'e en position g\'en\'erale
doivent encore \^etre transform\'es par une action \`a trois
param\`etres. Mais d'apr\`es le Th\'eor\`eme~35, p.~391\footnote{\,
%%%%%%%%%%%%%%%%%%%%%%%-------DEBUT--------%%%%%%%%%%%%%%%%%%%%%%%%%%%
Ce théorème établit que si un groupe {\em réel} continu dans un espace
à $n > 2$ variables $x_1, \dots, x_n$ laisse invariante une équation
réelle:
\[
\sum_{k,\,\nu}^{1\dots n}\,
f_{k\nu}(x_1,\dots,x_n)\,dx_k\,dx_\nu
=
0
\]
dont le déterminant ne s'annule pas identiquement, et si en outre, il
est constitué de telle sorte qu'il transforme de la manière la plus
générale possible les $\infty^{ n-1}$ éléments linéaires qui passent
par un tout point fixé en position générale, alors ce groupe possède
soit $\frac{ (n+1)(n+2)}{ 1 \cdot 2}$, soit $\frac{ n(n+1)}{ 1 \cdot
2}$, soit $\frac{ n ( n+1)}{ 1 \cdot 2}$ paramètres.
Dans chacun des trois cas, les groupes sont décrit par
un système de générateurs explicites. 
}, %%%%%%%%%%%%%%%%%%%%%%%%-----FIN-----%%%%%%%%%%%%%%%%%%%%%%%%%%%%%%%
%%%\Fill 
chaque groupe r\'eel \`a six param\`etres de cette esp\`ece est
semblable, par une transformation ponctuelle r\'eelle, ou bien au
groupe des mouvements euclidiens, ou bien au groupe projectif r\'eel
\`a six param\`etres qui laisse invariants le volume ainsi qu'une
conique non-d\'eg\'en\'er\'ee infiniment \'eloign\'ee, ou bien encore
au groupe projectif \`a six param\`etres d'une surface du second
degr\'e non-d\'eg\'en\'er\'ee, qui est repr\'esent\'ee par une
\'equation r\'eelle entre $x$, $y$ et $z$; ici cependant, on doit
encore distinguer trois cas, selon que cette surface du second degr\'e
est imaginaire, ou r\'eelle et non r\'egl\'ee, ou r\'eelle et
r\'egl\'ee.
%%%\Mathematiques 
%%%[[coniques infiniment \'eloign\'ee~?]]
%%%\Mathematiques 
%%%[[classification des quadriques r\'eelles~?]]

\bigskip
\`A pr\'esent, nous recherchons tous les groupes 
r\'eels \`a six param\`etres
qui sont semblables\footnote{\, ---\,\,par une transformation
ponctuelle complexe\,\,---} \`a l'un des groupes: [1] \dots \ [4].

Comme chacun des groupes: [1] \dots \ [4] laisse invariante seulement
la famille: $x = \text{\rm const.}$, $y = \text{\rm const.}$, mais
nulle autre famille de $\infty^2$ courbes, chaque groupe r\'eel \`a
six param\`etres: $X_1 f \dots X_6 f$ qui est semblable \`a l'un
d'entre eux, laisse invariante une et seulement une famille de
$\infty^2$ courbes. De plus, cette famille doit n\'ecessairement
\^etre r\'eelle, car si elle \'etait imaginaire, la famille des
courbes imaginaires conjugu\'ees resterait aussi invariante par le
groupe: $X_1 f \dots X_6 f$, puisqu'il est r\'eel, et il y aurait
alors deux familles invariantes distinctes, ce qui ne peut pas \^etre
le cas. Par une transformation ponctuelle r\'eelle, nous pouvons
maintenant transformer cette unique famille existante de $\infty^2$
courbes en la famille: $x = \text{\rm const.}$, $y = \text{\rm
const.}$, et par l\`a, nous pouvons\,\,---\,\,d\`es le
d\'ebut\,\,---\,\,rapporter le groupe r\'eel: $X_1 f \dots X_6 f$
\`a la forme:
\[
\aligned
X_kf
=
\xi_k(x,y)\,p
&
+
\eta_k(x,y)\,q
+
\zeta_k(x,y,z)\,r
\\
&
{\scriptstyle{(k\,=\,1,\,2\,\cdots\,6)}},
\endaligned
\] 
o\`u naturellement les $\xi_k$, $\eta_k$, $\zeta_k$ sont
des fonctions r\'eelles de leurs arguments.

Comme pr\'ec\'edemment ({\em voir} p.~\pageref{412a}), on obtient
maintenant aussi les six transformations infinit\'esimales:
\[
\overline{X}_kf
=
\xi_k(x,y)\,p
+
\eta_k(x,y)\,q
\ \ \ \ \ \ \ \ \
{\scriptstyle{(k\,=\,1\,\cdots\,6)}},
\]
qui forment \'evidemment un groupe r\'eel dans l'espace des deux
variables $x, y$; et gr\^ace au raisonnement de la
page~\pageref{412b}, on comprend imm\'ediatement que ce groupe
r\'eduit doit poss\`eder au moins cinq param\`etres et que, si on
l'interpr\`ete comme un groupe agissant sur un plan, alors chaque
famille r\'eelle ou imaginaire de $\infty^1$ courbes qu'il laisse
invariante doit \^etre transform\'ee par une action \`a trois
param\`etres.

Si maintenant le groupe r\'eduit: $\overline{ X}_1 f \dots
\overline{ X}_6 f$ ne laisse invariante 
absolument aucune famille de $\infty^1$ courbes dans le plan, alors
d'apr\`es les pages~370 sq.,
%%%\Fill
il est semblable, {\em via}\, une transformation ponctuelle r\'eelle,
\`a l'un des deux groupes: (I), p.~\pageref{413} et \thetag{ V},
p.~\pageref{422}.  Si, d'un autre c\^ot\'e, il laisse invariante
seulement une famille de $\infty^1$ courbes, alors cette famille est
n\'ecessairement r\'eelle et par suite, sur la base des
d\'eveloppements des pages~378 sq., 
%%%\Fill 
le groupe peut \^etre
rapport\'e, {\em via}\, une transformation ponctuelle r\'eelle, \`a
l'une des trois formes:
\thetag{ II}, \thetag{ III}, p.~\pageref{413} et \thetag{ VI},
p.~\pageref{422}. Si enfin le groupe: $\overline{ X}_1 f \dots
\overline{ X}_6f$ laisse invariantes deux familles distinctes de
$\infty^1$ courbes, alors il peut, lorsque ces deux familles sont
r\'eelles, recevoir la forme~\thetag{ IV}, p.~\pageref{413}, et
lorsqu'au contraire elles sont imaginaires conjugu\'ees, la forme:
\def\theequation{VII}\begin{equation}
\left\{
\aligned
&
p,\ \
q,\ \
xp+yq,\ \
yp-xq
\\
&
(x^2-y^2)\,p
+
2xyq,\ \
2xyp+
(y^2-x^2)\,q,
\endaligned\right.
\end{equation}
(\cf p.~380) 
%%%\Fill
et ceci, les deux fois, {\em via}\, une transformation
ponctuelle r\'eelle.

Si le groupe r\'eduit: $\overline{ X}_1 f \dots \overline{ X}_6f$ a
l'une des formes \thetag{ I} \dots\: \thetag{ VI}, alors en conduisant
exactement les m\^emes calculs que pr\'ec\'edemment ({\em voir} de la
p.~\pageref{414} \`a la p.~\pageref{425}), on d\'etermine toutes les
formes possibles du groupe: $X_1 f \dots X_6f$, et alors, comme on
s'en convainc facilement, on a seulement besoin, lors de
l'effectuation de ces calculs, de transformations ponctuelles qui
respectent pas \`a pas la r\'ealit\'e du groupe. Par l\`a, nous
voyons que les formes~\thetag{ I} et~\thetag{ II} du groupe r\'eduit
ne conduisent \`a aucun groupe r\'eel poss\'edant les qualit\'es ici
requises, et nous voyons aussi que chaque groupe r\'eel: $X_1 f
\dots X_6 f$ qui satisfait nos exigences et dont le groupe r\'eduit
poss\`ede l'une des formes~\thetag{ III}, \dots\:, \thetag{ VI} est
semblable, {\em via}\, une transformation ponctuelle r\'eelle, \`a l'un des
groupes: [2]\footnote{\, C'est le bon ordre.}, [1], [3], [4]. Et
naturellement, le param\`etre $c$ qui appara\^{\i}t dans les deux
groupes [1] et [2] doit avoir une valeur r\'eelle.

Enfin, si le groupe r\'eduit: $\overline{ X}_1 f \dots \overline{
X}_6 f$ a la forme~\thetag{ VII}, alors le groupe: $X_1 f \dots X_6
f$ est semblable, {\em via}\, une transformation imaginaire, \`a l'un des
groupes de la forme~[1], car le groupe~\thetag{ IV}, p.~\pageref{413}
provient du groupe~\thetag{ VII}, lorsqu'on introduit comme nouvelles
variables $x + iy$ et $x - iy$ \`a la place de $x$ et de $y$. 
%%%({\em voir}\, p.~336). 
%%%\Fill

Nous ne voulons n\'eanmoins faire aucun usage de
cette circonstance, car il est pr\'ef\'erable de d\'eterminer
directement tous les groupes r\'eels \`a six param\`etres qui
satisfont nos exigences et dont les groupes r\'eduits 
ont la forme~\thetag{
VII}.

\bigskip

Ainsi, nous cherchons maintenant tous les groupes r\'eels \`a six
param\`etres de la forme:
\def\theequation{56}\begin{equation}
\left\{
\aligned
p+\varphi_1\,r,\ \
q+\varphi_2\,r,\ \
xp+yq+\varphi_3\,r,\ \
yp-xq+\varphi_4\,r
\\
(x^2-y^2)\,p+2xyq+\varphi_5\,r,\ \
2xyp+(y^2-x^2)\,q+\varphi_6\,r
\endaligned\right.
\end{equation}
relativement auxquels deux points ont un et un seul invariant, tandis
qu'un nombre de points sup\'erieur \`a deux n'a aucun invariant
essentiel; ici, $\varphi_1 \dots \varphi_6$ sont naturellement des
fonctions de $x, y, z$.

Exactement comme aux pages~\pageref{414} sq., nous voyons\footnote{\,
Dans un premier temps, on redresse les deux premières transformations
en $p$ et en $q$.  Le crochet de la troisième avec $p$ et avec $q$
montre que $\varphi_3$ ne dépend que de $z$. Si $\varphi_3 \equiv 0$,
le travail est fini: il suffit de prendre $a = 0$.  Si $\varphi_3 \not
\equiv 0$, en introduisant $\int \, \frac{ 1}{ \varphi_3} \, dz$ 
comme nouveau $z$ (et en relocalisant bien sûr les
considérations autour d'un point où 
$\varphi_3$ ne s'annule pas), on redresse
la troisième transformation en $xp + yq + r$, c'est-à-dire
que $a = 1$. On 
réunit les deux cas en introduisant une constante $a$, qui
sera de toute façon ramenée à être égale à $0$ ou à $1$
à la fin de la démonstration.  
} %%%%%%%%%%%%%%%%%%%%%%%%%%%%%%%%%%%%%%%%%%%%%%%%%%%%%%%%%%% 
\deutsch{sehen wir ein}
que les trois premi\`eres transformations infinit\'esimales~\thetag{
56} peuvent \^etre rapport\'ees \`a la forme: 
\[
p,\ \
q,\ \
xp+yq+ar,
\]
o\`u $a$ d\'esigne une constante r\'eelle, {\em via}\, une
transformation ponctuelle, et pour pr\'eciser, {\em via}\, une
transformation ponctuelle r\'eelle.

En calculant le crochet de: $yp - xq + \varphi_4 \, r$ avec $p$ et
avec $q$, nous voyons que $\varphi_4$ est ind\'ependant de $x$ et de
$y$. De plus, nous avons:
\[
\big[
xp+yq+ar,\ \
yp-xq+\varphi_4(z)\,r
\big]
=
a\,\varphi_4'(z)\,r,
\]
et donc, si $a$ ne s'annule pas, alors $\varphi_4$ est aussi
ind\'ependant de $z$; si au contraire $a = 0$ et $\varphi_4 \neq 0$,
nous introduisons la fonction r\'eelle:
\[
\int\,
\frac{dz}{\varphi_4(z)}
\]
comme nouvelle variable $z$ et nous parvenons par cela \`a faire que
$\varphi_4 = 1$. Par cons\'equent, nous pouvons r\'eunir ces deux cas
en rapportant les quatre premi\`eres transformations
infinit\'esimales~\thetag{ 56} \`a la forme:
\[
p,\ \
q,\ \
xp+yq+ar,\ \
yp-xq+br,
\]
o\`u $b$ repr\'esente \'egalement une constante r\'eelle.

\`A pr\'esent, calculons les crochets:
\[
\aligned
\big[
p,\ \
(x^2-y^2)\,p+2xyq+\varphi_5\,r
\big]
&
=
2(xp+yq)
+
\frac{\partial\varphi_5}{\partial x}\,r
\\
\big[
q,\ \
(x^2-y^2)\,p+2xyq+\varphi_5\,r
\big]
&
=
-2(yp-xq)
+
\frac{\partial\varphi_5}{\partial y}\,r
\\
\big[
xp+yq+ar,\ \
(x^2-y^2)\,p+2xyq+\varphi_5\,r
\big]
&
=
(x^2-y^2)\,p
+
2xyq
+
\\
&
+
\Big\{
x\,\frac{\partial\varphi_5}{\partial x}
+
y\,\frac{\partial\varphi_5}{\partial y}
+
a\,\frac{\partial\varphi_5}{\partial z}
\Big\}\,r
\\
\big[
yp-xq+br,\ \
(x^2-y^2)\,p+2xyq+\varphi_5\,r
\big]
&
=
2xyp
+
(y^2-x^2)\,q
+
\\
&
+
\Big\{
y\,\frac{\partial\varphi_5}{\partial x}
-
x\,\frac{\partial\varphi_5}{\partial y}
+
b\,\frac{\partial\varphi_5}{\partial z}
\Big\}\,r
\\
\big[
yp-xq+br,\ \
2xyp+(y^2-x^2)\,q+\varphi_6\,r
\big]
&
=
-(x^2-y^2)\,p
-
2xyq
+
\\
&
+
\Big\{
y\,\frac{\partial\varphi_6}{\partial x}
-
x\,\frac{\partial\varphi_6}{\partial y}
+
b\,\frac{\partial\varphi_6}{\partial z}
\Big\}\,r
\endaligned
\]
d'o\`u nous obtenons: 
\def\theequation{57}\begin{equation}
\left\{
\aligned
\frac{\partial\varphi_5}{\partial x}
=
2a,
\ \ \ \ \
\frac{\partial\varphi_5}{\partial y}
&
=
-2b
\\
x\,\frac{\partial\varphi_5}{\partial x}
+
y\,\frac{\partial\varphi_5}{\partial y}
+
a\,\frac{\partial\varphi_5}{\partial z}
&
=
\varphi_5
\\
y\,\frac{\partial\varphi_5}{\partial x}
-
x\,\frac{\partial\varphi_5}{\partial y}
+
b\,\frac{\partial\varphi_5}{\partial z}
&
=
\varphi_6
\\
y\,\frac{\partial\varphi_6}{\partial x}
-
x\,\frac{\partial\varphi_6}{\partial y}
+
b\,\frac{\partial\varphi_6}{\partial z}
&
=
-\varphi_5.
\endaligned\right.
\end{equation}
De l\`a, on tire tout d'abord:
\[
\varphi_5
=
2(ax-by)
+
\omega(z),
\ \ \ \ \ \ \ \
a\,\omega'(z)
=
\omega(z)
=
a^2\omega''(z)
\]
et par cons\'equent: 
\[
\varphi_6
=
2(bx+ay)
+
b\,\omega'(z);
\]
mais si nous ins\'erons cette valeur de $\varphi_6$ dans la derni\`ere
des \'equations~\thetag{ 57}, viennent alors deux \'equations:
\[
b^2\,\omega''(z)
=
-\omega(z)
=
-a^2\,\omega''(z),
\]
qui, pour des valeurs r\'eelles de $a$ et de $b$, ne
peuvent alors valoir que lorsque $\omega ( z)$ s'annule\footnote{\,
On notera l'organisation fine des calculs: le terme $- \omega ( z)$,
dont la fin de la démonstration conclut qu'il
s'annule identiquement,  
avait déjà été placé au centre plus haut, 
et maintenant, dans cette dernière équation,
il est clair, puisque $b^2$ et $-a^2$ sont
de signe opposé, que $- \omega ( z)$ au centre
est à la fois de signe positif et de signe négatif, 
donc nul. 
}. %%%%%%%%%%%%%%%%%%%%%%%%%%%%%%%%%%%%%%%%%%%%%%%%%%%%%%%%%%%

Ainsi, nous sommes parvenus au groupe:
\def\theequation{58}\begin{equation}
\label{eq-58-p-432}
\boxed{
\aligned
p,\ \
q,\ \
xp+yq+ar,\ \
yp-xq+br
\\
(x^2-y^2)\,p+2xyq+2(ax-by)\,r
\\
2xyp+(y^2-x^2)\,q+2(bx+ay)\,r
\endaligned
}\,\,.
\end{equation}
Ici, les constantes $a$ et $b$ ne doivent naturellement pas s'annuler
toutes les deux, parce que sinon,\label{432} le groupe serait
intransitif; si $b=0$, nous pouvons toujours parvenir \`a faire que $a
= 1$ en introduisant une nouvelle variable $z$, et si au contraire,
$b$ ne s'annule pas, on peut toujours faire que $b$ soit simplement
\'egal \`a $1$.

Comme on le trouve ais\'ement, relativement au groupe~\thetag{ 58},
deux points $x_1, y_1, z_1$ et $x_2, y_2, z_2$ ont un et un seul
invariant, \`a savoir celui-ci:
\[
\label{invariant-58}
z_1+z_2
-
a\,\log
\big\{
(x_2-x_1)^2
+
(y_2-y_1)^2
\big\}
+
2\,b\,\arctan\,
\frac{y_2-y_1}{x_2-x_1}.
\] 
En raisonnant comme pr\'ec\'edemment (\cf par exemple
p.~\pageref{419}), on pourrait aussi se convaincre facilement que si
$a$ et $b$ ne s'annulent pas tous deux, trois points et plus n'ont
aucun invariant essentiel; toutefois, nous pr\'ef\'erons d\'emontrer
cela en indiquant une transformation ponctuelle imaginaire qui envoie
le groupe~\thetag{ 58} sur le groupe~[1] p.~\pageref{426}, puisque
nous savons d\'ej\`a du groupe~[1] qu'il satisfait toutes nos
exigences.

Si nous introduisons dans le groupe~\thetag{ 58} les nouvelles
variables:
\[
x_1
=
x+iy,
\ \ \ \ \ 
y_1
=
x-iy,
\ \ \ \ \
z_1
=
z,
\]
nous obtenons le groupe: 
\[
\aligned
p_1+q_1,\ \
i(p_1-q_1),\ \
x_1p_1+y_1q_1+ar_1,\ \
i(y_1q_1-x_1p_1)+br_1
\\
x_1^2p_1+y_1^2q_1+\big\{(a+ib)\,x_1+(a-ib)\,y_1\big\}\,r_1
\ \ \ \ \ \ \ \ \ \ \ \ \ \ \ \
\\
i(y_1^2q_1-x_1^2p_1)+i\big\{(a-ib)\,y_1-(a+ib)\,x_1\big\}\,r_1.
\ \ \ \ \ \ \ \ \ \ \ \ 
\endaligned
\]
Maintenant, comme $a+ib \neq 0$, nous pouvons 
encore introduire $\frac{ 2z_1}{ a +
ib}$ comme nouvelle variable $z_1$ et obtenir ainsi un groupe qui a
exactement la forme~[1]. Ainsi, le groupe~\thetag{ 58} est semblable
au groupe~[1] {\em via}\, la transformation imaginaire:
\[
x_1
=
x+iy,
\ \ \ \ \ \
y_1
=
x-iy,
\ \ \ \ \ \
z_1
=
{\textstyle{\frac{2z}{a+ib}}},
\]
et en fait, aux param\`etres $a$ et $b$ du groupe~\thetag{ 58}
correspond alors la valeur:
\def\theequation{59}\begin{equation}
{\textstyle{
c\,
=\,
\frac{a-ib}{a+ib}
}}
\end{equation}
du param\`etre $c$ du groupe~[1].

Avec ceci, on a d\'emontr\'e que le groupe r\'eel~\thetag{ 58}
satisfait toutes nos exigences, et on a maintenant trouv\'e en toute
g\'en\'eralit\'e tous les groupes r\'eels qui les satisfont.

\medskip
{\bf Th\'eor\`eme~37.}\label{Theorem-37-p-433}
{\em Si un groupe continu fini r\'eel de l'espace des $x, y, z$ est
donn\'e de telle sorte que relativement \`a lui, deux points ont un et
un seul invariant et que les invariants d'un nombre de points
sup\'erieur \`a deux peuvent tous s'exprimer au moyen des invariants
des paires de points, alors ce groupe poss\`ede six param\`etres et il
est semblable \deutsch{\"ahnlich}, {\rm via} une transformation
ponctuelle r\'eelle de l'espace, \`a l'un des onze groupes suivants}:

\begin{small}
\[\label{433-sq}
{\bf 1}\ \ \ \ \
\boxed{\ \
p,\ \
q,\ \
r,\ \
xq-yp,\ \
yr-zq,\ \
zp-xr\ \
}
\]

\[
{\bf 2}\ \ \ \ \
\boxed{\ \
p,\ \
q,\ \
r,\ \
xq-yp,\ \
yr+zq,\ \
zp+xr\ \
}
\]

\[
{\bf 3}\ \ \ \ \
\boxed{\ \
p+xU,\ \
q+yU,\ \
r+zU,\ \
xq-yp,\ \
yr-zq,\ \
zp-xr\ \
}
\]

\[
{\bf 4}\ \ \ \ \
\boxed{\ \
p-xU,\ \
q-yU,\ \
r-zU,\ \
xq-yp,\ \
yr-zq,\ \
zp-xr\ \
}
\]

\[
{\bf 5}\ \ \ \ \
\boxed{\ \
p-xU,\ \
q-yU,\ \
r+zU,\ \
xq-yp,\ \
yr+zq,\ \
zp+xr\ \
}
\]

\[
{\bf 6}\ \ \ \ \
\boxed{\ \
\aligned
p,\ \
q,\ \
xp+yq+cr,\ \
yp-xq+r
\ \ \ \ \ \ \ \ \ \ \ \ \ \ 
\ \ \ \ \ \ \ \ \ \ \ \ \ \ 
\\
(x^2-y^2)p+2xyq+2(cx-y)r,\ \
2xyp+(y^2-x^2)q+2(x+cy)r\ \
\endaligned
}
\]

\[\label{432-1}
{\bf 7}\ \ \ \ \
\boxed{\ \
\aligned
p,\ \
q,\ \
xp+yq+r,\ \
yp-xq
\ \ \ \ \ \ \ \ \ \ \ \ \ \ 
\ \ \ \ \ \ \ \ \ \ \ \ \ \ 
\\
(x^2-y^2)p+2xyq+2xr,\ \
2xyp+(y^2-x^2)q+2yr\ \
\endaligned
}
\]

\[
{\bf 8}\ \ \ \ \
\boxed{\ \
\aligned
p,\ \
q,\ \
xp+r,\ \
yq+cr,\ \
x^2p+2xr,\ \
y^2q+2cyr\ \
\\
{\scriptstyle{(c\,\neq\,0)}}
\ \ \ \ \ \ \ \ \ \ \ \ \
\ \ \ \ \ \ \ \ \ \ \ \ \
\ \ \ \ \ \ \ \ \ \ \ \ \
\endaligned
}
\]

\[
{\bf 9}\ \ \ \ \
\boxed{\ \
\aligned
p,\ \
q,\ \
xq+r,\ \
x^2q+2xr,\ \
xp+yq+cr\ \
\\
x^2p+2xyq+2(y+cx)r
\ \ \ \ \ \ \ \ \ \ \ \ \ \ \
\endaligned
}
\]

\[
{\bf 10}\ \ \ \ \
\boxed{\ \
p-yr,\ \
q+xr,\ \
r,\ \
xq,\ \
xp-yq,\ \
yp\ \
}
\]

\[
{\bf 11}\ \ \ \ \
\boxed{\ \
\aligned
p,\ \
q,\ \
r,\ \
xq+yr,\ \
2xp+yq\ \
\\
x^2p+xyq+{\textstyle{\frac{1}{2}}}y^2r
\ \ \ \ \ \ \ \ \ \ \
\endaligned
}
\]
\end{small}

\noindent
{\em Ici, le param\`etre r\'eel $c$ est \`a chaque fois essentiel et ne
peut pas \^etre supprim\'e.}

\medskip

Pour faciliter la compr\'ehension de ce tableau, nous ajoutons encore
les remarques suivantes.

Pour abr\'eger, on \'ecrit chaque fois $U$ \`a la place de $xp + yq +
zr$.

Le groupe~{\bf 2} laisse invariante la conique qui d\'ecoupe la
sph\`ere r\'eelle: $x^2 + y^2 - z^2 = 0$ sur le plan \`a l'infini, et
ce groupe provient du groupe~{\bf 1} constitu\'e des mouvements
euclidiens lorsqu'on remplace $z$ par $iz$.
%%%\Mathematiques

Le groupe~{\bf 4} laisse invariante la surface r\'eelle du second
degr\'e non r\'egl\'ee: $x^2 + y^2 + z^2 = 1$ et le groupe~{\bf 5},
la surface r\'eelle r\'egl\'ee: $x^2 + y^2 - z^2 = 1$. Ces groupes
sont obtenus \`a partir du groupe~{\bf 3} de la surface: $x^2 + y^2 +
z^2 + 1 = 0$, lorsqu'on introduit soit $ix, iy, iz$, soit $ix, iy,
z$ comme nouvelles variables, et l'on peut alors tr\`es simplement
indiquer quel est l'invariant de deux points qui leur correspond,
gr\^ace aux pages~\pageref{410} sq.

\`A la place du groupe~\thetag{ 58} se trouvent dans notre tableau les
deux groupes~{\bf 6} et~{\bf 7}, parce que, d'apr\`es la
page~\pageref{432}, l'un des deux param\`etres $a$ et $b$ peut
toujours \^etre rendu simplement \'egal \`a $1$; 
cependant, on doit faire
attention au fait que dans le groupe~{\bf 6}, la valeur du param\`etre
$c$ doit \^etre limit\'ee aux seules valeurs qui sont $\geqslant 0$,
car les deux groupes qui correspondent aux param\`etres $+ c$ et $- c$
sont semblables, comme on s'en convainc imm\'ediatement, lorsqu'on
intervertit $x$ et $y$ et qu'on remplace $z$ par $-z$.

Les groupes~{\bf 8} \dots \ {\bf 11} sont les groupes [1] \dots \:
[4] de la page~\pageref{426}; on doit encore mentionner
sp\'ecialement que dans le groupe~{\bf 8}, on a seulement besoin de
donner \`a $c$ des valeurs telles que $c^2 \leqslant 1$, car par la
transformation:
\[
x_1=x,
\ \ \ \ \ \
y_1=y,\ \ \ \ 
z_1=
\textstyle{\frac{z}{c}},
\]
le groupe s'envoie sur un groupe de la m\^eme forme, dont le
param\`etre poss\`ede la valeur $\frac{ 1}{ c}$.

\HEAD{Groupes de $R_3$ pour lesquels deux points ont un seul
invariant.}{Division\,\,V.\,\,\,Chapitre\,\,20.\,\,\,\S\,\,90.}

\sectiondritterV{\sf\S\,\,90.}
\label{S-90}
\setcounter{footnote}{0}

En quelques mots, nous voulons maintenant expliquer comment on peut
r\'esoudre dans le plan le probl\`eme que nous venons de r\'esoudre
pour l'espace trois fois \'etendu.

Si un groupe continu fini du plan doit \^etre constitu\'e de telle
sorte que relativement \`a son action, deux points ont un et un seul
invariant, tandis qu'un nombre de points sup\'erieur \`a deux n'a pas
d'invariant essentiel, alors on v\'erifie en premier lieu, comme au
\S~85, que le groupe est transitif, qu'il doit poss\`eder trois
param\`etres et qu'il ne peut pas contenir deux transformations
infinit\'esimales ayant les m\^emes courbes int\'egrales.

Inversement, il est clair que relativement \`a tout groupe transitif
$G$ du plan \`a trois param\`etres, deux points ont toujours un et un
seul invariant\footnote{ Soit $X_1, X_2, X_3$ les trois
transformations infinitésimales d'un tel groupe du plan. Un invariant
quelconque $J = J ( x_1, y_1; \, x_2, y_2)$ entre deux points
satisfait les équations $0 \equiv X_k^{ (1)} J + X_k^{ ( 2)} J$, pour
$k=1, 2, 3$. Ces $3$ équations (intégrables au sens de Frobenius) sur
un espace de dimension $4$ sont automatiquement indépendantes dès que
le groupe est transitif.  Aussi existe-t-il exactement $4 - 3 = 1$
invariant $J ( x_1, y_1; \, x_2, y_2)$.
}. %%%%%%%%%%%%%%%%%%%%%%%%%%%%%%%%%%%%%%%%%%%%%%%%%%%%%%%%%%%
Si de plus $G$ ne contient pas deux transformations
infinit\'esimales ayant les m\^emes courbes int\'egrales, on peut
d\'emontrer qu'un nombre de points sup\'erieur \`a deux n'a jamais
d'invariant essentiel\footnote{\, 
La Proposition~4 p.~\pageref{Satz-4} traitait déjà 
le cas de la dimension $3$. 
}. %%%%%%%%%%%%%%%%%%%%%%%%%%%%%%%%%%%%%%%%%%%%%%%%%%%%%%%%%%%

En effet, si l'on fixe un point $P_1$ en position g\'en\'erale, il
reste encore un sous-groupe \`a un param\`etre de $G$ par l'action
duquel les autres points du plan se meuvent sur $\infty^1$ courbes,
\`a savoir sur les pseudocercles de centre $P_1$ appartenant \`a $G$.
S'il y avait seulement $\infty^1$ pseudocercles, ces pseudocercles ne
d\'ependraient pas de leur centre, et par suite, la transformation
infinit\'esimale de $G$ qui laisse invariant le point $P_1$ aurait les
m\^emes courbes int\'egrales que la transformation infinit\'esimale
ind\'ependante d'elle
qui laisse invariant un autre point quelconque en position
g\'en\'erale. Mais comme cela ne doit pas se
produire, nous pouvons en conclure qu'il y a au minimum $\infty^2$
pseudocercles. Si maintenant on fixe, outre le point $P_1$, encore un
deuxi\`eme point $P_2$ en position g\'en\'erale, alors les deux
pseudocercles de centre $P_1$ et $P_2$, qui passent par un
troisi\`eme point en
position g\'en\'erale $P_3$, se coupent en $P_3$, et mis \`a part
$P_3$, ils se coupent au maximum en des points 
isol\'es\footnote{\,
Analyticité ...
}; %%%%%%%%%%%%%%%%%%%%%%%%%%%%%%%%%%%%%%%%%%%%%%%%%%%%%%%%%%%
par
cons\'equent, d\`es que $P_1$ et $P_2$ sont fix\'es, chaque autre
point $P$ du plan reste en g\'en\'eral au repos, et pour
pr\'eciser, \`a cause des deux invariants qui appartiennent aux deux
paires de points $P_1, P$ et $P_2, P$. Il en d\'ecoule que, sous les
hypoth\`eses pos\'ees, tous les invariants d'un nombre de points
sup\'erieur \`a trois peuvent \^etre exprim\'es au moyen des
invariants des paires de points (\cf aussi p.~\pageref{409} sq.).

Si donc nous cherchons tous les groupes du plan qui poss\`edent la
propri\'et\'e connue concernant les invariants de plusieurs points, et
si tout d'abord, nous ne tenons pas compte de la condition de
r\'ealit\'e, nous avons seulement besoin de rechercher, parmi les
groupes \`a trois param\`etres qui sont rassembl\'es 
\`a la page~57\footnote{\,
%%%%%%%%%%%%%%%%%%%%%%%-------DEBUT--------%%%%%%%%%%%%%%%%%%%%%%%%%%%
\`A cet endroit sont listés tous les groupes locaux
à un, deux, trois ou quatre paramètres qui agissent
en dimension deux. 
}, %%%%%%%%%%%%%%%%%%%%%%%%-----FIN-----%%%%%%%%%%%%%%%%%%%%%%%%%%%%%%%
%%%\Fill
ceux qui sont transitifs et qui n'ont jamais deux transformations
infinit\'esimales poss\'edant les m\^emes courbes int\'egrales.
Nous trouvons de cette fa\c con seulement les
quatre groupes suivants: 
\def\theequation{60}\begin{equation}
\left\{
\begin{array}{c}
\boxed{
\aligned
&
\ \
p,\ \
q,\ \
xp+cyq\ \
\\
&
\ \ \ \ \ \ \ \ \ \ \ \ 
{\scriptstyle{(c\,\neq\,0)}}
\endaligned
}
\\
\\
\boxed{
\ \
p+x^2p+xyq,\ \
q+xyp+y^2q,\ \
yp-xq\ \
}
\\
\\
\boxed{
\ \
xq,\ \
xp-yq,\ \
yp\ \
}
\\
\\
\boxed{
\ \
p,\ \
q,\ \
xp+(x+y)\,q\ \
}\,.
\end{array}
\right.
\end{equation}
Ici, le groupe:
\[
p,\ \
xp+yq,\ \
x^2p+(2xy+y^2)\,q
\]
est remplac\'e par le groupe projectif de la conique: $x^2 + y^2 + 1
= 0$ (\cf p.~70, p.~76 et p.~88). Pareillement, le groupe:\, $p,
\ \ 2xp+ yq, \ \ x^2 p + xy q$ est remplac\'e par un groupe projectif
qui lui est semblable (\cf p.~95).

Les invariants des deux points $x_1, \, y_1$ et $x_2, \, y_2$ pour les
groupes~\thetag{ 60}, \'ecrits l'un apr\`es l'autre, ont l'expression
suivante:
\def\theequation{61}\begin{equation}
\left\{
\aligned
&
c\log(x_2-x_1)^2-\log(y_2-y_1)^2,
\ \ \ \ \
\textstyle{\frac{
(x_2-x_1)^2
+
(y_2-y_1)^2
+
(x_1y_2-x_2y_1)^2}{
(1+x_1x_2+y_1y_2)^2
}}
\\
&
x_1y_2-x_2y_1,
\ \ \ \ \ \ \ \ \ \ \
(x_2-x_1)\,
e^{-\frac{y_2-y_1}{x_2-x_1}}.
\endaligned\right.
\end{equation}

%%%\Fill
Si l'on veut avoir tous les groupes r\'eels du plan qui poss\`edent les
caract\'eristiques requises, on doit ajouter aux groupes~\thetag{ 60}
encore les deux suivants (\cf pp.~370 sq. et 
p.~\pageref{428}):
\def\theequation{62}\begin{equation}
\left\{
\begin{array}{c}
\boxed{
\ \
p,\ \,
q,\ \
yp-xq+c(xp+yq)\ \
}
\\ 
\\
\boxed{
\ \
p-x^2p-xyq,\ \
q-xyp-y^2q,\ \
yp-xq\ \
}\,.
\end{array}
\right.
\end{equation}

%%%\Mathematiques

\renewcommand{\thefootnote}{\fnsymbol{footnote}}
\noindent
Le premier d'entre eux fournit l'invariant suivant\footnote[1]{\,
Monsieur de Helmholtz s'est en tout cas pr\'eoccup\'e, d\'ej\`a en
1868, d'une recherche dont le but peut s'\'enoncer, avec notre
mani\`ere actuelle de nous exprimer, comme suit: d\'eterminer les
groupes du plan relativement auxquels deux points ont un et un seul
invariant, tandis qu'un nombre de points sup\'erieur \`a deux n'a pas
d'invariant essentiel. \`A cette \'epoque, il avait d\'ej\`a
remarqu\'e que dans le plan, on peut imaginer une famille de
mouvements relativement auxquels deux points poss\`edent
l'invariant~\thetag{ 63}.
}:
\def\theequation{63}\begin{equation}
\big\{
(x_2-x_1)^2
+
(y_2-y_1)^2
\big\}\,
e^{2c\arctan\,
\frac{y_2-y_1}{x_2-x_1}};
\end{equation}
le second provient du deuxi\`eme 
des groupes~\thetag{ 60}, lorsqu'on
introduit $ix$, $iy$ comme nouvelles variables \`a la place de $x, y$
et par l\`a, l'invariant qui lui appartient peut \^etre imm\'ediatement
indiqu\'e\footnote[2]{\,
D\'ej\`a dans les ann\'ees 1874--76, Lie a d\'etermin\'e tous les
groupes \`a trois param\`etres du plan, sans toutefois ajouter
l'\'ecriture de tous les calculs jusque dans les d\'etails
(G\"ott. Nachr. v. 1874; norwegisches Archiv 1878). En l'ann\'ee
1884, il a donn\'e ensuite une \'enum\'eration d\'etaill\'ee des
diff\'erentes formes normales auxquelles tous les groupes de ce genre
peuvent \^etre rapport\'ees. Plus tard (en 1887) Monsieur Poincar\'e
s'est pr\'eoccup\'e des fondements de la g\'eom\'etrie du plan, et
dans son travail relatif \`a ce sujet (Bull. de la Soc. Math.,
vol. 15), il a aussi pris en consid\'eration les recherches relativement
anciennes de Lie sur la th\'eorie des groupes, tandis qu'il ne
connaissait ni le travail de Lie de 1884, ni la premi\`ere
communication de Lie sur les fondements de la g\'eom\'etrie (1886),
pas plus que l'\'etude helmholtzienne de l'ann\'ee 1868. Les
conclusions int\'eressantes auxquelles Monsieur Poincar\'e est parvenu
r\'esultent en v\'erit\'e imm\'ediatement des recherches 
anciennes de Lie.
}.
\renewcommand{\thefootnote}{\arabic{footnote}}

%%%%%%%%%%%%%%%%%%%%%%%%%%%%%%%%%%%%%%%%%%%%%%%%%%%%%%%%%%%%%%%%%%%%%

\newpage

% 33   :   438--470

\setcounter{footnote}{0}

$\:$
\bigskip\bigskip\bigskip

\centerline{\Large Chapitre~21.}
\label{Chapitre-21}
\thispagestyle{empty}

\bigskip

\noindent
\centerline{\large\bf
Critique des recherches helmholtziennes.
}

\bigskip\medskip
Monsieur de Helmholtz a apport\'e, dans son travail d\'ej\`a cit\'e
\`a la page~\pageref{396}, un traitement du probl\`eme de
Riemann-Helmholtz que nous voulons maintenant examiner de plus pr\`es.

Nous commen\c cons (dans le \S~91) par restituer la teneur des axiomes
qu'il a pos\'es. Ensuite (dans les \S~92 et~93), nous \'etablissons
quelles contenus ces axiomes peuvent recevoir lorsqu'on se sert de la
mani\`ere dont on s'exprime dans la th\'eorie des groupes. Dans le
\S~94 nous fournissons une critique des conclusions que Monsieur de
Helmholtz a tir\'ees de ses axiomes au sujet de l'espace \'etendu dans
trois directions, et nous montrons qu'{\em il est parvenu \`a son
r\'esultat final en pr\'esupposant tacitement toute une s\'erie
d'hypoth\`eses qui sont absolument incorrectes}.

Si l'on veut supprimer ces d\'efauts tr\`es substantiels des
d\'eveloppements helmholtziens, il se pr\'esente deux voies
diff\'erentes: on peut soit laisser compl\`etement de c\^ot\'e les
axiomes de M. de Helmhotz et seulement envisager ses calculs, soit
partir de ses axiomes comme point de d\'epart, sans prendre en
consid\'eration ses calculs.

Nous nous engageons dans la premi\`ere voie au \S~96, \label{438} o\`u
nous montrons {\em quels axiomes doivent \^etre pos\'es, afin de
parvenir au but en effectuant des calculs qui sont analogues, par
leurs id\'ees fondamentales, aux calculs helmholtziens.}

Dans le \S~96, nous empruntons la deuxi\`eme voie, et nous exposerons
s\'epar\'ement quelles conclusions la th\'eorie des groupes permet de
tirer \`a partir des axiomes de Helmhotz. On obtiendra alors comme
r\'esultat que ces axiomes (si on les interpr\`ete d'une mani\`ere
assez voisine) suffisent certes \`a caract\'eriser les mouvements
euclidiens et non-euclidiens de l'espace \`a trois dimensions, mais
qu'en tout cas, l'un des axiomes de Helmholtz, \`a savoir l'axiome de
monodromie, est superflu.

Le r\'esultat ainsi obtenu soul\`eve la question relativement proche:
savoir si l'on ne peut pas se passer aussi d'autres parties des
axiomes de Helmholtz. La r\'eponse \`a cette question sera donn\'ee
dans le Chapitre~23; on montrera que {\em les axiomes de Helmholtz
restants renferment aussi des \'el\'ements superflus}.

\HEAD{Critique des recherches helmholtziennes.}{
Division\,\,V.\,\,\,Chapitre\,\,21.\,\,\,\S\,\,91.}

\sectiondritterV{\sf\S\,\,91.
\\
Les axiomes helmholtziens.}
\label{S-91}
\setcounter{footnote}{0}

\renewcommand{\thefootnote}{\fnsymbol{footnote}}
Au fondement de sa recherche, Monsieur de Helmholtz place les axiomes
suivants, ou, hypoth\`eses, comme il les appelle\footnote[1]{\,
G\"ott. Nachr. 1868, p.~197 sq.}:
\renewcommand{\thefootnote}{\arabic{footnote}}

\medskip

{\small

<<\,{\bf I.}~L'espace \`a $n$ dimensions est une vari\'et\'e $n$ fois
\'etendue, c'est-\`a-dire que l'\'el\'ement individuel d\'efini en
celui-ci, le point, peut \^etre d\'etermin\'e par la mesure de
variables quantitatives continues et ind\'ependantes les unes des
autres (coordonn\'ees), dont le nombre est \'egal \`a $n$. Chaque
mouvement d'un point est ainsi accompagn\'e par une modification
continue d'au moins une des coordonn\'ees. S'il devait y avoir des
exceptions dans lesquelles, soit la modification deviendrait
discontinue, soit, malgr\'e le mouvement, absolument aucune
modification de toutes les coordonn\'ees ne se produirait, alors ces
exceptions devraient \^etre limit\'ees \`a certains lieux d\'efinis
par une ou plusieurs \'equations (donc \`a des points, \`a des lignes,
\`a des surfaces, \etc), et ces lieux peuvent initialement rester
exclus de la recherche.

<<\,On doit remarquer que par la continuit\'e du changement au cours
du mouvement, on ne signifie pas seulement que toutes les valeurs
interm\'ediaires jusqu'\`a la valeur finale sont parcourues par la
quantit\'e qui se diff\'erencie d'elle-m\^eme\footnote{\, Les
symétries par rapport à un hyperplan, telles que $( x_1, x_2, \dots,
x_n) \longmapsto ( -x_1, x_2, \dots, x_n)$, ne sont pas admises, car
elles ne sont pas reliées à la transformation identique par une
famille {\em continue} de mouvements puisqu'elles renversent
l'orientation. }, mais aussi que les quotients diff\'erentiels
existent, c'est-\`a-dire que le rapport des changements associ\'es des
coordonn\'ees se rapproche d'un rapport fixe, lorsque les grandeurs de
ces changements diminuent progressivement.

\medskip

<<\,{\bf II.}~On 
pr\'esuppose l'existence de corps mobiles mais rigides en
eux-m\^emes, c'est-\`a-dire de syst\`emes de points, comme cela est
n\'ecessaire pour pouvoir effectuer la comparaison des grandeurs
spatiales par congruence. Mais comme nous ne sommes pas encore
autoris\'es \`a pr\'esupposer une m\'ethode sp\'eciale pour la mesure
des grandeurs spatiales, la d\'efinition d'un corps rigide, ici, ne
peut \^etre que la suivante: {\sl entre les $2 n$ coordonn\'ees de
chaque paire de points qui appartient \`a un corps rigide en
lui-m\^eme, il existe une \'equation ind\'ependante du mouvement de
ce dernier, qui reste la m\^eme pour toutes les paires de points
congruentes}.

<<\,Sont {\sl congruentes}\, les paires de points qui peuvent \^etre
amen\'ees \`a co\"{\i}ncider simultan\'ement, ou l'une apr\`es
l'autre, avec la m\^eme paire de points de l'espace.

\medskip

\label{439}
<<\,{\bf III.}~{\em 
On pr\'esuppose une mobilit\'e parfaitement libre des
corps rigides}, c'est-\`a-dire, on suppose que chaque point en eux
peut \^etre d\'eplac\'e contin\^ument vers le lieu de chaque autre
point, tant que son mouvement n'est pas restreint par les \'equations
qui existent entre lui et les autres points restants du syst\`eme
rigide auquel il appartient.

<<\,Le premier point d'un syst\`eme rigide en lui-m\^eme est donc
absolument mobile. Lorsque ce point est fix\'e, il existe une
\'equation pour le deuxi\`eme point et l'une de ses coordonn\'ees
devient une fonction des $(n-1)$ coordonn\'ees restantes. Apr\`es que
le deuxi\`eme point est aussi fix\'e, il existe deux \'equations pour
le troisi\`eme point, \etc Par cons\'equent, ce sont au total
$\frac{ n ( n+1)}{ 1\, \cdot \, 2}$ quantit\'es qui sont requises pour
la d\'etermination de la position d'un corps rigide en lui-m\^eme.

<<\,Il d\'ecoule de cette supposition et de celle \'enonc\'ee en~II,
que, {\em si deux syst\`emes de points rigides en eux-m\^emes $A$ et
$B$ peuvent \^etre, dans une premi\`ere position initiale de $A$,
rapport\'es \`a une congruence de points correspondants, alors, pour
toute autre position de $A$, ils doivent pouvoir \^etre rapport\'es
\`a la congruence de tous les m\^emes points qui \'etaient congruents
auparavant}. Cela veut dire, en d'autres termes, que la congruence de
deux objets spatiaux ne d\'epend pas de leur situation initiale, ou
encore que toutes les parties de l'espace, lorsqu'on fait abstraction
de leurs bornes, sont congruentes les unes aux autres.

\medskip

<<\,{\bf IV.}~On doit finalement attribuer encore \`a l'espace une
propri\'et\'e qui est analogue \`a la {\sl monodromie}\, des fonctions
d'une variable complexe, et qui exprime en elle-m\^eme que deux corps
congruents sont encore \`a nouveau congruents apr\`es que l'un d'eux a
subi une r\'evolution compl\`ete autour d'un axe de rotation
arbitraire. Une {\sl rotation} est caract\'eris\'ee analytiquement en
ceci qu'un certain nombre de points du corps en mouvement conservent
leurs coordonn\'ees inchang\'ees au cours du mouvement; un {\sl
retour en arri\`ere}\, du mouvement est caract\'eris\'e par le fait
que les complexes de valeurs des coordonn\'ees qui se transformaient
auparavant contin\^ument l'un dans l'autre sont parcourus en sens
inverse. Nous pouvons donc exprimer le fait concern\'e comme suit:
{\em Lorsqu'un corps rigide tourne autour de $(n-1)$ de ses points et
que ces points sont choisis de telle sorte que sa position ne d\'epend
plus que d'une variable ind\'ependante, alors la rotation sans retour
en arri\`ere le reconduit finalement\`a la situation initiale dont il
est parti}.\,>>

}

\medskip

Les axiomes de Monsieur de Helmholtz ici reproduits attribuent
certaines propri\'et\'es aux mouvements de l'espace $n$ fois \'etendu,
et il s'agit essentiellement maintenant de d\'eterminer tous les
syst\`emes possibles de mouvements pour lesquels les propri\'et\'es
indiqu\'ees se manifestent. Afin de pouvoir employer notre th\'eorie
des groupes pour r\'esoudre ce probl\`eme, nous devons avant toute
chose montrer que nous avons au fond affaire ici \`a un probl\`eme de
la th\'eorie des groupes. C'est ce qui se va se passer dans le
prochain paragraphe.

\HEAD{Critique des recherches helmholtziennes.}{
Division\,\,V.\,\,\,Chapitre\,\,21.\,\,\,\S\,\,92.}

\sectiondritterV{\sf\S\,\,92.
\\
Cons\'equences des axiomes helmholtziens.}
\label{S-92}
\setcounter{footnote}{0}

Le {\em premier}\, des axiomes helmholtziens exprime simplement que
les mouvements continus sont possibles et il d\'etermine ce qui doit
\^etre entendu en g\'en\'eral par <<\,mouvement continu\,>>.

Lorsqu'on consid\`ere un mouvement de l'espace $n$ fois \'etendu, on
s'imagine de la mani\`ere la plus commode deux espaces $n$ fois
\'etendus contenus l'un dans l'autre, dont l'un est fixe, et l'autre
mobile\footnote{\, En dimension deux, l'espace fixe sera une grande
région rectangulaire et planaire, tandis que l'espace mobile,
superposé tel une nappe liquide, glissera de par toutes ses parties,
d'un seul tenant, avec des franges libres. }; la nature particuli\`ere
du mouvement consid\'er\'e sp\'ecifie alors de quelle mani\`ere les
points individuels de l'espace en mouvement modifient leur situation
\`a l'int\'erieur de l'espace fixe. Maintenant, l'Axiome~I demande que
chaque mouvement soit accompagn\'e d'une modification continue des
coordonn\'ees du point qui se d\'eplace. Par cons\'equent, si nous
nous imaginons un mouvement quelconque qui commence pour le temps $t =
0$ et qui est continuel pendant un certain temps $t$, et si nous
admettons qu'un point quelconque de l'espace mobile a les
coordonn\'ees $x_1, \dots, x_n$ au temps $t = 0$ par rapport \`a
l'espace fixe, et qu'il a les coordonn\'ees $\mathfrak{ x}_1, \dots,
\mathfrak{ x}_n$ au temps $t$, alors notre mouvement sera
repr\'esent\'e par des \'equations de la forme:
\def\theequation{1}\begin{equation}
\mathfrak{x}_\nu
=
F_\nu\big(
x_1,\dots,x_n,t
\big)
\ \ \ \ \ \ \ \ \ \ \
{\scriptstyle{(\nu\,=\,1\,\cdots\,n)}},
\end{equation}
qui pour $t = 0$ se r\'eduisent aux \'equations: $\mathfrak{ x}_\nu =
x_\nu$; avec cela, les $F_\nu$ sont des fonctions r\'eelles de leurs
arguments.

Ainsi, nous voyons que chaque mouvement continu de l'espace $n$ fois
\'etendu fournit une famille continue de $\infty^1$ transformations
ponctuelles r\'eelles, et pour pr\'eciser, une famille dans
laquelle est contenue la transformation identique. \label{440}

Pour ce qui touche \`a la nature des fonctions $F_\nu$, remarquons que
Monsieur de Helmholtz pr\'esuppose en tout cas l'existence des
quotients diff\'erentiels du premier ordre par rapport \`a $x$ et \`a
$t$; cela d\'ecoule d\'ej\`a de l'Axiome~I, bien que ce ne soit pas
\'enonc\'e avec toute la certitude d\'esirable, mais cela devient
indubitable lorsqu'on envisage les calculs situ\'es aux 
pages~202--206 de son travail. \`A cet
endroit-l\`a, Monsieur de Helmholtz utilise par ailleurs aussi
l'existence de certains {\em quotients diff\'erentiels du second
ordre}, lorsqu'il diff\'erentie notamment 
d'abord par rapport \`a $x$ et ensuite
par rapport au param\`etre pr\'esent, ce qui correspondrait au fait
que les quotients diff\'erentiels des $F_\nu$ par rapport \`a $x$ 
sont diff\'erentiables de leur c\^ot\'e par rapport \`a $t$.

On doit en outre observer que Monsieur de Helmholtz exclut de la
consid\'eration les circonstances o\`u la continuit\'e est
interrompue. Nous pouvons exprimer cela plus rigoureusement, en
\'etablissant, comme on le fait par ailleurs couramment, que l'on
doit se restreindre dans toute la recherche \`a une portion limit\'ee
de l'espace, \`a l'int\'erieur de laquelle toutes les fonctions qui se
pr\'esentent et leurs premiers quotients diff\'erentiels sont
continus; il s'agit l\`a de fixer les concepts d'une mani\`ere qui
est exactement la m\^eme que celle vers laquelle nous avons d\'ej\`a
converg\'e\footnote{\, %%%%%%%%%%%%%%%%%%%%%%%%%%%%%%%%%%%%% 
\`A ce sujet, \voir~les principes
de pensée p.~\pageref{Prologue} sq.
} %%%%%%%%%%%%%%%%%%%%%%%%%%%%%%%%%%%%%%%%%%%%%%%%%%%%%%%%%%%
depuis toujours dans nos recherches g\'en\'erales sur les
groupes continus.

\bigskip

Dans son {\em deuxi\`eme} axiome, Monsieur de Helmholtz caract\'erise
plus pr\'ecis\'ement les mouvements de l'espace $n$ fois \'etendu,
lorsqu'il indique comment deux points quelconques de l'espace mobile
mentionn\'e ci-dessus se comportent l'un
par rapport \`a l'autre au cours des
diff\'erents mouvements. Naturellement, les exigences de l'Axiome~II,
et \'egalement celles des axiomes suivants, se r\'ef\`erent seulement
aux points qui se trouvent \`a l'int\'erieur d'une telle portion
limit\'ee de l'espace $n$ fois \'etendu, et qui y demeurent aussi au
cours du mouvement.

Consid\'erons notre espace mobile dans une situation quelconque \`a
l'int\'erieur de l'espace fixe et envisageons deux points quelconques
de l'espace mobile, qui ont les coordonn\'ees: $x_1^0 \dots x_n^0$,
$y_1^0 \dots y_n^0$ par rapport \`a l'espace fixe; les quantit\'es:
$x_1^0 \dots x_n^0$, $y_1^0 \dots y_n^0$ poss\`edent alors des valeurs
num\'eriques d\'etermin\'ees, que nous pouvons cependant nous imaginer
comme \'etant choisies de mani\`ere absolument quelconque. En outre,
soit $x_1 \dots x_n$, $y_1 \dots y_n$ les coordonn\'ees des deux
points en question lorsqu'ils sont dans une autre position quelconque
de l'espace mobile. Le deuxi\`eme axiome helmholtzien demande alors
qu'entre les coordonn\'ees $x_1 \dots x_n$, $y_1 \dots y_n$, il existe
une \'equation qui est ind\'ependante de tous les mouvements, donc une
\'equation qui est aussi satisfaite, dans toute nouvelle position de
notre paire de points, par les coordonn\'ees de ses deux points. De
plus, dans l'axiome, il appara\^{\i}t explicitement qu'entre les $2n$
coordonn\'ees de toute paire de points qui appartient \`a un corps
rigide, une telle \'equation doit exister. En cela r\'eside le fait
que l'\'equation entre $x_1 \dots x_n$, $y_1 \dots y_n$ doit toujours
\^etre une \'equation v\'eritable, quelle que soit la mani\`ere dont
la position initiale: $x_1^0 \dots x_n^0$, $y_1^0 \dots y_n^0$ des
deux points peut \^etre choisie\footnote{\, On s'attend à ce que
l'équation mentionnée dépende des quatre variables $x$, $y$, $x^0$ et
$y^0$, comme c'est le cas pour l'équation $\vert \! \vert x - y
\vert\! \vert = \vert \! \vert x^0 - y^0 \vert \! \vert$ 
qui exprime la conservation de la norme euclidienne $\vert \! \vert z
\vert \! \vert :=
\big( z_1^2 + \cdots + z_n^2\big)^{ 1/ 2}$ à 
travers tous les déplacements. }.

Dans l'Axiome~II, il est demand\'e en plus que l'\'equation en
question soit la m\^eme pour toutes les paires de points congruentes,
donc pour toutes les paires de points qui peuvent \^etre envoy\'ees
l'une sur l'autre par les mouvements; par cons\'equent, elle ne
d\'epend, en dehors de $x_1 \dots x_n$, $y_1 \dots y_n$, que des
coordonn\'ees quelconques d'une paire de points congruents que l'on
peut choisir librement, par exemple de $x_1^0 \dots x_n^0$, $y_1^0
\dots y_n^0$, donc cette \'equation peut \^etre rapport\'ee \`a la
forme
\footnote{\, Plus bas, l'équation~\thetag{7} ramènera
encore plus précisément cette relation syncrétique à la forme
symétrique et résolue $\Omega ( x, y) = \Omega ( x^0, y^0)$. }:
\def\theequation{2}\begin{equation}
\Phi
\big(
x_1^0\dots x_n^0,\
y_1^0\dots y_n^0,\ \
x_1\dots x_n,\
y_1\dots y_n
\big)
=
0.
\end{equation}
En particulier, cette \'equation doit encore \^etre satisfaite lorsque
la paire de points $x_\nu$, $y_\nu$ co\"{\i}ncide avec la paire de
points $x_\nu^0$, $y_\nu^0$, et comme on peut donner aux quantit\'es
$x_\nu^0$, $y_\nu^0$ des valeurs num\'eriques quelconques, et que donc
par suite il ne peut pas exister de relation seulement entre: $x_1^0 \dots
x_n^0$, $y_1^0 \dots y_n^0$, il en d\'ecoule que l'\'equation~\thetag{
2} doit se r\'eduire \`a une identit\'e lorsqu'on fait la
susbstitution:
\[
x_\nu
=
x_\nu^0,
\ \ \ \ \ \ \
y_\nu
=
y_\nu^0
\ \ \ \ \ \ \ \ \ 
{\scriptstyle{(\nu\,=\,1\,\cdots\,n)}}.
\]
De l\`a, il suit en m\^eme temps que~\thetag{ 2} ne peut pas \^etre
libre de toutes les $2n$ quantit\'es $x_\nu^0$, $y_\nu^0$.

\`A pr\'esent, imaginons-nous que la paire de points $x_\nu$, $y_\nu$
est soumise \`a un mouvement continu qui est repr\'esent\'e par
l'\'equation~\thetag{ 1}; elle est ainsi envoy\'ee sur une nouvelle
paire de points, dont les points poss\`edent les coordonn\'ees:
\def\theequation{3}\begin{equation}
\left\{
\aligned
\mathfrak{x}_\nu
&
=
F_\nu
\big(
x_1\dots x_n,\,t
\big)
\\
\mathfrak{y}_\nu
&
=
F_\nu
\big(
y_1\dots y_n,\,t
\big)
\\
&
\ \ 
{\scriptstyle{(\nu\,=\,1\,\cdots\,n)}}
\endaligned\right.
\end{equation}
dans l'espace fixe. Puisque cette nouvelle paire de points doit
satisfaire la m\^eme \'equation que la paire de points $x_\nu$,
$y_\nu$, il en d\'ecoule que l'\'equation:
\def\theequation{2'}\begin{equation}
\Phi
\big(
x_1^0\dots x_n^0,\
y_1^0\dots y_n^0,\ \
\mathfrak{x}_1\dots \mathfrak{x}_n,\
\mathfrak{y}_1\dots \mathfrak{y}_n
\big)
=
0
\end{equation}
est valide, quelle que soit la valeur de $t$. Si nous effectuons la
substitution~\thetag{ 3} dans cette \'equation, nous obtenons une
\'equation entre $x_1 \dots x_n$, $y_1 \dots y_n$ et $t$ qui doit
\^etre satisfaite pour toutes les valeurs de $t$. Si nous nous
souvenons ici que $x_1 \dots x_n$, $y_1 \dots y_n$ sont li\'es par la
relation~\thetag{ 2} et que cette relation doit \^etre ind\'ependante
de tous les mouvement, alors nous réalisons que par la
substitution~\thetag{ 3}, l'\'equation~\thetag{ 2'} doit se
transformer en une \'equation qui est \'equivalente \`a~\thetag{ 2}.
Si tel n'\'etait pas le cas, alors par la substitution~\thetag{ 3},
l'\'equation~\thetag{ 2'} ne se transformerait pas en une \'equation
qui est une cons\'equence de~\thetag{ 2}, donc l'\'equation~\thetag{
2} entre les coordonn\'ees des paires de points mobiles ne serait
manifestement pas ind\'ependante du mouvement, mais au contraire, elle
changerait au cours du mouvement.

Puisque l'\'equation~\thetag{ 2} ne peut pas \^etre libre des $2n$
quantit\'es $x_\nu^0$, $y_\nu^0$, nous sommes autoris\'es \`a nous
imaginer qu'elle est r\'esolue par rapport \`a l'une de ces
quantit\'es, par exemple, par rapport \`a $y_n^0$:
\def\theequation{4}\begin{equation}
y_n^0
=
\varphi
\big(
x_1^0\dots x_n^0,\
y_1^0\dots y_{n-1}^0,\ \
x_1\dots x_n,\
y_1\dots y_n
\big).
\end{equation}
D'apr\`es ce qui a \'et\'e dit \`a l'instant, on obtient \`a pr\'esent
que l'\'equation:
\def\theequation{4'}\begin{equation}
y_n^0
=
\varphi
\big(
x_1^0\dots x_n^0,\
y_1^0\dots y_{n-1}^0,\ \
\mathfrak{x}_1\dots\mathfrak{x}_n,\
\mathfrak{y}_1\dots\mathfrak{y}_n
\big)
\end{equation}
se transforme en~\thetag{ 4} par la substitution~\thetag{ 3}, et
puisque nous pouvons attribuer aux quantit\'es $x_\nu^0$, $y_\nu^0$
toutes les valeurs num\'eriques quelconques, ceci doit valoir aussi
quelles que soient les valeurs que peuvent avoir $x_\nu^0$, $y_\nu^0$,
c'est-\`a-dire que l'\'equation:
\def\theequation{5}\begin{equation}
\left\{
\aligned
&
\ \ \ \ \
\varphi
\big(
x_1^0\dots x_n^0,\
y_1^0\dots y_{n-1}^0,\ \
\mathfrak{x}_1\dots\mathfrak{x}_n,\
\mathfrak{y}_1\dots\mathfrak{y}_n
\big)
=
\\
&
=
\varphi
\big(
x_1^0\dots x_n^0,\
y_1^0\dots y_{n-1}^0,\ \
x_1\dots x_n,\
y_1\dots y_n
\big)
\endaligned\right.
\end{equation}
doit se transformer en une identit\'e par la substitution~\thetag{ 3}.

Tout ceci reste encore vrai aussi, lorsque l'on attribue aux
quantit\'es: $x_1^0 \dots x_n^0$, $y_1^0 \dots y_{ n - 1}^0$ des valeurs
num\'eriques
d\'etermin\'ees: $\alpha_1 \dots \alpha_n$, $\beta_1 \dots \beta_{ n
- 1}$; si donc nous posons:
\[
\varphi
\big(
\alpha_1\dots \alpha_n,\
\beta_1\dots\beta_{n-1},\ \
x_1\dots x_n,\
y_1\dots y_n
\big)
=
\Omega
\big(
x_1\dots x_n,\
y_1\dots y_n
\big)
\]
nous obtenons \`a partir de~\thetag{ 5} l'\'equation:
\def\theequation{6}\begin{equation}
\Omega
\big(
\mathfrak{x}_1\dots\mathfrak{x}_n,\
\mathfrak{y}_1\dots\mathfrak{y}_n
\big)
=
\Omega
\big(
x_1\dots x_n,\
y_1\dots y_n
\big),
\end{equation}
qui se transforme \'egalement en une identit\'e par la
substitution~\thetag{ 3}. Nous trouvons ainsi que les deux points $x_1
\dots x_n$, $y_1 \dots y_n$ poss\`edent l'invariant: $\Omega
\big( x_1 \dots x_n, \ y_1 \dots y_n \big)$ relativement aux
$\infty^1$ transformations~\thetag{ 1}. En m\^eme temps, il en
d\'ecoule que l'\'equation~\thetag{ 2}, ou l'\'equation~\thetag{ 4}
qui lui est \'equivalente, peut \^etre rapport\'ee \`a la
forme\footnote{\, Cette dernière équation déduite, aussi bien que
l'équation abstraite~\thetag{ 7} postulée p.~\pageref{eq-7-p-10},
réfère sémantiquement à la conservation ${\rm dist} ( x, y) = {\rm
dist} ( x^0, y^0)$ d'une <<\,distance\,>> en un certain sens
généralisé; les résultats de la recherche le préciseront. }:
\def\theequation{7}\begin{equation}
\Omega
\big(
x_1\dots x_n,\
y_1\dots y_n
\big)
=
\Omega
\big(
x_1^0\dots x_n^0,\
y_1^0\dots y_n^0
\big).
\end{equation}

Les \'equations~\thetag{ 1} repr\'esentent un mouvement quelconque
parmi les mouvements continus qui sont admissibles, et ce, de la
mani\`ere la plus g\'en\'erale possible; les m\^emes consid\'erations
que celles \'emises auparavant montrent alors que la fonction $\Omega
\big( x_1 \dots x_n, \ y_1 \dots y_n \big)$ des coordonn\'ees de la
paire de points $x_1 \dots x_n, \ y_1 \dots y_n$ conserve sa valeur
num\'erique \`a travers tous les mouvements dont est susceptible cette
paire de points; et que donc, en prenant pour base le deuxi\`eme
axiome helmholtzien, chaque paire de points appartenant
\`a l'espace mobile mentionn\'e pr\'ec\'edemment poss\`ede un {\em
invariant}\, relativement \`a tous les mouvements. En d'autres termes:

\medskip

{\em Soit:
\[
\mathfrak{x}_\nu
=
F_\nu
\big(
x_1\dots x_n,\ t
\big)
\ \ \ \ \ \ \ \ \ \
{\scriptstyle{(\nu\,=\,1\cdots\,n)}}
\]
la famille de $\infty^1$ transformations qui est d\'etermin\'ee par un
mouvement continu quelconque de l'espace $n$ fois \'etendu, alors deux
points: $x_1 \dots x_n$, $y_1 \dots y_n$ poss\`edent toujours un
invariant:
\[
\Omega
\big(
x_1\dots x_n,\
y_1\dots y_n
\big)
\]
relativement \`a la famille~\thetag{ 1}; et pour pr\'eciser, les deux
points ont cet invariant relativement \`a toute famille de $\infty^1$
transformations, qui est sp\'ecifi\'ee comme \'etant l'un des
mouvements continus admissibles de l'espace le plus g\'en\'eral
possible.}

\medskip

Toutefois, on doit souligner que l'existence d'un tel invariant est
seulement une cons\'equence du deuxi\`eme axiome helmholtzien, {\em et
qu'en revanche, l'exigence qu'il doit exister un tel invariant ne peut
pas \^etre substitu\'ee compl\`etement \`a cet axiome}.

En effet, le deuxi\`eme axiome demande, comme nous l'avons vu plus
haut, qu'entre les coordonn\'ees de deux points distincts l'un de
l'autre dans notre espace mobile, il existe une \'equation
ind\'ependante du mouvement, et pour pr\'eciser, une v\'eritable
\'equation\footnote{\, 
En toute rigueur, Helmholtz a effectivement oublié de requérir que
l'équation invariante soit non triviale, sain réflexe logique dont nul
géomètre ne doit manquer. }, donc une \'equation qui n'est pas sans
signification pour une paire de points individuelle. Si maintenant
deux points: $x_1
\dots x_n$ et $y_1 \dots y_n$ de l'espace mobile ont l'invariant:
$\Omega ( x, y)$ relativement \`a tous les mouvements, alors il
s'ensuit\,\,---\,\,\`a vrai dire en toute
g\'en\'eralit\'e\,\,---\,\,qu'entre leurs $2 n$ coordonn\'ees, il
existe une \'equation ind\'ependante de tous les mouvements, \`a
savoir l'\'equation:
\def\theequation{7}\begin{equation}
\Omega
\big(
x_1\dots x_n;\
y_1\dots y_n
\big)
=
\Omega
\big(
x_1^0\dots x_n^0;\
y_1^0\dots y_n^0
\big),
\end{equation}
dans laquelle les deux syst\`emes de valeurs distincts l'un de
l'autre: $x_1^0 \dots x_n^0$ et: $y_1^0 \dots y_n^0$, d\'esignent
les coordonn\'ees de la paire de points: $x_1 \dots x_n$, $y_1 \dots
y_n$ dans une situation initiale quelconque bien d\'efinie de l'espace
mobile; dans certains cas exceptionnels, il peut cependant se produire
que~\thetag{ 7} ne repr\'esente pas d'\'equation v\'eritable entre
$x_1 \dots x_n$, $y_1 \dots y_n$, par exemple lorsque le membre de
droite de~\thetag{ 7} prend une forme absolument nulle pour certains
syst\`emes de valeurs: $x_1^0 \dots x_n^0$, $y_1^0 \dots y_n^0$.

Il s'ensuit de l\`a qu'afin d'\'epuiser compl\`etement le contenu du
deuxi\`eme axiome helmholtzien, nous ne pouvons pas nous contenter,
pour toute paire de points, de requ\'erir l'existence d'un invariant,
mais nous devons encore ajouter une deuxi\`eme exigence, qui peut
s'\'enoncer de la mani\`ere suivante:

{\em
L'invariant $\Omega \big( x_1 \dots x_n, \ y_1 \dots y_n \big)$ de la
paire de points: $x_1 \dots x_n$, $y_1 \dots y_n$ doit \^etre
constitu\'e de telle sorte que l'\'equation~\thetag{ 7} soit toujours
une \'equation v\'eritable entre $x_1 \dots x_n$ et $y_1 \dots y_n$,
quelles que soient les diff\'erentes valeurs num\'eriques qu'on puisse
donner \`a $x_1^0 \dots x_n^0$, $y_1^0 \dots y_n^0$.}

Il va de soi que cette exigence doit \^etre satisfaite seulement pour
les syst\`emes de valeurs: $x_1^0 \dots x_n^0$, $y_1^0 \dots y_n^0$
qui appartiennent \`a la r\'egion limit\'ee de l'espace $n$ fois \'etendu
qui a \'et\'e mentionn\'eee pr\'ec\'edemment.

\bigskip

Avant que nous passions \`a la consid\'eration du 
troisi\`eme axiome helmholtzien, nous voulons auparavant
mentionner certaines 
cons\'equences qui se laissent d\'eduire de l'existence
de l'invariant $\Omega ( x, y)$.

Supponsons que les \'equations: 
\def\theequation{8}\begin{equation}
x_\nu'
=
F_\nu
\big(
x_1\dots x_n;\ t
\big)
\ \ \ \ \ \ \ \ \ \ \ \
{\scriptstyle{(\nu\,=\,1\,\cdots\,n)}}
\end{equation}
et
\def\theequation{9}\begin{equation}
x_\nu''
=
\Phi_\nu
\big(
x_1'\dots x_n';\ \tau
\big)
\ \ \ \ \ \ \ \ \ \ \ \
{\scriptstyle{(\nu\,=\,1\,\cdots\,n)}}
\end{equation}
repr\'esentent deux mouvements continus quelconques de l'espace
$n$ fois \'etendu. Si nous nous imaginons d'abord le
mouvement~\thetag{ 8} s'accomplissant pendant le temps $t$, et ensuite
le mouvement~\thetag{ 9} pendant le temps $\tau$, alors le point:
$x_1 \dots x_n$ se trouvera finalement dans la position: $x_1 ''
\dots x_n ''$ qui est d\'efinie par les \'equations:
\def\theequation{10}\begin{equation}
x_\nu''
=
\Phi_\nu
\big(
F_1(x,t)\dots F_n(x,t);\ 
\tau
\big)
\ \ \ \ \ \ \ \ \ \ \ \
{\scriptstyle{(\nu\,=\,1\,\cdots\,n)}}.
\end{equation}
Maintenant, comme ces \'equations~\thetag{ 10} repr\'esentent \'evidemment
une transformation, on en d\'eduit que par l'effectuation de deux
mouvements l'un \`a la suite de l'autre, on obtient toujours une
certaine transformation de l'espace. La m\^eme chose vaut
naturellement lorsqu'on effectue un nombre quelconque de tels
mouvements l'un \`a la suite de l'autre.

Si nous nous souvenons maintenant que les deux points: $x_1 \dots
x_n$, $y_1 \dots y_n$ ont l'invariant $\Omega ( x, y)$ non seulement
relativement \`a la transformation~\thetag{ 8}, mais aussi
relativement \`a la transformation~\thetag{ 9}, nous voyons alors
imm\'ediatement qu'ils poss\`edent cet invariant relativement \`a
toute transformation~\thetag{ 10}, et plus g\'en\'eralement,
relativement \`a toute transformation qui est obtenue par
accomplissement de plusieurs mouvements l'un \`a la suite de l'autre.
Par cons\'equent:

{\em
Si l'on effectue un mouvement continu, ou plusieurs mouvements
de cette sorte
l'un \`a la suite de l'autre, alors on obtient toujours une
transformation relativement \`a laquelle deux points quelconques:
$x_\nu$ et $y_\nu$ ont l'invariant $\Omega ( x, y)$.}

Au moyen du mouvement continu le plus g\'en\'eral possible qui soit
admissible, une certaine famille de transformations de l'espace est
donc sp\'ecifi\'ee de telle sorte que, relativement \`a cette
famille, deux points $x_\nu$ et $y_\nu$ ont l'invariant $\Omega ( x,
y)$. Nous ne pouvons pas dire pour l'instant de quelle nature
particuli\`ere est faite cette famille, puisque c'est le
troisi\`eme axiome helmholtzien qui donne la premi\`ere information
l\`a-dessus. On ne peut d'ailleurs m\^eme pas encore conclure du
contenu du deuxi\`eme aximome que deux points doivent avoir {\em
seulement}\, un invariant.

\bigskip

Le {\em troisi\`eme}\, axiome de Monsieur de Helmholtz ({\em voir}\,
p.~\pageref{439}, lignes~1 \`a 11 ci-dessus) est constitu\'e de deux
parties. \label{446}

La premi\`ere partie (l.~1 \`a 5) peut \^etre, si l'on tient compte de
ce qui a \'et\'e dit auparavant, formul\'ee de la mani\`ere suivante:
chaque point de l'espace mobile doit pouvoir \^etre d\'eplac\'e
contin\^ument dans le lieu de chaque autre point de cet espace, tant
qu'il n'est pas restreint par le fait que les invariants de toutes les
paires de points de l'espace mobile auquel il appartient doivent
conserver leur valeur num\'erique au cours du mouvement.

En cela r\'eside le fait que chaque invariant que poss\`ede un
syst\`eme quelconque de points: $P_1, P_2, P_3 \dots$ relativement
\`a tous les mouvements, doit se laisser exprimer au moyen des
invariants des paires de points qui sont contenues dans le syst\`eme.
Si en effet $J$ est un invariant quelconque que le syst\`eme de
points: $P_1, P_2, P_3 \dots$ a relativement \`a tous les mouvements,
alors la fonction $J$ conserve sa valeur num\'erique au cours de tous
les mouvements. Si maintenant $J$ ne se laissait pas exprimer au moyen
des invariants des paires de points: $P_1, P_2$; $P_1, P_3$; $P_2,
P_3$; $\dots$, alors la mobilit\'e du syst\`eme de points: $P_1,
P_2, P_3 \dots$ ne serait pas seulement limit\'ee par le fait que les
invariants de ces paires de points-l\`a devraient conserver leur valeur
num\'erique au cours du mouvement, mais elle serait aussi limit\'ee
par la condition ind\'ependante par rapport \`a celles-l\`a que $J$
doit conserver toujours sa valeur num\'erique; mais cela
contredirait l'axiome pos\'e ci-dessus.

Il suit alors de la premi\`ere partie du troisi\`eme axiome
helmholtzien qu'un point n'a absolument aucun invariant relativement
\`a tous les mouvements et que trois points (ou plus) ne poss\`edent,
relativement \`a tous les mouvements, que les invariants qui
s'expriment au moyen des invariants des paires de points qui sont
contenues en eux; autrement dit: {\em un point individuel n'a en
g\'en\'eral aucun invariant; deux points ont en tout cas un invariant
relativement \`a tous les mouvements possibles; mais au contraire,
trois points (ou plus) n'ont aucun invariant {\small\sf essentiel}}. Mais
comme nous ne savons cependant pas encore si cette propri\'et\'e des
mouvements peut \^etre substitu\'ee enti\`erement \`a l'exigence qui
est pos\'ee dans cette premi\`ere partie-l\`a, nous devons comme
auparavant ajouter l'exigence que {\em chaque point $P$ de l'espace
mobile est limit\'e dans son mouvement seulement par les invariants
qu'il a avec les autres points de l'espace mobile}.

\bigskip

Il reste maintenant encore \`a discuter de la deuxi\`eme partie du
troisi\`eme axiome helmholtzien ({\em voir} p.~\pageref{439},
l.~6--11 ci-dessus).

La deuxi\`eme partie de l'Axiome~III semble au premier coup d'{\oe}il
ne contenir que des cons\'equences de la premi\`ere partie, et c'est
visiblement ainsi que Monsieur de Helmholtz a voulu l'entendre.
Il en va cependant tout autrement. Certes les lignes 6--9 expriment
seulement des faits qui sont cons\'equences imm\'ediates des
hypoth\`eses pos\'ees pr\'ec\'edemment, en particulier les mots:
<<\,le premier point \dots \ mobile\,>> ne sont qu'une autre mani\`ere
d'exprimer le fait qu'un point individuel ne poss\`ede aucun
invariant. Mais dans les lignes 10 et 11 se glisse une nouvelle
supposition qui n'est pas cons\'equence des pr\'ec\'edentes.

Monsieur de Helmholtz s'imagine en effet qu'un certain nombre, disons
$m$, de points: $P_1, P_2, \dots P_m$ de l'espace mobile sont
fix\'es. D'apr\`es l'Axiome~II, il existe alors pour les coordonn\'ees
de chaque autre point $P$ de l'espace mobile certaines \'equations,
qui expriment que les invariants des $m$ paires de points: $P, P_1$;
$P, P_2$; $\dots$; $P, P_m$ conservent leur valeur num\'erique au
cours de tous les mouvements encore possibles. Mais les hypoth\`eses qui
pr\'ec\`edent ne disent rien sur la nature et sur le nombre de ces
\'equations, ce au sujet de quoi certaines suppositions qui sont
faites tacitement aux lignes~11 et~12 donnent un premier
\'eclaircissement. \label{447}

\renewcommand{\thefootnote}{\fnsymbol{footnote}}
Monsieur de Helmholtz demande dans ces circonstances la chose
suivante, que nous pouvons exprimer ainsi: si les $m$ points $P_1,
P_2, \dots P_m$ sont fix\'es, il doit exister entre les $n$
coordonn\'ees de chaque autre point $P$ exactement $m$ (mais pas plus)
\'equations, et {\em ces \'equations doivent g\'en\'eralement \^etre
ind\'ependantes les unes des autres}, c'est-\`a-dire, qu'elles doivent
\^etre ind\'ependantes les unes des autres aussi longtemps que $P_1
\dots P_m$ sont des points mutuellement en position
g\'en\'erale\footnote[1]{\, Ainsi doivent \^etre compris les mots des
lignes 6--11 de l'Axiome~III ci-dessus. 
Monsieur de Helmholtz dit ici, certes
d'une mani\`ere qui n'est pas explicite, que lorsque le premier point
d'un syst\`eme rigide en lui-m\^eme est fix\'e, il doit exister une et
une seule \'equation pour chaque autre point; et puisqu'il ajoute
alors: <<\,l'une de ses coordonn\'ees devient une fonction des
$(n-1)$ coordonn\'ees restantes\,>>, il est clair qu'il exclut
l'existence de deux \'equations pour le deuxi\`eme point. Il r\'esulte
\'egalement de mani\`ere distincte de tout cela qu'apr\`es fixation
de $m$ points, alors pour tout autre point, il doit exister exactement
$m$ \'equations qui sont en g\'en\'eral ind\'ependantes les unes des
autres.}.
\renewcommand{\thefootnote}{\arabic{footnote}}

Il suit de l\`a tout d'abord qu'apr\`es fixation d'un point $P_1$, il
existe une et une seule \'equation pour tout autre point $P$; et
comme un point individuel ne doit avoir aucun invariant, il en
d\'ecoule que {\em deux points ont un et un seul invariant}.

En outre, les exigences helmholtziennes montrent qu'apr\`es fixation
de $m$ points: $P_1 \dots P_m$ qui sont mutuellement en position
g\'en\'erale, un mouvement continu est encore toujours possible, tant
que $m$ a l'une des valeurs $1, 2, \dots, n-1$, mais qu'au contraire,
dans le cas $m = n$, aucun mouvement continu n'est plus possible, et
qu'apr\`es fixation de $n$ points tels, tous les points de l'espace
demeurent plut\^ot g\'en\'eralement au repos. Si en effet les $m$
points: $P_1 \dots P_m$ sont fix\'es, il existe alors entre les $n$
coordonn\'ees de tout autre point $P$ en position g\'en\'erale $m$
\'equations ind\'ependantes les unes des autres. Et puisque, d'apr\`es
ce qui pr\'ec\`ede, lorsqu'on tient compte des hypoth\`eses admises,
ce sont les seules conditions auxquelles la mobilit\'e de $P$ est
soumise, il en d\'ecoule que $P$ peut encore \^etre envoy\'e sur
chaque autre point $P'$ dont les coordonn\'ees satisfont ces $m$
\'equations-l\`a, et qui est li\'e \`a $P$ par une s\'erie continue de
tels points. Si donc en particulier $m$ a la valeur $n$, alors le
point $P$ ne peut plus effectuer aucun mouvement continu, mais il doit
au contraire rester au repos.

\bigskip

Nous avons vu plus haut qu'en accomplissant un nombre quelconque de
mouvements continus l'un apr\`es l'autre, on obtient toujours une
transformation ponctuelle de $R_3$ parfaitement d\'etermin\'ee.
Nous pouvons maintenant donner des explications plus pr\'ecises
sur la famille des transformations qui sont engendr\'ees de cette
mani\`ere.

Consid\'erons notre espace mobile dans une situation quelconque \`a
l'int\'erieur de l'espace fixe. Soit $n$ points $P_1 \dots P_n$ qui
sont mutuellement en position g\'en\'erale dans l'espace mobile et
$\mathfrak{ P}_1 \dots \mathfrak{ P}_n$ les points de l'espace fixe
avec lesquels ils co\"{\i}ncident exactement. Si ensuite $P$ est un
autre point quelconque de l'espace mobile, alors $P$ a aussi une
position enti\`erement d\'etermin\'ee \`a l'int\'erieur de l'espace
fixe, car il n'existe plus aucun mouvement continu tel que $P_1 \dots
P_n$ conservent compl\`etement leur position: $\mathfrak{ P}_1 \dots
\mathfrak{ P}_n$.

Maintenant, imaginons-nous que l'espace mobile est soumis \`a un
nombre quelconque de mouvements continus. Nous allons voir que la
transformation la plus g\'en\'erale de $R_n$ que l'on peut obtenir de
cette mani\`ere ne d\'epend que d'un nombre fini de param\`etres
arbitraires.

En effet, par un mouvement continu, nous pouvons tout d'abord
transf\'erer le point $P_1$ de la position $\mathfrak{ P}_1$ vers tout
autre point $\mathfrak{ P}_1'$, et par cons\'equent, la position la
plus g\'en\'erale de $\mathfrak{ P}_1'$ d\'epend de $n$ param\`etres
arbitraires. Si l'on choisit $\mathfrak{ P}_1'$ fix\'e, alors $P_2$
peut encore \^etre transf\'er\'e vers tout autre point $\mathfrak{
P}_2 '$ qui satisfait une certaine \'equation, donc $\mathfrak{ P}_2'$
d\'epend encore, lorsqu'on a choisi $\mathfrak{ P}_1'$, de $n-1$
param\`etres, et ainsi de suite. En bref, par des mouvements continus,
nous pouvons parvenir \`a ce que $P_1 \dots P_n$ re\c coivent les
position nouvelles $\mathfrak{ P }_1 ' \dots \mathfrak{ P }_n'$, o\`u
le syst\`eme de points: $\mathfrak{ P }_1 ' \dots \mathfrak{ P }_n'$
d\'epend de:
\[
n+(n-1)+(n-2)+\cdots+1
=
{\textstyle\frac{n(n+1)}{1\,\cdot\,2}}
\]
param\`etres. Mais avec cela, toutes les possibilit\'es sont aussi
\'epuis\'ees, car aussit\^ot que $P_1 \dots P_n$ re\c
coivent la nouvelle position: $\mathfrak{ P}_1 ' \dots \mathfrak{
P}_n'$, tout autre point $P$ de l'espace mobile re\c coit en m\^eme
temps une nouvelle position $\mathfrak{ P}'$ compl\`etement
d\'etermin\'ee, puisqu'il n'y a plus aucun mouvement continu au cours
duquel $\mathfrak{ P}_1 ' \dots \mathfrak{ P}_n'$ demeurent
enti\`erement au repos.

%%%\Fill 
%%%[[Frames.]]

En cela r\'eside le fait que, par le transfert du syst\`eme de points
$P_1 \dots P_n$ depuis la position initiale: $\mathfrak{ P}_1 \dots
\mathfrak{ P}_n$ vers la nouvelle position: $\mathfrak{ P}_1 ' \dots
\mathfrak{ P}_n'$, une transformation ponctuelle enti\`erement
d\'etermin\'ee est d\'efinie de mani\`ere unique; quelle que soit la
fa\c con dont ce transfert puisse \^etre effectu\'e par une s\'erie de
mouvement continus se succ\'edant l'un apr\`es l'autre, on obtient
toujours la m\^eme transformation ponctuelle quand on effectue ces
mouvements l'un apr\`es l'autre. Maintenant, puisqu'il y a exactement
$\frac{ 1}{ 2}\, n ( n+1)$ param\`etres arbitraires \`a disposition
pour le choix du syst\`eme de points $\mathfrak{ P}_1 ' \dots
\mathfrak{ P}_n'$, lorsqu'on s'imagine le syst\`eme de points $P_1
\dots P_n$ d\'eplac\'e de toutes les mani\`eres possibles, il en
d\'ecoule que l'ensemble de toutes les transformations ponctuelles,
qui peuvent \^etre obtenues par un nombre quelconque de mouvements
continus de l'espace mobile \`a partir d'une position initiale
d\'etermin\'ee, constitue une famille ayant exactement $\frac{ 1}{
2}\, n ( n+1)$ param\`etres essentiels.
%%%\Mathematiques
%%%[[Est-ce que cela d\'epend du frame?]]
%%%\Mathematiques 
%%%[[$\{ g : d ( gx, gy) =  d( x, y) \}$ est \'evidemment un groupe]]

Remarquons ensuite que cette famille de transformations ponctuelles
est ind\'ependante du choix de la position initiale de l'espace
mobile, parce que, quelle que soit la mani\`ere dont on choisit la
position initiale, il y a toujours $n$ points $\overline{ P}_1 \dots
\overline{ P}_n$ de l'espace mobile qui co\"{\i}ncident avec les
points: $\mathfrak{ P}_1 \dots \mathfrak{ P}_n$ de l'espace fixe.
{\em En cela r\'eside le fait que la famille de transformations
ponctuelles ainsi d\'efinie \`a l'instant constitue un} 
{\small\sf groupe}.
En effet, si au moyen d'une transformation quelconque de notre
famille, nous amenons l'espace mobile de la position $R$ \`a la
position $R'$, et si ensuite, au moyen d'une autre transformation de
notre famille, nous l'amenons de la position $R'$ \`a la position
$R''$, alors il y a toujours une transformation
enti\`erement d\'etermin\'ee appartenant \`a la famille, gr\^ace \`a
laquelle $R$ se transforme en $R''$.

%%%\Fill 
%%%[[Gr\^ace au frame.]]

Le groupe que l'on trouve de cette mani\`ere est s\^urement
transitif, car il peut envoyer tout point de l'espace sur tout
autre point. De plus il est facile de voir que ses transformations
sont ordonn\'ees ensemble par paires de transformations
inverses. En effet, si $S$ est la transformation de notre groupe qui
envoie $\mathfrak{ P}_1 \dots \mathfrak{ P}_n$ sur $\mathfrak{ P}_1
' \dots \mathfrak{ P}_n'$, alors il est toujours possible, au moyen
d'un certain nombre de mouvements continus, de parvenir \`a ce que les
points de l'espace mobile, qui co\"{\i}ncident avec $\mathfrak{ P}_1'
\dots \mathfrak{ P}_n'$ dans une situation quelconque de cet espace,
arrivent finalement \`a la position: $\mathfrak{ P}_1 \dots
\mathfrak{ P}_n$. Ainsi $S^{ -1}$ fait toujours partie de notre
groupe en m\^eme temps que $S$.

%%%\Fill
%%%[[Commentaire sur $S^{ - 1}$.]]

\label{449}
Enfin, on peut aussi d\'emontrer que notre groupe, que nous voulons
appeler bri\`evement $g$, est continu. En effet, s'il ne l'\'etait
pas, d'apr\`es le Tome~I, 
Chap.~18\footnote{\,
%%%%%%%%%%%%%%%%%%%%%%%-------DEBUT--------%%%%%%%%%%%%%%%%%%%%%%%%%%%
Ce chapitre est consacré à l'étude des groupes de
transformations comportant plusiseurs composantes
connexes. 
}, %%%%%%%%%%%%%%%%%%%%%%%%-----FIN-----%%%%%%%%%%%%%%%%%%%%%%%%%%%%%%%
%%%\Fill
il serait compos\'e d'une s\'erie
de familles continues s\'epar\'ees, dont chacune contiendrait $\frac{
1}{ 2} \, n ( n+1)$ param\`etres et parmi ces familles, il y en aurait
une qui est continue et qui constituerait un groupe $\mathfrak{ g}$
\`a $\frac{ 1}{ 2} \, n ( n+1)$ param\`etres contenant des
transformations inverses par paires. Ce groupe $\mathfrak{ g}$ serait
le seul auquel appartient la transformation identique, parmi les
familles continues concern\'ees. Mais maintenant, chaque famille
continue de $\infty^1$ transformations qui est d\'etermin\'ee par un
mouvement continu contient la transformation identique ({\em voir}\,
p.~\pageref{440}) et elle fait donc partie du groupe $\mathfrak{ g}$;
par cons\'equent, le groupe $\mathfrak{ g}$ comprend aussi toutes les
transformations qui proviennent de l'accomplissement d'un nombre
quelconque de mouvements continus l'un apr\`es l'autre. Il en
d\'ecoule que $g$ est contenu dans $\mathfrak{ g}$ et comme
$\mathfrak{ g}$ \'etait un sous-groupe de $g$, nous pouvons conclure
que $g$ co\"{\i}ncide avec $\mathfrak{ g}$, c'est-\`a-dire que $g$ est
effectivement continu.

Avec les hypoth\`eses qu'a faites de Monsieur de Helmholtz,
l'\'enonc\'e suivant est donc valide:

{\em Lorsqu'on effectue un nombre quelconque de mouvements continus de
l'espace les uns \`a la suite des autres, parmi ceux qui sont les plus
admissibles possibles, on obtient un groupe continu
fini transitif $g$ de
transformations r\'eelles qui renferme exactement $\frac{ 1}{ 2} \, n
( n+1)$ param\`etres et dont les transformations sont inverses l'une
de l'autre par paires.}

%%%\Fill 
%%%[[diff\'erence conceptuelle entre groupe et mouvement.]]

Il va alors de soi que relativement \`a ce groupe, deux points ont un
et un seul invariant, et que $s > 2$ points n'ont pas d'invariant
essentiel. En effet, si deux points poss\'edaient par exemple plus
d'un invariant relativement au groupe, ils poss\`ederaient alors
\'evidemment aussi cet invariant relativement \`a tous les mouvements
continus, alors qu'ils ne doivent avoir qu'un seul invariant
relativement \`a ces mouvements.

\bigskip

Gr\^ace \`a ce qui pr\'ec\`ede, il a \'et\'e d\'emontr\'e que nous
avons affaire \`a un groupe continu fini et que relativement \`a ce
groupe, deux points ont un et un seul invariant, tandis que $s > 2$
points n'ont pas d'invariant essentiel. Nous avons d\'ej\`a \'etudi\'e
dans le Chapitre~20 les groupes de cette esp\`ece, et bien que nous
nous soyons limit\'es, \`a ce moment-l\`a, \`a l'espace de dimension
trois, on peut n\'eanmoins
appliquer imm\'ediatement au cas de l'espace $n$ fois
\'etendu au moins une partie des raisonnements conduits \`a cet
endroit-l\`a, \`a savoir les d\'eveloppements des pages
\pageref{399}--\pageref{404}. Nous reconnaissons ainsi que chaque
groupe ayant la constitution d\'efinie \`a l'instant poss\`ede les
propri\'et\'es suivantes: si un point en position g\'en\'erale est
fix\'e, alors tout autre point en position g\'en\'erale peut recevoir
encore $\infty^{ n- 1}$ positions diff\'erentes; si on fixe deux
points qui sont mutuellement en position g\'en\'erale, alors tout
autre troisi\`eme point en position g\'en\'erale peut encore recevoir
$\infty^{ n - 2}$ positions diff\'erentes, \etc; pour \^etre
bref, nous trouvons que la finitude du groupe et que les hypoth\`eses
faites sur les invariants de deux points (ou plus) qui ont \'et\'e
indiqu\'ees \`a la page~\pageref{447} se d\'eduisent des exigences de
Monsieur de Helmholtz, \`a savoir des exigences qu'apr\`es fixation
d'un point, il existe une et une seule \'equation pour tout autre
point, et ainsi de suite. Ainsi, nous n'avons plus besoin de
mentionner sp\'ecialement ces exigences.

\bigskip

Finalement, on peut encore mentionner que le groupe continu fini $g$
qui est d\'etermin\'e par les mouvements continus de l'espace se trouve
\^etre li\'e encore d'une mani\`ere simple \`a un autre groupe. Nous
entendons le groupe $g'$ de toutes les transformations ponctuelles
pour lesquelles deux points quelconques $x_\nu$, $y_\nu$ ont
l'invariant $\Omega ( x, y)$.

Il est \`a l'avance clair qu'il existe un tel groupe $g'$, car, quelle
que soit la mani\`ere dont nous pouvons choisir la fonction $\Omega (x,
y)$, il y a toujours des transformations relativement auxquelles deux
points ont l'invariant $\Omega ( x, y)$, et aussi, l'ensemble de ces
transformations constitue toujours un certain groupe, qui \`a vrai dire
peut \'eventuellement se r\'eduire \`a la transformation identique.

Dans notre cas, il est maintenant facile de voir que la transformation
la plus g\'en\'erale de $g'$ contient exactement $\frac{ 1}{ 2} \, n (
n+1)$ param\`etres arbitraires, ce qui d\'ecoule tout simplement de la
propri\'et\'e particuli\`ere qu'a l'invariant $\Omega ( x, y)$, sous les
hypoth\`eses ici pos\'ees. Par cons\'equent, il est clair que $g$ est le
plus grand groupe continu contenu dans $g'$; si $g'$ devait \^etre
lui-m\^eme continu, alors il co\"{\i}nciderait naturellement avec $g$.

\HEAD{Critique des recherches helmholtziennes.}{
Division\,\,V.\,\,\,Chapitre\,\,21.\,\,\,\S\,\,93.}

\sectiondritterV{\sf\S\,\,93.
\\
Formulation des axiomes helmholtziens
\\
en termes de th\'eorie des groupes.}
\label{S-93}
\setcounter{footnote}{0}

Nous allons maintenant rassembler les r\'esultats des pr\'ec\'edents
paragraphes. Afin de pouvoir nous exprimer d'une mani\`ere plus
commode, nous voulons toutefois appeler bri\`evement le groupe d\'efini \`a
la page~\pageref{449} un {\sl groupe des mouvements de $R_n$}. En
conformit\'e avec cela, nous appelons bri\`evement chaque transformation
de ce groupe un {\sl mouvement}, de telle sorte que par
<<\,mouvement\,>> nous entendons toujours une transformation qui
transforme l'espace mobile d'une situation vers une autre situation.
Ce que nous avons appel\'e jusqu'\`a pr\'esent un {\sl mouvement continu}
est maintenant 
simplement une famille continue de $\infty^1$ mouvements, dans
laquelle est contenue la transformation identique.

Cette mani\`ere de s'exprimer correspond manifestement \`a la
mani\`ere maintenant habituelle de s'exprimer dont nous avons
d\'ej\`a us\'e auparavant, lorsque nous parlions
du groupe des mouvements euclidiens, sachant aussi
que nous entendions chaque fois par <<\,mouvement\,>>
une transformation de ce groupe.

\`A pr\'esent, en collectant les r\'esultats des 
paragraphes pr\'ec\'edents, nous pouvons dire que les
trois premiers axiomes helmholtziens ont la 
m\^eme signification que les exigences suivantes:

\medskip

{\bf A)} {\em Chaque point de l'espace $n$ fois \'etendu peut \^etre
d\'etermin\'e par $n$ coordonn\'ees: $x_1 \dots x_n$.}

\medskip

{\bf B)} 
{\em Un groupe continu fini qui est r\'eel
et transitif et que nous appelons {\em groupe des mouvements:} 
\def\theequation{11}\begin{equation}
x_\nu'
=
f_\nu
\big(
x_1\dots x_n;\,
a_1\dots a_r
\big),
\quad\quad\quad
{\scriptstyle{(\nu\,=\,1\,\cdots\,n)}}
\end{equation}
est d\'etermin\'e par tous les mouvements possibles dans l'espace. Le
nombre $r$ de param\`etres de ce groupe poss\`ede la valeur $\frac{ 1}{
2}\, n ( n+1)$ et les transformations du groupe sont inverses l'une de
l'autre par paires. Les fonctions $f_1 \dots f_n$ sont non seulement
diff\'erentiables par rapport 
aux $x$, mais aussi par rapport aux $a$ et les
quotients diff\'erentiels par rapport aux $x$ sont de leur c\^ot\'e \`a
nouveau diff\'erentiables par rapport aux $a$.}

\medskip

{\bf C)} 
{\em Relativement au groupe~\thetag{ 11}, deux points ont un et un
\label{452}
seul invariant, et $s > 2$ points n'ont aucun invariant essentiel.
Si: $J \big(x_1 \dots x_n; \, y_1 \dots y_n \big)$ est l'invariant
des deux points $x_\nu$ et $y_\nu$, alors on doit pouvoir d\'elimiter \`a
l'int\'erieur de l'espace $n$ fois \'etendu une certaine r\'egion $n$ fois
\'etendue de telle sorte que les relations:
\def\theequation{12}\begin{equation}
J
\big(
x_1\dots x_n;\,y_1\dots y_n
\big)
=
J
\big(
x_1^0\dots x_n^0;\,y_1^0\dots y_n^0
\big)
\end{equation}
fournissent toujours une \'equation v\'eritable entre: $x_1 \dots x_n$,
$y_1 \dots y_n$, quels que soient les points distincts l'un de l'autre
que l'on puisse aussi choisir
pour $x_1^0 \dots x_n^0$, $y_1^0 \dots y_n^0$ dans la r\'egion.}

\medskip

{\bf D)} 
{\em \`A l'int\'erieur de la r\'egion d\'efinie \`a l'instant, tous les
points sont parfaitement et librement mobiles par le groupe~\thetag{
11}, tant qu'ils ne sont pas li\'es par les invariants que les paires
individuelles de points ont relativement au groupe. Si donc par
exemple: $x_1^0 \dots x_n^0$; $y_1^0 \dots y_n^0$ sont deux points
distincts quelconques de cette r\'egion-l\`a, alors, aussit\^ot que le
premier d'entre eux est fix\'e, le second peut occuper encore toutes les
positions: $y_1 \dots y_n$ qui satisfont les \'equations:
\def\theequation{13}\begin{equation}
J
\big(
x_1^0\dots x_n^0;\,y_1\dots y_n
\big)
=
J
\big(
x_1^0\dots x_n^0;\,y_1^0\dots y_n^0
\big),
\end{equation}
o\`u il est suppos\'e que depuis $y_1^0 \dots y_n^0$ vers $y_1 \dots y_n$,
une transition continue qui ne comporte que des syst\`emes de valeurs
r\'eelles satisfaisant~\thetag{ 13}, est possible.}

\`A cela s'ajoute maintenant encore le quatri\`eme axiome helmholtzien,
que nous n'avons encore pas pris en consid\'eration jusqu'\`a pr\'esent,
celui qu'on appelle souvent l'{\em Axiome de monodromie}. Nous
pouvons maintenant lui donner la version suivante:

\medskip

{\bf E)}
\label{condition-E}
{\em Si, \`a l'int\'erieur de la r\'egion introduite ci-dessus, on
fixe $n-1$ points mutuellement en position g\'en\'erale qui sont
choisis de telle sorte qu'ils ne restent au repos simultan\'ement que
par l'action d'un groupe \`a un param\`etre de mouvements, et si $X f$
est la transformation infinit\'esimale de ce groupe \`a un
param\`etre, alors les \'equations finies correspondantes:
\def\theequation{14}\begin{equation}
x_\nu'
=
x_\nu
+
{\textstyle{\frac{t}{1}}}\,X\,x_\nu
+
{\textstyle{\frac{t^2}{1\,\cdot\,2}}}\,XX\,x_\nu
+
\cdots
\quad\quad\quad
{\scriptstyle{(\nu\,=\,1\,\cdots\,n)}}
\end{equation}
doivent \^etre constitu\'ees de telle sorte que, si $t$ cro\^{\i}t
toujours \`a partir de la valeur nulle, tous les autres points: $x_1
\dots x_n$ de la r\'egion reviennent finalement tous en m\^eme temps \`a
leur position initiale. Bri\`evement, il est demand\'e que les
\'equations~\thetag{ 14} repr\'esentent un mouvement
qui a une p\'eriode r\'eelle.}

\HEAD{Critique des recherches helmholtziennes.}{
Division\,\,V.\,\,\,Chapitre\,\,21.\,\,\,\S\,\,94.}

\sectiondritterV{\sf\S\,\,94.
\\
Critique des conclusions que Monsieur de Helmholtz
\\
tire de ses axiomes.}
\label{S-94}
\setcounter{footnote}{0}

Dans les pr\'ec\'edents paragraphes, nous avons donn\'e aux axiomes
helmholtziens une version qui fait ressortir clairement le caract\`ere
<<\,th\'eorie des groupes\,>> qu'a le probl\`eme dans son ensemble.
Nous voulons maintenant examiner d'une mani\`ere critique les
cons\'equences que Monsieur de Helmholtz a tir\'ees de ses axiomes.
Afin de pouvoir effectuer cela de la mani\`ere la plus commode
possible, nous traduisons tout d'abord ces cons\'equences dans le
langage de la th\'eorie des groupes, autant qu'il est possible. Comme
Monsieur de Helmholtz s'est restreint dans sa recherche \`a l'espace
de dimension trois, nous faisons naturellement de m\^eme.

\bigskip

\renewcommand{\thefootnote}{\fnsymbol{footnote}}
Dans l'espace \`a trois dimensions chaque groupe de mouvements qui
satisfait les exigences helmholtziennes est \`a six param\`etres.
Monsieur de Helmholtz consid\`ere maintenant\footnote[1]{\,
G\"ott. Nachr. 1868, p.~202 sq. La mani\`ere dont il construit ces
mouvements est remarquable. Il s\'electionne un mouvement
d\'etermin\'e $S$, par lequel un point d\'etermin\'e $P_1$ est
envoy\'e vers une nouvelle position $P$, et d'un autre c\^ot\'e, il
s'imagine disposer du mouvement le plus g\'en\'eral contenu dans le
groupe, par lequel $P_1$ est \'egalement envoy\'e sur $P$. Sous ces
hypoth\`eses, $S^{ - 1} T$ est visiblement le mouvement le plus
g\'en\'eral appartenant au groupe au cours duquel $P$ reste au
repos. Avec cela, il est d\'emontr\'e non seulement qu'il y a de tels
mouvements, mais aussi que les \'equations des mouvements concern\'es
peuvent \^etre rapport\'ees \`a une forme telle qu'elles peuvent
s'appliquer aux points qui sont infiniments voisins de
$P$. Naturellement, Monsieur de Helmholtz ne calcule pas du tout de
mani\`ere symbolique avec les transformations; chez lui, le concept
de {\sl famille de transformations} ne se rencontre pas m\^eme une
seule fois explicitement. En utilisant les concepts et les mani\`eres
de s'exprimer de la th\'eorie des groupes, nous donnons en g\'en\'eral
une forme plus rigoureuse aux d\'eveloppements helmholtziens.}, comme
nous pouvons l'exprimer, tous les mouvements contenus dans le groupe,
\renewcommand{\thefootnote}{\arabic{footnote}}
qui laissent invariant un point d\'etermin\'e. Puisque le groupe des
mouvements est transitif, l'ensemble de tous les mouvements qui
laissent au repos un point d\'etermin\'e, constitue un groupe \`a
trois param\`etres. Si, par souci de simplicit\'e, nous nous imaginons
que le point invariant est choisi comme \'etant l'origine des
coordonn\'ees, alors les \'equations de ce groupe \`a trois
param\`etres ont la forme:
\def\theequation{15}\begin{equation}
\left\{
\aligned
x'
&
=
\lambda_1\,x+\lambda_2\,y+\lambda_3\,z+\cdots
\\
y'
&
=
\mu_1\,x+\mu_2\,y+\mu_3\,z+\cdots
\\
z'
&
=
\nu_1\,x+\nu_2\,y+\nu_3\,z+\cdots,
\endaligned\right.
\end{equation}
o\`u les constantes $\lambda, \mu, \nu$ et les termes d'ordre
sup\'erieur qui ont \'et\'e supprim\'es ne d\'ependent que de trois
param\`etres arbitraires et o\`u le d\'eterminant des $\lambda, \mu,
\nu$ ne s'annule s\^urement pas identiquement.

Cependant, Monsieur de Helmholtz ne consid\`ere pas le groupe~\thetag{
15} lui-m\^eme, mais il \'etudie seulement de quelle mani\`ere les
points infiniment proches de l'origine des coordonn\'ees sont
transform\'es \label{454} par le groupe~\thetag{ 15}, autrement dit:
il se restreint \`a la consid\'eration des transformations:
\def\theequation{16}\begin{equation}
\left\{
\aligned
dx'
&
=
\lambda_1\,dx+\lambda_2\,dy+\lambda_3\,dz
\\
dy'
&
=
\mu_1\,dx+\mu_2\,dy+\mu_3\,dz
\\
dz'
&
=
\nu_1\,dx+\nu_2\,dy+\nu_3\,dz,
\endaligned\right.
\end{equation}
que l'on obtient lorsqu'on diff\'erentie les \'equations~\thetag{ 15}
et que l'on pose ensuite: $x = y = z = 0$. Il est clair que ces
transformations en les variables: $dx, dy, dz$ forment un groupe, et
pour pr\'eciser, nul autre groupe que le groupe lin\'eaire homog\`ene
qui est attach\'e \deutsch{zugeordnet} \`a l'origine des coordonn\'ees
par tous les mouvements du groupe, et qui indique notamment comment
sont transform\'es les \'el\'ements lin\'eaires: $dx : dy : dz$
passant par l'origine des coordonn\'ees, d\`es qu'on l'a fix\'ee ({\em
voir} Tome~I, Th\'eor\`eme~109, 
p.~603\footnote{\,
%%%%%%%%%%%%%%%%%%%%%%%-------DEBUT--------%%%%%%%%%%%%%%%%%%%%%%%%%%%
Ce théorème énonce que chaque groupe $X_1f, \dots, X_r f$ de l'espace
$x_1, \dots, x_n$ associe à tout point $x_1^0, \dots, x_n^0$ fixé en
position générale un groupe linéaire homogène de $R_n$ parfaitement
déterminé, lequel indique de quelle manière sont transformés les
éléments linéaires passant par ce point.
} %%%%%%%%%%%%%%%%%%%%%%%%-----FIN-----%%%%%%%%%%%%%%%%%%%%%%%%%%%%%%%
).

%%%\Fill

Jusqu'\`a ce point, les d\'eveloppements de Monsieur de Helmholtz sont
irr\'eprochables. Mais maintenant, il admet tacitement et sans un mot
de justification que tous ses axiomes, qu'il a pos\'es au sujet des
mouvements encore possibles apr\`es fixation d'un point, peuvent aussi
s'appliquer aux points qui sont infiniment voisins du point fix\'e, et
donc que s'ils s'appliquent \`a des points
finiment \'eloign\'es les uns des
autres, il s'ensuit qu'ils s'appliquent en m\^eme temps \`a des points
infiniment voisins. Pour l'exprimer plus pr\'ecis\'ement: il
s'imagine le groupe lin\'eaire homog\`ene:
\def\theequation{16'}\begin{equation}
\left\{
\aligned
x'
&
=
\lambda_1\,x+\lambda_2\,y+\lambda_3\,z
\\
y'
&
=
\mu_1\,x+\mu_2\,y+\mu_3\,z
\\
z'
&
=
\nu_1\,x+\nu_2\,y+\nu_3\,z
\endaligned\right.
\end{equation} 
comme un groupe de mouvements qui laisse invariante l'origine des
coordonn\'ees, et il pose \`a l'avance que le groupe~\thetag{ 16'}
satisfait alors toujours ses axiomes, lorsque le groupe~\thetag{ 15}
les satisfait. C'est sur cette supposition que reposent tous ses
d\'eveloppements subs\'equents.

Nous voulons d'abord montrer que cette supposition a la m\^eme
signification qu'une autre supposition, que l'on peut exprimer de
mani\`ere appropri\'ee, et ensuite, {\em au moyen d'une s\'erie
d'exemples, nous allons montrer clairement l'irrecevabilit\'e de cette
supposition dans son int\'egralit\'e}.

\label{455} Le groupe~\thetag{ 15} est engendr\'e par trois
transformations infinit\'esimales de la forme:
\def\theequation{17}\begin{equation}
\left\{
\aligned
&
\quad 
(\alpha_{k1}\,x+\alpha_{k2}\,y+\alpha_{k3}\,z+\cdots)\,p+
\\
&
+
(\beta_{k1}\,x+\beta_{k2}\,y+\beta_{k3}\,z+\cdots)\,q+
\\
&
+
(\gamma_{k_1}\,x+\gamma_{k2}\,y+\gamma_{k3}\,z+\cdots)\,r
\\
&
\quad\quad\quad\quad\quad
{\scriptstyle{(k\,=\,1,\,2,\,3)}},
\endaligned\right.
\end{equation}
o\`u les $\alpha, \beta, \gamma$ d\'esignent des constantes et o\`u les
termes supprim\'es sont d'un ordre sup\'erieur par rapport \`a $x, y, z$.
Le groupe \`a six param\`etres de tous les mouvements contient encore,
hormis les transformations infinit\'esimales~\thetag{ 17}, trois
transformations de la forme:
\def\theequation{18}\begin{equation} 
p+\cdots,\quad\quad
q+\cdots,\quad\quad
r+\cdots.
\end{equation}
D'un autre c\^ot\'e, \renewcommand{\thefootnote}{\fnsymbol{footnote}}
le groupe~\thetag{ 16'} est le groupe r\'eduit associ\'e au
groupe~\thetag{ 17}, et ses transformations
infinit\'esimales\footnote[1]{\, Nous ne voulons pas passer sous
silence le fait que Monsieur de Helmholtz op\`ere avec les transformations
infinit\'esimales du groupe {\em lin\'eaire homog\`ene}~\thetag{ 16'},
et on peut m\^eme dire qu'il consid\`ere, quoique de mani\`ere
inconsciente, les groupes \`a un param\`etre qui sont engendr\'es par
certaines de ces transformations infinit\'esimales. N\'eanmoins, on
ne se trouve chez lui en aucune fa\c con le concept g\'en\'eral de
transformation infinit\'esimale, et encore moins le concept
g\'en\'eral de groupe \`a un param\`etre.}:
\def\theequation{19}\begin{equation}
\left\{
\aligned
L_kf
=
(\alpha_{k1}\,x+
&
\alpha_{k2}\,y+\alpha_{k3}\,z)\,p
+
(\beta_{k1}\,x+\beta_{k2}\,y+\beta_{k3}\,z)\,q
+
\\
&
+
(\gamma_{k1}\,x+\gamma_{k2}\,y+\gamma_{k3}\,z)\,r
\\
&
\quad\quad\quad
{\scriptstyle{(k\,=\,1,\,2,\,3)}}
\endaligned\right.
\end{equation}
proviennent de~\thetag{ 17} lorsqu'on supprime tous les termes d'ordre
sup\'erieur, et donc elles ne sont pas n\'ecessairement
ind\'ependantes les unes des autres.
\renewcommand{\thefootnote}{\arabic{footnote}}
Remarquons en outre
(\cf le Tome~I, p.~606), qu'en tenant compte des hypoth\`eses
pos\'ees, les transformations infinit\'esimales:
\def\theequation{20}\begin{equation}
p,\quad
q,\quad
r,\quad
L_1f,\quad
L_2f,\quad
L_3f,
\end{equation}
engendrent aussi un 
groupe\footnote{\,
%%%%%%%%%%%%%%%%%%%%%%%-------DEBUT--------%%%%%%%%%%%%%%%%%%%%%%%%%%% 
\label{justification-groupe-reduit}
Voici la justification. 
Dans des coordonnées $(x_1, \dots, x_n)$, soit
un groupe linéaire homogène quelconque constitué
d'un certain nombre $m \leqslant n^2$ de
générateurs infinitésimaux: 
\[
L_\nu
=
{\textstyle{\sum_{j=1}^n}}\,
\big(
{\textstyle{\sum_{i=1}^n}}\,
\alpha_{\nu ij}\,x_i
\big)
{\textstyle{\frac{\partial}{\partial x_j}}}
\ \ \ \ \ \ \ \ \ \ \ \ \
{\scriptstyle{(\nu\,=\,1\,\cdots\,m)}}
\]
à coefficients linéaires, donc par hypothèse
stable par crochets de Lie. 
Alors la collection: 
\[
p_1,\,\dots,\,p_n,\ \ \ \ \
L_1,\,\dots,\,L_m
\]
obtenue en leur adjoignant toutes les transformations infinitésimales
du groupe (transitif et commutatif) des translations est elle aussi
stable par crochets (donc engendre un groupe au sens de Lie), puisque:
\[
\big[p_k,\,L_\nu\big]
=
\big[
{\textstyle{\frac{\partial}{\partial x_k}}},\,\,
{\textstyle{\sum_{j=1}^n}}\,
\big(
{\textstyle{\sum_{i=1}^n}}\,
\alpha_{\nu ij}\,x_i
\big)
{\textstyle{\frac{\partial}{\partial x_j}}}
\big]
=
{\textstyle{\sum_{j=1}^n}}\,
\alpha_{\nu kj}\,
{\textstyle{\frac{\partial}{\partial x_j}}}
=
{\textstyle{\sum_{j=1}^n}}\,
\alpha_{\nu kj}\,p_j.
\]
} %%%%%%%%%%%%%%%%%%%%%%%%-----FIN-----%%%%%%%%%%%%%%%%%%%%%%%%%%%%%%% 
transitif, qui toutefois n'est pas n\'ecessairement \`a six
param\`etres. Ce groupe~\thetag{ 20} n'est autre que le groupe
r\'eduit relatif au groupe~\thetag{ 17}, \thetag{ 18} de tous les
mouvements.

L'hypoth\`ese introduite tacitement par Monsieur de Helmholtz a
maintenant simplement le sens que le groupe~\thetag{ 19} satisfait
alors chaque axiome ayant \'et\'e pos\'e au sujet de tous les
mouvements qui sont encore possibles apr\`es fixation d'un point,
lorsque le groupe~\thetag{ 17} satisfait ces axiomes. Mais d'un autre
c\^ot\'e, il est clair que le groupe~\thetag{ 17}, \thetag{ 18} de
tous les mouvements satisfait tous les axiomes helmholtziens si et
seulement si le groupe~\thetag{ 17} satisfait encore les axiomes
restants qui s'appliquent apr\`es fixation de l'origine des
coordonn\'ees. Enfin, il est \'evident qu'entre les deux
groupes~\thetag{ 20} et~\thetag{ 19}, on a exactement la m\^eme
relation qu'entre les groupes~\thetag{ 17}, \thetag{ 18} et \thetag{
17}. Par cons\'equent, ce qui a \'et\'e admis par Monsieur de
Helmholtz revient simplement \`a supposer que {\em le groupe
r\'eduit~\thetag{ 20} satisfait alors toujours chaque exigence qu'il a
pos\'ee, lorsque le groupe initial~\thetag{ 17}, \thetag{ 18} les
satisfait}.

Ceci, Monsieur de Helmholtz l'a admis tacitement, sans la moindre
indication de d\'emonstration, et sans faire
aucune allusion au fait que cela
n\'ecessite \`a vrai dire une d\'emonstration. 
\renewcommand{\thefootnote}{\fnsymbol{footnote}} 
Par une s\'erie d'exemples\footnote[1]{\, Lie a communiqu\'e ces exemples
pour la premi\`ere fois dans les Comptes Rendus de 1892 (Vol.~114,
p.~463).}, nous voulons montrer que ce qui a \'et\'e admis n'est pas
fond\'e.
\renewcommand{\thefootnote}{\arabic{footnote}} 
Nous allons trouver qu'un groupe transitif \`a six
param\`etres de $R_3$ peut tr\`es bien satisfaire certaines des
exigences helmholtziennes, sans que le groupe r\'eduit qui lui est
associ\'e les satisfasse; d'un autre c\^ot\'e, nous allons voir
qu'\`a un groupe transitif \`a six param\`etres qui ne satisfait pas
certaines des exigences helmholtziennes peut tr\`es bien \^etre
associ\'e un groupe r\'eduit qui satisfait les exigences en question.

\bigskip

Comme nous l'avons vu, une partie des exigences helmholtziennes se
ram\`ene \`a ce que deux points ont un et un seul invariant
relativement \`a tous les mouvements du groupe. Maintenant,
relativement au groupe transitif \`a six param\`etres:
\def\theequation{21}\begin{equation}
q,\quad
xq+r,\quad
x^2q+2xr,\quad
x^3q+3x^2r,\quad
x^4q+4x^3r,\quad
p
\end{equation}
les deux points: $x_1, y_1, z_1$ et $x_2, y_2, z_2$ ont le
seul\footnote{\, La Proposition~2 p.~\pageref{Satz-2} s'applique à
l'origine, qui est un point de position générale. }
invariant\footnote{\, Effectivement, $x_2 - x_1$ est annihilé
identiquement par $X_k^{ (1)} + X_k^{ (2)}$, pour $k=1, \dots, 6$,
ce qui est trivial pour 
les cinq premiers opérateurs différentiels, le sixième
s'écrivant $p_1 + p_2$, d'où $(p_1 + p_2) ( x_2 - x_1) = 1 - 1 = 0$. 
}: $x_2 - x_1$, tandis que relativement au groupe r\'eduit associ\'e:
\def\theequation{21'}\begin{equation}
q,\quad
r,\quad
xr,\quad
p,
\end{equation}
ils ont au contraire les deux invariants\footnote{\, D'après
l'équation~\thetag{ 2} p.~\pageref{eq-2-p-400}, une fonction
quelconque $J = J ( x_1, y_1, z_1, x_2, y_2, z_2)$ est invariante si
et seulement si elle est annihilée identiquement par $Y_l^{ (1)} +
Y_l^{ (2)}$, pour $l=1, \dots, 4$. L'annihilation par $p_1 + p_2$, par
$q_1 + q_2$ et par $r_1 + r_2$ équivaut à ce que $J = J ( x_2 - x_1,
y_2 - y_1, z_2 - z_1)$. Enfin, l'annihilation par $x_1 r_1 + x_2 r_2$
implique que $J = J ( x_2 - x_1, y_2 - y_1)$. La vérification
(similaire) des deux affirmations qui suivent immédiatement est
laissée au lecteur. }:
$x_2 - x_1$ et
$y_2 - y_1$. D'un autre c\^ot\'e, relativement au groupe transitif
\`a six param\`etres:
\def\theequation{22}\begin{equation}\label{457-0}
q,\quad
xq+r,\quad
x^2q+2xr,\quad
x^3q+3x^2r,\quad
p,\quad
xp-zr,
\end{equation}
ces deux points n'ont absolument aucun invariant, tandis que,
relativement au groupe r\'eduit associ\'e:
\def\theequation{22'}\begin{equation}
q,\quad
r,\quad
xr,\quad
p,\quad
xp-zr
\end{equation}
ils ont l'invariant: $y_2 - y_1$.

De l\`a, il r\'esulte que les deux groupes~\thetag{ 17}, 
\thetag{ 18} et \thetag{ 20} ne satisfont pas en g\'en\'eral
simultan\'ement les exigences qui ont \'et\'e admises.

\bigskip

Consid\'erons ensuite le groupe:
\def\theequation{23}\begin{equation}
\left\{
\aligned
&
q,\quad
p,\quad
xq+r,\quad
x^2q+2xr,\quad
xp+yq+cr
\\
&
\quad\quad\quad\quad
x^2p+2xyq+2(cx+y)\,r,
\endaligned\right.
\end{equation}
que nous avons d\'ej\`a rencontr\'e \`a la
page~\pageref{418}. D'apr\`es la page~\pageref{419-0}, deux points ont
un et un seul invariant relativement \`a ce groupe tandis que $s > 2$
points n'ont pas d'invariant essentiel, donc ce groupe satisfait
certaines exigences qui d\'ecoulent des axiomes helmholtziens (\cf
p.~\pageref{452}). Cependant, le groupe r\'eduit
associ\'e\footnote{\, Pour obtenir trois transformations
infinitésimales linéairement indépendantes
fixant l'origine, une combinaison linéaire préalable est requise:
retrancher de la cinquième transformation la troisième
multipliée par $c$, ce qui donne $xp + yq - cxq$. }:
\def\theequation{23'}\begin{equation}
q,\quad
p,\quad
r,\quad
xr,\quad
xp+yq-cxq,\quad
yr
\end{equation}
ne satisfait pas ces exigences, car il contient des transformations
infinit\'esimales, et m\^eme trois, \`a savoir: $r, \ xr, \ yr$, qui ont
des courbes int\'egrales en commun, alors que, d'apr\`es la Proposition~1,
p.~\pageref{Satz-1}, cela ne peut pas se produire lorsque deux points
doivent avoir un et un seul invariant, alors que $s > 2$ points ne
doivent pas avoir d'invariant essentiel.

Par cons\'equent, les exigences en question ne sont pas toujours
n\'ecessairement satisfaites par le groupe~\thetag{ 20} lorsqu'elles le
sont par le groupe~\thetag{ 17}, \thetag{ 18}.

\bigskip

Enfin, l'axiome de monodromie peut tr\`es bien appartenir aussi aux
exigences qui sont satisfaites par le groupe~\thetag{ 17}, \thetag{
18}, sans que celles-ci soient satisfaites par le groupe~\thetag{ 20}.

Le groupe:
\def\theequation{24}\begin{equation}\label{457}
\left\{
\aligned
&
p,\quad
q,\quad
xp+yq+r,\quad
yp-xq
\\
&
\quad
(x^2-y^2)\,p+2xyq+2xr
\\
&
\quad
2xyp+(y^2-x^2)\,q+2yr,
\endaligned\right.
\end{equation}
que nous avons d\'ej\`a trouv\'e \`a la
page~\pageref{eq-58-p-432}\footnote{\, à savoir le groupe~\thetag{ 58}
avec $a = 1$ et $b = 0$, ce qui a donné le groupe~7
p.~\pageref{432-1}.}, en fournit la d\'emonstration. Ce groupe {\em
extr\^emement curieux} \deutsch{äusserst merkwürdige Gruppe}
satisfait en effet toutes les exigences de l'axiome helmholtzien de
monodromie, tandis que le groupe r\'eduit associ\'e:
\def\theequation{24'}\begin{equation}
p,\quad
q,\quad
r,\quad
yp-xq,\quad
xr,\quad
yr
\end{equation}
ne les satisfait pas.

Il est facile de voir que le groupe~\thetag{ 24'} ne satisfait pas
l'axiome de monodromie. Si nous fixons en effet deux points qui sont
mutuellement en position g\'en\'erale, par exemple l'origine des
coordonn\'ees et le point $x_0, y_0, z_0$, o\`u $x_0$ et $y_0$ ne
s'annulent pas, alors les mouvements encore possibles forment un
groupe \`a un param\`etre qui est engendr\'e par la transformation
infinit\'esimale\footnote{\,
\`A un facteur multiplicatif près, 
c'est la seule transformation infinitésimale qui s'annule en $(0, 0,
0)$ et en $(x_0, y_0, z_0)$. }: $(x_0 y - y_0 x)\, r$, et dont les
\'equations finies sont par cons\'equent de la forme:
\[
x'=x,\quad\quad
y'=y,\quad\quad
z'
=
z+t(x_0y-y_0x).
\]
Mais alors ici, le point $x, y, z$ ne retourne manifestement jamais
\`a sa position initiale, quand $t$ cro\^{\i}t continuellement \`a
partir de z\'ero.

D'autre part, afin de nous persuader que le groupe~\thetag{ 24}
satisfait l'axiome de monodromie, nous devons examiner d'un peu plus
pr\`es les mouvements de ce groupe.

Si l'on fixe un point $x_0, y_0, z_0$ relativement \`a l'action du
groupe~\thetag{ 24}, alors chaque autre point $x, y, z$ se meut en
g\'en\'eral d'une mani\`ere compl\`etement libre sur la
pseudosph\`ere\footnote{\, L'unique invariant d'une paire de points a
déjà été calculé p.~\pageref{invariant-58} (faire $b = 0$). }:
\def\theequation{25}\begin{equation}
\big\{
(x-x_0)^2+(y-y_0)^2
\big\}\,
e^{-z}
=
\text{\rm const.}
\end{equation}
centr\'ee en ce point: $x_0, y_0, z_0$; seuls les points de la droite:
$x = x_0$, $y = y_0$, qui constitue \'evidemment une pseudosph\`ere
par elle-m\^eme, font exception, car ils restent enti\`erement au
repos.

Si donc nous voulons fixer deux points distincts qui ne restent
simultan\'ement au repos que par l'action d'un sous-groupe \`a un
param\`etre de~\thetag{ 24}, alors nous devons choisir ces points de
telle sorte que leur droite de liaison
\deutsch{Verbindungslinie}
ne soit pas parall\`ele \`a l'axe des $z$.

Deux tels points sont par exemple l'origine des coordonn\'ees et un
point quelconque $x_0, y_0, z_0$ pour lequel $x_0$ et $y_0$ ne
s'annulent pas tous deux. Les \'equations finies du groupe \`a un
param\`etre qui laisse invariants ces deux points se d\'eterminent
gr\^ace au syst\`eme simultan\'e\footnote{\, On détermine à l'avance
une certaine combinaison linéaire $W = \alpha\, T + \beta\, U + \gamma
\, V$ des trois dernières transformations infinitésimales~\thetag{ 24}
(qui s'annulent déjà à l'origine):
\[
W
=
\big[
\alpha y+\beta(x^2-y^2)+2\gamma xy
\big]\,p
+
\big[
-\alpha x+2\beta xy+\gamma(y^2-x^2)
\big]\,q
+
\big[
2\beta x+2\gamma y
\big]\,r
\]
afin qu'elle s'annule aussi en $(x_0, y_0, z_0)$.
\'Egaler à zéro les trois coefficients de
$p$, $q$, $r$ dans $W \big\vert_{ (x_0, y_0, z_0)}$ donne alors
l'unique (à un facteur multiplicatif près) solution:
\[
\alpha
=
y_0,\ \ \ \ \ \ \ \ \ \
\beta
=
{\textstyle{\frac{y_0}{x_0^2+y_0^2}}},\ \ \ \ \ \ \ \ \ \
\gamma
=
-{\textstyle{\frac{x_0}{x_0^2+y_0^2}}}.
\]
Le système considéré de trois équations différentielles ordinaires
représente alors les courbes intégrales de $W$, d'inconnues $\big( x'
(t), y' (t), z' (t) \big)$ valant $(x, y, z)$ pour $t = 0$.
}:
\def\theequation{26}\begin{equation}
\left\{
\aligned
\frac{dx'}{dt}
&
=
\quad
y'
+
\frac{y_0({x'}^2-{y'}^2)}{x_0^2+y_0^2}
-
\frac{2x_0x'y'}{x_0^2+y_0^2}
\\
\frac{dy'}{dt}
&
=
-
x'
+
\frac{2y_0x'y'}{x_0^2+y_0^2}
+
\frac{x_0({x'}^2-{y'}^2)}{x_0^2+y_0^2}
\\
\frac{dz'}{dt}
&
=
\frac{2(y_0x'-x_0y')}{x_0^2+y_0^2},
\endaligned\right.
\end{equation}
avec les conditions initiales bien connues: $[ x']_{ t= 0} = x$, \etc
\`A partir de~\thetag{ 26}, on trouve\footnote{\, 
Il s'agit d'intégrer l'équations différentielle $\frac{ dw}{ dt} = - i
\big\{ w - w^2 / \alpha \big\}$ d'inconnue complexe $w = w ( t) \in
\C$, où $\alpha := x_0 + iy_0$ est une constante, c'est-à-dire $\frac{
\alpha \, dw}{ \alpha w - w^2} = - i dt$, ou encore, après décomposition
en éléments simples et séparation des quotients différentiels:
\[
dw
\big(
{\textstyle{\frac{1}{w}}}
-
{\textstyle{\frac{1}{w-\alpha}}}
\big)
=
-
idt.
\]
Une intégration donne $\log \frac{ w}{ w - \alpha} = - it + {\rm
const.}$, puis par exponentiation $w = - \frac{ c \, \alpha}{ e^{ it}
- c}$, où $c$ est une constante que l'on détermine par la condition
initiale $w
\big\vert_{ t = 0} = x + iy$, c'est-à-dire $c = \frac{ - x - iy}{ x_0
- x + i ( y_0 - y)}$. Enfin, on remplace $c$ dans l'expression de $w$.
Très exceptionnellement, le texte allemand comporte ici une
coquille: au dénominateur, entre accolades, il est écrit $x - x_0 +
i ( y - y_0)$. }:
\[
\frac{d}{dt}\,
(x'+iy')
=
-i
\Big\{
x'+iy'
-
\frac{(x'+iy')^2}{x_0+iy_0}
\Big\},
\]
et de l\`a, par int\'egration: 
\def\theequation{27}\begin{equation}
x'+iy'
=
\frac{(x_0+iy_0)(x+iy)}{
\{x_0-x+i(y_0-y)\}\,e^{it}+x+iy}.
\end{equation}

\noindent
D'un autre c\^ot\'e, puisque l'on sait que le point $x, y, z$ se meut
sur la pseudosph\`ere passant par lui et qui a pour centre l'origine des
coordonn\'ees, on obtient encore entre $x, y, z$ et $x', y', z'$
l'\'equation:
\[ ({x'}^2+{y'}^2)\,e^{-z'} = (x^2+y^2)\,e^{-z},
\]
d'o\`u on d\'eduit imm\'ediatement: 
\def\theequation{28}\begin{equation}
z'
=
z
+
{\rm log}\,\frac{{x'}^2+{y'}^2}{x^2+y^2}.
\end{equation}
Ainsi, on a compl\`etement d\'etermin\'e les transformations finies du
groupe \`a un param\`etre dont il est question dans le discours. Il
serait facile de les indiquer sous une forme r\'eelle, mais pour notre
objectif, ce n'est pas du tout n\'ecessaire, car les
\'equations~\thetag{ 27} et~\thetag{ 28} montrent d\'ej\`a que notre
groupe~\thetag{ 28} satisfait l'axiome de monodromie. Si en effet on
soumet l'espace au mouvement continu qui est sp\'ecifi\'e par les
\'equations~\thetag{ 27}, \thetag{ 28},
\renewcommand{\thefootnote}{\fnsymbol{footnote}} et si on laisse la
variable $t$ parcourir toutes les valeurs r\'eelles entre $0$ et
$2\pi$, alors tous les points retournent
simultan\'ement\footnote[1]{\, Les $\infty^2$ courbes int\'egrales
qui, \`a cette occasion, sont parcourues par les points de l'espace,
sont les courbes obtenues par intersection des deux familles de
pseudosph\`eres:
\[
(x^2+y^2)\,e^{-z}
=
\text{\rm const.}
\quad\quad
\{(x-x_0)^2+(y-y_0)^2\}\,e^{-z}
=
\text{\rm const.}
\]
Leurs projections sur le plan des $x, y$ sont alors en g\'en\'eral des
cercles, mais ce sont aussi exceptionnellement
des lignes droites.
}
\`a leur position initiale pour $t = 2\pi$.
\renewcommand{\thefootnote}{\arabic{footnote}}

Ainsi, nous voyons que le groupe~\thetag{ 24} satisfait en fait
l'axiome de monodromie, alors que le groupe r\'eduit associ\'e~\thetag{
24'} ne le satisfait pas.

Pour cette raison encore, le groupe~\thetag{ 24} est alors aussi
particuli\`erement curieux, parce qu'il montre de la mani\`ere la plus
convaincante qui soit, combien il est impossible de d\'eduire quelque
chose au sujet du comportement de points infiniment voisins \`a partir
du comportement de points qui sont finiment \'eloign\'es les uns des
autres, et combien il est p\'erilleux d'extrapoler
\`a ces points-l\`a, comme l'a fait
Monsieur de Helmholtz, les axiomes qui ont \'et\'e pos\'es pour les
points finiment \'eloign\'es les uns des autres.

Si en effet l'on fixe un point pour l'action du groupe~\thetag{ 24},
par exemple l'origine des coordonn\'ees, alors chaque autre point se
meut en g\'en\'eral, comme nous l'avons vu plus haut, sur l'une des
$\infty^1$ surfaces invariantes:
\[
(x^2+y^2)\,e^{-z}
=
\text{\rm const}.
\]
Comme on s'en convainc facilement, les points de chacune de ces
surfaces sont transform\'es par l'action d'un groupe \`a trois
param\`etres qui a la m\^eme structure\footnote{\,
\deutsch{gleichzusammengesetzt ist}: possède 
les mêmes constantes de structure dans les crochets de Lie. En effet,
les mouvements encore possibles par le sous-groupe de~\thetag{ 24} qui
laisse fixe l'origine sont engendrés par les trois générateurs
infinitésimaux:
\[
P
:=
(x^2-y^2)p+2xyq+2xr,\ \ \ \ \ \ \ \ \
Q
:=
2xyp+(y^2-x^2)q+2yr
\ \ \ \ \ \ \
\text{\rm et}
\ \ \ \ \ \ \ 
yp-xq.
\]
Ils satisfont bien les mêmes relations de commutation: $\big[ P,
\, yp - xq \big] = -Q$; $\big[ Q, \, yp - xq \big] = P$ et: $\big[ P,
\, Q \big] = 0$ que les trois générateurs infinitésimaux: $p$, $q$,
$yp - xq$ des mouvements euclidiens du plan. } que le {\em groupe des
mouvements euclidiens} d'un plan, et qui lui est de surcro\^{\i}t
semblable par une transformation ponctuelle r\'eelle.

Les points qui sont infiniment voisins de l'origine des coordonn\'ees
se comportent d'une mani\`ere totalement diff\'erente. Afin de se
constituer une image de la fa\c con dont sont transform\'es ces
points, lorsque l'origine des coordonn\'ees est fix\'ee, on peut au
mieux envisager le groupe par lequel sont transform\'es les $\infty^2$
\'el\'ements lin\'eaires passant par l'origine des coordonn\'ees, car
chaque point qui est infiniment voisin de l'origine des coordonn\'ees
d\'etermine un tel \'el\'ement lin\'eaire. Si on entend par $x', y',
z'$ les coordonn\'ees homog\`enes d'un \'el\'ement lin\'eaire, alors
ce groupe est:
\[
y'p'-x'q',
\quad\quad
x'r',
\quad\quad
y'r',
\]
et on v\'erifie alors imm\'ediatement que la vari\'et\'e des
$\infty^2$ \'el\'ements lin\'eaires est transform\'ee par 
l'action d'un groupe projectif qui est {\em dualistique du groupe des
mouvements euclidiens d'un plan}.

%%%\Mathematiques 
%%%[[{\sc dualistisch}~??]]

Ainsi, on voit que, {\em lors du passage des points finiment
\'eloign\'es aux points infiniment voisins, un saut plus total peut se
produire, et que les points infiniment voisins ob\'eissent dans de
telles circonstances
\`a des lois tout autres que celles auxquelles sont soumis les points
finiment \'eloign\'es les uns des autres}.

\bigskip

\label{non-complet-Helmholtz}
Gr\^ace aux exemples pr\'ec\'edents, il a \'et\'e suffisamment
d\'emontr\'e que la supposition que Monsieur de Helmholtz a introduite
tacitement et qui a \'et\'e d\'ecrite plus pr\'ecis\'ement aux
pages~\pageref{455} sq. est erron\'ee. Et maintenant, comme ses
consid\'erations ult\'erieures prennent enti\`erement cette
supposition comme point de d\'epart et n'ont force de preuve que sur
la base de cette supposition, nous parvenons donc au r\'esultat que
Monsieur de Helmholtz n'a pas d\'emontr\'e l'assertion qu'il \'enonce
\`a la fin de son travail, \`a savoir: il n'a pas d\'emontr\'e que
ses axiomes suffisent \`a caract\'eriser les mouvements euclidiens et
non-euclidiens.

Apr\`es nous \^etre convaincus de cette fa\c con que les d\'eveloppements
de Monsieur de Helmholtz ne constituent pas une d\'emonstration, nous
allons nous rattacher tout d'abord 
dans le \S~95 aux hypoth\`eses que Monsieur de
Helmholtz a pos\'ees au cours de ses calculs.

\HEAD{Critique des recherches helmholtziennes.}{
Division\,\,V.\,\,\,Chapitre\,\,21.\,\,\,\S\,\,95.}

\sectiondritterV{\sf\S\,\,95.
\\
Consid\'erations se rattachant aux calculs helmholtziens.}
\label{S-95}
\setcounter{footnote}{0}

Dans le paragraphe pr\'ec\'edent, nous avons vu que Monsieur de Helmholtz
applique sans plus de fa\c cons ses axiomes \`a des points infiniment
voisins. Dans la supposition erron\'ee que cela soit autoris\'e
r\'eside la faiblesse des d\'eveloppements helmholtziens; l'introduction
de cette supposition \^ote toute force de preuve \`a sa r\'eflexion.

Mais on peut \'eviter cette erreur, si l'on dispose d\`es le d\'ebut
les axiomes helmholtziens de telle sorte qu'ils se r\'ef\`erent
simplement \`a des points infiniment voisins; on peut alors toujours
s'arranger pour que les calculs helmholtziens parviennent
effectivement au but en prenant pour base les axiomes formul\'es de
cette fa\c con. Comme cela est possible de diverses mani\`eres et
puisque cela ne vaut pas la peine de discuter \`a fond les
diff\'erentes possibilit\'es, nous nous contenterons de ce qui suit:
nous poserons un syst\`eme d'axiomes qui se r\'ef\`ere \`a des points
infiniment voisins et qui, sans concorder avec les hypoth\`eses que
Monsieur de Helmholtz a introduites tacitement dans ses calculs, sont
cependant tr\`es analogues avec celles-ci. Ensuite, nous
d\'emontrerons {\em par des calculs qui ne s'\'ecartent pas dans leur
principe de ceux de Monsieur de Helmholtz}, que ce syst\`eme d'axiomes
suffit pour la caract\'erisation des mouvements euclidiens et
non-euclidiens de l'espace ordinaire.

Les axiomes que nous posons s'\'enoncent ainsi\footnote{\,
%%%%%%%%%%%%%%%%%%%%%%%-------DEBUT--------%%%%%%%%%%%%%%%%%%%%%%%%%%%
La finitude du nombre des paramètres du groupe n'est pas 
explicitement demandé. 
}: %%%%%%%%%%%%%%%%%%%%%%%%-----FIN-----%%%%%%%%%%%%%%%%%%%%%%%%%%%%%%%

\medskip

{\bf I)} {\em L'espace trois fois \'etendu est une vari\'et\'e num\'erique.}

\medskip

{\bf II)} {\em Les mouvements de cet espace forment un groupe continu
r\'eel de transformations ponctuelles.}

\medskip

\label{III-infinitesimal}
{\bf III)} {\em Si l'on fixe un point r\'eel en position g\'en\'erale,
alors le groupe lin\'eaire homog\`ene, qui d\'etermine de quelle
mani\`ere les $\infty^2$ \'el\'ements lin\'eaires r\'eels qui passent
par ce point sont transform\'es, poss\`ede exactement trois
param\`etres.}

\medskip

\label{IV-periodique}
{\bf IV)} {\em Chaque sous-groupe r\'eel \`a un param\`etre du groupe
lin\'eaire homog\`ene mentionn\'e \`a l'instant est constitu\'e de
telle sorte que par son action, tous les \'el\'ements lin\'eaires
r\'eels qui ne restent pas au repos d\'ecrivent un c\^one r\'eel
qu'ils parcourent contin\^ument, \`a l'int\'erieur duquel ils
retournent finalement et simultan\'ement, mais sans revenir en
arri\`ere, \`a leur situation initiale; ou bien, pour exprimer cela
plus pr\'ecis\'ement: si:
\[
\aligned
x_1'
&
=
\alpha_1(t)\,x'+\alpha_2(t)\,y'+\alpha_3(t)\,z'
\\
y_1'
&
=
\beta_1(t)\,x'+\beta_2(t)\,y'+\beta_3(t)\,z'
\\
z_1'
&
=
\gamma_1(t)\,x'+\gamma_2(t)\,y'+\gamma_3(t)\,z'
\endaligned
\]
sont les \'equations finies d'un tel sous-groupe \`a un param\`etre sous
leur forme canonique, il se produit finalement le cas, lorsque la
variable r\'eelle $t$ cro\^{\i}t continuellement \`a partir de $0$ et pour
une certaine valeur finie positive de $t$, que: $x_1' : y_1' : z_1 '$
est proportionnel \`a: $x' : y ' : z'$.}

\medskip

%%%\Fill 
%%%[[m\^emes coordonn\'ees homog\`enes, 
%%%m\^eme point de l'espace projectif.]]

Nous allons montrer que ces axiomes suffisent enti\`erement \`a
caract\'eriser les mouvements euclidiens et non-euclidiens.

Tout d'abord, nous devons d\'eterminer quelle forme poss\`ede le
groupe lin\'eaire homog\`ene mentionn\'e dans les axiomes.

Soit: $x_1 ', x_2 ', x_3'$ les coordonn\'ees de l'\'el\'ement
lin\'eaire passant par un point fix\'e en position g\'en\'erale;
d'apr\`es les hypoth\`eses pos\'ees, le groupe lin\'eaire homog\`ene
$g$ attach\'e \`a ce point contient alors trois transformations
infinit\'esimales ind\'ependantes de la forme:
\def\theequation{29}\begin{equation}
\sum_{\mu\,\nu}^{1,2,3}\,
\alpha_{k\mu\nu}\,x_\mu'\,p_\nu'
\quad\quad\quad\quad
{\scriptstyle{(k\,=\,1,\,2,\,3)}}.
\end{equation}
Parmi ces transformations infinit\'esimales, il peut y en avoir au
plus une\footnote{\, 
Dans le groupe linéaire en dimension $n$ quelconque, 
seule l'homothétie de générateur
infinitésimal $x_1' p_1' + \cdots + x_n' p_n'$ 
laisse fixe toutes les éléments linéaires. } 
qui fixe chaque \'el\'ement lin\'eaire: $x_1 ' : x_2 ' : x_3
'$, et par cons\'equent il est certain que la vari\'et\'e deux fois
\'etendue des $\infty^2$ \'el\'ements lin\'eaires: $x_1' : x_2 ' :
x_3 '$ est transform\'ee par $g$ {\em via}\, l'action d'un groupe
projectif $\mathfrak{ g}$, \label{462} qui poss\`ede soit trois, soit
deux param\`etres.

Si nous rapportons projectivement notre vari\'et\'e des $\infty^2$
\'el\'ements lin\'eaires aux points r\'eels d'un plan, alors
$\mathfrak{ g}$ se pr\'esente comme un groupe projectif $\mathfrak{
g'}$ du plan ayant deux ou trois param\`etres. Chaque sous-groupe
r\'eel \`a un param\`etre de $\mathfrak{ g}'$ est alors constitu\'e de
telle sorte que dans la forme {\em canonique} de ses \'equations
finies:
\[
\mathfrak{x}'
=
\varphi(\mathfrak{x},\mathfrak{y};\,t),
\quad\quad\quad\quad
\mathfrak{y}'=
\psi(\mathfrak{x},\mathfrak{y};\,t),
\]
les quantit\'es $\mathfrak{ x}'$ et $\mathfrak{ y}'$ sont des
fonctions p\'eriodiques de la variable r\'eelle $t$. Tous les points
r\'eels qui ne restent pas invariants par l'action d'un tel groupe \`a
un param\`etre, se meuvent alors sur des courbes r\'eelles qu'ils
parcourent contin\^ument dans toute leur extension, et pour
pr\'eciser, de telle mani\`ere que, sans revenir en arri\`ere, ils
reviennent finalement et simultan\'ement \`a leur position initiale;
si la droite \`a l'infini ne devait pas rester invariante par l'action
du groupe \`a un param\`etre en question, alors un point situ\'e dans
la r\'egion finie pourrait naturellement traverser aussi l'infini au
cours de son mouvement.

Maintenant, nous allons tout d'abord rechercher tous les groupes
r\'eels projectifs \`a un param\`etre du plan qui ont la constitution
d\'ecrite \`a l'instant. Nous les connaissons\footnote{\, D'après les
développements substantiels du Tome~III qui précèdent cette
Division~V. }, donc il ne nous sera pas difficile de trouver tous les
groupes projectifs r\'eels du plan
\`a deux ou \`a trois param\`etres qui contiennent seulement des
sous-groupes \`a un param\`etre de cette nature. Enfin, \`a partir de
l\`a, nous pourrons conclure quelle forme poss\`ede le groupe
lin\'eaire homog\`ene~\thetag{ 29}.

%%%\Fill 
%%%[[107 + 384.]]

D'apr\`es les pages~107 et~384, tout groupe projectif r\'eel \`a un
param\`etre du plan peut \^etre rapport\'e, {\em via}\, une
transformation projective r\'eelle de ce plan, \`a l'une des sept
formes suivantes:
\def\theequation{30}\begin{equation}
\label{ch21-eq30}
\left\{
\aligned
&
\quad\quad\quad
\mathfrak{p}+\mathfrak{y}\,\mathfrak{q};
\quad
\mathfrak{p}+\mathfrak{x}\,\mathfrak{q};
\quad
\mathfrak{y}\,\mathfrak{q};
\quad
\mathfrak{q};
\\
&
\mathfrak{x}\,\mathfrak{p}+\mathfrak{c}\,\mathfrak{y}\,\mathfrak{q}
\quad
{\scriptstyle{(\mathfrak{c}\,\neq\,0,\,1)}};
\quad
\mathfrak{y}\,\mathfrak{p}-\mathfrak{x}\,\mathfrak{q}
+
\mathfrak{c}\,(\mathfrak{x}\,\mathfrak{p}+\mathfrak{y}\mathfrak{q})
\quad
{\scriptstyle{(\mathfrak{c}\,\neq\,0)}};
\\
&
\quad\quad\quad\quad\quad
\mathfrak{y}\,\mathfrak{p}-\mathfrak{x}\,\mathfrak{q}.
\endaligned\right.
\end{equation}
Maintenant, quels sont ceux, parmi ces groupes \`a un param\`etre, qui
ont ici la constitution que nous demandons~?

Les cinq premiers ne l'ont certainement pas. En effet, par l'action de
chacun d'entre eux, au moins une droite r\'eelle reste au 
repos\footnote{\, Pour les
groupes numéro 1, 3, 4 et 5, une telle droite
fixe est, respectivement $\{ \mathfrak{ y} = 0\}$, 
$\{ \mathfrak{ x} = 0\}$, $\{ \mathfrak{ x} = 0\}$
et $\{ \mathfrak{ y} = 0\}$. Pour le groupe
numéro 2, c'est la droite à l'infini qui
est invariante. On vérifie cela en effectuant
le changement de carte projective
$\mathfrak{ x}' = \frac{ 1}{ \mathfrak{ x}}$, 
$\mathfrak{ y}' = \frac{ \mathfrak{ y}}{ \mathfrak{ x}}$
qui transforme comme suit les champs de vecteurs basiques:
\[
{\textstyle{\frac{\partial}{\partial\mathfrak{x}}}}
=
-\mathfrak{x}'\mathfrak{x}'
{\textstyle{\frac{\partial}{\partial\mathfrak{x}'}}}
-
\mathfrak{x}'\mathfrak{y}'
{\textstyle{\frac{\partial}{\partial\mathfrak{y}'}}}
\ \ \ \ \ \ \ \ \
\text{\rm et}
\ \ \ \ \ \ \ \ \
{\textstyle{\frac{\partial}{\partial\mathfrak{y}}}}
=
\mathfrak{x}'
{\textstyle{\frac{\partial}{\partial\mathfrak{y}'}}},
\ \ \ \ \ \ \ \ \
\text{\rm d'où}
\ \ \ \ \ \ \ \ \
{\textstyle{\frac{\partial}{\partial\mathfrak{x}}}}
+
\mathfrak{x} 
{\textstyle{\frac{\partial}{\partial\mathfrak{y}}}}
=
-\mathfrak{x}'\mathfrak{x}'
{\textstyle{\frac{\partial}{\partial\mathfrak{x}'}}}
+
(1-\mathfrak{x}'\mathfrak{y}')
{\textstyle{\frac{\partial}{\partial\mathfrak{y}'}}}, 
\]
ce qui montre que la droite $\{ \mathfrak{ x}' = 0 \}$ reste
invariante dans ce système de coordonnées.  Dans les cinq cas, en
restriction à chaque droite invariante respective, le groupe est
équivalent ou bien à $\partial_x$ (le point à l'infini est fixe) ou
bien à $x \partial_x$ (l'origine et le point à l'infini sont tous deux
fixes). }, dont les
points sont transform\'es par l'action d'un groupe \`a un param\`etre,
et pour pr\'eciser, de telle mani\`ere que sur cette droite, soit deux
points {\em r\'eels}\, s\'epar\'es, soit deux points {\em r\'eels}\,
qui co\"{\i}ncident, conservent leur position. Il d\'ecoule de l\`a
qu'un point r\'eel d'une telle droite qui ne reste pas au repos peut
certes se mouvoir librement en g\'en\'eral sur la droite, mais qu'il
n'est pas en \'etat de parcourir contin\^ument la droite de telle fa\c
con qu'il retourne finalement \`a sa position initiale sans revenir en
arri\`ere; en effet, il ne peut pas franchir les points invariants de
la droite.

Le groupe \`a un param\`etre: 
\[
\mathfrak{y}\,\mathfrak{p}-\mathfrak{x}\,\mathfrak{q}
+
\mathfrak{c}\,(\mathfrak{x}\,\mathfrak{p}+\mathfrak{y}\,\mathfrak{q})
\quad\quad\quad\quad
{\scriptstyle{(\mathfrak{c}\,\neq\,0)}}
\]
satisfait tout aussi peu notre exigence. En effet, par son action,
chaque point r\'eel non invariant d\'ecrit une spirale 
logarithmique\footnote{\, 
En coordonnées polaires, ce générateur 
s'écrit $\partial_\theta + \mathfrak{ c}\, r \, \partial_r$; 
ses courbes intégrales satisfont
$\frac{ d \theta}{ dt} = 1$ et
$\frac{ dr }{ dt} = \mathfrak{ c}\, r$, d'où
$\theta = \theta_0 + t$ et
$r = r_0 e^{ \mathfrak{ c}\, t}$: ce sont des
spirales {\em non} périodiques lorsque $\mathfrak{ c} \neq
0$. 
};%%%%%%%%%%%%%%%%%%%%%%%%%%%%%%%%%%%%%%%%%%%%%%%%%%%%%%%%%%%%%%%%%
il
parcourt alors cette spirale contin\^ument sans revenir en arri\`ere,
donc, \'evidemment, il ne retourne pas \`a sa position initiale.

Par cons\'equent, le groupe \`a un param\`etre: 
\def\theequation{31}\begin{equation}
\mathfrak{y}\,\mathfrak{p}-\mathfrak{x}\,\mathfrak{q},
\end{equation}
par l'action duquel tout point r\'eel d\'ecrit un cercle, est le seul
qui satisfait notre exigence parmi les groupes~\thetag{ 30}.

Si on consid\`ere maintenant les diff\'erents types de groupes
projectifs r\'eels du plan \`a deux et \`a trois param\`etres ({\em
voir}\, pp.~106 sq. et p.~384), alors on v\'erifie imm\'ediatement
que presque chacun d'entre eux contient un sous-groupe r\'eel \`a un
param\`etre qui poss\`ede l'une des six premi\`eres formes~\thetag{
30}, ou toutefois pour le moins, qui peut \^etre ramen\'e \`a l'une
de ces formes par une transformation projective r\'eelle. Le seul
groupe projectif r\'eel \`a deux ou \`a trois param\`etres qui ne
contient aucun sous-groupe \`a un param\`etre tel que les six
premiers de~\thetag{ 30}, est le groupe projectif r\'eel \`a trois
param\`etres:
\def\theequation{32}\begin{equation}
\mathfrak{p}
+
\mathfrak{x}\,(\mathfrak{x}\,\mathfrak{p}
+\mathfrak{y}\,\mathfrak{q}),
\quad\quad
\mathfrak{q}
+
\mathfrak{y}\,(\mathfrak{x}\,\mathfrak{p}
+\mathfrak{y}\,\mathfrak{q}),
\quad\quad
\mathfrak{y}\,\mathfrak{p}-\mathfrak{x}\,\mathfrak{q}
\end{equation}
de la conique imaginaire: $\mathfrak{ x}^2 + \mathfrak{ y}^2 + 1 = 0$.

Il d\'ecoule de l\`a que, sous les hypoth\`eses pos\'ees, le groupe
$\mathfrak{ g}$ d\'efini \`a la page~\pageref{462} est toujours \`a trois
param\`etres et qu'il peut \^etre ramen\'e, {\em via}\, une transformation
projective r\'eelle, \`a la forme~\thetag{ 32}.

En revenant maintenant au groupe lin\'eaire homog\`ene~\thetag{ 29}, on
v\'erifie imm\'ediatement (\cf p.~110) que par une transformation
lin\'eaire homog\`ene r\'eelle, celui-ci re\c coit la forme:
\[
\aligned
&
x_\mu'\,p_\nu'-x_\nu'\,p_\mu'
+
\alpha_{\mu\nu}\,
(x_1'\,p_1'+x_2'\,p_2'+x_3'\,p_3')
\\
&
\quad\quad\quad
{\scriptstyle{(\mu,\,\nu\,=\,1,\,2,\,3;\,\,\,\,
\alpha_{\mu\nu}\,+\,\alpha_{\nu\mu}\,=\,0)}}.
\endaligned
\]
En calculant les crochets de ces transformations par paires, on
d\'eduit enfin que tous les $\alpha_{ \mu \nu}$ 
s'annulent\footnote{\,
%%%%%%%%%%%%%%%%%%%%%%%-------DEBUT--------%%%%%%%%%%%%%%%%%%%%%%%%%%%
\label{verification-crochets-rotations}
Posons $R_{ \mu \nu} ' := x_\mu' p_\nu' - x_\nu' p_\mu'$ (générateur
infinitésimal des rotations autour de la droite $\{ x_\mu' = x_\nu' =
0\}$) et $D' := x_1' p_1' + x_2 ' p_2 ' + x_3 ' p_3'$ (générateur
infinitésimal des homothéties de centre l'origine).  On a les
relations de commutation $\big[ R_{ 12}', \, R_{ 23}' \big] =
R_{13}'$, $\big[ R_{ 12}', \, R_{ 13}' \big] = - R_{ 23}'$, $\big[ R_{
13}', \, R_{ 23}' \big] = - R_{ 12}'$ et aussi $\big[ R_{ \mu \nu}' ,
\, D' \big] = 0$.  Il en découle que le crochet de Lie:
\[
\aligned
\big[R_{12}'+\alpha_{12}D',\,
R_{23}'+\alpha_{23}\,D'\big]
&
=
\big[R_{12}',\,R_{23}'\big]
+
\alpha_{23}\big[R_{12}',\,D'\big]
+
\alpha_{12}\,\big[D',\,R_{23}'\big]
+
\alpha_{12}\alpha_{23}\big[D',\,D'\big]
\\
&
=
R_{13}'
\\
&
=
\lambda
\big(R_{13}'+\alpha_{13}D'
\big)
\endaligned
\]
qui doit être combinaison linéaire des trois générateurs
infinitésimaux $R_{ \mu\nu'} +
\alpha_{ \mu\nu} D'$ ne
peut l'être que lorsque $\alpha_{ 13} = 0$, avec $\lambda = 1$. De
même, $\alpha_{ 12} = 0$ et $\alpha_{ 23} = 0$. 
}. %%%%%%%%%%%%%%%%%%%%%%%%-----FIN-----%%%%%%%%%%%%%%%%%%%%%%%%%%%%%%%

Avec cela, nous sommes parvenus au r\'esultat que, sous les
hypoth\`eses pos\'ees, le groupe lin\'eaire homog\`ene r\'eel~\thetag{
29} peut toujours \^etre suppos\'e de la forme:
\[
x_1'\,p_2'-x_2'\,p_1',
\quad\quad
x_2'\,p_3'-x_3'\,p_2',
\quad\quad
x_3'\,p_1'-x_1'\,p_3'.
\]
Mais il d\'ecoule de l\`a que chaque groupe r\'eel de transformations
qui satisfait nos Axiomes~III et~IV dans le voisinage de chaque point
r\'eel: $x_1^0, x_2^0, x_3^0$ en position g\'en\'erale contient trois
transformations infinit\'esimales du premier ordre en les $x_\nu -
x_\nu^0$, que l'on peut s'imaginer rapport\'ees \`a la forme:
\[
(x_\mu-x_\mu^0)\,p_\nu
-
(x_\nu-x_\nu^0)\,p_\mu
+
\cdots
\quad\quad\quad
{\scriptstyle{(\mu,\,\nu\,=1,\,2,\,3;\,\,\,\,\mu\,<\,\nu)}},
\]
alors qu'au contraire, il ne contient pas d'autres transformations
infinit\'esimales du premier ordre, et en particulier aucune de la
forme: 
\[
(x_1-x_1^0)\,p_1
+
(x_2-x_2^0)\,p_2
+
(x_3-x_3^0)\,p_3
+
\cdots.
\]
Mais avec cela, la d\'etermination de tous ces groupes est ramen\'ee
au probl\`eme qui a d\'ej\`a \'et\'e r\'esolu au \S~84 (p.~365
sq.). 
%%%\Fill 
%%%[[Explorer ces pages.]]
Ainsi, nous pouvons tout d'abord en conclure que les groupes
concern\'es sont finis, et pour pr\'eciser, qu'ils comportent six
param\`etres. En outre, on obtient qu'ils sont semblables, {\em via}\,
une transformation ponctuelle r\'eelle, soit au groupe des mouvements
euclidiens, soit au groupe projectif r\'eel d'une des deux surfaces du
second degr\'e:
\[
x_1^2+x_2^2+x_3^2+1
=
0,
\quad\quad\quad
x_1^2+x_2^2+x_3^2-1
=
0.
\]

Par cons\'equent, chaque groupe r\'eel de l'espace trois fois \'etendu
qui satisfait nos Axiomes~III et~IV peut \^etre transform\'e, {\em
via}\, une transformation ponctuelle r\'eelle, soit en le groupe des
mouvements euclidiens, soit en l'un des deux groupes de mouvements
non-euclidiens. En d'autres termes: nos Axiomes~I \dots\, IV suffisent
effectivement \`a caract\'eriser ces trois sortes de mouvements.

\HEAD{Critique des recherches helmholtziennes.}{
Division\,\,V.\,\,\,Chapitre\,\,21.\,\,\,\S\,\,96.}

\sectiondritterV{\sf\S\,\,96.
\\
Quelles conclusions peut-on tirer des axiomes helmholtziens?}
\label{S-96}
\setcounter{footnote}{0}

Nous voulons \`a pr\'esent faire compl\`etement abstraction des
recherches que Monsieur de Helmholtz lui-m\^eme a engag\'ees sur la
base de ses axiomes et nous voulons \'etablir, comme nous l'avons
d\'ej\`a annonc\'e a la page~\pageref{438}, quelles conclusions
peuvent \^etre tir\'ees des axiomes helmholtziens en eux-m\^emes.
Nous nous limitons donc \`a nouveau \`a l'espace ordinaire trois fois
\'etendu. \label{465}

\bigskip

Dans le \S~93, nous avons indiqu\'e comment les axiomes helmholtziens
peuvent \^etre formul\'es, quand on les interpr\`ete de la fa\c con
qui a \'et\'e expos\'ee au \S~92 et lorsqu'on emploie de surcro\^{\i}t
les concepts et les mani\`eres de s'exprimer de la th\'eorie des
groupes. Il nous reste maintenant \`a voir quels sont les groupes de
l'espace trois fois \'etendu qui satisfont les exigences pos\'ees dans
le \S~93.

Chacun des groupes requis est r\'eel, fini et continu, et relativement
\`a son action, deux points ont un et un seul invariant, tandis qu'un
nombre de points sup\'erieur \`a deux n'a pas d'invariant
essentiel. Mais puisque, d'apr\`es le Th\'eor\`eme~37
p.~\pageref{Theorem-37-p-433}, 
chaque groupe de cette esp\`ece est semblable, {\em
via}\, une transformation ponctuelle r\'eelle, \`a l'un des groupes
qui sont rassembl\'es \`a la page~\pageref{433-sq}, il ne nous reste
alors qu'\`a exclure, parmi ces groupes-l\`a, ceux qui ne
correspondent pas aux exigences restantes du \S~93.

Ces exigences restantes reviennent toutes \`a ce que nos groupes, \`a
l'int\'erieur d'une certaine r\'egion finie de l'espace trois fois
\'etendu, doivent poss\'eder certaines propri\'et\'es. Puisque
l'origine des coordonn\'ees est un point en position g\'en\'erale pour
tous les groupes de la page~\pageref{433-sq}, nous pouvons clairement
admettre que la r\'egion en question est constitu\'ee de tous les
points r\'eels qui se trouvent dans un certain voisinage de l'origine
des coordonn\'ees. Par cons\'equent, nous devons seulement rechercher
encore, parmi les groupes de la page~\pageref{433-sq}, quels sont ceux
qui poss\`edent les propri\'et\'es suivantes:

\medskip

{\small\sf Premi\`erement}, 
l'invariant: $J \big( x_1, y_1, z_1 ; \, x_2, y_2,
z_2 \big)$ de deux points doit \^etre constitu\'e de telle sorte que la
relation:
\def\theequation{33}\begin{equation}
J
\big(
x_1,y_1,z_1:\,
x_2,y_2,z_2
\big)
=
J
\big(
x_1^0,y_1^0,z_1^0;\,
x_2^0,y_2^0,z_2^0
\big),
\end{equation}
pour chaque paire de syst\`emes de valeurs distincts l'un de l'autre:
$x_1^0, y_1^0, z_1^0$ et $x_2^0, y_2^0, z_2^0$, fournit une v\'eritable
\'equation entre $x_1, y_1, z_1$, $x_2, y_2, z_2$, 
dans un certain voisinage de l'origine des coordonn\'ees.

\medskip

{\small\sf Deuxi\`emement}, 
apr\`es fixation de l'origine des coordonn\'ees,
chaque point quelconque $x_0, y_0, z_0$ qui est situ\'e dans un certain
voisinage de l'origine des coordonn\'ees, doit pouvoir encore se
transformer en tous les points $x, y, z$ de ce voisinage, qui
satisfont l'\'equation:
\def\theequation{34}\begin{equation}
J
\big(
x,y,z;\,
0,0,0
\big)
=
J
\big(
x_0,y_0,z_0;\,
0,0,0
\big),
\end{equation}
et \`a travers laquelle un transfert continu entre tous les points 
de\renewcommand{\thefootnote}{\fnsymbol{footnote}}
cette sorte est g\'en\'eralement possible\footnote[1]{\, \`A vrai dire,
dans le \S~93, on demande encore plus, mais nous verrons que cette
exigence suffit d\'ej\`a compl\`etement.}.

\renewcommand{\thefootnote}{\arabic{footnote}}

%%%\Fill 
%%%[[une seule orbite.]]

%%%\Fill 
%%%[[{\bf SUBTIL:} la transitivit\'e en plus dans D) \'etait cach\'ee dans l'existence d'un seul invariant $J$, et dans la première phrase qui demande la libre mobilité tant que les points ne sont pas liés avec d'autres points,  donc la libre mobilit\'e demand\'ee en Axiome D) \S~93 est inutile.]]

\medskip

{\small\sf Troisi\`emement}, 
l'axiome de monodromie doit \^etre satisfait.

\medskip

Nous recherchons tout d'abord, pour chaque groupe individuel de la
page~\pageref{433-sq}, si la premi\`ere de ces exigences est
satisfaite, et lorsqu'elle l'est, nous nous demandons alors si la
deuxi\`eme est satisfaite; et nous ne prendrons en consid\'eration
l'axiome de monodromie que dans les derni\`eres lignes. Comme par
ailleurs le groupe des mouvements euclidiens et les deux groupes de
mouvements non-euclidiens satisfont \'evidemment toutes nos exigences,
d\`es le d\'ebut, nous les laissons simplement de c\^ot\'e.

\bigskip

Le fait que les deux groupes:
\def\theequation{35}\begin{equation}
p,\quad
q,\quad
r,\quad
xq-yp,\quad
yr+zq,\quad
zp+xr
\end{equation}
et:
\def\theequation{36}\begin{equation}
p-xU,\quad
q-yU,\quad
r+zU,\quad
xq-yp,\quad
yr+zq,\quad
zp+xr
\end{equation}
satisfont notre premi\`ere exigence est \`a port\'ee de main, puisque par
exemple, pour le premier d'entre eux, l'\'equation~\thetag{33} s'\'enonce
de la mani\`ere suivante:
\[
\left\{
\aligned
&
\quad
(x_2-x_1)^2+(y_2-y_1)^2-(z_2-z_1)^2
=
\\
&
=
(x_2^0-x_1^0)^2+(y_2^0-y_1^0)^2-(z_2^0-z_1^0)^2,
\endaligned\right.
\]
et ceci est toujours une v\'eritable \'equation entre $x_1, y_1, z_1$,
$x_2, y_2, z_2$. N\'eanmoins, ces deux groupes ne satisfont pas la
deuxi\`eme exigence. En effet, l'\'equation~\thetag{ 34} correspondante
a, pour les deux groupes, la forme:
\def\theequation{37}\begin{equation}
x^2+y^2-z^2
=
x_0^2+y_0^2-z_0^2.
\end{equation}
Et maintenant, apr\`es fixation de l'origine des coordonn\'ees, le
point $x_0, y_0, z_0$ peut \`a vrai dire en g\'en\'eral \^etre
transform\'e encore en tous les points qui satisfont
l'\'equation~\thetag{ 37}; mais pour tous les points $x_0, y_0, z_0$
qui se trouvent sur la conique du second degr\'e: $x_0^2 + y_0^2 -
z_0^2 = 0$, cela n'est plus valable; effectivement, pour ces points,
l'\'equation~\thetag{ 37} re\c coit la forme:
\[
x^2+y^2-z^2
=
0,
\]
et il est clair que chaque point de cette esp\`ece peut se transformer
en tous les points qui satisfont cette \'equation, \`a l'exception de
l'origine des coordonn\'ees, car, sous les hypoth\`eses pos\'ees, elle
reste au repos.

Nous voyons donc que les groupes~\thetag{35} et~\thetag{ 36} sont
d\'ej\`a exclus, sans qu'on utilise l'axiome de monodromie.

\bigskip

Pour le groupe:
\def\theequation{38}\begin{equation}
p,\quad
q,\quad
xp+r,\quad
yq+cr,\quad
x^2p+2xr,\quad
y^2q+2cyr\quad\quad
{\scriptstyle{(c\,\neq\,0)}},
\end{equation}
on peut \'ecrire l'\'equation~\thetag{33} 
correspondante
sous la forme:
\[
\frac{(x_2-x_1)(y_2-y_1)^c}{e^{\frac{1}{2}(z_1+z_2)}}
=
\frac{(x_2^0-x_1^0)(y_2^0-y_1^0)^c}{e^{\frac{1}{2}(z_1^0+z_2^0)}}.
\]
Si $c$ est n\'egatif, alors cette \'equation n'a plus de sens aussit\^ot
que l'on pose $x_2^0 = x_1^0$ et $y_2^0 = y_1^0$;

%%%\Mathematiques: 
%%%[[{\sc Ist $c$ negativ, so wird diese Gleichung zu einer Identit\"at, sobald man $x_2^0 = x_1^0$ und $y_2^0 = y_1^0$ setz}~??]]

\noindent
si au contraire, $c$ est positif, elle fournit alors, pour toute paire
de syst\`emes de valeurs $x_1^0, y_1^0, z_1^0$ et $x_2^0, y_2^0,
z_2^0$ distincts l'un de l'autre, une v\'eritable \'equation entre
$x_1, y_1, z_1$, $x_2, y_2, z_2$. Le groupe~\thetag{ 38} satisfait
donc notre premi\`ere exigence seulement lorsque $c$ est positif.

Mais il ne satisfait pas notre deuxi\`eme exigence, m\^eme dans le
cas: $c > 0$. En effet, pour notre groupe, l'\'equation~\thetag{ 34}
re\c coit la forme:
\[
\frac{x\cdot y^c}{e^{\frac{1}{2}\,z}}
=
\frac{x_0\cdot y_0^c}{e^{\frac{1}{2}\,z_0}},
\]
et apr\`es fixation de l'origine des coordonn\'ees, le point $x_0, y_0,
z_0$ devrait pouvoir occuper toutes les positions qui satisfont cette
\'equation dans un certain voisinage de l'origine des coordonn\'ees.
C'est \`a vrai dire effectivement le cas pour un point $x_0, y_0, z_0$
en position g\'en\'erale. Cependant, en choisissant par exemple $x_0 =
y_0 = 0$, le point $x_0, y_0, z_0$ devrait encore pouvoir \^etre
transform\'e en tous les points se trouvant dans un certain
voisinage de l'origine des coordonn\'ees qui satisfont l'\'equation:
$x y^c = 0$, mais cela est impossible, car apr\`es fixation de l'origine
des coordonn\'ees, tous les points de la droite: $x = y = 0$ restent en
g\'en\'eral au repos ({\em voir}\, p.~\pageref{426}).

%%%\Fill 
%%%[[on l'a d\'ej\`a vu~!]]

Par cons\'equent, le groupe~\thetag{ 38} est
lui aussi exclu, sans qu'on soit oblig\'e de
se pr\'eoccuper de l'axiome de monodromie.

\bigskip

Qui plus est, notre premi\`ere exigence n'est m\^eme encore jamais
satisfaite par le groupe:
\def\theequation{39}\begin{equation}
\left\{
\aligned
&
p,\quad
q,\quad
xp+yq+ar,\quad
yp-xq+r,\quad
\\
&
(x^2-y^2)\,p+2xyq+2(ax-y)\,r
\\
&
2xyp+(y^2-x^2)\,q+2(x+ay)\,r,
\endaligned\right.
\end{equation}
puisque l'\'equation~\thetag{ 33} correspondante ne fournit en effet
pas toujours une v\'eritable \'equation entre $x_1, y_1, z_1$, $x_2,
y_2, z_2$, lorsqu'on pose: $x_2^0 = x_1^0$, $y_2^0 = y_2^0$.

%%%\Mathematiques
%%%[[Retrouver l'expression de l'invariant.]]

\bigskip

Il en va autrement pour le groupe:
\def\theequation{40}\begin{equation}
\left\{
\aligned
&
p,\quad
q,\quad
xp+yq+r,\quad
yp-xq
\\
&
(x^2-y^2)\,p+2xyq+2xr
\\
&
2xyp+(y^2-x^2)\,q+2yr. 
\endaligned\right.
\end{equation}
Pour celui-ci, l'\'equation~\thetag{ 33} s'\'ecrit en effet:
\[
\aligned
&
\quad \
\big\{
(x_2-x_1)^2+(y_2-y_1)^2
\big\}\,
e^{-(z_1+z_2)}
=
\\
&
=
\big\{
(x_2^0-x_1^0)^2+(y_2^0-y_1^0)^2
\big\}\,
e^{-(z_1^0+z_2^0)},
\endaligned
\]
et ceci constitue toujours une \'equation v\'eritable entre $x_1, y_1,
z_1$, $x_2, y_2, z_2$. Notre premi\`ere exigence est donc satisfaite;
mais la deuxi\`eme ne l'est pas.

En effet, l'\'equation~\thetag{ 34} s'\'ecrit maintenant:
\[
(x^2+y^2)\,e^{-z}
=
(x_0^2+y_0^2)\,e^{-z_0},
\]
et si notre deuxi\`eme exigence devait \^etre satisfaite, apr\`es
fixation de l'origine des coordonn\'ees, alors par exemple, chaque
point $x_0, y_0, z_0$ pour lequel $x_0 = y_0 = 0$ devrait pouvoir
\^etre transform\'e encore en tous les autres points qui satisfont
l'\'equation: $x^2 + y^2 = 0$ ou encore, puisqu'il s'agit de
quantit\'es r\'eelles, les deux \'equations: $x = y = 0$. Mais cela
est impossible, car, en m\^eme temps que l'origine des coordonn\'ees,
tous les points de la droite: $x = y = 0$ restent g\'en\'eralement au
repos.

%%%\Fill 
%%%[[V\'erifier.]]

Par cons\'equent, les groupes~\thetag{ 39} et~\thetag{ 40} 
n'entrent pas non plus en ligne de compte.

\bigskip

Les deux groupes: 
\def\theequation{41}\begin{equation}
\left\{
\aligned
&
p,\quad
q,\quad
xq+r,\quad
xp+yq+cr,\quad
x^2q+2xr,
\\
&
\quad\quad\quad\quad
x^2p+2xyq+2(y+cx)\,r
\endaligned\right.
\end{equation}
et: 
\def\theequation{42}\begin{equation}
p,\quad
q,\quad
r,\quad
2xp+yq,\quad
xq+yr,\quad
x^2p+xyq+{\textstyle{\frac{1}{2}}}\,y^2r
\end{equation}
ne satisfont pas non plus notre premi\`ere exigence; dans les deux
cas en effet, l'\'equation~\thetag{ 33} correspondante ne fournit
aucune \'equation v\'eritable entre $x_1, y_1, z_1$, $x_2, y_2, z_2$,
aussit\^ot que l'on pose $x_2^0 = x_1^0$, $y_2^0 = y_1^0$.

\bigskip

Enfin, le groupe: 
\def\theequation{43}\begin{equation}
r,\quad
p-yr,\quad
q+xr,\quad
xq,\quad
xp-yq,\quad
yp
\end{equation}
satisfait certes notre premi\`ere exigence, mais pas la deuxi\`eme.
L'\'equation~\thetag{ 34} correspondante s'\'ecrit en effet: $z =
z_0$, et par cons\'equent, apr\`es fixation de l'origine des
coordonn\'ees, chaque autre point $x_0, y_0, z_0$ devrait pouvoir se
transformer encore en tous les autres points qui se trouvent sur le
plan: $z = z_0$. Les points pour lesquels $x_0 = y_0 = 0$ montrent
que cela n'est pas le cas, car ils restent tous au repos apr\`es
fixation de l'origine des coordonn\'ees.

\bigskip

Par ce qui pr\'ec\`ede, on a d\'emontr\'e que {\em le groupe des
mouvements euclidiens et les deux groupes de mouvements non-euclidiens
sont les seuls groupes qui satisfont les exigences \'enonc\'ees au
\S~93; ce sont les seuls groupes de cette esp\`ece, m\^eme lorsqu'on
laisse compl\`etement de c\^ot\'e l'axiome de monodromie}.

Nous voyons donc que les axiomes helmholtziens sont suffisants pour la
caract\'erisation des mouvements euclidiens et non-euclidiens, {\em si
on les interpr\`ete de la m\^eme mani\`ere qu'aux \S~92 et~93}; mais
nous voyons en m\^eme temps qu'en prenant pour base cette
interpr\'etation, l'axiome de monodromie est de toute fa\c con
superflu. Mais bien qu'il se soit av\'er\'e que l'on peut se passer de
l'un des axiomes de Helmholtz, se pr\'esente juste apr\`es la question
de savoir si on ne peut pas se passer aussi encore d'autres parties de
ces axiomes. Dans le Chapitre~23 nous montrerons qu'il faut r\'epondre
par l'affirmative \`a cette question.

Mais il y a encore un point qui m\'erite d'\^etre pris en
consid\'eration. Dans les \S~92 et~93, nous avons interpr\'et\'e les
axiomes helmholtziens d'une certaine mani\`ere. Nous sommes loin
d'affirmer que cette interpr\'etation qui est la n\^otre soit la seule
possible, et nous dirons plut\^ot sans d\'etour que les axiomes en
question autorisent des interpr\'etations diverses. Cela r\'eside dans
la nature des choses: Monsieur de Helmholtz, qui n'avait pas \`a sa
disposition les concepts et la langue sp\'ecialis\'ee instruite de la
th\'eorie des groupes, a d\^u s'aider de longues circonlocutions par
lesquelles, cependant, il n'a pas pu prendre en compte toutes les
exceptions imaginables; ainsi, il est tout \`a fait compr\'ehensible
que ses axiomes apparaissent avoir d'autant plus de sens vari\'es
qu'on les examine plus soigneusement. Maintenant, bien que nous
n'affirmions pas avoir donn\'e la seule interpr\'etation possible des
axiomes helmholtziens, nous croyons toutefois que nous leur avons
attribu\'e l'interpr\'etation la plus int\'eressante qui soit, celle
qui est en accord avec leur teneur; nous y avons m\^eme d\'ej\`a trop
introduit maint lecteur.

Nous ne voulons pas rentrer plus avant dans la consid\'eration des
diverses conceptions que l'on peut avoir au sujet des axiomes
helmholtziens; toutefois, nous voulons encore en discuter
bri\`evement.

En effet, la mani\`ere dont on peut interpr\'eter ces axiomes d\'epend
pour l'essentiel du fait de savoir si on lit en eux qu'ils doivent
\^etre valables sans exception \`a l'int\'erieur d'une certaine
r\'egion, ou si l'on admet qu'ils doivent \^etre valables seulement en
g\'en\'eral, donc pour des points qui sont mutuellement en position
g\'en\'erale. Dans les paragraphes~92 et~93, nous nous sommes
d\'ecid\'es pour la premi\`ere conception, mais la deuxi\`eme
poss\`ede aussi une certaine l\'egitimit\'e, au moins en ce qui
concerne le troisi\`eme axiome helmholtzien, l'axiome sur la libre
mobilit\'e. Afin de montrer quelle influence cette conception a sur
la port\'ee des axiomes helmholtziens, nous voulons finalement
rechercher encore ce qui advient lorsque les axiomes helmholtziens
sont satisfaits seulement par des points qui sont mutuellement en
position g\'en\'erale. Nous allons voir que ces axiomes ne suffisent
alors pas \`a caract\'eriser les mouvements euclidiens et
non-euclidiens, m\^eme si on admet en plus l'axiome de monodromie.

Si un groupe r\'eel doit satisfaire les axiomes helmholtziens
seulement en ce qui concerne des points qui sont mutuellement en
position g\'en\'erale, alors il doit, quand on fait d'abord
abstraction de l'axiome de monodromie, poss\'eder les propri\'et\'es
suivantes: si un point en position g\'en\'erale est fix\'e, tout
autre point en position g\'en\'erale se meut librement sur une
surface; si deux points qui sont mutuellement en position
g\'en\'erale sont fix\'es, tout autre troisi\`eme point en position
g\'en\'erale se meut librement sur une courbe; si trois points qui
sont mutuellement en position g\'en\'erale sont fix\'es, tous les
points restent au repos. Il est clair que les groupes rassembl\'es \`a
la page~\pageref{433-sq} satisfont toutes ces exigences et qu'ils sont
les seuls \`a les satisfaire.

Si nous admettons en plus l'axiome de monodromie, alors tout un nombre
des groupes dont il s'agit seront exclus, parce qu'ils ne satisfont
pas l'axiome de monodromie; mais il reste alors non seulement les
mouvements euclidiens et non-euclidiens, mais en outre encore un
groupe, \`a savoir le groupe~{\bf 7} \`a la page~\pageref{433-sq}. On
a d\'emontr\'e aux pages~\pageref{457} sq. qu'il satisfait compl\`etement
l'axiome de monodromie helmholtzien. Par cons\'equent, sous les
hypoth\`eses pos\'ees, l'axiome de monodromie ne suffit pas \`a
exclure aussi ce groupe.

En d\'efinitive \deutsch{Alles
in Allem}, nous pouvons \'enoncer bri\`evement le r\'esultat
du pr\'esent paragraphe ainsi:

\medskip

{\em Les mouvements euclidiens et non-euclidiens sont compl\`etement
caract\'eris\'es par les axiomes helmholtziens, {\small\sf lorsqu'on
interpr\`ete ces axiomes de telle sorte qu'ils doivent \^etre
satisfaits sans exception par tous les points \`a l'int\'erieur d'une
certaine r\'egion; et quand on les interp\`ete de cette mani\`ere,
l'axiome de monodromie est alors superflu}. Si au contraire on demande
seulement que les axiomes helmholtziens soient satisfaits par des
points qui sont mutuellement en position g\'en\'erale, alors ces
axiomes ne sont pas suffisants pour caract\'eriser les mouvements
euclidiens et non-euclidiens, et il ne suffisent d'ailleurs m\^eme pas
lorsqu'on admet en plus l'axiome de monodromie.}

\bigskip

Nous avons vu que les axiomes helmholtziens suffisent certes pour la
caract\'erisation des mouvements euclidiens et non euclidiens,
lorsqu'on les interp\`ete d'une certaine mani\`ere, mais qu'ils
contiennent des \'el\'ements superflus quand on les interpr\`ete
ainsi, et par cons\'equent, ils ne repr\'esentent pas une v\'eritable
solution du probl\`eme de Riemann-Helmholtz. Dans les deux chapitres
suivants, nous allons maintenant traiter le probl\`eme de
Riemann-Helmholtz d'un point de vue nouveau, et nous allons le
r\'esoudre de deux fa\c cons diff\'erentes, \`a savoir, dans le
Chapitre~22, en \'etablissant certains axiomes qui se r\'ef\`erent \`a
des points infiniment voisins, et dans le Chapitre~23, {\em via}\, des
axiomes qui se r\'ef\`erent \`a des points finiment \'eloign\'es.

%%%%%%%%%%%%%%%%%%%%%%%%%%%%%%%%%%%%%%%%%%%%%%%%%%%%%%%%%%%%%%%%%%%%%

\newpage

% 27   :   471--497

\setcounter{footnote}{0}

$\:$
\bigskip\bigskip\bigskip

\centerline{\Large Chapitre~22.}
\label{Chapitre-22}
\thispagestyle{empty}

\bigskip

\noindent
{\large\bf
Premi\`ere solution du probl\`eme de Riemann-Helmholtz.
}

\bigskip\medskip

Le probl\`eme de Riemann-Helmholtz, comme nous l'avons formul\'e \`a
la page~\pageref{397}, demande d'indiquer des propri\'et\'es qui sont
communes \`a la famille de mouvements euclidiens et aux deux familles
de mouvements non-euclidiens et gr\^ace auxquelles ces trois familles
se distinguent de toutes les autres familles possibles de mouvements.

Comme ce probl\`eme est en fait ind\'etermin\'e, nous le sp\'ecifions
plus pr\'ecis\'ement d\`es le d\'ebut par l'exigence que l'ensemble
des propri\'et\'es \`a indiquer doivent non seulement \^etre
suffisantes pour caract\'eriser nos trois familles, mais encore \^etre
n\'ecessaires, c'est-\`a-dire que ces trois familles-l\`a de
mouvements ne doivent pas seulement \^etre caract\'eris\'ees par une
fraction des propri\'et\'es en question. Mais cependant, m\^eme sous
cette version d\'etermin\'ee, le probl\`eme admet encore des solutions
tr\`es diverses, puisque dans le choix des propri\'et\'es requises, on
n'est par ailleurs soumis \`a aucune limitation
suppl\'ementaire. \label{471}

\renewcommand{\thefootnote}{\fnsymbol{footnote}}
Parmi les propri\'et\'es qui sont communes \`a nos trois familles de
mouvements, il s'en trouve une qui nous touche de pr\`es, \`a savoir
que chacune des trois familles constitue un groupe r\'eel continu,
dont les transformations sont inverses l'une de l'autre par paires.
Nous supposons alors depuis le d\'ebut cette propri\'et\'e comme
quelque chose de donn\'e\footnote[1]{\, Nous ne voulons pas passer
sous silence le fait qu'il est tout \`a fait possible d'\'etablir
d'autres exigences dont d\'ecoulent cette propri\'et\'e de type
<<\,groupe\,>>. Mais on peut bien s\^ur soutenir qu'il est conforme
\`a la nature des choses d'\'etablir d\`es le d\'ebut ces exigences
quant \`a la nature des groupes. Dans le concept de mouvement se
trouve en effet le fait qu'on peut accomplir deux mouvements l'un
apr\`es l'autre, et il d\'ecoule de l\`a imm\'ediatement que
l'ensemble des transformations qui repr\'esentent des changements de
lieu possibles par mouvement forme un groupe. Nous ajoutons
seulement encore l'hypoth\`ese que le groupe est continu et qu'il est
constitu\'e de transformations inverses l'une de l'autre. Du reste,
dans le Chapitre~21, nous avons d\'ej\`a donn\'e une solution qui ne
fait pas usage de cette derni\`ere propri\'et\'e des groupes.
%%%\Mathematiques
%%%[{\sc Gruppeneigenschaft}: au singulier.
%%%Affirmation incompr\'ehensible~? Vieille conjecture
%%%de Lie~?]
},
et nous devons maintenant simplement demander encore:
\renewcommand{\thefootnote}{\arabic{footnote}} {\em Par quelles
propri\'et\'es peut-on distinguer le groupe des mouvements euclidiens
et les deux groupes de mouvements non-euclidiens parmi tous les
groupes r\'eels continus poss\'edant des transformations inverses par
paire}~? Pour all\'eger le discours, nous supprimerons toujours dans
la suite les mots: <<\,poss\'edant des transformations inverses par
paires\,>>, et ces mots seront ainsi \`a chaque fois sous-entendus.

Dans le pr\'esent chapitre nous allons r\'epondre \`a cette question
d'une mani\`ere telle que nous nous limiterons aux propri\'et\'es qui
se r\'ef\`erent \`a des points infiniment voisins. Nous allons
introduire un nouveau concept, \`a savoir le concept de {\sl libre
mobilit\'e dans l'infinit\'esimal}, et nous montrerons que chaque
groupe r\'eel continu agissant sur
un espace dont la dimension est strictement sup\'erieure
\`a deux qui poss\`ede la libre mobilit\'e dans l'infinit\'esimal
en un point de position g\'en\'erale est semblable, {\em via}\, une
transformation ponctuelle r\'eelle de l'espace en question, soit au
groupe des mouvements euclidiens, soit \`a l'un des deux groupes de
mouvements non-euclidiens de l'espace
en question. Par cons\'equent, il sera alors
d\'emontr\'e que les propri\'et\'es que nous aurons r\'eunies dans le
concept de libre mobilit\'e dans l'infnit\'esimal,
\renewcommand{\thefootnote}{\fnsymbol{footnote}} seront des
propri\'et\'es caract\'eristiques telles que celles que nous
recherchons\footnote[1]{\, Ce qui suit est une refonte de l'\'etude
que Lie a conduite aux pages~284 \`a 321 du {\em Leipziger Berichte}
de 1890.}.
\renewcommand{\thefootnote}{\arabic{footnote}}

Dans le \S~97 nous prendrons tout d'abord en consid\'eration le plan,
dans lequel, \`a vrai dire, le concept de libre mobilit\'e dans
l'infinit\'esimal ne suffit encore pas pour caract\'eriser les
mouvements euclidiens et non-euclidiens. Dans le \S~98, nous
traiterons de l'espace ordinaire trois fois \'etendu, et dans le \S~99,
de l'espace \`a un nombre quelconque de dimensions. Dans le \S~100
enfin, nous \'etudierons bri\`evement les liens qui existent entre nos
d\'eveloppements et ceux de Riemann.

\HEAD{Première solution au problème de Riemann-Helmholtz.}{
Division\,\,V.\,\,\,Chapitre\,\,22.\,\,\,\S\,\,97.}

\sectiondritterV{\sf\S\,\,97.}
\label{S-97}
\setcounter{footnote}{0}

Soit $\mathfrak{ G}$ un groupe r\'eel continu d'un espace quelconque
et soit $F$ une figure quelconque dans cet espace. Alors, si
$\mathfrak{ G}$ contient un sous-groupe continu \`a au moins un
param\`etre par l'action duquel la figure $F$ reste invariante, nous
voulons dire qu'apr\`es fixation de $F$, un mouvement continu est
encore possible par le groupe $\mathfrak{ G}$. Si au contraire le
plus grand sous-groupe continu de $\mathfrak{ G}$ qui laisse $F$
invariante est seulement constitu\'e de la transformation identique,
alors nous disons qu'apr\`es fixation de $F$, plus aucun mouvement
continu n'est possible.

En utilisant une telle mani\`ere de s'exprimer, nous d\'efinissons
maintenant ce que nous entendons par <<\,libre mobilit\'e dans
l'infinit\'esimal\,>> pour un plan.

\medskip

{\bf D\'efinition.} {\em Un groupe r\'eel continu du plan poss\`ede la
{\small\sf libre mobilit\'e dans l'infinit\'esimal} en un point r\'eel $P$
lorsque les conditions suivantes sont remplies: Apr\`es fixation de $P$,
un mouvement continu doit encore \^etre possible; au contraire,
aucun mouvement ne doit plus \^etre possible, aussit\^ot que l'on a
fix\'e, outre $P$, un \'el\'ement lin\'eaire quelconque passant par
$P$.}

\medskip

Nous cherchons maintenant tous les groupes r\'eels continus du plan
qui poss\`edent la libre mobilit\'e dans l'infinit\'esimal en un point
r\'eel {\em de position g\'en\'erale}.

\bigskip

Soit $G$ un groupe ayant les qualit\'es requises et soit $P$ un point
r\'eel de position g\'en\'erale, en lequel $G$ poss\`ede la libre
mobilit\'e dans l'infinit\'esimal. Si ensuite $g$ est le plus grand
sous-groupe continu de $G$ par lequel le point $P$ reste au repos,
chaque \'el\'ement lin\'eaire r\'eel passant par $P$ doit pouvoir tourner
contin\^ument autour de $P$ par l'action de $g$; car si tel n'\'etait
pas le cas, $g$ laisserait au moins un de ces \'el\'ements lin\'eaires
au repos, et $g$ serait alors le plus grand sous-groupe continu de $G$
par lequel le point $P$ et l'\'el\'ement lin\'eaire en question restent
invariants, et $g$ devrait par cons\'equent se r\'eduire \`a la
transformation identique, alors que nous avons toutefois suppos\'e
qu'apr\`es fixation de $P$, un mouvement continu \'etait encore possible.

\renewcommand{\thefootnote}{\fnsymbol{footnote}}
Il d\'ecoule de l\`a imm\'ediatement que $G$ n'est pas seulement
transitif, mais qu'il est m\^eme aussi r\'eel-primitif\,\footnote[1]{\, Par
souci de bri\`evet\'e, nous appelons {\sl r\'eel-primitif} tout groupe
r\'eel qui est primitif, d\`es qu'on se restreint au domaine r\'eel,
bien qu'il puisse \^etre imprimitif dans le domaine complexe 
(\cf~p.~363).
}.
\renewcommand{\thefootnote}{\arabic{footnote}}
Si en effet il \'etait intransitif ou s'il \'etait r\'eel-imprimitif,
il laisserait invariante une famille r\'eelle de courbes: $\varphi (
x, y) = \text{\rm const.}$, et alors par cons\'equent, en m\^eme temps
que chaque point r\'eel $x, y$ fix\'e en position g\'en\'erale, 
l'\'el\'ement lin\'eaire r\'eel: $dx : dy$ passant par ce
point, qui est d\'etermin\'e par l'\'equation: $\varphi_x \, dx +
\varphi_y \, dy = 0$, resterait simultan\'ement au repos; mais
comme, apr\`es fixation d'un point $P$ qui se trouve bien s\^ur aussi en
position g\'en\'erale, aucun \'el\'ement lin\'eaire r\'eel passant par
ce point ne doit rester au repos, nous parvenons ainsi \`a une
contradiction.

Il est par ailleurs facile de v\'erifier que notre groupe est \`a
trois param\`etres. En effet, si on fixe le point $P$ et un
\'el\'ement lin\'eaire r\'eel quelconque passant par ce point, alors
aucun mouvement continu n'est encore possible, et il ne reste donc plus
aucun param\`etre arbitraire dans le groupe. Comme en outre $P$, en
tant que point de position g\'en\'erale, peut prendre $\infty^2$
positions diff\'erentes \`a cause de la transitivit\'e du groupe, et
comme, apr\`es fixation de $P$, les \'el\'ements lin\'eaires r\'eels
passant par $P$ peuvent encore tourner tout autour de $P$, la
fixation du point et de l'\'el\'ement lin\'eaire fournit donc
exactement {\em trois}\, conditions pour les param\`etres de $G$, et
par cons\'equent, le nombre de ces param\`etres est simplement \'egal
\`a {\em trois}. Et maintenant, puisque nous avons d\'ej\`a
d\'etermin\'e aux pages~370 sq. tous les groupes r\'eels \`a trois
param\`etres du plan qui sont r\'eels-primitifs, nous obtenons le

\medskip

{\bf Th\'eor\`eme~38.}
{\em Si un groupe r\'eel continu du plan poss\`ede la libre mobilit\'e
dans l'infinit\'esimal en un point r\'eel de position g\'en\'erale,
alors il est \`a trois param\`etres et il est semblable, {\rm via} une
transformation ponctuelle r\'eelle de ce plan, soit au groupe r\'eel
continu projectif \`a trois param\`etres de la 
conique}\footnote{\,
%%%%%%%%%%%%%%%%%%%%%%%-------DEBUT--------%%%%%%%%%%%%%%%%%%%%%%%%%%%
Trois générateurs infinitésimaux sont: 
$p + x^2 p + xyq$, $q + xyp + y^2 q$ et
$yp - xp$; ils annulent en effet l'équation 
de la conique en restriction à la conique. 
}: %%%%%%%%%%%%%%%%%%%%%%%%-----FIN-----%%%%%%%%%%%%%%%%%%%%%%%%%%%%%%%
\[
x^2+y^2+1
=
0,
\]
{\em soit au groupe r\'eel continu projectif de la conique}\footnote{\,
%%%%%%%%%%%%%%%%%%%%%%%-------DEBUT--------%%%%%%%%%%%%%%%%%%%%%%%%%%%
Trois générateurs infinitésimaux sont: 
$-p + x^2 p + xyq$, $-q + xyp + y^2 q$ et
$yp - xp$. 
}: %%%%%%%%%%%%%%%%%%%%%%%%-----FIN-----%%%%%%%%%%%%%%%%%%%%%%%%%%%%%%%
\[
x^2+y^2-1
=
0,
\]
{\em soit \`a l'un des groupes projectifs \`a trois param\`etres de la
forme}:
\[
p,\quad
q,\quad
yp-xq+c\,(xp+yq),
\]
{\em o\`u $c$ d\'esigne une constante r\'eelle}.

\medskip 

En ce qui concerne les groupes que l'on a trouv\'es ici, il est avant
tout remarquable qu'ils poss\`edent tous la libre mobilit\'e dans
l'infinit\'esimal en tous les points r\'eels d'une certaine partie du
plan, alors que nous avons demand\'e seulement qu'ils poss\`edent
cette propri\'et\'e en {\em un}\, point r\'eel de position
g\'en\'erale\footnote{\,
%%%%%%%%%%%%%%%%%%%%%%%-------DEBUT--------%%%%%%%%%%%%%%%%%%%%%%%%%%%
De toute façon, par analyticité, dès qu'elle est satisfaite
en un seul point, 
la libre mobilité
dans l'infinitésimal est alors satisfaite en tous les points
qui n'appartiennent pas à un certain sous-ensemble analytique
fermé de codimension $\geqslant 1$, donc notamment en tous
les points d'un certain ouvert dense. 
}. %%%%%%%%%%%%%%%%%%%%%%%%-----FIN-----%%%%%%%%%%%%%%%%%%%%%%%%%%%%%%%
Ce ph\'enom\`ene trouve son origine dans le fait que
chaque groupe qui poss\`ede la libre mobilit\'e dans l'infinit\'esimal
en un point r\'eel de position g\'en\'erale poss\`ede aussi cette
propri\'et\'e en tous les points r\'eels qui se trouvent dans un
certain voisinage du point en question.

D'un autre c\^ot\'e, les groupes que l'on a trouv\'es montrent que
dans le plan, la libre mobilit\'e dans l'infinit\'esimal ne suffit 
encore pas pour caract\'eriser le groupe des mouvements euclidiens et les
deux groupes de mouvements non-euclidiens, car en dehors de ces
trois groupes, chaque groupe de la forme\footnote{\,
%%%%%%%%%%%%%%%%%%%%%%%-------DEBUT--------%%%%%%%%%%%%%%%%%%%%%%%%%%%
Les courbes intégrales du troisième champ sont des spirales
s'enroulant autour de l'origine, qui dégénèrent en cercles
pour $c = 0$. 
}: %%%%%%%%%%%%%%%%%%%%%%%%-----FIN-----%%%%%%%%%%%%%%%%%%%%%%%%%%%%%%%
\[
p,\quad
q,\quad
yp-xq+c\,(xp+yq)
\quad\quad\quad\quad
{\scriptstyle{(c\,\gtrless\,0)}}
\]
poss\`ede effectivement encore la libre mobilit\'e dans
l'infinit\'esimal en tout point r\'eel qui n'est pas infiniment
\'eloign\'e. Si l'on veut exclure ces groupes, afin de ne conserver
que les mouvements euclidiens et non-euclidiens, on doit donc
\'etablir un axiome sp\'ecial. Cependant, cette observation ne
co\"{\i}ncide pas avec la remarque de Monsieur de Helmholtz, d'apr\`es
laquelle dans le plan, le groupe euclidien et les deux groupes
non-euclidiens ne sont pas les seuls relativement auxquels deux
points ont un et un seul invariant, tandis qu'un nombre de points
sup\'erieur \`a deux n'a pas d'invariant essentiel.

%%%\Fill
%%%[[Trouver la r\'ef\'erence dans Helmholtz et commenter.]]

\bigskip

Nous voulons maintenant compl\'eter les r\'esultats acquis en
recherchant aussi tous les groupes r\'eels continus {\em projectifs}
du plan qui poss\`edent la libre mobilit\'e dans l'infinit\'esimal en
un point r\'eel de position g\'en\'erale.
La connaissance de ces groupes nous sera utile plus tard, lorsque nous
envisagerons l'espace trois fois \'etendu.

Exactement comme plus haut, nous v\'erifions tout d'abord que chaque
groupe projectif du plan qui satisfait nos exigences est transitif et
poss\`ede trois param\`etres. D'apr\`es notre \'enum\'eration des
groupes projectifs du plan \`a trois param\`etres ({\em voir}\,
p.~106), on obtient \`a nouveau que chacun des groupes demand\'es
laisse invariant soit une conique, soit un point.

Si une conique reste invariante, alors son \'equation doit \^etre
r\'eelle, car sinon la conique imaginaire conjugu\'ee resterait aussi
invariante et notre groupe ne pourrait pas alors \^etre \`a trois
param\`etres.

%%%\Mathematiques
%%%[[Classification des coniques imaginaires et r\'eelles~??]]

\noindent
Et maintenant, puisque chaque conique dont l'\'equation est r\'eelle
peut \^etre ramen\'ee, {\em via}\, une transformation projective
r\'eelle, \`a l'une des deux formes: $x^2 + y^2 \pm 1 = 0$, ce
cas-l\`a est d\'ej\`a r\'egl\'e.

Si d'un autre c\^ot\'e un point reste invariant, alors ce point ne
peut pas \^etre r\'eel\footnote{\,
%%%%%%%%%%%%%%%%%%%%%%%-------DEBUT--------%%%%%%%%%%%%%%%%%%%%%%%%%%%
Les transformations sont projectives, donc la droite passant par le
point universellement fixe et le point $P$ qu'on fixe en position
g\'en\'erale serait stabilis\'ee, ce qui contredirait la libre
mobilité dans l'infinitésimal en $P$.
}, %%%%%%%%%%%%%%%%%%%%%%%%-----FIN-----%%%%%%%%%%%%%%%%%%%%%%%%%%%%%%%
car sinon, si un point r\'eel en position
g\'en\'erale \'etait fix\'e, chaque \'el\'ement lin\'eaire passant par
ce point ne pourrait pas tourner autour de lui, et donc il n'y
aurait aucun point r\'eel de position g\'en\'erale en lequel
il poss\`ederait la libre mobilit\'e dans l'infinit\'esimal.
Par cons\'equent, le point invariant est imaginaire et le point
imaginaire conjugu\'e reste simultan\'ement au repos, ainsi que la
droite de liaison r\'eelle passant par les deux.

Si nous d\'epla\c cons alors les deux points appropriés gr\^ace \`a une
transformation r\'eelle projective en les deux points du cercle
%%%\Mathematiques
%%%[[{\sc Kreispunkte}: Quels points~? Quel cercle~?]]
%%%\Mathematiques
%%%[[Pourquoi ce groupe-l\`a \`a quatre param\`etres~?]]
notre groupe se transforme en un sous-groupe \`a trois param\`etres
du groupe suivant \`a quatre param\`etres:
\def\theequation{1}\begin{equation}
p,\quad
q,\quad
yp-xq,\quad
xp+yq,
\end{equation}
et nous devons par cons\'equent seulement rechercher encore tous les
sous-groupes \`a trois param\`etres de ce groupe, qui satisfont notre
exigence.

Si un sous-groupe r\'eel $G$ de~\thetag{ 1} \`a trois param\`etres
poss\`ede la libre mobilit\'e dans l'infinit\'esimal en un point $x_0,
y_0$ de position g\'en\'erale, alors, quand on fixe ce point, les
$\infty^1$ \'el\'ements lin\'eaires: $dx : dy$ qui passent par ce
point sont transform\'es par l'action d'un groupe \`a un
param\`etre. Si nous nous imaginons maintenant que le point: $x_0,
y_0$ est d\'eplac\'e sur l'origine des coordonn\'ees gr\^ace \`a une
transformation r\'eelle du groupe \`a deux param\`etres: $p, q$, alors
le groupe $G$ re\c coit une nouvelle forme dans laquelle
%%%\Mathematiques 
%%%[[Faire le calcul~??]]
il y a en tout cas une transformation infinit\'esimale de la forme:
\def\theequation{2}\begin{equation}
yp-xq+c\,(xp+yq),
\end{equation}
%%%\Mathematiques 
%%%[[Pourquoi cette forme~??]]
car: $xp + yq$ laisse g\'en\'eralement au repos chaque \'el\'ement
lin\'eaire passant par l'origine des coordonn\'ees. De plus, comme
groupe transitif, $G$ contient encore deux transformations
infinit\'esimales de la forme\footnote{\,
%%%%%%%%%%%%%%%%%%%%%%%-------DEBUT--------%%%%%%%%%%%%%%%%%%%%%%%%%%%
---\,\,après combinaison linéaire {\em via}~\thetag{2}\,---
}: %%%%%%%%%%%%%%%%%%%%%%%%-----FIN-----%%%%%%%%%%%%%%%%%%%%%%%%%%%%%%%
\[
p+\alpha\,(xp+yq),
\quad\quad\quad
q+\beta\,(xp+yq).
\]
En calculant maintenant les crochets avec~\thetag{ 2}, on obtient: 
\[
-q+c\,p,
\quad\quad\quad
p+c\,q,
\]
et comme $c$ est r\'eel, et que $G$ ne peut contenir que trois
transformations infinit\'esimales ind\'ependantes, on en
d\'eduit que $\alpha$ et $\beta$ s'annulent tous deux, et que par
cons\'equent $G$ poss\`ede la forme:
\[
p,\quad
q,\quad
yp-xq+c\,(xp+yq).
\]
Ainsi, la proposition suivante est valide.

\medskip

{\bf Proposition~1.}
{\em Si un groupe projectif r\'eel continu du plan poss\`ede la libre
mobilit\'e dans l'infinit\'esimal en un point r\'eel de position
g\'en\'erale, alors il est \`a trois param\`etres et il est semblable,
{\rm via} une transformation projective\footnote{\,
%%%%%%%%%%%%%%%%%%%%%%%-------DEBUT--------%%%%%%%%%%%%%%%%%%%%%%%%%%%
Comparé au Théorème~38, tout est projectif, maintenant. 
} %%%%%%%%%%%%%%%%%%%%%%%%-----FIN-----%%%%%%%%%%%%%%%%%%%%%%%%%%%%%%%
r\'eelle du plan, \`a l'un des
groupes projectifs indiqu\'es dans le Th\'eor\`eme~38.}

\medskip

Parmi les groupes projectifs rassembl\'es dans le Th\'eor\`eme~38, il
y en a deux qui ne poss\`edent pas la libre mobilit\'e dans
l'infinit\'esimal en certains points r\'eels, \`a savoir
premi\`erement le groupe de la conique: $x^2 + y^2 - 1 = 0$, et
deuxi\`emement le groupe: $p, \, q, \, yp - xq + c\, (xp+ yq)$. Pour
le premier, la libre mobilit\'e dans l'infinit\'esimal n'a pas lieu en
tous les points r\'eels de la conique et aussi en tous les points
de la conique qui ne sont par r\'eels\footnote{\,
%%%%%%%%%%%%%%%%%%%%%%%-------DEBUT--------%%%%%%%%%%%%%%%%%%%%%%%%%%%
En effet, puisque la conique (le cercle) est
stabilisé, les directions tangentes sont elles aussi
fixées, dans le domaine réel aussi bien que dans le
domaine complexe. 
}, %%%%%%%%%%%%%%%%%%%%%%%%-----FIN-----%%%%%%%%%%%%%%%%%%%%%%%%%%%%%%% 
et pour le deuxi\`eme, en tous les points r\'eels de la droite \`a
l'infini. Au contraire, il n'y a aucun point exceptionnel 
tel\footnote{\,
%%%%%%%%%%%%%%%%%%%%%%%-------DEBUT--------%%%%%%%%%%%%%%%%%%%%%%%%%%%
En effet, ce n'est qu'en les points, tous purement
complexes, de la conique qu'on n'a pas
libre mobilité dans l'infinitésimal. 
} %%%%%%%%%%%%%%%%%%%%%%%%-----FIN-----%%%%%%%%%%%%%%%%%%%%%%%%%%%%%%%

\medskip
{\bf Th\'eor\`eme~39.}
{\em Si un groupe projectif r\'eel continu du plan poss\`ede {\small\sf 
sans
exception} la libre mobilit\'e dans l'infinit\'esimal en tous les
points r\'eels, alors il est transitif \`a trois param\`etres et il
peut \^etre transform\'e, {\rm via} une transformation r\'eelle
projective, en le groupe projectif r\'eel \`a trois param\`etres de la
conique: $x^2 + y^2 + 1 = 0$.}

\HEAD{Première solution au problème de Riemann-Helmholtz.}{
Division\,\,V.\,\,\,Chapitre\,\,22.\,\,\,\S\,\,98.}

\sectiondritterV{\sf\S\,\,98.}
\label{S-98}
\setcounter{footnote}{0}

Les recherches des pr\'ec\'edents paragraphes doivent maintenant
\^etre accomplies pour l'espace ordinaire trois fois \'etendu. Nous
d\'efinissons donc tout d'abord ce que nous entendons ici par 
<<\,libre
mobilit\'e dans l'infinit\'esimal\,>>.

\medskip

{\bf D\'efinition.} {\em Un groupe r\'eel continu de l'espace trois
fois \'etendu poss\`ede la {\small\sf libre mobilit\'e dans
l'infinit\'esimal} en un point r\'eel $P$ 
lorsque les conditions suivantes
sont remplies: Si l'on fixe le point $P$ et un \'el\'ement lin\'eaire
quelconque passant par ce point, alors un mouvement continu doit
toujours \^etre possible; si au contraire on fixe encore, outre $P$
et cet \'el\'ement lin\'eaire-l\`a, un \'el\'ement de surface
quelconque, qui passe par les deux, alors aucun mouvement continu ne
doit \^etre encore possible.}

\medskip

Nous cherchons maintenant tous les groupes r\'eels continus de $R_3$
qui poss\`edent la libre mobilit\'e dans l'infinit\'esimal en un point
r\'eel de position g\'en\'erale.

Soit $G$ un groupe ayant la qualit\'e requise et soit $P$ le point
r\'eel de position g\'en\'erale en lequel $G$ poss\`ede la libre
mobilit\'e dans l'infinit\'esimal. Si on fixe ensuite $P$ et en outre
aussi un \'el\'ement lin\'eaire quelconque passant par lui, alors
chaque \'el\'ement de surface passant par les deux doit encore pouvoir
tourner contin\^ument autour de l'\'el\'ement lin\'eaire en question,
sinon aucun mouvement continu ne serait encore possible apr\`es
fixation de l'\'el\'ement lin\'eaire concern\'e, en
contradiction avec notre hypoth\`ese.

Se laisse d\'eduire de l\`a que $G$ est transitif. \label{482}
En effet, si $\infty^1$ surfaces r\'eelles: $\varphi ( x_1, x_2, x_3)
= \text{\rm const.}$ restaient invariantes, alors $P$ en commun avec
le plan tangent \`a la surface: $\varphi = \text{\rm const.}$ passant
par lui d\'eterminerait un \'el\'ement de surface qui ne pourrait
tourner autour d'aucun \'el\'ement lin\'eaire contenu
en lui. D'un autre c\^ot\'e, si $\infty^2$ courbes
r\'eelles restaient invariantes, alors on pourrait les arranger
d'une infinit\'e de mani\`eres en une famille de $\infty^1$ surfaces
invariantes, et l'on reviendrait donc \`a la situation dont on a
d\'emontr\'e \`a l'instant qu'elle est impossible.

\`A travers notre point $P$ passent $\infty^2$ \'el\'ements
lin\'eaires: $dx_1 : dx_2 : dx_3$, qui forment une vari\'et\'e plane
deux fois \'etendue. Si l'on fixe $P$, alors cette vari\'et\'e est
transform\'ee par un groupe r\'eel projectif $\mathfrak{ g}$ qui est
manifestement constitu\'e de telle sorte qu'apr\`es fixation d'un
\'el\'ement lin\'eaire r\'eel quelconque, chaque \'el\'ement de surface le
contenant peut tourner autour de l'\'el\'ement lin\'eaire, tandis
qu'absolument aucun mouvement continu n'est encore possible,
aussit\^ot qu'on fixe un \'el\'ement lin\'eaire r\'eel quelconque et
un \'el\'ement de surface quelconque le 
contenant\footnote{\,
%%%%%%%%%%%%%%%%%%%%%%%-------DEBUT--------%%%%%%%%%%%%%%%%%%%%%%%%%%%
Fixer ou bien seulement un élément-plan ou bien 
seulement un élément-ligne laisse bien entendu encore des degrés
de liberté. 
}. %%%%%%%%%%%%%%%%%%%%%%%%-----FIN-----%%%%%%%%%%%%%%%%%%%%%%%%%%%%%%%
Mais en rapportant maintenant les $\infty^2$ \'el\'ements lin\'eaires
passant par $P$ de mani\`ere projective aux points r\'eels d'un plan,
le groupe $\mathfrak{ g}$ se transforme alors en un groupe r\'eel
projectif $\mathfrak{ g'}$ de ce plan qui est isomorphe-holo\'edrique
\`a $\mathfrak{ g}$, et comme les \'el\'ements de surface r\'eels
passant par $P$ correspondent en m\^eme temps aux \'el\'ements
lin\'eaires r\'eels de ce plan, nous reconnaissons alors
imm\'ediatement que le groupe $\mathfrak{ g}'$ poss\`ede la libre
mobilit\'e dans l'infinit\'esimal sans exception en tous les points de
ce plan. D'apr\`es le Th\'eor\`eme~39, il d\'ecoule de l\`a que
$\mathfrak{ g}'$ est transitif \`a trois param\`etres et qu'il peut
\^etre transform\'e, {\em via}\, une transformation r\'eelle
projective du plan, en le groupe de la conique: $x^2 + y^2 + 1 =
0$. Par cons\'equent le groupe $\mathfrak{ g}'$ est aussi transitif
\`a trois param\`etres et il laisse invariante une conique imaginaire
d'\'el\'ements lin\'eaires dont nous pouvons nous imaginer
l'\'equation rapport\'ee \`a la forme:
\def\theequation{3}\begin{equation}
dx_1^2+dx_2^2+dx_3^2
=
0,
\end{equation}
{\em via}\, une transformation lin\'eaire homog\`ene en $dx_1, dx_2,
dx_3$.

En appliquant cela au groupe $G$, on obtient avant tout que $G$ est
\`a six param\`etres. Si en effet l'on fixe le point $P$, ainsi qu'un
\'el\'ement lin\'eaire r\'eel passant par $P$ et un \'el\'ement de
surface r\'eel quelconque qui les contient tous deux, alors les
param\`etres de $G$ seront ainsi soumis \`a exactement $3 + 2 + 1 = 6$
conditions, ce qui d\'ecoule de la transitivit\'e de $G$ et de la
propri\'et\'e de $\mathfrak{ g}$ qu'on a d\'ecrite plus haut. Mais
alors, comme d'apr\`es les hypoth\`eses qu'on a justement pos\'ees,
plus aucun mouvement continu n'est encore possible, le nombre de ces
param\`etres doit \^etre pr\'ecis\'ement \'egal \`a six.

En outre, il est clair qu'apr\`es fixation de $P$, le groupe transitif
\`a six param\`etres $G$ transforme les $\infty^2$ \'el\'ements
lin\'eaires r\'eels par une
action \`a trois param\`etres, et pour
pr\'eciser, de telle sorte que reste invariante une conique imaginaire
du second degr\'e constitu\'ee d'\'el\'ements lin\'eaires. Si donc
nous nous imaginons $P$ transf\'er\'e sur l'origine des coordonn\'ees
{\em via}\, une transformation ponctuelle r\'eelle, et si ensuite on
produit une transformation lin\'eaire homog\`ene r\'eelle en $x_1,
x_2, x_3$ pour que la conique d'\'el\'ements lin\'eaires rattach\'ee
\`a $P$ soit repr\'esent\'ee par l'\'equation~\thetag{ 3}, alors $G$
re\c coit une nouvelle forme telle que dans le voisinage de l'origine
des coordonn\'ees, il contient six transformations infinit\'esimales
ind\'ependantes de
la forme:
\[
\aligned
&
\quad\quad\quad
p_1+\cdots,\quad
p_2+\cdots,\quad
p_3+\cdots,
\\
&
x_\mu\,p_\nu
-
x_\nu\,p_\mu
+
\alpha_{\mu\nu}\,(x_1p_1+x_2p_2+x_3p_3)
+
\cdots
\\
&
\quad\quad\quad\quad\quad\quad
{\scriptstyle{(\mu,\,\nu\,=\,1,\,2,\,3;\ \ \
\mu\,<\,\nu)}}.
\endaligned
\]
Comme $G$ est \`a six param\`etres, les constantes r\'eelles $\alpha_{
\mu \nu}$ doivent ici toutes s'annuler ensemble, comme on s'en
convainc ais\'ement en calculant les crochets.

Le groupe $G$ appartient donc aux groupes r\'eels que nous avons
d\'etermin\'es aux pages~385 sq.; comme il est \`a six param\`etres,
nous obtenons le

%%%\Mathematiques
%%%[[\'Etudier et r\'esumer ces pages~!!]]

\medskip

{\bf Th\'eor\`eme~40.}
\label{Theorem-40}
{\em Si un groupe r\'eel continu de l'espace trois fois \'etendu
poss\`ede la libre mobilit\'e dans l'infinit\'esimal en un point de
position g\'en\'erale, alors il est transitif \`a six param\`etres et
il est équivalent, {\rm via} une transformation ponctuelle r\'eelle de
cet espace, soit au groupe des mouvements euclidiens\footnote{\,
%%%%%%%%%%%%%%%%%%%%%%%-------DEBUT--------%%%%%%%%%%%%%%%%%%%%%%%%%%%
Six générateurs sont: $p_1$, $p_2$, $p_3$ et $x_2p_1 - x_1p_2$,
$x_3p_1 - x_1p_3$, $x_3p_2 - x_2p_3$.
}, %%%%%%%%%%%%%%%%%%%%%%%%-----FIN-----%%%%%%%%%%%%%%%%%%%%%%%%%%%%%%%
soit \`a l'un
des deux groupes de mouvements non-euclidiens de cet espace, et par
cons\'equent dans le deuxi\`eme cas ou bien au groupe projectif r\'eel
continu \`a six param\`etres par lequel la surface imaginaire\footnote{\,
%%%%%%%%%%%%%%%%%%%%%%%-------DEBUT--------%%%%%%%%%%%%%%%%%%%%%%%%%%%
Six générateurs dans l'un 
et l'autre cas sont: $p_1$, $p_2$, $p_3$ et 
$\pm p_1 + x_1( x_1p_1 + x_2p_2 + x_3p_3)$, 
$\pm p_2 + x_2( x_1p_1 + x_2p_2 + x_3p_3)$, 
$\pm p_3 + x_3( x_1p_1 + x_2p_2 + x_3p_3)$.
}: %%%%%%%%%%%%%%%%%%%%%%%%-----FIN-----%%%%%%%%%%%%%%%%%%%%%%%%%%%%%%% 
$x_1^2
+ x_2^2 + x_3^2 + 1 = 0$ reste invariante, ou bien au groupe projectif
r\'eel continu de la surface r\'eelle non r\'egl\'ee: 
$x_1^2 + x_2^2 + x_3^2 - 1 = 0$.}

\medskip

Nous voyons ainsi que dans l'espace trois fois \'etendu, la libre
mobilit\'e dans l'infinit\'esimal suffit parfaitement pour
caract\'eriser les mouvements euclidiens et non-euclidiens, ce qui
offre un nouvel exemple du fait que l'espace trois fois \'etendu se
distingue essentiellement de l'espace deux fois \'etendu.

\bigskip

Nous pouvons maintenant indiquer aussi imm\'ediatement tous les
groupes projectifs r\'eels continus de l'espace ordinaire, qui
poss\`edent la libre mobilit\'e dans l'infinit\'esimal en un point de
position g\'en\'erale.

Chaque groupe projectif $\mathfrak{ G}$ de cette esp\`ece est en effet
semblable, {\em via}\, une transformation ponctuelle r\'eelle, \`a
l'un des trois groupes du th\'eor\`eme d\'emontr\'e \`a l'instant.
Mais d'apr\`es le Th\'eor\`eme~19, p.~292, chacun de ces groupes peut
\^etre transform\'e \`a nouveau en un groupe projectif seulement au
moyen d'une transformation projective, donc on obtient que la
transformation ponctuelle r\'eelle par laquelle le groupe $\mathfrak{
G}$ est semblable \`a l'un de nos trois groupes, est projective.
Enfin, si nous tenons encore compte du fait que l'un des groupes du
Th\'eor\`eme~40 poss\`ede la libre mobilit\'e dans l'infinit\'esimal
en chaque point r\'eel 
de l'espace sans exception, tandis que les deux autres
groupes ne poss\`edent cette propri\'et\'e que dans une portion de
l'espace, nous pouvons \'enoncer notre r\'esultat de la mani\`ere
suivante:

\medskip

{\bf Th\'eor\`eme~41.}
{\em Si un groupe projectif r\'eel continu de l'espace ordinaire trois
fois \'etendu poss\`ede la libre mobilit\'e dans l'infinit\'esimal en
tous les points r\'eels de cet espace sans exception, alors ce groupe
est transitif \`a six param\`etres et il est constitu\'e de toutes les
transformations projectives r\'eelles par l'action desquelles reste
invariante une surface imaginaire non-d\'eg\'en\'er\'ee du second
degr\'e qui est repr\'esent\'ee par une \'equation r\'eelle.}

\medskip

{\bf Proposition~2.}
{\em Si un groupe projectif r\'eel continu de l'espace ordinaire trois
fois \'etendu ne poss\`ede pas la libre mobilit\'e dans
l'infinit\'esimal en tous les points r\'eels de cet espace, mais la
poss\`ede seulement en les points r\'eels d'une certaine r\'egion,
alors ce groupe est transitif \`a six param\`etres, et pour pr\'eciser,
il est semblable, {\rm via} une transformation projective r\'eelle,
soit au groupe projectif r\'eel continu d'une surface r\'eelle du
second degr\'e non r\'egl\'ee, soit au groupe des mouvements
euclidiens.}

\HEAD{Première solution au problème de Riemann-Helmholtz.}{
Division\,\,V.\,\,\,Chapitre\,\,22.\,\,\,\S\,\,99.}

\sectiondritterV{\sf\S\,\,99.}
\label{S-99}
\setcounter{footnote}{0}

Avant que nous passions \`a la g\'en\'eralisation des \'enonc\'es du
pr\'ec\'edent paragraphe \`a l'espace ayant nombre quelconque de
dimensions, nous devons en quelques mots clarifier certaines
mani\`eres de s'exprimer, dont nous ferons usage (\cf aussi la
Remarque p.~338).

Par chaque point d'un espace $n$ fois \'etendu $x_1 \dots x_n$ passent
$\infty^{ n -1}$ \'el\'ements lin\'eaires: $dx_1 : \cdots : dx_n$ qui
forment une vari\'et\'e projective $(n-1)$-fois \'etendue.
Maintenant, dans l'espace ordinaire trois fois \'etendu, nous avons
qualifi\'e par le nom: <<\,\'El\'ement de surface\,>> chaque touffe
unie \deutsch{ebene B\"uschel} de $\infty^1$
tels \'el\'ements lin\'eaires, et dans
$R_n$, nous voulons plut\^ot dire pour cela: {\sl \'El\'ement d'une
vari\'et\'e \`a deux dimensions}, ou plus bri\`evement: {\sl
$M_2$-\'el\'ement}. Pareillement, chaque botte unie \deutsch{ebene
B\"undel} de $\infty^2$ \'el\'ements lin\'eaires devra \^etre
d\'esign\'ee comme {\sl \'El\'ement d'une vari\'et\'e \`a trois
dimensions}, ou plus bri\`evement comme {\sl $M_3$-\'el\'ement}, et
ainsi de suite. Enfin, nous voulons nous r\'eserver le droit
d'employer aussi l'appellation: {\sl $M_1$-\'el\'ement} pour les
\'el\'ements lin\'eaires.

\`A pr\'esent nous pouvons d\'efinir ce que nous
entendons par libre mobilit\'e dans l'infinit\'esimal
dans $R_n$.

\medskip

{\bf D\'efinition.}
{\em Un groupe r\'eel continu de transformations ponctuelles de $R_n$
poss\`ede la {\small\sf 
libre mobilit\'e dans l'infinit\'esimal} en un point
r\'eel $P$ lorsque les conditions suivantes sont remplies: Si l'on fixe le
point $P$, puis un $M_1$-\'el\'ement r\'eel quelconque passant par
lui, puis un $M_2$-\'el\'ement r\'eel quelconque passant par les deux,
puis un $M_3$-\'el\'ement r\'eel passant par les trois, \etc, et
tout de suite enfin un $M_q$-\'el\'ement r\'eel quelconque, qui passe
par tous les \'el\'ements ant\'erieurs, alors un mouvement continu
doit \^etre possible aussi longtemps et seulement aussi longtemps que
$q < n-1$.}

\medskip

Nous sommes tout pr\`es de pr\'esumer que les \'enonc\'es que nous
avons trouv\'es dans le \S~98 pour l'espace \`a trois dimensions se
laissent transf\'erer \`a chaque espace \`a $n > 3$ dimensions. En
d'autres termes, nous sommes sur le point de pr\'esumer que les
\'enonc\'es suivants sont valides:

\medskip

{\bf Th\'eor\`eme~42.}\label{481}
{\em Si un groupe r\'eel continu de $R_n$ $(n \geqslant 3)$ poss\`ede
la libre mobilit\'e dans l'infinit\'esimal en un point r\'eel de
position g\'en\'erale, alors il est transitif \`a $\frac{ 1}{ 2}\, n (
n+1)$ param\`etres
et il est semblable, {\rm via} une transformation ponctuelle
r\'eelle de $R_n$, soit au groupe des mouvements euclidiens, soit \`a
l'un des deux groupes de mouvements non-euclidiens, donc dans le
deuxi\`eme cas ou bien au groupe projectif r\'eel continu de la
vari\'et\'e:}
\[
x_1^2+x_2^2+\cdots+x_n^2+1
=
0
\]
ou bien au groupe projectif r\'eel continu de la vari\'et\'e: 
\[
x_1^2+x_2^2+\cdots+x_n^2-1
=
0.
\]

\medskip

{\bf Th\'eor\`eme~43.}
{\em Si un groupe {\small\sf projectif} r\'eel continu de $R_n$ $(n
\geqslant 3)$ poss\`ede la libre mobilit\'e dans l'infnit\'esimal en
tous les points r\'eels de cet espace sans exception, alors ce groupe
est transitif \`a $\frac{ 1}{ 2}\, n ( n+1)$ param\`etres et il est
constitu\'e de toutes les transformations projectives r\'eelles par
l'action desquelles reste invariante une surface imaginaire
non-d\'eg\'en\'er\'ee du second degr\'e qui est repr\'esent\'ee par
une \'equation r\'eelle.}

\medskip

{\bf Proposition~3.}
{\em Si un groupe projectif r\'eel continu de $R_n$ $(n \geqslant 3)$
ne poss\`ede pas la libre mobilit\'e dans l'infinit\'esimal en tous
les points r\'eels de cet espace, mais la poss\`ede seulement en les
points d'une certaine r\'egion, alors ce groupe est transitif \`a
$\frac{ 1}{ 2}\, n ( n+1)$ param\`etres, et pour pr\'eciser, il est
semblable, {\rm via} une transformation projective r\'eelle, soit au
groupe d'une vari\'et\'e r\'eelle du second degr\'e non r\'egl\'ee,
soit au groupe des mouvements euclidiens de cet espace.}

\medskip

Nous allons maintenant montrer que les choses se comportent
effectivement ainsi. \`A cette fin, nous avons \'evidemment seulement
besoin de montrer qu'en admettant la validit\'e de nos \'enonc\'es
dans un espace \`a $n \geqslant 3$ dimensions, leur justesse en
d\'ecoule toujours dans l'espace \`a $n+1$ dimensions; car puisque
nos \'enonc\'es sont d\'ej\`a reconnus vrais pour $n = 3$, leur
justesse sera alors d\'emontr\'ee pour tout $n \geqslant 3$.

\bigskip

Nous posons donc \`a l'avance que les Th\'eor\`eme~42 et~43 ainsi que
la Proposition~3 sont d\'ej\`a d\'emontr\'es pour l'espace \`a $n
\geqslant 3$ dimensions. Sous ces hypoth\`eses nous cherchons d'abord
tous les groupes r\'eels continus de $R_{ n+1}$ qui poss\`edent la
libre mobilit\'e dans l'infinit\'esimal en un point r\'eel quelconque
de position g\'en\'erale.

Soit $G$ un groupe r\'eel continu de $R_{ n+1}$ qui poss\`ede la libre
mobilit\'e dans l'infinit\'esimal en un point r\'eel $P$ de position
g\'en\'erale. Si ensuite on fixe le point $P$, puis un
$M_1$-\'el\'ement r\'eel quelconque passant par lui, puis un
$M_2$-\'el\'ement r\'eel quelconque passant par les deux, \etc,
et finalement un $M_{ n - 1}$-\'el\'ement r\'eel quelconque passant
par tous les \'el\'ements pr\'ec\'edents, alors chaque
$M_n$-\'el\'ement r\'eel qui passe par le $M_{ n - 1}$-\'el\'ement
fix\'e peut encore tourner contin\^ument autour de cet
\'el\'ement; car s'il ne le pouvait pas, aucun mouvement continu ne
serait encore possible apr\`es fixation des $n-1$ \'el\'ement
sus-nomm\'es.

Il d\'ecoule de l\`a que $G$ est transitif. Si en effet il \'etait
intransitif, alors $R_{ n+1}$ se d\'ecomposerait d'une mani\`ere, ou
bien m\^eme d'une infinit\'e de mani\`eres, en $\infty^1$ vari\'et\'es
r\'eelles invariantes $n$ fois \'etendues, donc il passerait par le
point $P$, qui est en position g\'en\'erale, au minimum un
$M_n$-\'el\'ement r\'eel qui resterait invariant en m\^eme temps que
$P$, en cons\'equence de quoi il ne pourrait pas se tourner
contin\^ument autour d'aucun
$M_{ n - 1}$-\'el\'ement r\'eel contenu en
lui.

Par ailleurs, souvenons-nous \`a ce sujet que les $\infty^n$
$M_1$-\'el\'ements r\'eels passant par $P$ forment une vari\'et\'e
lisse $n$ fois \'etendue, qui\,\,---\,\,le point $P$ \'etant
fix\'e\,\,---\,\,est transform\'ee par l'action d'un groupe projectif
r\'eel $\mathfrak{ g}$. Imaginons-nous les $M_1$-\'el\'ements de cette
vari\'et\'e rapport\'es projectivement
aux points r\'eels d'un espace lisse $n$ fois
\'etendu $R_n'$; alors chaque $M_1$-\'el\'ement passant par $P$ se
transforme en un point de $R_n'$; \`a chaque $M_2$-\'el\'ement
passant par $P$ correspond clairement et de mani\`ere r\'eversible un
$M_1'$-\'el\'ement r\'eel de $R_n'$; en g\'en\'eral, \`a chaque
$M_q$-\'el\'ement r\'eel passant par $P$ correspond clairement et de
mani\`ere r\'eversible un $M_{ q - 1}'$-\'el\'ement de $R_n'$; et
enfin, au groupe $\mathfrak{ g}$ correspond un groupe projectif r\'eel
$\mathfrak{ g}'$ de $R_n'$ qui lui est isomorphe-holo\'edrique. Si
maintenant nous fixons dans $R_n'$ un point, puis un
$M_1'$-\'el\'ement r\'eel quelconque passant par ce point, puis un
$M_2'$-\'el\'ement r\'eel quelconque passant par les deux, \etc,
et finalement un $M_{ n - 2}'$-\'el\'ement r\'eel quelconque qui passe
par tous les \'el\'ements pr\'ec\'edents, alors chaque $M_{ n -
1}'$-\'el\'ement r\'eel qui passe par le point et par tous ces
\'el\'ements doit pouvoir se tourner encore autour du $M_{ n -
2}'$-\'el\'ement; sinon, l'exigence de libre mobilit\'e dans
l'infinit\'esimal au point $P$ ne serait pas satisfaite par notre
groupe $G$.
Par cons\'equent, le groupe projectif r\'eel $\mathfrak{ g}'$
de $R_n'$ poss\`ede la libre mobilit\'e dans
l'infinit\'esimal sans exception en chaque point 
r\'eel de $R_n'$.

Mais puisque d'apr\`es notre hypoth\`ese, le Th\'eor\`eme~43 est valide
dans l'espace $n$ fois \'etendu $R_n'$, le groupe $\mathfrak{ g}'$ est
alors \`a $\frac{ 1}{ 2}\, n ( n+1)$ param\`etres et peut \^etre
rapport\'e, {\em via} une transformation projective r\'eelle de
$R_n'$, \`a une forme telle qu'il laisse invariante la vari\'et\'e
imaginaire du second degr\'e:
\[
{x_1'}^2+{x_2'}^2+\cdots+{x_n'}^2+1
=
0.
\]
Par cons\'equent, le groupe $\mathfrak{ g}$ est aussi \`a $\frac{ 1}{
2}\, n ( n+1)$ param\`etres et il peut \^etre rapport\'e, {\em via}
une transformation lin\'eaire homog\`ene r\'eelle des
diff\'erentielles: $dx_1 \dots dx_{ n +1}$, \`a une forme telle qu'il
laisse invariante la vari\'et\'e imaginaire des $M_1$-\'el\'ements qui
est repr\'esent\'ee par l'\'equation:
\def\theequation{4}\begin{equation}
dx_1^2+\cdots+dx_n^2+dx_{n+1}^2
=
0.
\end{equation}

Si nous revenons maintenant au groupe $G$, nous reconnaissons
qu'apr\`es fixation de $P$, les $\infty^n$ $M_1$-\'el\'ements r\'eels
sont transform\'es de $\frac{ 1}{ 2}\, n ( n+1)$ mani\`eres
diff\'erentes par ce groupe, et pour pr\'eciser, de telle sorte qu'une
certaine vari\'et\'e imaginaire du second degr\'e constitu\'ee de
$M_1$-\'el\'ements reste invariante. Si en particulier nous nous
imaginons que le point $P$ est transf\'er\'e sur l'origine des
coordonn\'ees par une transformation ponctuelle r\'eelle, nous pouvons
toujours parvenir, gr\^ace \`a une transformation lin\'eaire
homog\`ene r\'eelle en les variables $x_1 \dots x_{ n+1}$, \`a ce que
cette vari\'et\'e imaginaire soit repr\'esent\'ee par
l'\'equation~\thetag{ 4}.

Enfin, on doit encore remarquer que $G$ est \`a $\frac{ 1}{ 2}\, (n + 1)
( n+2)$ param\`etres. Si nous fixons en effet le point $P$, alors les
param\`etres de $G$ sont soumis \`a $n+1$ conditions, puis, si nous
fixons tous les $\infty^n$ $M_1$-\'el\'ements r\'eels passant par $P$, les
param\`etres de $G$ sont soumis \`a $\frac{ 1}{ 2}\, n ( n+1)$
conditions suppl\'ementaires, car ces $\infty^n$ $M_1$-\'el\'ements se
transforment, apr\`es fixation de $P$, pr\'ecis\'ement de $\frac{ 1}{
2} \, n ( n+1)$ mani\`eres diff\'erentes. Et maintenant, comme
aussit\^ot que sont fix\'es $P$ ainsi que tous les $M_1$-\'el\'ements
r\'eels passant par lui, tous les $M_2$-, $M_3$-, \dots
$M_n$-\'el\'ements r\'eels passant par $P$ demeurent en m\^eme temps
au repos, en cons\'equence de quoi plus aucun mouvement continu n'est
encore possible, alors le nombre de param\`etres de $G$ doit \^etre
simplement \'egal \`a:
\[
n+1
+
{\textstyle{
\frac{n(n+1)}{1\,\cdot\,2}}}
=
{\textstyle{
\frac{(n+1)(n+2)}{1\,\cdot\,2}}}.
\]

\`A pr\'esent, il est facile de d\'eterminer $G$.

Si en effet nous nous imaginons les variables $x_1 \dots x_{ n+1}$
choisies comme ci-dessus, donc en particulier de telle sorte que $P$
soit l'origine des coordonn\'ees, alors, dans le voisinage de
l'origine des coordonn\'ees, le groupe $G$ contient visiblement
$\frac{ 1}{ 2}\, n ( n+1)$ transformations infinit\'esimales 
ind\'ependantes de la
forme:
\[
\aligned
p_1+\cdots,\
\cdots,\quad
&
p_{n+1}+\cdots,
\quad\quad
x_\mu\,p_\nu-x_\nu\,p_\mu
+
\alpha_{\mu\nu}\,
\sum_1^{n+1}\,
x_\tau\,p_\tau
+\cdots,
\\
&
\quad
{\scriptstyle{(\mu,\,\nu\,=\,1\,\cdots\,n+1;\ \ \
\mu\,<\,\nu)}},
\endaligned
\]
o\`u les $\alpha_{ \mu \nu}$ sont des constantes r\'eelles. Comme $G$
a exactement $\frac{ 1}{ 2} \, (n+1) (n+2)$ param\`etres, on obtient
de plus imm\'ediatement en calculant les crochets que tous les
$\alpha_{ \mu \nu}$ s'annulent. Par cons\'equent, le groupe $G$
appartient aux groupes r\'eels que nous avons d\'etermin\'es aux pages
385~sq., et pour pr\'eciser, on obtient que $G$ peut \^etre
transform\'e {\em via}\, une transformation ponctuelle r\'eelle de
$R_{ n+1}$, soit en le groupe des mouvements euclidiens de cet espace,
soit en le groupe projectif r\'eel continu de l'une des deux
vari\'et\'es du second degr\'e: 
\[
x_1^2+\cdots+x_{n+1}^2\pm 1
=
0.
\]

\medskip

{\em Ainsi, on a d\'emontr\'e que, aussit\^ot que le Th\'eor\`eme~43
est valable dans un espace de dimension $n \geqslant 3$, le
Th\'eor\`eme~42 est toujours vrai dans un espace \`a $n+1$
dimensions.}

\medskip

Mais si le Th\'eor\`eme~42 vaut dans un espace \`a $n+1$ dimensions,
alors le Th\'eor\`eme~43 et la Proposition~3 valent en m\^eme temps
dans cet espace.

En effet, chaque groupe projectif r\'eel continu $\mathfrak{ G}$ de
$R_{ n+1}$, qui poss\`ede la libre mobilit\'e dans l'infinit\'esimal
en un point r\'eel de position g\'en\'erale est semblable, sous les
hypoth\`eses qu'on a justement pos\'ees, {\em via}\, une
transformation ponctuelle r\'eelle, soit au groupe des mouvements
euclidiens de $R_{ n+1}$, soit \`a l'un des deux groupes de mouvements
non-euclidiens de cet espace. Mais d'apr\`es le Th\'eor\`eme~19,
p.~292, cette transformation ponctuelle r\'eelle est n\'ecessairement
projective, et par cons\'equent, nous parvenons tout simplement, avec
les hypoth\`eses que nous avons pos\'ees, au Th\'eor\`eme~43 et \`a la
Proposition~3.

\medskip

{\em Avec cela, on a d\'emontr\'e que, aussit\^ot que les
Th\'eor\`emes~42 et~43 et que la Proposition~3 valent dans un espace
\`a $n \geqslant 3$ dimensions, il sont aussi vrais dans l'espace \`a
$n+1$ dimensions. Puisque nous les avons d\'ej\`a d\'emontr\'es dans
l'espace ordinaire trois fois \'etendu, ils sont donc valides en toute
g\'en\'eralit\'e.}

\bigskip

{\em Si donc $n \geqslant 3$, alors le groupe des mouvements
euclidiens et les deux groupes de mouvements non-euclidiens de $R_n$
sont compl\`etement caract\'eris\'es par l'exigence qu'il doivent
avoir la libre mobilit\'e dans l'infinit\'esimal en un point de
position g\'en\'erale.}

\medskip

Mentionnons encore que les d\'eveloppements qui pr\'ec\`edent ne
valent pas seulement pour les groupes r\'eels dont les transformations
finies sont repr\'esent\'ees par des \'equations analytiques, mais
aussi pour les groupes dont les \'equations finies sont
non-analytiques, pourvu seulement que ces \'equations autorisent un
certain nombre de diff\'erentiations par rapport aux variables et par
rapport aux param\`etres ({\em voir} p.~365 et p.~366).

\HEAD{Première solution au problème de Riemann-Helmholtz.}{
Division\,\,V.\,\,\,Chapitre\,\,22.\,\,\,\S\,\,100.}

\sectiondritterV{\sf\S\,\,100.
\\
Sur le discours d'habilitation de Riemann.}
\label{S-100}
\setcounter{footnote}{0}

Les recherches sur la libre mobilit\'e dans l'infinit\'esimal que nous
venons de conduire sont en relation avec certaines r\'eflexions que
Riemann \`a indiqu\'ees sans son discours d'habilitation. Nous
voulons donc g\'en\'eralement entrer plus avant dans le contenu de ce
discours, d'autant plus qu'il a plusieurs points de contact avec
toutes les recherches r\'ecentes sur les fondements de la
g\'eom\'etrie.

\label{Riemann-jamais-facile}
Riemann n'est jamais facile \`a lire, mais sa c\'el\`ebre allocution
<<\,{\em Sur les hypoth\`eses qui servent de fondement \`a la
g\'eom\'etrie}\,>> pr\'esente des difficult\'es de compr\'ehension
d'un type tout \`a fait sp\'ecial. Riemann a d\^u pr\'esenter sa
soutenance devant une assembl\'ee qui ne consistait qu'en partie de
math\'ematiciens, et pour cette raison, il s'est efforc\'e d'\^etre le
plus g\'en\'eral possible, afin d'\^etre intelligible; mais de mani\`ere
sous-jacente, la nettet\'e et la d\'etermination de l'expression a
subi des dommages en de tr\`es nombreux endroits, et la cons\'equence
en est que maintenant, le math\'ematicien est justement assez
souvent dans le doute quant \`a ce que Riemann a effectivement voulu
dire.

Cette insuffisance se fait particuli\`erement sensible dans la partie
de l'expos\'e o\`u Riemann parle de la mobilit\'e des figures d'un
espace $n$ fois \'etendu; Riemann s'aide d'expressions qui certes
peuvent sembler compr\'ehensibles \`a des non-math\'ematiciens, mais
dont le sens v\'eridique ne peut qu'\^etre devin\'e par le
math\'ematicien.

Malheureusement, \`a ce jour, personne n'a encore jusqu'\`a pr\'esent
examin\'e le contenu du discours de Riemann sous toutes ses facettes
aussi pr\'ecis\'ement que nous l'avons effectu\'e dans le Chapitre~21,
dont le contenu est celui de l'\'etude helmholtzienne <<\,{\em Sur les
faits qui se trouvent au fondement de la g\'eom\'etrie}\,>>. Certes,
on a r\'etabli la th\'eorie riemanienne de la mesure de courbure d'un
$R_n$ d'apr\`es les indications de Riemann, et on s'est convaincu que
la longueur d'un \'el\'ement courbe peu v\'eritablement, aussit\^ot
que la mesure de courbure est partout constante, \^etre rapport\'ee
\`a la forme indiqu\'ee par Riemann. Mais justement, les
consid\'erations de Riemann sur la mobilit\'e des figures semblent
avoir peu ou pas du tout attir\'e l'attention, bien que ces
consid\'erations soient pr\'ecis\'ement d'une importance toute
particuli\`ere pour la conception de Riemann au sujet de la question
sur les fondements de la g\'eom\'etrie. D'un autre c\^ot\'e, il y a
aussi encore par ailleurs une s\'erie de passages dans le discours de
Riemann dont un \'eclaircissement urgent est souhaitable
\deutsch{dringend zu wünschen ist}. 

Il serait tr\`es m\'eritoire que quelqu'un se soum\^{\i}t \`a l'effort
de suivre pas \`a pas l'encha\^{\i}nement des id\'ees de Riemann et
r\'epond\^{\i}t, autant que possible, aux diverses questions qu'on
peut alors s'imposer
\deutsch{die sich da aufdrängen}. 
Nous n'avons toutefois pas l'intention de
r\'ealiser cela, et nous nous contenterons de quelques indications,
qui nous l'esp\'erons, pourront contribuer \`a stimuler un tel travail
\`a fond des pens\'ees riemanniennes.

\bigskip

Comme nous l'avons d\'ej\`a mentionn\'e auparavant, Riemann s'imagine
en premier lieu que les points de l'espace $n$ fois \'etendu sont
d\'etermin\'es par $n$ coordonn\'ees $x_1 \dots x_n$. Nous avons dit
\`a ce moment-l\`a ({\em voir}\, p.~\pageref{394}), que Riemann a
cherch\'e \`a démontrer cela. Mais peut-\^etre avons-nous fait ainsi
du tort \`a Riemann.
%%%\Fill 
%%%[[C'est le moins qu'on puisse dire~! Faire absolument une note.]]

\renewcommand{\thefootnote}{\fnsymbol{footnote}}
On peut aussi interpr\'eter cela en disant que Riemann consid\'erait
que la d\'eterminabilit\'e des points par des coordonn\'ees est une
chose {\em qui va de soi}. Monsieur de Helmholtz semble admettre que
Riemann a \'etabli l'hypoth\`ese en question en tant
qu'axiome\,\footnote[1]{\, {\em voir}\, G\"ott. Nachr. 1868, p.~198.}.

Afin de pouvoir effectuer des d\'eterminations de mesure \`a
l'int\'erieur d'un espace $n$ fois \'etendu, Riemann \'etablit comme
l'hypoth\`ese la plus simple que chaque ligne, donc chaque vari\'et\'e
une fois \'etendue, est mesurable au moyen de chaque
autre\footnote[2]{\, {\em voir}\, {\sc Riemann}, ges. Werke, 1. Ausg.,
p.~258~sq.}. Pour une r\'ealisation effective du mesurage
\deutsch{Messung}, il est
ensuite n\'ecessaire de poss\'eder une expression pour la longueur
d'un \'el\'ement courbe, donc une expression pour la longueur d'un
morceau infiniment petit de ligne entre deux points infiniment voisins
$x_\nu$ et $x_\nu + d x_\nu$. Il pose \`a l'avance que cette longueur
d'un \'el\'ement courbe est une fonction homog\`ene du premier degr\'e
en $dx_1 \dots dx_n$, qui reste inchang\'ee lorsque tous les $dx_\nu$
changent de signe, et qui en outre d\'epend encore de $x_1 \dots x_n$.
\renewcommand{\thefootnote}{\arabic{footnote}}

Ce qui importe maintenant pour cela, c'est de d\'eterminer plus
pr\'ecis\'ement la forme de la longueur d'un \'el\'ement courbe. \`A
cette fin, Riemann admet sans aucun mot de justification que deux
points quelconques de $R_n$ qui ne sont pas infiniment voisins ont
aussi une distance tout \`a fait d\'etermin\'ee l'un par rapport \`a
l'autre, \`a savoir il pose \`a l'avance et sans plus de fa\c cons
l'existence d'une fonction:
\[
\Omega\big(
x_1\dots x_n;\quad
x_1^0\dots x_n^0
\big),
\]
qui a la m\^eme valeur pour tous les points $x_1 \dots x_n$ qui
<<\,{\em sont aussi lointainement distants} ({\small\sf gleich weit
abstehen})\,>> du point $x_1^0 \dots x_n^0$. On ne voit pas si Riemann
a voulu faire une nouvelle hypoth\`ese avec cela, ou si il a voulu
signifier que l'existence d'une telle fonction $\Omega$ d\'ecoule de
l'existence d'une longueur d'\'el\'ement courbe ayant la constitution
indiqu\'ee.  En tout cas, cette derni\`ere possibilit\'e n'est pas
imm\'ediatement claire, car il est certain que l'existence d'une telle
fonction $\Omega$ d\'ecoule de l'existence d'une longueur
d'\'el\'ement courbe seulement lorqu'entre deux points, une ligne la
plus courte possible est en m\^eme temps d\'etermin\'ee au moyen d'une
telle longueur d'\'el\'ement courbe; il devrait donc \^etre
d\'emontr\'e d'abord que l'existence de telles lignes les plus courtes
possibles, d\'ecoule des hypoth\`eses qui ont
\'et\'e faites auparavant sur la longueur d'un \'el\'ement courbe.

%%%\Fill
%%%[[longueur int\'egr\'ee: localement, \c ca marche; globalement: Hopf-Rinow;]]

\bigskip

{\small
On peut se convaincre directement que l'existence d'une longueur
d'\'el\'ement courbe:\label{487}
\[
\omega
\big(
x_1\dots x_n;\ \
dx_1\dots dx_n
\big),
\]
qui est une fonction homog\`ene {\em quelconque}\, du premier degr\'e
par rapport \`a $dx_1 \dots dx_n$, ne permet 
pas toujours d'en tirer
l'existence d'une fonction distance:
\[
\Omega
\big(
x_1\dots x_n;\ \
y_1\dots y_n
\big)
\]
entre deux points finiment \'eloign\'es l'un de l'autre.

Soit en effet $\Gamma$ le plus grand groupe continu de transformations
ponctuelles r\'eelles de $R_n$ par lesquelles la longueur $\omega ( x,
dx)$ d'un \'el\'ement courbe reste invariante. Alors nous pouvons
manifestement exprimer aussi l'existence de cette longueur d'un
\'el\'ement courbe en disant: les deux points infiniment voisins
$x_\nu$ et $x_\nu + d x_\nu$ ont l'invariant: $\omega ( x, dx)$
relativement \`a $\Gamma$. De m\^eme, l'existence d'une fonction
distance $\Omega ( x, y)$ reviendrait \`a ce que deux
points $x_\nu$ et $y_\nu$ qui sont finiment \'eloign\'es l'un de
l'autre ont l'invariant $\Omega ( x,y)$ relativement au groupe
$\Gamma$.

Mais maintenant, deux points infiniment voisins $x_\nu$ et $x_\nu +
dx_\nu$ peuvent tr\`es bien avoir un invariant de la forme $\omega (
x, dx)$ relativement \`a un groupe, sans que les deux points $x_\nu$
et $y_\nu$ aient un invariant. Cela est encore montr\'e par
le groupe~\thetag{ 22} de $R_3$ 
mentionn\'e \`a la page~\pageref{457-0}:
\[
q,\quad
xq+r,\quad
x^2q+2xr,\quad
x^3q+3x^2r,\quad
p,\quad
xp-zr.
\]
En effet, par son action les deux points
infiniment voisins $x, y, z$ et $x + dx$, $y + dy$, $z + dz$ ont
l'invariant: $dy - z\,dx$, alors que deux points finiment
\'eloign\'es l'un de l'autre n'ont aucun invariant.

%%%\Mathematiques
%%%[[Check $dy-zdx$.]]

\bigskip

On peut bien s\^ur encore penser que les choses se passent tout
autrement si, avec Riemann, on ajoute l'hypoth\`ese que la longueur
$\omega ( x, dx)$ d'un \'el\'ement courbe ne modifie pas sa valeur,
lorsque tous les $dx_\nu$ changent de signe. Dans ce cas, on ne peut
rien d\'eduire de l'exemple donn\'e \`a l'instant, car dans ces
circonstances, l'invariant: $dy - z \, dx$ ne satisfait pas
l'exigence de Riemann. Nous voulons cependant ne pas donner suite aux
questions sugg\'er\'ees par cela.

}

\bigskip

Poursuivons notre rapport sur l'encha\^{\i}nement des
id\'ees de Riemann.

Riemann pose donc \`a l'avance qu'il y a une fonction: 
\[
\Omega\big(
x_1\dots x_n;\ \
x_1^0\dots x_n^0
\big)
\]
qui a la m\^eme valeur pour tous les points qui sont aussi
lointainement distants de $x_1^0 \dots x_n^0$. Comme cela semble
\^etre le cas, le concept g\'en\'eral de fonction distance dans une
vari\'et\'e $n$ fois \'etendue appara\^{\i}t ici pour la premi\`ere
fois.

\label{Omega-metrique}
Au sujet de la fonction $\Omega$, Riemann suppose qu'en tant que
fonction de $x_1 \dots x_n$, dans un voisinage du syst\`eme de
valeurs: $x_1^0 \dots x_n^0$, elle cro\^{\i}t de tous les c\^ot\'es,
et par cons\'equent, elle a un minimum pour: $x_\nu = x_\nu^0$; il
admet de plus que ses quotients diff\'erentiels du premier et du
deuxi\`eme ordre sont tous finis pour $x_\nu = x_\nu^0$. Sous ces
hypoth\`eses, la premi\`ere diff\'erentielle de $\Omega$ doit
s'annuler pour $x_\nu = x_\nu^0$ et la seconde:
\[
d^2\Omega
=
\sum_{\mu\,\nu}^{1\dots n}\,
\frac{\partial^2\Omega}{\partial x_\mu\,\partial x_\nu}\,\,
dx_\mu\,dx_\nu
\]
ne doit devenir n\'egative pour aucun syst\`eme de valeurs: $dx_1
\dots dx_n$, lorsqu'on substitue: $x_\nu = x_\nu^0$.

Riemann se restreint \`a la consid\'eration du cas o\`u la
diff\'erentielle $d^2 \Omega$ se transforme, pour $x_\nu = x_\nu^0$,
en une forme quadratique constamment positive par rapport \`a $dx_1
\dots dx_n$.
Le carr\'e de la longueur $ds$ d'un \'el\'ement courbe appartenant au
point $x_1^0 \dots x_n^0$ ne peut alors se distinguer de
l'expression:
\[
\sum_{\mu\,\nu}^{1\dots n}\,
\left[
\frac{\partial^2\Omega}{\partial x_\mu\,\partial x_\nu}
\right]_{x=x^0}\,
dx_\mu\,dx_\nu
\]
que par un facteur constant d\'ependant de $x_1^0 \dots x_n^0$, et par
cons\'equent la longueur d'un \'el\'ement courbe en un point
quelconque $x_1\dots x_n$
est d\'etermin\'ee par une \'equation de la forme:
\def\theequation{5}\begin{equation}
ds^2
=
\sum_{\mu\,\nu}^{1\dots n}\,
\alpha_{\mu\,\nu}\,
(x_1\dots x_n)\,
dx_\mu\,dx_\nu,
\end{equation}
o\`u les $\alpha_{\mu\nu}$ sont des fonctions r\'eelles dont le
d\'eterminant ne s'annule pas identiquement, et o\`u le membre de
droite est une fonction constamment positive de $dx_1 \dots dx_n$ pour
des valeurs g\'en\'erales de $x_1 \dots x_n$.

Nous avons restitu\'e la d\'erivation de l'expression de la longueur
d'un \'el\'ement courbe d'une mani\`ere d\'etaill\'ee, parce qu'elle
est d'un int\'er\^et particulier, notamment {\em via}\, l'introduction
d'une fonction distance.

L'expression trouv\'ee pour la longueur d'un \'el\'ement courbe permet
d'utiliser le calcul des variations, et on obtient qu'entre deux
points qui se trouvent \`a l'int\'erieur d'une certaine r\'egion, une
et une seule ligne la plus courte est possible
et qu'elle se trouve
enti\`erement \`a cette r\'egion. De plus, il est clair que par
chaque point $x_1 \dots x_n$, il passe une 
unique ligne la plus courte
contenant un \'el\'ement lin\'eaire: $dx_1 : \dots : dx_n$ donn\'e
passant par ce point.

Riemann consid\`ere maintenant un $M_2$-\'el\'ement quelconque passant
par un point quelconque $x_1 \dots x_n$, et il s'imagine les lignes
les plus courtes qui passent par chacun des $\infty^1$ \'el\'ements
lin\'eaires de ce $M_2$-\'el\'ement. La vari\'et\'e deux fois
\'etendue ainsi engendr\'ee est uniquement d\'etermin\'ee par le
$M_2$-\'el\'ement choisi, et sa mesure gaussienne de courbure au point
$x_1 \dots x_n$ a une valeur enti\`erement d\'etermin\'ee. La valeur
de cette mesure gaussienne de courbure, qui d\'epend seulement de la
forme indiqu\'ee~\thetag{ 5} de la longueur d'un \'el\'ement courbe,
Riemann l'appelle la <<\,{\sl mesure de courbure} \deutsch{Mass der
Kr\"ummung}\,>> que l'espace $n$ fois \'etendu, muni de la
longueur~\thetag{ 5} d'un \'el\'ement courbe, possède au point $x_1 \dots
x_n$ dans la direction du $M_2$-\'el\'ement choisi. Il affirme de
surcro\^{\i}t que r\'eciproquement, la longueur~\thetag{ 5} d'un
\'el\'ement courbe d'un espace $n$ fois \'etendu est en g\'en\'eral
compl\`etement d\'etermin\'ee, aussit\^ot que la mesure de courbure
est donn\'ee en chaque point dans les directions de $\frac{ 1}{ 2}\,
n(n-1)$ $M_2$-\'el\'ements qui sont mutuellement en position
g\'en\'erale.

Parmi les vari\'et\'es $n$ fois \'etendues qui ont une longueur
d'\'el\'ement courbe de la forme~\thetag{ 5} mentionn\'ee ci-dessus,
Riemann consid\`ere en particulier toutes celles pour lesquelles la
mesure de courbure est constante en tous lieux, donc pour lesquelles
cette mesure de courbure a la m\^eme valeur en tous les points et en
tous les $M_2$-\'el\'ements passant par ces points. Il dit ce qui
suit au sujet de ces vari\'et\'es:

\medskip

\renewcommand{\thefootnote}{\fnsymbol{footnote}}
\label{Riemann-dehnung}
<<\,Le caract\`ere commun de ces vari\'et\'es dont la mesure de
courbure est constante, peut aussi \^etre exprim\'e en disant que les
figures peuvent se mouvoir\footnote[1]{\, Ici \`a vrai dire, Riemann
aurait m\^eme d\^u ajouter le mot <<\,librement\,>>.} sans
\'elargissement \deutsch{Dehnung} en elles.
\renewcommand{\thefootnote}{\arabic{footnote}} Car il est \'evident
que les figures en elles ne pourraient pas coulisser
\deutsch{verschiebar sein}
et pivoter \deutsch{drehbar sein} librement, si
la mesure de courbure n'\'etait pas la m\^eme en chaque point et dans
toutes les directions. Mais d'autre part, les rapports m\'etriques de
la vari\'et\'e sont compl\`etement d\'etermin\'es par la mesure de
courbure; donc les rapports m\'etriques autour d'un point et dans
toutes les directions sont exactement les m\^emes qu'autour d'un autre
point, et par cons\'equent, \`a partir de ce premier point, les
m\^emes constructions peuvent \^etre transf\'er\'ees, d'o\`u il
s'ensuit que, dans les vari\'et\'es dont la mesure de courbure est
constante, on peut donner aux figures chaque position quelconque.\,>>

En ces termes, il est exprim\'e que l'exigence d'apr\`es laquelle la
mesure de courbure doit \^etre partout constante, a la m\^eme
signification que certaines exigences concernant la mobilit\'e des
figures. Si nous pla\c cons cette derni\`ere exigence au tout
d\'ebut, nous pouvons restituer de la mani\`ere suivante
l'encha\^{\i}nement des id\'ees de Riemann d'une mani\`ere quelque peu
plus pr\'ecise, quoique non absolument pr\'ecise:

\medskip

\label{EL-RH}
{\em Riemann cherche, parmi les vari\'et\'es dont la longueur d'un
\'el\'ement courbe a la forme~\thetag{ 5}, toutes celles dans
lesquelles les figures peuvent occuper chaque position quelconque,
c'est-\`a-dire, dans lesquelles les figures peuvent coulisser et
tourner, sans subir d'\'elargissement. Il parvient \`a ce r\'esultat que
les vari\'et\'es dont la mesure de courbure est constante en tous
lieux sont les seules dans lesquelles les figures sont mobiles de
cette mani\`ere.}

\medskip

Mais maintenant, que doit-on entendre l\`a-dessous, 
quand on dit que les figures
peuvent occuper chaque position quelconque, ou 
encore qu'elles doivent
pouvoir coulisser et tourner librement~?

Tout d'abord, nous voulons \'etablir ce que cela veut dire que les
figures doivent \^etre mobiles sans \'elargissement et en g\'en\'eral.

\label{490}
Quand Riemann parle d'un mouvement des figures, il s'imagine
indubitablement que ce mouvement est continu, ce que montrent
notamment les mots <<\,coulisser\,>> et <<\,tourner\,>> dont il fait
usage.  Mais maintenant, chaque mouvement continu am\`ene avec lui une
modification continue des cordonn\'ees du point qui se meut, et
fournit ainsi une famille continue de $\infty^1$ transformations
r\'eelles:
\def\theequation{6}\begin{equation}
x_\nu'
=
f_\nu(x_1\dots x_n;\ \ 
t)
\ \ \ \ \ \ \ \ \ \ \ \ \
{\scriptstyle{(\nu\,=\,1\,\cdots\,n)}}
\end{equation}
de la vari\'et\'e $x_1 \dots x_n$, et pour pr\'eciser, une famille
dans laquelle se trouve la transformation identique (\cf
p.~\pageref{440}). Si les figures ne doivent subir aucun
\'elargissement au cours de ce mouvement, alors la longueur de chaque
ligne doit rester inchang\'ee, autrement dit: les $\infty^1$
transformations~\thetag{ 6} doivent laisser invariante la
longueur~\thetag{ 5} d'un \'el\'ement lin\'eaire. Inversement, il est
\'evident que chaque famille~\thetag{ 6} de $\infty^1$ transformations
qui laisse invariante la longueur~\thetag{ 5} d'un \'el\'ement
lin\'eaire et qui contient la transformation identique, repr\'esente
un mouvement continu relativement auquel les figures ne subissent
aucun \'elargissement. Par cons\'equent, l'exigence que le
mouvement des figures soit en g\'en\'eral possible sans
\'elargissement,
revient \`a ce qu'il doive exister au moins une famille continue
de $\infty^1$ transformations contenant la transformation identique
qui laisse invariante la longueur~\thetag{ 5} d'un \'el\'ement
lin\'eaire.

Maintenant, l'ensemble de toutes les transformations r\'eelles,
relativement auxquelles la longueur~\thetag{ 5} d'un \'el\'ement
lin\'eaire reste invariante, constitue certainement un groupe
$\mathfrak{ G}$, et ce groupe $\mathfrak{ G}$ contient toujours un
certain sous-groupe continu $G$, qui n'est contenu dans aucun
sous-groupe continu plus grand de $\mathfrak{ G}$. Si ensuite~\thetag{
6} est une famille de $\infty^1$ transformations r\'eelles continues
qui contient la transformation identique et qui appartient au groupe
$G$, alors cette famille d\'etermine visiblement un mouvement continu
au cours duquel les figures ne subissent aucun \'elargissement.
Inversement, si~\thetag{ 6} est un mouvement continu au cours duquel
les figures ne subissent aucun \'elargissement, alors la
famille~\thetag{ 6} appartient \'evidemment au groupe $G$. Il
d\'ecoule de l\`a que $G$ est compl\`etement d\'etermin\'e par 
l'ensemble de tous les 
mouvements au cours desquels les figures ne subissent
aucun \'elargissement; en m\^eme temps, il est clair que les
transformations de $G$ d\'eterminent toutes les positions que l'on
peut donner aux figures par des mouvements continus sans
\'elargissement.

Et maintenant, qu'est-ce que cela veut dire, lorsque Riemann demande
que les figures puissent {\em coulisser et tourner librement
sans \'elargissement}~?

Manifestement, le mot <<\,librement\,>> est si ind\'etermin\'e,
qu'on pourrait en tirer n'importe quelle conclusion. 
On en est donc r\'eduit \`a indiquer des pr\'esomptions.

Nous croyons qu'avec la {\em possibilit\'e libre de coulisser} 
\deutsch{belibige Verschiebbarkeit} 
pour les figures, Riemann s'est imagin\'e
qu'au cours du mouvement dont une figure est susceptible sans subir
d'\'elargissement, un point de la figure s\'electionn\'e
arbitrairement peut \^etre envoy\'e sur chaque autre point de
l'espace. Avec cette interpr\'etation, l'exigence de possibilit\'e
libre de coulisser revient \'evidemment \`a ce que le groupe d\'efini
plus haut doive \^etre transitif.

D'un autre c\^ot\'e, par <<\,rotations\,>> 
\deutsch{Drehungen},
Riemann peut en tout cas seulement entendre les mouvements qui sont
encore possibles apr\`es fixation d'un point quelconque. L'exigence
de {\em possibilit\'e libre de tourner} 
\deutsch{beliebigen Drehbarkeit}
pour les figures, nous la comprenons 
donc maintenant en disant que le
mouvement le plus g\'en\'eral qui est encore possible sans
\'elargissement apr\`es fixation de l'origine des coordonn\'ees,
doit d\'ependre d'autant de param\`etres qu'{\em autorise
l'invariance de l'expression}:
\[
\sum_{\mu\,\nu}^{1\dots n}\,
\alpha_{\mu\nu}\,(0\dots 0)\,
dx_\mu\,dx_\nu,
\]
de telle sorte que le groupe transitif $G$ doive contenir le plus
grand nombre possible de param\`etres.

Afin de nous rendre compte de la port\'ee de cette exigence, 
nous devons rechercher quels sont les mouvements qui sont encore
possibles sans \'elargissement apr\`es fixation d'un point r\'eel
en position g\'en\'erale.

Pour des raisons de simplicit\'e, nous supposons que l'origine des
coordonn\'ees est un point de position g\'en\'erale, et nous nous
imaginons de plus que la longueur de l'\'el\'ement courbe
est rapport\'ee, {\em via}\, une transformation lin\'eaire
homog\`ene de $x$, \`a une forme telle que pour:
\[
x_1=x_2=\cdots=x_n=0,
\]
elle prend la forme:
\def\theequation{7}\begin{equation}
dx_1^2+dx_2^2+\cdots+dx_n^2.
\end{equation}
Chaque mouvement qui est 
encore \'eventuellement
possible sans \'elargissement apr\`es fixation
de l'origine des coordonn\'ees est alors repr\'esent\'e par une
transformation de la valeur:
\def\theequation{8}\begin{equation}
x_k'
=
\sum_1^n\,\alpha_{k\nu}\,x_\nu
+\cdots
\ \ \ \ \ \ \ \ \ \ \ \ \
{\scriptstyle{(k\,=\,1\,\cdots\,n)}},
\end{equation}
par l'action de 
laquelle la longueur~\thetag{ 5} d'un \'el\'ement courbe reste
invariante; \`a ce sujet, le d\'eterminant des $\alpha_{ k \nu}$ ne
doit pas s'annuler et les termes supprim\'es sont du deuxi\`eme ordre
en $x$, et d'ordre sup\'erieur.

Sous les hypoth\`eses pos\'ees il est clair que la transformation
r\'eduite issue de~\thetag{ 8}:
\def\theequation{9}\begin{equation}
x_k'
=
\sum_1^n\,\alpha_{k\nu}\,x_\nu
\ \ \ \ \ \ \ \ \ \ \ \ \
{\scriptstyle{(k\,=\,1\,\cdots\,n)}}
\end{equation}
doit laisser invariante l'expression diff\'erentielle~\thetag{ 7}. En
outre, il est facile de voir que la transformation~\thetag{ 8} est
compl\`etement d\'etermin\'ee, aussit\^ot qu'on conna\^{\i}t la
transformation r\'eduite associ\'ee~\thetag{ 9}. La transformation
r\'eduite associ\'ee~\thetag{ 9} indique en effet comment sont
transform\'es les $\infty^{ n- 1}$ \'el\'ements lin\'eaires passant
par l'origine des coordonn\'ees. Maintenant, une unique ligne la plus
courte est dirig\'ee par chaque \'el\'ement lin\'eaire $dx_1 : \cdots
: dx_n$ passant par l'origine des coordonn\'ees et chaque point $x_1
\dots x_n$ de cette ligne la plus courte est compl\`etement
d\'etermin\'e, lorsque, mis \`a part l'\'el\'ement lin\'eaire en
question, on conna\^{\i}t aussi sa distance $r$ \`a l'origine des
coordonn\'ees. Si donc l'on sait que par la transformation~\thetag{ 8}
qui laisse invariante la longueur~\thetag{ 5} d'un \'el\'ement courbe,
l'\'el\'ement lin\'eaire: $dx_1 : \cdots : dx_n$ est transform\'e en
l'\'el\'ement lin\'eaire: $dx_1 ' : \cdots : dx_n'$, alors en m\^eme
temps, la nouvelle position $x_1 ' \dots x_n'$ que le point $x_1 \dots
x_n$ re\c coit par la transformation~\thetag{ 8}, est elle aussi
d\'etermin\'ee, car en effet, $x_1 ' \dots x_n'$ est \'evidemment le
point qui se trouve \`a la distance $r$ de l'origine des coordonn\'ees
sur la ligne la plus courte dirig\'ee par l'\'el\'ement lin\'eaire:
$dx_1 ' : \cdots : dx_n'$.

Nous voyons \`a partir de l\`a que chaque mouvement sans
\'elargissement qui est encore possible apr\`es fixation de l'origine
des coordonn\'ees, est parfaitement d\'etermin\'e par la fa\c con dont
il d\'eplace les $\infty^{ n-1}$ \'el\'ements lin\'eaires passant par
l'origine des coordonn\'ees. Si donc nous demandons que le mouvement
le plus g\'en\'eral
encore possible sans \'elargissement apr\`es fixation de l'origine
des coordonn\'ees d\'epende du plus grand nombre de param\`etres,
alors cela revient \`a demander qu'apr\`es fixation par le groupe $G$
de l'origine des coordonn\'ees, les \'el\'ements lin\'eaires passant
par ce point sont transform\'es par un groupe projectif ayant le plus
grand nombre possible de param\`etres.

Si nous rassemblons les r\'esultats acquis, nous pouvons dire: {\em
Lorsque Riemann demande qu'il existe une longueur d'\'el\'ement
courbe constamment positive de la forme~\thetag{ 5}
et que les figures de l'espace
$x_1 \dots x_n$ puissent coulisser et tourner librement
sans \'elargissement, il demande
avec cela exactement la m\^eme chose que lorsqu'on exige qu'une
longueur d'\'el\'ement courbe~\thetag{ 5} reste invariante par un
groupe continu $G$ qui premi\`erement, est transitif et
deuxi\`emement, transforme de la mani\`ere la plus g\'en\'erale
possible les $\infty^{ n-1}$ \'el\'ements lin\'eaires passant par
chaque point r\'eel fix\'e en position g\'en\'erale}
(\cf pp.~355 sq.).

Si nous relions maintenant les d\'eveloppements des pages
353~sq. avec ceux des pages 385~sq., nous reconnaissons
imm\'ediatement que chaque groupe r\'eel continu $G$ ayant la
constitution ici demand\'ee peut \^etre transform\'e, {\em via}\, une
transformation ponctuelle r\'eelle, soit en le groupe des mouvements
euclidiens de $R_n$, soit en l'un des deux groupes de mouvements
non-euclidiens. Les exigences que Riemann \'enonce au sujet de la
mobilit\'e des figures suffisent donc, en commun avec l'exigence qu'il
existe une longueur d'\'el\'ement courbe, pour caract\'eriser ces
trois groupes de mouvements.

Nous ins\'erons \`a pr\'esent encore quelques remarques au sujet de la
solution que Riemann a lui-m\^eme donn\'ee \`a son probl\`eme.

Riemann dit: <<\,Car il est \'evident que les figures ne
pourraient pas coulisser et tourner librement, si la mesure de
courbure n'\'etait pas la m\^eme en chaque point et dans toutes les
directions\,>>. De la constance de la mesure de courbure, il d\'eduit
maintenant que les rapports m\'etriques sont, dans toutes les
directions autour d'un point, exactement les m\^emes qu'autour d'un
autre point,
et de l\`a il d\'ecoule \`a nouveau que les figures peuvent coulisser et
tourner librement sans \'elargissement de la mani\`ere qui a \'et\'e
d\'ecrite plus haut.

Maintenant, l'exigence que la mesure de courbure soit constante en
tous les points et dans toutes les directions de $M_2$-\'el\'ements
passant par ces points, est visiblement n\'ecessaire, lorsque le
groupe $G$ doit \^etre transitif et lorsqu'il doit \^etre en outre
constitu\'e de telle sorte qu'apr\`es fixation d'un point r\'eel en
position g\'en\'erale, chaque $M_2$-\'el\'ement r\'eel passant par ce
point peut encore \^etre transform\'e en tout autre. Inversement, si
le groupe est constitu\'e de telle sorte qu'apr\`es fixation d'un
point r\'eel en position g\'en\'erale, chaque $M_2$-\'el\'ement r\'eel
peut encore \^etre transform\'e en tout autre, on en d\'eduit
ais\'ement que le groupe est transitif
%%%\Fill
%%%[[Implique toujours transitif.]]
et que la mesure de courbure a la m\^eme valeur en tous les points et
dans toutes les directions de $M_2$-\'el\'ements passant par ces
points. Mais il d\'ecoule de l\`a, d'apr\`es Riemann, que les figures
peuvent coulisser et tourner librement sans \'elargissement 
%%%\Mathematiques
%%%[[Je ne devine pas le raisonnement~! Lie admet quelque chose de Riemann que Riemann n'a pas d\'etaill\'e dans son discours.]]
et \`a nouveau, cela revient \`a ce qu'apr\`es fixation d'un point
r\'eel de position g\'en\'erale, les \'el\'ements lin\'eaires passant
par ce point sont transform\'es par $G$ de la mani\`ere la plus
g\'en\'erale possible.

Par cons\'equent, la proposition suivante d\'ecoule des
d\'eveloppements riemanniens:

\smallskip

{\em Si $G$ est le plus grand groupe r\'eel continu qui laisse
invariante une expression diff\'erentielle constamment positive de la
forme~\thetag{ 5}, et si ce groupe est constitu\'e de telle sorte
qu'apr\`es fixation d'un point r\'eel en position g\'en\'erale, chaque
$M_2$-\'el\'ement r\'eel passant par ce point peut \^etre transform\'e
en tout autre, alors il transforme les $\infty^{ n - 1}$ \'el\'ements
lin\'eaires passant par ce point de la mani\`ere la plus g\'en\'erale
possible.}

%%%\Mathematiques
%%%[[Encore une fois ici, je n'ai pas saisi le raisonnement math\'ematique. Comment assure-t-on que les \'el\'ements lin\'eaires se transforment par le groupe maximal du $ds^2 \vert_P$, seulement \`a partir du transfert des $M_2$-\'el\'ements~? \c Ca doit provenir du raisonnement sur les g\'eod\'esiques qui \'etoilent-prolongent une transformation infinit\'esimale locale.]]

%%%\Mathematiques
%%%[[Existe-t-il des cas non maximalement isotropes~?]]

\bigskip

Riemann pose \`a l'avance que le carr\'e de la longueur d'un
\'el\'ement courbe est une expression diff\'erentielle du second ordre
constamment positive et il ajoute l'exigence que les figures puissent
coulisser et tourner librement sans subir d'\'elargissement. Cette
derni\`ere exigence a, si on l'interpr\`ete de la mani\`ere dont nous
l'avons envisag\'ee, la m\^eme signification que l'exigence que
l'expression diff\'erentielle en question admette un groupe r\'eel
continu qui est transitif et qui transforme de la mani\`ere la plus
g\'en\'erale possible les \'el\'ements lin\'eaires r\'eels passant par
le point r\'eel de position g\'en\'erale qui a \'et\'e fix\'e.

Il est clair que le groupe r\'eel $G$, sous ces hypoth\`eses,
poss\`ede la libre mobilit\'e dans l'infinit\'esimal en chaque point
de position g\'en\'erale.
%%%\Mathematiques
%%%[[C'est elliptique: il faut analyser le groupe orthogonal pour v\'erifier que la suite des \'el\'ements lin\'eaires embo\^{\i}t\'es permet des mobilit\'es successives jusqu'au blocage.]]
D'un autre c\^ot\'e, dans les \S\S~97--99, nous avons d\'emontr\'e que
chaque groupe r\'eel $\mathfrak{ G}$ de $R_n$ $(n > 2)$, qui poss\`ede
la libre mobilit\'e dans l'infinit\'esimal en un point r\'eel de
position g\'en\'erale peut \^etre transform\'e, {\em via}\, une
transformation ponctuelle r\'eelle soit en le groupe des mouvements
euclidiens, soit en l'un des deux groupes de mouvements
non-euclidiens. Par cons\'equent, les deux exigences riemanniennes:
Existence d'une longueur quadratique d'\'el\'ement courbe constamment
positive, et possibilit\'e pour les figures de coulisser et de tourner
librement sans \'elargissement, ces deux exigences suffisent pour la
caract\'erisation des mouvements euclidiens et non-euclidiens, pourvu
qu'on interpr\`ete la possibilit\'e de coulisser et de tourner
librement comme cela a \'et\'e fait plus haut. Mais nous voyons en
m\^eme temps que {\em l'exigence de libre mobilit\'e dans
l'infinit\'esimal, en commun avec l'exigence que les mouvements
doivent constituer un groupe, ces deux exigences peuvent apparemment
\^etre compl\`etement substitu\'ees aux exigences poursuivies par
Riemann}.

\bigskip

{\small 

Nous allons encore montrer que l'on peut d\'eduire l'existence d'une
longueur quadratique d'\'el\'ement courbe constamment positive
directement \`a partir de l'exigence de libre mobilit\'e dans
l'infinit\'esimal, sans s'appuyer sur les d\'eveloppements des
\S\S~97--99.

Tout d'abord, on peut montrer d'une mani\`ere essentiellement plus
simple que pr\'ec\'edemment que chaque groupe {\em projectif}\, r\'eel
continu, qui poss\`ede la libre mobilit\'e dans l'infinit\'esimal en
tout les points r\'eels de l'espace, co\"{\i}ncide avec le groupe
projectif r\'eel d'une surface imaginaire du second degr\'e, dont
l'\'equation est r\'eelle.

Pour le plan, 
nous avons d\'emontr\'e cela dans le \S~97.
Afin de d\'emontrer cette proposition en g\'en\'eral, nous supposons
que nous l'avons d\'emontr\'ee dans l'espace \`a $n - 1 \geqslant 2$
dimensions et nous d\'emontrons ensuite qu'elle est aussi valide dans
l'espace \`a $n$ dimensions.

Si donc la proposition est d\'emontr\'ee pour l'espace \`a $n-1$
dimension, il en d\'ecoule imm\'ediatement que chaque groupe projectif
r\'eel $\gamma$ de $R_n$ qui poss\`ede la libre mobilit\'e dans
l'infinit\'esimal en un point de position g\'en\'erale, a $\frac{
1}{2}\, n ( n+1)$ param\`etres et qu'il peut \^etre ramen\'e, {\em
via}\, une transformation projective r\'eelle, \`a la forme:
\[
\aligned
&
p_\mu
+
\sum_{j\,k}^{1\dots n}\,
\alpha_{\mu\,j\,k}\,x_j\,p_k
+
\sum_\tau^{1\dots n}\,
\beta_{\mu\,\tau}\,x_\tau\,U
\\
&
\quad
x_\mu\,p_\nu-x_\nu\,p_\mu
+
\sum_\tau^{1\dots n}\,
\gamma_{\mu\,\nu\,\tau}\,x_\tau\,U
\\
&
\quad\quad
{\scriptstyle{(\mu,\,\nu\,=\,1\,\cdots\,n;\ \ \
\gamma_{\mu\,\nu\,\tau}\,+\,\gamma_{\nu\,\mu\,\tau}\,=\,0)}}.
\endaligned
\]
Mais il s'ensuit de l\`a par des calculs, qui ne diff\`erent
pratiquement pas de ceux des pages~290 sq., que $\gamma$ peut
recevoir, par une transformation projective r\'eelle, la forme:
\[
p_\mu+cx_\mu\,U,
\quad\quad
x_\mu\,p_\nu-x_\nu\,p_\mu
\quad\quad\quad
{\scriptstyle{(\mu,\,\nu\,=\,1\cdots\,n)}}.
\]
Si maintenant en particulier $\gamma$ doit poss\'eder la libre
mobilit\'e dans l'infinit\'esimal en tous les points r\'eels de
l'espace, alors $c$ doit \^etre positif, et par suite dans ce
cas $\gamma$ consiste, apr\`es un choix appropri\'e des variables, en
toutes les transformations projectives r\'eelles qui laissent
invariante la vari\'et\'e:
\[
x_1^2+\cdots+x_n^2+1
=
0.
\]

Ainsi, nous avons d\'emontr\'e directement que notre proposition vaut
pour $R_n$ si elle vaut pour $R_{ n-1}$ et puisqu'elle est vraie dans le
plan, elle est donc g\'en\'eralement d\'emontr\'ee.

Consid\'erons maintenant un groupe quelconque de $R_n$ qui, dans le
voisinage d'un point r\'eel en position g\'en\'erale, poss\`ede la
libre mobilit\'e dans l'infinit\'esimal. On obtient ensuite de la
m\^eme mani\`ere qu'aux pages~\pageref{482} sq. que $G$ est transitif,
qu'il a $\frac{ 1}{ 2}\, n ( n+1)$ param\`etres et qu'il peut \^etre
rapport\'e \`a la forme:
\[
p_\nu+\cdots,
\quad\quad
x_\mu\,p_\nu-x_\nu\,p_\mu
+\cdots
\quad\quad\quad
{\scriptstyle{(\mu,\,\nu\,=\,1\cdots\,n)}},
\]
{\em via}\, une transformation r\'eelle.

Il d\'ecoule de l\`a qu'apr\`es fixation de l'origine des
coordonn\'ees, les diff\'erentielles $dx_1 \dots dx_n$ sont
transform\'ees de telle sorte que l'expression diff\'erentielle:
\def\theequation{10}\begin{equation}
dx_1^2+\cdots+dx_n^2
\end{equation}
reste invariante. Maintenant, comme l'origine des coordonn\'ees est un
point de position g\'en\'erale et comme notre groupe est transitif,
nous obtenons donc g\'en\'eralement pour chaque point de position
g\'en\'erale une expression diff\'erentielle de ce type qui 
lui est attachée \deutsch{zugeordnet}
et pour pr\'eciser, nous trouvons l'expression diff\'erentielle
attach\'ee au point $x_1^0 \dots x_n^0$, lorsque nous soumettons
l'expression~\thetag{ 10} \`a une transformation quelconque:
\[
x_\nu'
=
f_\nu(x_1\dots x_n)
\quad\quad\quad\quad
{\scriptstyle{(\nu\,=\,1\,\cdots\,n)}},
\]
qui envoie l'origine des coordonn\'ees sur le point $x_1^0 \dots
x_n^0$, et quand nous effectuons apr\`es 
la substitution: $x_\nu ' = x_\nu^0$ dans les coefficients de
l'expression ainsi obtenue.

Les expressions diff\'erentielles qui sont attach\'ees
de cette mani\`ere aux points de l'espace, \label{496} 
sont transformées
par les
transformations de notre groupe exactement comme les points
eux-m\^emes, 
et de l\`a, il d\'ecoule \`a nouveau, que notre groupe laisse
invariante une certaine expression diff\'erentielle de la forme:
\[
\sum_{k\,\nu}^{1\dots n}\,
\alpha_{k\,\nu}\,
(x_1\dots x_n)\,dx_k\,dx_\nu,
\]
qui rev\^et la forme~\thetag{ 10} \`a l'origine des coordonn\'ees, et
qui par cons\'equent se transforme, pour tous les points r\'eels $x_1
\dots x_n$ de position g\'en\'erale, en une forme quadratique
constamment positive par rapport \`a $dx_1 \dots dx_n$.

\renewcommand{\thefootnote}{\fnsymbol{footnote}}
Ainsi, sans utiliser le Th\'eor\`eme~42, p.~\pageref{481}, on a
d\'emontr\'e que chaque groupe r\'eel qui poss\`ede la libre
mobilit\'e dans l'infinit\'esimal en un point de position
g\'en\'erale, laisse invariante une expression diff\'erentielle du
second degr\'e constamment positive, et on a donc d\'emontr\'e qu'avec
l'aide de l'axiome de libre mobilit\'e dans l'infinit\'esimal, on peut
d\'eduire l'axiome de Riemann sur la longueur d'un \'el\'ement courbe,
sans m\^eme pr\'esenter tous les groupes qui poss\`edent la libre
mobilit\'e dans l'infinit\'esimal\footnote[1]{\,
Lie a d\'ej\`a signal\'e cela en 1890 dans le 
Leiziger Berichten ({\em voir}\, la Remarque p.~292).}.
\renewcommand{\thefootnote}{\arabic{footnote}}

\bigskip

Pour terminer, nous voulons encore prononcer quelques mots au sujet
d'un passage de la soutenance orale de Riemann, qui nous para\^{\i}t
inintelligible dans la version d\'ej\`a pr\'esent\'ee
\deutsch{die uns in der vorliegenden Fassung unverst\"andlich ercheint}.

\noindent
Nous voulons dire le passage \`a la page~265 de ses {\em {\OE}uvres
compl\`etes} ($1^{\text{\rm \`ere}}$ \'edition), o\`u Riemann applique
ses recherches sur la d\'etermination des rapports m\'etriques d'une
grandeur $n$ fois \'etendue dans l'espace.

Le <<\,premi\`erement\,>> dans la phrase <<\,{\em Elles peuvent
premi\`erement s'exprimer} \etc\,>> pr\'esente avant tout des
difficult\'es. Ce <<\,premi\`erement\,>> signale un
<<\,deuxi\`emement\,>> \`a venir, qui vient aussi, car il vient dans
la phrase: <<\,{\em Si l'on suppose deuxi\`emement} \etc\,>>.
Mais ce <<\,deuxi\`emement\,>> ne correspond pas du tout au
<<\,premi\`erement\,>>. En effet, tandis qu'apr\`es le
<<\,premi\`erement\,>>, des conditions sont indiqu\'ees qui sont
n\'ecessaires et suffisantes pour la d\'etermination des rapports
m\'etriques de l'espace (euclidien), lorsque certaines hypoth\`eses
sont faites, Riemann ne donne pas apr\`es le <<\,deuxi\`emement\,>>,
comme on devrait s'y attendre v\'eritablement, une autre version de
ces conditions, mais il introduit une hypoth\`ese tout \`a fait
nouvelle. Les mots: <<\,{\em Si l'on suppose deuxi\`emement} 
\etc\,>> ne correspondent donc pas au <<\,premi\`erement\,>>, mais au
membre de phrase: <<\,{\em lorsqu'on admet comme hypoth\`eses \dots\,
dans ses parties infinit\'esimales}\,>>. On voit donc que le
<<\,premi\`erement\,>> n'a rien \`a voir avec le
<<\,deuxi\`emement\,>>,
\renewcommand{\thefootnote}{\fnsymbol{footnote}} et qu'il est au fond
tout \`a fait superflu. En fait, le passage tout entier devient
passablement clair\footnote[1]{\, \`A vrai dire, il resterait toujours
\`a \'etablir ce que les mots <<\,{\em Enfin on pourrait
troisi\`emement} \etc\,>> veulent dire, et nous n'y sommes pas
parvenus. Du reste, il se trouve par ailleurs aussi encore des
passages dans la soutenance orale de Riemann qui nous paraissent pour
le moins inintelligibles.} si on supprime
\renewcommand{\thefootnote}{\arabic{footnote}} le
<<\,premi\`erement\,>> entre les mots <<\,dans ses parties
infinit\'esimales\,>> et <<\,Elles peuvent premi\`erement
s'exprimer\,>>, si on supprime l'alin\'ea, et si enfin, on ajoute le
mot <<\,euclidien\,>>, comme cela a \'et\'e indiqu\'e plus haut.

D'apr\`es une communication que nous devons \`a la bont\'e de Monsieur
H. Weber, l'\'editeur des {\oe}uvres de Riemann, il n'y a aucun doute
que Riemann a v\'eritablement \'ecrit le <<\,premi\`erement\,>>; la
possibilit\'e d'une autre version est donc exclue. En cons\'equence
de cela, il ne reste aucune autre possibilit\'e que d'admettre qu'ici
se pr\'esente une inadvertance stylistique de Riemann

}

%%%%%%%%%%%%%%%%%%%%%%%%%%%%%%%%%%%%%%%%%%%%%%%%%%%%%%%%%%%%%%%%%%%%%

\newpage

% 24   :   498--521

\setcounter{footnote}{0}

$\:$
\bigskip\bigskip\bigskip

\centerline{\Large Chapitre~23.}
\label{Chapitre-23}
\thispagestyle{empty}

\bigskip

\noindent
{\large\bf
Deuxi\`eme solution du probl\`eme de Riemann-Helmholtz.}

\bigskip\medskip

Souvenons-nous des remarques effectu\'ees tout au d\'ebut de
l'introduction du pr\'ec\'edent chapitre (p.~\pageref{471} sq.). Nous
prenons maintenant ici en consid\'eration le point de vue annonc\'e
l\`a-bas, \ie nous demandons \`a nouveau quelles propri\'et\'es
sont communes au groupe des mouvements euclidiens ainsi qu'aux deux
groupes de mouvements non-euclidiens et qui font que ces trois groupes
sont remarquables parmi tous les groupes continus poss\'edant des
transformations inverses l'une de l'autre par paires. Tandis que dans
le pr\'ec\'edent chapitre, nous avons consid\'er\'e des propri\'et\'es
caract\'eristiques qui se r\'ef\`erent seulement \`a des points
infiniment voisins, nous voulons maintenant utiliser seulement des
propri\'et\'es qui concernent les rapports mutuels de points finiment
\'eloign\'es les uns des autres.

Chez Riemann, la diff\'erence entre les axiomes qui se rapportent aux
points infiniment voisins et ceux qui se rapportent aux points
finiment \'eloign\'es les uns des autres n'est pas encore faite. Les
axiomes utilis\'es par Riemann, qui \`a vrai dire ne sont pas
formul\'es express\'ement par lui, ou en tout cas ne sont cependant
formul\'es ni pr\'ecis\'ement, ni distinctement, se r\'ef\`erent en
partie \`a une esp\`ece de points, en partie \`a l'autre. L'exigence
de Riemann sur l'existence d'une longueur d'\'el\'ement courbe se
r\'ef\`ere \`a des points infiniment proches. Mais maintenant, comme
nous l'avons d\'ej\`a mis en \'evidence aux pages~\pageref{487} sq.,
l'existence d'une fonction distance entre deux points finiment
\'eloign\'es l'un de l'autre ne d\'ecoule pas encore de l'existence
d'une longueur d'\'el\'ement courbe, tant qu'on ne sait rien de plus
pr\'ecis sur la nature de la longueur d'un \'el\'ement courbe. Tout de
m\^eme, Riemann suppose aussi sans plus de
justification l'existence d'une telle fonction distance, et c'est une
hypoth\`ese qui se r\'ef\`ere \`a des points finiment \'eloign\'es les
uns des autres. Enfin, les exigences de Riemann\,\,---\,\,exprim\'ees
\`a vrai dire de mani\`ere tr\`es ind\'etermin\'ee\,\,---\,\, sur la
mobilit\'e des figures sans \'elargissement se r\'ef\`erent
\'evidemment aussi \`a des points finiment \'eloign\'es les uns des
autres, mais, \`a cause de la nature particuli\`ere de la longueur
d'un \'el\'ement courbe, que Riemann d\'eduit de la fonction distance
qu'il a suppos\'ee, ces exigences se laissent aussi imm\'ediatement
appliquer \`a des points infiniment voisins (\cf
pp.~\pageref{490} sq.).

Les axiomes de Monsieur de Helmholtz se r\'ef\`erent, dans la version
qu'il leur a donn\'ee au d\'ebut de son travail, \`a des points
finiment \'eloign\'es les uns des autres; cependant, il les applique
non seulement \`a de tels points, mais encore aussi \`a des points
infiniment voisins, un proc\'ed\'e dont nous avons suffisamment
montr\'e l'inadmissibilit\'e dans le Chapitre~21. La diff\'erence de
nature entre les deux genres d'axiomes, \ie entre ceux relatifs
aux points finiment \'eloign\'es et ceux relatifs aux points
infiniment voisins, appara\^{\i}t, dans notre exposition, pour la
premi\`ere fois distinctement \`a l'expression.
\`A vrai dire,\renewcommand{\thefootnote}{\fnsymbol{footnote}} Riemann
fait des allusions qui peuvent conduire \`a pr\'esumer qu'il a eu des
pens\'ees similaires, mais ses expressions\footnote[1]{\, {\em Voir}\,
ses {\em {\OE}uvres Compl\`etes}, $1^{\text{\rm \`ere}}$ \'Edition,
p.~267.} sont fix\'ees de mani\`ere si g\'en\'erale qu'on ne ne peut
pas trancher quant \`a ce qu'il a voulu dire.
\renewcommand{\thefootnote}{\arabic{footnote}}

\bigskip

Ce que nous avons appel\'e notre {\em premi\`ere}\, solution du
probl\`eme de Riemann-Helmholtz se r\'ef\`ere seulement \`a des points
infiniment voisins, ce que l'on peut exprimer de mani\`ere plus
pr\'ecise comme suit: cette solution r\'esout le probl\`eme de
d\'eterminer tous les groupes qui sont d\'efinis par certaines
propri\'et\'es dans l'infinit\'esimal. \`A vrai dire, il faut aboutir
\`a ce que l'exigence, d'apr\`es laquelle les mouvements doivent
former un groupe, se r\'ef\`ere \`a des points finiment \'eloign\'es
les uns des autres, et une telle exigence ne se laisse en g\'en\'eral
pas \'eluder, puisqu'il s'agit n\'ecessairement, lorsqu'on consid\`ere
les mouvements, de changements de lieu finis.

Notre {\em deuxi\`eme}\, solution traite un probl\`eme qui se
r\'ef\`ere seulement \`a des points finiment \'eloign\'es les uns des
autres; le concept de points infiniments voisins intervient ici pour
autant que nous demandons que les points s\'epar\'es restent aussi
s\'epar\'es, \label{499} et donc que deux points finiment \'eloign\'es
les uns des autres restent finiment \'eloign\'es au cours de tous les
mouvements,\renewcommand{\thefootnote}{\fnsymbol{footnote}} et ne sont
jamais transform\'es en des points inifiniment voisins\footnote[2]{\,
La solution du probl\`eme de Riemann-Helmholtz que nous avons donn\'ee
dans le Chapitre~21 sur la base des axiomes de Helmholtz partage le
m\^eme caract\`ere.
}.
\renewcommand{\thefootnote}{\arabic{footnote}}

Si, pour le probl\`eme de Riemann-Helmholtz, on pose au fondement des
axiomes qui se r\'ef\`erent \`a des points finiment \'eloign\'es les
uns des autres, alors la solution est beaucoup plus difficile que
lorsqu'on utilise des axiomes au sujet de points infiniment voisins.
C'est sur cette derni\`ere base que notre {\em premi\`ere}\, solution
({\em voir}\, le Chapitre~22) a \'et\'e achev\'ee d\'efinitivement, y
compris pour l'espace \`a $n$ dimensions; mais au contraire, {\em
notre deuxi\`eme solution}, qui est expos\'ee dans le pr\'esent
chapitre, {\em est achev\'ee d\'efinitivement, au moins en un certain
sens, seulement pour l'espace, \`a vrai dire le plus important, de
dimension trois}. Pour l'espace \`a $n$ dimensions, nous montrons
seulement qu'il est suffisant de poser certains axiomes qui se
rapportent {\em \`a des points finiment \'eloign\'es les uns des
autres}, et qui en tout cas requi\`erent moins que les axiomes
helmholtziens; mais nous ne pr\'etendons toutefois pas que pour $n >
3$ ces axiomes ne contiennent pas d'\'el\'ements superflus; nous
estimons m\^eme qu'il est vraisemblable qu'ils en contiennent.

\bigskip

Dans le \S~101, nous d\'emontrons en premier lieu une proposition
g\'en\'erale sur les groupes relativement auxquels deux points
poss\`edent un invariant. Dans le \S~102, nous traitons ensuite
l'espace \`a trois dimensions, et \`a cette occasion, la proposition
d\'emontr\'ee dans le \S~101 sera de quelque utilit\'e. Dans le
\S~103 enfin, nous r\'esolvons le cas d'un espace \`a quatre
dimensions et nous indiquons encore pour terminer comment on peut
parvenir au but dans les espaces \`a un nombre de dimensions
sup\'erieur \`a quatre.

\HEAD{Deuxième solution au problème de Riemann-Helmholtz.}{
Division\,\,V.\,\,\,Chapitre\,\,23.\,\,\,\S\,\,101.}

\sectiondritterV{\sf\S\,\,101.}
\label{S-101}
\setcounter{footnote}{0}

Pendant le compte rendu critique du travail de Helmholtz, nous avons
d\'ej\`a mentionn\'e que l'on doit \^etre extr\^emement
pr\'ecautionneux lorsqu'on veut d\'eduire quelque chose au sujet du
comportement de points infiniment voisins
\`a partir du comportement
de points finiment \'eloign\'es les uns des autres. Toutefois, nous
ne disons pas par l\`a qu'on ne peut rien d\'eduire au sujet du
comportement de ces points-l\`a \`a partir du comportement de ces
points-ci. En effet, le th\'eor\`eme suivant est vrai en tout cas.

\medskip

{\bf Th\'eor\`eme~44.}
{\em Si, relativement \`a un groupe continu de l'espace $n$ fois
\'etendu qui est constitu\'e de transformations infinit\'esimales,
deux points finiment \'eloign\'es les uns des autres $x_1 \dots x_n$
et $y_1 \dots y_n$ ont un invariant, alors aussi, deux points
infiniment voisins $x_\nu$ et $x_\nu + dx_\nu$ poss\`edent {\small\sf au
moins un} invariant relativement au groupe, ou pour s'exprimer de
mani\`ere plus pr\'ecise~{\rm :} le groupe laisse invariante au
minimum une expression de la forme}:
\[
J
\big(
x_1\dots x_n,\ \
dx_1\dots dx_n
\big).
\]

\medskip

Par souci de simplicit\'e, nous d\'emontrons tout d'abord 
ce th\'eor\`eme pour les groupes finis continus.

Soit:
\[
X_kf
=
\sum_{\nu=1}^n\,
\xi_{k\nu}(x_1\dots x_n)\,
\frac{\partial f}{\partial x_\nu}
\ \ \ \ \ \ \ \ \ \ \ \ \ 
{\scriptstyle{(k\,=\,1\,\cdots\,r)}}
\] 
un groupe quelconque \`a $r$ param\`etres; on peut poser:
\[\label{500}
Y_kf
=
\sum_{\nu=1}^n\,
\xi_{k\nu}(y_1\dots y_n)\,
\frac{\partial f}{\partial y_\nu}.
\]
Si ensuite les deux points $x_\nu$ et $y_\nu$ doivent avoir un
invariant $\Omega ( x, y)$, alors il est n\'ecessaire et suffisant que les
$r$ \'equations:
\def\theequation{1}\begin{equation}
X_kf+Y_kf
=
0
\ \ \ \ \ \ \ \ \ \ \ \ \ 
{\scriptstyle{(k\,=\,1\,\cdots\,r)}}
\end{equation}
en les $2n$ variables $x_\nu$ et $y_\nu$ aient une solution en commun,
ou, ce qui revient au m\^eme, que 
tous les $2n \times 2n$ d\'eterminants de
la matrice:
\def\theequation{2}\begin{equation}
\left\vert
\aligned
\xi_{k1}(x)\dots\,\xi_{kn}(
&
x)\ \ \
\xi_{k1}(y)\dots\,\xi_{kn}(y)\
\\
&
{\scriptstyle{(k\,=\,1,\,2\,\cdots\,r)}}
\endaligned
\right\vert
\end{equation}
s'annulent identiquement.

Si maintenant cette condition est remplie et si nous posons: $y_\nu =
x_\nu + dx_\nu$, on obtient imm\'ediatement que les $2n \times 2n$
d\'eterminants de la matrice:
\def\theequation{2'}\begin{equation}
\left\vert
\aligned
\xi_{k1}(x)\dots\,\xi_{kn}(
&
x)\ \ \
d\xi_{k1}(x)\dots\,d\xi_{kn}(x)\
\\
&
{\scriptstyle{(k\,=\,1\,\cdots\,r)}}
\endaligned
\right\vert
\end{equation}
s'annulent identiquement. Par cons\'equent, le groupe
en les $2n$ variables $x_\nu$, $x_\nu'$ provenant de
$X_1 f \dots X_r f$ par prolongation ({\em voir} le Tome~I, pp.~524
sq.):
\[
X_kf
+
\sum_{\nu=1}^n\,
\bigg\{
\sum_{\tau=1}^n\,
\frac{\partial\xi_{k\nu}(x)}{\partial x_\tau}\,x_\tau'
\bigg\}\,
\frac{\partial f}{\partial x_\nu'}
\ \ \ \ \ \ \ \ \ \ \ \ \ 
{\scriptstyle{(k\,=\,1\,\cdots\,r)}}
\]
est s\^urement intransitif, et
il laisse invariante au minimum une fonction: $\omega ( x_1 \dots
x_n, \ x_1' \dots x_n')$. Du reste, il est facile de voir que l'on
peut toujours choisir cette fonction $\omega$ de telle sorte qu'elle
n'est pas ind\'ependante de tous les $x'$. Si en effet les $r$
\'equations~\thetag{ 1} n'ont pas de solution commune ind\'ependante
de $y_1 \dots y_n$, alors les \'equations:
\def\theequation{3}\begin{equation}
X_kf
+
\sum_{\nu=1}^n\,
\bigg\{
\sum_{\tau=1}^n\,
\frac{\partial\xi_{k\nu}(x)}{\partial x_\tau}\,x_\tau'
\bigg\}\,
\frac{\partial f}{\partial x_\nu'}
=
0
\ \ \ \ \ \ \ \ \ \ \ \ \ 
{\scriptstyle{(k\,=\,1\,\cdots\,r)}}
\end{equation}
n'ont visiblement 
%%%\Fill
pas non plus de solution commune ind\'ependante de $x_1 ' \dots x_n'$.
D'autre part, si les \'equations~\thetag{ 1} ont une solution: $\Omega
( x_1 \dots x_n)$ qui est ind\'ependante de tous les $y$, alors elles
ont aussi la solution: $\Omega ( y_1 \dots y_n)$ ind\'ependante de
tous les $x$, et alors il est visible que l'expression:
\[
\sum_{\tau=1}^n\,
\frac{\partial\Omega(x_1\dots x_n)}{\partial x_\tau}\,
x_\tau'
\]
est une solution commune des \'equations~\thetag{ 3}
qui n'est pas ind\'ependante de tous les $x'$.

Il d\'ecoule de l\`a que les deux points infiniment voisins
$x_\nu$ et $x_\nu + dx_\nu$ ont toujours un invariant 
relativement au groupe $X_1 f \dots X_r f$, 
aussit\^ot que les deux points $x_\nu$ et $y_\nu$ en 
ont un, et la v\'eracit\'e de notre Th\'eor\`eme~44 est
donc d\'emontr\'ee pour tous les groupes continus finis.

Jusqu'\`a pr\'esent, pour la d\'emonstration du Th\'eor\`eme~44, nous
nous sommes limit\'es aux groupes continus finis, mais il n'est pas
difficile de voir que notre d\'emonstration s'applique en
g\'en\'eral \`a tous les groupes continus finis ou infinis qui sont
constitu\'es de transformations infinit\'esimales.

Si en effet $Xf$ est la transformation infinit\'esimale g\'en\'erale
d'un groupe et si $Yf$ a la m\^eme signification que ci-dessus, alors
les deux points $x_\nu$ et $y_\nu$ ont un invariant relativement \`a
ce groupe si et seulement si, dans la collection de toutes les
\'equations:
\[
Xf+Yf
=
0
\]
ne sont pas contenues plus que $2n-1$ \'equations qui sont
ind\'ependantes les unes des autres. Et si cette condition est
remplie, alors dans la collection de toutes les \'equations:
\[
Xf
+
\sum_{\nu=1}^n\,
\bigg\{
\sum_{\tau=1}^n\,
\frac{\partial\xi_{\nu}(x)}{\partial x_\tau}\,x_\tau'
\bigg\}\,
\frac{\partial f}{\partial x_\nu'}
=
0,
\]
sont contenues au maximum $2n - 1$ \'equations qui sont
ind\'ependantes les unes des autres, et par suite, les deux points
$x_\nu$ et $x_\nu + dx_\nu$ ont s\^urement un invariant.

Avec cela, notre th\'eor\`eme est d\'emontr\'e en g\'en\'eral.

\bigskip

On sait d\'ej\`a que les deux points $x_\nu$ et $y_\nu$ finiment
\'eloign\'es l'un de l'autre ont, relativement \`a un groupe continu
donn\'e, l'invariant $\Omega ( x, y)$, et ceci soul\`eve donc la
question de savoir si l'on ne peut pas d\'eduire de $\Omega$ un
invariant pour deux points infiniment voisins $x_\nu$ et $x_\nu +
dx_nu$; car ces deux points ont un invariant, comme il suit du
Th\'eor\`eme~44. Nous ne voulons pas nous engager en plus pour
r\'epondre \`a cette question, mais seulement mentionner que
l'invariant cherch\'e pour les deux points $x_\nu$ et $x_\nu + dx_\nu$
peut toujours \^etre indiqu\'e simplement quand l'expression $\Omega (
x, x+ dx)$ est d\'eveloppable en s\'erie enti\`ere ordinaire par
rapport \`a $dx_1 \dots dx_n$; on obtient en effet alors l'invariant
d\'esir\'e en recherchant la fonction homog\`ene du plus bas degr\'e
par rapport \`a $dx_1 \dots dx_n$ qui est contenue dans le
d\'eveloppement de $\Omega ( x, x + dx)$ par rapport aux puissances de
$dx_1 \dots dx_n$. Mais quand $\Omega ( x, x + dx)$ n'est pas
d\'eveloppable en s\'erie enti\`ere ordinaire par rapport \`a $dx_1
\dots dx_n$, on ne peut pas d\'eduire imm\'ediatement l'invariant
entre $x_\nu$ et $x_\nu + dx_\nu$ \`a partir de l'expression de
$\Omega ( x, x+ dx)$.

\bigskip

\label{502}
Il est souhaitable de poss\'eder un crit\`ere simple par lequel on
peut v\'erifier si deux points infiniment voisins ont, ou n'ont pas un
invariant relativement \`a un groupe. Si en effet on s'est convaincu
qu'ils n'ont pas d'invariant, alors d'apr\`es le Th\'eor\`eme~44, il
est imm\'ediatement clair que deux points finiment \'eloign\'es l'un
de l'autre n'ont pas non plus d'invariant. Nous allons montrer pour
cela comment on peut parvenir \`a un tel crit\`ere.

Consid\'erons un groupe continu quelconque $G$, fini ou infini, qui
est constitu\'e de transformations infinit\'esimales. Parmi les
transformations infinit\'esimales de $G$, qui laissent invariant un
point $x_1^0 \dots x_n^0$ de position g\'en\'erale, il y en a un
certain nombre, 
disons exactement $m \leqslant n^2$,
ind\'ependantes les unes des autres:
\def\theequation{4}\begin{equation}
\sum_{\mu\,\nu}^{1\dots n}\,
\alpha_{k\,\mu\,\nu}\,
(x_\mu-x_\mu^0)\,
\frac{\partial f}{\partial x_\nu}
+
\cdots
\ \ \ \ \ \ \ \ \ \ \ \
{\scriptstyle{(k\,=\,1\,\cdots\,m)}},
\end{equation}
qui sont du premier ordre en les $x_\mu - x_\mu^0$ et \`a partir
desquelles on ne peut d\'eduire, par combinaison lin\'eaire, aucune
transformation du second ordre, ou d'un ordre sup\'erieur, en les
$x_\mu - x_\mu^0$. Les transformations infinit\'esimales~\thetag{ 4}
sont alors visiblement constitu\'ees de telle sorte que les termes du
premier ordre, dans chaque transformation infinit\'esimale de $G$ qui
fixe le point $x_1^0 \dots x_n^0$, se laissent exprimer par
combinaison lin\'eaire \`a partir des termes du premier ordre des
transformations~\thetag{ 4}.

Il d\'ecoule de l\`a que les transformations infinit\'esimales:
\def\theequation{5}\begin{equation}
A_kf
=
\sum_{\mu\,\nu}^{1\dots n}\,
\alpha_{k\,\mu\,\nu}\,x_\mu'\,
\frac{\partial f}{\partial x_\nu'}
\ \ \ \ \ \ \ \ \ \ \ \
{\scriptstyle{(k\,=\,1\,\cdots\,m)}}
\end{equation}
en les variables $x_1 ' \dots x_n'$ engendrent de leur c\^ot\'e un
groupe lin\'eaire homog\`ene \label{503} \`a $m$ param\`etres, qui est
enti\`erement d\'etermin\'e par le groupe $G$ et le point $x_1^0 \dots
x_n^0$. Ce groupe lin\'eaire homog\`ene, dont l'existence nous est
d\'ej\`a connue dans le cas d'un groupe continu fini, d'apr\`es le
Tome~I, p.~599 sq., a une signification tr\`es simple: il indique
notamment de quelle mani\`ere les points:
\[
x_\nu^0+dx_\nu
=
x_\nu^0+x_\nu'\,dt
\ \ \ \ \ \ \ \ \ \ \ \
{\scriptstyle{(k\,=\,1\,\cdots\,m)}}
\]
infiniment voisins de $x_1^0 \dots x_n^0$ sont transform\'es par le
groupe $G$, aussit\^ot que l'on fixe le point $x_1^0 \dots x_n^0$. On
v\'erifie cela imm\'ediatement, si l'on pense que les termes
du\renewcommand{\thefootnote}{\fnsymbol{footnote}} premier ordre dans
les transformations infinit\'esimales~\thetag{ 4} sont les seuls par
lesquelles sont effectivement
transform\'es\footnote[1]{\, Jusqu'\`a pr\'esent, nous
avons habituellement interpr\'et\'e $x_1 ' \dots x_n'$ comme les
coordonn\'ees homog\`enes des $\infty^{ n-1}$ \'el\'ements lin\'eaires
passant par le point $x_1^0 \dots x_n^0$, mais cependant,
l'interpr\'etation du groupe~\thetag{ 5} que nous indiquons ici,
nous l'avons d\'ej\`a mentionn\'ee aussi (\cf p.~\pageref{454}).
}
\renewcommand{\thefootnote}{\arabic{footnote}}
les points infiniment voisins du point $x_1^0 \dots x_n^0$.

Si maintenant les deux points infiniment voisins $x_\nu$ et $x_\nu +
dx_\nu$ ont, relativement au groupe $G$, l'invariant $\omega ( x_1
\dots x_n, \ dx_1 \dots dx_n)$, alors il est clair qu'apr\`es fixation
du point $x_\nu^0$, le point infiniment voisin: $x_\nu^0 + dx_\nu$
est transform\'e de telle sorte que l'expression: $\omega ( x_1^0
\dots x_n^0, \ dx_1 \dots dx_n)$ reste invariante, et donc dans ce
cas, $\omega ( x_1^0 \dots x_n^0, \ x_1' \dots x_n')$ est un invariant
du groupe~\thetag{ 5}.

Si d'un autre c\^ot\'e le groupe~\thetag{ 5} poss\`ede un invariant:
$\varphi ( x_1' \dots x_n')$, il s'ensuit que deux points infiniment
voisins quelconques ont un invariant relativement au groupe $G$. En
effet, si $G$ est intransitif et si $\chi (x_1 \dots x_n)$ est un
invariant quelconque de $G$, alors l'expression:
\[
\sum_{\nu=1}^n\,
\frac{\partial\chi(x)}{\partial x_\nu}\,dx_\nu
\]
est visiblement un invariant des deux points infiniment voisins
$x_\nu$ et $x_\nu +
dx_\nu$. Mais si $G$ est transitif, alors en s'imaginant que
toutes les transformations de $G$ s'exercent sur l'expression:
\[
\varphi(dx_1\dots dx_n),
\]
on obtient une expression compl\`etement d\'etermin\'ee de cette
esp\`ece qui est attach\'ee \`a chaque point de l'espace, et toutes ces
expressions sont \'echangeables l'une avec l'autre par les
transformations de $G$, ce qui ne veut toutefois rien dire d'autre
qu'il y a une expression:
\[
\omega(x_1\dots x_n,\
dx_1\dots dx_n)
\]
invariante par $G$, et par cons\'equent que deux points infiniment
voisins quelconques ont un invariant relativement \`a $G$ (\cf
p.~\pageref{496}).
%%%\Mathematiques 
%%%[[Analyse de l'isotropie.]]

Nous voyons \`a partir de l\`a que deux points infiniment
voisins quelconques ont un invariant relativement \`a $G$
si et seulement si le groupe lin\'eaire homog\`ene~\thetag{ 5}
d\'efini \`a l'instant en les variables $x_1 ' \dots x_n'$ poss\`ede
un invariant. Nous \'enon\c cons donc ce r\'esulat
de la mani\`ere suivante: 

\medskip
{\bf Proposition~1.}
{\em Si l'on souhaite d\'ecider si deux points infiniment voisins
quelconques ont un invariant, relativement \`a un groupe continu qui
est constitu\'e de transformations infinit\'esimales, alors il n'est
pas n\'ecessaire de conna\^{\i}tre les transformations
infinit\'esimales du groupe lui-m\^eme, mais il suffit de
conna\^{\i}tre les termes du premier ordre en les $x_\nu - x_\nu^0$,
dans les transformations infinit\'esimales du groupe qui laissent
invariant un point $x_1^0 \dots x_n^0$ en position g\'en\'erale.}

\medskip

Il est d'une importance particuli\`ere de pouvoir d\'ecider
si deux points infiniment voisins $x_\nu$ et $x_\nu + dx_\nu$
poss\`edent un invariant $\omega ( x, dx)$
qui est homog\`ene du premier ordre
en les $dx_\nu$. S'ils ont en effet un tel invariant, on peut
le consid\'erer comme longueur d'un \'el\'ement courbe,
et l'on trouve par int\'egration que chaque courbe a une certaine
longueur.

Il n'est pas difficile d'\'etablir si un tel invariant existe
ou non; \`a cette fin, on doit \'evidemment constater simplement
si le groupe~\thetag{ 5} poss\`ede un invariant qui est homog\`ene
du premier ordre en les $x_\nu'$. Mais cela est tr\`es facile.

Si la transformation:
\[
Uf
=
\sum_{\nu=1}^n\,x_\nu'\,
\frac{\partial f}{\partial x_\nu'}
\]
se laisse exprimer par combinaison lin\'eaire \`a partir des
transformations infinit\'esimales~\thetag{ 5}, alors tous les
invariants \'eventuels de~\thetag{ 5} sont homog\`enes d'ordre z\'ero,
et il n'y en aura s\^urement aucun qui est homog\`ene du premier
ordre. Si au contraire $Uf$ ne se peut pas s'exprimer par combinaison
lin\'eaire \`a partir des transformations~\thetag{ 5} et si le
groupe~\thetag{ 5} a en g\'en\'eral un invariant, alors il en a
toujours aussi un qui est homog\`ene du premier ordre; car, 
sous les hypoth\`eses pos\'ees, les \'equations: 
\def\theequation{5'}\begin{equation}
A_kf
=
\sum_{\mu\,\nu}^{1\dots n}\,
\alpha_{k\,\mu\,\nu}\,x_\mu'\,
\frac{\partial f}{\partial x_\nu'}
=
0
\ \ \ \ \ \ \ \ \ \ \ \
{\scriptstyle{(k\,=\,1\,\cdots\,m)}}
\end{equation}
ont certainement une solution commune. Maintenant, puisque l'\'equation
$Uf = 0$ n'est pas cons\'equence de~\thetag{ 5'} et puisque les $m$
crochets $[ A_k, U ]$ s'annulent identiquement, il est clair que les
\'equations:
\[
A_1f
=
0,\dots,\
A_mf
=
0,
\ \ \ 
Uf
=
f
\]
ont \'egalement une solution commune.

On voit donc qu'\`a partir de l'existence d'un invariant de deux
points \`a distance finie on ne peut pas d\'eduire avec certitude
l'existence d'une longueur d'\'el\'ement courbe, {\em et par
cons\'equent, qu'\`a deux points finiment \'eloign\'es l'un de l'autre
peut tr\`es bien \^etre associ\'ee une fonction distance, sans que les
courbes aient une longueur}. Si en effet deux points finiment
\'eloign\'es l'un de l'autre ont un invariant, il d\'ecoule du
Th\'eor\`eme~44 que deux points infiniment voisins $x_\nu$ et $x_\nu +
dx_\nu$ ont au minimum un invariant, mais cependant, l'invariant en
question peut \^etre homog\`ene d'ordre z\'ero par rapport \`a tous
les $dx_\nu$, alors que les courbes n'ont toutefois une longueur que
lorsqu'il existe un invariant $\omega ( x, dx)$ qui est homog\`ene du
premier ordre en les $dx_\nu$.

Le groupe du plan \`a trois param\`etres:
\[
p,\quad
q,\quad
xp+yq
\]
fournit un exemple pour cela. Relativement \`a lui, deux points
finiment \'eloign\'es l'un de l'autre ont un et un seul invariant, \`a
savoir:
\[
\frac{y_2-y_1}{x_2-x_1},
\]
et pareillement, deux points infiniment voisins ont un et un seul
invariant, \`a savoir:
\[
\frac{dy}{dx}.
\]
Ainsi, bien que deux points finiment \'eloign\'es l'un de l'autre
poss\`edent une fonction distance, les courbes n'ont cependant pas de
longueur relativement \`a ce groupe.

Inversement, nous avons vu pr\'ec\'edemment qu'\`a partir de
l'existence d'un invariant $\omega ( x, dx)$ qui est homog\`ene du
premier ordre en $dx$, on ne peut en aucune fa\c con d\'eduire
l'existence d'une fonction distance entre deux points finiment
\'eloign\'es l'un de l'autre ({\em voir}\, p.~\pageref{487} sq.).
Nous pouvons maintenant exprimer aussi ce dernier fait de la mani\`ere
suivante: {\em M\^eme si les courbes ont une longueur, ce n'est
pourtant pas pour cette raison que deux points finiment \'eloign\'es
l'un de l'autre ont une fonction distance, ou qu'ils poss\`edent 
une ligne de liaison la plus courte}.

\bigskip

Les consid\'erations pr\'ec\'edentes sur les relations entre les
invariants de points finiment \'eloign\'es les uns des autres et les
invariants de points infiniment voisins peuvent encore \^etre
compl\'et\'ees d'une mani\`ere essentielle; on peut ainsi les
g\'en\'eraliser aussi en introduisant trois points ou plus, en
\'epuisant les diff\'erentielles d'ordre deux et d'ordre sup\'erieur
en les $x_\nu$, et ainsi de suite. Nous nous r\'eservons le droit
de traiter d'une mani\`ere
d\'etaill\'ee de ces questions
\`a une autre occasion.

\bigskip

Nous nous tournons maintenant vers le probl\`eme v\'eritable du
pr\'esent chapitre, \ie vers la solution du probl\`eme de
Riemann-Helmholtz, et pour pr\'eciser, nous nous limitons tout
d'abord, comme cela a \'et\'e annonc\'e plus haut, \`a l'espace trois
fois \'etendu.

\HEAD{Deuxième solution au problème de Riemann-Helmholtz.}{
Division\,\,V.\,\,\,Chapitre\,\,23.\,\,\,\S\,\,102.}

\sectiondritterV{\sf\S\,\,102.}
\label{S-102}
\setcounter{footnote}{0}

Nous affirmons que les mouvements euclidiens et non-euclidiens de
$R_3$ sont enti\`erement caract\'eris\'es, lorsque l'on pose
les axiomes suivants:

\medskip

{\bf I)} {\em $R_3$ est une vari\'et\'e num\'erique.}\label{506}

\medskip

{\bf II)} {\em Les mouvements de $R_3$ forment un groupe r\'eel
continu, qui est engendr\'e par des transformations infinit\'esimales.}

\medskip

{\bf III)} {\em Si l'on fixe un point r\'eel quelconque: $y_1^0,
y_2^0, y_3^0$ en position g\'en\'erale, alors tous les autres points
r\'eels: $x_1, x_2, x_3$, 
en lesquels
un autre point r\'eel: $x_1^0, x_2^0, x_3^0$ peut 
encore \^etre transform\'e, satisfont une
\'equation r\'eelle de la forme:
\def\theequation{6}\begin{equation}
W
\big(
y_1^0,y_2^0,y_3^0\,;\ 
x_1^0,x_2^0,x_3^0\,;\
x_1,x_2,x_3
\big)
=
0,
\end{equation}
qui n'est pas r\'ealis\'ee pour: $x_1 = y_1^0$, $x_2 = y_2^0$, $x_3 =
y_3^0$, et qui repr\'esente {\rm (}en g\'en\'eral{\rm )} une surface
r\'eelle passant par le point: $x_1^0, x_2^0, x_3^0$.}

\medskip

{\bf IV)} {\em Autour du point: $y_1^0, y_2^0, y_3^0$, une r\'egion
finie trois fois \'etendue peut \^etre d\'elimit\'ee de telle sorte
qu'apr\`es fixation du point: $y_1^0, y_2^0, y_3^0$, chaque autre
point r\'eel: $x_1^0, x_2^0, x_3^0$ de la r\'egion peut encore \^etre
transform\'e contin\^ument en chaque autre point r\'eel appartenant
\`a la r\'egion qui satisfait l'\'equation~\thetag{ 6} et pouvant être
reli\'e au point: $x_1^0, x_2^0, x_3^0$ par une s\'erie de points
continue et irr\'eductible.}

%%%\Mathematiques
%%%[[{\sc irreducible}~??]]
%%%\Fill
%%%[[Am\'eliorer la traduction]]

Si l'on raye les deux mots mis entre parenth\`eses dans l'Axiome~III,
ces axiomes suffisent s\^urement aussi; ce qui suit le montre, quand
on supprime partout les deux mots mis entre parenth\`ese. Si on
conserve dans l'Axiome~III et dans la suite les mots mis entre
parenth\`ese, alors la suite montre, si nous ne nous trompons pas, que
les axiomes sont encore aussi satisfaits.

%%%\Mathematiques
%%%[[Doute-t-il~??]]

\bigskip

Par $G$, nous voulons entendre un groupe continu fini ou infini
quelconque qui est engendr\'e par des transformations
infinit\'esimales et qui satisfait les Axiomes~III et~IV.

Puisque les \'equations~\thetag{ 6} sont identiquement r\'ealis\'ees
pour: $x_1 = x_1^0$, $x_2 = x_2^0$, $x_3 = x_3^0$, il y a toujours,
parmi l'ensemble de tous les points r\'eels: $x_1, x_2, x_3$ qui
satisfont l'\'equation~\thetag{ 6}, une famille continue non divis\'ee
en plusieurs parties, dans laquelle est contenu le point: $x_1^0,
x_2^0, x_3^0$. D'apr\`es l'Axiome~III, cette famille de points forme
(en g\'en\'eral) une surface passant par le point: $x_1^0, x_2^0,
x_3^0$ (il est cependant imaginable aussi qu'elle se r\'eduise \`a une
courbe passant par ce point, voire \`a ce point lui-m\^eme). Quoi
qu'il en soit, nous voulons l'appeler une {\sl pseudosph\`ere}\,
appartenant au groupe $G$, et pour pr\'eciser, une pseudosph\`ere de
centre: $y_1^0, y_2^0, y_3^0$ (\cf p.~\pageref{402-0}).

Notre Axiome~III exprime maintenant visiblement qu'{\em une
pseudosph\`ere de centre: $y_1^0, y_2^0, y_3^0$ ne passe jamais par
son centre}. En cela r\'eside le fait que deux points s\'epar\'es
restent aussi s\'epar\'es, quand toutes les transformations de $G$
agissent, et donc que deux points finiment \'eloign\'es l'un de
l'autre ne peuvent jamais \^etre transform\'es en deux points
infiniment voisins (\cf p.~\pageref{499}).

Mais notre Axiome~IV peut \`a pr\'esent \^etre exprim\'e de la
mani\`ere suivante: autour du point: $y_1^0, y_2^0, y_3^0$, une
r\'egion finie trois fois \'etendue peut \^etre d\'elimit\'ee de telle
sorte qu'apr\`es fixation du point: $y_1^0, y_2^0, y_3^0$, chaque
autre point r\'eel de la r\'egion peut se mouvoir de mani\`ere
compl\`etement libre sur la pseudosph\`ere centr\'ee en le point:
$y_1^0, y_2^0, y_3^0$ qui passe par lui.

\medskip

{\small
Avant de passer \`a la d\'etermination de tous les groupes $G$ qui
satisfont nos axiomes, nous voulous encore comparer en peu de mots nos
axiomes aux axiomes helmholtziens.

Dans ses axiomes, Monsieur de Helmholtz demande r\'eellement plus que
nous, m\^eme quand on fait compl\`etement 
abstraction de son axiome de monodromie.
Comme cela a \'et\'e d\'emontr\'e auparavant, il d\'ecoule en effet de
ses trois premiers axiomes que les mouvements forment 
en g\'en\'eral
un groupe
continu {\em fini}, qui est engendr\'e par des transformations
infinit\'esimales; toutefois, nous demandons simplement dans notre
Axiome~II que les mouvements forment un groupe continu engendr\'e par
des transformations infinit\'esimales, donc nous n'excluons pas depuis
le d\'ebut la possibilit\'e que le groupe soit {\em infini}. D'autre
part, nos Axiomes~III et~IV exigent essentiellement moins que le
troisi\`eme axiome helmholtzien.

En effet, alors que Monsieur de Helmholtz exige que chaque point soit
mobile d'une mani\`ere parfaitement libre, tant qu'il n'est pas
restreint par les \'equations qui existent entre lui et les points
mobiles restants, nous demandons seulement qu'apr\`es fixation d'{\em
un}\, point, chaque autre point soit parfaitement libre, tant qu'il
n'est pas limit\'e par les \'equations existant entre lui et le point
fix\'e. \`A vrai dire, nous ajoutons encore l'exigence apparemment
nouvelle qu'une pseudosph\`ere ne passe jamais par son centre, mais
cette exigence n'est rien de plus qu'une version pr\'ecise de ce que
Monsieur de Helmholtz a aussi demand\'e implicitement. En effet, ce
dernier exige en particulier qu'apr\`es fixation du point: $y_1^0,
y_2^0, y_3^0$, chaque autre point puisse se mouvoir d'une mani\`ere
compl\`etement libre sur la pseudosph\`ere passant par lui et
centr\'ee en le point: $y_1^0, y_2^0, y_3^0$. Si maintenant une des
pseudosph\`eres en question passait par son centre, alors, apr\`es
fixation de: $y_1^0, y_2^0, y_3^0$, aucun point $P$ de cette
pseudosph\`ere ne pourrait se mouvoir de mani\`ere parfaitement libre
sur sa pseudosph\`ere, car parmi les transformations
non-d\'eg\'en\'er\'ees de $G$ qui laissent au repos le point: $y_1^0,
y_2^0, y_3^0$, il n'y en a aucune qui transforme le point $P$ en le
point invariant: $y_1^0, y_2^0, y_3^0$.

}

\medskip

Venons-en maintenant \`a la d\'etermination de tous les groupes qui
satisfont nos axiomes.

Nous d\'emontrons en premier lieu que chaque groupe de cette esp\`ece
{\em est transitif}\, et que {\em deux points finiment \'eloign\'es
l'un de l'autre ont un et un seul invariant relativement \`a $G$}.

Si $G$ \'etait intransitif, il laisserait s\^urement invariantes
$\infty^1$ surfaces r\'eelles dont les \'equations pourraient \^etre
rapport\'ees, au moyen d'une transformation ponctuelle r\'eelle, \`a
la forme: $x_3 = \text{\rm const.}$ Si maintenant nous fixions un
point r\'eel: $y_1^0, y_2^0, y_3^0$ en position g\'en\'erale, alors
chaque autre point: $x_1^0, x_2^0, x_3^0$ pourrait \'evidemment
encore occuper seulement les positions qui satisfont l'\'equation:
$x_3 = x_3^0$, et par cons\'equent, l'\'equation~\thetag{ 6} serait
n\'ecessairement de la forme: $x_3 = x_3^0$. Donc la pseudosph\`ere
centr\'ee en le point: $y_1^0, y_2^0, y_3^0$ serait repr\'esent\'ee
par l'\'equation: $x_3 = \text{\rm const.}$, \`a la suite de quoi une
de ces pseudosph\`eres passerait par son centre, ce qui est
exclu. Ainsi, $G$ ne peut pas \^etre intransitif, mais il doit au
contraire \^etre {\em transitif}.

En outre, il d\'ecoule de l'Axiome~IV qu'apr\`es fixation d'un point
r\'eel en position g\'en\'erale, chaque autre point r\'eel peut
occuper (en g\'en\'eral) encore exactement $\infty^2$ positions
diff\'erentes. Ainsi, il est clair que par l'action de $G$, chaque
paire de points de $R_3$ ne peut pas \^etre transform\'ee en chaque
autre paire de points, mais plut\^ot: une paire de points re\c coit,
{\em via}\, l'action de $G$, au maximum $\infty^5$ positions
diff\'erentes. Par cons\'equent, si $Xf$ est une transformation
infinit\'esimale quelconque de $G$, et si $Yf$ a la m\^eme
signification qu'\`a la page~\pageref{500}, alors il n'y a s\^urement,
parmi l'ensemble des \'equations: $Xf + Yf = 0$, pas plus que cinq
\'equations qui sont ind\'ependantes les unes des autres, donc ces
\'equations ont au minimum une solution en commun: $\Omega ( x_1,
x_2, x_3; \ y_1, y_2, y_3)$, c'est-\`a-dire que les deux points
$x_\nu$ et $y_\nu$ ont en tout cas un invariant relativement \`a $G$,
\`a savoir l'invariant $\Omega ( x, y)$.

De la transitivit\'e de $G$ il d\'ecoule maintenant imm\'ediatement
que les deux points $x_\nu$ et $y_\nu$ ont seulement un invariant
relativement \`a $G$; en effet, s'ils avaient deux ou m\^eme trois
invariants, alors chaque paire de points occuperait au maximum
$\infty^4$ positions diff\'erentes par l'action de $G$, et par
cons\'equent, apr\`es fixation d'un point en position g\'en\'erale,
chaque autre point pourrait se mouvoir au plus sur une courbe, ce qui
se trouve \^etre en contradiction avec l'Axiome~IV.

Ainsi, $G$ est r\'eellement transitif et deux points $x_\nu$ et
$y_\nu$ finiment \'eloign\'es l'un de l'autre ont, relativement \`a
$G$, seulement un invariant: $\Omega ( x, y)$. Nous pouvons conclure
de l\`a que l'\'equation~\thetag{ 6} peut \^etre ramen\'ee \`a la
forme:
\def\theequation{6'}\begin{equation}
\Omega
\big(
x_1,x_2,x_3\,;\
y_1^0,y_2^0,y_3^0
\big)
=\Omega
\big(
x_1^0,x_2^0,x_3^0\,;\
y_1^0,y_2^0,y_3^0
\big).
\end{equation}

\bigskip

Nous d\'emontrons maintenant que notre groupe est {\em r\'eel-primitif}.

\smallskip

(Puisque l'ensemble de toutes les pseudosph\`eres de centre: $y_1^0,
y_2^0, y_3^0$ est repr\'esent\'e par {\em une}\, \'equation et
puisque, parmi ces pseudosph\`eres, se trouvent $\infty^1$ surfaces
r\'eelles, il est clair que seulement un nombre discret de
pseudosph\`eres de centre: $y_1^0, y_2^0, y_3^0$ peut se r\'eduire
\`a des courbes ou \`a des points.) Comme (de plus) aucune des
pseudosph\`eres de centre: $y_1^0, y_2^0, y_3^0$ ne peut passer par
son centre, il para\^{\i}t clair qu'apr\`es fixation du point:
$y_1^0, y_2^0, y_3^0$, aucune courbe ou surface passant par ce point
ne peut rester au repos,
%%%\Mathematiques 
%%%[[Pas clair pour moi.]]
et par cons\'equent, $G$ est effectivement r\'eel-primitif.

\`A pr\'esent, il est facile aussi de d\'emontrer que $G$ est fini et
pour pr\'eciser, qu'il a au plus {\em six param\`etres}.

Tout d'abord, il est en effet certain qu'au moins $\infty^2$
pseudosph\`eres diff\'erentes appartiennent \`a $G$. S'il n'y avait en
effet que $\infty^1$ pseudosph\`eres, alors il ne passerait par chaque
point qu'une pseudosph\`ere et chaque point serait le centre de la
pseudosph\`ere passant par lui;
%%%\Mathematiques 
%%%[[Pas clair pour moi.]]
mais ceci ne doit pas \^etre.

Ensuite, soit $P_1$, $P_2$, $P_3$ et $P$ quatre points r\'eels
quelconques qui sont mutuellement en position g\'en\'erale. Puisqu'il
y a $\infty^2$ pseudosph\`eres diff\'erentes, les deux pseudosph\`eres
centr\'ees en $P_1$ et en $P_2$ passant par $P$ se coupent
n\'ecessairement en une courbe r\'eelle $C$ passant par $P$. Si
maintenant $C$ se trouvait aussi dans la pseudosph\`ere centr\'ee en
$P_3$ passant par $P$, alors toutes les pseudosph\`eres passant par
$P$ se couperaient g\'en\'eralement en la courbe $C$; si donc l'on
fixait $P$, la courbe $C$ passant par $P$ resterait elle aussi
n\'ecessairement au repos.
%%%\Fill
Mais nous avons d\'emontr\'e \`a l'instant que cela est impossible. On
obtient donc que les trois pseudosph\`eres de centres $P_1$, $P_2$ et
$P_3$ qui passent par $P$ ont seulement le point $P$ en commun, mais
n'ont pas en commun de courbe r\'eelle passant par ce point. Si
maintenant nous fixons les trois points $P_1$, $P_2$ et $P_3$, alors
$P$ peut en tout cas se mouvoir seulement sur la
vari\'et\'e-intersection 
\deutsch{Schnittmannigfaltigkeit} des trois
pseudosph\`eres passant par $P$ centr\'ees en $P_1$, $P_2$ et $P_3$,
et comme cette vari\'et\'e-intersection ne consiste qu'en le point
$P$, ce point $P$ doit rester au repos. En d'autres termes: d\`es
qu'on fixe trois points r\'eels qui sont mutuellement en position
g\'en\'erale, tous les points de $R_3$ restent g\'en\'eralement au
repos. Maintenant, puisque $G$ est transitif et puisqu'apr\`es
fixation d'un point, chaque autre point r\'eel en position
g\'en\'erale peut encore d\'ecrire une surface, la fixation de trois
point revient \`a imposer au plus six conditions. Donc le groupe $G$
est fini, et pour pr\'eciser, il a au plus six param\`etres.

Si nous fixons un point r\'eel $P$ en position g\'en\'erale, alors la
vari\'et\'e projective des $\infty^2$ \'el\'ements lin\'eaires passant
par $P$ se transforme par l'action d'un groupe projectif r\'eel $g$,
qui a visiblement au plus trois param\`etres. Mais maintenant,
d'apr\`es notre d\'etermination de tous les groupes projectifs r\'eels
du plan ({\em voir}\, p.~106 sq. et p.~380 sq.), nous obtenons que
chaque groupe projectif r\'eel du plan qui a moins de trois
param\`etres laisse au repos au moins un point r\'eel. Si donc $g$
avait moins de trois param\`etres, il laisserait au repos au minimum
un \'el\'ement lin\'eaire r\'eel passant par $P$, et il d\'ecoulerait
de l\`a que $G$ serait r\'eel-imprimitif, alors qu'il doit cependant
\^etre r\'eel-primitif. Par cons\'equent $g$ a exactement trois
param\`etres, \`a la suite de quoi $G$ a exactement six param\`etres.
Ainsi:

\medskip

{\em Si un groupe r\'eel $G$ satisfait nos Axiomes~III et~IV, alors il
est r\'eel-primitif \`a six param\`etres; en outre, deux points
finiment \'eloign\'es l'un de l'autre ont, relativement \`a
son action,
seulement un invariant; et par ailleurs, $G$ est constitu\'e de telle
sorte que les $\infty^2$ \'el\'ements lin\'eaires r\'eels passant par
chaque point r\'eel fix\'e en position g\'en\'erale sont transform\'es
par l'action d'un groupe \`a trois param\`etres.}

\medskip

Avant toute chose, nous devons maintenant \'etablir quelles sont les
diff\'erentes formes que peut prendre le groupe $g$ d\'efini
ci-dessus. Nous savons que $g$ a trois param\`etres et qu'il ne laisse
invariant aucun \'el\'ement lin\'eaire. Si donc nous rapportons
projectivement les $\infty^2$ \'el\'ements lin\'eaires r\'eels passant
par $P$ aux points r\'eels d'un plan, alors $g$ se transforme en un
groupe projectif r\'eel $\mathfrak{ g}$ \`a trois param\`etres de ce
plan, et pour pr\'eciser, en un groupe par l'action duquel aucun point
r\'eel ne reste invariant. Mais d'apr\`es les pages~106 et~384, chaque
tel groupe $\mathfrak{ g}$ 
est semblable, {\em via}\, une transformation projective
r\'eelle du plan, \label{511} soit au groupe projectif r\'eel d'une conique
non-d\'eg\'en\'er\'ee qui est repr\'esent\'ee par une \'equation
r\'eelle, soit au groupe:
\[
\mathfrak{p},\quad
\mathfrak{q},\quad
\mathfrak{y}\,\mathfrak{p}
-
\mathfrak{x}\,\mathfrak{q}
+
\mathfrak{c}\,
(\mathfrak{x}\,\mathfrak{p}+\mathfrak{y}\,\mathfrak{q}).
\]
Par cons\'equent, $g$ laisse invariante soit une conique
non-d\'eg\'en\'er\'ee du second degr\'e form\'ee d'\'el\'ements
lin\'eaires qui est repr\'esent\'ee par une \'equation r\'eelle, soit
il laisse invariant un \'el\'ement de surface r\'eel et dans cet
\'el\'ement, deux \'el\'ements lin\'eaires conjugu\'es.

Lorsque le deuxi\`eme cas se produit, \`a chaque point r\'eel en
position g\'en\'erale est associ\'e par $G$ un \'el\'ement de surface
r\'eel invariant passant par ce point, donc $G$ laisse invariante une
\'equation de Pfaff r\'eelle, qui naturellement, \`a cause de la
primitivit\'e de $G$, ne doit pas \^etre int\'egrable. Si nous nous
imaginons cette \'equation rapport\'ee \`a la forme: $dz - ydx$ par
une transformation ponctuelle r\'eelle de $R_3$, alors $G$ se
transforme en un groupe r\'eel \`a six param\`etres $G'$, par l'action
duquel l'\'equation: $dz - y dx = 0$ reste invariante. De plus, nous
pouvons faire en sorte, {\em via}\, une transformation ponctuelle
r\'eelle \`a travers laquelle l'\'equation de Pfaff en question reste
invariante, que l'origine des coordonn\'ees: $x = y = z = 0$ soit un
point de position g\'en\'erale relativement \`a l'action de $G'$ ({\em
voir}\, le Tome~II, p.~402 sq.).

Si nous fixons l'origine des coordonn\'ees, alors l'\'el\'ement de
surface: $dz = 0$ reste invariant, ainsi que deux \'el\'ements
lin\'eaires imaginaires conjugu\'es contenus en lui; nous pouvons
toujours, {\em via}\, une transformation r\'eelle \`a travers laquelle
l'origine des coordonn\'ees et l'\'equation: $dz - ydx = 0$ restent
invariantes, transformer ces deux \'el\'ements lin\'eaires en les
\'el\'ements lin\'eaires:
\def\theequation{7}\begin{equation}
dz 
=
dx+idy
=0,
\quad\quad\quad\quad
d\bar z
=
dx-idy
=
0.
\end{equation}
Maintenant, puisque $G'$ est transitif, et puisque, apr\`es fixation
de l'origine des coordonn\'ees, les \'el\'ements lin\'eaires passant
par l'origine sont transform\'es par l'action d'un groupe \`a trois
param\`etres, $G'$ doit contenir trois transformations
infinit\'esimales du premier ordre en les $x, y, z$ dans le voisinage
de l'origine des coordonn\'ees, mais au contraire, aucune
transformation du second ordre, ou d'un ordre sup\'erieur. Mais par
ailleurs, nous connaissons toutes les transformations
infinit\'esimales qui laissent invariante l'\'equation: $dz - y dx =
0$ et qui sont du premier ordre en les $x, y, z$ ({\em voir}\, le
Tome~II, p.~405), et nous reconnaissons imm\'ediatement que parmi ces
transformations, il y en a seulement quatre par lesquelles les deux
\'el\'ements lin\'eaires~\thetag{ 7} restent invariants et \`a partir
desquelles on ne peut d\'eduire, par combinaison lin\'eaire, aucune
transformation du second ordre ou d'un ordre sup\'erieur, et ce sont
les quatre:
\[
yp-xq+\cdots,
\quad\quad
xp+yq+2zr+\cdots,
\quad\quad
zp+\cdots,
\quad\quad
zq+\cdots,
\]
o\`u les termes supprim\'es sont d'ordre deux ou sup\'erieur par
rapport \`a $x, y, z$. Si nous tenons compte du fait que notre groupe
$G'$ contient pr\'ecis\'ement trois transformations infinit\'esimales
ind\'ependantes du premier ordre et que celles-ci ne doivent laisser
au repos aucun \'el\'ement lin\'eaire passant par l'origine des
coordonn\'ees, alors nous trouvons facilement que les transformations
infinit\'esimales de $G'$, qui laissent au repos l'origine des
coordonn\'ees, ont la forme:
\[
yp-xq+c\,(xp+yq+2zr)+\cdots,
\quad\quad
zp+\cdots,
\quad\quad
zq+\cdots.
\]
En outre, en tant que groupe transitif, $G'$ contient naturellement
encore trois transformations infinit\'esimales d'ordre z\'ero de la
forme:
\def\theequation{8}\begin{equation}
p+\cdots,
\quad\quad\quad
q+\cdots,
\quad\quad\quad
r+\cdots.
\end{equation}

Le groupe lin\'eaire homog\`ene, que $G'$ associe \`a l'origine des
coordonn\'ees (\cf\, p.~\pageref{503}), est \'evidemment de la
forme:
\def\theequation{9}\begin{equation}
y'p'-x'q'+c\,(x'p'+y'q'+2z'r'),
\quad\quad
z'p',
\quad\quad
z'q'.
\end{equation}
Mais maintenant, deux points finiment \'eloign\'es l'un de l'autre
doivent avoir un invariant par rapport \`a $G'$, donc d'apr\`es les
pages~\pageref{502} sq., le groupe lin\'eaire homog\`ene~\thetag{ 9}
doit poss\'eder un invariant, ce qui ne peut manifestement se produire
que si $c$ s'annule, et par cons\'equent on obtient qu'en dehors des
transformations~\thetag{ 8}, $G'$ en contient encore trois de la
forme:
\def\theequation{10}\begin{equation}
yp-xq+\cdots,
\quad\quad
zp+\cdots,
\quad\quad
zq+\cdots,
\end{equation}
et chaque transformation infinit\'esimale de $G'$ peut \^etre
exprim\'ee par combinaison lin\'eaire \`a partir de ces six-l\`a.

Rappelons maintenant \`a ce sujet que les transformations
infinit\'esimales de $G$ peuvent aussi \^etre interpr\'et\'ees comme
transformations de contact infinit\'esimales du plan $x, z$,
%%%\Mathematiques
%%%[[How~??]]
et qu'\`a chacune de ces transformations de contact infinit\'esimales
est associ\'ee une fonction caract\'eristique par laquelle elle est
parfaitement d\'etermin\'ee. En particulier, les fonctions
caract\'eristiques des transformations~\thetag{ 10} s'\'enoncent de la
mani\`ere suivante:
\[
x^2+y^2+\cdots,
\quad\quad
yz+\beta_3(x,y)+\cdots,
xz+\alpha_3(x,y)+\cdots,
\]
o\`u $\alpha_3$ et $\beta_3$ sont certaines fonctions enti\`erement
homog\`enes de degr\'e trois par rapport \`a $x$ et \`a $y$, et o\`u
les termes supprim\'es sont \`a chaque fois d'un rang en $x, y, z$ qui
est sup\'erieur \`a celui des termes qui sont \'ecrits ({\em voir}\,
le Tome~II, p.~526 sq. et p.~532 sq.). Mais maintenant, $G'$ doit
n\'ecessairement contenir aussi la transformation infinit\'esimale
dont la fonction caract\'eristique poss\`ede la forme:
\def\theequation{11}\begin{equation}
\big\{
yz+\beta_3+\cdots,\ \
xz+\alpha_3+\cdots
\big\}
=
z^2
+
\lambda_2(x,y)\,z
+
\lambda_4(x,y)
+
\cdots,
\end{equation}
%%%\Mathematiques
%%%[[Note sur le crochet de Poisson.]]
o\`u nous entendons par $\lambda_2$ et par $\lambda_4$ des fonctions
enti\`erement homog\`enes de degr\'e 2 et 4 par rapport \`a $x$ et $y$
({\em loc. cit.}, p.~321 et p.~526 sq.). La seule transformation
infinit\'esimale dont la fonction caract\'eristique poss\`ede la
forme~\thetag{ 11} serait du second ordre en $x, y, z$, alors que $G'$
ne contient cependant absolument aucune transformation
infinit\'esimale du second ordre. Par cons\'equent, nous parvenons au
r\'esultat qu'il n'y a pas du tout de groupe $G'$ \`a six param\`etres
ayant la constitution 
consid\'er\'ee ici, et donc que le deuxi\`eme des deux
cas distingu\'es \`a la page~\pageref{511} ne peut pas du tout se
r\'ealiser.

\medskip

{\small
Si l'on veut \'eviter les calculs avec les fonctions
caract\'eristiques, on peut aussi proc\'eder de la mani\`ere
suivante: $G'$ laisse invariante l'\'equation de Pfaff: $dz - ydx =
0$ et, quand on fixe un point r\'eel en position g\'en\'erale, alors
deux \'el\'ements lin\'eaires imaginaires conjugu\'es restent au repos
dans l'\'el\'ement lin\'eaire que l'\'equation: $dz - ydx = 0$
attache \`a ce point. Il d\'ecoule de l\`a que $G'$, interpr\'et\'e
comme groupe de transformations de contact du plan $x, z$, laisse
invariantes deux \'equations diff\'erentielles ordinaires conjugu\'ees
du second ordre.
%%%\Mathematiques
%%%[[?? A l'air tr\`es int\'eressant~!!]]
Ainsi, gr\^ace \`a une transformation de contact (imaginaire) du plan,
nous pouvons parvenir \`a ce que les courbes int\'egrales de l'une de
ces deux \'equations diff\'erentielles soient transform\'ees en les
points du plan. Nous pouvons donc nous supposer les variables $x, y,
z$ choisies depuis le d\'ebut de telle sorte que $G'$, interpr\'et\'e
comme groupe de transformations de contact du plan $x, z$, ne consiste
qu'en des transformations ponctuelles prolong\'ees, et qui en outre
laisse invariante encore une \'equation diff\'erentielle ordinaire du
second ordre. Mais maintenant, nous avons vu \`a la page~76 que chaque
groupe continu fini de transformations ponctuelles du plan qui laisse
invariante une \'equation diff\'erentielle ordinaire du second ordre
est semblable, {\em via}\, une transformation ponctuelle, \`a un
groupe projectif du plan. Par cons\'equent, nous pouvons admettre que
$G'$, lorsqu'en g\'en\'eral il existe, provient par prolongation d'un
groupe projectif du plan \`a six param\`etres; avec cela, nous n'avons
naturellement pas besoin de consid\'erer comme diff\'erents l'un de
l'autre deux groupes projectifs qui, \`a l'int\'erieur du groupe
projectif g\'en\'eral, 
sont conjugués l'un à l'autre 
\deutsch{gleichberechtigt sind}, ou sont
dualistiques l'un de l'autre;
nous pouvons donc supposer que $G'$ provient par prolongation
du groupe lin\'eaire g\'en\'eral: 
\[
p,\quad
r,\quad
xp,\quad
zp,\quad
xr,\quad
zr
\]
du plan, et ainsi, que dans les
variables choisies $x, y, z$, il poss\`ede la forme:
\[
p,\quad
r,\quad
xp-yq,\quad
zp-y^2q,\quad
xr+q,\quad
zr+yq.
\]
Mais nous avons d\'ej\`a d\'emontr\'e aux pages~\pageref{415} sq. au
sujet de ce groupe que relativement \`a lui, deux points finiment
\'eloign\'es l'un de l'autre n'ont aucun invariant. Ainsi, on obtient
\`a nouveau qu'un groupe $G'$ de la nature ici demand\'ee n'existe pas
du tout.

}

\medskip

Il reste encore \`a traiter le premier des deux cas distingu\'es \`a
la page~\pageref{511}.

Si ce cas se produit, les $\infty^2$ \'el\'ements lin\'eaires passant
par chaque point r\'eel fix\'e en position g\'en\'erale sont
transform\'es par l'action d'un groupe \`a trois param\`etres, et pour
pr\'eciser, de telle fa\c con qu'un c\^one non-d\'eg\'en\'er\'e du
second degr\'e constitu\'e d'\'el\'ements lin\'eaires reste invariant.
Par cons\'equent, $G$ laisse invariante une \'equation r\'eelle de la
forme:
\[
\sum_{\mu\,\nu}^{1,2,3}\,
\alpha_{\mu\nu}\,
(x_1,x_2,x_3)\,dx_\mu\,dx_\nu
=
0,
\]
dont le d\'eterminant ne s'annule pas. Maintenant, puisque notre $G$
est 
\`a six param\`etres, on obtient, gr\^ace au Th\'eor\`eme~35
p.~391, qu'il 
peut \^etre transform\'e, {\em via}\, une transformation
ponctuelle r\'eelle de $R_3$, soit en le groupe des mouvements
euclidiens, soit en l'un des deux groupes des mouvements
non-euclidiens, soit en l'un des deux groupes suivants:
\[
\aligned
p_1,\ \
&
p_2,\ \ 
p_3,\quad
x_1p_2-x_2p_1,\quad
x_2p_3+x_3p_2,\quad
x_3p_1+x_1p_3;
\\
&
\quad
p_1+x_1U_1,\quad
p_2+x_2U_2,\quad
p_3-x_3U_3,
\\
&
x_1p_2-x_2p_1,\quad
x_2p_3+x_3p_2,\quad
x_3p_1+x_1p_3.
\endaligned
\]
Cependant, les deux derniers groupes ne satisfont pas notre
Axiome~III, car pour ces deux-l\`a l'origine des coordonn\'ees est un
point en position g\'en\'erale et une des pseudosph\`eres ayant pour
centre l'origine des coordonn\'ees, est la conique:
\[
x_1^2+x_2^2-x_3^2
=
0,
\]
mais cette pseudosph\`ere passe par son centre, ce qui est justement
exclu par l'Axiome~III.

\medskip

{\em Ainsi on a d\'emontr\'e que les mouvements euclidiens
et non-euclidiens sont effectivement enti\`erement caract\'eris\'es 
par les axiomes pos\'es \`a la page~\pageref{506}.}

\HEAD{Deuxième solution au problème de Riemann-Helmholtz.}{
Division\,\,V.\,\,\,Chapitre\,\,23.\,\,\,\S\,\,103.}

\sectiondritterV{\sf\S\,\,103.}
\label{S-103}
\setcounter{footnote}{0}

Nous allons maintenant montrer que les mouvements euclidiens
et non-euclidiens de $R_4$ sont enti\`erement caract\'eris\'es, 
lorsqu'on pose les axiomes suivants.

\medskip

{\bf I)} {\em $R_4$ est une vari\'et\'e num\'erique.} \label{515}

\medskip

{\bf II)} {\em Les mouvements de $R_4$ forment un groupe r\'eel continu, 
qui est engendr\'e par des transformations infinit\'esimales.}

\medskip

{\bf III)} {\em Si l'on fixe un point r\'eel quelconque: $y_1^0 \dots
y_4^0$ en position g\'en\'erale, alors tous les autres points
r\'eels: $x_1, x_2, x_3$ en lesquels un autre point r\'eel: $x_1^0
\dots x_4^0$ peut encore \^etre transform\'e, satisfont une \'equation
r\'eelle de la forme:
\def\theequation{12}\begin{equation}
W
\big(
y_1^0\dots y_4^0\,;\
x_1^0\dots x_4^0\,;\
x_1\dots x_4
\big)
=0,
\end{equation}
qui n'est pas r\'ealis\'ee pour: $x_1 = y_1^0, \dots, x_4 = y_4^0$,
et qui repr\'esente en g\'en\'eral une vari\'et\'e trois fois
\'etendue passant par le point: $x_1^0 \dots x_4^0$.}

%%%\Fill 
%%%[[Maintien ici du en g\'en\'eral.]]

\medskip

{\bf IV)} {\em Autour du point: $y_1^0 \dots y_4^0$, une r\'egion
finie quatre fois \'etendue peut \^etre d\'elimit\'ee de telle sorte
que les conditions suivantes sont satisfaites: Si l'on fixe le
point: $y_1^0 \dots y_4^0$, alors chaque autre point r\'eel: $x_1^0
\dots x_4^0$ de la r\'egion peut encore \^etre transform\'e
contin\^ument en tous les points r\'eels: $x_1 \dots x_4$ qui
satisfont l'\'equation~\thetag{ 12}. Si l'on fixe, outre le point:
$y_1^0 \dots y_4^0$, encore un deuxi\`eme point r\'eel: $z_1^0 \dots
z_4^0$ de la r\'egion, alors chaque autre point r\'eel: $x_1^0 \dots
x_4^0$ de la r\'egion peut encore \^etre transform\'e contin\^ument en
tous les points r\'eels: $x_1 \dots x_4$ de la r\'egion qui satisfont
les deux \'equations: \thetag{ 12} et:
\def\theequation{13}\begin{equation}
W
\big(
z_1^0\dots z_4^0\,;\
x_1^0\dots x_4^0\,;\
x_1\dots x_4
\big)
=
0.
\end{equation}
Avec cela, il est \`a chaque fois suppos\'e que les deux points:
$x_1^0 \dots x_4^0$ et $x_1 \dots x_4$ sont reli\'es l'un \`a l'autre
par une s\'erie continue et irr\'eductible de points de cette sorte.}

\medskip

Nous faisons \`a nouveau remarquer express\'ement \`a ce sujet que ces
axiomes contiennent peut-\^etre certains 
\'el\'ements superflus, bien
qu'ils demandent moins que les axiomes helmholtziens. En effet,
Monsieur de Helmholtz demande en plus que lorsque trois points sont
fix\'es, chaque quatri\`eme point est encore compl\`ement libre de se
mouvoir, autant que l'autorisent les \'equations par lesquelles il est
li\'e avec les points fix\'es; \`a cet effet, c'est m\^eme son axiome
de monodromie qui se pr\'esente.

\bigskip

Soit \`a nouveau $G$ un groupe quelconque qui satisfait nos
axiomes. La vari\'et\'e passant par le point: $x_1^0 \dots x_4^0$ qui
est d\'efinie par l'\'equation~\thetag{ 12}, nous l'appelerons
naturellement encore une pseudosph\`ere de centre: $y_1^0 \dots
y_4^0$ relative \`a $G$. Notre Axiome~III exprime ensuite qu'en
g\'en\'eral, les pseudosph\`eres sont des vari\'et\'es r\'eelles trois
fois \'etendues et qu'une pseudosph\`ere ne passe jamais par son
centre.

\smallskip

On peut maintenant d\'emontrer, exactement comme dans le pr\'ec\'edent
paragraphe, que $G$ est transitif, et que deux points finiment
\'eloign\'es l'un de l'autre: $x_1 \dots x_4$ et: $y_1 \dots y_4$
ont un et un seul invariant: $\Omega ( x, y)$ relativement \`a $G$.
Il d\'ecoule alors imm\'ediatement de l\`a que l'\'equation~\thetag{
12} peut \^etre ramen\'ee \`a la forme:
\def\theequation{14}\begin{equation}
\Omega
\big(
x_1\dots x_4\,;\
y_1^0\dots y_4^0
\big)
=
\Omega
\big(
x_1^0\dots x_4^0\,;\
y_1^0\dots y_4^0
\big).
\end{equation} 
Nous ne devons donc pas nous attarder plus longtemps sur ce sujet.

\smallskip

De m\^eme, la d\'emonstration que $G$ est {\em r\'eel-primitif}\, se
pr\'esente presque exactement comme au \S~102. Parmi les
pseudosph\`eres de centre: $y_1^0 \dots y_4^0$, il y en a seulement
un nombtre discret qui ne sont pas des vari\'et\'es r\'eelles \`a
trois dimensions, mais qui sont seulement deux fois ou une fois
\'etendues, voire m\^eme se r\'eduisent \`a un point. Comme par
ailleurs aucune pseudosph\`ere ne passe par son centre, alors apr\`es
fixation du point: $y_1^0 \dots y_4^0$, il est certain qu'aucune
vari\'et\'e ponctuelle passant par ce point ne peut rester au repos.
Par cons\'equent, $G$ doit \^etre {\em r\'eel-primitif}.

%%%\Mathematiques
%%%[[Donner plus d'explications.]]

\smallskip

Enfin, exactement comme au \S~102, on peut aussi d\'emontrer que $G$
est fini, et pour pr\'eciser, qu'il a au plus six param\`etres.

\smallskip

Si en effet $P_1 \dots P_4$ et $P$ sont cinq points mutuellement en
position g\'en\'erale, alors les quatre pseudosph\`eres de centres
$P_1 \dots P_4$ passant par $P$ ne peuvent se couper qu'en $P$; car
si elles se coupaient en une vari\'et\'e $M$ passant par $P$ qui ne
consistait pas seulement en le point $P$, alors toutes les
pseudosph\`eres passant par $P$ se couperaient g\'en\'eralement en la
vari\'et\'e $M$, par suite de quoi $M$ resterait aussi en m\^eme temps
au repos apr\`es fixation de $P$, ce qui est exclu. Si maintenant
nous fixons les quatre points $P_1 \dots P_4$, alors $P$ ne peut
visiblement se mouvoir que sur l'intersection des pseudosph\`eres
centr\'ees en $P_1 \dots P_4$ qui passent par $P$, et puisque cette
intersection consiste seulement en le point $P$ lui-m\^eme, $P$ doit
rester au repos, et parce que $P$ est un point en position
g\'en\'erale, chaque point de l'espace reste en m\^eme temps
g\'en\'eralement au repos. Sous les hypoth\`eses pos\'ees, la fixation
de $P_1 \dots P_4$ n\'ecessite au plus: $4 + 3 + 2 + 1 = 10$
conditions, donc on obtient que $G$ est fini, et pour pr\'eciser,
qu'il a {\em au plus dix param\`etres}.

\smallskip

Choisissons maintenant, parmi les pseudosph\`eres de centre $P_1$, une
quelconque d'entre elles en position g\'en\'erale, et appelons-la
$K_1$. Si $P_2$ est un point quelconque de $K_1$, il y a $\infty^1$
pseudosph\`eres de centre $P_2$, mais dont aucune ne co\"{\i}ncide
avec $K_1$, car une pseudosph\`ere ne contient jamais son centre.

Si nous fixons $P_1$, alors $P_2$ peut se mouvoir d'une mani\`ere
compl\`etement libre sur la pseudosph\`ere $K_1$, qui, sous les
hypoth\`eses pos\'ees, est certainement une vari\'et\'e r\'eelle trois
fois \'etendue. Si, hormis $P_1$, nous fixons aussi encore $P_2$,
alors chaque autre point $P_3$ de $K_1$ peut seulement se mouvoir sur
la vari\'et\'e $M'$ qui est d\'ecoup\'ee dans $K_1$ par la
pseudosph\`ere de centre $P_2$ qui passe par $P_3$;
et d'apr\`es l'Axiome~IV, $P_3$ peut se mouvoir d'une mani\`ere
compl\`etement libre sur cette vari\'et\'e. Nous pouvons ajouter que,
sous les hypoth\`eses pos\'ees, $M'$ est en g\'en\'eral une
vari\'et\'e r\'eelle deux fois \'etendue.

Si nous d\'eterminons les points de $K_1$ par trois
coordonn\'ees: $\mathfrak{ x}_1, \mathfrak{ x}_2, \mathfrak{ x}_3$ et
si nous appelons: $\mathfrak{ y}_1^0, \mathfrak{ y}_2^0, \mathfrak{
y}_3^0$ les coordonn\'ees de $P_2$, et: $\mathfrak{ x}_1^0,
\mathfrak{ x}_2^0, \mathfrak{ x}_3^0$ celles de $P_3$, nous pouvons
aussi exprimer tout cela comme suit: Si $P_1$ est fix\'e, alors les
points: $\mathfrak{ x}_1, \mathfrak{ x}_2, \mathfrak{ x}_3$ de la
pseudosph\`ere $K_1$ sont transform\'es de telle sorte qu'apr\`es
fixation d'un point r\'eel: $\mathfrak{ y}_1^0, \mathfrak{ y}_2^0,
\mathfrak{ y}_3^0$ de $K_1$, 
chaque autre point r\'eel $\mathfrak{ x}_1^0,
\mathfrak{ x}_2^0, \mathfrak{ x}_3^0$ de $K_1$ peut
encore \^etre envoy\'e sur
tous les points r\'eels: $\mathfrak{ x}_1, \mathfrak{ x}_2,
\mathfrak{ x}_3$ de $K_1$ qui satisfont une \'equation de la forme:
\[
\mathfrak{W}
\big(\mathfrak{ y}_1^0, \mathfrak{ y}_2^0,
\mathfrak{ y}_3^0\,;\
\mathfrak{ x}_1^0, \mathfrak{ x}_2^0, \mathfrak{ x}_3^0\,;\
\mathfrak{ x}_1, \mathfrak{ x}_2, \mathfrak{ x}_3
\big)
=
0;
\]
en g\'en\'eral, cette \'equation d\'etermine une vari\'et\'e r\'eelle
deux fois \'etendue se trouvant sur $K_1$. Et maintenant, lorsqu'on
fixe $P_1$, les points de la vari\'et\'e trois fois \'etendue $K_1$
sont visiblement transform\'es par un groupe $G_1$; d'apr\`es ce qui
a justement \'et\'e dit, il s'ensuit de plus que ce groupe satisfait
tous les axiomes pos\'es dans le paragraphe pr\'ec\'edent, par suite
de quoi nous pouvons conclure que $G_1$ a six param\`etres et qu'il
peut \^etre transform\'e en le groupe des mouvements euclidiens, ou en
les deux groupes de mouvements non-euclidiens de $R_3$, {\em via}\, une
transformation ponctuelle r\'eelle en les variables: $\mathfrak{
x}_1, \mathfrak{ x}_2, \mathfrak{ x}_3$. D'autre part, nous avons vu
plus haut que $G$ est transitif et qu'il a au plus dix param\`etres,
donc le sous-groupe de $G$ par lequel $P_1$ reste invariant n'a
certainement pas plus que six param\`etres. Par cons\'equent, on
obtient que $G$ poss\`ede exactement dix param\`etres, et qu'apr\`es
fixation de $P_1$, les points de $R_4$ sont transform\'es par un
groupe \`a six param\`etres; en m\^eme temps, les points de chaque
pseudosph\`ere de centre $P_1$ qui est g\'en\'eralement situ\'ee sont
transform\'es par un groupe isomorphe-holo\'edrique \`a six
param\`etres, qui est semblable, par une transformation ponctuelle
r\'eelle de $R_3$, soit au groupe des mouvements euclidiens, soit \`a
l'un des deux groupes de mouvements non-euclidiens.

%%%\Mathematiques
%%%[[flots transversaux aux pseudosph\`eres~??]]

Si nous fixons $P_1$ et $P_2$, alors les points de $K_1$ sont encore
transform\'es par l'action d'un groupe \`a trois param\`etres, et il
en va de m\^eme pour les points de chaque autre pseudosph\`ere de
centre $P_1$ qui est g\'en\'eralement situ\'ee. Les groupes r\'eels
\`a trois param\`etres en question sont naturellement
isomorphes-holo\'edriques l'un avec l'autre. Mais si un sous-groupe
r\'eel \`a trois param\`etres des mouvements euclidiens ou
non-euclidiens de $R_3$ est isomorphe-holo\'edrique avec un
sous-groupe \`a trois param\`etres qui laisse invariant un point
r\'eel, alors il est manifestement lui-m\^eme le sous-groupe \`a trois
param\`etres qui laisse invariant un certain point r\'eel (\cf\,
p.~385 et Chap.~10). Par cons\'equent, lorsque $P_1$ et $P_2$ sont
fix\'es, sur chaque pseudosph\`ere de centre $P_1$ qui est
g\'en\'eralement situ\'ee, un certain point r\'eel doit en g\'en\'eral
rester au repos, et pour pr\'eciser, les points en question doivent
visiblement constituer une courbe qui passe par $P_2$. La sym\'etrie
montre qu'apr\`es fixation de $P_1$ et de $P_2$, il passe aussi en
m\^eme temps par $P_1$ une courbe continue dont les points restent
enti\`erement au repos.

Si donc nous fixons $P_1$ et $P_2$, alors il passe par $P_1$ ainsi que
par $P_2$ une courbe continue dont les points restent au repos. On
peut s'imaginer que ces deux courbes co\"{\i}ncident, mais il est
aussi possible qu'elles diff\`erent l'une de l'autre. Nous voulons
laisser en suspens la question de savoir si le deuxi\`eme cas peut
r\'eellement se produire. Afin de pouvoir nous exprimer d'une
mani\`ere commode, nous voulons appeler {\sl pseudodroite}\,
l'ensemble des deux courbes. Ainsi nous pouvons dire: Pour chaque
paire de points $P_1$ et $P_2$ de $R_4$, une pseudodroite est
d\'etermin\'ee qui passe par $P_1$ et par $P_2$; si on fixe $P_1$ et
$P_2$, alors les points de cette pseudodroite restent enti\`erement au
repos.

Il est clair que chaque point de la pseudosph\`ere $K_1$ d\'etermine
une pseudodroite passant par $P_1$ et que les pseudodroites ainsi
d\'etermin\'ees passent en g\'en\'eral toutes par $P_1$; car si nous
fixons en m\^eme temps que $P_1$ un point $P$ en position g\'en\'erale
qui ne se trouve pas sur $K_1$, alors un point $P'$ sur $K_1$ reste
aussi au repos et la pseudodroite d\'etermin\'ee par $P_1$ et par $P$
co\"{\i}ncide donc avec celle qui est d\'etermin\'ee par $P_1$ et par
$P'$. Par cons\'equent, exactement $\infty^3$ pseudodroites passant par
$P_1$ sont en g\'en\'eral d\'etermin\'ees, mais puisque chacune
de ces pseudodroites consiste \'eventuellement en deux courbes, dont
l'une seulement passe par $P_1$ lui-m\^eme, il n'est pas certain
depuis le d\'ebut que la famille des courbes passant par 
$P_1$, qui est d\'etermin\'ee par ces $\infty^3$ pseudodroites, 
est constitu\'ee de $\infty^3$ courbes diff\'erentes.

Comme $\infty^3$ pseudodroites diff\'erentes passant par $P$ sont
d\'etermin\'ees et comme, d'autre part, $\infty^3$ \'el\'ements
lin\'eaires diff\'erents passent aussi par $P_1$, il y a parmi ces
\'el\'ements lin\'eaires un certain nombre, \ie $\infty^m$ ($0
\leqslant m \leqslant 3$) qui sont diff\'erents, et qui sont
constitu\'es de telle sorte que par chacun d'entre eux passent
certaines de ces pseudodroites-l\`a, alors qu'aucune telle
pseudodroite ne passe par les \'el\'ements lin\'eaires restants. Nous
allons d\'emontrer que le nombre $m$ en question est simplement \'egal
\`a trois.

En effet, si $P_1$ est fix\'e, alors les pseudodroites passant par
$P_1$ sont \'echangeables l'une avec l'autre,
%%%\Mathematiques
%%%[[Il y a une notion de commutation au sens des groupes, mais je ne la devine pas.]]
et donc la vari\'et\'e des $\infty^m$ \'el\'ements lin\'eaires
d\'efinie \`a l'instant reste au repos. Comme de plus \`a chaque point de
la pseudosph\`ere $K_1$ est associ\'ee l'une 
des pseudodroites passant par
$P_1$ et comme, apr\`es fixation de $P_1$, chaque point de $K_1$ peut
\^etre encore transform\'e en chaque autre, alors, apr\`es fixation de
$P_1$, chacune des pseudodroites d\'etermin\'ees
passant par $P_1$ peut aussi \^etre
transform\'ee en chaque autre et par cons\'equent aussi, chacun des
ces $\infty^m$ \'el\'ements lin\'eaires-l\`a peut \^etre transform\'e
en chaque autre. Il d\'ecoule de l\`a que par chacun de ces
$\infty^m$ \'el\'ements lin\'eaires passent $\infty^{ 3 - m}$
pseudodroites diff\'erentes, et que la pseudosph\`ere $K_1$ se
d\'ecompose en $\infty^m$ vari\'et\'es r\'eelles $(3 - m)$ fois
\'etendues, de telle sorte que chacun de nos $\infty^m$ \'el\'ements
lin\'eaires est associ\'e \`a l'une de ces vari\'et\'es; apr\`es
fixation de $P_1$, ces $\infty^m$ vari\'et\'es $(3 - m)$ fois \'etendues
sont \'echangeables l'une avec l'autre. 
%%%\Mathematiques
%%%[[M\^eme probl\`eme.]]
Si maintenant on avait $m = 1$ ou $= 2$, alors, apr\`es fixation de
$P_1$, les points de $K_1$ se transformeraient de mani\`ere
r\'eelle-imprimitive, ce qui n'est \'evidemment pas le cas. Si d'un
autre c\^ot\'e, on avait $m = 0$, alors un nombre discret
d'\'el\'ements lin\'eaires r\'eels passant par $P_1$ resterait au
repos en m\^eme temps que $P_1$, et notre groupe $G$ \`a dix
param\`etres serait donc r\'eel-imprimitif, alors qu'il doit cependant
\^etre r\'eel-primitif.

Par cons\'equent, le nombre $m$ d\'efini plus haut est r\'eellement
\'egal \`a $3$, en cons\'equence de quoi il passe en 
g\'en\'eral par $P_1$
une pseudodroite dirig\'ee par chaque \'el\'ement lin\'eaire r\'eel.

Ainsi, on a d\'emontr\'e qu'il existe, entre les $\infty^3$
\'el\'ements lin\'eaires en $P_1$ et les $\infty^3$ points r\'eels de
la pseudosph\`ere $K_1$, une relation qui en tout cas est univoque et
r\'eversible \deutsch{\small\sf eindeutig umkehrbar} \`a l'int\'erieur
d'une certaine r\'egion
et qui reste maintenue par
toutes les transformations de $G$ qui laissent invariant le point
$P_1$. Si nous nous rappelons
alors qu'apr\`es fixation de $P_1$, les 
$\infty^3$ \'el\'ements lin\'eaires passant par $P_1$ sont 
transform\'es par l'action d'un groupe projectif r\'eel $\mathfrak{
g}$, nous reconnaissons donc imm\'ediatement que ce groupe $\mathfrak{
g}$ est semblable au groupe $G_1$, par lequel les points de $K_1$ sont
transform\'es en m\^eme temps, et pour pr\'eciser, que $\mathfrak{ g}$
est semblable \`a $G_1$ {\em via}\, une transformation ponctuelle
r\'eelle de $R_3$. Mais comme $G_1$ \'etait de son c\^ot\'e
semblable, {\em via}\, une transformation ponctuelle r\'eelle
de $R_3$, soit au
groupe des mouvements euclidiens, soit \`a l'un des deux groupes de
mouvements non-euclidiens de cet espace, il en d\'ecoule donc, en
tenant compte du Th\'eor\`eme~19, p.~292, que $\mathfrak{ g}$ est \`a
six param\`etres et qu'il peut \^etre transform\'e au moyen d'une
transformation projective r\'eelle en l'un des trois groupes de
mouvements mentionn\'es.

\medskip

{\em Ainsi, notre groupe $G$ \`a dix param\`etres est constitu\'e de
telle sorte qu'apr\`es fixation d'un point r\'eel en position
g\'en\'erale, les $\infty^3$ \'el\'ements lin\'eaires r\'eels passant
par ce point sont transform\'es par l'action d'un groupe \`a six
param\`etres, et pour pr\'eciser, d'un groupe euclidien ou
non-euclidien.}

\medskip

Si les \'el\'ements lin\'eaires en question sont transform\'es par
l'action d'un groupe euclidien, alors, en m\^eme temps que chaque
point r\'eel en position g\'en\'erale, une botte 
\deutsch{\small\sf B\"undel}
r\'eelle de
$\infty^2$ \'el\'ements lin\'eaires passant par un tel point reste
invariante, 
%%%\Mathematiques
%%%[[What~??]]
donc $G$ laisse invariante une \'equation de Pfaff
r\'eelle:
\def\theequation{15}\begin{equation}
\sum_1^4\,\alpha_\nu(x_1\dots x_4)\,dx_\nu
=
0,
\end{equation}
qui ne doit naturellement pas \^etre int\'egrable, parce que sinon $G$
serait r\'eel-imprimitif. Mais maintenant, il est connu d'apr\`es la
th\'eorie du probl\`eme de Pfaff,
%%%\Mathematiques
%%%[[What~??]]
qu'\`a toute \'equation de Pfaff non int\'egrable \`a quatre variables
est associ\'ee un syst\`eme invariant simultan\'e. Ce syst\`eme
simultan\'e, qui est \'egalement r\'eel pour l'\'equation
r\'eelle~\thetag{ 15}, devrait naturellement rester invariant par $G$,
donc dans ce cas $G$ devrait \^etre r\'eel-imprimitif;
%%%\Mathematiques
%%%[[What~??]]
par cons\'equent, le cas o\`u les $\infty^3$ \'el\'ements lin\'eaires
passant par un point r\'eel fix\'e sont transform\'es de mani\`ere
euclidienne n'entre g\'en\'eralement pas en ligne de compte.

D'un autre c\^ot\'e, si les \'el\'ements lin\'eaires en question se
transforment de mani\`ere non-euclidienne, alors \`a chaque point
r\'eel en position g\'en\'erale et associ\'ee une conique du second
degr\'e r\'eelle ou imaginaire constitu\'ee d'\'el\'ements
lin\'eaires, qui, par un choix appropri\'e des variables, re\c coit ou
bien la forme: 
\def\theequation{16}\begin{equation}
dx_1^2+\cdots+dx_4^2
=
0,
\end{equation}
ou bien la forme: 
\def\theequation{17}\begin{equation}
dx_1^2+dx_2^2+dx_3^2-dx_4^2
=
0.
\end{equation}
Maintenant, puisque $G$ a dix param\`etres et que les \'el\'ements
lin\'eaires passant par chaque point r\'eel fix\'e se transforment par
l'action d'un groupe \`a six param\`etres, on obtient gr\^ace aux
d\'eveloppements des pages~385 sq., que $G$ peut \^etre transform\'e
{\em via}\, une transformation ponctuelle r\'eelle de $R_4$ soit en le
groupe des mouvements euclidiens, soit en l'un des deux groupes de
mouvements non-euclidiens\,\,---\,\,ces cas sont possibles, lorsque la
forme~\thetag{ 16} se produit\,\,---, ou bien qu'il est semblable, 
{\em via}\, une transformation ponctuelle r\'eelle de $R_4$ 
soit au groupe projectif d'une des deux vari\'et\'es: 
\[
x_1^2+x_2^2+x_3^2-x_4^2
=
\pm 1,
\]
soit au groupe: 
\[
p_1\dots p_4,\quad\quad
x_\mu\,p_\nu-x_\nu\,p_\mu,\quad\quad
x_\mu\,p_4+x_4\,p_\mu
\quad\quad\quad
{\scriptstyle{(\mu,\,\nu\,=\,1,\,2,\,3)}}
\]
---\,\,ces cas correspondent \`a la forme~\thetag{17}.
%%%\Mathematiques
%%%[[\`A comprendre.]]

Mais les trois derniers groupes
n'entrent pas du tout en ligne de compte, 
car pour chacun d'entre eux, il y a visiblement, parmi
les pseudosph\`eres ayant un centre donn\'e, toujours 
une pseudosph\`ere qui
passe par son centre. 
%%%\Mathematiques
%%%[[What??]]
{\em Par cons\'equent, il ne reste plus que les mouvements euclidiens
et non-euclidiens, et on a donc d\'emontr\'e que ces trois familles de
mouvements de $R_4$ sont compl\`etement caract\'eris\'ees par les
axiomes indiqu\'es \`a la page~\pageref{515}.}

\bigskip

\`A pr\'esent, nous voulons indiquer de quelle mani\`ere ces
consid\'erations se r\'ealisent pour les espaces de dimension
sup\'erieure \`a quatre.

Dans $R_n$ ($n > 2$) nous demandons pareillement que l'espace soit une
vari\'et\'e num\'erique et que les mouvements forment un groupe
continu engendr\'e par des transformations infinit\'esimales. Si un
point r\'eel: $y_1^0 \dots y_n^0$ en position g\'en\'erale est
fix\'e, alors tout autre point r\'eel: $x_1^0 \dots x_n^0$ peut
\^etre envoy\'e encore seulement sur les points r\'eels: $x_1 \dots
x_n$ qui satisfont une certaine \'equation:
\[
W
\big(
y_1^0\dots y_n^0\,;\
x_1^0\dots x_n^0\,;\
x_1\dots x_n
\big)
=
0.
\]
Avec cela, nous supposons que cette \'equation repr\'esente en
g\'en\'eral une vari\'et\'e r\'eelle $(n-1)$ fois \'etendue, et
qu'elle n'est pas satisfaite pour: $x_1 = y_1^0, \dots, x_n = y_n^0$.

La vari\'et\'e r\'eelle passant par le point: $x_1^0 \dots x_n^0$ qui
est d\'etermin\'ee par l'\'equation: $W = 0$, nous l'appelons bien
s\^ur une pseudosph\`ere de centre: $y_1^0 \dots y_n^0$. Notre
derni\`ere exigence exprime alors visiblement qu'une pseudosph\`ere ne
peut jamais passer par son centre.

Enfin, nous demandons encore que dans $R_n$, puisse \^etre
d\'elimit\'ee une r\'egion finie $n$ fois \'etendue, \`a l'int\'erieur
de laquelle les exigences suivantes soient satisfaites: Si l'on fixe
un point $P_1$ de la r\'egion, alors tout autre point de la r\'egion
doit pouvoir se mouvoir d'une mani\`ere enti\`erement libre sur la
pseudosph\`ere de centre $P_1$ passant par lui. Si l'on fixe $q$
points $P_1 \dots P_q$ de la r\'egion qui sont mutuellement en
position g\'en\'erale, alors, tant que $q$ est $< n-1$, chaque autre
point en position g\'en\'erale appartenant \`a la r\'egion doit
pouvoir se mouvoir d'une mani\`ere enti\`erement libre sur la
vari\'et\'e passant par $P$ qui est l'intersection des $q$
pseudosph\`eres de centres $P_1 \dots P_q$.

Nous affirmons que ces exigences suffisent pour caract\'eriser les
mouvements euclidiens et non-euclidiens. Mais nous voulons seulement
indiquer la d\'emonstration de cette assertion.

Nous supposons que notre assertion est d\'ej\`a d\'emontr\'ee pour
l'espace \`a $n-1$ dimensions et nous d\'emontrons, en prenant cette
hypoth\`ese pour base, que notre assertion est vraie aussi pour pour
l'espace \`a $n$ dimensions. Comme nous avons l'avons d\'ej\`a
d\'emontr\'ee dans les cas $n = 3$ et $n = 4$, sa validit\'e en
g\'en\'eral sera alors \'etablie.

Mais pour d\'emontrer que notre assertion est vraie dans $R_n$,
aussit\^ot qu'elle l'est dans $R_{ n - 1}$, nous proc\'edons comme
dans le cas $n = 4$. Nous montrons que chaque groupe $G$ qui
satisfait nos axiomes est transitif et que relativement \`a lui, deux
points finiment \'eloign\'es l'un de l'autre ont un et un seul
invariant; de plus, chaque groupe tel est fini, r\'eel-primitif, et
il contient au plus $\frac{ 1}{ 2}\, n ( n+1)$ param\`etres. Si \`a
pr\'esent nous fixons un point r\'eel $P_1$ en position g\'en\'erale,
alors on obtient facilement que les points de chaque pseudosph\`ere de
centre $P_1$ g\'en\'eralement situ\'ee sont transform\'es par un
groupe qui satisfait tous les axiomes dans $R_{ n - 1}$, et par
cons\'equent, gr\^ace \`a l'hypoth\`ese prise pour base, ce groupe
peut \^etre transform\'e, {\em via}\, une transformation ponctuelle de
$R_{ n - 1}$, soit en le groupe des mouvements euclidiens, soit en
l'un des deux groupes de mouvements non-euclidiens de cet espace.

Il d\'ecoule de l\`a tout d'abord que $G$ contient pr\'ecis\'ement
$\frac{ 1}{ 2}\, n ( n+1)$ param\`etres. De plus, exactement comme
dans le cas $n = 4$, on peut d\'emontrer que deux points de $R_n$
d\'eterminent une pseudodroite. Pour cela, on peut se baser sur les
recherches de Monsieur Werner, qui, \`a l'initiative de Lie, a
d\'etermin\'e les plus grands sous-groupes qui sont contenus dans le
groupe projectif d'une vari\'et\'e non-d\'eg\'en\'er\'ee du second
degr\'e dans l'espace \`a $n > 5$ dimensions.

%%%\Mathematiques
%%%[[Trouver l'article.]]

Apr\`es que l'existence des pseudodroites est connue, on montre \`a
nouveau, exactement comme dans le cas $n = 4$, que les $\infty^{ n -
1}$ \'el\'ements lin\'eaires r\'eels passant par chaque point r\'eel
fix\'e en position g\'en\'erale sont transform\'es par $G$ de
mani\`ere euclidienne ou non-euclidienne. Enfin, on d\'emontre
facilement qu'ils ne peuvent certainement pas \^etre transform\'es de
mani\`ere euclidienne, et gr\^ace \`a la consid\'eration des cas
restants, on obtient que $G$ peut \^etre transform\'e, {\em via}\, une
transformation ponctuelle r\'eelle de $R_n$, soit en le groupe des
mouvements euclidiens, soit en l'un des deux groupes de mouvements
non-euclidiens de $R_n$.

\bigskip

Nous voulons terminer le pr\'esent chapitre en rectifiant une petite
erreur qui s'est gliss\'ee aux pages~\pageref{446} sq. \`A cet
endroit-l\`a, nous avons en effet dit que le troisi\`eme axiome de
Monsieur de Helmholtz est constitu\'e de deux parties, dont la seconde
contient certaines exigences qui ne d\'ecoulent pas de celles qui sont
pos\'ees dans la premi\`ere, bien que, d'apr\`es sa version, cette
seconde partie semble contenir seulement des cons\'equences de la
premi\`ere. Mais en fait, les choses se passent tout autrement. {\em
Si en effet on comprend ce que nous avons appel\'e \`a ce moment-l\`a
la premi\`ere partie du troisi\`eme axiome helmholtzien} de telle
sorte que, quand un point $P_1$ est fix\'e, chaque autre point doit
pouvoir se mouvoir d'une mani\`ere enti\`erement libre sur la
pseudosph\`ere de centre $P_1$ qui passe par lui, et si l'on en tire
la conclusion qu'une pseudosph\`ere ne peut jamais passer par son
centre, alors \`a vrai dire, ce que nous avons appel\'e la seconde
partie du troisi\`eme axiome helmholtzien ne contient rien qui ne
d\'ecoule pas d\'ej\`a de la premi\`ere partie; cela est montr\'e par
des consid\'erations similaires \`a celles que nous avons
d\'evelopp\'ees dans le pr\'esent chapitre, et par exemple, dans$R_3$, afin de d\'emontrer que trois pseudosph\`eres dont les centres
sont mutuellement en position g\'en\'erale se coupent g\'en\'eralement
en un seul point. Et maintenant, puisqu'au cours de notre recherche
sur les axiomes helmholtziens, nous avons toujours interpr\'et\'e cela
de telle sorte qu'une pseudosph\`ere ne doit pas passer par son centre
(\cf~p.~\pageref{465}), notre division en deux parties du
troisi\`eme axiome effectu\'ee aux pages~\pageref{446} sq. est
erron\'ee. \`A vrai dire, il reste toujours des chances que Monsieur
de Helmholtz ait propos\'e les cons\'equences qu'il tire de son
troisi\`eme axiome, sans les justifier, bien que leur justification ne
soit pas du tout aussi simple que cela.

%%%%%%%%%%%%%%%%%%%%%%%%%%%%%%%%%%%%%%%%%%%%%%%%%%%%%%%%%%%%%%%%%%%%%

\newpage

\renewcommand{\chaptermark}{}

\HEAD{Bibliographie}{
Bibliographie}

\renewcommand{\bibname}
{\large\bf Bibliographie\hspace{-2.953cm}\large\bf}
\label{-Bibliographie}

\newpage

$\:$

\bigskip\bigskip\bigskip

\begin{center}
{\bf
Quatrième de couverture:
}\end{center}
\label{quatrieme-de-couverture}

\bigskip

\begin{center}
{\large\sf LE PROBL\`EME DE L'ESPACE}

\medskip

{\sf Sophus Lie, Friedrich Engel 

\medskip

et le problème de Riemann-Helmholtz}

\end{center}

\bigskip

Est-il possible de caractériser l'espace euclidien tridimensionnel qui
s'offre si immédiatement à l'intuition physique au moyen d'axiomes
mathématiques simples et naturels? Plus généralement, est-il possible
de caractériser les espaces de Bolyai-Lobatchevski\u{\i} à courbure
constante négative, ainsi que les espaces de Riemann à courbure
constante positive, à l'exclusion de toute autre géométrie contraire à
une intuition directe?

\`A une époque (1830--50) où l'émergence nécessaire des géométries
dites non-euclidiennes devenait incontestable, c'est Riemann qui a
soulevé cette question profonde et difficile dans son discours
d'habilitation (1854), sans chercher, toutefois, à la résoudre
complètement. Helmholtz (1868) l'interprétera en conceptualisant le
mouvement des corps dans l'espace et il tentera d'établir
rigoureusement que le caractère métrique et localement homogène d'un
espace se déduit d'axiomes de mobilité maximale pour des corps
rigides.

Mais il fallut attendre les travaux de Sophus Lie, et notamment la
{\em Theorie der Transformationsgruppen} 
(2100 pages, 1888--93) écrite en
collaboration avec Friedrich Engel, pour qu'une solution complète et
rigoureuse soit apportée à ce fascinant problème, à la fois au plan
local et au plan global. L'introduction historique, philosophique et
mathématique ainsi que la traduction que nous proposons ici aspirent à
faire connaître un aspect de l'{\oe}uvre monumentale de Sophus Lie qui
demeure essentiellement peu évoqué au sein de la philosophie
traditionnelle de la géométrie.

\bigskip\bigskip\noindent
{\bf Joël Merker}, agrégé de mathématiques et de philosophie,
spécialiste d'analyse et de géométrie à plusieurs variables réelles ou
complexes, chercheur au CNRS - Département de Mathématiques et
Applications, \'Ecole Normale Supérieure.

%%%%%%%%%%%%%%%%%%%%%%%%%%%%%%%%%%%%%%%%%%%%%%%%%%%%%%%%%%%%%%%%%%%%%

\vfill\end{document}